\documentclass[a4paper, 10pt, leqno,twoside]{renereport}

\usepackage{amsfonts}
\usepackage{amsmath}
\usepackage{amsthm}
\usepackage{amssymb}
\usepackage{amsthm}
\usepackage{changebar}
\usepackage[dvips]{graphicx}
\usepackage{fancyhdr}
\usepackage{geometry}
\usepackage[latin1]{inputenc}
\usepackage{makeidx}
\usepackage{mathtools}
\usepackage{mathrsfs}
\usepackage{nicefrac}
\usepackage[normalem]{ulem}
\usepackage{pifont}
\usepackage{pstricks}
\usepackage{stmaryrd}
\usepackage{t1enc}
\usepackage[T1]{fontenc}
\usepackage{times}
\usepackage{trsym}
\usepackage{wasysym}
\usepackage{yfonts} 
\usepackage[arrow, matrix, curve]{xy}
\renewcommand{\chaptermark}[1]{ \markboth {}{#1}}
\renewcommand{\sectionmark}[1]{\markboth{\uppercase {#1}}{}}
\fancypagestyle{plain}{}																					% damit auch "plain" Seiten fancy werden
\newcommand{\A}{{\mathbb A}}
\newcommand{\Z}{{\mathbb Z}}

\newcommand{\Q}{{\mathbb Q}}
\newcommand{\R}{{\mathbb R}}
\newcommand{\N}{{\mathbb N}}
\newcommand{\C}{{\mathbb C}}

\renewcommand{\H}{\mathbb{H}}

\newcommand{\id}{\operatorname{id}}

\newcommand{\coker}{\operatorname{coker}}
\newcommand{\Hom}{\operatorname{Hom}}

\newcommand{\Spec}{{\operatorname{Spec}\,}}

\newcommand{\Yr}{Y^{\natural}}
\newcommand{\pn}{\pi^{\natural}_n}

\newcommand{\var}{\mathrm{pol}_{{\mathrm{dR}},D^2 \cdot 1_{\{ \epsilon \}}-1_{E_{[D]}}}}
\newcommand{\varn}{\mathrm{pol}^n_{{\mathrm{dR}},D^2 \cdot 1_{\{ \epsilon \}}-1_{E_{[D]}}}}
\newcommand{\varD}{\mathrm{pol}_{{\mathrm{dR}},D^2 \cdot 1_{\{ \epsilon \}}-1_{E_{[D]}}}}
\newcommand{\varDn}{\mathrm{pol}^n_{{\mathrm{dR}},D^2 \cdot 1_{\{ \epsilon \}}-1_{E_{[D]}}}}
\newcommand{\varDk}{\mathrm{pol}^k_{{\mathrm{dR}},D^2 \cdot 1_{\{ \epsilon \}}-1_{E_{[D]}}}}

\newcommand{\Ext}{\mathrm{Ext}}
\newcommand{\Mod}{\mathrm{Mod}}
\newcommand{\Lie}{\mathrm{Lie}}
\newtheoremstyle{break1}% name  
{9pt}%      Space above, empty = 'usual value'  
{9pt}%      Space below  
{\itshape}% Body font  
{}%         Indent amount (empty = no indent, \parindent = para indent)  
{\bfseries}% Thm head font  
{}%        Punctuation after thm head  
{\newline}% Space after thm head: \newline = linebreak  
{}%         Thm head spec  
\theoremstyle{break1} 
\newtheorem{satz}{Satz}[section]
\newtheorem{lemma}[satz]{Lemma}
\newtheorem{Lemma}[satz]{Lemma} 
\newtheorem{theorem}[satz]{Theorem} 
 
\newtheorem{proposition}[satz]{Proposition}
\newtheorem{Korollar}[satz]{Corollary}

\newtheorem{corollary}[satz]{Corollary} 
\newtheorem{Proposition}[satz]{Proposition}
\newtheorem{Proposition/Definition} [satz] {Proposition/Definition}
\newtheorem{Lemma/Definition}[satz]{Lemma/Definition}
\theoremstyle{definition}

\newtheoremstyle{break2}% name  
{9pt}%      Space above, empty = 'usual value'  
{9pt}%      Space below  
{\normalfont}% Body font  
{}%         Indent amount (empty = no indent, \parindent = para indent)  
{\bfseries}% Thm head font  
{}%        Punctuation after thm head  
{\newline}% Space after thm head: \newline = linebreak  
{}%         Thm head spec  
\theoremstyle{break2}
\newtheorem{definition}[satz]{Definition}

\newtheorem{Remark}[satz]{Remark}
\newtheorem{Bemerkung/Definition}[satz]{Bemerkung/Definition}

\newtheorem{guess}[satz]{Guess}
\newtheorem{remark}[satz]{Remark}

\newtheorem{definition/notation}[satz]{Definition/Notation}
\newtheorem{remark/notation}[satz]{Remark/Notation}
\date{}
\title{\Huge{\textsf{The de Rham realization of the elliptic polylogarithm in families}}\\ \bigskip}
\author{DISSERTATION ZUR ERLANGUNG DES DOKTORGRADES\\
DER NATURWISSENSCHAFTEN (DR. RER. NAT.)\\
DER FAKULTÄT FÜR MATHEMATIK\\
\vspace{10mm}
DER UNIVERSITÄT REGENSBURG\\
vorgelegt von\\
René Achim Scheider\\
aus Deggendorf\\
im Jahr 2014\\
}

\begin{document}

\maketitle

\vspace*{\fill}
Promotionsgesuch eingereicht am: 21. Januar 2014\\
\newline
Die Arbeit wurde angeleitet von: Prof. Dr. Guido Kings\\
\newline
Prüfungsausschuss:\\
Prof. Dr. Helmut Abels (Vorsitzender)\\
Prof. Dr. Guido Kings (1. Gutachter)\\
Prof. Dr. Kenichi Bannai, Keio University, Japan (2. Gutachter)\\
Prof. Dr. Klaus Künnemann\\
Prof. Dr. Ulrich Bunke (Ersatzprüfer)
\thispagestyle{empty}
\renewcommand{\contentsname}{Contents}
\tableofcontents

\clearpage
\mbox{}
\thispagestyle{empty}
\clearpage

\
\vspace{1cm}
\section*{\LARGE Introduction}
\markboth{\uppercase{Introduction}}{}
\addcontentsline{toc}{chapter}{Introduction}
Polylogarithms in their various manifestations provide a key instrument for modern arithmetic geometry's quest to determine special values of $L$-functions resp. to investigate these within the context of algebraic $K$-theory and motivic cohomology.\\
Illustration of this principle is best woven into a general review of some major development steps that led to the concept of the polylogarithm as it appears in the present work.

\subsubsection{Classical polylogarithm}

As a starting point, let us consider the classical polylogarithm functions which are defined at first for $|z|<1$ via the power series
\[\mathrm{Li}_m(z)=\sum_{n=1}^{\infty}\frac{z^n}{n^m} \quad (m\geq 1)\]
and then extended to multivalued holomorphic functions on $\mathbb P^1(\C) \backslash \{ 0,1,\infty \}$ using the expressions
\[\mathrm{Li}_{m+1}(z)=\int_0^z\mathrm{Li}_m(t)\frac{\mathrm{d}t}{t} \ \ \ (m\geq 1), \quad \mathrm{Li}_1(z)=\int_0^z\frac{\mathrm{d}t}{1-t}.\vspace{0.30cm}\]
Let now $F$ be a number field of degree $n=r_1+2r_2$ and with discriminant $d_F$; write $\zeta_F(s)$ for its Dedekind zeta function. It is a highly nontrivial task and of fundamental arithmetic interest to find information about the special values $\zeta_F(m)$ for integers $m\geq 2$. The above polylogarithm functions have a striking relevance for this problem:\\
In the 1980's Zagier \cite{Za1} used volume computations from $3$-dimensional hyperbolic geometry to prove that the number $\zeta_F(2)\pi^{-2(r_1+r_2)}|d_F|^{\frac{1}{2}}$ is a rational linear combination of products of the Bloch-Wigner function $D$ evaluated at algebraic arguments; this last function can be imagined as a single valued version $D(z):\mathbb P^1(\C)\rightarrow \R$ of the dilogarithm $\mathrm{Li}_2(z)$.\\
Subsequently, a refined formulation of this result and a conjectural generalization for arbitrary $m\geq 2$ was given in the conceptual framework of $K$-theory (cf. \cite{Za3}, §1 and §8): Again, one introduces a single valued version $P_m(z): \mathbb P^1(\C) \rightarrow \R$ of $\mathrm{Li}_m(z)$. Functional equations satisfied by $P_m$ model the definition of a subquotient $\mathcal B_m(F)$ of $\Z[F^*]$ on which $P_m$ (and embeddings of $F$ into $\C$) gives rise to a map into euclidean space of the same dimension, say $j(m)$, as the target of the Borel regulator map on $K_{2m-1}(F)$. Zagier's conjecture then claims that (up to torsion) there is a canonical isomorphism between $K_{2m-1}(F)$ and $\mathcal B_m(F)$ under which these two maps coincide. Using Borel's theorem on the covolume of the regulator lattice this would imply $\Q^*$- equivalence of $\zeta_F(m)$ with $\pi^{mj(m)}|d_F|^{-\frac{1}{2}}$- times the determinant of a matrix with entries given by $P_m$ evaluated at $F$-algebraic arguments.\\
\newline
In \cite{Be-De} (part of) Zagier's conjecture was reinterpreted in motivic and Hodge-theoretic formalism. Important for the approach to polylogarithms is the observation that the monodromy and differential equations of the $\mathrm{Li}_m(z)$ permit "sheafifying" them into an inverse system of $\Q$-variations of mixed Hodge-Tate structures on $\mathbb P^1(\C) \backslash \{ 0,1,\infty \}$: the classical polylogarithm functions now appear as entries in the period matrices of these variations.\\
\newline
\markright{\uppercase{Introduction}}Beilinson also introduced the $\ell$-adic version of the classical polylogarithm pro-sheaf and moreover showed that the specialization of the Hodge-theoretic resp. $\ell$-adic polylogarithm to primitive roots of unity gives the regulator to absolute Hodge cohomology of his cyclotomic elements in motivic cohomology resp. the Deligne-Soulé elements in $\ell$-adic cohomology (for references cf. \cite{Hu-Wi}, 2).\\
\newline
Following ideas of Beilinson and Deligne, Huber and Wildeshaus \cite{Hu-Wi} revealed that the classes in absolute Hodge and $\ell$-adic cohomology defined by the classical polylogarithm come from a "universal motivic polylogarithm". A corollary of their motivic constructions is a compatibility for Beilinson's cyclotomic elements and the Deligne-Soulé elements, needed for the completion of Bloch's and Kato's proof of the Tamagawa number conjecture for the Riemann zeta function (modulo powers of $2$); for an alternative solution of this problem, also using the machinery of polylogarithms, see \cite{Hu-Ki}.\\
\newline
Bannai \cite{Ba1} developed the syntomic formalism required to transfer the construction of the classical polylogarithm pro-sheaf on the projective line minus three points to the category of filtered overconvergent $F$-isocrystals. He also described explicitly the so defined $p$-adic polylogarithm sheaves and their specialization to roots of unity (cf. also \cite{Ba2}), using $p$-adic polylogarithm functions which were defined by Coleman as analogues of the $\mathrm{Li}_m(z)$ and whose values at roots of unity are related to special values of Kubota-Leopoldt $p$-adic $L$-functions at positive integers. Analogous to the Hodge case these specializations are the image of the motivic Beilinson elements by the syntomic regulator.

\subsubsection{Elliptic polylogarithm}

The concept of the polylogarithm pro-sheaf for elliptic curves was created by Beilinson and Levin in the fundamental paper \cite{Be-Le}. The formalism introduced there is applicable for any reasonable theory of topological or mixed sheaves on a relative elliptic curve, and the elliptic polylogarithm appears as a pro-one-extension on the complement of an étale closed subscheme of the curve, characterized by a certain residue condition. Its specialization along torsion sections induces a collection of cohomology classes on the base, the so-called Eisenstein classes, which in \cite{Be-Le} are determined essentially by a computation of their residue at infinity in the modular case.\\
\newline
The period matrix of the $\Q$-Hodge elliptic polylogarithm is described in \cite{Be-Le} by elliptic polylogarithm functions which are $q$-averaged versions of the above $\mathrm{Li}_m(z)$ and studied extensively in \cite{Le1}. In the $\R$-Hodge case they use real analytic Eisenstein-Kronecker series by which one can further express the specialization to torsion sections (cf. also \cite{Wi2}, V). These also appear in the elliptic Zagier conjecture which predicts that the determinant of a matrix built by these functions applied to certain divisors gives a special $L$-value of the symmetric power of the curve (cf. \cite{Den2}, \cite{Go}, \cite{Wi1}).\\
A different approach to compute the real Hodge polylogarithm sheaves on a single complex elliptic curve can be found in \cite{Ba-Ko-Ts}, App. A. Its two crucial components are: the concrete knowledge of the underlying modules with integrable connection, by which the variations of mixed $\R$-Hodge structures are in fact determined (a manifestation of "rigidity", cf. e.g. \cite{Wi2}, III); it is available from the explicit description given in \cite{Ba-Ko-Ts} for the de Rham realization of the polylogarithm on an elliptic curve over a general subfield of $\C$; and second, based on this knowledge, the definition of the real structures, which is achieved by constructing multivalued meromorphic functions $D_{m,n},D_{m,n}^*$ that solve certain iterated differential equations and give rise to sections inducing these structures.\\
\newline
By \cite{Be-Le}, similar to the classical situation the absolute Hodge and $\ell$-adic cohomology classes given by the elliptic polylogarithm and the Eisenstein classes are the realizations under the respective regulator of a single "motivic elliptic polylogarithm" and associated "motivic Eisenstein classes".\\
It is convenient to anticipate that the same is true when working on abelian schemes (see below).\\
\newline
At this point arises the fundamental meaning of the polylogarithm for Beilinson- and Bloch-Kato-type conjectures: these essentially predict that the motivic cohomology of a smooth projective variety over a number field contains elements by whose images under the Deligne and $\ell$-adic regulator one can express the leading Taylor coefficient at zero of $L$-functions attached to the variety. Hence, the need emerges to construct classes in motivic cohomology whose regulators are accessible to explicit computations, and this is where the polylogarithm comes into play: the Eisenstein classes are generally expected to provide promising elements, and though a concrete description of their regulators is a major nontrivial task, calculations are faciliated by a number of convenient properties enjoyed by the polylogarithm sheaf (e.g. compatibility with base change, norm compatibility, rigidity).\\
\newline
A prime example in this context, illustrating the outlined philosophy, is Kings' \cite{Ki5} proof of the (weak) Bloch-Kato conjecture for elliptic curves over an imaginary quadratic field $K$ with CM by the ring of integers $\mathcal O_K$. It is based on considering the $\mathcal O_K$-linear subspace in motivic cohomology generated by an element which is constructed in Deninger's \cite{Den1} proof of the Beilinson conjecture for Hecke characters and which comes about by applying a variant of the Eisenstein symbol to a certain torsion divisor of the curve. The task imposed by the Bloch-Kato conjecture then consists in determining the $\ell$-adic regulator on this subspace which in turn is known from \cite{Hu-Ki} to express via the $\ell$-adic Eisenstein classes (in fact, the Eisenstein symbol on torsion points and the Eisenstein classes are expected to coincide up to a factor already on the motivic level). The explicit computation of these $\ell$-adic Eisenstein classes is deduced in \cite{Ki5} from a geometric construction of the $\ell$-adic elliptic polylogarithm as inverse limit over torsion points of $1$-motives, and it shows that they are describable in terms of elliptic units resp. elliptic Soulé elements. Results from Iwasawa theory concerning Soulé's elements finally enable to translate this description into a proof of the conjecture.\\
For an explicit treatment of the $\ell$-adic polylogarithm via Iwasawa theory and elliptic units cf. \cite{Ki3}.\\
\newline
The syntomic version of the elliptic polylogarithm was studied extensively by Bannai, Kobayashi and Tsuji \cite{Ba-Ko-Ts} for the situation of a single elliptic curve over an imaginary quadratic field $K$ with CM by $\mathcal O_K$ and a fixed Weierstraß model over $\mathcal O_K$ having good reduction above an unramified prime $p\geq 5$. The technical fundament for defining this $p$-adic elliptic polylogarithm is the theory of rigid syntomic cohomology and its relation to filtered overconvergent $F$-isocrystals as contained in \cite{Ba1}. The $p$-adic polylogarithm sheaves appear as filtered overconvergent $F$-isocrystals on the syntomic datum provided by the curve minus its zero section, and the corresponding cohomology classes are expected to be the image of the motivic elliptic polylogarithm classes by the syntomic regulator.\\
One of the main results in \cite{Ba-Ko-Ts} is the explicit construction of these sheaves. As in the $\R$-Hodge case the crucial ingredients for this are the knowledge of the de Rham realization of the polylogarithm sheaves, by which their syntomic realization is in fact determined ("rigidity"), and the definition of the Frobenius structures, the last achieved by constructing overconvergent functions $D^{(p)}_{m,n}$ that satisfy certain iterated differential equations imposed by the horizontality of the Frobenius isomorphisms.\\
Its motivic origin assumed, the relevance of this $p$-adic elliptic polylogarithm for the $p$-adic Beilinson conjecture comes from the further result of \cite{Ba-Ko-Ts} that it specializes along torsion sections to the so-called $p$-adic Eisenstein Kronecker numbers which in the case of ordinary reduction over $p$ are shown to be directly connected to special values of $p$-adic $L$-functions.\\
So far, a description of the syntomic elliptic polylogarithm for the relative situation, that is to say: for the universal family with level $N$ structure, has not yet been established. Nevertheless, Bannai and Kings \cite{Ba-Ki2} were able to determine the syntomic Eisenstein classes on the ordinary locus of the modular curve in terms of $p$-adic Eisenstein series constructed via a version of Katz's $p$-adic Eisenstein measure. A vital ingredient to prove this result is again the knowledge of the underlying de Rham datum by which the syntomic data are in fact determined. The de Rham Eisenstein classes in turn are shown to be given by certain holomorphic Eisenstein series with explicit formulas for their $q$-expansions; this last fact is proven in \cite{Ba-Ki2} by comparing the residues at the cusps of these modular forms with the residues of the de Rham Eisenstein classes which in turn are obtained by deducing the motivic residues from the étale residues computed in \cite{Be-Le} resp. \cite{Hu-Ki}.\\
For implications of the result of \cite{Ba-Ki2} for the $p$-adic Beilinson conjecture see \cite{Ba-Ki1} and \cite{Ni}.\\
\newline
Before turning closer to the de Rham realization of the elliptic polylogarithm and outlining the ambitions of this work, let us insert some brief remarks concerning more general geometric situations.

\subsubsection{Results in higher dimension and genus}

Wildeshaus \cite{Wi2} constructed and studied the Hodge and $\ell$-adic polylogarithm in the context of mixed Shimura varieties, from which one also obtains the definition for abelian schemes (cf. \cite{Ki4}) with associated Eisenstein classes (cf. \cite{Bl2}). As in the classical and elliptic case the abelian polylogarithm and Eisenstein classes in their realizations have a common motivic origin in $K$-theory (cf. \cite{Ki4}).\\
Blottière \cite{Bl2} described the Hodge polylogarithm for complex abelian schemes by proving that the associated extension of $\C$-pro-local systems - which determines the Hodge data ("rigidity") - can be expressed via "polylogarithmic currents" on the underlying $\mathscr{C}^{\infty}$-manifold; these currents had been constructed by Levin \cite{Le2} as higher-dimensional analogues of Eisenstein-Kronecker series. Drawing on this result, he also gave a (again topological) description of the Eisenstein classes for the situation of a Hilbert-Blumenthal family of abelian varieties and showed that the residue of these classes along the Baily-Borel boundary of the Hilbert modular variety is described by special values of $L$-functions of the defining totally real field (cf. \cite{Bl1}). For a different proof of the last result, resolving the residue computation more functorially by a systematic use of the topological polylogarithm, cf. \cite{Ki2}.\\
The definition of the polylogarithm for arbitrary relative curves (originally part of \cite{Be-Le}, but then excluded) can be found in \cite{Ki1}, where it is shown that the latter induces the polylogarithm of the asscociated Jacobian by taking cup-product with the fundamental class of the curve.\\
Apart from the mentioned results one can generally say that in dimension resp. genus greater than $1$ the polylogarithm sheaves and the Eisenstein classes are - e.g. as to explicit description or as to the latter's non-vanishing and relation to $L$-functions - little understood and subject to active research.

\subsubsection{De Rham realization of the elliptic polylogarithm: the approach of \cite{Ba-Ko-Ts}}

Starting point and main inspiration for this work is the explicit description given in \cite{Ba-Ko-Ts}, 1, for the de Rham realization of the polylogarithm on an elliptic curve $E$ defined over a subfield $F$ of $\C$.\\
\newline
In this situation, if $U:=E-[0], \mathcal H:=H^1_{\mathrm{dR}}(E/F)^\vee$ and $\mathcal H_U:=\mathcal H\otimes_F \mathcal O_U$, the polylogarithm classes $\mathrm{pol}^n_{\mathrm{dR}}$ are by definition the components of an inverse system
\[(\mathrm{pol}_{{\mathrm{dR}}}^n)_{n\geq 1} \in \lim_{n\geq 1} H^1_{\mathrm{dR}}(U/F,\mathcal H_U^\vee \otimes_{\mathcal O_U}\mathcal L_n)\]
which is characterized by a certain residue condition along the divisor given by the zero point $[0]$; here, $\mathcal L_n$ is our notation for the $n$-th logarithm sheaf of $E/F$ with a splitting for its fiber in $[0]$ fixed (in the sense of \cite{Ba-Ko-Ts}, Def. 1.22 and the subsequent explanations).\\
The first step in \cite{Ba-Ko-Ts} to approach a description of $\mathrm{pol}^n_{\mathrm{dR}}$ is an explicit construction of $\mathcal L_1$ on an open affine covering $\{ U_k \}_k$ of $E$, which is done by choosing differentials of the second kind $\{\omega^*,\omega \}$ that give an $F$-basis $\{\underline{\omega}^*,\underline{\omega} \}$ of $H^1_{\mathrm{dR}}(E/F)$ and by then glueing the free modules
\[\mathcal O_{U_k}\cdot \underline{e}_k\oplus \mathcal O_{U_k} \cdot \underline{\omega}^{*\vee}\oplus \mathcal O_{U_k}\cdot \underline{\omega}^\vee\]
to the desired extension $\mathcal L_1$ of $\mathcal O_E$ by $\mathcal H_E$, where the glueing maps, the effect of the connection on the generator $\underline{e}_k$ and the splitting in $[0]$ are defined essentially by using a suitable \v{C}ech cocycle for $\underline{\omega}^*$. This also implies a construction of $\mathcal L_n=\mathrm{Sym}^n_{\mathcal O_E}\mathcal L_1$, where the underlying sheaf on $U_k$ is written as
\[\bigoplus_{0\leq i+j\leq n \atop 0 \leq i,j}\mathcal O_{U_k}\cdot \frac{\underline{e}_k^{n-i-j} \underline{\omega}^{*\vee i}\underline{\omega}^{\vee j}}{(n-i-j)!}.\]
One can then deduce a similar description of $\mathcal L_n$ on $U$. Taking $\mathcal O_U$-linear tensor combinations of the sections $\underline{\omega}^*,\underline{\omega}$ with the thus obtained generators for $\mathcal  L_{n|U}$ and with the differentials $\omega^*,\omega$ defines elements of $\Gamma(U,\mathcal H_U^\vee \otimes_{\mathcal O_U}\mathcal L_n\otimes_{\mathcal O_U}\Omega^1_{U/F})$ and hence classes in $H^1_{\mathrm{dR}}(U/F,\mathcal H_U^\vee \otimes_{\mathcal O_U}\mathcal L_n)$.\\
\newline
The main result, Thm. 1.41, of \cite{Ba-Ko-Ts}, 1, then exhibits the class $\mathrm{pol}^n_{\mathrm{dR}}$ as such a combination, ingeniously constructing the coefficients $L_k \in \Gamma(U,\mathcal O_U)$ occurring in that combination as follows:\\
Writing $E(\C)$ as complex torus $\C/\Gamma$, associated to the line bundle $\mathcal O_{\C/\Gamma}([0])$ is a unique normalized canonical theta function $\theta(z)$ which one can express in terms of the Weierstraß sigma function. It gives rise to a meromorphic function in two variables, the so-called Kronecker theta function for $\Gamma$:
\[\Theta(z,w):=\frac{\theta(z+w)}{\theta(z)\theta(w)}.\]
Modifying $\Theta(z,w)$ by an exponential factor gives a function $\Xi(z,w)$ whose Laurent expansion around $w=0$ yields coefficient functions $L_k(z), k\geq 0$, which turn out to come from $F$-algebraic rational functions $L_k \in \Gamma(U,\mathcal O_U)$. These so-called "connection functions" $L_k$ are the mentioned coefficients that are used in \cite{Ba-Ko-Ts}, 1, for the indicated construction of $\mathrm{pol}^n_{\mathrm{dR}}$; they can in fact be expressed entirely in terms of the Laurent coefficient functions of $\Theta(z,w)$ expanded around $w=0$.\\
A technical inconvenience arises from the fact that in general the pole order of $L_k$ in $[0]$ is greater than one (indeed, if $k\geq 2$ this order is exactly $k$), such that the proof of the mentioned main result requires extensive \v{C}ech calculations to find representatives in logarithmic de Rham cohomology for the constructed cohomology class, which is necessary for performing the required residue computation.

\subsubsection{Outline of the work}

The purpose of this thesis is to establish a new geometric approach to the study of the de Rham realization of the polylogarithm. Our central result in this context shows how to construct the logarithm sheaves of rational abelian schemes from the birigidified Poincaré bundle with universal integrable connection on the product of the abelian scheme and the universal vectorial extension of its dual. This is done essentially by restricting the mentioned data of the Poincaré bundle along the infinitesimal neighborhoods of the zero section of the universal extension. Our perspective also permits a useful interpretation of fundamental formal properties of the logarithm sheaves within the standard theory of the Poincaré bundle. For the situation of a relative elliptic curve we present in addition a related viewpoint on the first logarithm extension in terms of $1$-motives.\\
Having developed in detail the outlined geometric understanding of the logarithm sheaves, we proceed to exploit it systematically for an investigation of the polylogarithm for the universal family of elliptic curves with level $N$ structure. To be more precise, the object in the focus of our study here is a slightly modified version of the usual small elliptic polylogarithm class that provides better access for explicit computations but still contains all the information about the de Rham Eisenstein classes. A main theorem of the work then gives an explicit analytic description for this variant of the polylogarithm via the coefficient functions appearing in the one-variable Laurent expansion of a meromorphic Jacobi form originally defined by Kronecker in the 19th century. Furthermore, using this result, we are able to determine the specialization of the modified polylogarithm along torsion sections concretely in terms of certain algebraic Eisenstein series. From this we regain in particular the already known expressions of the de Rham Eisenstein classes by algebraic modular forms.\\
Our conceptual approach via the Poincaré bundle additionally brings light into the so far rather obscure appearance of theta functions in the study of the elliptic polylogarithm.\\
Moreover and as a matter of future research, we expect our method to produce new insights also for the syntomic resp. higher-dimensional case.

\markright{\uppercase{Overview}}

\
\vspace{1cm}
\section*{\LARGE Overview}
\addcontentsline{toc}{chapter}{Overview}

Let us now discuss the contents of this work in somewhat more detail.\\
\newline
The purpose of the preliminary Chapter 0 on the one hand is to recall the most basic vocabulary required to develop the polylogarithmic formalism in its de Rham realization: this includes the notion of modules with connection, its incorporation in the functorial language of $\mathcal D$-modules as well as the elementary definitions and facts concerning de Rham cohomology. On the other hand, as the universal vectorial extension and the Poincaré bundle are the crucial objects for our geometric construction of the logarithm sheaves on an abelian scheme, we present a thorough account of these concepts, thereby also integrating the viewpoint of extensions and biextensions. Maybe our presentation is also of some use for the non-expert striving for a unified picture of the different facets of the Poincaré bundle.\\
\newline
In the subsequent Chapters 1 and 2 we fix the geometric setting of an abelian scheme $X\xrightarrow{\pi}S$ over a connected base $S$ which is smooth, separated and of finite type over $\Q$; we write $\epsilon$ for its zero section.\\
\newline
\sectionmark{\uppercase{Overview}}The major part of Chapter 1 is occupied with working out the basic formalism of the logarithm sheaves for the given situation $X/S/\Q$ as it expresses in the framework of de Rham cohomology. Of course, the formal structure of the definitions and proofs in this context can (and will) be extracted from their counterparts articulated in other realizations (we will mainly use \cite{Hu-Ki}, App. A, and \cite{Ki4}) and thus are always an elaboration of the condensed exposition given in the original source \cite{Be-Le}. Nevertheless, we have made the experience that a rigorous adjustment to the de Rham setting is at some points not at all immediate and requires supplying a number of additional details and arguments. We have therefore decided to give from the beginning on a thorough self-contained account with full proofs.\\
Specifically, we begin in 1.1 by defining the logarithm sheaves of $X/S/\Q$, written as $(\mathcal L_n,\nabla_n,\varphi_n)$ with $\nabla_n$ the integrable $\Q$-connection of the $\mathcal O_X$-vector bundle $\mathcal L_n$ and $\varphi_n: \prod_{k=0}^n \mathrm{Sym}^k_{\mathcal O_S}\mathcal H \simeq \epsilon^*\mathcal L_n$ the splitting of its zero fiber, where $\mathcal H:=H^1_{\mathrm{dR}}(X/S)^\vee$ is equipped with the dual of the Gauß-Manin connection relative $\Q$. Their relative de Rham cohomology sheaves are computed in 1.2. In 1.3 we elaborate in more detail the viewpoint on the logarithm sheaves as unipotent objects. For this purpose, we introduce a suitable notion of unipotent vector bundles with integrable connection for our situation $X/S/\Q$ and respective categories $U_n(X/S/\Q)$, where $n$ denotes the length of unipotency. We then prove the universal property of the logarithm sheaves which states that with $1^{(n)}:=\varphi_n(\frac{1}{n!})$ the pair $(\mathcal L_n, 1^{(n)})$ is (up to unique isomorphism) the unique pair consisting of an object in $U_n(X/S/\Q)$ and a global horizontal $S$-section of its zero fiber such that for any $\mathcal U$ in $U_n(X/S/\Q)$ the map
\[\pi_*\underline{\Hom}_{\mathcal D_{X/S}}(\mathcal L_n, \mathcal U) \rightarrow \epsilon^*\mathcal U, \quad f \mapsto \epsilon^*(f)\big(1^{(n)}\big)\]
is a horizontal isomorphism (Thm. 1.3.6). Subsequently, we show that the assignment $\mathcal U \mapsto \epsilon^*\mathcal U$ gives an equivalence of $U_n(X/S/\Q)$ with the category of $\mathcal O_S$-vector bundles with integrable $\Q$-connection that carry the structure of a $\prod_{k=0}^n \mathrm{Sym}^k_{\mathcal O_S}\mathcal H$-module with certain compatibilities (Thm. 1.3.13). Side corollaries are further non-evident informations about our categories of unipotent bundles (1.3.4).\\
In 1.4 we discuss the crucial invariance of the logarithm sheaves under isogenies. Finally, in 1.5, we let $X=E$ be an elliptic curve over $S$ and define the elliptic polylogarithm for $E/S/\Q$:
\[\mathrm{pol}_{\mathrm{dR}}=(\mathrm{pol}_{{\mathrm{dR}}}^n)_{n\geq 1} \in \lim_{n\geq 1} H^1_{\mathrm{dR}}(U/\Q,\mathcal H_U^\vee \otimes_{\mathcal O_U}\mathcal L_n),\]
where $U:=E-[0]$ and $\mathcal H_U:=\mathcal H \otimes_{\mathcal O_S} \mathcal O_U$, as well as the $D$-variant of the elliptic polylogarithm:
\[\var = \Big(\varn \Big)_{n\geq 0} \in \lim_{n\geq 0} H^1_{\mathrm{dR}}(U_D/\Q, \mathcal L_n),\]
where $D>1$ is a fixed integer and $U_D:=E-E[D]$. The investigation of this $D$-variant will be the main goal of Chapter 3. The idea to introduce such a better behaved modification of the polylogarithm and to use a formula connecting it with the latter in order to extract from it the Eisenstein classes can be found (for the $\ell$-adic setting) in \cite{Ki3}, 4, and will be adopted in this work.\\
\newline
The heart of Chapter 2 consists in establishing the already indicated geometric construction of the logarithm sheaves via the birigidified Poincaré bundle $(\mathcal P,r,s,\nabla_{\mathcal P})$ on $X\times_S \Yr$, where $\Yr$ is the universal extension of the dual abelian scheme $Y$ and $r$ resp. $s$ is the rigidification along $\Yr$ resp. $X$.\\
However, a basic difficulty arises when one aims at such a construction: the universal integrable connection $\nabla_{\mathcal P}$ of $\mathcal P$ being a relative connection (namely,  relative to $\Yr$), we can construct the logarithm sheaves only with their $\Q$-connections restricted relative $S$. Our remedy for this problem is contained in 2.1: there, we prove that given an extension of vector bundles on $X$ with integrable $S$-connection \markright{\uppercase{Overview}}
\[0 \rightarrow \mathcal H_X \rightarrow \mathcal L_1' \rightarrow \mathcal O_X \rightarrow 0,\]
mapping to the identity under the natural projection
\[\tag{1} \mathrm{Ext}^1_{\mathcal D_{X/S}}(\mathcal O_X,\mathcal H_X) \rightarrow \mathrm{Hom}_{\mathcal O_S}(\mathcal O_S,\mathcal H^\vee\otimes_{\mathcal O_S}\mathcal H),\]
and a $\mathcal O_S$-linear splitting $\varphi_1':\mathcal O_S\oplus \mathcal H \simeq \epsilon^*\mathcal L_1'$ for its pullback via $\epsilon$, then the $S$-connection on $\mathcal L_1'$ uniquely extends to an integrable $\Q$-connection such that the previous data become the first logarithm sheaf of $X/S/\Q$ (Prop. 2.1.4). Then, in 2.3, we proceed to construct such data from $(\mathcal P,r,s,\nabla_{\mathcal P})$.\\
To be more precise, if $\Yr_1$ is the first infinitesimal neighborhood of the zero section of $\Yr$ and $\mathcal P_1$ the restriction of $\mathcal P$ along $X\times_S \Yr_1\rightarrow X\times_S \Yr$, equipped with the induced integrable $\Yr_1$-connection, then by adjunction along the natural morphism $\iota_1:X \rightarrow X\times_S \Yr_1$ together with the rigidification $s$ and the crucial identification $\mathrm{Lie}(\Yr/S)^\vee \simeq \mathcal H$ one obtains a horizontal exact sequence
\[0 \rightarrow (\iota_1)_*\mathcal H_X \rightarrow \mathcal P_1 \rightarrow (\iota_1)_*\mathcal O_X \rightarrow 0\]
whose pushout along the projection $p_1:X\times_S \Yr_1\rightarrow X$ gives the exact sequence of vector bundles on $X$ with integrable $S$-connection
\[\tag{2} 0 \rightarrow \mathcal H_X \rightarrow (p_1)_*\mathcal P_1 \rightarrow \mathcal O_X \rightarrow 0.\]
On the other hand, it is easy to construct from the rigidification $r$ a natural $\mathcal O_S$-linear splitting
\[\mathcal O_S \oplus \mathcal H\simeq\epsilon^*(p_1)_*\mathcal P_1.\]
Our main result in Chapter 2 (Thm. 2.3.1 resp. Cor. 2.3.2) proves that $(2)$ maps to the identity in $(1)$, such that due to the above explanations we have achieved a construction of the first logarithm sheaf of $X/S/\Q$ (with the mentioned limitation concerning our knowledge of the absolute connection).\\
Though we always work with the higher logarithm sheaves as the symmetric powers of the first, it is natural to ask how also they can be obtained from the infinitesimal geometry of the data $(\mathcal P,r,s,\nabla_{\mathcal P})$. This is explained in detail in 2.4, where it also turns out that in our geometric perspective the comultiplication of the logarithm sheaves expresses via the biextension structure of the Poincaré bundle.\\
Our constructions have a natural and appealing formulation within the language of the Fourier-Mukai transformation for $\mathcal D$-modules on abelian schemes as formulated in \cite{Lau}. To illuminate this point in full clarity we extend in 2.2 the basic theory developed in \cite{Lau} by introducing (entirely in the spirit of Mukai \cite{Mu}) the notion of WIT-sheaves on $\Yr$, which permits to leave the derived categories and consider honest sheaves. We also define categories of unipotent sheaves $U_n(\Yr/S)$ and $U_n(X/S)$ on $\Yr$ and $X$, the last being analogous to the above $U_n(X/S/\Q)$ but forgetting $\Q$-structures, and prove in particular that Fourier transformation establishes an equivalence between them (Thm. 2.2.12). This and other results of 2.2 are of independent interest, but they were actually also our heuristic starting point for the interpretation of the logarithm sheaves by the Poincaré bundle: noting that $\mathcal L_n$ defines an object of $U_n(X/S)$, one immediately obtains from Thm. 2.2.12 that it must be the Fourier transform of a sheaf on $\Yr$ which in fact lives on the $n$-th infinitesimal neighborhood of $S$ in $\Yr$. And indeed, as is explained in 2.3 resp. 2.4, the standard exact sequences of $\mathcal O_X$-vector bundles with $S$-connection
\[0 \rightarrow \mathrm{Sym}^n_{\mathcal O_X}\mathcal H_X \rightarrow \mathcal L_n \rightarrow \mathcal L_{n-1} \rightarrow 0\]
are simply the Fourier transforms of the canonical exact sequences of WIT-sheaves of index $0$ on $\Yr$:
\[0 \rightarrow \mathcal J^n/\mathcal J^{n+1} \rightarrow \mathcal O_{\Yr}/\mathcal J^{n+1} \rightarrow \mathcal O_{\Yr}/\mathcal J^n \rightarrow 0,\]
where $\mathcal J$ denotes the augmentation ideal of the zero section of $\Yr$.\\
In 2.5 we reveal that in our viewpoint the invariance of the logarithm sheaves under isogenies is the manifestation of a symmetry isomorphism of the Poincaré bundle for the isogeny and its transpose map on the universal extension. This will also be relevant for explicit computations in Chapter 3.\\
Finally, in 2.6 we elaborate the $1$-motivic origin of the first logarithm extension $\mathcal Log^1$. For this one needs at first to be able to equip the de Rham realization of a $1$-motive over a base scheme in a natural way with an integrable connection, which is a nontrivial problem. In a recent work of Andreatta and Bertapelle \cite{An-Ber} such a "motivic Gauß-Manin connection" is constructed in full generality by using crystalline techniques. We apply their results in the following way: if $X=E$ is an elliptic curve and $E\times_S E$ is considered as $1$-motive over $E$ via the second projection, the Barsotti-Rosenlicht-Weil isomorphism and taking de Rham realizations with motivic Gauß-Manin connections produces a map
\[\tag{3} (E\times_S E)(E) \rightarrow \Ext^1_{\mathcal D_{E/\Q}}(\mathcal H_E^\vee,\mathcal O_E)\]
which is the negative of the "motivic de Rham-Manin map" investigated in \cite{An-Ber} (cf. Rem. 2.6.2). We then prove that $\mathcal Log^1$ is the extension dual to the image of the diagonal $\Delta_E$ under $(3)$ (Thm. 2.6.3). To achieve this we relate the extension $\mathcal Log^1$ to the more explicit "classical Manin map" of \cite{Co2} and use the latter's relation with the motivic de Rham-Manin map as established in \cite{An-Ber} (to dispose of that comparison is also the reason why we restrict to relative dimension one).\\
Explicating this further leads to another viewpoint on the relation between the first logarithm extension and the Poincaré bundle: namely, the $\mathbb G_{m,E\times_S (E^\vee)^\natural}$-torsor $P$ associated to the Poincaré bundle $\mathcal P$ on $E\times_S(E^\vee)^\natural$ naturally sits in an exact sequence of $E$-group schemes
\[0 \rightarrow \mathbb G_{m,E} \rightarrow P \rightarrow E\times_S (E^\vee)^\natural \rightarrow 0\]
whose associated sequence of Lie algebras relative $E$ reads as
\[\tag{4} 0 \rightarrow \mathcal O_E \rightarrow \Lie(P/E) \rightarrow \mathcal H_E^\vee \rightarrow 0.\]
One can equip all terms in $(4)$ with motivic Gauß-Manin connections relative $\Q$, and the dual of the thus obtained extension is $\mathcal Log^1$ (Cor. 2.6.16).\\
We wish to remark that the basic idea that one should be able to obtain the logarithm sheaves as the formal completion of the Poincaré bundle was pointed out to us by Guido Kings. An initial hint for a relation between these objects can be seen in the observation that the Kronecker theta function used in \cite{Ba-Ko-Ts} is a meromorphic section of the Poincaré bundle (cf. \cite{Ba-Ko}, 1.2).\\
\newline
Chapter 3 uses the geometric approach towards the logarithm sheaves developed in Chapter 2 to find a way of describing the $D$-variant of the polylogarithm for the universal elliptic curve with level $N$ structure and to derive from the latter also a description for the specialization along torsion sections.\\
In 3.1 we fix for a general elliptic curve $E/S$ with dual abelian scheme $\widehat{E}/S$ the principal polarization
\[E \xrightarrow{\sim} \widehat{E}\]
defined by the ample invertible sheaf $\mathcal O_E([0])$. The birigidified Poincaré bundle on $E\times_S E$ then is
\[\tag{5} (\mathcal M \otimes_{\mathcal O_{E\times_S E}}(\pi \times \pi)^*\epsilon^*\mathcal O_E([0]), \mathrm{can}, \mathrm{can}),\]
where $\mathcal M$ denotes the Mumford bundle for $\mathcal O_E([0])$:
\[\mathcal M:= \mu^* \mathcal O_E([0]) \otimes_{\mathcal O_{E\times_S E}} \mathrm{pr}_1^* \mathcal O_E([0])^{-1} \otimes_{\mathcal O_{E\times_S E}} \mathrm{pr}_2^*\mathcal O_E([0])^{-1},\]
and where $\mathrm{can}$ means the obvious canonical rigidification along the second resp. first factor of $E\times_SE$.\\
The Poincaré bundle on $E\times_S \widehat{E}^\natural$ and its rigidifications arise from $(5)$ by pullback via the natural map
\[\tag{6} E\times_S \widehat{E}^\natural \rightarrow E\times_S \widehat{E} \xrightarrow{\sim} E\times_S E.\] 
Large parts of our considerations will take place on the analytic side. As a convenient method to work with vector bundles on complex manifolds we use the yoga of automorphy matrices, i.e. we fix a trivialization for the pullback of the bundle to the universal covering, compute the matrix describing the effect of deck transformations on the chosen trivializing sections and then express sections of the bundle as vectors of holomorphic functions on the universal covering transforming under deck transformations with this automorphy matrix. The details are explained in 3.2. The advantage of this approach via the universal covering is that it avoids choosing open coverings and coordinate charts.\\
In 3.3 we introduce the function by whose inverse we will trivialize $\mathcal O_{E^{an}}([0])$ (componentwise) on the universal covering of the analytification $E^{an}$ of the universal elliptic curve $E$ with level $N$ structure ($N\geq 3$) over the modular curve $S$: this is the "elementary theta function"
\[\tag{7} \vartheta(z,\tau):=\exp\bigg[\frac{z^2}{2}\eta(1,\tau)\bigg]\cdot \sigma(z,\tau),\]
where $\eta(1,\tau)$ is the quasi-period (equivalently: is $-G_2(\tau)$). It differs from the theta function $\theta(z,\tau)$ used in \cite{Ba-Ko-Ts} by an exponential factor; the crucial point is that $\theta(z,\tau)$ does not vary holomorphically in both variables, whereas $\vartheta(z,\tau)$ does. What we do here is performing a shift from canonical to classical theta functions, explained in detail in the first part of 3.3; we only remark that for $\tau\in \H$ fixed the function $\vartheta(z;\tau)$ is the unique holomorphic function on $\C$ whose inverse induces the classical factor of automorphy for $\mathcal O_{\C/\Gamma_{\tau}}([0])$ and such that its derivative in $z=0$ is normalized to $1$.\\
The trivialization via (the inverse of) $(7)$ for $\mathcal O_{E^{an}}([0])$ on the universal covering of $E^{an}$ induces a trivialization for $\mathcal M^{an}$ on the universal covering of $E^{an}\times_{S^{an}}E^{an}$ (always meant componentwise): it is given by (the inverse of) the "fundamental meromorphic Jacobi form"
\[\tag{8} J(z,w,\tau):=\frac{\vartheta(z+w,\tau)}{\vartheta(z,\tau)\vartheta(w,\tau)}=\exp[zw\cdot \eta(1,\tau)]\cdot \frac{\sigma(z+w,\tau)}{\sigma(z,\tau)\sigma(w,\tau)}.\]
This function is an exponential term times the Kronecker theta function $\Theta(z,w,\tau)$ of \cite{Ba-Ko-Ts}, it is $2\pi i$-times the meromorphic Jacobi form $F(2\pi i z,2\pi i w, \tau)$ used in \cite{Za2} to construct a generating function for the period polynomials of all Hecke eigenforms for the full modular group, it appears in \cite{Le1} to relate the Debye elliptic polylogarithm functions to Eisenstein functions and series, and it equals the function introduced in \cite{Le-Ra}, 2.2, to describe the relative nilpotent de Rham fundamental torsor for a pointed family of elliptic curves. Its definition goes back to Kronecker (cf. \cite{Le-Ra}, 2.2.1).\\
In the rest of 3.3 we investigate some of the analytic properties of the functions $(7)$ and $(8)$ that will become important for us; in particular, we examine the coefficient functions in the Laurent expansion around $w=0$ of $J(z,w,\tau)$ and reveal their connection to modular forms (Thm. 3.3.16).\\ 
From $(7)$ we also get an induced trivialization on the universal covering for the second factor in $(5)$:
\[(\pi^{an} \times \pi^{an})^*(\epsilon^{an})^*\mathcal O_{E^{an}}([0])=(\pi^{an}\times \pi^{an})^*\omega_{E^{an}/S^{an}}^\vee\]
which coincides with the trivialization induced by the dual of the canonical differential form.\\
In sum, the function $(7)$ provides us with a trivialization for $(5)$ on the universal covering; by taking pullback via $(6)$ we obtain a trivializing section $t$ on the universal covering for the Poincaré bundle $\mathcal P^{an}$ on $E^{an}\times_{S^{an}}(\widehat{E}^\natural)^{an}$ whose factor of automorphy we compute. We then express in this language the birigidification of $\mathcal P^{an}$ and give an explicit formula for $\nabla_{\mathcal P}^{an}$. The details are explained in 3.4.\\
With this explicit knowledge of $(\mathcal P^{an},r^{an},s^{an},\nabla_{\mathcal P}^{an})$ and with the main result of Chapter 2 (which constructs the first logarithm sheaf of $E/S/\Q$ from the data $(\mathcal P,r,s,\nabla_{\mathcal P})$) we can proceed in 3.5 to give a description of the analytified logarithm sheaves $(\mathcal L_n^{an}, \nabla_n^{an},\varphi_n^{an})$ on $E^{an}$.\\
Writing $\mathrm{pr}: \C\times \H \rightarrow E^{an}$ for the projection of the universal covering ($E^{an}$ now means a fixed connected component) we construct from $t$ a global section $e$ of $\mathrm{pr}^*\mathcal L_1^{an}$ which splits the exact sequence
\[0 \rightarrow \mathrm{pr}^*\mathcal H^{an}_{E^{an}} \rightarrow \mathrm{pr}^*\mathcal L_1^{an} \rightarrow \mathrm{pr}^*\mathcal O_{E^{an}} \rightarrow 0.\]
Trivializing $\mathrm{pr}^*\mathcal H^{an}_{E^{an}}$ by the basic sections $\{f, g \}$ defined by $\{\eta^\vee, \omega^\vee\}$ we obtain a trivialization
\[\mathrm{pr}^*\mathcal L_1^{an} = \mathcal O_{\C\times \H}\cdot e \oplus \mathcal O_{\C\times \H} \cdot f \oplus \mathcal O_{\C\times \H}\cdot g\]
and induced ones for the other logarithm sheaves (note the analogy with \cite{Ba-Ko-Ts}, Cor. 1.28):
\[\mathrm{pr}^*\mathcal L_n^{an} = \bigoplus_{0\leq i+j\leq n \atop 0 \leq i,j}\mathcal O_{\C\times \H}\cdot \frac{e^{n-i-j} f^ig^j}{(n-i-j)!}.\]
Fixing these, we derive from our knowledge of $(\mathcal P^{an},r^{an},s^{an},\nabla_{\mathcal P}^{an})$ explicit formulas for the induced automorphy matrices (Prop. 3.5.2 and Prop. 3.5.6), for the restrictions of the connections $\nabla_n^{an}$ relative $S^{an}$ (Prop. 3.5.3 and Prop. 3.5.7), for the splittings $\varphi_n^{an}$ (Prop. 3.5.8) and for the pullback of sections of $\mathcal L_n^{an}$ resp. $\Omega^1_{E^{an}}\otimes_{\mathcal O_{E^{an}}}\mathcal L_n^{an}$ along torsion sections (Cor. 3.5.13 resp. Prop. 3.5.14).\\
At this point one piece is still missing for a complete knowledge of the analytified logarithm sheaves: a description for their absolute connections $\nabla_n^{an}$. We solve this problem by characterizing $\nabla_n^{an}$ analogously as in the algebraic situation (Prop. 3.5.20 resp. Prop. 2.1.4) and then prove that a reasonable guess of a formula for $\nabla_n^{an}$ satisfies all required conditions of this characterization (Thm. 3.5.21).\\
With these preparations the next goal is to establish a concrete description for the system
\[\tag{9} (\varD)^{an}=\Big(\varDn\Big)^{an}_{n\geq0} \in \lim_{n\geq 0} H^1_{\mathrm{dR}}(U_D^{an}, \mathcal L_n^{an})\]
which is defined as the image of $\varD$ under the natural analytification map
\[\tag{10} \lim_{n\geq 0} H^1_{\mathrm{dR}}(U_D/\Q,\mathcal L_n) \hookrightarrow \lim_{n\geq 0} H^1_{\mathrm{dR}}(U_D^{an},\mathcal L_n^{an}).\]
We approach this problem by proving that $(9)$ is characterized analytically in the same way as $\varD$ is algebraically, i.e. it is the unique element in $\lim_{n\geq 0} H^1_{\mathrm{dR}}(U_D^{an},\mathcal L_n^{an})$ having a certain residue under an analytic residue map (Thm. 3.7.5). The reason why such a characterization is possible is the regularity of the logarithm sheaves (Prop. 3.7.1); this regularity is also responsible for the injectivity of $(10)$, hence for the fact that $\varD$ is determined by its analytification.\\
To construct the cohomology classes $(9)$ we return to the fundamental meromorphic Jacobi form.\\
Namely, we consider the functions $s_k^D(z,\tau)$ defined by the Laurent expansion
\[D^2 \cdot J(z,-w,\tau)-D\cdot J\Big(Dz,-\frac{w}{D},\tau \Big)=s^D_0(z,\tau)+s^D_1(z,\tau)w+...\]
They are meromorphic functions on $\C \times \H$, holomorphic on $\mathrm{pr}^{-1}(U_{D}^{an})$, and satisfy (cf. 3.3.3):

\[{\small
\tag{11}
\begin{split} &s^D_k \ \textrm{has at worst simple poles along} \ z=m\tau+n \ (m,n\in \Z, \tau \in \H), \textrm{with residue} \ (D^2-1)\cdot \frac{(2\pi i m)^k}{k!},\\
&\textrm{and along} \ z=\frac{m}{D}\tau+\frac{n}{D} \ (\textrm{with} \ D \ \textrm{not simultaneously dividing} \ m \ \textrm{and} \ n),  \textrm{with residue} \ -\frac{(2\pi i \frac{m}{D})^k}{k!}.
\end{split}
}
\]
In 3.6 we construct from the functions $s_k^D(z,\tau)$ a certain vector $p_n^D(z,\tau)$ of functions (cf. $(3.6.2)$) and show in laborious calculations that this vector defines an element of $\Gamma(U_D^{an},\Omega^1_{E^{an}}\otimes_{\mathcal O_{E^{an}}} \mathcal L_n^{an})$ which goes to zero in the de Rham complex (Thm. 3.6.2); here, we trivialize $\mathcal L_n^{an}$ on the universal covering of each component of $E^{an}$ as above and $\Omega^1_{E^{an}}$ by $\{\mathrm{d}z,\mathrm{d}\tau \}$. The $p_n^D(z,\tau)$ are compatible for the transition maps of the logarithm sheaves and thus induce an inverse system
\[p^D=(p^D_n)_{n\geq 0} \in \lim_{n\geq 0} H^1_{\mathrm{dR}}(U_D^{an},\mathcal L_n^{an}).\]
Our main result about the $D$-variant of the polylogarithm (Thm. 3.8.3) then is the equality
\[\tag{12} (\varD)^{an}=p^D.\]
Its proof is rather technical, but the crucial points are the already mentioned characterization of the left side and a computation of the residue of the right side using $(11)$; here, it is of enormous convenience that the $s_k^D$ have at worst simple poles and hence already define a logarithmic de Rham class.\\
From now on assume $(D,N)=1$ and that $a,b$ are two integers not simultaneously divisible by $N$. Via the Drinfeld basis $(e_1,e_2) \in E[N](S)$ for $E[N]$ one obtains the $N$-torsion section
\[t_{a,b}:=ae_1+be_2: S \rightarrow U_D \subseteq E.\]
We consider the specialization of $\varDn (n\geq 0)$ along $t_{a,b}$, in the following sense: let
\[\big(t_{a,b}^*(\varDn)\big)^{(n)} \in H^1_{\mathrm{dR}}(S/\Q,\mathrm{Sym}^n_{\mathcal O_S} H^1_{\mathrm{dR}}(E/S))\]
be the de Rham cohomology class received by pulling back
\[\varDn \in H^1_{\mathrm{dR}}(U_D/\Q, \mathcal L_n)\]
along $t_{a,b}$, by using the (horizontal) identifications
\[t_{a,b}^*\mathcal L_n \simeq \epsilon^*\mathcal L_n \simeq \prod_{k=0}^n \mathrm{Sym}^k_{\mathcal O_S}\mathcal H\]
that come from the invariance of $\mathcal L_n$ under $N$-multiplication (cf. 1.4.2) and from the splitting $\varphi_n$, by then taking the $n$-th component and by finally identifying
\[\tag{13} \mathrm{Sym}^n_{\mathcal O_S}\mathcal H \simeq \mathrm{Sym}^n_{\mathcal O_S} H^1_{\mathrm{dR}}(E/S)\]
via the following choice of the Poincaré duality isomorphism:
\[H^1_{\mathrm{dR}}(E/S) \xrightarrow{\sim} \mathcal H, \quad x \mapsto \{\ y\mapsto \mathrm{tr}(x\cup y)\}.\]
On the other hand, the Hodge filtration and Kodaira-Spencer map induce a canonical homomorphism
\[\tag{14} \Gamma\big(S, \omega_{E/S}^{\otimes(n+2)}\big) \rightarrow H^1_{\mathrm{dR}}(S/\Q,\mathrm{Sym}^n_{\mathcal O_S}H^1_{\mathrm{dR}}(E/S))\]
defined on the space of weakly holomorphic algebraic modular forms of weight $n+2$ and level $N$.\\
In our main result about the specialization of the $D$-variant of the elliptic polylogarithm (Thm. 3.8.15) we show that $\big(t_{a,b}^*(\varDn)\big)^{(n)}$ is the image under $(14)$ of the algebraic modular form
\[\begin{cases}
-{_D}F^{(2)}_{\frac{a}{N},\frac{b}{N}} \ \ &\textrm{if} \ \ n=0
\vspace{0.9mm}\\
\frac{(-1)^{n}}{n!} \cdot {_D}F^{(n+2)}_{\frac{a}{N},\frac{b}{N}} \ \ &\textrm{if} \ \ n> 0,
\end{cases}
\]
where in general for $k\geq 1$:
\[_DF^{(k)}_{\frac{a}{N},\frac{b}{N}}:=D^2F^{(k)}_{\frac{a}{N},\frac{b}{N}}-D^{2-k}F^{(k)}_{\frac{Da}{N},\frac{Db}{N}}.\]
The $F^{(k)}_{\frac{a}{N},\frac{b}{N}}$ resp. $F^{(k)}_{\frac{Da}{N},\frac{Db}{N}}$ in turn are algebraic modular forms of weight $k$ and level $N$ constructed as in Ch. I, 3, of Kato's work \cite{Ka}, where they are used to define the Euler system of zeta elements in the space of modular forms, related to operator-valued zeta functions via a period map. Essentially, they are given as averaged sum of algebraic Eisenstein series won by specializing along torsion sections certain iterated derivations of the logarithmic derivative of Kato-Siegel functions (at least for $k\neq 2$, otherwise one specializes an algebraic Weierstraß $\wp$-element). Their analytic expressions as holomorphic functions in $\tau$ can be found in 3.3.4, observing the explanations in Rem. 3.8.12 and Rem. 3.8.13.\\ 
For the proof of the theorem we first resolve the problem on the analytic side, using the crucial result $(12)$ and the fact that we can explicitly compute the analytic specialization of the section $p_n^D(z,\tau)$ along $t_{a,b}^{an}$ (Thm. 3.6.5). From this analytic result we can then in fact deduce the algebraic statement. The detailed strategy of proof is explained at the beginning of 3.8.2.\\
With the already indicated relation between the polylogarithm and its $D$-variant the previous theorem yields in particular a formula for the de Rham Eisenstein classes at $t_{a,b}$, the latter defined as
\[\mathrm{Eis}^n(t_{a,b}):=-N^{n-1}\cdot\big(\mathrm{contr}_n(t_{a,b}^*\mathrm{pol}^{n+1}_{\mathrm{dR}})\big)^{(n)} \in H^1_{\mathrm{dR}}(S/\Q,\mathrm{Sym}^n_{\mathcal O_S}\mathcal H),\]
where $\mathrm{contr}_n: H^1_{\mathrm{dR}}\Big(S/\Q,\mathcal H^\vee\otimes_{\mathcal O_S}\prod_{k=0}^{n+1}\mathrm{Sym}^k_{\mathcal O_S}\mathcal H\Big) \rightarrow H^1_{\mathrm{dR}}\Big(S/\Q,\prod_{k=0}^n\mathrm{Sym}^k_{\mathcal O_S}\mathcal H\Big)$ is a certain contraction map. Namely, we show the equalities
\[\begin{cases}
\mathrm{Eis}^0(t_{a,b})=-N^{-1}\cdot\big(-F^{(2)}_{\frac{a}{N},\frac{b}{N}}\big) \quad &\textrm{in} \ \ H^1_{\mathrm{dR}}(S/\Q)
\vspace{0.9mm}\\
\mathrm{Eis}^n(t_{a,b})=-N^{n-1} \frac{(-1)^{n}}{n!}\cdot F^{(n+2)}_{\frac{a}{N},\frac{b}{N}} \quad &\textrm{in} \ \ H^1_{\mathrm{dR}}(S/\Q,\mathrm{Sym}^n_{\mathcal O_S}\mathcal H), \ n>0, 
\end{cases}\]
where here $F^{(2)}_{\frac{a}{N},\frac{b}{N}}$ resp. $F^{(n+2)}_{\frac{a}{N},\frac{b}{N}}$ means the element of $H^1_{\mathrm{dR}}(S/\Q)$ resp. $H^1_{\mathrm{dR}}(S/\Q,\mathrm{Sym}^n_{\mathcal O_S}\mathcal H)$ induced by the algebraic modular form $F^{(2)}_{\frac{a}{N},\frac{b}{N}}$ resp. $F^{(n+2)}_{\frac{a}{N},\frac{b}{N}}$ via $(14)$ resp. via $(14)$ and $(13)$.\\
As was already explained in more detail during the Introduction, a very different way to determine the de Rham Eisenstein classes can be found in \cite{Ba-Ki2}.

\
\vspace{1cm}
\section*{\LARGE Acknowledgements}
\addcontentsline{toc}{chapter}{Acknowledgements}
\markright{\uppercase{Acknowledgements}}
I would like to express my profound gratitude to my advisor Guido Kings for his constant support, encouragement and optimism during my work on this thesis. He has never been tired to share his intuition with me and has provided me with innumerable useful suggestions.\\
It is a pleasure to further thank my colleagues Sandra Eisenreich, Philipp Graf, Johannes Sprang and Georg Tamme: each of them has devoted on more than one occasion time and energy to valuable discussions and has helped me straighten out details of the work.\\
Part of this research was done during a semestral stay at the Università degli Studi di Padova. I heartily thank my host Bruno Chiarellotto as well as Alessandra Bertapelle for great hospitality and for integration into the stimulating atmosphere of the Dipartimento di Matematica. There, I also had the opportunity to present a first draft of the results of this work, and I thank all participants of my lectures for valuable comments and questions. I am in particular obliged to Fabrizio Andreatta, Francesco Baldassarri, Bruno Chiarellotto, Adrian Iovita, Matteo Longo and Shanwen Wang for fruitful discussions. Special thanks go to Alessandra Bertapelle for helping me understand the details of the work \cite{An-Ber} and for her patient support in my search of the lost sign.\\
Moreover, I would like to express my gratitude to Kenichi Bannai, Minhyong Kim, Lars Kindler, Bernd Schober and Tobias Sitte  for helpful discussions and for kindly answering my questions.\\
I thank Tobias Sitte and Georg Tamme for patiently helping me with the fine-tuning of the \LaTeX - file.\\
Since the beginning of the last year I have been supported in my research by a scholarship of the Studienstiftung des deutschen Volkes, which I gratefully acknowledge.

\newpage
\setcounter{chapter}{-1}
\chapter{Preliminaries and notation}
\section{Abelian schemes: duality theory, universal vectorial extension and Poincaré bundle}
\sectionmark{Abelian schemes: duality theory, universal vectorial extension...}
\markright{\uppercase{Preliminaries and notation}}

\subsection{Introduction of the basic objects}
Given an abelian scheme, the universal vectorial extension of its dual and the associated birigidified Poincaré bundle with universal integrable connection will be key instruments for this work.\\
We therefore begin with an adequate review of all of these notions, striving to be as self-contained and detailed as seems possible without going beyond the scope of preliminary remarks. Our presentation consists in compiling and supplementing scattered material from the literature, whereby we mention as our main sources \cite{Ch-Fa}, Ch. I, 1, \cite{Lau}, $(1.1)$ and $(2.1)$-$(2.2)$, and \cite{Maz-Mes}, Ch. I.

\subsubsection{Algebraic equivalence to zero and rigidifications}
If $A$ is an abelian scheme over an arbitrary base scheme $B$ we denote by
\[\begin{split} \pi_A \colon A \rightarrow B\\
\mu_A \colon A\times_B A \rightarrow A\\
\epsilon_A \colon B \rightarrow A\\
(-1)_A \colon A\rightarrow A\\
\mathrm{pr}_{1,A}, \ \mathrm{pr}_{2,A} \colon A \times_B A \rightarrow A \end{split}\]
the structure map, the multiplication map, the zero section, the inverse map and the two projections.

\begin{definition}
A line bundle $\mathcal L$ on $A$ is \underline{algebraically equivalent to zero} if the line bundle on $A \times_B A$ given by
\[\mu_A^* \mathcal L \otimes_{\mathcal O_{A\times_B A}} \mathrm{pr}_{1,A}^* \mathcal L ^{-1} \otimes_{\mathcal O_{A\times_B A}} \mathrm{pr}_{2,A}^*\mathcal L^{-1}\]
is \underline{trivial over $B$}, i.e. if it is isomorphic to $(\pi_A \times \pi_A)^* \mathcal M$ for some line bundle $\mathcal M$ on $B$.\\
Here, $\pi_A\times \pi_A: A\times_B A \rightarrow B$ is the canonical map and $\mathcal L^{-1}=\underline{\Hom}_{\mathcal O_A}(\mathcal L,\mathcal O_A)$ is the dual of $\mathcal L$.
\end{definition}
\begin{remark}
It is easy to check that algebraic equivalence to zero is compatible with base change in the following sense: if we have a morphism $B'\rightarrow B$, set $A':=A\times_B B'$ and consider the cartesian square
\begin{equation*}
\begin{xy}
\xymatrix{
A' \ar[r]\ar[d]  & B' \ar[d]\\
A \ar[r]^{\pi_A} & B}
\end{xy}
\end{equation*}
then pulling back line bundles along $A' \rightarrow A$ preserves the property of being algebraically equivalent to zero. Of course, we regard $A'$ as abelian scheme over $B'$ in the natural way.\\
It is moreover clear that the property of being algebraically equivalent to zero is stable under the formation of tensor product and dual, and that the line bundle $\mathcal O_A$ has this property.
\end{remark}
\begin{definition}
(i) Let $\mathcal L$ be a line bundle on $A$. By a \underline{$B$-rigidification} of $\mathcal L$ we mean an isomorphism of $\mathcal O_B$-modules
\[\alpha: \mathcal O_B \xrightarrow{\sim} \epsilon_A^*\mathcal L.\]
(ii) Let $(\mathcal L_1,\alpha_1)$ and $(\mathcal L_2,\alpha_2)$ be two $B$-rigidified line bundles on $A$. An \underline{isomorphism} between them is an isomorphism $\varphi: \mathcal L_1 \xrightarrow{\sim} \mathcal L_2$ of the line bundles which is compatible with the rigidifications in the obvious sense, i.e. $\epsilon_A^*(\varphi)$ becomes the identity on $\mathcal O_B$ when using the isomorphisms $\alpha_1$ and $\alpha_2$.\\
(iii) The \underline{tensor product} of two $B$-rigidified line bundles $(\mathcal L_1,\alpha_1)$ and $(\mathcal L_2,\alpha_2)$ on $A$ is the pair $(\mathcal L_1 \otimes_{\mathcal O_A} \mathcal L_2, \alpha_1 \otimes \alpha_2)$, where $\alpha_1 \otimes \alpha_2$ means the obvious induced $B$-rigidification of $\mathcal L_1\otimes_{\mathcal O_A} \mathcal L_2$.\\
(iv) The \underline{inverse} of a $B$-rigidified line bundle $(\mathcal L,\alpha)$ on $A$ is the pair $(\mathcal L^{-1},\alpha^{-1})$, where $\alpha^{-1}$ is the $B$-rigidification of $\mathcal L^{-1}$ naturally induced by dualizing $\alpha$.\\
(v) The line bundle $\mathcal O_A$ together with its canonical $B$-rigidification will be written $(\mathcal O_A,\mathrm{can})$.
\end{definition}

\begin{remark}
Given the situation of Rem. 0.1.2 we have a commutative (in fact cartesian) diagram
\begin{equation*}
\begin{xy}
\xymatrix{
A' \ar[d]  & B' \ar[l]_{ \ \epsilon_{A'}} \ar[d]\\
A & B \ar[l]_{ \ \epsilon_A}}
\end{xy}
\end{equation*}
which shows that the pullback of a $B$-rigidified line bundle $\mathcal L$ on $A$ along the projection $A'\rightarrow A$ is naturally equipped with an induced $B'$-rigidification.
\end{remark}

\textit{\textbf{From now on we fix an abelian scheme $X$ of relative dimension $g$ over a locally noetherian base $S$.\\
We will write $\pi$, $\mu$, $\epsilon$, $\mathrm{pr}_1$, $\mathrm{pr}_2$ instead of $\pi_X$, $\mu_X$, $\epsilon_X$, $\mathrm{pr}_{1,X}$, $\mathrm{pr}_{2,X}$, but we keep the notation $(-1)_X$.\\
For $S$-schemes $T$ we will often use the abbreviation $X_T$ for the base extension $X\times_S T$ and view it as abelian scheme over $T$ in the natural way}}.

\subsubsection{Dual abelian scheme and Poincaré bundle}
Consider the \underline{dual functor of $X/S$}, i.e. the contravariant commutative group-functor on the category of all $S$-schemes given by
\begin{equation*}
T \mapsto \mathrm{Pic}^0(X_T/T):=\textrm{\{Isomorphism classes of pairs} \ (\mathcal L, \alpha)\},
\end{equation*}
where $\mathcal L$ is a line bundle on $X_T$ which is algebraically equivalent to zero and where $\alpha$ is a $T$-rigidification of $\mathcal L$. The group law of $\mathrm{Pic}^0(X_T/T)$ is defined by taking the tensor product of representatives; its neutral element is the class of $(\mathcal O_{X_T}, \mathrm{can})$, and the inverse of the class of $(\mathcal L,\alpha)$ is represented by $(\mathcal L^{-1},\alpha^{-1})$. The assignment $T \mapsto \mathrm{Pic}^0(X_T/T)$ becomes contravariant functorial by means of Rem. 0.1.2 and Rem. 0.1.4.

\begin{lemma}
Let $\mathcal L$ be an arbitrary line bundle on $X_T$. Then any automorphism of $\mathcal L$ which restricts to the identity on $\epsilon_{X_T}^*\mathcal L$ is already the identity on $\mathcal L$.\\
In particular, there are no nontrivial automorphisms of a pair $(\mathcal L, \alpha) \in \mathrm{Pic}^0(X_T/T)$, i.e. if
\[(\mathcal L, \alpha) \simeq (\mathcal L, \alpha)\]
is an isomorphism, then it must be the identity.
\begin{proof}
It suffices to show only the first claim, which is done by a well-known standard argument (cf. e.g. \cite{Kl}, Lemma 9.2.10); in view of a later spot of the work we here recall it explicitly.\\
Namely, let $\varphi: \mathcal L \xrightarrow{\sim} \mathcal L$ be an automorphism of a line bundle $\mathcal L$ on $X_T$ with $\epsilon_{X_T}^*(\varphi)=\id$ on $\epsilon_{X_T}^*\mathcal L$. Note that $\varphi$ belongs to
\[\tag{\textbf{0.1.1}} \Gamma(T,\mathcal O_T) \xrightarrow{\sim} \Gamma(X_T,\mathcal O_{X_T}) \xrightarrow{\sim} \Hom_{\mathcal O_{X_T}}(\mathcal L, \mathcal L) \ni \varphi,\]
where the first arrow comes from the natural map $\mathcal O_T \rightarrow (\pi_{X_T})_*\mathcal O_{X_T}$ which is an isomorphism because $\mathcal O_S \xrightarrow{\sim} \pi_* \mathcal O_X$ holds universally for the abelian scheme $X/S$ (cf. \cite{Maz-Mes}, Ch. I, $(1.9)$). The second isomorphism in $(0.1.1)$ is defined by scalar multiplication.\\
The chain of identifications $(0.1.1)$ says that $\varphi$ is given by multiplication with a unit $u$ of $\Gamma(T,\mathcal O_T)$. To determine this unit observe that the composition
\[\Gamma(T,\mathcal O_T) \xrightarrow{\sim} \Gamma(X_T,\mathcal O_{X_T}) \rightarrow \Gamma(T,\mathcal O_T)\]
is the identity: here, the first arrow is as in $(0.1.1)$ and the second arrow comes from the map on structure sheaves defined by the zero section $\epsilon_{X_T}: T \rightarrow X_T$. This implies that the induced isomorphism of line bundles on $T$:
\[\epsilon_{X_T}^*(\varphi): \epsilon_{X_T}^*\mathcal L \xrightarrow{\sim} \epsilon_{X_T}^*\mathcal L\]
is given by $u$-multiplication. But by assumption $\epsilon_{X_T}^*(\varphi)=\id$, hence $u=1$ and the claim follows.
\end{proof}
\end{lemma}
We have the following deep theorem which can be deduced from a more general representability result in the theory of algebraic spaces (cf. the discussion in \cite{Ch-Fa}, Ch. I, p. 2-7).
\begin{theorem}
The dual functor of $X/S$
\begin{equation*}
T \mapsto \mathrm{Pic}^0(X_T/T)
\end{equation*}
is representable by an abelian scheme $Y$ of relative dimension $g$ over $S$.
\end{theorem}
\begin{definition}
The abelian $S$-scheme of Thm. 0.1.6 is called the \underline{dual abelian scheme of $X$}.\\
We denote by $\pi_Y, \ \mu_Y, \ \epsilon_Y $ the structure map, the multiplication map and the zero section of $Y/S$.
\end{definition}
The scheme $Y$ of Thm. 0.1.6 comes with a distinguished isomorphism class in $\mathrm{Pic}^0(X\times_S Y/Y)$ such that $Y$ together with this class forms a universal object for the dual functor. Each two representatives of this distinguished isomorphism class are uniquely isomorphic, as follows from Lemma 0.1.5.\\
From now on we fix such a representative and denote it by $(\mathcal P^0,r^0)$.
\begin{remark}
We explain the birigidification of $\mathcal P^0$.\\
By definition, $\mathcal P^0$ is a line bundle on the abelian $Y$-scheme $X \times_S Y$ which is algebraically equivalent to zero and $r^0$ is a $Y$-rigidification of $\mathcal P^0$, i.e. an isomorphism
\[r^0: \mathcal O_Y \simeq (\epsilon \times \id_Y)^*\mathcal P^0,\]
where $\epsilon \times \id_Y$ defines the zero section of $X\times_S Y/Y$:
\[Y \simeq S\times_S Y \xrightarrow{\epsilon \times \id_Y} X \times_S Y.\]
The commutative diagram of group homomorphisms
\begin{equation*}
\begin{xy}
\xymatrix@C-0.3cm{
\Hom_S(Y,Y) \ar@{-}[r]^{\sim \quad } \ar[d]_{\circ \epsilon_Y}  & \mathrm{Pic}^0(X \times_S Y/Y) \ar[d]^{(\mathrm{id}_X\times \epsilon_Y)^*}\\
\Hom_S(S,Y) \ar@{-}[r]^{ \ \sim} & \mathrm{Pic}^0(X/S)}
\end{xy}
\end{equation*}
and Lemma 0.1.5 show that there is a unique isomorphism $(\mathrm{id}_X \times \epsilon_Y)^*(\mathcal P^0,r^0) \simeq (\mathcal O_X,\mathrm{can})$ in $\mathrm{Pic}^0(X/S)$. In other words, there is a unique trivialization
\[s^0:(\id_X \times \epsilon_Y)^*\mathcal P^0 \simeq \mathcal O_X\]
of $\mathcal P^0$ along the map
\[X \simeq X\times_S S \xrightarrow{\mathrm{id}_X \times \epsilon_Y} X\times_S Y\]
such that the restriction of $s^0$ along $\epsilon: S \rightarrow X$ coincides with the restricition of $r^0$ along $\epsilon_Y: S \rightarrow Y$.\\
It is the existence of the two compatible rigidifications $r^0$ and $s^0$ which is meant by the common parlance that $\mathcal P^0$ is \underline{birigidified}, and we may write $(\mathcal P^0,r^0,s^0)$ to stress this fact.
\end{remark}
\begin{definition}
We call $(\mathcal P^0,r^0,s^0)$ the \underline{birigidified Poincaré bundle on $X \times_S Y$.}
\end{definition}
Finally, let us mention the phenomenon of \underline{biduality}:\\
For this one first recognizes that $\mathcal P^0$ is algebraically equivalent to zero not only with respect to the base scheme $Y$ - which is true by definition - but also with respect to $X$ (cf. \cite{SGA7-I}, exp. VII, Ex. 2.9.5 and Rem. 2.9.6; the substantial ingredient is the rigidity theorem for abelian schemes).\\
If we denote (only for a moment) by $Z$ the dual abelian scheme of $Y$ and by $\mathcal Q^0$ the birigidified Poincaré bundle on $Y \times_S Z$ there is then a unique $S$-morphism $\iota: X\rightarrow Z$ inducing an isomorphism
\[\tag{\textbf{0.1.2}} (\mathrm{id}_Y \times \iota)^*\mathcal Q^0 \simeq \sigma^*\mathcal P^0\]
such that the $Z$-rigidification of $\mathcal Q^0$ induces the $X$-rigidification $s^0$ of $\sigma^*\mathcal P^0$.\\
Here, $\sigma: Y \times_S X \xrightarrow{\sim} X \times_S Y$ denotes the shift automorphism and with the last $s^0$ we more precisely mean the $X$-rigidification of $\sigma^*\mathcal P^0$ naturally induced by $s^0$ noting that $\sigma \circ (\epsilon_Y \times \id_X)=\id_X \times \epsilon_Y$.\\
In shorter words, $\iota$ is the map corresponding to $(\sigma^*\mathcal P^0,s^0)$ under $\Hom_S(X,Z)\simeq \mathrm{Pic}^0(Y\times_S X/X)$.

\begin{theorem}
The map $\iota: X \rightarrow Z$ defined above is an isomorphism of abelian schemes.
\end{theorem}
\begin{proof}
Cf. \cite{Bo-Lü-Ray}, 8.4, Thm. 5 (b).\footnote{The hypothesis made there that $X/S$ is projective is of course not necessary for the argument: in the reference this assumption has the sole purpose to guarantee representability of the dual functor and thus is superfluous in view of Thm. 0.1.6.}
\end{proof}
It is easy to see that $(0.1.2)$ also respects the $Y$-rigidifications of both sides, and hence altogether we may identify $Z$ with $X$ as abelian schemes and $\mathcal Q^0$ with $\sigma^*\mathcal P^0$ as birigidified line bundles.
\begin{corollary}
The dual functor of $Y/S$ is represented by the abelian scheme $X$ and $(\sigma^*\mathcal P^0,s^0,r^0)$ is the birigidified Poincaré bundle on $Y \times_S X$. \qquad \qed
\end{corollary}

\subsubsection{Universal vectorial extension and Poincaré bundle}
In the following, we use the notion of an integrable connection and some associated standard constructions; for detailed explanations we refer to the review of algebraic connections given in 0.2.1.\\
\newline
Consider the contravariant commutative group-functor on the category of all $S$-schemes given by
\begin{equation*}
T \mapsto \mathrm{Pic}^\natural(X_T/T):=\textrm{\{Isomorphism classes of triples} \ (\mathcal L, \alpha, \nabla_{\mathcal L})\},
\end{equation*}
where $\mathcal L$ and $\alpha$ are as in the definition of $\mathrm{Pic}^0(X_T/T)$ and where $\nabla_{\mathcal L}: \mathcal L \rightarrow \Omega^1_{X_T/T} \otimes_{\mathcal O_{X_T}} \mathcal L$ is an integrable $T$-connection on $\mathcal L$. We call it the \underline{$\natural$-dual functor} of $X/S$.\\
For the precise definition of this functor one takes into account the following points:\\
First, an isomorphism between two such triples is by definition an isomorphism of the line bundles which respects the $T$-rigidifications and the connections.\\
Second, the group law of $\mathrm{Pic}^\natural(X_T/T)$ is defined by taking the tensor product of representing triples: this means that one forms the tensor product of the line bundles, of the $T$-rigidifications and of the integrable $T$-connections. The neutral element then is given by the class of $(\mathcal O_{X_T}, \mathrm{can}, \mathrm{d})$ with $\mathrm{d}: \mathcal O_{X_T} \rightarrow \Omega^1_{X_T/T}$ denoting the exterior derivative. The inverse of the class of $(\mathcal L, \alpha, \nabla_{\mathcal L})$ is represented by $(\mathcal L^{-1}, \alpha^{-1}, \nabla_{\mathcal L}^{-1})$, where $\nabla^{-1}_{\mathcal L}$ is the dual connection of $\nabla_{\mathcal L}$.\\
Finally, the assignment $T \mapsto \mathrm{Pic}^\natural(X_T/T)$ is contravariant functorial via Rem. 0.1.2, Rem. 0.1.4 and pullback of an integrable connection along the respective commutative diagram as in 0.2.1 (v).\\
\newline
From Lemma 0.1.5 it follows a fortiori that there are no nontrivial automorphisms of a triple $(\mathcal L, \alpha, \nabla_{\mathcal L})$.\\
\newline
In the course of the work we will exclusively be concerned with the case that the base scheme $S$ is of characteristic zero, i.e. a $\mathbb Q$-scheme. In this case we can remove the requirement of algebraic equivalence to zero in the definition of $\mathrm{Pic}^\natural(X_T/T)$ by the following more general statement:
\begin{lemma}
If $S$ is of characteristic zero, $T$ an $S$-scheme and $(\mathcal L,\nabla_{\mathcal L})$ is a line bundle on $X_T$ equipped with a (not necessarily integrable) $T$-connection, then $\mathcal L$ is algebraically equivalent to zero.
\end{lemma}
\begin{proof}
It is a fact that a line bundle on $X_T$ is algebraically equivalent to zero already if its restriction to each geometric fiber of $X_T/T$ is algebraically equivalent to zero.\footnote{There seems to be no place in the official literature where this is proven. It is nevertheless an easy formal consequence of the deeper fact that the dual functor of $X/S$ as we defined it is an open subfunctor of the relative Picard functor of $X/S$.}
We may hence check the claim after pullback of $\mathcal L$ to the fiber over a geometric point $\Spec(\Omega)$ of $T$. This pullback has an induced connection relative $\Spec(\Omega)$, given by the pullback of $\nabla_{\mathcal L}$ along the occurring fiber product diagram (cf. 0.2.1 (v)). We are thus reduced to show that a line bundle with connection on an abelian variety over an algebraically closed field of characteristic zero is algebraically equivalent to zero. This is well-known, cf. e.g. \cite{Brio}, Prop. 2.18.
\end{proof}

Before stating the key theorem about the $\natural$-dual functor we remark that a commutative group scheme over a base $B$ will often be freely identified with its associated abelian $fppf$-sheaf on the category of all $B$-schemes and that a sequence of homomorphisms of commutative group schemes is said to be exact if the corresponding sequence of abelian $fppf$-sheaves is exact. Moreover, for a commutative $B$-group scheme $G$ we write $\Lie(G/B)$ to denote the Lie algebra of $G/B$ which in most cases will be treated as the $\mathcal O_B$-module given by $(\epsilon_G^*\Omega^1_{G/B})^\vee$, where $\epsilon_G$ is the zero section of $G/B$; for its equivalent definition as a functor and for basic facts about Lie algebras cf. \cite{Li-Lo-Ray}, 1. Finally, for a $\mathcal O_B$-vector bundle $\mathcal E$ the notation $\mathbb V(\mathcal E)$ means the geometric vector bundle over $B$ associated with $\mathcal E$ (cf. \cite{Gö-We}, Ch. 11, $(11.4)$).\\
\newline
Returning to our fixed situation of an abelian scheme $X$ over a locally noetherian base $S$ we have the following fundamental theorem about the $\natural$-dual functor; it is a recapitulation in inverse order of the results proven in \cite{Maz-Mes}, Ch. I, §1-§4.

\begin{theorem}
The $\natural$-dual functor of $X/S$
\begin{equation*}
T \mapsto \mathrm{Pic}^\natural(X_T/T)
\end{equation*}
is representable by a $S$-group scheme $\Yr$ sitting in a natural short exact sequence of $S$-group schemes
\[\tag{\textbf{0.1.3}} 0 \rightarrow \mathbb V(\Lie(X/S)) \rightarrow \Yr \rightarrow Y \rightarrow 0,\]
given in $T$-rational points (for an $S$-scheme $T$) by the exact sequence of abelian groups
\[\tag{\textbf{0.1.4}}0 \rightarrow H^0(X_T,\Omega^1_{X_T/T}) \rightarrow \mathrm{Pic}^\natural(X_T/T) \rightarrow \mathrm{Pic}^0(X_T/T),\]
where the injection maps a form $\omega$ to the class of $(\mathcal O_{X_T},\mathrm{can},  \mathrm{d}+\omega)$ and the subsequent arrow is defined by "forgetting the connection", i.e. by  $(\mathcal L, \alpha,\nabla_{\mathcal L}) \mapsto (\mathcal L, \alpha)$ on representatives.\\
The scheme $\Yr$ in particular is of finite type, separated and smooth of relative dimension $2g$ over $S$ with geometrically integral fibers.\\
Moreover, there exists a canonical isomorphism of $\mathcal O_S$-vector bundles
\[\tag{\textbf{0.1.5}} \Lie(\Yr/S) \simeq H^1_{\mathrm{dR}}(X/S)\]
such that the sequence of Lie algebras of $(0.1.3)$ identifies with the exact sequence ("Hodge filtration")
\[0\rightarrow \pi_*\Omega^1_{X/S}\rightarrow H^1_{\mathrm{dR}}(X/S) \rightarrow R^1\pi_*\mathcal O_X \rightarrow 0\]
induced by the degeneration of the Hodge-de Rham spectral sequence for $X/S$ on the first sheet.\footnote{General remarks about algebraic de Rham cohomology can be found in 0.2.2. For detailed information about the de Rham cohomology of an abelian scheme cf. \cite{Bert-Br-Mes}, 2.5.}\\
The exact sequence $(0.1.3)$ parametrizes all extensions of $Y$ by $S$-vector groups $W$ via pushout along a unique $S$-vector group homomorphism $\mathbb V(\Lie(X/S)) \rightarrow W$.
\end{theorem}

\begin{definition}
The $S$-group scheme $\Yr$ of Thm. 0.1.13 is called the \underline{universal (vectorial) extension of $Y$}.\\
We denote by $\pi^{\natural} , \ \mu^{\natural}, \ \epsilon^{\natural}, \ (-1)^\natural$ the structure map, the multiplication map, the zero section and the inverse map of $\Yr/S$.
\end{definition}

\begin{remark}
(i) Of course, the term "universal vectorial extension" refers to the property of the extension $(0.1.3)$ described in the last part of Thm. 0.1.13. This property is used to define more generally the universal vectorial extension of an arbitrary $1$-motive over $S$ (cf. 2.6.1). Its existence is always guaranteed by in fact purely formal arguments, and in the case of the $1$-motive given by $Y$ the preceding theorem thus provides us with an explicit interpretation of this object in terms of the $\natural$-dual functor of $X/S$.\vspace{1mm}\\
(ii) There are further equivalent viewpoints on $\Yr$ by means of rigidified extensions (cf. \cite{Maz-Mes}, Ch. I, §2) or by Grothendieck's $\natural$-extensions (cf. ibid., Ch. I, §3-§4 resp. the account given in 0.1.3). For its classical expression via differential forms of the third kind in the case of an abelian variety as well as for its interpretation within the theory of generalized Picard varieties resp. jacobians cf. the references in the introduction of \cite{Maz-Mes}, \cite{Col}, I.7, and \cite{Co4}.
\end{remark}

$\Yr$ comes with a distinguished isomorphism class in $\mathrm{Pic}^\natural(X \times_S \Yr/\Yr)$ such that $\Yr$ together with this class forms a universal object for the $\natural$-dual functor. We now want to fix a representative for this universal isomorphism class by using the $Y$-rigidified Poincaré bundle $(\mathcal P^0,r^0)$ on $X\times_S Y$.\\
Write $(\mathcal P,r)$ for the pullback of $(\mathcal P^0,r^0)$ along the map $X\times_S \Yr\rightarrow X\times_S Y$ induced by the canonical arrow $\Yr \rightarrow Y$ of $(0.1.3)$: it is a $\Yr$-rigidified line bundle on the abelian $\Yr$-scheme $X\times_S \Yr$ which is algebraically equivalent to zero. We then have:
\begin{lemma}
There is a unique integrable $\Yr$-connection $\nabla_{\mathcal P}$ on $\mathcal P$ such that $(\mathcal P,r,\nabla_{\mathcal P})$ represents the universal isomorphism class in $\mathrm{Pic}^\natural(X\times_S \Yr/\Yr)$.
\end{lemma}
\begin{proof}
Fix some representative $(\mathcal Q,\nu, \nabla_\mathcal Q)$ of this universal class. By definition of $Y$ there is a unique $S$-morphism $f:\Yr \rightarrow Y$ such that the pullback of $(\mathcal P^0,r^0)$ under $\id_X \times f: X \times_S \Yr \rightarrow X\times_S Y$ is isomorphic to $(\mathcal Q,\nu)$. We first claim that $f$ is identical to the canonical map in $(0.1.3)$: by Thm. 0.1.13 we need to check that for all $S$-schemes $T$ the map $f$ in $T$-rational points
\[f(T): \mathrm{Pic}^\natural(X\times_S T/T) \simeq \Yr(T) \rightarrow Y(T) \simeq \mathrm{Pic}^0(X\times_S T/T)\]
is given as in $(0.1.4)$, i.e. "by forgetting the connection". This in turn is straightforwardly seen by using the universality of $(\mathcal Q,\nu, \nabla_\mathcal Q)$ and $(\mathcal P^0,r^0)$ together with $(\mathrm{id}_X \times f)^*(\mathcal P^0,r^0) \simeq (\mathcal Q,\nu)$.\\
We thus have an isomorphism
\[\tag{\textbf{0.1.6}} (\mathcal P,r)=(\mathrm{id}_X\times f)^*(\mathcal P^0,r^0) \simeq (\mathcal Q,\nu).\]
Now we endow the left side of $(0.1.6)$ with the connection $\nabla_{\mathcal P}$ induced by $\nabla_\mathcal Q$. Clearly, $(\mathcal P,r,\nabla_{\mathcal P})$ is a representative of the universal class in $\mathrm{Pic}^\natural(X\times_S \Yr/\Yr)$ because $(\mathcal Q,\nu, \nabla_\mathcal Q)$ is. Finally, the uniqueness part is an application of Lemma 0.1.5 (note that $(\mathcal P,r)$ is assumed to be fixed).
\end{proof}
\begin{remark}
The commutative diagram of group homomorphisms
\begin{equation*}
\begin{xy}
\xymatrix@C-0.3cm{
\Hom_S(\Yr,\Yr) \ar@{-}[r]^{\sim \quad } \ar[d]_{\circ \epsilon^\natural}  & \mathrm{Pic}^\natural(X \times_S \Yr/\Yr) \ar[d]^{(\mathrm{id}_X\times \epsilon^\natural)^*}\\
\Hom_S(S,\Yr) \ar@{-}[r]^{ \ \sim} & \mathrm{Pic}^\natural(X/S)}
\end{xy}
\end{equation*}
and Lemma 0.1.5 show that there is a unique isomorphism $(\mathrm{id}_X \times \epsilon^\natural)^*(\mathcal P,r,\nabla_{\mathcal P}) \simeq (\mathcal O_X,\mathrm{can},\mathrm{d})$ in $\mathrm{Pic}^\natural(X/S)$. In other words, there is a unique trivialization
\[s:(\mathrm{id}_X \times \epsilon^\natural)^*(\mathcal P,\nabla_{\mathcal P}) \simeq (\mathcal O_X,\mathrm{d})\]
of $(\mathcal P,\nabla_{\mathcal P})$ along
\[X \simeq X\times_S S \xrightarrow{\mathrm{id}_X \times \epsilon^\natural} X\times_S \Yr\]
with the property that (on the level of line bundles) the restriction of $s$ along $\epsilon: S \rightarrow X$ coincides with the restriction of $r$ along $\epsilon^\natural: S \rightarrow \Yr$. One can further check that the trivialization $s$ on the level of line bundles is just the one induced by $s^0$ in the obvious sense.\footnote{The reason is that both give an isomorphism of $(\mathrm{id}_X \times \epsilon^\natural)^*(\mathcal P,r)$ with $(\mathcal O_X,\mathrm{can})$: for $s$ this holds by definition and for the trivialization induced by $s^0$ because $s^0$ and $r^0$ are compatible (cf. Rem. 0.1.8) and $r$ is induced by $r^0$. Now use Lemma 0.1.5 to conclude the desired equality.}\\
We may express these data simultaneously by writing $(\mathcal P,r,s,\nabla_{\mathcal P})$.
\end{remark}
\begin{definition}
We call $(\mathcal P,r,s,\nabla_{\mathcal P})$ the \underline{birigidified Poincaré bundle with universal integrable $\Yr$-connection}\\
\underline{on $X\times_S \Yr$.}
\end{definition}
It is of course reasonable to speak about rigidifying a line bundle also along the zero section of a group scheme, and in this sense we apply the term "rigidification" to $s$.\\
\newline
The reason why we decide to write $(\mathcal P,r,s,\nabla_{\mathcal P})$ instead of the more consequential $(P^\natural, r^\natural,s^\natural, \nabla_{\mathcal P^{\natural}})$ is for convenience of notation: playing a crucial role in this work, these objects will occur a huge number of times and will moreover need to be furnished with other super- and subscripts.

\subsection{Extension and biextension structures: generalities}
The language of $\mathbb G_m$-(bi-)extensions and $\natural$-structures on them permits to develop another perspective on the objects introduced in the preceding subsection. We therefore
insert a review of some of its basic vocabulary, essentially summarizing material from the sources \cite{SGA7-I}, exp. VII, \cite{De2}, 10.2, \cite{Ber}, 3, and \cite{Po}, II, 10.2-10.3.

\subsubsection{Extensions and $\natural$-extensions}

We fix an arbitrary base scheme $T$ and a commutative $T$-group scheme $H$ with multiplication map $m: H\times_T H \rightarrow H$ and projections $p_1, p_2: H\times_T H \rightarrow H$.\\
Let us further write $q_1,q_2,q_3: H\times_T H \times_T H \rightarrow H$ for the three projections of the triple product.
\begin{definition}
A \underline{(commutative) $\mathbb G_{m,T}$-extension of $H$} is a line bundle $\mathcal L$ on $H$ together with an isomorphism of line bundles on $H\times_T H$:
\[\tag{\textbf{0.1.7}} m^*\mathcal L \simeq p_1^*\mathcal L \otimes_{\mathcal O_{H\times_T H}} p_2^*\mathcal L\]
such that the following two diagrams of isomorphisms of line bundles on $H\times_T H \times_T H$ resp. on $H\times_T H$ (expressing "associativity" and "commutativity") commute; note that we leave away the index $H\times_T H \times_T H$ resp. $H \times_T H$ in the tensor products.
\begin{equation*} \tag{\textbf{0.1.8}} \begin{split}
\begin{xy}
\xymatrix{
(q_1+q_2+q_3)^*\mathcal L \ar[r]^{(a) \ \ } \ar[d]_{(b)}  & (q_1+q_2)^*\mathcal L \otimes q_3^*\mathcal L \ar[d]^{(d) \otimes \id}\\
q_1^*\mathcal L \otimes (q_2+q_3)^*\mathcal L \ar[r]^{\id\otimes (c)} & q_1^*\mathcal L \otimes q_2^*\mathcal L \otimes q_3^*\mathcal L}
\end{xy}
\end{split}
\end{equation*}
resp.
\begin{equation*} \tag{\textbf{0.1.9}} \begin{split}
\begin{xy}
\xymatrix{
m^*\mathcal L \ar[r]\ar[d]_{\id}  & p_1^*\mathcal L \otimes p_2^*\mathcal L \ar[d]^{\mathrm{can}}\\
m^*\mathcal L \ar[r]& p_2^*\mathcal L \otimes p_1^*\mathcal L}
\end{xy}
\end{split}
\end{equation*}
The arrows in $(0.1.8)$ are induced by pulling back the isomorphism $(0.1.7)$ along the maps\\
$H\times_T H \times_T H \rightarrow H\times_T H$ defined by
\[\begin{split} (a')&: \ (x,y,z) \mapsto (x+y,z) \quad (b'): \ (x,y,z) \mapsto (x,y+z)\\
(c') &: \ (x,y,z) \mapsto (y,z) \qquad \ \ \ (d'): \ (x,y,z) \mapsto (x,y). \qquad \end{split}\]
\end{definition}

The upper arrow of $(0.1.9)$ is $(0.1.7)$ and the lower one comes from $(0.1.7)$ by pullback via the shift map $H\times_T H \rightarrow H \times_T H$, defined by $(x,y) \mapsto (y,x)$.

\begin{remark}
Let $\mathcal L$ be a $\mathbb G_{m,T}$-extension of $H$ and write $L$ for the $\mathbb G_{m,H}$-torsor associated with the line bundle $\mathcal L$. By a purely formal procedure the datum $(0.1.7)$ induces on $L$ canonically the structure of a commutative $T$-group scheme fitting into an exact sequence of $T$-group schemes
\[\tag{\textbf{0.1.10}} 0\rightarrow \mathbb G_{m,T} \rightarrow L \rightarrow H \rightarrow 0.\]
Conversely, an extension as in $(0.1.10)$ makes $L$ into a $\mathbb G_{m,H}$-torsor whose associated line bundle $\mathcal L$ on $H$ naturally carries the structure of a $\mathbb G_{m,T}$-extension of $H$.\\
For the explicit (entirely formal) constructions and for more details cf. \cite{SGA7-I}, exp. VII, 1.1-1.2.
\end{remark}

Keeping the above notations let us  further assume that $H$ is smooth over $T$.

\begin{definition}
A \underline{$\natural$-$\mathbb G_{m,T}$-extension of $H$} is a $\mathbb G_{m,T}$-extension $\mathcal L$ of $H$ together with a (automatically integrable) $T$-connection on $\mathcal L$ such that the isomorphism of $(0.1.7)$:
\[m^*\mathcal L \simeq p_1^*\mathcal L \otimes_{\mathcal O_{H\times_T H}} p_2^*\mathcal L\]
is horizontal for the induced $T$-connections on both sides.\\
Write $\Ext^\natural(H,\mathbb G_{m,T})$ for the group of isomorphism classes of $\natural$-$\mathbb G_{m,T}$-extensions of $H$, where an isomorphism is of course an isomorphism of the line bundles which is compatible with the respective isomorphisms $(0.1.7)$ and with the connections; the group structure is induced by taking the tensor product of line bundles, connections and isomorphisms $(0.1.7)$.
\end{definition}

\subsubsection{Biextensions and $\natural$-structures}
We fix a base scheme $T$ and two commutative $T$-group schemes $H_1,H_2$ with multiplication maps $m_{H_1}, m_{H_2}$. We write $p_{ij}$ for the projection of a triple product to the $i$-th and $j$-th factor.
\begin{definition}
A \underline{(commutative) $\mathbb G_{m,T}$-biextension of $H_1\times_T H_2$} is a line bundle $\mathcal L$ on $H_1\times_T H_2$ together with isomorphisms of line bundles on $H_1\times_T H_1 \times_T H_2$ and $H_1\times_T H_2\times_T H_2$:
\[\tag{\textbf{0.1.11}} p_{13}^*\mathcal L \otimes p_{23}^*\mathcal L \simeq (m_{H_1} \times \mathrm{id}_{H_2})^*\mathcal L\]
and
\[\tag{\textbf{0.1.12}} p_{12}^*\mathcal L \otimes p_{13}^*\mathcal L \simeq (\mathrm{id}_{H_1} \times m_{H_2})^*\mathcal L\]
such that certain five induced diagrams of isomorphisms of line bundles on the fiber products
\[\bullet \ H_1 \times_T H_1 \times_T H_1 \times_T H_2, \quad H_1\times_T H_1 \times_T H_2,\]
\[ \bullet \ H_1 \times_T H_2 \times_T H_2 \times_T H_2, \quad H_1 \times_T H_2 \times_T H_2,\]
\[\bullet \ H_1 \times_T H_1 \times_T H_2 \times_T H_2 \qquad \qquad \qquad \qquad \quad\]
are commutative.\footnote{These diagrams are intuitive to write down but would cost further notation and space, so we don't explicate them here and refer to the diagrams $(2.0.5)$, $(2.0.6)$, $(2.0.8)$, $(2.0.9)$, $(2.1.1)$ in \cite{SGA7-I}, exp. VII, 2.0 and 2.1, for the precise requirements.}
\end{definition}

\begin{remark}
Write $H_{1_{H_2}}$ resp. $H_{2_{H_1}}$ for the product $H_1\times_T H_2$, considered as group scheme relative $H_2$ resp. $H_1$.\\
If $\mathcal L$ is a $\mathbb G_{m,T}$-biextension of $H_1\times_T H_2$, then $\mathcal L$ has canonically an induced structure as a $\mathbb G_{m,H_2}$-extension of $H_{1_{H_2}}$ and as a $\mathbb G_{m,H_1}$-extension of $H_{2_{H_1}}$. One gets these structures precisely from the isomorphisms $(0.1.11)$ and $(0.1.12)$ (the two associativity and commutativity conditions are exactly the meaning of the first four commutative diagrams of Def. 0.1.22). They are compatible in the sense of \cite{SGA7-I}, exp. VII, Def. 2.1 (this is the meaning of the fifth commutative diagram of Def. 0.1.22).\\
One thus sees that the datum of such two compatible extension structures for a line bundle $\mathcal L$ on $H_1\times_T H_2$ is simply tantamount to have on $\mathcal L$ the structure of a $\mathbb G_{m,T}$-biextension of $H_1\times_T H_2$.\\
Note that with Rem. 0.1.20 we then obtain exact sequences of group schemes relative $H_2$ resp. $H_1$:
\[\tag{\textbf{0.1.13}} 0 \rightarrow \mathbb G_{m,H_2} \rightarrow L \rightarrow H_{1_{H_2}} \rightarrow 0,\]
\[\tag{\textbf{0.1.14}}0 \rightarrow \mathbb G_{m,H_1} \rightarrow L \rightarrow H_{2_{H_1}} \rightarrow 0.\]
\end{remark}

Keeping the above notations we further assume that $H_1$ and $H_2$ are smooth over $T$.

\begin{definition}
Let $\mathcal L$ be a $\mathbb G_{m,T}$-biextension of $H_1\times_T H_2$.\\
(i) A \underline{$\natural$-$1$-structure} on $\mathcal L$ is the datum of a (automatically integrable) $H_2$-connection on $\mathcal L$ such that $(0.1.11)$ is horizontal for the induced $H_2$-connections (i.e. $\mathcal L$ becomes a $\natural$-$\mathbb G_{m,H_2}$-extension of $H_{1_{H_2}}$) and such that $(0.1.12)$ is horizontal for the induced $(H_2\times_T H_2)$-connections.\\
(ii) By working in (i) relative to $H_1$ instead of $H_2$ one obtains the notion of a \underline{$\natural$-$2$-structure} on $\mathcal L$.\\
(iii) A \underline{$\natural$-structure} on $\mathcal L$ is the datum of a $\natural$-$1$-structure and of a $\natural$-$2$-structure on $\mathcal L$.
\end{definition}
By the canonical decomposition $\Omega^1_{H_1\times_T H_2/H_2}\oplus \Omega^1_{H_1\times_T H_2/H_1} \simeq \Omega^1_{H_1\times_T H_2/T}$ one easily sees that a  $\natural$-structure on $\mathcal L$ is equivalent to the datum of a $T$-connection on $\mathcal L$ such that $(0.1.11)$ and $(0.1.12)$ are horizontal for the induced $T$-connections.

\subsection{Extension and biextension structures: applications}
Let us return to the objects we have defined in 0.1.1 starting from the abelian scheme $X$ over the locally noetherian base $S$. The terminology developed in 0.1.2 will now provide us with a different interpretation of the dual abelian scheme $Y$ and its universal vectorial extension $\Yr$ as parametrizing schemes for the $\mathbb G_m$- and $\natural$-$\mathbb G_m$-extensions of $X$. We finally apply the concept of biextension to the Poincaré bundle and, following \cite{Ber}, 4, resp. \cite{De2}, 10.2, explain the construction of Deligne's pairing for the first de Rham cohomology of $X$ and $Y$.\\
All of these viewpoints play a particularly basic and natural role when working within the framework of $1$-motives, which we will do extensively in 2.6.

\subsubsection{Extensions and the dual abelian scheme}
Let $(\mathcal L, \alpha)$ be a line bundle on $X_T$ which is algebraically equivalent to zero together with a $T$-rigidification. The existence of this $T$-rigidification easily implies that any line bundle $\mathcal M$ on $T$ as in Def. 0.0.1 must be trivial, and hence we obtain some isomorphism of line bundles on $X_T \times_T X_T$:
\[\tag{\textbf{0.1.15}} \mu_T^* \mathcal L \simeq \mathrm{pr}_{1,T}^* \mathcal L \otimes_{\mathcal O_{X_T\times_T X_T}} \mathrm{pr}_{2,T}^* \mathcal L.\]
If we view $X_T\times_T X_T$ as abelian $T$-scheme both sides of $(0.1.15)$ carry a natural $T$-rigidification (induced by $\alpha$ in the obvious way), and there is a unique isomorphism as in $(0.1.15)$ respecting them.\footnote{This is straightforward: first, choose any isomorphism as in $(0.1.15)$. Then it can be changed by a constant in $\Gamma(T,\mathcal O_T)^*=\Gamma(X_T,\mathcal O_{X_T})^*$ such that it becomes compatible with the rigidifications: namely, restrict the chosen isomorphism $(0.1.15)$ along the zero section of $X_T\times_T X_T$ and use the rigidifications to obtain an automorphism of $\mathcal O_T$ which is given by an element of $\Gamma(T,\mathcal O_T)^*$. The inverse of this element is the desired constant, as is easily checked. Finally, that an isomorphism as in $(0.1.15)$ which respects the rigidifications is unique follows from Lemma 0.1.5 applied to the abelian $T$-scheme $X_T\times_T X_T$.}

\begin{lemma}
With this choice of $(0.1.15)$ we obtain on $\mathcal L$ the structure of a (commutative) $\mathbb G_{m,T}$-extension of $X_T$.
\end{lemma}
\begin{proof}
We need to show that the diagrams $(0.1.8)$ resp. $(0.1.9)$ expressing associativity and commutativity are commutative. All line bundles occurring in $(0.1.8)$ resp. $(0.1.9)$ carry a natural $T$-rigidification (induced by $\alpha$), and by the choice of $(0.1.15)$ and the definition of the arrows in $(0.1.8)$ resp. $(0.1.9)$ these rigidifications are respected by all isomorphisms in $(0.1.8)$ resp. $(0.1.9)$. Now use Lemma 0.1.5, applied to the abelian $T$-scheme $X_T\times_T X_T \times_T X_T$ resp. $X_T\times_T X_T$.
\end{proof}
Observing moreover Rem. 0.1.20 we see that the pair $(\mathcal L,\alpha)$ thus yields in a natural way an extension
\[0 \rightarrow \mathbb G_{m,T} \rightarrow L \rightarrow X_T \rightarrow 0,\]
where $L$ stands as usual for the $\mathbb G_{m,X_T}$-torsor associated with the line bundle $\mathcal L$.\\
Conversely, given an exact sequence of abelian $fppf$-sheaves
\[0 \rightarrow \mathbb G_{m,T} \rightarrow L \rightarrow X_T \rightarrow 0\]
write $\mathcal L$ for the line bundle on $X_T$ corresponding to $L$; then $\mathcal L$ naturally becomes a $\mathbb G_{m,T}$-extension of $X_T$ (cf. Rem. 0.1.20), and the associated isomorphism $(0.1.7)$ implies that $\mathcal L$ is algebraically equivalent to zero; moreover, restricting $(0.1.7)$ along the zero section of $X_T\times_T X_T$ induces a $T$-rigidification $\alpha$ of $\mathcal L$ such that in sum we have obtained a pair $(\mathcal L,\alpha)$ as before.
\begin{theorem}["Barsotti-Rosenlicht-Weil formula"]
The assignments described above induce mutually inverse isomorphisms of abelian groups
\[\tag{\textbf{0.1.16}} \mathrm{Pic}^0(X_T/T) \simeq \Ext^1_{fppf}(X_T, \mathbb G_{m,T})\]
which is functorial in $T$. In other words, the dual abelian scheme $Y$ of $X/S$ represents the functor on the category of $S$-schemes
\[T \mapsto \Ext^1_{fppf}(X_T, \mathbb G_{m,T}).\]
In particular, we obtain an isomorphism of abelian $fppf$-sheaves
\[\tag{\textbf{0.1.17}} Y \simeq \underline{\Ext}^1_{fppf}(X,\mathbb G_{m,S}).\]
\end{theorem}
\begin{proof}
That we obtain an induced isomorphism of groups $\Ext^1_{fppf}(X_T,\mathbb G_{m,T}) \xrightarrow {\sim} \mathrm{Pic}^0(X_T/T)$ follows from \cite{Oo}, Ch. III, Thm. 18.1; observe that the additional hypotheses made there can be completely removed: cf. the discussion in \cite{Jo2}, p. 7, footnote 1.\\
It is easy to see that assigning to a pair $(\mathcal L,\alpha)$ an $fppf$-extension of $X_T$ by $\mathbb G_{m,T}$ as outlined above induces a well-defined map $\mathrm{Pic}^0(X_T/T) \rightarrow \Ext^1_{fppf}(X_T,\mathbb G_{m,T})$ which becomes the identity on $\mathrm{Pic}^0(X_T/T)$ when further composed with the previous isomorphism of groups. This suffices to establish the isomorphism $(0.1.16)$. The remaining claims are clear.
\end{proof}

\subsubsection{$\natural$-extensions and the universal vectorial extension}
Recall that representatives of the abelian group $\Ext^\natural(X_T,\mathbb G_{m,T})$ are $\mathbb G_{m,T}$-extensions $\mathcal L$ of $X_T$ together with an integrable $T$-connection $\nabla_{\mathcal L}$ such that the isomorphism defining the extension structure
\[\mu_T^* \mathcal L \simeq \mathrm{pr}_{1,T}^* \mathcal L \otimes_{\mathcal O_{X_T\times_T X_T}} \mathrm{pr}_{2,T}^* \mathcal L\]
is horizontal for the induced integrable $T$-connections. Pullback of the preceding isomorphism along the zero section of $X_T \times_T X_T$ induces in the natural way a $T$-rigidification $\alpha$ of $\mathcal L$, and the triple $(\mathcal L, \alpha, \nabla_{\mathcal L})$ yields a well-defined class in $\mathrm{Pic}^\natural(X_T/T)$.\\
In this way we obtain a homomorphism
\[\Ext^\natural(X_T,\mathbb G_{m,T}) \rightarrow \mathrm{Pic}^\natural(X_T/T),\]
and an application of Lemma 0.1.5 (similar as in the proof of Lemma 0.1.25) shows its injectivity.\\
It is in fact also surjective: if a class $(\mathcal L,\alpha, \nabla_{\mathcal L})$ of $\mathrm{Pic}^\natural(X_T/T)$ is given, then the construction and isomorphy of $(0.1.16)$ implies that there is a $\mathbb G_{m,T}$-extension of $X_T$ which - together with its induced $T$-rigidification - is (uniquely) isomorphic to $(\mathcal L,\alpha)$. Now, using this isomorphism, endow the line bundle of that extension with the connection induced by $\nabla_{\mathcal L}$. One can then show that the isomorphism defining the extension structure is in fact horizontal (cf. the argument in $\cite{Maz-Mes}$, Ch. I, proof of Prop. 4.2.1), which thus shows the desired surjectivity.\\
\newline
In sum, this establishes the following result (cf. \cite{Maz-Mes}, Ch. I, $(4.2)$):
\begin{theorem}
The universal vectorial extension $\Yr$ of $Y$ represents the functor on the category of $S$-schemes
\[T \mapsto \Ext^\natural(X_T, \mathbb G_{m,T}).\] \qquad \qed
\end{theorem}

\subsubsection{Applications to the Poincaré bundle}
Consider the birigidified Poincaré bundle $(\mathcal P^0,r^0,s^0)$ on $X\times_S Y$ and recall that it is compatibly rigidified and algebraically equivalent to zero with respect to each of the two factors of $X\times_S Y$ (cf. Rem. 0.1.8, Def. 0.1.9 and the subsequent discussion). Applying two times the construction of Lemma 0.1.25 (the second time with $X$ replaced by $Y$) we obtain on $\mathcal P^0$ naturally  the structure of a $\mathbb G_{m,Y}$- and of a $\mathbb G_{m,X}$-extension of $X \times_S Y$. Denoting by $P^0$ the $\mathbb G_{m,X\times_S Y}$-torsor associated with $\mathcal P^0$ we may express this (by Rem. 0.1.23) via associated exact sequences
\[\tag{\textbf{0.1.18}} 0 \rightarrow \mathbb G_{m,Y} \rightarrow P^0 \rightarrow X_Y\rightarrow 0,\]
\[\tag{\textbf{0.1.19}} 0 \rightarrow \mathbb G_{m,X} \rightarrow P^0 \rightarrow Y_X\rightarrow 0.\]
Their extension classes are the universal ones under the isomorphisms of functors given by $(0.1.17)$:
\[Y \simeq \underline{\Ext}^1_{fppf}(X,\mathbb G_{m,S}),\]
\[X \simeq \underline{\Ext}^1_{fppf}(Y,\mathbb G_{m,S}),\]
where for the second we have made use of biduality (cf. Cor. 0.1.11).\\
The compatibility of $r^0$ and $s^0$ allow to show (as in the proof of Lemma 0.1.25) that the two extension structures on $\mathcal P^0$ are compatible in the sense of Rem. 0.1.23, and thus $\mathcal P^0$ is a $\mathbb G_{m,S}$-biextension of $X\times_S Y$. It is often referred to as the \underline{canonical or universal $\mathbb G_{m,S}$-biextension of $X\times_S Y$}.\\
\newline
We have furthermore introduced the birigidified Poincaré bundle with universal integrable $\Yr$-connection $(\mathcal P,r,s,\nabla_\mathcal P)$ on $X\times_S \Yr$ (cf. Def. 0.1.18). Pullback of the $\mathbb G_{m,S}$-biextension structure of $\mathcal P^0$ makes $\mathcal P$ a $\mathbb G_{m,S}$-biextension of $X\times_S \Yr$ whose corresponding isomorphisms write as
\[\tag{\textbf{0.1.20}} p_{13}^*\mathcal P \otimes p_{23}^*\mathcal P \simeq (\mu \times \mathrm{id}_{\Yr})^*\mathcal P \quad \textrm{on} \ X\times_S X \times_S \Yr,\]
\[\tag{\textbf{0.1.21}} p_{12}^*\mathcal P \otimes p_{13}^*\mathcal P \simeq (\mathrm{id}_X \times \mu^{\natural})^*\mathcal P \quad \textrm{on} \ X\times_S \Yr \times_S \Yr.\]
One checks that $(0.1.20)$ is horizontal for the integrable $\Yr$-connections induced by $\nabla_\mathcal P$ on both sides (use the argument of \cite{Maz-Mes}, Ch. I, proof of Prop. 4.2.1, applied to the abelian $\Yr$-scheme $X\times_S \Yr$); in particular, $(\mathcal P,\nabla_\mathcal P)$ becomes what in Def. 0.1.21 we have called a $\natural$-$\mathbb G_{m,\Yr}$-extension of $X\times_S \Yr$, and its class in $\Ext^\natural(X\times_S \Yr,\mathbb G_{m,\Yr})$ is the universal one with respect to Thm. 0.1.27. Moreover, the isomorphism $(0.1.21)$ is horizontal for the integrable $(\Yr\times_S \Yr)$-connections induced by $\nabla_\mathcal P$ on both sides (cf. \cite{Ber}, proof of Prop. 3.9, or \cite{De2}, proof of Prop. $(10.2.7.4)$).\\
In the terminology of Def. 0.1.24 the horizontality of the two isomorphisms $(0.1.20)$ and $(0.1.21)$ defining the $\mathbb G_{m,S}$-biextension structure of $\mathcal P$ says that $\nabla_\mathcal P$ induces on $\mathcal P$ a $\natural$-$1$-structure.\\
Finally, note that by Rem. 0.1.23 the extension structures on $\mathcal P$ provide us with exact sequences
\[\tag{\textbf{0.1.22}} 0 \rightarrow \mathbb G_{m,\Yr} \rightarrow P \rightarrow X_{\Yr}\rightarrow 0,\]
\[\tag{\textbf{0.1.23}} 0 \rightarrow \mathbb G_{m,X} \rightarrow P \rightarrow \Yr_X\rightarrow 0.\]

\subsubsection{Deligne's pairing}
The objects introduced in this section can be used to construct a natural duality between the first de Rham cohomology sheaves of an abelian scheme and its dual.\\
Let us denote by $X^\natural$ the universal vectorial extension of $X$.\\
Observing biduality (cf. Cor. 0.1.11) and then proceeding analogously as for $\mathcal P$ (cf. Lemma 0.1.16 and the above arguments) one obtains that the pullback of $\mathcal P^0$ along the canonical homomorphism $X^\natural \times_S Y \rightarrow X \times_S Y$ is a birigidified $\mathbb G_{m,S}$-biextension of $X^\natural \times_S Y$ which carries a distinguished integrable $X^\natural$-connection defining a $\natural$-$2$-structure on it. One can then take its further pullback along $X^\natural\times_S \Yr \rightarrow X^\natural \times_S Y$ and perform the analogous pullback for $\mathcal P$ and its $\natural$-$1$-structure.\\
Hence, writing $^\natural\mathcal P^\natural$ for the pullback of $\mathcal P^0$ along $X^\natural\times_S Y^\natural \rightarrow X\times_S Y$ we obtain a $\natural$-structure on ${^\natural}\mathcal P^\natural$ which expresses as a $S$-connection
\[\nabla_{{^\natural}\mathcal P^\natural}: {^\natural}\mathcal P^\natural \rightarrow \Omega^1_{X^\natural \times_S \Yr/S} \otimes_{\mathcal O_{X^\natural\times_S \Yr}} {^\natural}\mathcal P^\natural.\]
The global $2$-form on $X^\natural \times_S \Yr$ defined by the curvature of $\nabla_{{^\natural}\mathcal P^\natural}$ is checked to be invariant (cf. \cite{Co3}, p. 636), thus giving rise to an alternating $\mathcal O_S$-bilinear form
\[\Lie(X^\natural \times_S \Yr/S) \oplus \Lie(X^\natural \times_S \Yr/S) \rightarrow \mathcal O_S\]
which in view of the identifications provided by $(0.1.5)$ writes as
\[R:(H^1_{\mathrm{dR}}(Y/S) \oplus H^1_{\mathrm{dR}}(X/S)) \oplus (H^1_{\mathrm{dR}}(Y/S)\oplus H^1_{\mathrm{dR}}(X/S)) \rightarrow \mathcal O_S.\]
The integrability of the relative connections summing to $\nabla_{{^\natural}\mathcal P^\natural}$ implies that $R$ vanishes on $H^1_{\mathrm{dR}}(Y/S)$ and on $H^1_{\mathrm{dR}}(X/S)$. It is hence of the form
\[R((v,v'), (w,w'))=\Phi(v \otimes w')-\Phi(w \otimes v')\]
with a $\mathcal O_S$-linear pairing
\[\Phi:H^1_{\mathrm{dR}}(Y/S) \otimes H^1_{\mathrm{dR}}(X/S) \rightarrow \mathcal O_S\]
which in terms of $R$ is given by the formula
\[\Phi(v\otimes w)=R((v,0), (0,w)).\]
The pairing $\Phi$ is called \underline{Deligne's pairing}. That it is in fact perfect is a fundamental result which was proven under additional assumptions on the base scheme $S$ by Deligne and Coleman (cf. \cite{De2} and \cite{Co1}) and for arbitrary $S$ in the more general situation of Cartier dual $1$-motives by Bertapelle in \cite{Ber}. We will come back to the motivic nature of Deligne's pairing in 2.6.1.

\section{Algebraic connections, de Rham cohomology and $\mathcal D$-modules}
\markright{\uppercase{Preliminaries and notation}}
\subsection{Connections}
The algebraic theory of modules with connection is a basic tool used throughout this work, and we here recall various of its elementary notions and constructions which will be needed in the future.\\
\newline
Let $f:X \rightarrow S$ be a smooth morphism of schemes.\\
\newline
(i) If $\mathcal M$ is a $\mathcal O_X$-module a \underline{connection relative $S$ or $S$-connection} on $\mathcal M$ is a homomorphism of abelian sheaves
\[\nabla: \mathcal M \rightarrow \Omega^1_{X/S}\otimes_{\mathcal O_X} \mathcal M\]
satisfying the Leibniz rule. For each $i \geq 1$ one has a $f^{-1}\mathcal O_S$-linear map defined on local sections by
\[\nabla^i:\Omega^i_{X/S}\otimes_{\mathcal O_X} \mathcal M \rightarrow \Omega^{i+1}_{X/S}\otimes_{\mathcal O_X} \mathcal M, \quad \omega \otimes m \mapsto \mathrm{d}\omega \otimes m +(-1)^i\cdot \omega \wedge \nabla(m).\]
The connection $\nabla$ is called \underline{integrable} if its $\mathcal O_X$-linear \underline{curvature homomorphism}
\[K(\nabla):=\nabla^1\circ \nabla: \mathcal M \rightarrow \Omega^2_{X/S}\otimes_{\mathcal O_X} \mathcal M\]
vanishes, and in this case the maps $\nabla^i$ extend $\nabla$ to a $f^{-1}\mathcal O_S$-linear complex starting in degree zero
\[\Omega^{\bullet}_{X/S} (\mathcal M): \quad \lbrack \mathcal M \rightarrow \Omega^1_{X/S} \otimes_{\mathcal O_X} \mathcal M \rightarrow \Omega^2_{X/S} \otimes_{\mathcal O_X} \mathcal M \rightarrow...\rbrack,\]
called the \underline{de Rham complex} of $(\mathcal M,\nabla)$.\\
\newline
(ii) The datum of a $S$-connection on $\mathcal M$ is equivalent to a $\mathcal O_X$-linear map
\[\Delta: \underline{\mathrm {Der}}_{f^{-1}\mathcal O_S}(\mathcal O_X) \rightarrow \underline{\mathrm {End}}_{f^{-1}\mathcal O_S}(\mathcal M)\]
from the sheaf of $f^{-1}\mathcal O_S$-derivations on $\mathcal O_X$ into the sheaf of $f^{-1}\mathcal O_S$-linear endomorphisms on $\mathcal M$, satisfying the Leibniz rule; note that the $\mathcal O_X$-module structures of $\underline{\mathrm {Der}}_{f^{-1}\mathcal O_S}(\mathcal O_X)$ and $\underline{\mathrm {End}}_{f^{-1}\mathcal O_S}(\mathcal M)$ come from (outer) multiplication by sections of $\mathcal O_X$.\\
Observing that as $\mathcal O_X$-module $\underline{\mathrm {Der}}_{f^{-1}\mathcal O_S}(\mathcal O_X)$ identifies with the dual of $\Omega^1_{X/S}$ the mentioned equivalence is induced by composing a given $S$-connection $\nabla$ on $\mathcal M$ with sections $\theta \in \Gamma(U,\underline{\mathrm {Der}}_{f^{-1}\mathcal O_S}(\mathcal O_X))$:
\[\Delta(\theta): \mathcal M_{|U} \xrightarrow{\nabla} \Omega^1_{U/S}\otimes_{\mathcal O_U} \mathcal M_{|U} \xrightarrow{\theta \otimes \id} \mathcal M_{|U}.\]
The inverse assignment comes from choosing open subsets $U \subseteq X$ where a relative local coordinate system $\{ x_i,\partial_{x_i} \}_{1\leq i \leq n}$ is available (cf. 0.2.3 for this notion) and checking that the local formulas
\[\nabla(m)=\sum_{i=1}^ndx_i \otimes \Delta(\partial_{x_i})(m)\]
induce a well-defined $S$-connection globally on $\mathcal M$.\\
Under this correspondence the integrability condition of $\nabla$ translates into compatibility of $\Delta$ with the natural $f^{-1}\mathcal O_S$-Lie algebra structures of $\underline{\mathrm {Der}}_{f^{-1}\mathcal O_S}(\mathcal O_X)$ and $\underline{\mathrm {End}}_{f^{-1}\mathcal O_S}(\mathcal M)$.\\
\newline
(iii) A \underline{morphism} between two $\mathcal O_X$-modules with (integrable) $S$-connection is a $\mathcal O_X$-linear homomorphism that is compatible with the connections; such a map is also called \underline{horizontal}.\\
\newline
(iv) If $\mathcal M$ and $\mathcal N$ are $\mathcal O_X$-modules with (integrable) $S$-connections $\nabla_{\mathcal M}$ and $\nabla_{\mathcal N}$, then the tensor product $\mathcal M \otimes_{\mathcal O_X}\mathcal N$ can be equipped with a (integrable) $S$-connection
\[\nabla_{\mathcal M}\otimes \nabla_{\mathcal N}: \mathcal M \otimes_{\mathcal O_X}\mathcal N \rightarrow \Omega^1_{X/S}\otimes_{\mathcal O_X} (\mathcal M \otimes_{\mathcal O_X}\mathcal N),\]
the \underline{tensor product connection} of $\nabla_{\mathcal M}$ and $\nabla_{\mathcal N}$. Using associativity and commutativity of the tensor product it is induced by the following rule on local sections:
\[(\nabla_{\mathcal M}\otimes \nabla_{\mathcal N})(m \otimes n) \mapsto \nabla_{\mathcal M}(m)\otimes n + m \otimes \nabla_{\mathcal N}(n).\]
We can further endow the internal $\Hom$-sheaf $\underline{\Hom}_{\mathcal O_X}(\mathcal M,\mathcal N)$ with a (integrable) $S$-connection
\[\nabla_{\mathcal M,\mathcal N}: \underline{\Hom}_{\mathcal O_X}(\mathcal M,\mathcal N) \rightarrow \Omega^1_{X/S}\otimes_{\mathcal O_X} \underline{\Hom}_{\mathcal O_X}(\mathcal M,\mathcal N),\]
the \underline{internal $\Hom$-connection}. Using the canonical identification
\[\Omega^1_{X/S}\otimes_{\mathcal O_X} \underline{\Hom}_{\mathcal O_X}(\mathcal M,\mathcal N) \simeq \underline{\Hom}_{\mathcal O_X}(\mathcal M,\Omega^1_{X/S}\otimes_{\mathcal O_X}\mathcal N)\]
it is induced by the following rule on local sections:
\[\nabla_{\mathcal M,\mathcal N}(\varphi)(m)=\nabla_{\mathcal N}(\varphi(m))-(\id\otimes \varphi)(\nabla_{\mathcal M}(m)).\]
In the case that $\mathcal M$ is a vector bundle and $\mathcal N=\mathcal O_X$, equipped with its canonical integrable $S$-connection (exterior derivation), the resulting connection on the dual bundle $\mathcal M^\vee=\underline{\Hom}_{\mathcal O_X}(\mathcal M, \mathcal O_X)$ is called the \underline{dual connection} of $\nabla_{\mathcal M}$. One checks that its tensor product with $\nabla_{\mathcal N}$ becomes the connection $\nabla_{\mathcal M,\mathcal N}$ under the natural identification $\underline{\Hom}_{\mathcal O_X}(\mathcal M, \mathcal O_X) \otimes_{\mathcal O_X}\mathcal N \simeq \underline{\Hom}_{\mathcal O_X}(\mathcal M, \mathcal N)$.\\
If $\mathcal M$ is a line bundle the tensor product of $\nabla_{\mathcal M}$ with its dual connection becomes exterior derivation under the canonical isomorphism $\mathcal M \otimes_{\mathcal O_X}\mathcal M^\vee \simeq \mathcal O_X$.\\
\newline
Next, assume that we are given a commutative (not necessarily cartesian) diagram
\begin{equation*} \tag{\textbf{0.2.1}} \begin{split}
\begin{xy}
\xymatrix{
X' \ar[d]^{g} \ar[r] & S' \ar[d]^{h}\\
X \ar[r]^{f} & S}
\end{xy}
\end{split}
\end{equation*}
with $X/S$ and $X'/S'$ smooth.\\
\newline
(v) For a $\mathcal O_X$-module $\mathcal M$ with (integrable) $S$-connection we define its \underline{pullback along $(0.2.1)$}.\\
Namely, we may canonically equip the (module) pullback $g^*\mathcal M$ with a (integrable) $S'$-connection by the following formula on local sections:
\[g^*\mathcal M=\mathcal O_{X'}\otimes_{g^{-1}\mathcal O_X}g^{-1}\mathcal M \rightarrow \Omega^1_{X'/S'}\otimes_{\mathcal O_{X'}}g^*\mathcal M,\]
\[\qquad \qquad \qquad \qquad \quad x' \otimes g^{-1}(m) \mapsto \mathrm{d}x'\otimes g^*(m)+x'\cdot (\mathrm{can} \otimes \id) (g^*(\nabla(m))),\]
where $\mathrm{can} \otimes \id: g^*\Omega^1_{X/S}\otimes_{\mathcal O_{X'}}g^*\mathcal M \rightarrow \Omega^1_{X'/S'} \otimes_{\mathcal O_{X'}} g^*\mathcal M$ is the obvious canonical map.\\
\newline
By contrast to conventions sometimes made in the literature (cf. e.g. \cite{Kat2}, $(1.1.4)$) we won't work with an additional notation like $(g,h)^*\mathcal M$ when referring to this construction but will simply write $g^*\mathcal M$ as it will be transparent or mentioned explicitly along which diagram we form the pullback.\\
\newline
We record the special case $S'=S$ and $h=\id$, which says that we have a notion of pullback for modules with (integrable) $S$-connection along a morphism of smooth $S$-schemes.\\
As further special case one sees that for a smooth map $f:X \rightarrow S$ as above and any $\mathcal O_S$-module $\mathcal M$ the pullback $f^*\mathcal M$ carries a \underline{canonical integrable $S$-connection}, given by the rule
\[\mathrm{d}\otimes \id:\mathcal O_X \otimes_{f^{-1}\mathcal O_S} f^{-1}\mathcal M\rightarrow \Omega^1_{X/S} \otimes_{f^{-1}\mathcal O_S}f^{-1}\mathcal M.\]
(vi) Assume in addition that $(0.2.1)$ is cartesian. For a $\mathcal O_{X'}$-module $\mathcal M'$ with (integrable) $S'$-connection
\[\nabla: \mathcal M' \rightarrow \Omega^1_{X'/S'} \otimes_{\mathcal O_{X'}}\mathcal M'\]
we then have the notion of its \underline{higher direct images along $(0.2.1)$}.\\
Namely, applying the $i$-th higher direct image functor for abelian sheaves to $\nabla$ and identifying
\[R^ig_*(\Omega^1_{X'/S'} \otimes_{\mathcal O_{X'}}\mathcal M') \simeq R^ig_*(g^*\Omega^1_{X/S} \otimes_{\mathcal O_{X'}}\mathcal M') \simeq \Omega^1_{X/S} \otimes_{\mathcal O_X} R^ig_*\mathcal M'\]
yields a homomorphism of abelian sheaves $R^ig_*\mathcal M' \rightarrow \Omega^1_{X/S} \otimes_{\mathcal O_X} R^ig_* \mathcal M'$ which is a (integrable) $S$-connection on the $\mathcal O_X$-module $R^ig_*\mathcal M'$.\\
If we consider $g_*$ as a left exact functor from the abelian category of $\mathcal O_{X'}$-modules with integrable $S'$-connection to the same category with $X',S'$ replaced by $X,S$, then the value of its $i$-th right derived functor at $\mathcal M'$ coincides with our previous construction, as one can easily show.\\
In terms of functors on derived categories it equals the $i$-th cohomology sheaf of $(g,h)_+(\mathcal M')$, where $(g,h)_+=Rg_*$ is the triangulated functor of \cite{Lau}, $(3.3)$ - to be precise, we here assume all schemes locally noetherian, the map $g$ quasi-compact and $\mathcal M'$ quasi-coherent as $\mathcal O_{X'}$-module.\\
\newline
(vii) At some places of the work it will be convenient to adapt \underline{Grothendieck's viewpoint on connections}.\\
For this let $\Delta^1_{X/S}$ be the first infinitesimal neighborhood of the diagonal immersion $X \rightarrow X\times_S X$ and $p_1,p_2: \Delta^1_{X/S} \rightarrow X$ the maps induced by the projections of $X\times_S X$. Hence, the two compositions
\[X \xrightarrow{i} \Delta^1_{X/S} \xymatrix{
  \ar@<1ex>[r]^{p_1}\ar[r]_{p_2} & X}\]
both are the identity, where $i$ is the natural nilpotent closed immersion of square zero. To give a $S$-connection on the $\mathcal O_X$-module $\mathcal M$ then is equivalent to give an isomorphism 
\[p_1^*\mathcal M \simeq p_2^*\mathcal M\]
of $\mathcal O_{\Delta^1_{X/S}}$-modules which becomes the identity when restricted to $X$ (cf. \cite{Bert-Og}, §2, Prop. 2.9).\\
As we won't use it explicitly we only remark that the translation of integrability into this setting expresses as a cocycle condition for certain isomorphisms, which is understood best by using the concept of "stratifications" (cf. ibid., §4, Def. 4.3 and Thm. 4.8, resp. §2, Thm. 2.15, for $\Q$-schemes).

\subsection{De Rham cohomology}
Algebraic de Rham cohomology is naturally associated to modules with integrable connection and provides the cohomological framework in which most of this work will take place.\\
We briefly review its basic definitions and fix some standard notation. Some more details, especially about the Gauß-Manin connection whose construction we won't repeat in detail here, can be found in \cite{Kat2}, $(2.0)$-$(3.3)$\footnote{In this reference all modules with integrable connections are assumed as quasi-coherent, an assumption which is superfluous for the definitions and facts we provide here (cf. also the second mentioned source).} and \cite{Har}, Ch. III, 4.\\
\newline
Let $f:X \rightarrow S$ be a smooth morphism of schemes.\\
\newline
For a $\mathcal O_X$-module $\mathcal M$ with integrable $S$-connection and $i \geq 0$ we define a $\mathcal O_S$-module
\begin{equation*}
H^i_{\mathrm{dR}}(X/S,\mathcal M):=\mathbb R^i f_*(\Omega^{\bullet}_{X/S}(\mathcal M)),
\end{equation*}
where as in 0.2.1 the notation $\Omega^{\bullet}_{X/S} (\mathcal M)$ means the $f^{-1}\mathcal O_S$-linear de Rham complex of $\mathcal M$.\\
We call $H^i_{\mathrm{dR}}(X/S,\mathcal M)$ the \underline{$i$-th de Rham cohomology sheaf of $X/S$ with coefficients in $\mathcal M$}.\\
As already done here we usually suppress in our notation the integrable connection underlying $\mathcal M$.\\
For $\mathcal M:=\mathcal O_X$ with its canonical integrable $S$-connection (exterior derivation) we abbreviate
\[H^i_{\mathrm{dR}}(X/S):= H^i_{\mathrm{dR}}(X/S,\mathcal O_X)\]
and call $H^i_{\mathrm{dR}}(X/S)$ the \underline{$i$-th de Rham cohomology sheaf of $X/S$}.\\
\newline
Assuming additionally that $S$ is locally noetherian and $f$ quasi-compact the spectral sequence of hyperderived functors (use \cite{Huy}, Rem. 2.67)
\[\tag{\textbf{0.2.2}} E_1^{p,q}=R^q f_*(\Omega^p_{X/S} \otimes_{\mathcal O_X} \mathcal M) \Rightarrow E^{p+q}=H^{p+q}_{\mathrm{dR}}(X/S,\mathcal M)\]
and standard cohomological results imply that all $H^i_{\mathrm{dR}}(X/S,\mathcal M)$ are quasi-coherent $\mathcal O_S$-modules if $\mathcal M$ is quasi-coherent, and that they are coherent $\mathcal O_S$-modules if $f$ is proper and $\mathcal M$ coherent.\\
\newline
Now assume that $S$ is smooth over another scheme $T$.
\begin{equation*}
\begin{xy}
\xymatrix@C-0.3cm{
X \ar[rr]^{f} \ar[dr]& & S \ar[dl]\\
& T &}
\end{xy}
\end{equation*}
If $\mathcal M$ is a $\mathcal O_X$-module with integrable $T$-connection, then via the natural map
\[\Omega^1_{X/T}\otimes_{\mathcal O_X}\mathcal M \xrightarrow{\mathrm{can}\otimes \id} \Omega^1_{X/S}\otimes_{\mathcal O_X}\mathcal M\]
we obtain an induced integrable $S$-connection on $\mathcal M$ and associated de Rham cohomology sheaves $H^i_{\mathrm{dR}}(X/S,\mathcal M)$ as above. Each of them can then be equipped with a canonical integrable $T$-connection, its \underline{Gauß-Manin connection relative $T$}. The idea is to construct the whole de Rham complex of this connection as the row $E_1^{\bullet, i}$ in the spectral sequence associated to the naturally finitely filtered complex $\Omega^{\bullet}_{X/T}(\mathcal M)$ and the $0$-th hyperderived functor of $f_*$.

\subsection{$\mathcal D$-modules}
Giving a module with integrable connection is equivalent to giving a left module over the sheaf $\mathcal D$ of differential operators for the considered geometric situation. Therefore, when working with integrable connections we will often freely switch into the language of $\mathcal D$-modules, where a convenient functorial formalism on the level of derived categories is available.\\
\newline
As it was the case for modules with integrable connection we will consider $\mathcal D$-modules for a relative situation. In the present subsection it thus becomes at first necessary to briefly outline how we will deal with the lack of references in the literature for a theory of relative $\mathcal D$-modules. We then introduce the sheaf of relative differential operators for a smooth map of schemes and explain that the left modules over it correspond to modules with integrable connection. The definition of the functorial machinery in a relative situation, for which we refer to \cite{Lau}, Rappels $(3.3.1)$, happens completely analogously as in the absolute case (cf. e.g. \cite{Ho-Ta-Tan}, Ch. 1); we therefore content ourselves with introducing explicitly only the inverse and direct image functor for relative $\mathcal D$-modules, which are in fact the two derived functors mainly used in the course of this work. Finally, as these results will be frequently needed, we record the statement of Kashiwara's equivalence and the availability of the localization sequence in the relative case.\\
\newline
\textit{\textbf{In this subsection we assume that all occurring schemes are separated $\Q$-schemes.}}

\subsubsection{General remarks}

The literature known to us is almost exclusively restricted to the study of $\mathcal D$-modules on smooth varieties over the complex numbers. But, as was already said, in the course of this work we will want to use $\mathcal D$-module formalism for the situation of a scheme lying smooth over an arbitrary field of characteristic zero or, more generally, over a base scheme defined over $\Q$.\\
The fundamental functorial setup for such a general case is (very briefly) outlined in \cite{Lau}, Rappels $(3.3.1)$, and this is in fact the only source we know concerning relative $\mathcal D$-modules.\\
Nevertheless, the basic formal results about $\mathcal D$-modules that we will need entirely carry over from the absolute to the relative situation. When making use of such a result in the progress of the work our policy will be to cite the absolute statement, referring to the source \cite{Ho-Ta-Tan}, and briefly say that resp. (if it is not immediately clear) why no problems arise in the transition to the relative situation.

\subsubsection{Basic concepts}
Let $S$ be a noetherian scheme and $f:X \rightarrow S$ a morphism which is of finite type and smooth of relative dimension $n$.\\
\newline
The \underline{sheaf of differential operators of $X/S$}, denoted $\mathcal D_{X/S}$, is the $f^{-1}\mathcal O_S$-subalgebra of $\underline{\mathrm {End}}_{f^{-1}\mathcal O_S}(\mathcal O_X)$ generated by $\underline{\mathrm {Der}}_{f^{-1}\mathcal O_S}(\mathcal O_X)$ and $\mathcal O_X$, where $\mathcal O_X$ is included in $\underline{\mathrm {End}}_{f^{-1}\mathcal O_S}(\mathcal O_X)$ by multiplication.\\
\newline
$\mathcal D_{X/S}$ is a sheaf of $\mathcal O_X$-algebras which is non-commutative as soon as $n \geq 1$. In the same way as in the absolute case (cf. \cite{Ho-Ta-Tan}, 1.1) one has the following local description of $\mathcal D_{X/S}$:\\
As $X/S$ is smooth of relative dimension $n$ we can choose $X$-locally, say on $U \subseteq X$ open affine, sections $\{x_1,...,x_n \}$ of $\mathcal O_U$ such that $\{\mathrm{d}x_1,...,\mathrm{d}x_n \}$ forms a $\mathcal O_U$-basis for $\Omega^1_{U/S}$. We then have
\[\mathcal D_{U/S}=\bigoplus_{\alpha \in \N_0^n}\mathcal O_U \partial^{\alpha}_x \quad \textit{with} \ \ \partial^{\alpha}_x:=\partial^{\alpha_1}_{x_1}...\partial^{\alpha_n}_{x_n},\]
where $\{\partial_{x_1}, ..., \partial_{x_n} \}$ is the dual basis of $\{\mathrm{d}x_1,...,\mathrm{d}x_n \}$: recall that $\underline{\mathrm {Der}}_{f^{-1}\mathcal O_S}(\mathcal O_X)$ as $\mathcal O_X$-module identifies with the dual of $\Omega^1_{X/S}$. In particular, $\mathcal D_{X/S}$ is a locally free $\mathcal O_X$-module and thus quasi-coherent.\\
\newline
Our above choice of sections $\{x_i,\partial_{x_i} \}_{1\leq i \leq n}$ is the immediate generalization of what in the absolute case is called "local coordinate system" (cf. \cite{Ho-Ta-Tan}, 1.1 resp. A.5). The existence of such a \underline{relative local coordinate system for $X/S$} together with the characteristic zero assumption is the main reason why many of the basic formal results of the absolute situation carry over to relative $\mathcal D$-modules.\\
\newline
By a \underline{(left) $\mathcal D$-module for $X/S$} we understand a sheaf of left modules over the ring sheaf $\mathcal D_{X/S}$.\\
The relation to the earlier defined modules with integrable connection (cf. 0.2.1) is as follows:\\
For a $\mathcal O_X$-module $\mathcal M$ the datum of an integrable $S$-connection on $\mathcal M$ is equivalent to the datum of a left $\mathcal D_{X/S}$-module structure on $\mathcal M$ which is compatible with the $\mathcal O_X$-module structure.\\
To see this recall from 0.2.1 (ii) that an integrable $S$-connection on $\mathcal M$ tantamounts to a map
\[\Delta: \underline{\mathrm {Der}}_{f^{-1}\mathcal O_S}(\mathcal O_X) \rightarrow \underline{\mathrm {End}}_{f^{-1}\mathcal O_S}(\mathcal M)\]
which is $\mathcal O_X$-linear, satisfies the Leibniz rule and is compatible with the natural $f^{-1}\mathcal O_S$-Lie algebra structures of $\underline{\mathrm {Der}}_{f^{-1}\mathcal O_S}(\mathcal O_X)$ and $\underline{\mathrm {End}}_{f^{-1}\mathcal O_S}(\mathcal M)$.\\
The equivalence with a $\mathcal D_{X/S}$-left module structure on $\mathcal M$ compatible with the given $\mathcal O_X$-module structure then is established by the formula
\[\theta \cdot m = \nabla(\theta)(m), \quad \textit{where} \ \ \theta \in \underline{\mathrm {Der}}_{f^{-1}\mathcal O_S}(\mathcal O_X), m \in \mathcal M.\]
The verification is trivial (cf. also \cite{Ho-Ta-Tan}, Lemma 1.2.1).\\
\newline
Under this correspondence the $\mathcal D_{X/S}$-linear homomorphisms between two (left) $\mathcal D_{X/S}$-modules identify with the horizontal $\mathcal O_X$-linear homomorphisms between them. Note furthermore that $\mathcal O_X$ with its canonical $\mathcal D_{X/S}$-module structure becomes equipped with the integrable $S$-connection given by exterior derivation.
Finally, observe that in the special case $X=S$ we simply have $\mathcal D_{X/S}=\mathcal O_X$.

\subsubsection{Derived Categories}
We write $\Mod(\mathcal D_{X/S})$ for the abelian category of (left) $\mathcal D$-modules for $X/S$ and $\Mod_{\mathrm{qc}}(\mathcal D_{X/S})$ for the full subcategory of $\Mod(\mathcal D_{X/S})$ consisting of those $\mathcal D_{X/S}$-modules which are quasi-coherent over $\mathcal O_X$.\footnote{For the notion of quasi-coherence for sheaves of modules on a ringed space cf. \cite{EGAI}, Ch. 0, $(5.1.3)$. Furthermore, we remark that for a $\mathcal D_{X/S}$-module the conditions to be quasi-coherent over $\mathcal D_{X/S}$ and to be quasi-coherent over $\mathcal O_X$ coincide: use \cite{EGAI}, Ch. I, Prop. $(2.2.4)$, and the already seen quasi-coherence of $\mathcal D_{X/S}$ over $\mathcal O_X$.} It is a thick abelian subcategory of $\Mod(\mathcal D_{X/S})$. \footnote{This follows from \cite{EGAI}, Ch. I, Cor. $(2.2.2)$.}\\
When considering the derived category of $\Mod(\mathcal D_{X/S})$ we use the notation
\[D^{\#}(\mathcal D_{X/S}):=D^{\#}(\Mod(\mathcal D_{X/S})),\]
where $\#$ is one of the boundedness conditions $\emptyset, +,-,b$.\\
We write $D^{\#}_{\mathrm{qc}}(\mathcal D_{X/S})$ for the full triangulated subcategory of $D^{\#}(\mathcal D_{X/S})$ consisting of those complexes with cohomology sheaves in $\Mod_{\mathrm{qc}}(\mathcal D_{X/S})$.\\
\newline
Finally, we use the occasion to insert analogous notational conventions for modules over the structure sheaf of a scheme\footnote{Of course, for this short general insertion we deactivate the assumption that we are working with separated $\Q$-schemes.}:\\
If $(B,\mathcal O_B)$ is a scheme we write $\Mod(\mathcal O_B)$ for the abelian category of sheaves of modules over $\mathcal O_B$ and $\Mod_{\mathrm{qc}}(\mathcal O_B)$ for the full subcategory consisting of those $\mathcal O_B$-modules which are quasi-coherent over $\mathcal O_B$. It is a thick abelian subcategory of $\Mod(\mathcal O_B)$.\footnote{As before, this follows from \cite{EGAI}, Ch. I, Cor. $(2.2.2)$.}\\
Concerning the derived category of $\Mod(\mathcal O_B)$ we write
\[D^{\#}(\mathcal O_B):=D^{\#}(\Mod(\mathcal O_B)),\]
where again $\#$ is one of the boundedness conditions $\emptyset, +,-,b$.\\
We denote with $D^{\#}_{\mathrm{qc}}(\mathcal O_B)$ the full triangulated subcategory of $D^{\#}(\mathcal O_B)$ consisting of those complexes with cohomology sheaves in $\Mod_{\mathrm{qc}}(\mathcal O_B)$.

\subsubsection{Inverse and direct image for $\mathcal D$-modules}
We briefly record the two most important triangulated functors on derived categories of $\mathcal D$-modules, namely the inverse and direct image functor. We will adopt the notations of \cite{Lau}, Rappels $(3.3.1)$, where one can also find the definitions of the other basic functors for relative $\mathcal D$-modules. But observe that compared with \cite{Lau} we impose more geometric conditions on the involved schemes, essentially regularity and finite-dimensionality: in our view, these assumptions seem necessary to ensure that the various triangulated functors are well-defined between bounded derived categories with quasi-coherent cohomology; we don't further elaborate on this technical issue here.\footnote{Let us only mention as a main reason that one needs the existence of bounded flat resolutions for $\mathcal D$-modules, which is guaranteed if the weak global dimensions of the stalks of the considered sheaf of relative differential operators are bounded; under our additional assumptions one can indeed bound them, namely by the sum of the dimension of the scheme with the relative dimension of the considered smooth morphism - this works similarly as in\cite{Ho-Ta-Tan}, p. 26 and Prop. 1.4.6.}\\
\newline
Let $S$ be a noetherian, regular and finite-dimensional scheme and let $X,Y$ be schemes which are of finite type and smooth of relative dimensions $d_{X/S}, d_{Y/S}$ over $S$; note that $X,Y$ then have the same geometric properties as $S$.\\
Assume that we are given an $S$-morphism $g: X \rightarrow Y$.
\begin{equation*}
\begin{xy}
\xymatrix@C-0.3cm{
X \ar[rr]^{g} \ar[dr]& & Y \ar[dl]\\
& S &}
\end{xy}
\end{equation*}
The \underline{inverse image functor} associated with $g$ is the triangulated functor
\[g^!: D^b_{\mathrm{qc}}(\mathcal D_{Y/S}) \rightarrow D^b_{\mathrm{qc}}(\mathcal D_{X/S})\]
given by
\[g^!\mathcal N^{\bullet}=\mathcal D_{(X\rightarrow Y)/S} \otimes^L_{g^{-1}\mathcal D_{Y/S}}g^{-1}\mathcal N^{\bullet}[d_{X/S}-d_{Y/S}]\]
for $\mathcal N^{\bullet}$ an object of $D^b_{\mathrm{qc}}(\mathcal D_{Y/S})$.\\
\newline
The \underline{direct image functor} associated with $g$ is the triangulated functor
\[g_+: D^b_{\mathrm{qc}}(\mathcal D_{X/S}) \rightarrow D^b_{\mathrm{qc}}(\mathcal D_{Y/S})\]
given by
\[g_+\mathcal M^{\bullet}=Rg_*(\mathcal D_{(Y \leftarrow X)/S}\otimes^L_{\mathcal D_{X/S}} \mathcal M^{\bullet})\]
for $\mathcal M^{\bullet}$ an object of $D^b_{\mathrm{qc}}(\mathcal D_{X/S})$.\\
\newline
Here, $\mathcal D_{(X\rightarrow Y)/S}$ and $\mathcal D_{(Y \leftarrow X)/S}$ are the two \underline{transfer bimodules} associated with the morphism $g$. For their definition cf. \cite{Lau}, $(3.3.1)$ or \cite{Ho-Ta-Tan}, Def. 1.3.1, Def. 1.3.3 and Lemma 1.3.4.\\
\newline
We will frequently and tacitly use the following two easy facts about the inverse image functor:\\
First, if $\mathcal N$ is a locally free $\mathcal O_Y$-module with integrable $S$-connection, considered as object of $D^b_{\mathrm{qc}}(\mathcal D_{Y/S})$ in the natural way, then we have canonically
\[g^!\mathcal N \simeq g^*\mathcal N[d_{X/S}-d_{Y/S}] \quad \textrm{in} \ D^b_{\mathrm{qc}}(\mathcal D_{X/S}),\]
where $g^*\mathcal N$ is equipped with the integrable $S$-connection given by pullback (cf. 0.2.1 (v)).\footnote{To verify this claim one first shows as in \cite{Ho-Ta-Tan}, proof of Prop. 1.5.8, that on the level of $\mathcal O$-modules the functor $g^!$ is derived pullback shifted by $[d_{X/S}-d_{Y/S}]$. One then uses the local freeness of $\mathcal N$ to see that $g^!\mathcal N[d_{Y/S}-d_{X/S}]$ is concentrated in degree zero. Finally, as in \cite{Ho-Ta-Tan}, p. 22, one checks that the zeroth cohomology of $g^!\mathcal N[d_{Y/S}-d_{X/S}]$ canonically identifies with $g^*\mathcal N$ as $\mathcal D_{X/S}$-module.}\\
Second, if $g$ is an open immersion, then the functor $g^!$ is given by the natural restriction $g^{-1}$ to the open subset $X$ of $Y$, and we simply write $\mathcal N_{|X}$ for this restriction.\\
\newline
As we will cite results from the source \cite{Ho-Ta-Tan} it should be remarked that in their notation the functor $g^!$ is written as $g^\dagger$ and the functor $g_+$ as $\int_g$ (cf. ibid., p. 33 and p. 40).

\subsubsection{Kashiwara's equivalence and localization sequence}
Going through the proofs of \cite{Ho-Ta-Tan}, Thm. 1.6.1 and Cor. 1.6.2, one checks that they can be modified to yield the following relative version of Kashiwara's equivalence:\\
\newline
Let $S$ be a noetherian, regular and finite-dimensional scheme and let $X,Z$ be schemes which are of finite type and smooth of some fixed relative dimensions over $S$.\\
Assume that we have a closed immersion $i: Z \rightarrow X$ of a fixed codimension $c$ and write $j$ for the associated open immersion $j: U:=X \backslash Z \rightarrow X$; here, as in \cite{Mi}, Ch. VI, §5, we say that $Z$ has codimension $c$ in $X$ if for any $s \in S$ the fiber $Z_s$ has pure codimension $c$ in $X_s$.
\begin{equation*}
\begin{xy}
\xymatrix{
\qquad \qquad & Z \ar[dr] \ar[r]^{i} & X  \ar[d]^{f} & U \ar[dl] \ar[l]_{j}\\
\qquad \qquad &  & S & &}
\end{xy}
\end{equation*}
\underline{Kashiwara's equivalence} in its relative form says that then the functor $H^0 i_+$ induces an equivalence
\[\tag{\textbf{0.2.3}} H^0 i_+: \Mod_{\mathrm{qc}}(\mathcal D_{Z/S}) \xrightarrow{\sim} \Mod_{\mathrm{qc}}^Z(\mathcal D_{X/S}),\]
and that the functor $i_+$ induces an equivalence
\[\tag{\textbf{0.2.4}} i_+:D^b_{\mathrm{qc}}(\mathcal D_{Z/S}) \xrightarrow{\sim} D^{b,Z}_{\mathrm{qc}}(\mathcal D_{X/S}).\]
Here, the superscripts mean the full subcategory of $\Mod_{\mathrm{qc}}(\mathcal D_{X/S})$ consisting of those modules whose support is contained in $Z$ resp. the full triangulated subcategory of $D^b_{\mathrm{qc}}(\mathcal D_{X/S})$ consisting of those complexes whose cohomology sheaves have support contained in $Z$.\\
The quasi-inverse of $(0.2.3)$ resp. $(0.2.4)$ is induced by the functor $H^0i^!$ resp. $i^!$.\\
\newline
In particular, from $(0.2.4)$ one deduces (along the same lines as in \cite{Ho-Ta-Tan}, Prop. 1.7.1) that for every $\mathcal M^{\bullet} \in D^b_{\mathrm{qc}}(\mathcal D_{X/S})$ there is a canonical distinguished triangle in $D^b_{\mathrm{qc}}(\mathcal D_{X/S})$
\[\tag{\textbf{0.2.5}}i_+ i^! \mathcal M^{\bullet} \rightarrow \mathcal M^{\bullet} \rightarrow j_+\mathcal M^{\bullet}_{|U}\]
which we will refer to as the \underline{localization sequence} or the \underline{canonical distinguished triangle} for $\mathcal M^{\bullet}$.
\begin{remark}
(i) The essential point to prove the equivalences $(0.2.3)$ and $(0.2.4)$ is the following:\\
Our geometric assumptions imply that around points of $X\cap Z$ the map $f$ factorizes locally on $X$ via an étale morphism into $\Spec (\mathcal O_S [T_1,...T_n])$ such that $Z$ is the inverse image of the closed subscheme defined by $T_{n-c+1}=...=T_n=0$. Here, $c$ is again the codimension of $Z$ in $X$, and $n$ is the relative dimension of $X/S$. The étale morphism comes from local sections $\{x_1,...,x_n \}$ of $\mathcal O_X$ such that $\{\mathrm{d}x_1,...\mathrm{d}x_n \}$ is locally a basis of $\Omega^1_{X/S}$. By dualizing we obtain a local basis $\{\partial_{x_1},...\partial_{x_n} \}$ of the tangent sheaf of $X/S$, and the set $\{x_i,\partial_{x_i} \}_{1\leq i \leq n}$ then is a relative local coordinate system for $X/S$, directly generalizing the (absolute) "local coordinate system" of \cite{Ho-Ta-Tan}, 1.1 resp. A.5. With this observation it is straightforward to generalize the proof of Kashiwara's equivalence (and of many other basic formal results for absolute $\mathcal D$-modules) to the relative situation.\vspace{1mm}\\
(ii) In the course of the work we will need the localization sequence $(0.2.5)$ almost exclusively in the case that $S$ is the spectrum of a field of characteristic zero.
\end{remark}

\newpage
\markright{\uppercase{Preliminaries and notation}}
\section{Summary of the basic notation}

It seems convenient to collect in a list the most important and steadily used pieces of the notation introduced so far.
All further symbols appearing in the progress of this work and not ranging among common standard notation will be introduced in situ.\\
\newline
$\bullet$ If $(B,\mathcal O_B)$ is a scheme:\\
\newline
$\Mod(\mathcal O_B)$: the category of sheaves of $\mathcal O_B$-modules\\
$\Mod_{\mathrm{qc}}(\mathcal O_B)$: the category of quasi-coherent sheaves of $\mathcal O_B$-modules\\
$D^{\#}(\mathcal O_B)$: the derived category of $\Mod(\mathcal O_B)$ with boundedness condition $\#=\emptyset,+,-,b$\\
$D^{\#}_{\mathrm{qc}}(\mathcal O_B)$: the full triangulated subcategory of $D^{\#}(\mathcal O_B)$ of complexes with cohomology in $\Mod_{\mathrm{qc}}(\mathcal O_B)$\\
$\mathcal F^\vee:=\underline{\Hom}_{\mathcal O_B}(\mathcal F, \mathcal O_B)$ for a $\mathcal O_B$-module $\mathcal F$\\
$\mathcal L^{-1}:=\mathcal L^\vee$ for a $\mathcal O_B$-line bundle $\mathcal L$ \\
$\Lie(G/B)$: the Lie algebra of $G/B$ for a commutative group scheme $G$ over $B$ \\
$\mathbb V(\mathcal E)$: the geometric vector bundle associated with $\mathcal E$ for a $\mathcal O_B$-vector bundle $\mathcal E$\\
\newline
$\bullet$ If $f:X \rightarrow S$ is a smooth map of schemes and $\mathcal M$ is a $\mathcal O_X$-module with integrable $S$-connection:\\
\newline
$\underline{\mathrm {Der}}_{f^{-1}\mathcal O_S}(\mathcal O_X)$: the sheaf of $f^{-1}\mathcal O_S$-derivations on $\mathcal O_X$\\
$\underline{\mathrm {End}}_{f^{-1}\mathcal O_S}(\mathcal M)$: the sheaf of $f^{-1}\mathcal O_S$-linear endomorphisms on $\mathcal M$\\
$\Omega^{\bullet}_{X/S}(\mathcal M)$: the de Rham complex of $\mathcal M$ (starting in degree zero)\\
$H^i_{\mathrm{dR}}(X/S,\mathcal M)$: the $i$-th de Rham cohomology sheaf of $X/S$ with coefficients in $\mathcal M$\\
$H^i_{\mathrm{dR}}(X/S)$: the $i$-th de Rham cohomology sheaf of $X/S$\\
\newline
For the next two points all schemes are assumed to be separated $\Q$-schemes.\\
\newline
$\bullet$ If $X$ is of finite type and smooth (of a fixed relative dimension) over a noetherian scheme $S$:\\
\newline
$\mathcal D_{X/S}$: the sheaf of differential operators of $X$ relative $S$\\
$\Mod(\mathcal D_{X/S})$: the category of sheaves of left $\mathcal D_{X/S}$-modules\\
$\Mod_{\mathrm{qc}}(\mathcal D_{X/S})$: the category of sheaves of left $\mathcal D_{X/S}$-modules which are quasi-coherent over $\mathcal O_X$\\
$D^{\#}(\mathcal D_{X/S})$: the derived category of $\Mod(\mathcal D_{X/S})$ with boundedness condition $\#=\emptyset,+,-,b$\\
$D^{\#}_{\mathrm{qc}}(\mathcal D_{X/S})$: the full triangulated subcategory of $D^{\#}(\mathcal D_{X/S})$ of complexes with cohomology in $\Mod_{\mathrm{qc}}(\mathcal D_{X/S})$\\
\newline
$\bullet$ If $X, Y$ are of finite type and smooth (of fixed relative dimensions) over a noetherian, regular and finite-dimensional scheme $S$, and if $g: X\rightarrow Y$ is a $S$-morphism:\\
\newline
$g^!:D^b_{\mathrm{qc}}(\mathcal D_{Y/S}) \rightarrow D^b_{\mathrm{qc}}(\mathcal D_{X/S})$: the inverse image functor associated with $g$\\
$g_+:D^b_{\mathrm{qc}}(\mathcal D_{X/S}) \rightarrow D^b_{\mathrm{qc}}(\mathcal D_{Y/S})$: the direct image functor associated with $g$\\
$\mathcal D_{(X\rightarrow Y)/S}$: the first transfer bimodule associated with $g$\\
$\mathcal D_{(Y \leftarrow X)/S}$: the second transfer bimodule associated with $g$\\
\newline
$\bullet$ If $S$ is a locally noetherian scheme and $X$ an abelian scheme of relative dimension $g$ over $S$:\\
\newline
$\pi: X \rightarrow S$: the structure map of $X/S$\\
$\mu: X\times_S X \rightarrow X$: the multiplication map of $X/S$\\
$\epsilon: S \rightarrow X$: the zero section of $X/S$\\
$(-1)_X: X\rightarrow X$: the inverse map of $X/S$\\
$\mathrm{pr}_1, \mathrm{pr}_2: X\times_S X \rightarrow X$: the two projection maps\\
$T\mapsto \mathrm{Pic}^0(X\times_S T/T)$: the dual functor of $X/S$ on the category of all $S$-schemes\\
$Y$: the dual abelian scheme of $X$\\
$(\mathcal P^0,r^0)$: the fixed representative for the universal isomorphism class in $\mathrm{Pic}^0(X\times_S Y/Y)$\\
$(\mathcal P^0,r^0,s^0)$: the birigidified Poincaré bundle on $X\times_S Y$\\
$T\mapsto \mathrm{Pic}^\natural(X\times_S T/T)$: the $\natural$-dual functor of $X/S$ on the category of all $S$-schemes\\
$\Yr$: the universal vectorial extension of $Y$\\
$\pi^\natural: \Yr \rightarrow S$: the structure map of $\Yr/S$\\
$\mu^\natural: \Yr \times_S \Yr \rightarrow \Yr$: the multiplication map of $\Yr/S$\\
$\epsilon^\natural: S \rightarrow \Yr$: the zero section of $\Yr/S$\\
$(-1)^\natural: \Yr \rightarrow \Yr$: the inverse map of $\Yr/S$\\
$(\mathcal P,r,s,\nabla_\mathcal P)$: the birigidified Poincaré bundle with universal integrable $\Yr$-connection on $X\times_S \Yr$
\markright{\uppercase{Preliminaries and notation}}

\newpage
\markright{\uppercase{The formalism of the logarithm sheaves and the elliptic...}}
\chapter{The formalism of the logarithm sheaves and the elliptic polylogarithm}
For the whole chapter we let $X$ be an abelian scheme of relative dimension $g\geq 1$ over a connected scheme $S$ which is smooth, separated and of finite type over $\Spec(\Q)$.\\
As fixed in 0.1 the structure morphism resp. zero section of $X/S$ will be denoted by $\pi$ resp. $\epsilon$.
\begin{equation*}
\begin{xy}
\xymatrix@C-0.3cm{
X \ar[rr]^{\pi} \ar[dr]& & S \ar@ /_ 0.6cm/[ll]_{\epsilon}\ar[dl]\\
& \Spec(\Q) &}
\end{xy}
\end{equation*}

\subsubsection{Further remarks about de Rham cohomology}
We briefly recall some additional information about the de Rham cohomology of $X/S$ which will be tacitly used in the future. For more details and for proofs we refer to \cite{Bert-Br-Mes}, 2.5.\\
Basic notations and facts concerning algebraic de Rham cohomology can be found in 0.2.2.\\
\newline
For the abelian scheme $X/S$ each cohomology sheaf $H^i_{\mathrm{dR}}(X/S)$ is a vector bundle on $S$ commuting with arbitrary base change. In the lowest and top degree we have the following description:\\
The degeneration of the Hodge-de Rham spectral sequence
\[E^{p,q}_1=R^q\pi_*\Omega^p_{X/S} \Rightarrow E^{p+q}=H^{p+q}_{\mathrm{dR}}(X/S)\]
at the first sheet and the canonical isomorphism $\mathcal O_S \xrightarrow{\sim}\pi_*\mathcal O_X$ provides a natural isomorphism
\[H^0_{\mathrm{dR}}(X/S) \simeq \mathcal O_S.\]
The same degeneration yields the identification $R^g\pi_*\Omega^g_{X/S}\xrightarrow{\sim} H^{2g}_{\mathrm{dR}}(X/S)$ whose inverse can be composed with the Grothendieck trace map (cf. \cite{Con1}, Ch. I, 1.1) to give the trace isomorphism
\[\mathrm{tr}:H^{2g}_{\mathrm{dR}}(X/S) \xrightarrow{\sim} \mathcal O_S.\]
Moreover, the cup product defines the structure of an alternating graded $\mathcal O_S$-algebra on $\bigoplus_{i=0}^{2g}H^i_{\mathrm{dR}}(X/S)$, and it induces for each $i=0,...,2g$ an isomorphism
\[\bigwedge^i H^1_{\mathrm{dR}}(X/S) \xrightarrow{\sim} H^i_{\mathrm{dR}}(X/S).\]
From this we clearly obtain a Poincaré duality identification
\[H^i_{\mathrm{dR}}(X/S) \simeq H^{2g-i}_{{\mathrm{dR}}}(X/S)^\vee, \quad x \mapsto \{\ y\mapsto \mathrm{tr}(x\cup y)\}.\]
Let us finally point out that if all sheaves $H^i_{\mathrm{dR}}(X/S)$ are equipped with their Gauß-Manin connection relative $\Spec(\Q)$, if duals resp. exterior powers are endowed with the naturally induced connection, and if $\mathcal O_S$ carries its canonical connection (exterior derivative), then all of the previous four isomorphisms become horizontal, as one would expect.\footnote{For the identification $H^0_{\mathrm{dR}}(X/S) \simeq \mathcal O_S$ this is obvious by \cite{Kat2}, $(3.1.0)$. The horizontality of the trace isomorphism is a well-known fact which will also come out as a side result later in 1.2.3. The claim for the fourth identification then becomes clear if one further observes that the cup product $H^i_{\mathrm{dR}}(X/S) \otimes_{\mathcal O_S} H^{2g-i}_{{\mathrm{dR}}}(X/S) \rightarrow H^{2g}_{\mathrm{dR}}(X/S)$ is horizontal. This in turn, as well as the horizontality of the third of the above isomorphisms, is easily deduced from \cite{Kat-Od}, $(11)$.}

\section{The definition of the logarithm sheaves}

We define the de Rham realization of the logarithm sheaves for our fixed geometric situation $X/S/\Q$ and introduce some elementary features like the transition maps, the natural unipotent filtration and a first remark about compatibility with base change. The general proceeding is analogous to \cite{Hu-Ki}, Def. A. 1.3, where the $\ell$-adic setting is considered, with the difference that we need to fix a splitting for the pullback of the first logarithm sheaf along the zero section as in our case it is not unique.\\
\newline
Let us write $\mathcal H$ for the dual of the $\mathcal O_S$-vector bundle $H^1_{\mathrm{dR}}(X/S)$ and endow it with the integrable $\Q$-connection given by the dual of the Gauß-Manin connection on $H^1_{\mathrm{dR}}(X/S)$.\\
The $\mathcal O_X$-vector bundle
\[\mathcal H_X:=\pi^*\mathcal H\]
then carries the pullback connection relative $\Spec(\Q)$ (cf. 0.2.1 (v)).
\subsubsection{Definition of the logarithm sheaves}
The Leray spectral sequence in de Rham cohomology (cf. \cite{Kat2}, $(3.3.0)$) for $\mathcal H_X$ reads as
\[E_2^{p,q}=H^p_{\mathrm{dR}}(S/\Q,H^q_{\mathrm{dR}}(X/S)\otimes_{\mathcal O_S}\mathcal H) \Rightarrow E^{p+q}=H_{\mathrm{dR}}^{p+q}(X/\Q,\mathcal H_X),\]
where $H^q_{\mathrm{dR}}(X/S)\otimes_{\mathcal O_S}\mathcal H$ carries the tensor product connection.\\
Under the standard canonical identifications
\[E_2^{p,q} \simeq \Ext^p_{\Mod_{\mathrm{qc}}(\mathcal D_{S/\Q})}(\mathcal O_S, H^q_{\mathrm{dR}}(X/S)\otimes_{\mathcal O_S}\mathcal H) \ \ , \ \ E^{p+q} \simeq \Ext^{p+q}_{\Mod_{\mathrm{qc}}(\mathcal D_{X/\Q})}(\mathcal O_X, \mathcal H_X)\]
it can be viewed as $\mathcal H_X$ plugged into the spectral sequence for the composition of functors
\[\Hom_{\mathcal D_{S/\Q}}(\mathcal O_S,-) \circ H^0_{\mathrm{dR}}(X/S,-): \Mod_{\mathrm{qc}}(\mathcal D_{X/\Q})\rightarrow \Mod_{\mathrm{qc}}(\mathcal D_{S/\Q})\rightarrow (\Q\textrm{-vector spaces}),\]
noting that this composition equals $\Hom_{\mathcal D_{X/\Q}}(\mathcal O_X,-)$.\\
\newline
The associated five term exact sequence yields the $\Q$-linear exact sequence
\[\tag{\textbf{1.1.1}} 0\rightarrow \Ext^1_{\mathcal D_{S/\Q}}(\mathcal O_S, \mathcal H) \xrightarrow{\pi^*} \Ext^1_{\mathcal D_{X/\Q}}(\mathcal O_X, \mathcal H_X) \rightarrow \Hom_{\mathcal D_{S/\Q}}(\mathcal O_S, \mathcal H^\vee \otimes_{\mathcal O_S} \mathcal H) \rightarrow 0,\]
where we have used the existence of the section $\epsilon$ to deduce injectivity of the map $E_2^{2,0} \rightarrow E^2$ and hence the surjectivity in $(1.1.1)$. That the injection in $(1.1.1)$ is given by pullback along $\pi$ is straightforwardly checked as well as the following description of the occurring projection:\\
For a class in $\Ext^1_{\mathcal D_{X/\Q}}(\mathcal O_X, \mathcal H_X)$, represented by a $\mathcal D_{X/\Q}$-linear extension
\[0 \rightarrow \mathcal H_X \rightarrow \mathcal M \rightarrow \mathcal O_X \rightarrow 0,\]
the first boundary map in the long exact sequence for the derived functors of $H^0_{\mathrm{dR}}(X/S,-)$ (cf. \cite{Kat2}, Rem. $(3.1)$) is a horizontal map
\[\mathcal O_S \rightarrow \mathcal H^\vee \otimes_{\mathcal O_S} \mathcal H,\]
and this is the image of our class in $\Hom_{\mathcal D_{S/\Q}}(\mathcal O_S, \mathcal H^\vee \otimes_{\mathcal O_S} \mathcal H)$.\\
We also remark that we will often tacitly identify $\Hom_{\mathcal D_{S/\Q}}(\mathcal O_S, \mathcal H^\vee \otimes_{\mathcal O_S} \mathcal H) \simeq \Hom_{\mathcal D_{S/\Q}}(\mathcal H, \mathcal H)$.\\
\newline
Now observe that the map $\pi^*$ in $(1.1.1)$ has a retraction defined by $\epsilon^*$, and that hence our exact sequence $(1.1.1)$ splits. With this we are already prepared to introduce the logarithm sheaves.
\markright{\uppercase{The formalism of the logarithm sheaves and the elliptic...}}

\begin{definition}
(i) The class in $\Ext^1_{\mathcal D_{X/\Q}}(\mathcal O_X, \mathcal H_X)$ mapping to the identity under the projection of $(1.1.1)$ and to zero under the retraction $\epsilon^*$ is called the \underline{first logarithm extension class of $X/S/\Q$} and written as $\mathcal Log^1$.\\
(ii) Suppose we are given a fixed pair consisting of an exact sequence of $\mathcal D_{X/\Q}$-modules
\[\tag{\textbf{1.1.2}} 0 \rightarrow \mathcal H_X \rightarrow \mathcal L_1 \rightarrow \mathcal O_X \rightarrow 0\]
representing the first logarithm extension class together with the choice of a $\mathcal D_{S/\Q}$-linear splitting
\[\varphi_1: \mathcal O_S \oplus \mathcal H \simeq \epsilon^* \mathcal L_1\]
for the pullback of $(1.1.2)$ along $\epsilon$.\\
Then the vector bundle $\mathcal L_1$ on $X$ together with its integrable $\Q$-connection, denoted by $\nabla_1$, and with the splitting $\varphi_1$ will be called the \underline{first logarithm sheaf of $X/S/\Q$} and denoted by $(\mathcal L_1, \nabla_1, \varphi_1)$.
\end{definition}

\begin{remark}
One can show that an extension representing $\mathcal Log^1$ has a nontrivial automorphism if and only if the sheaf $\mathcal H$ has a nonzero global horizontal section; if we however additionally require compatibility with a fixed splitting for the pullback of the extension along $\epsilon$, then any such automorphism is the identity. Such facts will be picked up and proven slightly more generally in Lemma 2.1.2 and Rem. 2.1.3.
\end{remark}

From now on let $(\mathcal L_1,\nabla_1,\varphi_1)$ be as in Def. 1.1.1 (ii).\\
\newline
For each $n\geq 1$ we define a vector bundle $\mathcal L_n$ on $X$ by $\mathcal L_n:=\mathrm{Sym}^n_{\mathcal O_X}\mathcal L_1$. It is equipped with the induced integrable $\Q$-connection, denoted by $\nabla_n$, and $\varphi_1$ induces the (horizontal) decomposition
\[\varphi_n: \prod_{k=0}^n \mathrm{Sym}^k_{\mathcal O_S}\mathcal H \simeq \epsilon^*\mathcal L_n.\]
For $n=0$ we let $\mathcal L_n:=\mathcal O_X$, equipped with its canonical integrable $\Q$-connection, and we let $\varphi_0$ be the natural isomorphism $\varphi_0: \mathcal O_S \simeq \epsilon^*\mathcal O_X$.
\begin{definition}
For each $n\geq 0$ the so defined triple $(\mathcal L_n,\nabla_n,\varphi_n)$ is called \underline{$n$-th logarithm sheaf of $X/S/\Q$}.
\end{definition}
\begin{remark}
We remark that although the connection $\nabla_n$ and the splitting $\varphi_n$ will frequently be suppressed in the notation these data always remain fixed; they belong to the definition of the logarithm sheaves.
\end{remark}
\subsubsection{The transition maps}
Let us denote for a moment with $\mathrm{pr}: \mathcal L_1 \rightarrow \mathcal O_X$ the projection in $(1.1.2)$.\\
For each $n\geq 0$ there is a natural exact sequence of $\mathcal D_{X/\Q}$-modules
\[\tag{\textbf{1.1.3}} 0 \rightarrow \mathrm{Sym}^{n+1}_{\mathcal O_X}\mathcal H_X \rightarrow \mathcal L_{n+1} \rightarrow \mathcal L_n \rightarrow 0.\]
Here, the first (nontrivial) arrow is induced by the map $\mathcal H_X \rightarrow \mathcal L_1$ of $(1.1.2)$, and the second is defined to be the composition
\[\mathrm{Sym}_{\mathcal O_X}^{n+1}\mathcal L_1 \rightarrow \mathrm{Sym}_{\mathcal O_X}^{n+1}(\mathcal L_1 \oplus \mathcal O_X) \rightarrow \mathrm{Sym}^n_{\mathcal O_X}\mathcal L_1\]
of the morphism on symmetric powers related to $\id \oplus \mathrm{pr}: \mathcal L_1 \rightarrow \mathcal L_1 \oplus \mathcal O_X$ with the map coming from the decomposition of the symmetric power of a direct sum. Exactness of $(1.1.3)$ is readily checked.\\
\newline
We briefly record how the previous transition maps of the logarithm sheaves express after pullback along the zero section:
\begin{Lemma}
Let $n\geq 0$. If we pull back the transition map $\mathcal L_{n+1} \rightarrow \mathcal L_n$ of $(1.1.3)$ via $\epsilon$ and use the splittings $\varphi_{n+1}$ and $\varphi_n$ then the induced $\mathcal D_{S/\Q}$-linear map
\[\qquad \prod_{k=0}^{n+1} \mathrm{Sym}^k_{\mathcal O_S}\mathcal H \rightarrow \prod_{k=0}^n \mathrm{Sym}^k_{\mathcal O_S}\mathcal H
\]
is given on sections explicitly by
\[\qquad (s_0,s_1,...,s_{n+1}) \mapsto ((n+1)\cdot s_0, n \cdot s_1,...,s_n), \quad \textrm{where}\ s_k \in \mathrm{Sym}^k_{\mathcal O_S}\mathcal H, \ k=0,...,n+1.
\]
\end{Lemma}
\begin{proof}
The proof is a straightforward calculation with the definitions.
\end{proof}
\subsubsection{The unipotent filtration}
For the computation of their de Rham cohomology and the formulation of their universal property it is crucial to consider the logarithm sheaves as unipotent objects. We will expose this in more detail later (cf. 1.3) and for now content ourselves with the following basic observation:\\
\newline
Namely, the exact sequence
\[0 \rightarrow \mathcal H_X \rightarrow \mathcal L_1 \rightarrow \mathcal O_X \rightarrow 0\]
of $(1.1.2)$ implies that for each $n\geq 0$ we have a natural filtration
\[\tag{\textbf{1.1.4}} \mathcal L_n=A^0\mathcal L_n \supseteq A^1\mathcal L_n \supseteq...\supseteq A^n\mathcal L_n \supseteq A^{n+1}\mathcal L_n=0\]
of $\mathcal L_n$ by subvector bundles stable under $\nabla_n$ and with quotients given by
\[\tag{\textbf{1.1.5}} \qquad A^i\mathcal L_n/A^{i+1}\mathcal L_n \simeq \mathrm{Sym}^i_{\mathcal O_X}\mathcal H_X \simeq \pi^*\mathrm{Sym}^i_{\mathcal O_S}\mathcal H, \quad i=0,...,n.\]
For this one sets
\[\tag{\textbf{1.1.6}} A^i\mathcal L_n:=\mathrm{im}(\mathrm{Sym}^i_{\mathcal O_X}\mathcal H_X \otimes_{\mathcal O_X} \mathrm{Sym}^{n-i}_{\mathcal O_X}\mathcal L_1 \xrightarrow{\mathrm{mult}} \mathrm{Sym}^n_{\mathcal O_X}\mathcal L_1),\]
which defines the filtration in $(1.1.4)$. The (horizontal) isomorphisms
\[\tag{\textbf{1.1.7}} \mathrm{Sym}^i_{\mathcal O_X}\mathcal H_X \xrightarrow{\sim} A^i\mathcal L_n/A^{i+1}\mathcal L_n\]
are first defined locally by choosing a local section $s$ of $\mathcal L_1$ mapping to $1$ in $(1.1.2)$ and then setting
\[s^i \mapsto s^i\cdot \underbrace{s\cdot ... \cdot s}_{(n-i)-times} \ \textrm{mod} \ A^{i+1}\mathcal L_n \ , \ \textrm{for local sections} \ s^i \ \textrm{of} \ \mathrm{Sym}^i_{\mathcal O_X}\mathcal H_X.\]
One can check that this is independent of the chosen splitting of $(1.1.2)$ and that the so defined global map $(1.1.7)$ is indeed an isomorphism respecting the connections.
\begin{remark}
The filtration we have just constructed can equivalently be described as follows:\\
Set $A^0\mathcal L_n:=\mathcal L_n$ and
\begin{equation*}
A^i\mathcal L_n:=\ker(\mathcal L_n \rightarrow \mathcal L_{i-1}), \quad i=1,...n+1,
\end{equation*}
where the arrow is given by iterated composition of the transition maps. These $A^i\mathcal L_n$ are exactly those defined in $(1.1.6)$.\\
Furthermore, consider $\frac{1}{(n-i)!}$-times the projection $\mathcal L_n \rightarrow \mathcal L_i$. In view of the exact sequence $(1.1.3)$ this map induces a surjective arrow $A^i\mathcal L_n \rightarrow \mathrm{Sym}^i_{\mathcal O_X}\mathcal H_X$ with kernel $A^{i+1}\mathcal L_n$. The obtained isomorphism
\[\qquad A^i\mathcal L_n/A^{i+1}\mathcal L_n \xrightarrow{\sim} \mathrm{Sym}^i_{\mathcal O_X}\mathcal H_X\]
then is precisely the inverse of $(1.1.7)$.
\end{remark}
\subsubsection{Compatibility with base change}
The logarithm sheaves as defined above behave naturally under base change in the following sense:\\
Let $S'$ be another connected scheme which is smooth, separated and of finite type over $\Spec(\Q)$ and $f: S' \rightarrow S$ a $\Q$-morphism. In the induced cartesian diagram
\begin{equation*} \tag{\textbf{1.1.8}}\begin{split}
\begin{xy}
\xymatrix{
X' \ar[r]^{\pi'} \ar[d]_{g} & S' \ar[d]^{f} \\
X \ar[r]^{\pi} & S}
\end{xy}
\end{split}
\end{equation*}
we view $X'$ as abelian $S'$-scheme. Using the canonical horizontal isomorphism
\[f^* H^1_{\mathrm{dR}}(X/S) \xrightarrow{\sim} H^1_{\mathrm{dR}}(X'/S')\]
we obtain by pullback of $(1.1.2)$ resp. of $\varphi_1$ via $g$ resp. via $f$ a pair as in Def. 1.1.1 (ii) for the abelian scheme $X'/S'$. This is straightforward to check. We thus see (cf. Prop. 1.4.7 for another viewpoint):
\begin{lemma}
The pullback of $(\mathcal L_n,\nabla_n,\varphi_n)$ along the arrows of the cartesian diagram $(1.1.8)$ induces the datum of the $n$-th logarithm sheaf for $X'/S'/\Q$ \qquad \qed
\end{lemma}

\markright{\uppercase{The formalism of the logarithm sheaves and the elliptic...}}
\section{The de Rham cohomology of the logarithm sheaves}
\markright{\uppercase{The formalism of the logarithm sheaves and the elliptic...}}
\subsection{The computation of $H^i_{\mathrm{dR}}(X/S,\mathcal L_n)$}
We calculate the de Rham cohomology sheaves
\[H^i_{\mathrm{dR}}(X/S,\mathcal L_n)\]
for $n \geq 1$ and $0 \leq i \leq 2g$. They are equipped with the Gauß-Manin connection relative $\Spec(\Q)$.\\
\newline
The formal frame of the subsequent arguments is essentially drawn from \cite{Ki4}, Prop. 1.1.3 (a), where the higher direct images of the logarithm pro-sheaf in the $\ell$-adic setting are computed. We need to adjust the formalism to our situation of de Rham realization and additionally evaluate the occurring spectral sequence in more detail to get the de Rham cohomology of each $\mathcal L_n$.\\
\newline
Recall (cf. \cite{Kat2}, Rem. $(3.1)$ or the proof of \cite{Har}, Ch. III, Prop. $(4.2)$) that the functor between abelian categories
\[ H^i_{\mathrm{dR}}(X/S,-): \Mod_{\mathrm{qc}}(\mathcal D_{X/\Q}) \rightarrow \Mod_{\mathrm{qc}}(\mathcal D_{S/\Q})\]
may be viewed as the $i$-th right derivation of the left exact functor
\[ H^0_{\mathrm{dR}}(X/S,-): \Mod_{\mathrm{qc}}(\mathcal D_{X/\Q}) \rightarrow \Mod_{\mathrm{qc}}(\mathcal D_{S/\Q}).\]
The filtration $A^\bullet \mathcal L_n$ of $(1.1.4)$ then yields (cf. \cite{EGAIII}, Ch. 0, 13.6) a spectral sequence in $\Mod_{\mathrm{qc}}(\mathcal D_{S/\Q})$
\[E_1^{p,q}=H^{p+q}_{\mathrm{dR}}(X/S, gr^p A^\bullet \mathcal L_n) \Rightarrow E^{p+q}=H^{p+q}_{\mathrm{dR}}(X/S, \mathcal L_n),\]
where with $(1.1.5)$ we see that $E_1^{p,q} \simeq H^{p+q}_{\mathrm{dR}}(X/S) \otimes_{\mathcal O_S} \mathrm{Sym}^p_{\mathcal O_S} \mathcal H$ for $0\leq p \leq n$ and zero otherwise.\\
\newline
The differential $d_1^{pq}: E_1^{p,q} \rightarrow E_1^{p+1,q}$ is given by the connecting morphism induced by applying $H_{\mathrm{dR}}^{p+q}(X/S,-)$ to the exact sequence
\[0 \rightarrow \mathrm{Sym}^{p+1}_{\mathcal O_X}\mathcal H_X \rightarrow A^p \mathcal L_n /A^{p+2}\mathcal L_n \rightarrow \mathrm{Sym}^p_{\mathcal O_X} \mathcal H_X \rightarrow 0.\]
This morphism in turn is equal to the composition

\[\tag{\textbf{1.2.1}}
\begin{split}
H^{p+q}_{\mathrm{dR}}(X/S) \otimes_{\mathcal O_S} \mathrm{Sym}^p_{\mathcal O_S}\mathcal H &\rightarrow H_{\mathrm{dR}}^{p+q+1}(X/S) \otimes_{\mathcal O_S} \mathcal H \otimes_{\mathcal O_S} \mathrm{Sym}^p_{\mathcal O_S}\mathcal H \\
&\rightarrow H_{\mathrm{dR}}^{p+q+1}(X/S) \otimes_{\mathcal O_S} \mathrm{Sym}^{p+1}_{\mathcal O_S}\mathcal H\end{split}\]in which the first arrow comes from the canonical map $\mathcal O_S \rightarrow \mathcal H^\vee \otimes_{\mathcal O_S} \mathcal H$ together with cup product and the second is given by multiplication: one deduces this straightforwardly with the same argument as in \cite{Ki4}, proof of Prop. 1.1.3 (a), observing the explicit construction of the filtration $(1.1.4)$-$(1.1.7)$ and the fact that (by definition) $\mathcal Log^1$ goes to the canonical map under the surjection of $(1.1.1)$.\\
\newline
Using the canonical (horizontal) identifications\footnote{For this one composes the isomorphism
\[H^i_{\mathrm{dR}}(X/S) \simeq \Bigg(\bigwedge^{2g-i} H^1_{\mathrm{dR}}(X/S)\Bigg)^\vee\]
from the beginning of this chapter with the usual natural identification
\[\Bigg(\bigwedge^{2g-i} H^1_{\mathrm{dR}}(X/S)\Bigg)^\vee\simeq \bigwedge^{2g-i} \mathcal H,\]
locally defined as in \cite{La}, Ch. XIX, §1, Prop. 1.5.}
\[H_{\mathrm{dR}}^i(X/S) \simeq \bigwedge^{2g-i} \mathcal H\]
the differential $d^{p,q}_1$ writes as a map
\[d^{p,q}_1: \bigwedge^{2g-p-q} \mathcal H \otimes_{\mathcal O_S} \mathrm{Sym}^p_{\mathcal O_S}\mathcal H \rightarrow \bigwedge^{2g-p-q-1} \mathcal H \otimes_{\mathcal O_S} \mathrm{Sym}^{p+1}_{\mathcal O_S}\mathcal H,\]
and the $q$-th line of the first sheet of the spectral sequence thus becomes a complex of the form
\[0 \rightarrow \bigwedge^{2g-q} \mathcal H \rightarrow \bigwedge^{2g-q-1} \mathcal H \otimes_{\mathcal O_S} \mathrm{Sym}^1_{\mathcal O_S}\mathcal H \rightarrow...\rightarrow \bigwedge^{2g-q-n} \mathcal H \otimes_{\mathcal O_S} \mathrm{Sym}^n_{\mathcal O_S}\mathcal H\rightarrow 0.\]
With the explicit knowledge of the maps in $(1.2.1)$ it is routinely verified that the preceding differentials are precisely the Koszul differentials in the Koszul complex
\[0 \rightarrow \bigwedge^{2g-q} \mathcal H \rightarrow \bigwedge^{2g-q-1} \mathcal H \otimes_{\mathcal O_S} \mathrm{Sym}^1_{\mathcal O_S}\mathcal H \rightarrow...\rightarrow \mathrm{Sym}_{\mathcal O_S}^{2g-q}\mathcal H\rightarrow 0\]
which is exact for all $q < 2g$.\footnote{We refer to \cite{La}, Ch. XXI, §4, Thm. 4.13 and Cor. 4.14, for the explicit formula of the Koszul differential and the acyclicity of the Koszul complex in the case of free modules of finite rank over a ring; this globalizes to our situation of a vector bundle. To get the Koszul exact sequence of above one sets $r=2g$ and $n=2g-q$ in ibid., Cor. 4.14.}\\
\newline
This admits a detailed evaluation of the second sheet of the spectral sequence. The result is:\\
\newline
$E_2^{n,-n} \simeq \mathrm{Sym}^n_{\mathcal O_S} \mathcal H$\\
$E_2^{n,q} \simeq \frac{H^{n+q}_{{\mathrm{dR}}}(X/S) \otimes_{\mathcal O_S} \mathrm{Sym}^n_{\mathcal O_S}\mathcal H}{\mathrm{im}(H^{n+q-1}_{{\mathrm{dR}}}(X/S) \otimes_{\mathcal O_S} \mathrm{Sym}^{n-1}_{\mathcal O_S}\mathcal H \rightarrow H^{n+q}_{{\mathrm{dR}}}(X/S) \otimes_{\mathcal O_S} \mathrm{Sym}^n_{\mathcal O_S}\mathcal H)} \quad \forall q \in ]-n;2g-n[$\\
$E_2^{0,2g} \simeq \mathcal O_S$\\
\newline
and zero for all other cases; the map appearing in the $E_2^{n,q}$-term is given as in $(1.2.1)$.\\
Of course, all the isomorphisms take place in the category $\Mod_{\mathrm{qc}}(\mathcal D_{S/\Q})$. \\
We obtain that the spectral sequence degenerates at $r=2$: note that there are only two nontrivial columns, namely $p=0$ and $p=n$, and that for $p=0$ we only have a non-zero entry for $q=2g$; thus the only possible non-zero differential for $r\geq 2$ could appear if $r=n$, namely $d_n^{0,2g}: E_n^{0,2g} \rightarrow E_n^{n,2g-n+1}$; but the last term is zero.\\
\newline
From this one deduces for each $n\geq 1$ the following

\begin{theorem}
(i) We have the following isomorphisms in $\Mod_{\mathrm{qc}}(\mathcal D_{S/\Q})$:\\
\newline
$\bullet \quad \mathrm{Sym}^n_{\mathcal O_S} \mathcal H \simeq H^0_{\mathrm{dR}}(X/S, \mathcal L_n)$,\\
\newline
induced by the inclusion $\mathrm{Sym}^n_{\mathcal O_X}\mathcal H_X \hookrightarrow \mathcal L_n$ given as in $(1.1.3)$.\\
\newline
$\bullet \quad \frac{H^i_{\mathrm{dR}}(X/S) \otimes_{\mathcal O_S} \mathrm{Sym}^n_{\mathcal O_S}\mathcal H}{\mathrm{im} (H^{i-1}_{{\mathrm{dR}}}(X/S) \otimes_{\mathcal O_S} \mathrm{Sym}^{n-1}_{\mathcal O_S}\mathcal H \rightarrow H^i_{\mathrm{dR}}(X/S) \otimes_{\mathcal O_S} \mathrm{Sym}^n_{\mathcal O_S} \mathcal H)} \simeq H_{\mathrm{dR}}^i(X/S, \mathcal L_n), \ \ \textrm{where} \ i\neq 0,2g$,\\
\newline
induced by the inclusion $\mathrm{Sym}^n_{\mathcal O_X}\mathcal H_X \hookrightarrow \mathcal L_n$; the lower map is defined as in $(1.2.1)$.\\
\newline
$\bullet \quad H_{\mathrm{dR}}^{2g}(X/S, \mathcal L_n) \simeq \mathcal O_S$,\\
\newline
induced by the projection $\mathcal L_n \rightarrow \mathcal L_{n-1} \rightarrow ... \rightarrow \mathcal O_X$ and the trace isomorphism for $H^{2g}_{\mathrm{dR}}(X/S)$.\\
\newline
(ii) For all $i<2g$ the transition map $\mathcal L_n \rightarrow \mathcal L_{n-1}$ induces the zero morphism
\[H^i_{\mathrm{dR}}(X/S, \mathcal L_n) \xrightarrow{0} H^i_{\mathrm{dR}}(X/S, \mathcal L_{n-1}),\]
and for $i=2g$ it induces an isomorphism
	\[H^{2g}_{\mathrm{dR}}(X/S, \mathcal L_n) \xrightarrow{\sim} H^{2g}_{\mathrm{dR}}(X/S, \mathcal L_{n-1})\]
which becomes the identity when taking into account the identifications with $\mathcal O_S$ in (i).\\
\qed
\end{theorem}

\begin{Remark}
All the sheaves $H^i_{\mathrm{dR}}(X/S, \mathcal L_n)$ are vector bundles on $S$: as the Gauß-Manin connection relative $\Spec(\Q)$ operates on them it suffices to check their coherence (cf. \cite{Bert-Og}, §2, Note 2.17). This in turn follows from the spectral sequence $(0.2.2)$ and the properness of $X/S$ (cf. 0.2.2).
\end{Remark}

\subsection{The computation of $H^i_{\mathrm{dR}}(U/S,\mathcal L_n)$}
We set $U:=X -\epsilon(S)$ and write $j: U \rightarrow X$ for the associated open immersion. The structure morphism of $U$ as an $S$-scheme will be denoted by $\pi_U$.
\[\begin{xy}
\xymatrix{
\qquad \qquad & S \ar[dr]_{\id} \ar[r]^{\epsilon} & X  \ar[d]^{\pi} & U \ar[dl]^{\pi_U} \ar[l]_{j}\\
\qquad \qquad &  & S & &}
\end{xy}\]

For an object $\mathcal E \in \Mod_{\mathrm{qc}}(\mathcal D_{X/\Q})$ we write $\mathcal E_{|U} \in \Mod_{\mathrm{qc}}(\mathcal D_{U/S})$ for its restriction to $U$ and $H^i_{\mathrm{dR}}(U/S, \mathcal E) \in \Mod_{\mathrm{qc}}(\mathcal D_{S/\Q})$ for the $i$-th de Rham cohomology sheaf of $\mathcal E_{|U}$, equipped with the Gauß-Manin connection relative $\Spec(\Q)$.\\
\newline
We now want to calculate
\[H^i_{\mathrm{dR}}(U/S,\mathcal L_n)\]
for $n \geq 0$ and $0 \leq i \leq 2g$.\\
\newline
For each $n \geq 0$ we have the canonical distinguished triangle in $D^b_{\mathrm{qc}}(\mathcal D_{X/\Q})$ (cf. $(0.2.5)$)
\[\epsilon_+ \epsilon^{!} \mathcal L_n \rightarrow \mathcal L_n \rightarrow j_+ \mathcal L_{n|U}.\]
Applying $\pi_+$ and observing $\epsilon^!\mathcal L_n \simeq \epsilon^* \mathcal L_n [-g]$ (cf. 0.2.3) gives the distinguished triangle in $D^b_{\mathrm{qc}}(\mathcal D_{S/\Q})$
\[\epsilon^*\mathcal L_n [-g] \rightarrow \pi_+ \mathcal L_n \rightarrow (\pi_U)_+ \mathcal L_{n|U}.\]The $k$-th cohomology sheaf of $\pi_+ \mathcal L_n$ resp. $(\pi_U)_+ \mathcal L_{n|U}$ is canonically isomorphic to $H^{k+g}_{{\mathrm{dR}}}(X/S, \mathcal L_n)$ resp. $H_{\mathrm{dR}}^{k+g}(U/S, \mathcal L_n)$ in $\Mod_{\mathrm{qc}}(\mathcal D_{S/\Q})$ (the proof is as in \cite{Dim-Ma-Sa-Sai}, Prop. 1.4).\\
The long exact cohomology sequence for the preceding distinguished triangle hence shows
\begin{lemma}
Let $n \geq 0$. Then:\\
\newline
(i) For $i \neq 2g-1,2g$ the canonical map $H_{\mathrm{dR}}^i(X/S, \mathcal L_n) \rightarrow H_{\mathrm{dR}}^i(U/S, \mathcal L_n)$ is an isomorphism.\\
\newline
(ii) One has an exact sequence of $\mathcal D_{S/\Q}$-modules
\[0 \rightarrow H_{\mathrm{dR}}^{2g-1}(X/S, \mathcal L_n) \xrightarrow{\mathrm{can}} H_{\mathrm{dR}}^{2g-1}(U/S, \mathcal L_n) \rightarrow \epsilon^*\mathcal L_n \xrightarrow{\sigma_n} H_{\mathrm{dR}}^{2g}(X/S, \mathcal L_n) \xrightarrow{\mathrm{can}} H_{\mathrm{dR}}^{2g}(U/S, \mathcal L_n) \rightarrow 0\]
\qed
\end{lemma}

\begin{definition}
For each $n\geq 0$ we write
\[\mathrm{Res}^n: H_{\mathrm{dR}}^{2g-1}(U/S, \mathcal L_n) \rightarrow \prod_{k=0}^n \mathrm{Sym}_{\mathcal O_S}^k \mathcal H\]
for the $\mathcal D_{S/\Q}$-linear arrow induced by the exact sequence of Lemma 1.2.3 (ii) and by the splitting
\[\varphi_n: \epsilon^* \mathcal L_n \simeq \prod_{k=0}^n \mathrm{Sym}_{\mathcal O_S}^k \mathcal H.\]
The claim of the following lemma implies that for $n\geq 1$ the map $\mathrm{Res}^n$ factors in the form
\[\mathrm{Res}^n: H_{\mathrm{dR}}^{2g-1}(U/S, \mathcal L_n) \rightarrow \prod_{k=1}^n \mathrm{Sym}_{\mathcal O_S}^k \mathcal H \subseteq \prod_{k=0}^n \mathrm{Sym}_{\mathcal O_S}^k \mathcal H,\]
and we write $\mathrm{Res}^n$ also for the first of these arrows.
\end{definition}

\begin{lemma}
Let $n\geq 0$. Then, under the identifications $\epsilon^* \mathcal L_n \simeq \prod_{k=0}^n \mathrm{Sym}_{\mathcal O_S}^k \mathcal H$ and $H_{\mathrm{dR}}^{2g}(X/S, \mathcal L_n) \simeq \mathcal O_S$ given by the splitting resp. by Thm. 1.2.1 (i), the map $\sigma_n$ appearing in Lemma 1.2.3 (ii) is $n!$-times the natural projection.
\end{lemma}
As the proof is a bit long and technical we postpone it to the next subsection.\\
\newline
With this result at hand we deduce from Lemma 1.2.3 (ii) that $H^{2g}_{\mathrm{dR}}(U/S, \mathcal L_n) = 0$ for all $n\geq 0$, that the map $H_{\mathrm{dR}}^{2g-1}(X/S) \xrightarrow{\mathrm{can}} H_{\mathrm{dR}}^{2g-1}(U/S)$ is an isomorphism and that for $n \geq 1$ there is an exact sequence of $\mathcal O_S$-vector bundles with integrable $\Q$-connection:
\[\tag{\textbf{1.2.2}} 0 \rightarrow H_{\mathrm{dR}}^{2g-1}(X/S, \mathcal L_n) \xrightarrow{\mathrm{can}} H_{\mathrm{dR}}^{2g-1}(U/S, \mathcal L_n) \xrightarrow{\mathrm{Res}^n} \prod_{k=1}^n \mathrm{Sym}^k_{\mathcal O_S} \mathcal H \rightarrow 0.\]
The long exact sequence for the derived functors of $H^0_{\mathrm{dR}}(U/S,-)$ applied to the restriction
\[0 \rightarrow \mathrm{Sym}^{n+1}_{\mathcal O_U} \mathcal H_U \rightarrow \mathcal L_{n+1|U} \rightarrow \mathcal L_{n|U} \rightarrow 0\]
of $(1.1.3)$ to $U$ and the vanishing of $H^{2g}_{\mathrm{dR}}(U/S, \mathrm{Sym}^{n+1}_{\mathcal O_U} \mathcal H_U) \simeq H^{2g}_{\mathrm{dR}}(U/S) \otimes_{\mathcal O_S}\mathrm{Sym}^{n+1}_{\mathcal O_S}\mathcal H$ shows:
\begin{lemma} 
For each $n\geq 0$ the transition map $H^{2g-1}_{\mathrm{dR}}(U/S,\mathcal L_{n+1}) \rightarrow H^{2g-1}_{\mathrm{dR}}(U/S,\mathcal L_n)$ is surjective. \\
\qed
\end{lemma} 
As the transition map $H^{2g-1}_{\mathrm{dR}}(X/S, \mathcal L_{n+1}) \rightarrow H^{2g-1}_{\mathrm{dR}}(X/S, \mathcal L_n)$ is zero (cf. Thm. 1.2.1 (ii)) the surjection of Lemma 1.2.6 vanishes on the subobject $H^{2g-1}_{\mathrm{dR}}(X/S, \mathcal L_{n+1})$, from which by $(1.2.2)$ (used for $n+1$ instead of $n$) we get an induced horizontal surjection of $\mathcal O_S$-vector bundles:
\[\tag{\textbf{1.2.3}} \prod_{k=1}^{n+1} \mathrm{Sym}^k_{\mathcal O_S}\mathcal H \rightarrow H^{2g-1}_{\mathrm{dR}}(U/S, \mathcal L_n).\]
\begin{lemma}
For each $n\geq 0$ the surjection of $(1.2.3)$ is an isomorphism.
\end{lemma}
\begin{proof}
As a surjection between two $\mathcal O_S$-vector bundles of the same rank is an isomorphism we only need to show that both involved bundles have the same $\mathcal O_S$-rank. In view of the exact sequence $(1.2.2)$ it is enough to show that $H^{2g-1}_{\mathrm{dR}}(X/S, \mathcal L_n)$ has the same rank as $\mathrm{Sym}_{\mathcal O_S}^{n+1}\mathcal H$. But these two sheaves are indeed isomorphic. We give two arguments for this fact:\\
First, look at the final terms of the de Rham cohomology sequence for the exact sequence
\[0 \rightarrow \mathrm{Sym}^{n+1}_{\mathcal O_X} \mathcal H_X \rightarrow \mathcal L_{n+1} \rightarrow \mathcal L_n \rightarrow 0.\]
of $(1.1.3)$. Then Thm. 1.1.1 (ii) implies that the connecting arrow
\[H^{2g-1}_{\mathrm{dR}}(X/S,\mathcal L_n)\rightarrow H^{2g}_{\mathrm{dR}}(X/S) \otimes_{\mathcal O_S} \mathrm{Sym}^{n+1}_{\mathcal O_S}\mathcal H \simeq \mathrm{Sym}^{n+1}_{\mathcal O_S}\mathcal H\]
is an isomorphism, which shows the claim.\\
Alternatively, it follows from Thm. 1.1.1 (i) that $H^{2g-1}_{\mathrm{dR}}(X/S,\mathcal L_n)$ identifies for $n\geq 1$ with
\[\frac{\mathcal H \otimes_{\mathcal O_S} \mathrm{Sym}^n_{\mathcal O_S} \mathcal H}{\mathrm{im}( \bigwedge^2 \mathcal H \otimes_{\mathcal O_S}\mathrm{Sym}_{\mathcal O_S}^{n-1} \mathcal H \rightarrow \mathcal H \otimes_{\mathcal O_S} \mathrm{Sym}^n_{\mathcal O_S} \mathcal H)},\]
where the map in the denominator is the Koszul differential and hence fits into an exact sequence
\[\bigwedge^2 \mathcal H \otimes_{\mathcal O_S}\mathrm{Sym}_{\mathcal O_S}^{n-1} \mathcal H \rightarrow \mathcal H \otimes_{\mathcal O_S} \mathrm{Sym}^n_{\mathcal O_S} \mathcal H \rightarrow \mathrm{Sym}_{\mathcal O_S}^{n+1} \mathcal H \rightarrow 0.\]
This shows the claim for $n\geq 1$, and for $n=0$ we know it from the beginning of this chapter.
\end{proof}
The proof of Lemma 1.2.7 has in particular shown the following supplement to Thm. 1.2.1:
\begin{Korollar}
For each $n\geq 1$ we have a horizontal isomorphism
\[\mathrm{Sym}_{\mathcal O_S}^{n+1} \mathcal H \xrightarrow{\sim} H^{2g-1}_{\mathrm{dR}}(X/S,\mathcal L_n),\]
determined by the commutative diagram
\begin{equation*}
\begin{xy}
\xymatrix@C-0.3cm{
 & & \mathcal H \otimes_{\mathcal O_S} \mathrm{Sym}_{\mathcal O_S}^n \mathcal H \ar[rd]^{\mathrm{mult}} \ar[rr] &  & H^{2g-1}_{\mathrm{dR}}(X/S, \mathcal L_n) \\
 & & &\mathrm{Sym}^{n+1}_{\mathcal O_S}\mathcal H \ar[ur]^{\sim} &}
\end{xy}
\end{equation*}
where the upper arrow is induced by the map $\mathrm{Sym}_{\mathcal O_X}^n \mathcal H_X \rightarrow \mathcal L_n$ in $(1.1.3)$ and the canonical identification $H^{2g-1}_{\mathrm{dR}}(X/S) \simeq \mathcal H$. \\
\qed
\end{Korollar}

Observing Lemma 1.1.5 in part (iv) we deduce in sum for every $n\geq 0$ the following

\begin{theorem}
(i) For each $i \neq 2g-1,2g$ the canonical map $H_{\mathrm{dR}}^i(X/S, \mathcal L_n) \rightarrow H_{\mathrm{dR}}^i(U/S, \mathcal L_n)$ is an isomorphism.\\
\newline
(ii) $H^{2g}_{\mathrm{dR}}(U/S, \mathcal L_n)=0.$\\
\newline
In particular, the transition map in de Rham cohomology
\[H_{\mathrm{dR}}^i(U/S, \mathcal L_{n+1}) \rightarrow H_{\mathrm{dR}}^i(U/S, \mathcal L_n)\]
is zero for all $i \neq 2g-1$.\\
\newline
(iii) The canonical map $H_{\mathrm{dR}}^{2g-1}(X/S) \rightarrow H_{\mathrm{dR}}^{2g-1}(U/S)$ is an isomorphism.\\
\newline
(iv) We have an isomorphism
	\[\prod_{k=1}^{n+1} \mathrm{Sym}^k_{\mathcal O_S}\mathcal H \xrightarrow{\sim} H^{2g-1}_{\mathrm{dR}}(U/S, \mathcal L_n)\]
determined by the commutative diagram
\begin{equation*}
\begin{xy}
\xymatrix@C-0.3cm{
 & & H^{2g-1}_{\mathrm{dR}}(U/S, \mathcal L_{n+1}) \ar[rd]^{\mathrm{Res}^{n+1}} \ar[rr] &  & H^{2g-1}_{\mathrm{dR}}(U/S, \mathcal L_n) \\
 & & & \prod_{k=1}^{n+1} \mathrm{Sym}^k_{\mathcal O_S}\mathcal H \ar[ur]^{\sim} &}
\end{xy}
\end{equation*}
Under these identifications the transition map in de Rham cohomology for $n\geq 1$
\[H_{\mathrm{dR}}^{2g-1}(U/S, \mathcal L_n) \rightarrow H_{\mathrm{dR}}^{2g-1}(U/S, \mathcal L_{n-1})\]
resp. the map $\mathrm{Res}^n$ induces an arrow
\[\prod_{k=1}^{n+1} \mathrm{Sym}^k_{\mathcal O_S}\mathcal H \rightarrow \prod_{k=1}^n \mathrm{Sym}^k_{\mathcal O_S}\mathcal H\]
which is given explicitly as
\[\qquad (h_1,h_2,...,h_{n+1}) \mapsto (n \cdot h_1,(n-1)\cdot h_2,...,h_n) \quad \textrm{with}\ h_k \in \mathrm{Sym}^k_{\mathcal O_S}\mathcal H, \ k=1,...,n+1.\]
\qed
\end{theorem}

\begin{remark}
The preceding theorem and Rem. 1.2.2 imply that the sheaves $H^i_{\mathrm{dR}}(U/S,\mathcal L_n)$ are all vector bundles.
\end{remark}

\subsection{Proof of Lemma 1.2.5}
We recall what we have to show:\\
Let $n\geq 0$ and consider the arrow in $D^b_{\mathrm{qc}}(\mathcal D_{X/\Q})$
\[\epsilon_+ \epsilon^* \mathcal L_n[-g] \rightarrow \mathcal L_n\]
appearing in the localization triangle for $\mathcal L_n$. By application of $\pi_+$ we get the arrow in $D^b_{\mathrm{qc}}(\mathcal D_{S/\Q})$
\[\tag{\textbf{1.2.4}} \epsilon^* \mathcal L_n [-g] \rightarrow \pi_+ \mathcal L_n,\]
and we need to verify that under the identifications
\[\epsilon^*\mathcal L_n \simeq \prod_{k=0}^n \mathrm{Sym}_{\mathcal O_S}^k \mathcal H, \quad H^g \pi_+ \mathcal L_n \simeq H^{2g}_{\mathrm{dR}}(X/S, \mathcal L_n) \xrightarrow{\sim} \mathcal O_S\]
we obtain $n!$-times the natural projection when taking $g$-th cohomology in $(1.2.4)$.\\
The isomorphism $H^{2g}_{\mathrm{dR}}(X/S,\mathcal L_n) \xrightarrow{\sim}\mathcal O_S$ is induced by the composition $\mathcal L_n \rightarrow \mathcal L_{n-1} \rightarrow ... \rightarrow \mathcal O_X$ together with the trace map $H^{2g}_{\mathrm{dR}}(X/S) \xrightarrow{\sim} \mathcal O_S$, and for the identification $H^g (\pi_+ \mathcal L_n) \simeq H^{2g}_{\mathrm{dR}}(X/S, \mathcal L_n)$ cf. \cite{Dim-Ma-Sa-Sai}, Prop. 1.4. To show the claim we don't need to assume the horizontality of the trace map; in fact, this will come out as a side result below.\\
\newline
\underline{Proof}:\\
We first explain that it suffices to verify the lemma for $n=0$:\\
The functoriality of the localization triangle induces a commutative diagram in $D^b_{\mathrm{qc}}(\mathcal D_{X/\Q})$
\begin{equation*}
\begin{xy}
\xymatrix{
\epsilon_+ \epsilon^* \mathcal L_n[-g] \ar[r]\ar[d] & \mathcal L_n \ar[d] \\
\epsilon_+ \mathcal O_S [-g] \ar[r] & \mathcal O_X}
\end{xy}
\end{equation*}
from which we obviously get commutative diagrams in $D^b_{\mathrm{qc}}(\mathcal D_{S/\Q})$ resp. $\Mod_{\mathrm{qc}}(\mathcal D_{S/\Q})$
\begin{equation*}
\begin{xy}
\xymatrix{
\epsilon^* \mathcal L_n \ar[r]\ar[d] & \pi_+ \mathcal L_n[g] \ar[d] & & \prod_{k=0}^n \mathrm{Sym}^k_{\mathcal O_S} \mathcal H \ar[r]\ar[d]_{n! \cdot \mathrm{can}} & H_{\mathrm{dR}}^{2g}(X/S, \mathcal L_n) \ar[d]^{\sim} \\
\mathcal O_S \ar[r] & \pi_+ \mathcal O_X[g] & & \mathcal O_S \ar[r] & H^{2g}_{\mathrm{dR}}(X/S)}
\end{xy}
\end{equation*}
In the right diagram the right vertical arrow is precisely the morphism in de Rham cohomology induced from $\mathcal L_n \rightarrow \mathcal L_{n-1} \rightarrow ... \rightarrow \mathcal O_X$ (and is an isomorphism by Thm. 1.2.1 (i)). The left vertical arrow is given by $n!$-times the canonical projection, as follows from Lemma 1.1.5.\\
Knowing the claim for $n=0$ means that we may add a lower commutative diagram in $\Mod_{\mathrm{qc}}(\mathcal O_S)$
\begin{equation*}\tag{\textbf{1.2.5}} \begin{split}
\begin{xy}
\xymatrix{
\prod_{k=0}^n \mathrm{Sym}^k_{\mathcal O_S} \mathcal H \ar[r]\ar[d]_{n! \cdot \mathrm{pr}} & H_{\mathrm{dR}}^{2g}(X/S, \mathcal L_n) \ar[d]^{\sim} \\
\mathcal O_S \ar[r] \ar[d]_{\id} & H^{2g}_{\mathrm{dR}}(X/S) \ar[d]^{\mathrm{tr}}\\
\mathcal O_S \ar[r]^{\id} & \mathcal O_S}
\end{xy}
\end{split}
\end{equation*}
such that also the frame diagram commutes. But this frame commutativity is the claim of the lemma, and thus we are indeed reduced to the case $n=0$.\\
\newline
It hence remains to show the commutativity of the lower diagram
\begin{equation*}\tag{\textbf{1.2.6}} \begin{split}
\begin{xy}
\xymatrix{
\mathcal O_S \ar[r]^{\psi \qquad} \ar[d]_{\id} & H^{2g}_{\mathrm{dR}}(X/S) \ar[d]^{\mathrm{tr}}\\
\mathcal O_S \ar[r]^{\id} & \mathcal O_S}
\end{xy}
\end{split}
\end{equation*}
of $\mathcal O_S$-linear maps in $(1.2.5)$. As the occurring arrow $\psi$ is (by construction) $\mathcal D_{S/\Q}$-linear, the commutativity of $(1.2.6)$ then clearly also implies the horizontality of the trace isomorphism.\\
\newline
The $\mathcal D_{S/\Q}$-linear map $\psi$ is given by applying the functor $\pi_+: D^b_{\mathrm{qc}}(\mathcal D_{X/\Q}) \rightarrow D^b_{\mathrm{qc}}(\mathcal D_{S/\Q})$ to the arrow $\epsilon_+\mathcal O_S[-g] \rightarrow \mathcal O_X$ of the localization triangle for $\mathcal O_X$, shifting the obtained map $\mathcal O_S[-g] \rightarrow \pi_+\mathcal O_X$ by $[g]$, taking $0$-th cohomology and identifying $H^0(\pi_+\mathcal O_X[g])\simeq H^g(\pi_+\mathcal O_X) \simeq H^{2g}_{\mathrm{dR}}(X/S)$.\\
\newline
If we replace $\Q$ by $S$ in the procedure of the previous passage and observe the identification
\[ \tag{\textbf{1.2.7}} \begin{split} H^g(\pi_+\mathcal O_X) = H^gR\pi_*(\mathcal D_{(S \leftarrow X)/S} \otimes^L_{\mathcal D_{X/S}}\mathcal O_X)=H^gR\pi_*(\Omega^g_{X/S}\otimes^L_{\mathcal D_{X/S}}\mathcal O_X)\\
\simeq H^gR\pi_*(((\Omega^{\bullet}_{X/S}\otimes_{\mathcal O_X} \mathcal D_{X/S})[g])\otimes ^L_{\mathcal D_{X/S}}\mathcal O_X) \simeq H^{2g}R\pi_*(\Omega^{\bullet}_{X/S})\simeq H^{2g}_{\mathrm{dR}}(X/S), \end{split}\]
then the obtained $\mathcal O_S$-linear arrow $\mathcal O_S \rightarrow H^{2g}_{\mathrm{dR}}(X/S)$ is again the map $\psi$ in $(1.2.6)$, only that we have forgotten the $\mathcal D_{S/\Q}$-linear structures (one can check this).\\
\newline
We now define arrows
\[Ad_{\epsilon}: \epsilon_+ \mathcal O_S \rightarrow \mathcal O_X[g] \quad \textrm{in} \ D^b_{\mathrm{qc}}(\mathcal D_{X/S})\]
resp.
\[Ad_{\pi}: \pi_+ \mathcal O_X[g] \rightarrow \mathcal O_S \quad \textrm{in} \ D^b_{\mathrm{qc}}(\mathcal O_S) \quad \]
as follows: $Ad_{\epsilon}$ comes from shifting by $[g]$ the map $\epsilon_+ \epsilon^!\mathcal O_X \rightarrow \mathcal O_X$ appearing in the localization triangle (relative $S$) for $\mathcal O_X$. Hence, if we apply $\pi_+$ to $Ad_{\epsilon}$, take $0$-th cohomology and identify $H^0(\pi_+\mathcal O_X[g]) \simeq H^g(\pi_+\mathcal O_X)\simeq H^{2g}_{\mathrm{dR}}(X/S)$ as we just did, then we get the morphism $\psi$.\\
The map $Ad_{\pi}$ has the important property that if we take $0$-th cohomology and use $H^0(\pi_+\mathcal O_X[g]) \simeq H^{2g}_{\mathrm{dR}}(X/S)$ it becomes the trace isomorphism; namely, it is constructed as follows:\\
As in $(1.2.7)$ we have canonically in $\Mod_{\mathrm{qc}}(\mathcal O_S)$
\[H^i(\pi_+ \mathcal O_X[g])\simeq H^{2g+i}_{{\mathrm{dR}}}(X/S).\]
From the degeneration of the Hodge-de Rham spectral sequence at the first sheet we further have
\[H^{2g}_{\mathrm{dR}}(X/S) \simeq R^g \pi_* \Omega^g_{X/S},\]
and composition with the Grothendieck trace isomorphism (cf. \cite{Con1}, Ch. I, 1.1)
\[R^g \pi_* \Omega^g_{X/S} \xrightarrow{\sim} \mathcal O_S\]
yields altogether a canonical $\mathcal O_S$-linear arrow \[H^0(\pi_+ \mathcal O_X[g]) \xrightarrow{\sim} \mathcal O_S.\]
As $\pi_+ \mathcal O_X [g]$ has no cohomology above $0$ there is a canonical morphism $\pi_+ \mathcal O_X [g]\rightarrow H^0(\pi_+ \mathcal O_X [g])$ in $D^b_{\mathrm{qc}}(\mathcal O_S)$ which is the identity when taking $0$-th cohomology. We define $Ad_{\pi}$ to be the composition
\[Ad_{\pi}: \pi_+ \mathcal O_X [g] \rightarrow H^0(\pi_+ \mathcal O_X[g]) \xrightarrow{\sim} \mathcal O_S \quad \textrm{in} \ D^b_{\mathrm{qc}}(\mathcal O_S),\]
and by construction it has the announced property to be the trace map when taking $0$-th cohomology and using the identification $H^0(\pi_+\mathcal O_X[g]) \simeq H^{2g}_{\mathrm{dR}}(X/S)$.\\
\newline
Now consider the following composition in $D^b_{\mathrm{qc}}(\mathcal O_S)$
\[\tag{\textbf{1.2.8}} \mathcal O_S \rightarrow \pi_+\mathcal O_X [g] \rightarrow \mathcal O_S\]
in which the first arrow is $\pi_+$ applied to $Ad_{\epsilon}$ and the second is $Ad_{\pi}$. If it is the identity then by taking $0$-th cohomology and using the property of $Ad_{\pi}$ explained above the remaining claim follows.\\
\newline
In the case that $S=\Spec(k)$, with $k$ a field of characteristic zero, the morphisms $Ad_{\epsilon}$ and $Ad_{\pi}$ coincide precisely with the adjunctions (of the same notation) constructed for morphisms between quasi-projective algebraic varieties of characteristic zero in \cite{Me}, Thm. (7.1)\footnote{That these adjunctions coincide with our constructions follows from their alternative description given in \cite{Me}, p. 95, combined with ibid., p. 69 and p. 72.}, and these adjunctions are functorial in the morphism (cf. \cite{Ho-Ta-Tan}, 2.7.2, where the notation is $Tr$ instead of $Ad$). This functoriality and $\pi \circ \epsilon=\id$ then imply that $(1.2.8)$ is indeed the identity.\\
\newline
For general $S$ the integrality\footnote{Our general assumptions imply that $S$ is a regular scheme and hence its local rings are integral domains. As it is moreover connected and (locally) noetherian we can readily conclude its integrality (cf. also \cite{Gö-We}, Ch. 3, Ex. 3.15 and Ex. 3.16 (b)).} of $S$ and a standard compatibility of the occurring arrows with base change (for the case of the Grothendieck trace map cf. \cite{Con1}, Ch. I, 1.1) reduce the commutativity of $(1.2.6)$ without problems to the situation of $S=\Spec(k)$. \qquad \qed

\section{Unipotent vector bundles with integrable connection}
\markright{\uppercase{The formalism of the logarithm sheaves and the elliptic...}}

\subsection{The notion of unipotency}
In 1.1 we saw that there exists a filtration
\[\mathcal L_n=A^0\mathcal L_n \supseteq A^1\mathcal L_n \supseteq...\supseteq A^n\mathcal L_n \supseteq A^{n+1}\mathcal L_n=0\]
of $\mathcal L_n$ by subvector bundles stable under $\nabla_n$ and with quotients given by
\[A^i\mathcal L_n/A^{i+1}\mathcal L_n \simeq \pi^*\mathrm{Sym}^i_{\mathcal O_S}\mathcal H, \quad i=0,...,n.\]
We now take this observation as a model to give a general definition of unipotency for vector bundles with integrable connection, adapted to our fixed geometric setting
\begin{equation*}
\begin{xy}
\xymatrix@C-0.3cm{
X \ar[rr]^{\pi} \ar[dr]& & S\ar[dl]\\
& \Spec(\Q) &}
\end{xy}
\end{equation*}
We further record the behaviour of unipotent bundles under some basic operations.\\
\newline
Let us denote by $\textit{VIC}(X/\Q)$ the category whose objects are the vector bundles on $X$ with integrable $\Q$-connection and whose morphisms are the $\mathcal O_X$-module homomorphisms respecting the connections. The fact that a coherent $\mathcal O_X$-module with integrable $\Q$-connection is already a vector bundle (cf. \cite{Bert-Og}, §2, Note 2.17) readily implies that $\textit{VIC}(X/\Q)$ is abelian.\\
Replacing $X$ by $S$ in what we just said defines the abelian category $\textit{VIC}(S/\Q)$.\\
By pullback via $X \xrightarrow{\pi} S$ we obtain an exact functor $\pi^*: \textit{VIC}(S/\Q) \rightarrow \textit{VIC}(X/\Q)$.

\begin{definition}
Let $n\geq 0$.\\
(i) An object $\mathcal U$ of $\textit{VIC}(X/\Q)$ is called \underline{unipotent of length $n$ for $X/S/\Q$} if there exists a filtration
\[\mathcal U=A^0\mathcal U \supseteq A^1\mathcal U \supseteq...\supseteq A^n\mathcal U \supseteq A^{n+1}\mathcal U=0\]
by subvector bundles stable under the connection of $\mathcal U$ such that for all $i=0,...,n$ there are objects $\mathcal Y_i$ of $\textit{VIC}(S/\Q)$ and isomorphisms in $\textit{VIC}(X/\Q)$:
\[A^i\mathcal U/A^{i+1}\mathcal U \simeq \pi^*\mathcal Y_i.\]
(ii) We write $U_n(X/S/\Q)$ for the full subcategory of $\textit{VIC}(X/\Q)$ consisting of those $\mathcal U$ in $\textit{VIC}(X/\Q)$ which are unipotent of length $n$ for $X/S/\Q$.\\
(iii) We write $U(X/S/\Q)$ for the full subcategory of $\textit{VIC}(X/\Q)$ consisting of those $\mathcal U$ in $\textit{VIC}(X/\Q)$ which are unipotent of some length for $X/S/\Q$. In other words, $U(X/S/\Q)$ is the union of the $U_n(X/S/\Q)$ for the canonical embeddings
\[U_0(X/S/\Q) \hookrightarrow U_1(X/S/\Q) \hookrightarrow U_2(X/S/\Q) \hookrightarrow... \hookrightarrow \textit{VIC}(X/\Q).\] 
Note that the zero vector bundle on $X$ with its unique $\Q$-connection is an object of each $U_n(X/S/\Q)$.\\
Note further that $U_0(X/S/\Q)$ is just the essential image of the functor $\pi^*: \textit{VIC}(S/\Q) \rightarrow \textit{VIC}(X/\Q)$.
\end{definition}
In particular, $\mathcal L_n$ with its integrable $\Q$-connection $\nabla_n$ becomes an object of $U_n(X/S/\Q)$.\\
\newline
The following lemma gives some first easy properties of unipotent bundles with integrable connection; for the dual and tensor product of modules with integrable connection cf. 0.2.1 (iv). 
\begin{lemma}
(i) If
\[0 \rightarrow \mathcal U' \rightarrow \mathcal U \rightarrow \mathcal U'' \rightarrow 0\]
is an exact sequence in $\textit{VIC}(X/\Q)$ with $\mathcal U'$ in $U_m(X/S/\Q)$ and $\mathcal U''$ in $U_n(X/S/\Q)$, then $\mathcal U$ is in $U_{m+n+1}(X/S/\Q)$.\\
\newline
(ii) If $\mathcal U$ is in $U_n(X/S/\Q)$, then its dual $\mathcal U^\vee$ is also in $U_n(X/S/\Q)$.\\
\newline
(iii) If $\mathcal V$ is in $U_m(X/S/\Q)$ and $\mathcal W$ is in $U_n(X/S/\Q)$, then the tensor product $\mathcal V\otimes_{\mathcal O_X} \mathcal W$ is in $U_{m+n}(X/S/\Q)$.
\end{lemma}
\begin{proof}
(i) is straightforward and (ii) follows easily from (i) by induction on the length of $\mathcal U$.\\
For (iii) define for each $k=0,...,m+n+1$:
\[A^k(\mathcal V \otimes_{\mathcal O_X}\mathcal W):=\sum_{i+j=k}(A^i\mathcal V \otimes_{\mathcal O_X} A^j \mathcal W) \ \ \textrm{in} \ \ \mathcal V \otimes_{\mathcal O_X} \mathcal W,\]
i.e.
\[A^k(\mathcal V \otimes_{\mathcal O_X}\mathcal W)=\mathrm{im} \Big(\bigoplus_{i+j=k}(A^i\mathcal V \otimes_{\mathcal O_X} A^j \mathcal W) \xrightarrow{\mathrm{can}} \mathcal V \otimes_{\mathcal O_X} \mathcal W \Big),\]
viewed as submodule of $\mathcal V \otimes_{\mathcal O_X}\mathcal W$.\\
We may endow these coherent subsheaves with the induced integrable $\Q$-connection such that they are vector bundles. There clearly is a chain of inclusions
\[\mathcal V \otimes_{\mathcal O_X} \mathcal W=A^0(\mathcal V \otimes_{\mathcal O_X} \mathcal W) \supseteq A^1(\mathcal V \otimes_{\mathcal O_X} \mathcal W) \supseteq...\supseteq A^{n+m}(\mathcal V \otimes_{\mathcal O_X} \mathcal W) \supseteq A^{m+n+1}(\mathcal V \otimes_{\mathcal O_X} \mathcal W)=0,\]
and it is not hard to check that the quotients are given for each $k=0,...,m+n$ by
\[A^k(\mathcal V \otimes_{\mathcal O_X} \mathcal W)/A^{k+1}(\mathcal V \otimes_{\mathcal O_X} \mathcal W) \simeq \bigoplus_{i+j=k}\Big((\mathcal A^i\mathcal V/A^{i+1}\mathcal V) \otimes_{\mathcal O_X} (A^j\mathcal W/A^{j+1}\mathcal W)\Big)\]
and hence obviously of the desired form.
\end{proof}
\subsection{The universal property of the logarithm sheaves}
The main goal is to show that the logarithm sheaves $(\mathcal L_n,\nabla_n) \in U_n(X/S/\Q)$ together with a distinguished section of their zero fiber are characterized by a universal property. This is well-known in other realizations (cf. \cite{Be-Le}, Prop. 1.2.6, or \cite{Hu-Ki}, Lemma A.2.3), from which we may extract the formal structure of the arguments; adjusting them properly for our case of de Rham realization and supplying all necessary details will be the task in what follows. The essential ingredient for the proof then consists in our knowledge of the de Rham cohomology of the logarithm sheaves (cf. 1.2.1).

\subsubsection{An auxiliary lemma}
To show the main theorem we will need the following duality result which is of independent interest.

\begin{lemma}
Let $\mathcal V$ be in $U_n(X/S/\Q)$ and $\mathcal Z$ in $\textit{VIC}(S/\Q)$. Then for each $0\leq i \leq 2g$ we have a canonical isomorphism in $\textit{VIC}(S/\Q)$
\[H^{2g-i}_{{\mathrm{dR}}}(X/S,\underline{\Hom}_{\mathcal O_X}(\mathcal V,\pi^*\mathcal Z)) \simeq \underline{\Hom}_{\mathcal O_S}(H^i_{\mathrm{dR}}(X/S,\mathcal V),\mathcal Z),\]
functorial in $\mathcal V$ and in $\mathcal Z$.
\end{lemma}
\begin{remark}
We explain how to consider both sides in the preceding lemma as objects of $\textit{VIC}(S/\Q)$:\\
On the left side observe that we endow the $\mathcal O_X$-vector bundle $\underline{\Hom}_{\mathcal O_X}(\mathcal V,\pi^*\mathcal Z)$ with the natural integrable $\Q$-connection (cf. 0.2.1 (iv)) which becomes the tensor product connection when identifying the bundle with $\mathcal V^\vee \otimes_{\mathcal O_X} \pi^*\mathcal Z$. The de Rham cohomology on the left side of the lemma is then equipped with the Gauß-Manin connection relative $\Spec(\Q)$ and indeed is a vector bundle on $S$ (by the same argument as in Rem. 1.2.2).\\
On the right side $H^i_{\mathrm{dR}}(X/S,\mathcal V)$ is a vector bundle on $S$ (cf. again Rem. 1.2.2) and equipped with the Gauß-Manin connection relative $\Spec(\Q)$. The internal $\Hom$ of the right side then clearly is a vector bundle on $S$ and carries its usual integrable $\Q$-connection.
\end{remark}
We now prove Lemma 1.3.3.
\begin{proof}
The left side canonically identifies with
\[H^{2g-i}_{{\mathrm{dR}}}(X/S,\underline{\Hom}_{\mathcal O_X}(\mathcal V,\pi^*\mathcal Z)) \simeq H^{2g-i}_{{\mathrm{dR}}}(X/S, \mathcal V^\vee \otimes_{\mathcal O_X} \pi^*\mathcal Z) \simeq H^{2g-i}_{{\mathrm{dR}}}(X/S,\mathcal V^\vee)\otimes_{\mathcal O_S}\mathcal Z\]
in $\textit{VIC}(S/\Q)$, whereas for the right side we have
\[\underline{\Hom}_{\mathcal O_S}(H^i_{\mathrm{dR}}(X/S,\mathcal V),\mathcal Z) \simeq H^i_{\mathrm{dR}}(X/S,\mathcal V)^\vee \otimes_{\mathcal O_S}\mathcal Z\]
in $\textit{VIC}(S/\Q)$, and hence it suffices to show that for each $i$ there is a functorial horizontal isomorphism
\[H^i_{\mathrm{dR}}(X/S,\mathcal V) \simeq H^{2g-i}_{{\mathrm{dR}}}(X/S,\mathcal V^\vee)^\vee.\]
For this consider the composition
\[H^i_{\mathrm{dR}}(X/S,\mathcal V) \otimes_{\mathcal O_S} H^{2g-i}_{{\mathrm{dR}}}(X/S,\mathcal V^\vee) \rightarrow H^{2g}_{\mathrm{dR}}(X/S,\mathcal V \otimes_{\mathcal O_X}\mathcal V^\vee) \rightarrow H^{2g}_{\mathrm{dR}}(X/S) \xrightarrow{\sim} \mathcal O_S,\]
where the first arrow is given by cup product, the second is induced by the canonical map $\mathcal V \otimes_{\mathcal O_X} \mathcal V^\vee \rightarrow \mathcal O_X$ and the last is the trace isomorphism. All three arrows are horizontal: for the second it is clear and for the two other maps cf. footnote 1 (the argument used there for the horizontality of the cup product carries over to the present situation of de Rham cohomology with coefficients).\\
We thus obtain an induced (and in fact functorial) map in $\textit{VIC}(S/\Q)$:
\[\tag{\textbf{1.3.1}} H^i_{\mathrm{dR}}(X/S,\mathcal V) \rightarrow H^{2g-i}_{{\mathrm{dR}}}(X/S,\mathcal V^\vee)^\vee,\]
and by what we already said it only remains to show that $(1.3.1)$ is an isomorphism.\\
Indeed, this isomorphism should be valid for arbitrary $\mathcal V$ in $\textit{VIC}(X/\Q)$, but we only need it for unipotent $\mathcal V$ where the arguments are easier: we proceed by induction over the length $n$ of $\mathcal V$.\\
For $n=0$ one writes $\mathcal V\simeq \pi^*\mathcal Y$ for some $\mathcal Y$ in $\textit{VIC}(S/\Q)$ and identifies the left side of $(1.3.1)$ with
\[H^i_{\mathrm{dR}}(X/S,\pi^*\mathcal Y) \simeq H^i_{\mathrm{dR}}(X/S)\otimes_{\mathcal O_S}\mathcal Y,\]
the right side with
\[H^{2g-i}_{{\mathrm{dR}}}(X/S,\pi^*\mathcal Y^\vee)^\vee \simeq H^{2g-i}_{{\mathrm{dR}}}(X/S)^\vee \otimes_{\mathcal O_S}\mathcal Y,\]
and uses the canonical isomorphism $H^i_{\mathrm{dR}}(X/S) \simeq H^{2g-i}_{{\mathrm{dR}}}(X/S)^\vee$ (given precisely by cup-product and the trace isomorphism) in order to deduce the claim for $n=0$.\\
If $n\geq 1$ we perform the induction step by considering an exact sequence in $\textit{VIC}(X/\Q)$
\[\tag{\textbf{1.3.2}} 0\rightarrow A^1\mathcal V \rightarrow \mathcal V \rightarrow \pi^* \mathcal X \rightarrow 0\]
in which $A^1\mathcal V$ is unipotent of length $n-1$ and $\pi^*\mathcal X$ is unipotent of length $0$ (where $\mathcal X$ is of course a suitable object of $\textit{VIC}(S/\Q)).$\\
There is a diagram of exact sequences of finite length in $\textit{VIC}(S/\Q)$:
\begin{equation*}
{\footnotesize
\begin{xy}
\xymatrix@C-0.3cm{
0 \ar[r] \ar[d]^{\sim} & H^0_{\mathrm{dR}}(X/S,A^1 \mathcal V) \ar[d]^{\sim} \ar[r]  & H^0_{\mathrm{dR}}(X/S, \mathcal V) \ar[d] \ar[r] & H^0_{\mathrm{dR}}(X/S,\pi^*\mathcal X) \ar[r] \ar[d]^{\sim} & H^1_{\mathrm{dR}}(X/S,A^1 \mathcal V) \ar[d]^{\sim} \ar[r] & ... \\
0 \ar[r] & H^{2g}_{\mathrm{dR}}(X/S,(A^1 \mathcal V)^\vee)^\vee \ar[r] & H^{2g}_{\mathrm{dR}}(X/S,\mathcal V^\vee)^\vee \ar[r] & H^{2g}_{\mathrm{dR}}(X/S,\pi^*\mathcal X^\vee)^\vee \ar[r] & H^{2g-1}_{\mathrm{dR}}(X/S,(A^1 \mathcal V)^\vee)^\vee \ar[r] & ...}
\end{xy}
}
\end{equation*}
The upper row is given by the long exact sequence of de Rham cohomology (cf. \cite{Kat2}, Rem. $(3.1)$) for $(1.3.2)$. The lower row comes from dualizing $(1.3.2)$, applying to this dual sequence the long exact sequence of de Rham cohomology and finally dualizing the obtained long exact sequence in $\textit{VIC}(S/\Q)$. The vertical arrows are the maps defined in $(1.3.1)$ for the various bundles, where for the indicated isomorphisms we have already taken into account the induction hypothesis.\\
The squares in which no connecting homomorphism is involved are commutative by the (easily seen) functoriality of $(1.3.1)$.
The squares in which the upper row is a connecting morphism
\[H^j_{\mathrm{dR}}(X/S,\pi^*\mathcal X) \rightarrow H^{j+1}_{{\mathrm{dR}}}(X/S,A^1 \mathcal V)\]
commute up to a sign as one can indeed verify.\footnote{This rests on a straightforwardly checked naturality property of the cup product: If
\[0 \rightarrow \mathcal F \rightarrow \mathcal G \rightarrow \mathcal H \rightarrow 0\]
is a short exact sequence in $\textit{VIC}(X/\Q)$, then for each $j$ the diagram
\begin{equation*}
\begin{xy}
\xymatrix@C-0.3cm{
H^j_{\mathrm{dR}}(X/S,\mathcal H) \ar[r] \ar[d] & H^{j+1}_{{\mathrm{dR}}}(X/S,\mathcal F)\ar[d]\\
\underline{\Hom}_{\mathcal O_S}(H^{2g-j}_{{\mathrm{dR}}}(X/S,\mathcal H^\vee),H^{2g}_{\mathrm{dR}}(X/S)) \ar[r] & \underline{\Hom}_{\mathcal O_S}(H^{2g-j-1}_{{\mathrm{dR}}}(X/S,\mathcal F^\vee),H^{2g}_{\mathrm{dR}}(X/S))}
\end{xy}
\end{equation*}
is commutative up to a sign. Here, the upper horizontal arrow is the connecting morphism on level $j$ for our sequence, and the lower horizontal map comes from applying $\underline{\Hom}_{\mathcal O_S}(-,H^{2g}_{\mathrm{dR}}(X/S))$ to the connecting morphism on level $2g-j-1$ of the dualized sequence. The vertical arrows are defined analogously to $(1.3.1)$, but without identifying $H^{2g}_{\mathrm{dR}}(X/S)\simeq \mathcal O_S$.
}\\
This obviously yields the claim by pursuing the big diagram successively until its end.
\end{proof}
The comments of the following remark will frequently (and often tacitly) be used in what follows.
\begin{remark}
(i) Let $\mathcal V$ and $\mathcal W$ be objects in $\textit{VIC}(X/\Q)$.\\
We then have a canonical identification of $\mathcal O_S$-modules
\[\tag{\textbf{1.3.3}} H^0_{\mathrm{dR}}(X/S,\underline{\Hom}_{\mathcal O_X}(\mathcal V, \mathcal W)) \simeq \pi_*\underline{\Hom}_{\mathcal D_{X/S}}(\mathcal V, \mathcal W)\]
which comes about by noting that
\[H^0_{\mathrm{dR}}(X/S,\underline{\Hom}_{\mathcal O_X}(\mathcal V, \mathcal W))=\pi_*\Big(\underline{\Hom}_{\mathcal O_X}(\mathcal V, \mathcal W)^{\nabla_{X/S}}\Big)=\pi_*\underline{\Hom}_{\mathcal D_{X/S}}(\mathcal V, \mathcal W).\]
The superscript in the middle term means that we take the subsheaf of $\underline{\Hom}_{\mathcal O_X}(\mathcal V, \mathcal W)$ consisting of those sections which are horizontal for the connection restricted relative $S$ (cf. \cite{Kat2}, Rem. $(3.1)$).\\
Observe that via the Gauß-Manin connection relative $\Spec(\Q)$ on the left side of $(1.3.3)$ we may and will equip $\pi_*\underline{\Hom}_{\mathcal D_{X/S}}(\mathcal V, \mathcal W)$ with the induced connection and thus view it as object of $\textit{VIC}(S/\Q)$.\vspace{1mm}\\
(ii) In particular, in the situation (and with the result of) Lemma 1.3.3 we obtain isomorphisms in $\textit{VIC}(S/\Q)$:
\[\tag{\textbf{1.3.4}} \pi_*\underline{\Hom}_{\mathcal D_{X/S}}(\mathcal V, \pi^*\mathcal Z) \simeq H^0_{\mathrm{dR}}(X/S,\underline{\Hom}_{\mathcal O_X}(\mathcal V, \pi^* \mathcal Z)) \simeq \underline{\Hom}_{\mathcal O_S}(H^{2g}_{\mathrm{dR}}(X/S,\mathcal V),\mathcal Z).\]
On sections this composition is given as follows: if $\mathcal V \rightarrow \pi^*\mathcal Z$ is a morphism which is horizontal for the connections restricted relative $S$, one applies $H^{2g}_{\mathrm{dR}}(X/S,-)$ to it and uses $H^{2g}_{\mathrm{dR}}(X/S,\pi^*\mathcal Z) \simeq H^{2g}_{\mathrm{dR}}(X/S) \otimes_{\mathcal O_S} \mathcal Z \simeq \mathcal Z$, where the last isomorphism is given by the trace map.
\end{remark}
\subsubsection{The universal property}
We now show how morphisms of $\mathcal L_n$ into another unipotent bundle $\mathcal U \in U_n(X/S/\Q)$ are parametrized via the fiber $\epsilon^*\mathcal U$, which is a direct expression of the universal property as we will then explain.\\
The reader acquainted with other realizations (cf. \cite{Be-Le}, Prop. 1.2.6, or \cite{Hu-Ki}, Lemma A.2.3) should - after a correct translation of the formalism - expect an isomorphism of the form
\[H^0_{\mathrm{dR}}(X/S, \underline{\Hom}_{\mathcal O_X}(\mathcal L_n, \mathcal U)) \simeq \epsilon^*\mathcal U\]
resp., by using the identification of $(1.3.3)$, an isomorphism
\[\pi_*\underline{\Hom}_{\mathcal D_{X/S}}(\mathcal L_n, \mathcal U) \simeq \epsilon^*\mathcal U.\]
This is indeed the case, and the precise statement in the case of de Rham realization is as follows:
\begin{theorem}
Let $n \geq 0$ and $\mathcal U$ be an object of $U_n(X/S/\Q)$.\\
\newline
(i) For each $k\geq n$ we then have an isomorphism in $\textit{VIC}(S/\Q)$
\[\tag{\textbf{1.3.5}} \pi_*\underline{\Hom}_{\mathcal D_{X/S}}(\mathcal L_k, \mathcal U) \xrightarrow{\sim} \epsilon^*\mathcal U,\]
functorial in $\mathcal U$ and compatible with the projections $\mathcal L_k \rightarrow \mathcal L_l$ of the logarithm sheaves for $k \geq l \geq n$.\\
It is defined on sections by
\[f \mapsto \epsilon^*(f) \Big(\frac{1}{k!}\Big),\]
where we use the identification
\[\varphi_k: \prod_{i=0}^k \mathrm{Sym}^i_{\mathcal O_S}\mathcal H \simeq \epsilon^*\mathcal L_k.\]
\\
(ii) A section $f: \mathcal L_k \rightarrow \mathcal U$ of the left side of $(1.3.5)$ is $\mathcal D_{X/\Q}$-linear (and not only $\mathcal D_{X/S}$-linear) if and only if its image $\epsilon^*(f)(\frac{1}{k!})$ under $(1.3.5)$ is a horizontal section of $\epsilon^*\mathcal U$.
\end{theorem}
\begin{proof}
(i) The formal frame of the proof is as in \cite{Hu-Ki}, Lemma A.2.3.\\
We proceed by induction over the length $n$ of $\mathcal U$.\\
For $n=0$ we have $\mathcal U \simeq \pi^*\mathcal Z$ with an object $\mathcal Z$ of $\textit{VIC}(S/\Q)$. By $(1.3.4)$ we then have a chain of isomorphisms in $\textit{VIC}(S/\Q)$:
\[\tag{\textbf{1.3.6}} \pi_*\underline{\Hom}_{\mathcal D_{X/S}}(\mathcal L_k, \pi^*\mathcal Z) \simeq H^0_{\mathrm{dR}}(X/S,\underline{\Hom}_{\mathcal O_X}(\mathcal L_k, \pi^* \mathcal Z)) \simeq \underline{\Hom}_{\mathcal O_S}(H^{2g}_{\mathrm{dR}}(X/S,\mathcal L_k),\mathcal Z).\]
Because of the identification $H^{2g}_{\mathrm{dR}}(X/S,\mathcal L_k) \xrightarrow{\sim} \mathcal O_S$, induced by the projection $\mathcal L_k \rightarrow \mathcal O_X$ and the trace map (cf. Thm. 1.2.1 (i)), we get from $(1.3.6)$ an isomorphism in $\textit{VIC}(S/\Q)$ of the desired form
\[\tag{\textbf{1.3.7}} \pi_*\underline{\Hom}_{\mathcal D_{X/S}}(\mathcal L_k, \pi^*\mathcal Z) \simeq \mathcal Z \simeq \epsilon^*\pi^*\mathcal Z.\]
We claim that for a section $f: \mathcal L_k \rightarrow \pi^*\mathcal Z$ of the left side it really acts as $f\mapsto \epsilon^*(f) (\frac{1}{k!})$. To see this we evoke the following diagram of $\mathcal O_S$-vector bundles:
\begin{equation*}
\begin{xy}
\xymatrix@C-0.3cm{
\prod_{i=0}^k \mathrm{Sym}^i_{\mathcal O_S}\mathcal H \ar@{-}[r]^{\qquad \sim} & \epsilon^*\mathcal L_k \ar[d]_{\epsilon^*(f)} \ar[r]^{\sigma_k \qquad }  & H^{2g}_{\mathrm{dR}}(X/S,\mathcal L_k) \ar[d] \ar[r]^{\sim \quad} & H^{2g}_{\mathrm{dR}}(X/S) \simeq \mathcal O_S \\
& \mathcal Z \simeq \epsilon^*\pi^*\mathcal Z \ar[r] & H^{2g}_{\mathrm{dR}}(X/S,\pi^*Z) \ar[r]^{\sim \qquad} & H^{2g}_{\mathrm{dR}}(X/S) \otimes_{\mathcal O_S}\mathcal Z \simeq \mathcal Z}
\end{xy}
\end{equation*}
Here, the vertical arrows of the square are the naturally induced ones and its horizontal arrows come about as follows: if $\mathcal G$ is $\mathcal L_k$ or $\pi^*\mathcal Z$, then one restricts its $\Q$-connection relative $S$, applies the functor $\pi_+: D^b_{\mathrm{qc}}(\mathcal D_{X/S}) \rightarrow D^b_{\mathrm{qc}}(\mathcal O_S)$ to the map $\epsilon_+\epsilon^*\mathcal G[-g]\rightarrow \mathcal G$ in its localization triangle relative $S$,  takes $g$-th cohomology and canonically identifies $H^g(\pi_+\mathcal G) \simeq H^{2g}_{\mathrm{dR}}(X/S,\mathcal G)$ as in $(1.2.7)$. The square commutes by functoriality of the localization triangle. The other maps in the diagram are clear.\\
By Rem. 1.3.5 (ii) the image of $f$ under $(1.3.6)$ is the composition of the right vertical arrow of the square with the two lower right arrows of the diagram. In order to get the desired image of $f$ under $(1.3.7)$ we need to precompose this last composition with the inverse of the upper right isomorphism. But the whole upper row of the diagram maps the section $1$ of $\mathcal O_S \subseteq \prod_{i=0}^k \mathrm{Sym}^i_{\mathcal O_S}\mathcal H$ to $k!$,\footnote{The proof of this is as for Lemma 1.2.5 (cf. 1.2.3), only that one works with $\Q$ replaced by $S$ from the beginning on.} and the whole lower row is the identity.\footnote{One can derive this pretty straightforwardly from the preceding assertion applied for $k=0$.}\\
With this the claim follows directly from the diagram, and we hence conclude the case $n=0$.\\
\newline
If $n\geq 1$ we find an exact sequence in $\textit{VIC}(X/\Q)$
\[0\rightarrow \pi^*\mathcal Y \rightarrow \mathcal U \rightarrow \mathcal U/\pi^*\mathcal Y \rightarrow 0\]
with $\mathcal U/\pi^*\mathcal Y \in U_{n-1}(X/S/\Q)$. We then obtain the exact sequence
\[\tag{\textbf{1.3.8}} 0 \rightarrow \underline{\Hom}_{\mathcal O_X}(\mathcal L_k, \pi^*\mathcal Y) \rightarrow \underline{\Hom}_{\mathcal O_X}(\mathcal L_k, \mathcal U) \rightarrow \underline{\Hom}_{\mathcal O_X}(\mathcal L_k, \mathcal U/\pi^*\mathcal Y) \rightarrow 0\]
in $\textit{VIC}(X/\Q)$ and the commutative diagram of exact sequences of $\mathcal O_S$-vector bundles
\begin{equation*}
{\footnotesize
\begin{xy}
\xymatrix@C-0.3cm{
0 \ar[r] & \pi_*\underline{\Hom}_{\mathcal D_{X/S}}(\mathcal L_k, \pi^*\mathcal Y) \ar[d]^{\sim} \ar[r]  & \pi_*\underline{\Hom}_{\mathcal D_{X/S}}(\mathcal L_k, \mathcal U) \ar[d] \ar[r] & \pi_*\underline{\Hom}_{\mathcal D_{X/S}}(\mathcal L_k, \mathcal U/\pi^*\mathcal Y) \ar[r]^{\delta_k\quad} \ar[d]^{\sim} & H^1_{\mathrm{dR}}(X/S,\underline{\Hom}_{\mathcal O_X}(\mathcal L_k, \pi^*\mathcal Y)) \ar[d]\\
0 \ar[r] & \mathcal Y \ar[r] & \epsilon^*\mathcal U \ar[r] & \epsilon^*(\mathcal U/\pi^*\mathcal Y) \ar[r] & 0}
\end{xy}
}
\end{equation*}
where we define the first three vertical arrows on sections by the rule $f \mapsto \epsilon^*(f) (\frac{1}{k!})$ and where we have already used the induction hypothesis. The upper row comes from applying the long exact sequence of de Rham cohomology to $(1.3.8)$ and usage of $(1.3.3)$.\\
To show that the second vertical arrow is an isomorphism it suffices to see that $\delta_k$ is the zero map.\\
For this we use the commutative diagram
\begin{equation*}
\begin{xy}
\xymatrix{
\pi_*\underline{\Hom}_{\mathcal D_{X/S}}(\mathcal L_{k-1}, \mathcal U/\pi^*\mathcal Y))  \ar[r]^{\quad \delta_{k-1} \qquad } \ar[d]^{\mathrm{can}} & H^1_{\mathrm{dR}}(X/S,\underline{\Hom}_{\mathcal O_X}(\mathcal L_{k-1}, \pi^*\mathcal Y)) \ar[d]^{\mathrm{can}}\\
\pi_*\underline{\Hom}_{\mathcal D_{X/S}}(\mathcal L_k, \mathcal U/\pi^*\mathcal Y)\ar[r]^{\delta_k \quad} & H^1_{\mathrm{dR}}(X/S,\underline{\Hom}_{\mathcal O_X}(\mathcal L_k, \pi^*\mathcal Y))} 
\end{xy}
\end{equation*}
and show that the right vertical arrow is zero. This yields the claim $\delta_k=0$ because the left vertical arrow becomes the identity if we identify both terms with $\epsilon^*(\mathcal U/\pi^*\mathcal Y)$ via the induction hypothesis (use Lemma 1.1.5). But according to Lemma 1.3.3 the right vertical map identifies with an arrow
\[\underline{\Hom}_{\mathcal O_S}(H^{2g-1}_{\mathrm{dR}}(X/S,\mathcal L_{k-1}), \mathcal Y) \rightarrow \underline{\Hom}_{\mathcal O_S}(H^{2g-1}_{\mathrm{dR}}(X/S,\mathcal L_k), \mathcal Y)\]
which is induced by the transition map $\mathcal L_k \rightarrow \mathcal L_{k-1}$. Thm. 1.2.1 (ii) then implies that this arrow is indeed the zero morphism.\\
\newline
We now know that the map of $(1.3.5)$:
\[\pi_*\underline{\Hom}_{\mathcal D_{X/S}}(\mathcal L_k, \mathcal U) \rightarrow \epsilon^*\mathcal U, \quad f \mapsto \epsilon^*(f)\Big(\frac{1}{k!}\Big)\]
is an isomorphism of $\mathcal O_S$-vector bundles. The left side carries an integrable $\Q$-connection via
\[H^0_{\mathrm{dR}}(X/S,\underline{\Hom}_{\mathcal O_X}(\mathcal L_k, \mathcal U)) \simeq \pi_*\underline{\Hom}_{\mathcal D_{X/S}}(\mathcal L_k, \mathcal U)\]
(cf. Rem. 1.3.5 (i)), and with the easy explicit knowledge of the Gauß-Manin connection in $0$-th cohomology (cf. \cite{Kat2}, Rem. $(3.1)$) one checks by hands that $(1.3.5)$ is indeed horizontal: when doing this the essential point is the horizontality of $\frac{1}{k!}$ for the pullback of the connection of $\mathcal L_k$ via $\epsilon$.\\
That $(1.3.5)$ is functorial in $\mathcal U$ is clear, and that it is compatible with the transition maps of the logarithm sheaves follows from Lemma 1.1.5. This finishes the proof of (i).\\
\newline
(ii) Consider $\pi_*\underline{\Hom}_{\mathcal D_{X/S}}(\mathcal L_k, \mathcal U)$ with its integrable $\Q$-connection induced by the identification
\[H^0_{\mathrm{dR}}(X/S, \underline{\Hom}_{\mathcal O_X}(\mathcal L_n, \mathcal U)) \simeq \pi_*\underline{\Hom}_{\mathcal D_{X/S}}(\mathcal L_k, \mathcal U)\]
of $(1.3.3)$. Let us compute its global horizontal $S$-sections: they are given by
\[H^0_{\mathrm{dR}}(S/\Q, H^0_{\mathrm{dR}}(X/S, \underline{\Hom}_{\mathcal O_X}(\mathcal L_n, \mathcal U))) = H^0_{\mathrm{dR}}(X/\Q,\underline{\Hom}_{\mathcal O_X}(\mathcal L_k,\mathcal U))=\Hom_{\mathcal D_{X/\Q}}(\mathcal L_k,\mathcal U),\]
and the analogous calculation holds on open subsets of $S$, showing that the subsheaf of horizontal sections of $\pi_*\underline{\Hom}_{\mathcal D_{X/S}}(\mathcal L_k, \mathcal U)$ is given by $\pi_*\underline{\Hom}_{\mathcal D_{X/\Q}}(\mathcal L_k, \mathcal U)$. As the isomorphism $(1.3.5)$ is horizontal by (i) the claim of (ii) follows.
\end{proof}

\begin{definition}
For each $n\geq 0$ we denote by $1^{(n)}$ the image of $\frac{1}{n!}$ under the (horizontal) splitting
\[\varphi_n: \prod_{k=0}^n \mathrm{Sym}^k_{\mathcal O_S}\mathcal H \simeq \epsilon^*\mathcal L_n,\]
such that $1^{(n)}$ is a global horizontal $S$-section of $\epsilon^*\mathcal L_n$.
\end{definition}
By Thm. 1.3.6 (i) we may view $(\mathcal L_n, 1^{(n)})$ as a pair consisting of an object $\mathcal L_n$ in $U_n(X/S/\Q)$ and a global horizontal $S$-section $1^{(n)}$ of $\epsilon^*\mathcal L_n$ with the property that for any $\mathcal U$ in $U_n(X/S/\Q)$ the map
\[\pi_*\underline{\Hom}_{\mathcal D_{X/S}}(\mathcal L_n, \mathcal U) \rightarrow \epsilon^*\mathcal U, \quad f \mapsto \epsilon^*(f)\big(1^{(n)}\big)\]
is a horizontal isomorphism. Using an analogous claim as in Thm. 1.3.6 (ii) it is routinely checked that a pair consisting of an object in $U_n(X/S/\Q)$ and a global horizontal $S$-section of the fiber in $\epsilon$ with the preceding property is unique up to unique isomorphism; here, isomorphism means an isomorphism in $\textit{VIC}(X/\Q)$ which respects the distinguished sections after pullback via $\epsilon$.\\
This is the manifestation of the \underline{universal property of the $n$-th logarithm sheaf.}
\begin{remark}
According to Lemma 1.1.5 the transition map $\mathcal L_{n+1} \rightarrow \mathcal L_n$ is then given by $1^{(n+1)} \mapsto 1^{(n)}$.
\end{remark}

\subsection{An equivalence of categories}
As a further fundamental result we show that passage to the zero fiber identifies $U_n(X/S/\Q)$ with a certain category $\mathcal C_n$ of $\mathcal O_S$-vector bundles with integrable $\Q$-connection carrying a (compatible) module structure over the sheaf of rings $\prod_{k=0}^n \mathrm{Sym}^k_{\mathcal O_S} \mathcal H$. This is the analogue in the de Rham realization of \cite{Be-Le}, 1.2.10 (v) (cf. also \cite{Hu-Ki}, Thm. A.2.5). Our goal for what follows is to give the proper definition of the category $\mathcal C_n$, to explicitly construct and give sense to the quasi-inverse of the mentioned fiber functor and to then prove the equivalence result in full detail. This is a somewhat laborious and technical task to do, but the theorem will be rather useful in our further study of unipotent bundles and the logarithm sheaves.

\subsubsection{The category $\boldsymbol{\mathcal C_n}$}
For all $n\geq 0$ consider $\prod_{k=0}^n \mathrm{Sym}^k_{\mathcal O_S} \mathcal H$ as a commutative sheaf of rings on $S$ with multiplication defined on components by
\[\tag{\textbf{1.3.9}} s^k \circ s^l:=\frac{(n-k)!(n-l)!}{(n-k-l)!} s^k \cdot s^l\]
if $s^k$ is a section of $\mathrm{Sym}^k_{\mathcal O_S}\mathcal H$ and $s^l$ is a section of $\mathrm{Sym}^l_{\mathcal O_S}\mathcal H$ such that $k+l \leq n$ resp. by $s^k \circ s^l=0$ in the case $k+l >n$. Note that on the right side of $(1.3.9)$ we mean multiplication in symmetric powers, hence $s^k \circ s^l$ is a section of $\mathrm{Sym}_{\mathcal O_S}^{k+l}\mathcal H$. We extend the multiplication to the whole product by linearity. The multiplicative identity is then given by $\frac{1}{n!}$.\\
We get an $\mathcal O_S$-algebra structure inducing the original $\mathcal O_S$-module structure on $\prod_{k=0}^n \mathrm{Sym}^k_{\mathcal O_S} \mathcal H$ by
\[\tag{\textbf{1.3.10}} \mathcal O_S \rightarrow \prod_{k=0}^n \mathrm{Sym}^k_{\mathcal O_S} \mathcal H, \quad s \mapsto \frac{s}{n!}.\]
\begin{definition}
For each $n\geq 0$ we write $\mathcal R^{(n)}$ to denote the so defined sheaf of commutative $\mathcal O_S$-algebras $\prod_{k=0}^n \mathrm{Sym}^k_{\mathcal O_S} \mathcal H$.
\end{definition}
Using Lemma 1.1.5 one sees that the map
\[\mathcal R^{(n+1)} \rightarrow \mathcal R^{(n)},\]
induced by $\mathcal L_{n+1} \rightarrow \mathcal L_n$ via the splittings $\varphi_{n+1}$ and $\varphi_n$, then becomes a morphism of $\mathcal O_S$-algebras.
\begin{remark}
As for any sheaf of quasi-coherent $\mathcal O_S$-algebras we may consider the spectrum of $\mathcal R^{(n)}$ which gives an affine morphism $\Spec(\mathcal R^{(n)}) \rightarrow S$. In our later geometric interpretation of the logarithm sheaves we will identify this spectrum with the $n$-th infinitesimal thickening $\Yr_n$ of the zero section of $\Yr$, where $\Yr$ is the universal vectorial extension of the dual abelian scheme.
\end{remark}
For each $\mathcal U$ in $U_n(X/S/\Q)$ the $\mathcal O_S$-vector bundle $\epsilon^*\mathcal U$ carries a structure of $\mathcal R^{(n)}$-module by 
\[\tag{\textbf{1.3.11}} \begin{split}
\Big(\prod_{k=0}^n \mathrm{Sym}^k_{\mathcal O_S} \mathcal H\Big) \otimes_{\mathcal O_S} \epsilon^*\mathcal U &\simeq \pi_*\underline{\Hom}_{\mathcal D_{X/S}}(\mathcal L_n,\mathcal L_n) \otimes_{\mathcal O_S} \pi_*\underline{\Hom}_{\mathcal D_{X/S}}(\mathcal L_n,\mathcal U) \\
&\xrightarrow{\mathrm{can}} \pi_*\underline{\Hom}_{\mathcal D_{X/S}}(\mathcal L_n,\mathcal U) \simeq \epsilon^*\mathcal U,
\end{split}\]
where the isomorphisms are due to Thm. 1.3.6 (i).\\
In other words: if $r$ is a section of $\mathcal R^{(n)}$ and $\xi$ is a section of $\epsilon^*\mathcal U$, then the multiplication is given by
\[r \cdot \xi=\epsilon^*(f)(r),\]
where $f: \mathcal L_n \rightarrow \mathcal U$ is the unique $\mathcal D_{X/S}$-linear arrow with
\[\epsilon^*(f)\Big(\frac{1}{n!}\Big)=\xi.\]
For $\mathcal U=\mathcal L_n$ the rule $(1.3.11)$ yields precisely the previously defined multiplication of $\mathcal R^{(n)}$.\\
\newline
The so defined $\mathcal R^{(n)}$-module structure on $\epsilon^*\mathcal U$ satisfies:\\
(i) It is compatible with the $\mathcal O_S$-module structure of $\epsilon^*\mathcal U$, i.e. the restriction of the multiplication via the arrow $\mathcal O_S \rightarrow \mathcal R^{(n)}$ of $(1.3.10)$ gives the original $\mathcal O_S$-multiplication on $\epsilon^*\mathcal U$.\\
(ii) The map $\mathcal R^{(n)} \otimes_{\mathcal O_S} \epsilon^*\mathcal U \rightarrow \epsilon^* \mathcal U$ is horizontal for the respective integrable $\Q$-connections.
\begin{definition}
We let $\mathcal C_n$ be the category whose objects are the $\mathcal O_S$-vector bundles $\mathcal E$ with integrable $\Q$-connection which carry the structure of a sheaf of $\mathcal R^{(n)}$-modules satisfying (i) and (ii) (with $\epsilon^*\mathcal U$ replaced by $\mathcal E$).\\
Morphisms in $\mathcal C_n$ are defined to be the $\mathcal R^{(n)}$-linear and horizontal sheaf homomorphisms.
\end{definition}
With what we have already said it is clear that we obtain a covariant functor
\[\tag{\textbf{1.3.12}} F_n: U_n(X/S/\Q) \rightarrow \mathcal C_n, \qquad \mathcal U \mapsto \epsilon^*\mathcal U,\]
and we set out to prove that $F_n$ is an equivalence of categories. This is the manifestation in the de Rham realization of \cite{Be-Le}, 1.2.10 (v) (cf. also \cite{Hu-Ki}, Thm. A.2.5).

\subsubsection{Some additional structures}
Before we prove the announced theorem we need to discuss some more formalities which will be helpful for a clean definition of the quasi-inverse of $F_n$.\\
\newline
Let $n\geq 0$. First, it is clear that
\[\prod_{k=0}^n \mathrm{Sym}^k_{\mathcal O_X}\mathcal H_X= \pi^* \mathcal R^{(n)}=\pi^{-1}\Big(\prod_{k=0}^n \mathrm{Sym}^k_{\mathcal O_S}\mathcal H\Big) \otimes_{\pi^{-1}\mathcal O_S}\mathcal O_X\]
becomes a sheaf of rings on $X$ with unit given by $\frac{1}{n!}$ and multiplication analogously as for $\mathcal R^{(n)}$; it becomes a sheaf of commutative $\mathcal O_X$-algebras, inducing the original $\mathcal O_X$-module structure, via
\[\tag{\textbf{1.3.13}} \mathcal O_X \rightarrow \pi^*\mathcal R^{(n)}, \quad t \mapsto \frac{t}{n!}.\]
Consider for each $0\leq k \leq n$ the map
\[\mathrm{Sym}^k_{\mathcal O_X}\mathcal H_X \otimes_{\mathcal O_X} \mathcal L_n \rightarrow \mathcal L_k \otimes_{\mathcal O_X}\mathcal L_{n-k} \rightarrow \mathcal L_n,\]
where the first arrow is induced by the inclusion $\mathcal H_X \hookrightarrow \mathcal L_1$ in $(1.1.2)$ together with projection and the second by multiplication in symmetric powers.\\
If we let each component $\mathrm{Sym}^k_{\mathcal O_X}\mathcal H_X$ act on $\mathcal L_n$ by $(n-k)!$-times this composition, then we get on $\mathcal L_n$ the structure of $\pi^*\mathcal R^{(n)}$-module\footnote{To check the associativity requires some basic combinatorics.} giving the earlier defined multiplication on $\mathcal R^{(n)}$ after pullback via $\epsilon$. In what follows we will consider $\mathcal L_n$ with the $\pi^*\mathcal R^{(n)}$-module structure just explained.\\
\newline
Furthermore, if $\mathcal E$ is a $\mathcal R^{(n)}$-module on $S$, the pullback $\pi^*\mathcal E$ is a $\pi^*\mathcal R^{(n)}$-module.\\
Hence, the tensor product $\pi^*\mathcal E \otimes_{\pi^*\mathcal R^{(n)}} \mathcal L_n$ makes sense and is itself a $\pi^*\mathcal R^{(n)}$-module.
\begin{lemma}
Assume that $\mathcal E$ is a $\mathcal R^{(n)}$-module on $S$ and a $\mathcal O_S$-vector bundle in the induced $\mathcal O_S$-module structure (i.e. via $(1.3.10)$). Then the $\pi^* \mathcal R^{(n)}$-module $\pi^*\mathcal E \otimes_{\pi^*\mathcal R^{(n)}} \mathcal L_n$ is coherent over $\mathcal O_X$ in the induced $\mathcal O_X$-module structure (i.e. via $(1.3.13)$).
\end{lemma}
\begin{proof}
We have a canonical epimorphism of $\mathcal O_X$-modules
\[\tag{\textbf{1.3.14}} \pi^*\mathcal E \otimes_{\mathcal O_X} \mathcal L_n \rightarrow \pi^*\mathcal E \otimes_{\pi^*\mathcal R^{(n)}} \mathcal L_n\]
induced by the map $\mathcal O_X \rightarrow \pi^*\mathcal R^{(n)}$.\\
If we know that the right side of $(1.3.14)$ is $\mathcal O_X$-quasi-coherent, then the claim follows from the fact that the left side is a vector bundle on $X$ (use \cite{Li}, Ch. 5, Prop. 1.11).\\
Note that $\mathcal L_n$ is quasi-coherent over $\pi^*\mathcal R^{(n)}$ and that $\mathcal E$ is quasi-coherent over $\mathcal R^{(n)}$, both by \cite{EGAI}, Ch. I, Prop. (2.2.4). Hence, $\pi^*\mathcal E$ is quasi-coherent over $\pi^*\mathcal R^{(n)}$, and as the same holds for $\mathcal L_n$ the $\pi^*\mathcal R^{(n)}$-quasi-coherence of $\pi^*\mathcal E \otimes_{\pi^*\mathcal R^{(n)}} \mathcal L_n$ follows (cf. the comment after the proof of ibid.). Another application of ibid. yields the $\mathcal O_X$-quasi-coherence of $\pi^*\mathcal E \otimes_{\pi^*\mathcal R^{(n)}} \mathcal L_n$, and by what we said above thus also its coherence.
\end{proof}
\subsubsection{The quasi-inverse functor}
We now construct a functor
\[G_n: \mathcal C_n \rightarrow U_n(X/S/\Q)\]
which will turn out to provide a quasi-inverse for the functor $F_n$ of $(1.3.12)$.\\
\newline
Assume that $\mathcal E$ is an object of the category $\mathcal C_n$ (cf. Def. 1.3.11). We first define an integrable $\Q$-connection on the coherent $\mathcal O_X$-module $\pi^*\mathcal E \otimes_{\pi^*\mathcal R^{(n)}} \mathcal L_n$ (cf. Lemma 1.3.12) as follows:\\
Let $\pi^*\nabla_{\mathcal E}$ be the pullback of the connection $\nabla_{\mathcal E}$ on $\mathcal E$ and $\nabla_n$ as usual the connection on $\mathcal L_n$. Then
\[\tag{\textbf{1.3.15}} \pi^*\mathcal E \otimes_{\pi^*\mathcal R^{(n)}} \mathcal L_n \rightarrow \Omega^1_{X/\Q} \otimes_{\mathcal O_X} (\pi^*\mathcal E \otimes_{\pi^*\mathcal R^{(n)}} \mathcal L_n)\]
is determined on local sections $f$ resp. $l$ of $\pi^*\mathcal E$ resp. $\mathcal L_n$ by
\[f\otimes l \mapsto (\pi^*\nabla_{\mathcal E})(f) \otimes l + \nabla_n(l)\otimes f,\]
where the first resp. second summand is to read in
\[(\Omega^1_{X/\Q} \otimes_{\mathcal O_X} \pi^*\mathcal E) \otimes_{\pi^*\mathcal R^{(n)}} \mathcal L_n\]
resp. in
\[(\Omega^1_{X/\Q} \otimes_{\mathcal O_X} \mathcal L_n) \otimes_{\pi^*\mathcal R^{(n)}} \pi^*\mathcal E,\]
both canonically isomorphic to the right side of $(1.3.15)$.\\
That $(1.3.15)$ is well-defined is straightforwardly checked by using that the multiplication maps
\[\pi^*\mathcal R^{(n)} \otimes_{\mathcal O_X} \mathcal L_n \rightarrow \mathcal L_n\]
and
\[\pi^*\mathcal R^{(n)} \otimes_{\mathcal O_X} \pi^*\mathcal E \rightarrow \pi^* \mathcal E\]
are horizontal, the first by definition of the $\pi^*\mathcal R^{(n)}$-module structure of $\mathcal L_n$ and the second by definition of $\mathcal C_n$. In this way we make $\pi^*\mathcal E \otimes_{\pi^*\mathcal R^{(n)}} \mathcal L_n$ an object of $\textit{VIC}(X/\Q)$ (it is coherent by Lemma 1.3.12 and carries an integrable $\Q$-connection, hence is indeed a vector bundle).\\
\newline
Next, we define a unipotent filtration of length $n$ on $\pi^*\mathcal E \otimes_{\pi^*\mathcal R^{(n)}} \mathcal L_n$:\\
For this recall from 1.1 the filtration making $\mathcal L_n$ into an object of $U_n(X/S/\Q)$:
\[\mathcal L_n=A^0\mathcal L_n \supseteq A^1\mathcal L_n \supseteq...\supseteq A^n\mathcal L_n \supseteq A^{n+1}\mathcal L_n=0\]
with
\[A^i\mathcal L_n:=\mathrm{im}(\mathrm{Sym}^i_{\mathcal O_X}\mathcal H_X \otimes_{\mathcal O_X} \mathrm{Sym}^{n-i}_{\mathcal O_X}\mathcal L_1 \xrightarrow{\mathrm{mult}} \mathrm{Sym}^n_{\mathcal O_X}\mathcal L_1)\]
and quotients given by
\[\qquad A^i\mathcal L_n/A^{i+1}\mathcal L_n \simeq \mathrm{Sym}^i_{\mathcal O_X}\mathcal H_X, \quad i=0,...,n.\]
In terms of the $\pi^*\mathcal R^{(n)}$-module structure of $\mathcal L_n$ we may write the filtration objects as
\[A^i\mathcal L_n=\Big(\prod_{k=i}^n \mathrm{Sym}^k_{\mathcal O_X}\mathcal H_X \Big)\cdot \mathcal L_n, \qquad i=0,...,n,\]
where $\prod_{k=i}^n \mathrm{Sym}^k_{\mathcal O_X}\mathcal H_X$ is considered as a sheaf of ideals in $\pi^*\mathcal R^{(n)}$. In this way the $A^i\mathcal L_n$ and with them also their quotients $\mathrm{Sym}^i_{\mathcal O_X}\mathcal H_X$ become $\pi^*\mathcal R^{(n)}$-modules such that via $(1.3.13)$ the original $\mathcal O_X$- module structure is induced. For each $i=0,...,n$ we then have the exact sequence of $\pi^*\mathcal R^{(n)}= \prod_{k=0}^n \mathrm{Sym}^k_{\mathcal O_X}\mathcal H_X$-modules
\[\tag{\textbf{1.3.16}} 0 \rightarrow A^{i+1} \mathcal L_n \rightarrow A^i \mathcal L_n \rightarrow \mathrm{Sym}^i_{\mathcal O_X}\mathcal H_X \rightarrow 0.\]
Note that the $\prod_{k=0}^n \mathrm{Sym}^k_{\mathcal O_X}\mathcal H_X$-module structure of $\mathrm{Sym}^i_{\mathcal O_X}\mathcal H_X$ is given by (usual) multiplication with $n!$-times the $\mathcal O_X$-component; this implies that the functor
\[(\pi^*\mathcal R^{(n)} \textrm{-modules}) \rightarrow (\pi^*\mathcal R^{(n)}\textrm{-modules}), \quad \mathcal M \mapsto \mathcal M \otimes_{\pi^*\mathcal R^{(n)}} \mathrm{Sym}^i_{\mathcal O_X}\mathcal H_X\]
is exact: observe that when considering $\mathcal M \otimes_{\pi^*\mathcal R^{(n)}} \mathrm{Sym}^i_{\mathcal O_X}\mathcal H_X$ as $\mathcal O_X$-module via $(1.3.13)$ it can be identified with $\mathcal M \otimes_{\mathcal O_X} \mathrm{Sym}^i_{\mathcal O_X}\mathcal H_X$, where $\mathcal M$ is an $\mathcal O_X$-module via $(1.3.13)$.\\
With this we see that $(1.3.16)$ remains exact after tensoring with $\pi^*\mathcal E$:
\[\tag{\textbf{1.3.17}} 0 \rightarrow \pi^*\mathcal E \otimes_{\pi^*\mathcal R^{(n)}} A^{i+1} \mathcal L_n \rightarrow \pi^*\mathcal E \otimes_{\pi^*\mathcal R^{(n)}}A^i \mathcal L_n \rightarrow \pi^*\mathcal E \otimes_{\pi^*\mathcal R^{(n)}}\mathrm{Sym}^i_{\mathcal O_X}\mathcal H_X \rightarrow 0.\]
The $\pi^*\mathcal E \otimes_{\pi^*\mathcal R^{(n)}} A^i\mathcal L_n$, viewed as $\mathcal O_X$-modules by $(1.3.13)$, thus define a filtration of length $n$ of $\pi^*\mathcal E \otimes_{\pi^*\mathcal R^{(n)}} \mathcal L_n$ by $\mathcal O_X$-submodules with quotients given by
\[\pi^*\mathcal E \otimes_{\pi^*\mathcal R^{(n)}} \mathrm{Sym}^i_{\mathcal O_X}\mathcal H_X \simeq \pi^*\mathcal E \otimes_{\mathcal O_X} \mathrm{Sym}^i_{\mathcal O_X}\mathcal H_X \simeq \pi^*(\mathcal E \otimes_{\mathcal O_S} \mathrm{Sym}^i_{\mathcal O_S}\mathcal H);\]
note that $\pi^*\mathcal E$ in the previous chain is an $\mathcal O_X$-module via $(1.3.13)$ and that this is exactly its original $\mathcal O_X$-structure because of property (i) in the definition of the category $\mathcal C_n$.\\
Exactly the same argument as in Lemma 1.3.12 shows that the $\pi^*\mathcal E \otimes_{\pi^*\mathcal R^{(n)}} A^i\mathcal L_n$ are $\mathcal O_X$-coherent. Moreover, they are stable under the integrable $\Q$-connection of $\pi^*\mathcal E \otimes_{\pi^*\mathcal R^{(n)}} \mathcal L_n$, hence locally free.\\
Altogether, it follows that the defined filtration makes $\pi^*\mathcal E \otimes_{\pi^*\mathcal R^{(n)}} \mathcal L_n$ an object of $U_n(X/S/\Q)$.\\
In this way we obtain a covariant functor
\[\tag{\textbf{1.3.18}} G_n: \mathcal C_n \rightarrow U_n(X/S/\Q), \quad \mathcal E \mapsto \pi^*\mathcal E \otimes_{\pi^*\mathcal R^{(n)}} \mathcal L_n.\]
\subsubsection{The equivalence result}
With all these preparations we are finally in the position to make sense of and prove the statement of
\begin{theorem}
The functor
\[F_n: U_n(X/S/\Q) \rightarrow \mathcal C_n, \quad \mathcal U \mapsto \epsilon^*\mathcal U\]
is an equivalence of categories with quasi-inverse
\[G_n:  \mathcal C_n \rightarrow U_n(X/S/\Q), \quad \mathcal E \mapsto \pi^*\mathcal E \otimes_{\pi^*\mathcal R^{(n)}} \mathcal L_n.\]
\end{theorem}
\begin{proof}
It is easy to check that $F_n \circ G_n \simeq \id$, hence it remains to show $G_n \circ F_n \simeq \id$.\\
For this we explicate the argument sketched very briefly in \cite{Hu-Ki}, proof of Thm. A.2.5.\\
Let us start with an object $\mathcal U$ of $U_n(X/S/\Q)$ and define an arrow in $\textit{VIC}(X/\Q)$
\[\tag{\textbf{1.3.19}} \pi^*\epsilon^* \mathcal U \otimes_{\pi^*\mathcal R^{(n)}} \mathcal L_n \rightarrow \mathcal U\]
as follows:\\
Note that $\pi^* \epsilon^*\mathcal U$ identifies by Thm. 1.3.6 (i) with $\pi^* \pi_* \underline{\Hom}_{\mathcal D_{X/S}}(\mathcal L_n, \mathcal U)$. Let us define an arrow
\[\tag{\textbf{1.3.20}} \pi^* \pi_* \underline{\Hom}_{\mathcal D_{X/S}}(\mathcal L_n, \mathcal U) \otimes_{\pi^*\mathcal R^{(n)}} \mathcal L_n \rightarrow \mathcal U,\]
or in other words
\[(\pi^{-1} \pi_* \underline{\Hom}_{\mathcal D_{X/S}}(\mathcal L_n, \mathcal U) \otimes_{\pi^{-1}\mathcal O_S}\mathcal O_X) \otimes_{\pi^*\mathcal R^{(n)}} \mathcal L_n \rightarrow \mathcal U.\]
Noting that $\pi$ is an open map we define the last morphism at the level of presheaves on open subsets $V \subseteq X$ by the rule
\[(\Hom_{\mathcal D_{\pi^{-1}(\pi(V))/\pi(V)}}(\mathcal L_n, \mathcal U) \otimes_{\mathcal O_S(\pi(V))}\mathcal O_X(V)) \otimes_{ \lbrack \mathcal R^{(n)}(\pi(V)) \otimes_{\mathcal O_S(\pi(V))}\mathcal O_X(V) \rbrack} \mathcal L_n(V) \rightarrow \mathcal U(V),\]
\[(f \otimes t) \otimes l \mapsto t \cdot (f_{|V})(l),\]
which is easily seen to be well-defined and which induces a $\mathcal O_X$-linear map $(1.3.20)$. It is a routine calculation to check that $(1.3.20)$ is horizontal for the $\Q$-connections of both sides. This finishes the definition of the arrow $(1.3.19)$ in $\textit{VIC}(X/\Q)$. What remains to show is that it is an isomorphism (naturality is clear). As $(1.3.19)$ induces (by construction) the identity on $\epsilon^*\mathcal U$ after pullback via $\epsilon$ and as its cokernel is a vector bundle (recall once more that $\textit{VIC}(X/\Q)$ is abelian) we may conclude the proof by the following general argument:
\end{proof}
\begin{lemma}
Let $\psi: \mathcal V \rightarrow \mathcal W$ be a morphism of vector bundles on $X$ such that the induced map $\epsilon^*(\psi): \epsilon^*\mathcal V \rightarrow \epsilon^*\mathcal W$
is an isomorphism and such that the cokernel of $\psi$ is a vector bundle. Then $\psi$ is an isomorphism.
\end{lemma}
\begin{proof}
At first, it is an easy application of the Nakayama lemma (together with \cite{EGAI}, Ch. I, Prop. (3.4.6)) to see that $\psi$ is an isomorphism if and only if for all $s \in S$ the induced map
\[\psi_s: \mathcal V_{|X_s} \rightarrow \mathcal W_{|X_s}\]
on the fiber $X_s$ over $s$ is an isomorphism. We may hence assume from the beginning that we are in the situation of an abelian variety $X$ over a field $k$. The hypothesis that $\psi$ induces an isomorphism after pullback to $\Spec(k)$ via $\epsilon$ immediately shows that the vector bundles $\mathcal V$ and $\mathcal W$ must have the same rank (which is constant as $X$ is connected). By a standard argument we thus only need to see that $\psi$ is an epimorphism. Now consider the pullback of $\coker(\psi)$ along $\epsilon$ to $\Spec(k)$ and use the hypothesis that $\epsilon^*(\psi)$ is an isomorphism together with the Nakayama lemma to see that the stalk of $\coker(\psi)$ vanishes in the zero point; as $\coker(\psi)$ is a vector bundle on $X$ it must hence be zero.
\end{proof}
\subsection{Some categorical structure results}
The preceding theorem can in particular be used to obtain non-trivial information about categories of unipotent vector bundles with integrable connection.
\begin{corollary}
For each $n\geq 0$ the category $U_n(X/S/\Q)$ is abelian.
\end{corollary}
\begin{proof}
One checks without problems that the category $\mathcal C_n$ is abelian. Now use Thm. 1.3.13.
\end{proof}
\begin{corollary}
The category $U(X/S/\Q)$ is abelian.
\end{corollary}
\begin{proof}
Note that $U(X/S/\Q)$ is a full subcategory of the abelian category $\textit{VIC}(X/\Q)$ and that then the only potentially nontrivial task is to see that for each morphism $\mathcal U \rightarrow \mathcal V$ of objects of $U(X/S/\Q)$ its kernel and cokernel is again in $U(X/S/\Q)$. But of course there is a suitable $n\geq 0$ such that $\mathcal U$ and $\mathcal V$ both are in $U_n(X/S/\Q)$. Now use Cor. 1.3.15 to conclude.
\end{proof}
Using the vocabulary of tensor categories (for which we refer to \cite{Sh}, Ch. 1, 1.1) we may view $\textit{VIC}(X/\Q)$ as rigid abelian tensor category (with unit object $(\mathcal O_X,\mathrm{d})$ and the usual tensor product resp. internal $\Hom$-objects). With Cor. 1.3.16 and Lemma 1.3.2 we may summarize the knowledge we have won about $U(X/S/\Q)$:
\begin{theorem}
The category $U(X/S/\Q)$ is a rigid abelian tensor subcategory of $\textit{VIC}(X/\Q)$ which is closed under extensions. \qquad \qed
\end{theorem}

\begin{remark}
Let us mention what we get in the special case that $S$ is the spectrum of a field $k$.\\
Our general assumptions on $S/\Spec(\Q)$ tantamount to requiring that $k$ is a number field.\\
As $\Spec(k) \rightarrow \Spec(\Q)$ then is étale we note that connections relative $\Spec(\Q)$ are the same as connections relative $\Spec(k)$. It is easy to see that $U(X/S/\Q)$ then becomes the category of vector bundles on $X$ with integrable $k$-connection which have a filtration of finite length by subbundles stable under the connection and quotients isomorphic to $(\mathcal O_X,\mathrm{d})$; write $U(X/k)$ for this category.\\
Denoting by $Vec_{fd}(k)$ the category of finite-dimensional $k$-vector spaces the functor
\[\epsilon^*: U(X/k) \rightarrow Vec_{fd}(k)\]
then is exact and faithful: exactness is obvious and faithfulness easily follows from Thm. 1.3.13 by choosing for two objects of $U(X/k)$ a common length of unipotency.\\
Because of the standard fact $k=\Gamma(X,\mathcal O_X)$ it is also clear that $\mathrm {End}_{U(X/k)}((\mathcal O_X,\mathrm{d}))=k$.\\
Together with Thm. 1.3.17 this shows that $U(X/k)$ is a (neutral) Tannakian category over $k$ with fiber functor given by $\epsilon^*$ (for the notion of a Tannakian category cf. \cite{Sh}, Ch. 1, Def. 1.1.7).
\end{remark}

\markright{\uppercase{The formalism of the logarithm sheaves and the elliptic...}}
\section{The invariance results for the logarithm sheaves}
\markright{\uppercase{The formalism of the logarithm sheaves and the elliptic...}}
\subsection{A technical preparation}
The zero fiber $\epsilon^* \mathcal U$ of a bundle $\mathcal U \in U_1(X/S/\Q)$ carries a module structure over $\mathcal O_S\oplus \mathcal H$ whose restriction to $\mathcal O_S$ coincides with the usual $\mathcal O_S$-multiplication on $\epsilon^*\mathcal U$ (cf. $(1.3.11)$ and $(1.3.12)$).\\
In addition to the explanations subsequent to $(1.3.11)$ the following auxiliary result gives another description of how the multiplication coming from the $\mathcal H$-component looks like.
\begin{lemma}
Let $\mathcal U$ be an object of $\textit{VIC}(X/\Q)$ and $\mathcal Y_0,\mathcal Y_1$ objects of $\textit{VIC}(S/\Q)$ sitting in an exact sequence of $\mathcal D_{X/\Q}$-modules
\[\tag{\textbf{1.4.1}} 0\rightarrow \pi^*\mathcal Y_1 \rightarrow \mathcal U \xrightarrow{p} \pi^*\mathcal Y_0 \rightarrow 0.\]
Denote by $\psi: \mathcal Y_0 \rightarrow \mathcal H^\vee \otimes_{\mathcal O_S}\mathcal Y_1$ the $\mathcal D_{S/\Q}$-linear map induced by the first edge morphism in the long exact sequence of de Rham cohomology relative $S$ for $(1.4.1)$; define a $\mathcal D_{S/\Q}$-linear map $\tau$ as the composition
\[\tau: \mathcal H \otimes_{\mathcal O_S} \mathcal Y_0 \xrightarrow{\id\otimes \psi} \mathcal H \otimes_{\mathcal O_S} \mathcal H^\vee \otimes_{\mathcal O_S} \mathcal Y_1 \xrightarrow{\mathrm{eval}\otimes \id} \mathcal Y_1.\]
On the other hand, by pullback of $(1.4.1)$ via $\epsilon$ one obtains the exact $\mathcal D_{S/\Q}$-linear sequence
\[\tag{\textbf{1.4.2}} 0\rightarrow \mathcal Y_1 \rightarrow \epsilon^*\mathcal U \xrightarrow{\epsilon^*(p)} \mathcal Y_0 \rightarrow 0.\]
Let $\xi \in \Gamma(S,\epsilon^*\mathcal U)$ and consider the section $(\epsilon^*(p))(\xi) \in \Gamma(S,\mathcal Y_0)$; together with $\tau$ it defines the $\mathcal O_S$-linear map
\[\eta: \mathcal H \rightarrow \mathcal Y_1, \quad s\mapsto \tau(s\otimes (\epsilon^*(p))(\xi)).\]
Then, if $f:\mathcal L_1\rightarrow \mathcal U$ is the unique $\mathcal D_{X/S}$-linear map with $\epsilon^*(f)(1)=\xi$ (cf. Thm. 1.3.6 (i)), the restriction of
\[\epsilon^*(f): \mathcal O_S \oplus \mathcal H \rightarrow \epsilon^*\mathcal U\]
to the direct summand $\mathcal H$ (which defines the multiplication of $\xi$ by sections of $\mathcal H$) is equal to $\eta$ composed with the inclusion of $(1.4.2)$.
\end{lemma}
\begin{proof}
We define a $\mathcal D_{X/S}$-linear map
\[g: \mathcal O_X \rightarrow \pi^*\mathcal Y_0, \quad 1 \mapsto \pi^*((\epsilon^*(p))(\xi))\]
and observe that we then have a commutative square
\begin{equation*}
\begin{xy}
\xymatrix{
\mathcal L_1 \ar[r] \ar[d]_{f} & \mathcal O_X \ar[d]^{g}\\
\mathcal U \ar[r]^{p \ \ } & \pi^*\mathcal Y_0}
\end{xy}
\end{equation*}
in which the upper map is given by the projection in $(1.1.2)$; the commutativity follows because both arrows $\mathcal L_1 \rightarrow \pi^*\mathcal Y_0$ are $\mathcal D_{X/S}$-linear with the property that after pullback via $\epsilon$ they send $1$ to $(\epsilon^*(p))(\xi)$, hence by Thm. 1.3.6 (i) must be equal.\\
The induced commutative square of $\mathcal O_S$-linear maps
\begin{equation*}
\begin{xy}
\xymatrix{
\mathcal \epsilon^*\mathcal L_1 \ar[r] \ar[d]_{\epsilon^*(f)} & \mathcal O_S \ar[d]^{\epsilon^*(g)}\\
\epsilon^* \mathcal U \ar[r]^{\epsilon^*(p)} & \mathcal Y_0}
\end{xy}
\end{equation*}
then permits a unique commutative continuation
\begin{equation*}
\begin{xy}
\xymatrix{
0 \ar[r] & \mathcal H \ar[d]_{h} \ar[r] & \epsilon^*\mathcal L_1  \ar[r] \ar[d]_{\epsilon^*(f)} & \mathcal O_S \ar[d]^{\epsilon^*(g)} \ar[r] & 0 \\
0 \ar[r] & \mathcal Y_1 \ar[r] & \epsilon^* \mathcal U \ar[r]^{\epsilon^*(p)} & \mathcal Y_0 \ar[r] & 0}
\end{xy}
\end{equation*}
such that the upper resp. lower row is given by pullback along $\epsilon$ of $(1.1.2)$ resp. by $(1.4.2)$.\\
We next consider the diagram of $\mathcal D_{X/S}$-linear maps
\begin{equation*}
\begin{xy}
\xymatrix{
0 \ar[r] & \mathcal H_X \ar[d]_{\pi^*(h)} \ar[r] & \mathcal L_1  \ar[r] \ar[d]_{f} & \mathcal O_X \ar[d]^{g} \ar[r] & 0 \\
0 \ar[r] &  \pi^*\mathcal Y_1 \ar[r] & \mathcal U \ar[r]^{p \ \ } & \pi^* \mathcal Y_0 \ar[r] & 0}
\end{xy}
\end{equation*}
with upper resp. lower row given by $(1.1.2)$ resp. by $(1.4.1)$. The claim is that it commutes.\\
With what we have already said above it only remains to prove the commutativity of the left square, which tantamounts to showing that a certain $\mathcal D_{X/S}$-linear arrow $\mathcal H_X \rightarrow \mathcal U$ (namely the difference of the two maps in the square) which is zero after pullback via $\epsilon$ is already zero. Taking into account that $X$ is integral\footnote{Use that $S$ is integral (cf. footnote 5) and \cite{Li}, Ch. 4, Prop. 3.8, to conclude that $X$ is integral.} and that $\mathcal U$ is a vector bundle one easily reduces the question to the situation $S=\Spec(k)$ with $k$ a field of characteristic zero. In this case \cite{Bert-Og}, §2, Prop. 2.16, yields that the considered map $\mathcal H_X \rightarrow \mathcal U$ is zero not only in the fiber, but already in the stalk of the zero point of $X$. As $X$ is integral and $\mathcal U$ is a vector bundle one can conclude from this that $\mathcal H_X \rightarrow \mathcal U$ is indeed zero. The desired commutativity is thus shown.\\
From the long exact sequence of de Rham cohomology (cf. \cite{Kat2}, $(2.0)$), applied to the preceding commutative diagram, we obtain the commutative square of $\mathcal O_S$-linear morphisms
\begin{equation*}
\begin{xy}
\xymatrix{
\mathcal O_S \ar[r]^{\mathrm{can} \quad \ \ } \ar[d] & \mathcal H^\vee \otimes_{\mathcal O_S}\mathcal H \ar[d]^{\id\otimes h}\\
\mathcal Y_0 \ar[r]^{\psi \quad \ \ \ } & \mathcal H^\vee \otimes_{\mathcal O_S}\mathcal Y_1}
\end{xy}
\end{equation*}
in which the left vertical arrow is given by $1\mapsto (\epsilon^*(p))(\xi)$ (by definition of $g$) and where the above horizontal arrow is the standard map (as follows from the definition of $\mathcal L_1$).\\
Tensoring with $\mathcal H$ and composing the horizontal arrows with $\mathrm{eval}\otimes \id$ we obtain the diagram
\begin{equation*}
\begin{xy}
\xymatrix{
\mathcal H \ar[r]^{\id \otimes \mathrm{can} \quad \quad \quad } \ar[d] & \mathcal H \otimes_{\mathcal O_S}\mathcal H^\vee \otimes_{\mathcal O_S} \mathcal H \ar[d]^{\id\otimes \id\otimes h}\ar[r]^{ \qquad \qquad \mathrm{eval}\otimes \id} & \mathcal H  \ar[d]^{h}\\
\mathcal H \otimes_{\mathcal O_S}\mathcal Y_0 \ar[r]^{\id\otimes \psi \quad \quad} &  \mathcal H\otimes_{\mathcal O_S}\mathcal H^\vee \otimes_{\mathcal O_S} \mathcal Y_1 \ar[r]^{\qquad \quad \mathrm{eval}\otimes \id} & \mathcal Y_1}
\end{xy}
\end{equation*}
in which the two small squares commute, implying commutativity of the whole. If we note that the upper horizontal composition is the identity, the lower horizontal composition is $\tau$ and the left vertical map is given by $s\mapsto s \otimes (\epsilon^*(p))(\xi)$, we get $h=\eta$. This implies the claim of the lemma.
\end{proof}

\subsection{The invariance results}
We now prove the fundamental fact that the logarithm sheaves of two abelian schemes become canonically identified under pullback by isogenies, from which we derive in particular the invariance of $\mathcal L_n$ under translation by torsion sections. This implies that the fiber of $\mathcal L_n$ in a torsion section is canonically isomorphic to its zero fiber, a property that will be important later when we define and compute the specialization of the polylogarithm along torsion sections.\\
We finally append a brief observation concerning compatibility of the logarithm sheaves under base change, supplementing the content of Lemma 1.1.7.

\subsubsection{Invariance under isogenies}
Let us asssume that $u: X \rightarrow X'$ is an isogeny\footnote{By an isogeny we mean a surjective and finite homomorphism over $S$; it is automatically flat and hence finite locally free. We thus have the notion of its degree $\mathrm{deg}(u)$, defined as the rank of the $\mathcal O_{X'}$-vector bundle $u_*\mathcal O_X$. A priori, $\mathrm{deg}(u)$ is a locally constant function on $X'$ and hence constant in our setting because $X'$ is connected ($X'$ is integral by the same argument as in footnote 10). Moreover, as we are in characteristic zero each isogeny is an étale morphism.} of abelian schemes over $S$ and use primed notation $\pi', \epsilon', \mathcal H', \mathcal L_n'$ etc. for the usual objects when they refer to the abelian scheme $X'/S$.\\
\newline
Pullback via $u$ of the exact sequence associated with $\mathcal L_1'$
\[0\rightarrow (\pi')^*\mathcal H' \rightarrow \mathcal L_1' \rightarrow \mathcal O_{X'} \rightarrow 0\]
gives an exact sequence of $\mathcal D_{X/\Q}$-modules
\[\tag{\textbf{1.4.3}} 0\rightarrow \pi^*\mathcal H' \rightarrow u^*\mathcal L_1' \rightarrow \mathcal O_X \rightarrow 0,\]
exhibiting $u^*\mathcal L_1'$ as object of $U_1(X/S/\Q)$. By Thm. 1.3.6 (i) we may determine a $\mathcal D_{X/S}$-linear map
\[\tag{\textbf{1.4.4}} f_1: \mathcal L_1 \rightarrow u^*\mathcal L_1'\]
by choosing a global $S$-section of $\epsilon^*u^*\mathcal L_1' \simeq (\epsilon')^*\mathcal L_1' \simeq \mathcal O_S \oplus \mathcal H'$, for which we take the section $1$.\\
By part (ii) of the same theorem we know that $f_1$ is even $\mathcal D_{X/\Q}$-linear.
\begin{theorem}
The canonical $\mathcal D_{X/\Q}$-linear map $f_1: \mathcal L_1 \rightarrow u^*\mathcal L_1'$ defined in $(1.4.4)$ is an isomorphism.\\
The induced morphism $\epsilon^*(f_1): \mathcal O_S \oplus \mathcal H \rightarrow \mathcal O_S \oplus \mathcal H'$ is given by
\[\epsilon^*(f_1)=\id \oplus (u^*_{{\mathrm{dR}}})^\vee,\]
where
\[u^*_{{\mathrm{dR}}}: H^1_{\mathrm{dR}}(X'/S) \rightarrow H^1_{\mathrm{dR}}(X/S)\]
denotes the canonical morphism on de Rham cohomology induced by $u$.
\end{theorem}
\begin{proof}
For the isomorphy of $f_1$ it suffices - by the equivalence of categories in Thm. 1.3.13 - to show that the map
\[\epsilon^*(f_1):\mathcal O_S \oplus \mathcal H \rightarrow \mathcal O_S \oplus \mathcal H',\]
induced by $f_1$ and the identification $\epsilon^*u^*\mathcal L_1' \simeq (\epsilon')^*\mathcal L_1'\simeq \mathcal O_S \oplus \mathcal H'$, is an isomorphism.\\
We will determine $\epsilon^*(f_1)$ explicitly, see that it is given as stated in the second claim of the theorem and then argue that this is indeed an isomorphism.\\
By definition $\epsilon^*(f_1)$ sends the global $S$-section $1$ of $\mathcal O_S$ to itself.\\
What it does on the direct summand $\mathcal H$ can be determined accurately with the help of Lemma 1.4.1:\\
The map $\eta$ in the claim of that lemma writes in the present situation as $\eta: \mathcal H \rightarrow \mathcal H', s\mapsto \tau(s\otimes 1)$, where $\tau$ is given by the composition
\[\tau:\mathcal H \otimes_{\mathcal O_S}\mathcal O_S \xrightarrow{\id\otimes \psi} \mathcal H \otimes_{\mathcal O_S}\mathcal H^\vee \otimes_{\mathcal O_S}\mathcal H' \xrightarrow{\mathrm{eval}\otimes \id} \mathcal H'.\]
But the map $\psi: \mathcal O_S \rightarrow \mathcal H^\vee \otimes_{\mathcal O_S}\mathcal H'$ was defined to be the first edge morphism of de Rham cohomology for the exact sequence $(1.4.3)$. One can check that $\psi$ is nothing else than the map associated with the canonical arrow
\[u^*_{{\mathrm{dR}}}: H^1_{\mathrm{dR}}(X'/S) \rightarrow H^1_{\mathrm{dR}}(X/S),\]
and hence $\eta: \mathcal H \rightarrow \mathcal H'$ is obviously equal to $(u^*_{{\mathrm{dR}}})^\vee$.\\
Altogether, by Lemma 1.4.1 we conclude that $\epsilon^*(f_1)$ is given on the summand $\mathcal H$ by $(u^*_{{\mathrm{dR}}})^\vee$ followed by the canonical inclusion, hence 
\[\epsilon^*(f_1)= \id \oplus (u^*_{{\mathrm{dR}}})^\vee: \mathcal O_S \oplus \mathcal H \rightarrow \mathcal O_S \oplus \mathcal H'.\]
But the map $(u^*_{{\mathrm{dR}}})^\vee$ - and hence $\epsilon^*(f_1)$ - is an isomorphism:\\
For this let $v: X'\rightarrow X$ denote the isogeny characterized by $u\circ v=[\mathrm{deg}(u)], v\circ u =[\mathrm{deg}(u)]$, where $[\mathrm{deg}(u)]$ means the multiplication map(s) by the degree $\mathrm{deg}(u)$ of $u$. Then the compositions
\[\mathcal H \xrightarrow{(u^*_{{\mathrm{dR}}})^\vee} \mathcal H' \xrightarrow{(v^*_{{\mathrm{dR}}})^\vee} \mathcal H\]
\[\mathcal H' \xrightarrow{(v^*_{{\mathrm{dR}}})^\vee} \mathcal H \xrightarrow{(u^*_{{\mathrm{dR}}})^\vee} \mathcal H'\]
are isomorphisms because they are given by multiplication with $\mathrm{deg}(u)$ which is invertible on $S$; this suffices to conclude.
\end{proof}
Taking for each $n\geq 0$ the $n$-th symmetric power of the isomorphism $f_1$ in Thm. 1.4.2 we obtain an induced $\mathcal D_{X/\Q}$-linear isomorphism
\[\mathcal L_n \xrightarrow{\sim} u^*\mathcal L_n'.\]
Its pullback via $\epsilon$ gives the $\mathcal D_{S/\Q}$-linear isomorphism
\[\prod_{k=0}^n\mathrm{Sym}^k_{\mathcal O_S}\mathcal H \simeq \epsilon^*\mathcal L_n \xrightarrow{\sim} \epsilon^*u^*\mathcal L_n' \simeq (\epsilon')^*\mathcal L_n' \simeq \prod_{k=0}^n\mathrm{Sym}^k_{\mathcal O_S}\mathcal H'\]
which is the $n$-th symmetric power of the map $\id \oplus (u^*_{{\mathrm{dR}}})^\vee: \mathcal O_S \oplus \mathcal H \xrightarrow{\sim} \mathcal O_S \oplus \mathcal H'$ in the theorem.\\
As this symmetric power sends the section $\frac{1}{n!}$ of $\mathcal O_S$ to itself we see in particular:
\begin{corollary}
The canonical $\mathcal D_{X/\Q}$-linear map
\[f_n: \mathcal L_n \rightarrow u^*\mathcal L_n',\]
characterized by the condition $1^{(n)}\mapsto (1^{(n)})'$ in the induced map
\[\epsilon^*(f_n): \epsilon^*\mathcal L_n \rightarrow \epsilon^*u^*\mathcal L_n' \simeq (\epsilon')^*\mathcal L_n',\]
(cf. Thm. 1.3.6 \footnote{Of course, $u^*\mathcal L_n'$ is an object of $U_n(X/S/\Q)$: simply take the pullback via $u$ of the natural unipotent filtration for $\mathcal L_n'$.}), is an isomorphism. Under the splittings of the logarithm sheaves the arrow
\[\epsilon^*(f_n): \prod_{k=0}^n\mathrm{Sym}^k_{\mathcal O_S}\mathcal H \rightarrow \prod_{k=0}^n\mathrm{Sym}^k_{\mathcal O_S}\mathcal H'\]
is given on each factor as the map
\[\mathrm{Sym}^k_{\mathcal O_S}\mathcal H \rightarrow \mathrm{Sym}^k_{\mathcal O_S}\mathcal H'\]
induced on symmetric powers by $(u_{{\mathrm{dR}}}^*)^\vee: \mathcal H \rightarrow \mathcal H'$. \qquad \qed
\end{corollary}
A special kind of isogenies are the $N$-multiplication endomorphisms of the abelian scheme $X$. We note this important case separatedly: 
\begin{corollary}
Let $N\neq 0$ be an integer and $[N]: X \rightarrow X$ the $N$-multiplication isogeny of $X$.\\
Then for each $n\geq 0$ the canonical $\mathcal D_{X/\Q}$-linear map
\[\mathcal L_n \rightarrow [N]^*\mathcal L_n,\]
characterized by $1^{(n)}\mapsto 1^{(n)}$ after pullback via $\epsilon$:
\[\epsilon^*\mathcal L_n \rightarrow \epsilon^*[N]^*\mathcal L_n \simeq \epsilon^*\mathcal L_n,\]
is an isomorphism. The induced morphism on the product of symmetric powers is given on each factor as the $N^k$-multiplication map\footnote{Note the standard fact (which one can also deduce from the later Prop. 2.5.1) that $[N]^*_{{\mathrm{dR}}}$ is multiplication by $N$.}
\[\mathrm{Sym}^k_{\mathcal O_S}\mathcal H \xrightarrow{\cdot N^k} \mathrm{Sym}^k_{\mathcal O_S}\mathcal H.\]
\qquad \qed
\end{corollary}

\subsubsection{Invariance under translation by torsion sections}

The preceding result implies another crucial invariance property of the logarithm sheaves:\\
For this let $t: S \rightarrow X$ be a torsion section, i.e. $t \in X(S)$ and $N\cdot t=0$ for some integer $N\neq0$. 
Denote with
\[T_t: X \rightarrow X\]
the translation by $t$: for each $S$-scheme $Z$ it is defined in $Z$-rational points as the map $X(Z) \rightarrow X(Z)$ given by $\id+t^Z$, where $t^Z$ means the element in $X(Z)$ naturally induced by $t \in X(S)$.\\
For any $N$ as above we have $[N]\circ T_t=[N]$, and combining this with the natural identification
\[\mathcal L_n \simeq [N]^*\mathcal L_n\]
of Cor. 1.4.4 we obtain the following isomorphism of $\mathcal D_{X/\Q}$-modules
\[\tag{\textbf{1.4.5}} T_t^*\mathcal L_n \simeq T_t^*[N]^*\mathcal L_n \simeq [N]^* \mathcal L_n \simeq  \mathcal L_n.\]
The composite isomorphism expresses the fundamental fact that $\mathcal L_n$ is invariant under translation with torsion sections; let us note that this invariance is indeed canonical:
\begin{lemma}
For each torsion section $t \in X(S)$ the $\mathcal D_{X/\Q}$-linear isomorphism
\[T_t^*\mathcal L_n \simeq \mathcal L_n\]
in $(1.4.5)$ is canonical, i.e. independent of the choice of the integer $N\neq 0$ annihilating $t$.\\
In particular, after pullback via $\epsilon$ we have a canonical $\mathcal D_{S/\Q}$-linear identification
\[t^*\mathcal L_n  \simeq \epsilon^*\mathcal L_n.\]
\end{lemma}
\begin{proof}
Let $M$ be another non-zero integer with $M\cdot t=0$. Then $MN$ has the same property, and we show that the isomorphism $T_t^*\mathcal L_n \simeq \mathcal L_n$ in $(1.4.5)$ constructed by using $N$ coincides with the one constructed by using $MN$; by exchanging $N$ and $M$ we thus conclude that the isomorphisms $(1.4.5)$ for $N$ and for $M$ coincide.\\
Consider the following diagram in which all arrows are $\mathcal D_{X/\Q}$-linear isomorphisms:
\begin{equation*}
\begin{xy}
\xymatrix@C-0.3cm{
T_t^*\mathcal L_n \ar@{-}[dd]_{\id} \ \ar@{-}[r] & T_t^*[N]^*\mathcal L_n \ar@{-}[r] \ar@{-}[d] & [N]^*\mathcal L_n \ar@{-}[r] \ar@{-}[d] &  \mathcal L_n \ar@{-}[dd]^{\id}\\
 & T_t^*[N]^*[M]^*\mathcal L_n \ar@{-}[d]^{\mathrm{can}} & [N]^*[M]^*\mathcal L_n  \ar@{-}[d]^{\mathrm{can}} & \\
T_t^*\mathcal L_n\ \ar@{-}[r] & T_t^*[MN]^*\mathcal L_n \ar@{-}[r]  & [MN]^*\mathcal L_n \ar@{-}[r] &  \mathcal L_n}
\end{xy}
\end{equation*}
The outer lines are $(1.4.5)$ for $N$ resp. $MN$, the arrows denoted with $\mathrm{can}$ are the obvious standard identifications and the two vertical arrows without label are the maps induced by the isomorphism $\mathcal L_n \simeq [M]^*\mathcal L_n$ of Cor. 1.4.4 by applying $T_t^*[N]^*$ resp. $[N]^*$.\\
The right square commutes: by Thm. 1.3.13 this may be checked after pullback via $\epsilon$, where it is straightforwardly seen. The left square comes from the right by applying the functor $T_t^*$, hence also commutes. The commutativity of the middle square easily boils down to the functoriality of the standard isomorphism $g^*f^*\mathcal F \simeq (f\circ g)^*\mathcal F$ for morphisms of schemes $f,g$ and sheaves of modules $\mathcal F$. We deduce that the whole diagram commutes, which is what we wanted to show.
\end{proof}

\begin{remark}
The isomorphism $t^*\mathcal L_n  \xrightarrow{\sim} \epsilon^*\mathcal L_n$ recorded in the preceding lemma can also be obtained as follows:\\
Apply $t^*$ to $\mathcal L_n \xrightarrow{\sim} [N]^*\mathcal L_n$ to get $t^*\mathcal L_n \xrightarrow{\sim} \epsilon^*\mathcal L_n$; then compose this isomorphism with the map $\epsilon^*\mathcal L_n \xrightarrow{\sim} \epsilon^*\mathcal L_n$ given by applying $\epsilon^*$ to $[N]^*\mathcal L_n \xrightarrow{\sim} \mathcal L_n$.
\end{remark}

\subsubsection{Invariance under base change}
We consider the situation
\begin{equation*}
\begin{xy}
\xymatrix{
X' \ar[r]^{\pi'} \ar[d]_{g} & S' \ar[d]^{f} \\
X \ar[r]^{\pi} & S}
\end{xy}
\end{equation*}
of $(1.1.8)$ and use the notations $\epsilon', \mathcal H', \mathcal L_n'$ for the usual objects when they refer to the abelian scheme $X'/S'$. Further, let us recall that the canonical (horizontal) map
\[\tag{\textbf{1.4.6}} f^* H^1_{\mathrm{dR}}(X/S) \rightarrow H^1_{\mathrm{dR}}(X'/S')\]
is an isomorphism.\\
It is clear that for each $n\geq 0$ the bundle $g^*\mathcal L_n$ defines an object of $U_n(X'/S'/\Q)$, and by Thm. 1.3.6 we may hence determine a $\mathcal D_{X'/\Q}$-linear map
\[\tag{\textbf{1.4.7}} \mathcal L_n' \rightarrow g^*\mathcal L_n\]
by choosing a global horizontal $S'$-section of $(\epsilon')^*g^*\mathcal L_n \simeq f^* \epsilon^* \mathcal L_n$, for which we take $f^*(1^{(n)})$.\\
In complete analogy to the arguments that led to Cor. 1.4.3 one arrives at  
\begin{proposition}
The canonical $\mathcal D_{X'/\Q}$-linear map $\mathcal L_n' \rightarrow g^*\mathcal L_n$ defined in $(1.4.7)$ is an isomorphism.\\
Under the splittings of the logarithm sheaves the associated morphism in the zero fibers
\[\prod_{k=0}^n\mathrm{Sym}^k_{\mathcal O_{S'}}\mathcal H' \rightarrow \prod_{k=0}^n\mathrm{Sym}^k_{\mathcal O_{S'}} f^* \mathcal H\]
is given on each factor as the map
\[\mathrm{Sym}^k_{\mathcal O_{S'}}\mathcal H' \rightarrow \mathrm{Sym}^k_{\mathcal O_{S'}} f^* \mathcal H\]
induced on symmetric powers by the dual of $(1.4.6)$. \qquad \qed
\end{proposition}

\section{The elliptic polylogarithm and its $D$-variant}
\markright{\uppercase{The formalism of the logarithm sheaves and the elliptic...}}
In the whole present section we assume that the relative dimension $g$ of $X/S$ equals $1$, i.e. that $X$ is an elliptic curve over $S$. For this situation we write $E$ instead of $X$.\\
\newline
As in 1.2.2 we denote by $U$ the open complement of the zero section of $E/S$.\\
If $T$ is an open subscheme of $E$ or a closed subscheme of $E$ which is smooth over $\Q$ the notation $\mathcal H_T$ means the $\mathcal O_T$-vector bundle with integrable $\Q$-connection given by the pullback of $\mathcal H$ to $T$.

\subsection{The elliptic polylogarithm}
We construct the elliptic polylogarithm associated to our fixed geometric setting $E/S/\Q$ as an inverse system of de Rham cohomology classes with components in $H^1_{\mathrm{dR}}(U/\Q,\mathcal H_U^\vee \otimes_{\mathcal O_U}\mathcal L_n)$, characterized by a certain natural residue property. It is the de Rham realization of the general formalism provided in 1.3.13 in the original source \cite{Be-Le}; cf. also \cite{Ba-Ko-Ts}, Def. 1.39.\\
\newline
Recall that in $(1.2.2)$ we obtained for each $n\geq 1$ the exact sequence of $\mathcal D_{S/\Q}$-modules
\[\tag{\textbf{1.5.1}} 0 \rightarrow H^1_{\mathrm{dR}}(E/S, \mathcal L_n) \xrightarrow{\mathrm{can}} H^1_{\mathrm{dR}}(U/S, \mathcal L_n) \xrightarrow{\mathrm{Res}^n} \prod_{k=1}^n \mathrm{Sym}^k_{\mathcal O_S} \mathcal H \rightarrow 0.\]
These sequences are compatible with respect to the morphisms induced by the transition maps of the logarithm sheaves. Note that by Lemma 1.1.5 the associated transition morphism
\[\prod_{k=1}^{n+1} \mathrm{Sym}^k_{\mathcal O_S} \mathcal H \rightarrow \prod_{k=1}^n \mathrm{Sym}^k_{\mathcal O_S} \mathcal H\]
is given explicitly as
\[(h_1, h_2,...,h_{n+1}) \mapsto (n \cdot h_1, (n-1)\cdot h_2,...,h_n) \quad \textrm{with}\ h_k \in \mathrm{Sym}^k_{\mathcal O_S}\mathcal H, \ k=1,...,n+1.\]
It follows that we have a well-defined system
\[\tag{\textbf{1.5.2}} \Big(\frac{1}{(n-1)!}\cdot \mathrm{id}_{\mathcal H}\Big)_{n\geq 1} \in \lim_{n\geq 1} \Hom_{\mathcal D_{S/\Q}}\Big(\mathcal H,\prod_{k=1}^n \mathrm{Sym}^k_{\mathcal O_S} \mathcal H \Big).\]

The exact sequence $(1.5.1)$ and Thm. 1.2.1 (ii) readily imply that the maps $\mathrm{Res}^n$ induce an isomorphism of $\Q$-vector spaces
\[\tag{\textbf{1.5.3}} \lim_{n\geq 1} \Hom_{\mathcal D_{S/\Q}}(\mathcal H, H^1_{\mathrm{dR}}(U/S, \mathcal L_n)) \xrightarrow{\sim} \lim_{n\geq 1} \Hom_{\mathcal D_{S/\Q}}\Big(\mathcal H,\prod_{k=1}^n \mathrm{Sym}^k_{\mathcal O_S} \mathcal H \Big).\]
We have Leray spectral sequences (cf. \cite{Kat2}, $(3.3.0)$)
\[E_2^{p,q}=H^p_{\mathrm{dR}}(S/\Q,\mathcal H^\vee \otimes_{\mathcal O_S} H^q_{\mathrm{dR}}(U/S,\mathcal L_n)) \Rightarrow E^{p+q}=H^{p+q}_{\mathrm{dR}}(U/\Q,\mathcal H_U^\vee \otimes_{\mathcal O_U}\mathcal L_n),\]
and one can conclude from Thm. 1.2.9 (ii) that the edge morphisms
\[H^1_{\mathrm{dR}}(U/\Q,\mathcal H_U^\vee \otimes_{\mathcal O_U}\mathcal L_n) \rightarrow H^0_{\mathrm{dR}}(S/\Q,\mathcal H^\vee \otimes_{\mathcal O_S} H^1_{\mathrm{dR}}(U/S,\mathcal L_n))\]
give an isomorphism
\[\tag{\textbf{1.5.4}} \lim_{n\geq 1}H^1_{\mathrm{dR}}(U/\Q,\mathcal H_U^\vee \otimes_{\mathcal O_U}\mathcal L_n) \xrightarrow{\sim} \lim_{n\geq 1} H^0_{\mathrm{dR}}(S/\Q,\mathcal H^\vee \otimes_{\mathcal O_S} H^1_{\mathrm{dR}}(U/S,\mathcal L_n)).\]
Combined with the natural identification
\[H^0_{\mathrm{dR}}(S/\Q,\mathcal H^\vee \otimes_{\mathcal O_S} H^1_{\mathrm{dR}}(U/S,\mathcal L_n)) \simeq \Hom_{\mathcal D_{S/\Q}}(\mathcal H, H^1_{\mathrm{dR}}(U/S,\mathcal L_n))\]
and with $(1.5.3)$ we obtain from $(1.5.4)$ the isomorphism
\[\tag{\textbf{1.5.5}} \lim_{n\geq 1} H^1_{\mathrm{dR}}(U/\Q,\mathcal H_U^\vee \otimes_{\mathcal O_U}\mathcal L_n) \xrightarrow{\sim} \lim_{n \geq 1} \Hom_{\mathcal D_{S/\Q}}\Big(\mathcal H, \prod_{k=1}^n \mathrm{Sym}^k_{\mathcal O_S}\mathcal H\Big).\]
\begin{definition}
(i) The \underline{elliptic polylogarithm cohomology system for $E/S/\Q$} is the system
\[\mathrm{pol}_{\mathrm{dR}}=(\mathrm{pol}_{{\mathrm{dR}}}^n)_{n\geq 1} \in \lim_{n\geq 1} H^1_{\mathrm{dR}}(U/\Q,\mathcal H_U^\vee \otimes_{\mathcal O_U}\mathcal L_n)\]
which under $(1.5.5)$ is mapped to the system $\Big(\frac{1}{(n-1)!}\cdot \mathrm{id}_{\mathcal H}\Big)_{n\geq 1}$ of $(1.5.2)$.\vspace{2.5mm}\\
(ii) The \underline{elliptic polylogarithm extension system for $E/S/\Q$} is the system
\[\mathrm{pol}=(\mathrm{pol}^n)_{n\geq 1}\in \lim_{n\geq 1} \Ext^1_{\mathcal D_{U/\Q}}(\mathcal H_U, \mathcal L_n)\]
corresponding to $\mathrm{pol}_{{\mathrm{dR}}}$ under the standard canonical identification
\[\lim_{n\geq 1} \Ext^1_{\mathcal D_{U/\Q}}(\mathcal H_U, \mathcal L_n) \simeq \lim_{n\geq 1} H^1_{\mathrm{dR}}(U/\Q,\mathcal H_U^\vee \otimes_{\mathcal O_U}\mathcal L_n).\]
\end{definition}

We append a remark which one can use to show that our definition of the elliptic polylogarithm coincides with the one given in \cite{Ba-Ko-Ts}, Def. 1.39, in the case that $S$ is the spectrum of a field.

\begin{remark}
We may define a morphism in $\Mod(\mathcal D_{S/\Q})$
\[\tag{\textbf{1.5.6}} \lim_{n\geq 1}\prod_{k=1}^n \mathrm{Sym}^k_{\mathcal O_S} \mathcal H \rightarrow \prod_{k=1}^{\infty} \mathrm{Sym}^k_{\mathcal O_S} \mathcal H\]
by requiring that the composition with the $n$-th projection ($n\geq 1$) is the chain of canonical maps
\[\lim_{n\geq 1}\prod_{k=1}^n \mathrm{Sym}^k_{\mathcal O_S} \mathcal H \rightarrow \prod_{k=1}^n \mathrm{Sym}^k_{\mathcal O_S} \mathcal H \rightarrow \mathrm{Sym}^n_{\mathcal O_S}\mathcal H.\]
By the already mentioned formula for the transition morphisms $\prod_{k=1}^{n+1} \mathrm{Sym}^k_{\mathcal O_S} \mathcal H \rightarrow \prod_{k=1}^n \mathrm{Sym}^k_{\mathcal O_S} \mathcal H$ it becomes clear that $(1.5.6)$ is an isomorphism. Then, under the composition
\[\begin{split} \lim_{n\geq 1} \Hom_{\mathcal D_{S/\Q}}\Big(\mathcal H,\prod_{k=1}^n \mathrm{Sym}^k_{\mathcal O_S} \mathcal H\Big) &\simeq \Hom_{\mathcal D_{S/\Q}}\Big(\mathcal H, \lim_{n\geq 1}\prod_{k=1}^n \mathrm{Sym}^k_{\mathcal O_S} \mathcal H\Big)\\
& \simeq \Hom_{\mathcal D_{S/\Q}}\Big(\mathcal H,\prod_{k=1}^{\infty} \mathrm{Sym}^k_{\mathcal O_S} \mathcal H\Big),\end{split}\]
the system $\Big(\frac{1}{(n-1)!}\cdot \mathrm{id}_{\mathcal H}\Big)_{n\geq 1}$ of $(1.5.2)$ maps to $\mathrm{id}_{\mathcal H} \in \Hom_{\mathcal D_{S/\Q}}(\mathcal H,\prod_{k=1}^{\infty} \mathrm{Sym}^k_{\mathcal O_S} \mathcal H)$.
\end{remark}

\subsection{The $D$-variant of the elliptic polylogarithm}
The elliptic polylogarithm $(\mathrm{pol}^n_{{\mathrm{dR}}})_{n\geq 1}$ was defined as a certain system of de Rham cohomology classes on $U=E - \epsilon(S)$ with coefficients in $\mathcal H_U^\vee \otimes_{\mathcal O_U}\mathcal L_n$.\\
Our main interest in the future, however, will focus on a variant of this object which we introduce now. For its construction we need to remove some more points of the curve, namely the $D$-torsion subscheme $E[D]$ for a fixed integer $D>1$, and then apply a formally analogous procedure as in 1.5.1. The outcome is an inverse system of de Rham cohomology classes on $E - E[D]$ with coefficients in $\mathcal L_n$ which again is determined by a prescribed residue along $E[D]$. In the formalism of \cite{Be-Le} it would be obtained essentially by pulling back the extension considered in ibid., 1.3.12, along the morphism induced by the section $D^2\cdot 1_{ \{\epsilon \}} - 1_{E[D]}$ to be defined below.\\
This "$D$-variant" of the elliptic polylogarithm turns out to provide a much better access for concrete description and computations than the original object. We remark that in explicit terms it was introduced and studied (for the $\ell$-adic realization) in the recent work \cite{Ki3}, where it serves as a crucial tool to derive the relation between the elliptic Soulé elements and $\ell$-adic Eisenstein classes (cf. ibid., 4).

\subsubsection{Some additional notation}

We fix an integer $D>1$ and write $[D]: E \rightarrow E$ for the $D$-multiplication endomorphism of $E$.\\
By the cartesian diagram
\begin{equation*}
\begin{xy}
\xymatrix{
E[D] \ar[r]^{ \ \ \pi_{E[D]}}\ar[d]_{i_D} & S \ar[d]^{\epsilon} \\
E \ar[r]^{[D]} & E}
\end{xy}
\end{equation*}
we define the closed subgroup scheme $E[D]$ of $E$ whose $T$-rational points (for an $S$-scheme $T$) are given by $\ker(E(T) \xrightarrow{\cdot D} E(T))$. The corresponding well-known properties of the morphism $[D]$ imply that $E[D]$ is a finite locally free $S$-group scheme of rank $D^2$, and as we are in characteristic zero it is moreover étale over $S$ (cf. also \cite{Kat-Maz}, Thm. 2.3.1).\\
We let $j_D: U_D \rightarrow E$ be the open immersion of $U_D:=E-E[D]$ and set $\pi_{U_D}:=\pi \circ j_D$.

\begin{equation*}
\begin{xy}
\xymatrix@C-0.3cm{
& E[D] \ar[dr]_{\pi_{E[D]}} \ar[r]^{i_D} & E  \ar[d]_{\pi} & U_D \ar[dl]^{\pi_{U_D}} \ar[l]_{j_D}\\
&  & S \ar[d] & &\\
&  & \Spec(\Q) & &}
\end{xy}
\end{equation*}

\subsubsection{Construction of the $D$-variant of the elliptic polylogarithm}

We start by using the machinery of the localization sequence analogously as at the beginning of 1.2.2.\\
\newline
Namely, for each $n \geq 0$ we have the canonical distinguished triangle in $D^b_{\mathrm{qc}}(\mathcal D_{E/\Q})$ (cf. $(0.2.5)$)
\[(i_D)_+ (i_D)^* \mathcal L_n[-1] \rightarrow \mathcal L_n \rightarrow (j_D)_+ \mathcal L_{n|U_D}\]
from which we obtain by applying the functor $\pi_+$ the distinguished triangle in $D^b_{\mathrm{qc}}(\mathcal D_{S/\Q})$
\[(\pi_{E[D]})_+i_D^*\mathcal L_n [-1] \rightarrow \pi_+ \mathcal L_n \rightarrow (\pi_{U_D})_+ \mathcal L_{n|U_D}.\]
Taking cohomology and using \cite{Dim-Ma-Sa-Sai}, Prop. 1.4, we obtain the exact sequence of vector bundles on $S$ with integrable $\Q$-connection
\[\textbf{(1.5.7)} \quad 0 \rightarrow H^1_{\mathrm{dR}}(E/S, \mathcal L_n) \xrightarrow{\mathrm{can}} H^1_{\mathrm{dR}}(U_D/S, \mathcal L_n) \xrightarrow{\mathrm{Res}^D_n} (\pi_{E[D]})_* i_D^*\mathcal L_n \xrightarrow{\sigma_n^D} H^2_{\mathrm{dR}}(E/S, \mathcal L_n) \rightarrow 0,\]
where all sheaves are equipped with their Gauß-Manin connection relative $\Spec(\Q)$. For $(\pi_{E[D]})_* i_D^*\mathcal L_n$ this means nothing else than applying $(\pi_{E[D]})_*$ and the projection formula to the (pullback) connection
\[i_D^*\mathcal L_n \rightarrow \Omega^1_{E[D]/\Q} \otimes_{\mathcal O_{E[D]}} i_D^*\mathcal L_n \simeq \pi_{E[D]}^*\Omega^1_{S/\Q} \otimes_{\mathcal O_{E[D]}} i_D^*\mathcal L_n.\]
For the outmost right zero in $(1.5.7)$ one observes $H^2_{\mathrm{dR}}(U_D/S,\mathcal L_n)=0$: use the spectral sequence $(0.2.2)$ and that $\pi_{U_D}: U_D \rightarrow S$ is affine and of relative dimension one.\\
\newline
We continue with constructions analogous to those which in 1.5.1 led to $(1.5.5)$.\\
\newline
Namely, from $(1.5.7)$ and Thm. 1.2.1 (ii) one derives that the maps $\mathrm{Res}^D_n$ induce an isomorphism
\[\tag{\textbf{1.5.8}} \lim_{n\geq 0}\Hom_{\mathcal D_{S/\Q}}(\mathcal O_S, H^1_{\mathrm{dR}}(U_D/S, \mathcal L_n)) \xrightarrow{\sim} \lim_{n\geq 0} \Hom_{\mathcal D_{S/\Q}}(\mathcal O_S, \ker(\sigma^D_n)).\]
Prolonging $(1.5.7)$ to the left one sees that for all $n\geq 0$ the canonical maps $H^0_{\mathrm{dR}}(E/S,\mathcal L_n) \rightarrow H^0_{\mathrm{dR}}(U_D/S,\mathcal L_n)$ are isomorphisms, which by application of Thm. 1.2.1 (ii) implies that the transition morphisms $H^0_{\mathrm{dR}}(U_D/S,\mathcal L_{n+1}) \rightarrow H^0_{\mathrm{dR}}(U_D/S,\mathcal L_n)$ are zero. From this one readily deduces that the edge morphisms in the Leray spectral sequences (cf. \cite{Kat2}, $(3.3.0)$)
\[E_2^{p,q}=H^p_{\mathrm{dR}}(S/\Q, H^q_{\mathrm{dR}}(U_D/S,\mathcal L_n)) \Rightarrow E^{p+q}=H^{p+q}_{\mathrm{dR}}(U_D/\Q,\mathcal L_n)\]
yield an isomorphism
\[\tag{\textbf{1.5.9}} \lim_{n\geq 0} H^1_{\mathrm{dR}}(U_D/\Q,\mathcal L_n) \xrightarrow{\sim} \lim_{n\geq 0} H^0_{\mathrm{dR}}(S/\Q, H^1_{\mathrm{dR}}(U_D/S,\mathcal L_n)).\]
Combined with the natural identification
\[H^0_{\mathrm{dR}}(S/\Q, H^1_{\mathrm{dR}}(U_D/S,\mathcal L_n)) \simeq \Hom_{\mathcal D_{S/\Q}}(\mathcal O_S, H^1_{\mathrm{dR}}(U_D/S,\mathcal L_n))\]
and with $(1.5.8)$ we obtain from $(1.5.9)$ the isomorphism
\[\tag{\textbf{1.5.10}} \lim_{n\geq 0} H^1_{\mathrm{dR}}(U_D/\Q,\mathcal L_n) \xrightarrow{\sim} \lim_{n\geq 0} \Hom_{\mathcal D_{S/\Q}}(\mathcal O_S, \ker(\sigma^D_n)).\]
We continue by defining a certain distinguished system in $\displaystyle\lim_{n\geq 0} \Hom_{\mathcal D_{S/\Q}}(\mathcal O_S, \ker(\sigma^D_n))$.\\
\newline
First, note that an element of $\Hom_{\mathcal D_{S/\Q}}(\mathcal O_S, \ker(\sigma^D_n))$ tantamounts to a section in $H^0(E[D],i_D^*\mathcal L_n)$ which is horizontal for the $\Q$-connection on $i_D^*\mathcal L_n$ and which goes to zero under the map $\sigma_n^D$ of $(1.5.7)$ in global $S$-sections:
\[\sigma_n^D(S):H^0(E[D],i_D^*\mathcal L_n) \rightarrow H^0(S,H^2_{\mathrm{dR}}(E/S,\mathcal L_n)).\]
We have an injection
\[\tag{\textbf{1.5.11}} H^0(E[D],\mathcal O_{E_{[D]}}) \hookrightarrow H^0(E[D],i_D^*\mathcal L_n),\]
coming about by taking global $E[D]$-sections in the chain
\[\tag{\textbf{1.5.12}} \mathcal O_{E_{[D]}} \hookrightarrow (\pi_{E[D]})^* \prod_{k=0}^n \mathrm{Sym}^k_{\mathcal O_S}\mathcal H \simeq (\pi_{E[D]})^*\epsilon^*\mathcal L_n \simeq i_D^*[D]^*\mathcal L_n \simeq i_D^*\mathcal L_n,\]
where the monomorphism is $\frac{1}{n!}$-times the obvious inclusion, the first isomorphism is the splitting of $\mathcal L_n$, the second uses the diagram defining $E[D]$ and the last comes from Cor. 1.4.4.\\
Note that we need the normalization by $\frac{1}{n!}$ to obtain from the injections of $(1.5.11)$ an induced map
\[\tag{\textbf{1.5.13}} H^0(E[D],\mathcal O_{E_{[D]}}) \hookrightarrow \lim_{n\geq 0}H^0(E[D],i_D^*\mathcal L_n)\]
into the inverse limit, as one can check with Lemma 1.1.5.\\
\newline
The zero section $\epsilon$ of $E/S$ induces the zero section of $E[D]/S$ which we likewise denote by $\epsilon$.\\
As $E[D]$ is étale and separated over $S$ it follows that $\epsilon: S\rightarrow E[D]$ is an open and closed immersion (cf. \cite{EGAIV}, Ch. IV, Cor. $(17.9.3)$). One obtains the schematic decomposition
\[E[D]=(E[D]\ - \{\epsilon\}) \amalg \{\epsilon\}.\]
We may hence define a section $1_{ \{\epsilon \}} \in H^0(E[D],\mathcal O_{E[D]})$ by determining it to be $0$ on $E[D]\ - \{\epsilon\}$ and to be $1$ on $\{\epsilon\}$; we further write $1_{E[D]}$ for the section $1\in H^0(E[D],\mathcal O_{E[D]})$.\\
We thus have a well-defined section
\[D^2\cdot 1_{ \{\epsilon \}} - 1_{E[D]} \in H^0(E[D],\mathcal O_{E[D]}).\]
Its image under $(1.5.11)$ gives a horizontal section in $H^0(E[D], i_D^*\mathcal L_n)$.\\
\newline
Let us now make the following assumption which will be proven in 1.5.4:
\begin{lemma}
Under the composition of $(1.5.11)$ with the map $\sigma_n^D(S)$ the section $D^2\cdot 1_{ \{\epsilon \}} - 1_{E[D]} \in H^0(E[D],\mathcal O_{E[D]})$ goes to zero in $H^0(S,H^2_{\mathrm{dR}}(E/S, \mathcal L_n))$.
\end{lemma}

Then, by what we said about how elements of $\Hom_{\mathcal D_{S/\Q}}(\mathcal O_S, \ker(\sigma^D_n))$ look like, we conclude:
\begin{lemma}
The image of $D^2\cdot 1_{ \{\epsilon \}} - 1_{E[D]}$ under the injection
\[ H^0(E[D],\mathcal O_{E_{[D]}}) \hookrightarrow \lim_{n\geq 0}H^0(E[D],i_D^*\mathcal L_n)\]
of $(1.5.13)$ is already contained in $\displaystyle\lim_{n\geq 0} \Hom_{\mathcal D_{S/\Q}}(\mathcal O_S, \ker(\sigma^D_n))$. \qquad \qed
\end{lemma}
\begin{definition}
(i) The \underline{$D$-variant of the elliptic polylogarithm cohomology system for $E/S/\Q$} is the system
\[\var = \Big(\varn \Big)_{n\geq 0} \in \lim_{n\geq 0} H^1_{\mathrm{dR}}(U_D/\Q, \mathcal L_n)\]
which under the isomorphism of $(1.5.10)$
\[\lim_{n\geq 0} H^1_{\mathrm{dR}}(U_D/\Q,\mathcal L_n) \xrightarrow{\sim} \lim_{n\geq 0} \Hom_{\mathcal D_{S/\Q}}(\mathcal O_S, \ker(\sigma^D_n))\]
is mapped to the image of $D^2\cdot 1_{ \{\epsilon \}} - 1_{E[D]}$ under the inclusion of $(1.5.13)$
\[H^0(E[D],\mathcal O_{E_{[D]}}) \hookrightarrow \lim_{n\geq 0}H^0(E[D],i_D^*\mathcal L_n).\]
Note that by Lemma 1.5.4 this image is contained in $\displaystyle\lim_{n\geq 0} \Hom_{\mathcal D_{S/\Q}}(\mathcal O_S, ker(\sigma^D_n))$.\vspace{2.5mm}\\
(ii) The \underline{$D$-variant of the elliptic polylogarithm extension system for $E/S/\Q$} is the system
\[\mathrm{pol}_{D^2\cdot 1_{ \{\epsilon \}} - 1_{E[D]}} = \Big(\mathrm{pol}^n_{D^2\cdot 1_{ \{\epsilon \}} - 1_{E[D]}} \Big)_{n\geq 0} \in \lim_{n\geq 0} \Ext^1_{\mathcal D_{U_D/\Q}}(\mathcal O_{U_D}, \mathcal L_n)\]
corresponding to $\var$ under the standard canonical identification
\[\lim_{n\geq 0}\Ext^1_{\mathcal D_{U_D/\Q}}(\mathcal O_{U_D}, \mathcal L_n) \simeq \lim_{n\geq 0} H^1_{\mathrm{dR}}(U_D/\Q, \mathcal L_n).\]
\end{definition}

\subsection{The relation between the elliptic polylogarithm and its $D$-variant}
We connect the elliptic polylogarithm classes $\mathrm{pol}_{{\mathrm{dR}}}^n \in H^1_{\mathrm{dR}}(U/\Q,\mathcal H^\vee_U \otimes_{\mathcal O_U}\mathcal L_n)$ with the $D$-variant classes $\varn \in H^1_{\mathrm{dR}}(U_D/\Q, \mathcal L_n)$.\\
To formulate their relation we apply an elementary change of coefficients
\[\mathrm{mult}_n: H^1_{\mathrm{dR}}(U_D/\Q,\mathcal L_n) \rightarrow H^1_{\mathrm{dR}}(U_D/\Q,\mathcal H^\vee_{U_D} \otimes_{\mathcal O_{U_D}}\mathcal L_n)\]
to $\varn$ whose result then has a direct expression in terms of $\mathrm{pol}_{\mathrm{dR}}^n$.\\
\newline
All of this is a direct translation of \cite{Ki3}, Prop. 4.3.3, to our setting, only adding a bit more details.
 
\subsubsection{A change of coefficients}

For each $n\geq 1$ we define
\[\mathrm{mult}_n: \mathcal L_n \rightarrow \mathcal H_E^\vee \otimes_{\mathcal O_E} \mathcal L_n\]
to be the composition
\[ \mathcal L_n \rightarrow \mathcal H_E^\vee \otimes_{\mathcal O_E} \mathcal H_E \otimes_{\mathcal O_E} \mathcal L_n \rightarrow \mathcal H_E^\vee \otimes_{\mathcal O_E} \mathcal L_1 \otimes_{\mathcal O_E} \mathcal L_n \rightarrow \mathcal H_E^\vee \otimes_{\mathcal O_E} \mathcal L_1 \otimes_{\mathcal O_E} \mathcal L_{n-1} \rightarrow \mathcal H_E^\vee \otimes_{\mathcal O_E} \mathcal L_n,\]
where the first arrow is given by the standard map, the second by the inclusion $\mathcal H_E \hookrightarrow \mathcal L_1$ (cf. $(1.1.2)$), the third by the transition map $\mathcal L_n \rightarrow \mathcal L_{n-1}$ and the last by multiplication in symmetric powers.\\
The morphisms $\mathrm{mult}_n$ are checked to be compatible with the transition maps of the logarithm sheaves and horizontal for the respective $\Q$-connections on $\mathcal L_n$ and $\mathcal H_E^\vee \otimes_{\mathcal O_E} \mathcal L_n$.\\
We get an induced homomorphism of $\Q$-vector spaces, also denoted by $\mathrm{mult}_n$:
\[\mathrm{mult}_n: H^1_{\mathrm{dR}}(U_D/\Q,\mathcal L_n) \rightarrow H^1_{\mathrm{dR}}(U_D/\Q,\mathcal H^\vee_{U_D} \otimes_{\mathcal O_{U_D}}\mathcal L_n),\]
compatible with respect to the transition maps of the logarithm sheaves.\\
Hence, by applying the maps $\mathrm{mult}_n$ to the $D$-variant of the elliptic polylogarithm, we obtain a system
\[\Big(\mathrm{mult}_n \Big(\varn \Big)\Big)_{n\geq 1} \in \lim_{n\geq 1} H^1_{\mathrm{dR}}(U_D/\Q,\mathcal H^\vee_{U_D} \otimes_{\mathcal O_{U_D}}\mathcal L_n)\]
which we will be able to relate with the elliptic polylogarithm of Def. 1.5.1 (i).

\begin{Remark}
Under the standard identification of de Rham- and $\Ext$-spaces $\mathrm{mult}_n$ writes as homomorphism
\[\mathrm{mult}_n: \Ext^1_{\mathcal D_{U_D/\Q}}(\mathcal O_{U_D}, \mathcal L_n) \rightarrow \Ext^1_{\mathcal D_{U_D/\Q}}(\mathcal H_{U_D},\mathcal L_n).\]
One can check that the preceding map is alternatively obtained as follows: tensor a given extension with $\mathcal H_{U_D}$ and take the pushout of the resulting extension along the map
\[\widetilde{\mathrm{mult}_n}: \mathcal H_{U_D} \otimes_{\mathcal O_{U_D}} \mathcal L_n  \rightarrow \mathcal L_n\]
which is defined as the composition of the (by now obvious) arrows
\[\mathcal H_{U_D} \otimes_{\mathcal O_{U_D}} \mathcal L_n \rightarrow \mathcal L_1 \otimes_{\mathcal O_{U_D}} \mathcal L_n \rightarrow \mathcal L_1 \otimes_{\mathcal O_{U_D}} \mathcal L_{n-1}\rightarrow \mathcal L_n.\]
\end{Remark}

\subsubsection{The comparison result}
Recall from Def. 1.5.1 (i) that the system of the elliptic polylogarithm
\[\mathrm{pol}_{{\mathrm{dR}}}=(\mathrm{pol}_{{\mathrm{dR}}}^n)_{n\geq 1} \in \lim_{n\geq 1} H^1_{\mathrm{dR}}(U/\Q,\mathcal H^\vee_U \otimes_{\mathcal O_U}\mathcal L_n)\]
is defined by mapping to the system $\Big(\frac{1}{(n-1)!}\cdot \mathrm{id}_{\mathcal H}\Big)_{n\geq 1}$ of $(1.5.2)$ under the isomorphism of $(1.5.5)$:
\[\lim_{n \geq 1} H^1_{\mathrm{dR}}(U/\Q,\mathcal H^\vee_U \otimes_{\mathcal O_U}\mathcal L_n) \xrightarrow{\sim} \lim_{n \geq 1} \Hom_{\mathcal D_{S/\Q}}\Big(\mathcal H, \prod_{k=1}^n\mathrm{Sym}^k_{\mathcal O_S}\mathcal H\Big).\]
We have $U_D \subseteq U$, an induced morphism $[D]: U_D \rightarrow U$, and hence for each $n\geq 1$ the classes
\[(\mathrm{pol}_{{\mathrm{dR}}}^n)_{|U_D} \in H^1_{\mathrm{dR}}(U_D/\Q,\mathcal H^\vee_{U_D} \otimes_{\mathcal O_{U_D}}\mathcal L_n),\]
\[[D]^*\mathrm{pol}_{{\mathrm{dR}}}^n \in H^1_{\mathrm{dR}}(U_D/\Q,\mathcal H^\vee_{U_D} \otimes_{\mathcal O_{U_D}}\mathcal L_n),\]
where for the second class we have used the isomorphism $[D]^*\mathcal L_n \simeq \mathcal L_n$ of Cor. 1.4.4 and the canonical identification $[D]^*\mathcal H_U^\vee \simeq \mathcal H^\vee_{U_D}$. Finally, for each $n \geq 1$ we also have the class
\[\mathrm{mult}_n\Big(\varn \Big) \in  H^1_{\mathrm{dR}}(U_D/\Q,\mathcal H^\vee_{U_D} \otimes_{\mathcal O_{U_D}}\mathcal L_n).\]
The relation between the previous three classes is given by

\begin{proposition}
For all $n\geq1$ we have the following equality in $H^1_{\mathrm{dR}}(U_D/\Q,\mathcal H^\vee_{U_D} \otimes_{\mathcal O_{U_D}}\mathcal L_n)$:
\[\mathrm{mult}_n \Big (\varn \Big)=D^2\cdot (\mathrm{pol}_{{\mathrm{dR}}}^n)_{|U_D} -D\cdot [D]^*\mathrm{pol}_{{\mathrm{dR}}}^n.\]
\end{proposition}

\begin{proof}
Note that the isomorphism of $(1.5.10)$ in particular establishes an injection
\[\lim_{n\geq 0} H^1_{\mathrm{dR}}(U_D/\Q,\mathcal L_n) \hookrightarrow \lim_{n\geq 0} \Hom_{\mathcal D_{S/\Q}}(\mathcal O_S, (\pi_{E[D]})_*i_D^*\mathcal L_n).\]
By an analogous procedure\footnote{To be precise, one replaces $\mathcal O_S$ by $\mathcal H$ in $(1.5.8)$ and $\mathcal L_n$ by $\mathcal H^\vee_{U_D}\otimes_{\mathcal O_{U_D}}\mathcal L_n$ in the subsequent Leray spectral sequence.} one obtains an injection
\[\begin{split} \lim_{n\geq 0}H^1_{\mathrm{dR}}(U_D/\Q,\mathcal H^\vee_{U_D} \otimes_{\mathcal O_{U_D}}\mathcal L_n) &\hookrightarrow \lim_{n\geq 0} \Hom_{\mathcal D_{S/\Q}}(\mathcal H,(\pi_{E[D]})_*i_D^*\mathcal L_n)\\
&\simeq \lim_{n\geq 0} \Hom_{\mathcal D_{E[D]/\Q}}(\mathcal H_{E[D]},i_D^*\mathcal L_n), \end{split}\]
where $i_D: E[D] \rightarrow E$ and $\pi_{E[D]}: E[D] \rightarrow S$ are still the morphisms fixed at the beginning of 1.5.2.\\
If we further use the identification
\[\tag{\textbf{1.5.14}} i_D^*\mathcal L_n \simeq \prod_{k=0}^n \mathrm{Sym}^k_{\mathcal O_{E[D]}}\mathcal H_{E[D]},\]
coming about by the chain of (by now obvious) maps
\[i_D^*\mathcal L_n \simeq i_D^*[D]^*\mathcal L_n \simeq \pi_{E[D]}^*\epsilon^*\mathcal L_n \simeq \pi_{E[D]}^*\prod_{k=0}^n\mathrm{Sym}^k_{\mathcal O_S}\mathcal H \simeq \prod_{k=0}^n \mathrm{Sym}^k_{\mathcal O_{E[D]}}\mathcal H_{E[D]},\]
then in sum we obtain an injection
\[\lim_{n \geq 0} H^1_{\mathrm{dR}}(U_D/\Q,\mathcal H^\vee_{U_D} \otimes_{\mathcal O_{U_D}}\mathcal L_n) \hookrightarrow \lim_{n \geq 0} \Hom_{\mathcal D_{E[D]/\Q}}\Big(\mathcal H_{E[D]},\prod_{k=0}^n \mathrm{Sym}^k_{\mathcal O_{E[D]}}\mathcal H_{E[D]}\Big)\]
which we will consider in components $n\geq 1$ only:

\[\tag{\textbf{1.5.15}} \lim_{n \geq 1} H^1_{\mathrm{dR}}(U_D/\Q,\mathcal H^\vee_{U_D} \otimes_{\mathcal O_{U_D}}\mathcal L_n) \hookrightarrow \lim_{n \geq 1} \Hom_{\mathcal D_{E[D]/\Q}}\Big(\mathcal H_{E[D]},\prod_{k=0}^n \mathrm{Sym}^k_{\mathcal O_{E[D]}}\mathcal H_{E[D]}\Big).\]
Note that the transition maps on the right side are given analogously as in Lemma 1.1.5.\\
\newline
We now show that the images of
\[\Big(\mathrm{mult}_n \Big(\varn \Big)\Big)_{n\geq 1}\]
and
\[(D^2\cdot (\mathrm{pol}^n_{{\mathrm{dR}}})_{|U_D}-D \cdot [D]^*\mathrm{pol}^n_{{\mathrm{dR}}})_{n\geq 1}\]
under the injection $(1.5.15)$ are equal, which then proves the theorem.\\
\newline
This follows by explicating the definitions of the three involved systems and all the identifications we have made to define them and the map $(1.5.15)$; it is the variety of these which makes the detailed verification of the desired equality rather lengthy. It seems reasonable if we here only record the main steps and results of this laborious task.\\
\newline
Recalling the decomposition $E[D]=(E[D]\ - \{\epsilon\}) \amalg \{\epsilon\}$ it follows by tracing back the definition of the polylogarithm system and of the map $(1.5.15)$ that the image of $(D^2 \cdot (\mathrm{pol}^n_{{\mathrm{dR}}})_{|U_D})_{n\geq 1}$ under $1.5.15)$ is given by $D\cdot \Big(D^2\cdot \frac{1}{(n-1)!}\cdot \mathrm{id}_{\mathcal H_{\{\epsilon \}}}\Big)_{n\geq 1}$. The additional factor $D$ appears because of the isomorphism $\mathcal L_n \xrightarrow{\sim} [D]^*\mathcal L_n$ involved in $(1.5.15)$ (cf. $(1.5.14)$) and because of Cor. 1.4.4.\\
\newline
In a similar way one checks that the image of $(D\cdot [D]^*\mathrm{pol}^n_{{\mathrm{dR}}})_{n\geq 1}$ is given by $\Big(D\cdot \frac{1}{(n-1)!}\cdot \mathrm{id}_{\mathcal H_{E[D]}}\Big)_{n\geq 1}$. This time no additional factor $D$ comes in because in the definition of the class $D\cdot [D]^*\mathrm{pol}^n_{{\mathrm{dR}}}$ the isomorphism $[D]^*\mathcal L_n \xrightarrow{\sim} \mathcal L_n$ is used and eliminates the inverse isomorphism involved in $(1.5.15)$.\\
\newline
Finally, the system $\Big(\mathrm{mult}_n \Big (\varn \Big)\Big)_{n\geq 1}$ is mapped under $(1.5.15)$ to
\[D\cdot \Big(\frac{1}{(n-1)!}\cdot (D^2\cdot \mathrm{id}_{\mathcal H_{\{\epsilon\}}}-\mathrm{id}_{\mathcal H_{E[D]}})\Big)_{n\geq 1} \in \lim_{n\geq 1}\Hom_{\mathcal D_{E[D]/\Q}}\Big(\mathcal H_{E[D]},\prod_{k=0}^n \mathrm{Sym}^k_{\mathcal O_{E[D]}}\mathcal H_{E[D]}\Big).\]
This can be verified as follows:\\
Observe that (by definition) the component $\var^n$ is mapped under the composition
\[\begin{split} &\lim_{n \geq 1} H^1_{\mathrm{dR}}(U_D/\Q,\mathcal L_n) \xrightarrow{\sim} \lim_{n\geq 1} \Hom_{\mathcal D_{S/\Q}}(\mathcal O_S, \ker(\sigma^D_n)) \\
\hookrightarrow &\lim_{n\geq 1} \Hom_{\mathcal D_{S/\Q}}(\mathcal O_S,(\pi_{E[D]})_*i_D^*\mathcal L_n)\simeq \lim_{n\geq 1} \Hom_{\mathcal D_{E[D]/\Q}}\Big(\mathcal O_{E[D]},\prod_{k=0}^n \mathrm{Sym}^k_{\mathcal O_{E[D]}}\mathcal H_{E[D]}\Big)\end{split}\]
to the section $\frac{1}{n!}\cdot (D^2\cdot 1_{\{\epsilon\}}-1_{{E[D]}})$ of $\prod_{k=0}^n \mathrm{Sym}^k_{\mathcal O_{E[D]}}\mathcal H_{E[D]}$.\\
A careful analysis shows that the previous map and $(1.5.15)$ fit into a commutative diagram
\begin{equation*}
\begin{xy}
\xymatrix@C-0.3cm{
\lim_{n \geq 1} H^1_{\mathrm{dR}}(U_D/\Q,\mathcal L_n) \ar[r] \ar[d]_{(\mathrm{mult}_n)_{n\geq 1}} & \lim_{n\geq 1} \Hom_{\mathcal D_{E[D]/\Q}}(\mathcal O_{E[D]},\prod_{k=0}^n \mathrm{Sym}^k_{\mathcal O_{E[D]}}\mathcal H_{E[D]}) \ar[d]\\
\lim_{n \geq 1} H^1_{\mathrm{dR}}(U_D/\Q,\mathcal H^\vee_{U_D} \otimes_{\mathcal O_{U_D}}\mathcal L_n) \ar[r] & \lim_{n \geq 1} \Hom_{\mathcal D_{E[D]/\Q}}(\mathcal H_{E[D]},\prod_{k=0}^n \mathrm{Sym}^k_{\mathcal O_{E[D]}}\mathcal H_{E[D]})}
\end{xy}
\end{equation*}
where the right vertical arrow is given under the canonical identification
\[\Hom_{\mathcal D_{E[D]/\Q}}\Big(\mathcal H_{E[D]},\prod_{k=0}^n \mathrm{Sym}^k_{\mathcal O_{E[D]}}\mathcal H_{E[D]}\Big) \simeq \Hom_{\mathcal D_{E[D]/\Q}}\Big(\mathcal O_{E[D]},\mathcal H^\vee_{E[D]}\otimes_{\mathcal O_{E[D]}}\prod_{k=0}^n \mathrm{Sym}^k_{\mathcal O_{E[D]}}\mathcal H_{E[D]}\Big)\]
by composition with $D$-times the following chain of maps:
\begin{align*}
&\prod_{k=0}^n \mathrm{Sym}^k_{\mathcal O_{E[D]}}\mathcal H_{E[D]} \rightarrow \mathcal H^\vee_{E[D]}\otimes_{\mathcal O_{E[D]}} \mathcal H_{E[D]} \otimes_{\mathcal O_{E[D]}} \prod_{k=0}^n \mathrm{Sym}^k_{\mathcal O_{E[D]}}\mathcal H_{E[D]}\\
\rightarrow &\mathcal H^\vee_{E[D]}\otimes_{\mathcal O_{E[D]}} \mathcal H_{E[D]} \otimes_{\mathcal O_{E[D]}} \prod_{k=0}^{n-1} \mathrm{Sym}^k_{\mathcal O_{E[D]}}\mathcal H_{E[D]} \rightarrow \mathcal H^\vee_{E[D]}\otimes_{\mathcal O_{E[D]}} \prod_{k=0}^n \mathrm{Sym}^k_{\mathcal O_{E[D]}}\mathcal H_{E[D]}.
\end{align*}
Here, the first arrow is the standard map, the second is given by the transition map of the symmetric power - hence acts on $\mathcal O_{E[D]}$ by $n$-multiplication (cf. Lemma 1.1.5) - and the third is multiplication in symmetric powers. The additional factor $D$ comes in essentially because an inclusion of the form $\mathcal H_E \hookrightarrow \mathcal L_1 \xrightarrow{\sim} [D]^*\mathcal L_1$ appears when one constructs (by using the maps $\mathrm{mult}_n$) the right vertical arrow of the above diagram such that it really commutes; then one applies again Cor. 1.4.4.\\
The image of the section $\frac{1}{n!}\cdot (D^2\cdot 1_{\{\epsilon\}}-1_{{E[D]}})$ under the previous chain multiplied with $D$ is indeed equal to $D\cdot \frac{1}{(n-1)!}\cdot (D^2\cdot \mathrm{id}_{\mathcal H_{\{\epsilon\}}}-\mathrm{id}_{\mathcal H_{E[D]}})$, as one readily checks. This clearly finishes the proof.
\end{proof}

\begin{remark}
In terms of extension systems Prop. 1.5.7 expresses as the following equality in $\Ext^1_{\mathcal D_{U_D/\Q}}(\mathcal H_{U_D},\mathcal L_n)$:
\[\mathrm{mult}_n \Big(\mathrm{pol}^n_{D^2\cdot 1_{ \{\epsilon \}} - 1_{E[D]}} \Big)=D^2\cdot (\mathrm{pol}^n)_{|U_D} -D\cdot [D]^*\mathrm{pol}^n,\]
where $\mathrm{mult}_n$ is the map of Rem. 1.5.6; analogously as before, to define the extension class of the right side one uses the identifications $[D]^*\mathcal H_U \simeq \mathcal H_{U_D}$ and $[D]^*\mathcal L_n \simeq \mathcal L_n$.
\end{remark}

\subsection{Proof of Lemma 1.5.3}
We have to show the following statement:\\
Let $n\geq 0$ and consider the composition
\[\tag{\textbf{1.5.16}} H^0(E[D],\mathcal O_{E_{[D]}}) \hookrightarrow H^0(E[D],i_D^*\mathcal L_n) \xrightarrow{\sigma_n^D(S)} H^0(S,H^2_{\mathrm{dR}}(E/S,\mathcal L_n))\]
of the inclusion $(1.5.11)$ with the morphism $\sigma^D_n$ of $(1.5.7)$ in global $S$-sections. Then the element $D^2\cdot 1_{ \{\epsilon \}} - 1_{E[D]}\in H^0(E[D],\mathcal O_{E[D]})$ is mapped to zero under $(1.5.16)$.\\
\newline
\underline{Proof:}\\
We will show that $D^2\cdot 1_{ \{\epsilon \}} - 1_{E[D]}$ is mapped to zero when we further compose $(1.5.16)$ with the identification $H^0(S,H^2_{\mathrm{dR}}(E/S,\mathcal L_n)) \xrightarrow{\sim} H^0(S,\mathcal O_S)$ coming from Thm. 1.2.1 (i):
\[\tag{\textbf{1.5.17}} H^0(E[D],\mathcal O_{E[D]}) \hookrightarrow H^0(E[D], i_D^*\mathcal L_n) \xrightarrow{\sigma_n^D(S)} H^0(S,H^2_{\mathrm{dR}}(E/S,\mathcal L_n)) \xrightarrow{\sim} H^0(S,\mathcal O_S).\]
\underline{First Step:}\\
We claim that the section $1_{E[D]} \in H^0(E[D],\mathcal O_{E[D]})$ maps to $D^2 \in H^0(S,\mathcal O_S)$ under $(1.5.17)$.\\
\newline
For this we begin by observing that the following diagram commutes:
\begin{equation*} \tag{\textbf{1.5.18}}\begin{split} \ \ 
\begin{xy}
\xymatrix{
H^0(E[D], i_D^*\mathcal L_n) \ar[d]^{\mathrm{can}} \ar[r]^{\sigma_n^D(S) \quad \ }  & H^0(S,H^2_{\mathrm{dR}}(E/S,\mathcal L_n)) \ar[r]^{\qquad \ \ \psi} \ar[d]^{\mathrm{can}} & H^0(S,\mathcal O_S) \ar[d]^{\id}\\
H^0(E[D],\mathcal O_{E[D]}) \ar[r]^{\sigma_0^D(S) \ \ } & H^0(S,H^2_{\mathrm{dR}}(E/S)) \ar[r]^{\qquad \mathrm{tr}(S) \ } \ar[r] & H^0(S,\mathcal O_S)}
\end{xy}
\end{split}
\end{equation*}
Here, $\psi$ is the last arrow of $(1.5.17)$ and $\mathrm{tr}(S)$ denotes the trace isomorphism in global $S$-sections. The morphisms $\mathrm{can}$ are induced by the composition $\mathcal L_n \rightarrow \mathcal O_E$ of transition maps. The two small diagrams commute by the definition of $\psi$ and functoriality of the procedure that led to $(1.5.7)$.\\
We may prolong the previous diagram on the left by the following commutative diagram:
\begin{equation*} \tag{\textbf{1.5.19}}\begin{split}\ \
\begin{xy}
\xymatrix@C-0.3cm{
 & H^0(E[D], i_D^*\mathcal L_n) \ar[dd]^{\mathrm{can}}\\
H^0(E[D],\mathcal O_{E[D]}) \ar[ur] \ar[dr]^{\id} & \\
 &  H^0(E[D],\mathcal O_{E[D]})}
\end{xy}
\end{split}
\end{equation*}
Here, the upper left arrow is the map $(1.5.11)$, and the diagram indeed commutes: when checking this one has to observe the factor $\frac{1}{n!}$ used in the definition of $(1.5.11)$ and that the projection $\prod_{k=0}^n \mathrm{Sym}^k_{\mathcal O_S}\mathcal H \rightarrow \mathcal O_S$ induced by $\mathcal L_n \rightarrow \mathcal O_E$ is given by $n!\cdot \mathrm{id}_{\mathcal O_S}$ according to Lemma 1.1.5.\\
From the commutativity of $(1.5.18)$ and $(1.5.19)$ we see that the claim of the first step reduces to verifying that $1_{E[D]}$ maps to $D^2$ under the composition
\[\tag{\textbf{1.5.20}} H^0(E[D],\mathcal O_{E[D]}) \xrightarrow{\sigma_0^D(S)} H^0(S,H^2_{\mathrm{dR}}(E/S)) \xrightarrow{\mathrm{tr}(S)} H^0(S,\mathcal O_S).\]
It suffices to show that this holds after transition to an étale covering $(S_i \rightarrow S)_i$ of $S$, as one checks without difficulties.\footnote{One uses that the assignment $T\mapsto H^0(T,\mathcal O_T)$ defines a sheaf on the étale site of $S$ (cf. \cite{Mi}, Ch. II, §1, Cor. 1.6) and that the maps in $(1.5.20)$ are compatible with base change over $S$ in the sense that when setting $E_i:=E\times_S S_i$ we have a commutative diagram
\begin{equation*} 
\begin{xy}
\xymatrix{
H^0(E[D],\mathcal O_{E[D]}) \ar[r]^{\sigma_0^D(S) \quad } \ar[d]_{\mathrm{can}} & H^0(S,H^2_{\mathrm{dR}}(E/S)) \ar[d]_{\mathrm{can}} \ar[r]^{ \quad \ \ \mathrm{tr}(S)} & H^0(S,\mathcal O_S)\ar[d]_{\mathrm{can}}\\
H^0(E_i[D],\mathcal O_{E_i[D]}) \ar[r]^{\sigma_0^D(S_i) \ \ } & H^0(S_i,H^2_{\mathrm{dR}}(E_i/S_i)) \ar[r]^{\qquad \mathrm{tr}(S_i) \ } \ar[r] & H^0(S_i,\mathcal O_{S_i})}
\end{xy}
\end{equation*}}
Further, we can choose an étale covering of $S$ with affine schemes $S_i$ and the property that over $S_i$ the divisor $E[D]$ becomes equal to a divisor of the form $[P^i_1]+...+[P^i_{D^2}]$ with sections $P^i_k \in (E\times_S S_i)(S_i)$, as follows from \cite{Kat-Maz}, Thm. 2.3.1, Prop. 1.10.12 and Thm. 1.10.1.\\ Altogether, we may hence assume from the beginning that $S$ is affine and that
\[E[D]=[P_1]+...+[P_{D^2}]\]
with sections $P_k \in E(S)$. Now consider again the composition $(1.5.20)$:\\
The image of $1_{E[D]}$ under $\sigma_0^D(S)$ is the fundamental class of the divisor $E[D]=[P_1]+...+[P_{D^2}]$ (essentially, this is \cite{Co2}, Lemma 1.5.1). As $\mathrm{tr}(S)$ maps this class precisely to $D^2$ (cf. \cite{Kat1}, 7.4\footnote{The projectiveness assumption made there is unnecessary.}) the claim of the first step of proof is shown.\\
\newline
\underline{Second Step:}\\
We claim that the section $1_{\{\epsilon\}}\in H^0(E[D],\mathcal O_{E[D]})$ maps to $1 \in H^0(S,\mathcal O_S)$ under $(1.5.17)$.\\
\newline
As in the first step one is reduced to show that the image of $1_{\{\epsilon\}}$ under $(1.5.20)$ equals $1\in H^0(S,\mathcal O_S)$.\\
To see this one considers the composition
\[\tag{\textbf{1.5.21}} H^0(S,\mathcal O_S) \rightarrow  H^0(E[D],\mathcal O_{E[D]}) \xrightarrow{\sigma_0^D(S)}  H^0(S,H^2_{\mathrm{dR}}(E/S))\]
in which the first arrow is given in the natural way by the decomposition $E[D]=(E[D]\ - \{\epsilon\}) \amalg \{\epsilon\}$, i.e. by extending functions to zero on $E[D]\ - \{\epsilon\}$. It is clear that $(1.5.21)$ is just the map $\sigma_0(S)$ of Lemma 1.2.3 (ii). From Lemma 1.2.5 it then follows that the composition
\[H^0(S,\mathcal O_S) \rightarrow H^0(E[D],\mathcal O_{E[D]}) \xrightarrow{\sigma_0^D(S)} H^0(S,H^2_{\mathrm{dR}}(E/S)) \xrightarrow{\mathrm{tr}(S)} H^0(S,\mathcal O_S)\]
is the identity. As the image of $1 \in H^0(S,\mathcal O_S)$ under the first of these arrows is $1_{\{\epsilon\}}\in H^0(E[D],\mathcal O_{E[D]})$ we deduce the claim of the second step of proof.\\
\newline
Combining the two steps of proof yields the statement of Lemma 1.5.3. \qquad \qed\\

\chapter{The logarithm sheaves and the Poincaré bundle}
We resume the general geometric setting
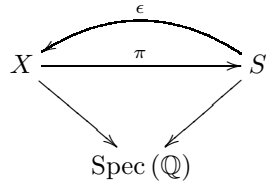
\begin{equation*}
\begin{xy}
\xymatrix@C-0.3cm{
X \ar[rr]^{\pi} \ar[dr]& & S \ar@ /_ 0.6cm/[ll]_{\epsilon}\ar[dl]\\
& \Spec(\Q) &}
\end{xy}
\end{equation*}
fixed at the beginning of the previous chapter, i.e. $X$ is an abelian scheme of relative dimension $g\geq 1$ over a connected scheme $S$ which is assumed smooth, separated and of finite type over $\Spec(\Q)$.\\
All notation introduced so far remains valid for the present chapter.

\section{A preliminary discussion}
\markright{\uppercase{The logarithm sheaves and the Poincaré bundle}}
We provide a clean articulation of the heuristic which will stand as a guiding principle behind all further progress of the present chapter. Pointedly formulated it says that we can fully reconstruct the first logarithm sheaf of $X/S/\Q$ if we know "the first logarithm sheaf of $X/S$". The last might be defined in complete analogy to the first logarithm sheaf of $X/S/\Q$, but by considering $\mathcal D_{X/S}$- resp. $\mathcal O_S$-linear structures instead of $\mathcal D_{X/\Q}$- resp. $\mathcal D_{S/\Q}$-linear structures.\\
In again more immediate terms: when searching for a (conceptual or explicit) description of the logarithm sheaves of $X/S/\Q$ we may at first safely forget that the connections are in fact absolute, i.e. $\Q$-connections, and instead try to find the adapted connections relative $S$: the reason is that these prolong uniquely to the desired $\Q$-connections, which will be the content of the crucial Prop. 2.1.4.\\
Let us finally remark already now that, of course, at a certain future stage it will become inevitable to leave the lines of this policy and to interpret conceptually resp. compute explicitly those completely abstract prolongations. This will happen essentially in 2.6 resp. (in the universal elliptic case) in 3.5.4.

\subsubsection{From absolute to relative connections}

In 1.1 we introduced the first logarithm sheaf of $X/S/\Q$ as a triple $(\mathcal L_1,\nabla_1,\varphi_1)$ consisting of a $\mathcal O_X$-vector bundle $\mathcal L_1$ with integrable $\Q$-connection $\nabla_1$, sitting in an exact sequence of $\mathcal D_{X/\Q}$-modules
\[\tag{\textbf{2.1.1}} 0 \rightarrow \mathcal H_X \rightarrow \mathcal L_1 \rightarrow \mathcal O_X \rightarrow 0\]
which represents the first logarithm extension class, together with the choice of a $\mathcal D_{S/\Q}$-linear splitting
\[\varphi_1: \mathcal O_S \oplus \mathcal H \simeq \epsilon^*\mathcal L_1\]
for the pullback of $(2.1.1)$ along the zero section $\epsilon$.\\
By the first logarithm extension class of $X/S/\Q$ we understood the class in $\Ext^1_{\mathcal D_{X/\Q}}(\mathcal O_X,\mathcal H_X)$ projecting to the identity and retracting to zero in the split exact sequence
\[\tag{\textbf{2.1.2}} 0\rightarrow \Ext^1_{\mathcal D_{S/\Q}}(\mathcal O_S, \mathcal H) \xrightarrow{\pi^*} \Ext^1_{\mathcal D_{X/\Q}}(\mathcal O_X, \mathcal H_X) \rightarrow \Hom_{\mathcal D_{S/\Q}}(\mathcal O_S, \mathcal H^\vee \otimes_{\mathcal O_S} \mathcal H) \rightarrow 0\]
which is part of the five term exact sequence associated with the Leray spectral sequence for $\mathcal H_X$.\\
Throughout, the vector bundle $\mathcal H_X=\pi^*\mathcal H$ was endowed with the pullback of the dual of the Gauß-Manin connection relative $\Spec(\Q)$ and $\mathcal O_X$ with the connection defined by exterior $\Q$-derivation.\\
\newline
We now consider the restriction relative $S$ of these connections such that $\mathcal H_X=\pi^*\mathcal H$ is equipped with its canonical $S$-connection (given by $\mathrm{d}\otimes \id$, cf. 0.2.1 (v)) and $\mathcal O_X$ with exterior $S$-derivation.\\
For $\Omega^{\bullet}_{X/S}(\mathcal H_X)$, the de Rham complex relative $S$ for $\mathcal H_X$, we have a Leray spectral sequence of hypercohomology (use \cite{Dim}, Thm. 1.3.19 (ii)):
\[E^{p,q}_2=H^p(S, H^q_{\mathrm{dR}}(X/S)\otimes_{\mathcal O_S}\mathcal H) \Rightarrow E^{p+q}=\mathbb H^{p+q}(X, \Omega^{\bullet}_{X/S}(\mathcal H_X)) \simeq \Ext^{p+q}_{\mathcal D_{X/S}}(\mathcal O_X,\mathcal H_X).\]
The beginning of its five term exact sequence and $(2.1.2)$ fit into the following commutative diagram with exact rows\footnote{For the surjectivity in the lower line observe that the five term exact sequence continues with the arrow $H^2(S,\mathcal H) \rightarrow \Ext^2_{\mathcal D_{X/S}}(\mathcal O_X,\mathcal H_X)$ which is in fact injective: note that (by compatibility of the Leray spectral sequences for $\Omega^{\bullet}_{X/S}(\mathcal H_X)$ and for $\mathcal H_X$) the composition of this arrow with the canonical morphism $\Ext^2_{\mathcal D_{X/S}}(\mathcal O_X,\mathcal H_X) \xrightarrow{\mathrm{can}} \Ext^2_{\mathcal O_X}(\mathcal O_X,\mathcal H_X) \simeq H^2(X,\mathcal H_X)$ is the map on cohomology induced by $\pi$ which in turn is injective because of the existence of $\epsilon$.} which are both split via the retraction induced by $\epsilon^*$:
\[
{\small
\tag{\textbf{2.1.3}}
\begin{split}
\xymatrix@C-0.3cm{
0 \ar[r] & \Ext^1_{\mathcal D_{S/\Q}}(\mathcal O_S, \mathcal H) \ar[d]_{\mathrm{can}} \ar[r]^-{\pi^*} & \Ext^1_{\mathcal D_{X/\Q}}(\mathcal O_X, \mathcal H_X) \ar[d]_{\mathrm{can}} \ar[r]& \Hom_{\mathcal D_{S/\Q}}(\mathcal O_S, \mathcal H^\vee \otimes_{\mathcal O_S} \mathcal H) \ar[d]_{\mathrm{can}} \ar[r] &0 \\
0 \ar[r] & \Ext^1_{\mathcal O_S}(\mathcal O_S, \mathcal H) \ar[r]^-{\pi^*} & \Ext^1_{\mathcal D_{X/S}}(\mathcal O_X, \mathcal H_X) \ar[r] & \Hom_{\mathcal O_S}(\mathcal O_S, \mathcal H^\vee \otimes_{\mathcal O_S} \mathcal H) \ar[r] & 0}
\end{split}
}
\]
Here, the vertical arrows are the forgetful ones, i.e. given by restricting $\Q$-connections to $S$-connections, the pullbacks via $\pi$ in the lower line are equipped with their canonical $S$-connections, and the projection in this line is given by mapping the class of a $\mathcal D_{X/S}$-linear extension
\[0\rightarrow \mathcal H_X \rightarrow \mathcal M \rightarrow \mathcal O_X \rightarrow 0\]
to the first boundary map $\mathcal O_S \rightarrow \mathcal H^\vee \otimes_{\mathcal O_S}\mathcal H$ in the long exact sequence for the derived functors of $H^0_{\mathrm{dR}}(X/S,-)$ (cf. \cite{Kat2}, $(2.0)$).

\subsubsection{Two auxiliary results}

We will need the following two lemmas.
\begin{lemma}
Suppose we are given two extensions of $\mathcal D_{X/S}$-modules
\begin{align*}
M&: \quad 0\rightarrow \mathcal H_X \xrightarrow{j_M} \mathcal M \xrightarrow{p_M} \mathcal  O_X \rightarrow 0\\
N&: \quad 0\rightarrow \mathcal H_X \xrightarrow{j_N} \mathcal N \xrightarrow{p_N} \mathcal  O_X \rightarrow 0
\end{align*}
with $\mathcal O_S$-linear splittings
\begin{align*}
\varphi_M &: \mathcal O_S \oplus \mathcal H \simeq \epsilon^*\mathcal M\\
\varphi_N &: \mathcal O_S \oplus \mathcal H \simeq \epsilon^*\mathcal N
\end{align*}
and the property that the classes of $M$ and $N$ in $\Ext^1_{\mathcal D_{X/S}}(\mathcal O_X,\mathcal H_X)$ are equal - e.g. if they both map to the identity under the lower projection of $(2.1.3)$.\\
Then there exists a unique isomorphism of $M$ and $N$ which respects the splittings.
\end{lemma}
\begin{proof} \ \\
\underline{Existence}:\\
As the classes of $M$ and $N$ in $\Ext^1_{\mathcal D_{X/S}}(\mathcal O_X,\mathcal H_X)$ are equal there exists a $\mathcal D_{X/S}$-linear isomorphism $\mathcal M \xrightarrow{f} \mathcal N$ which is compatible with the extension structures. For any $\mathcal O_S$-linear morphism $\mathcal O_S \xrightarrow{\mu} \mathcal H$ the map $f+j_N\circ \pi^*(\mu) \circ p_M: \mathcal M \rightarrow \mathcal N$ defines a $\mathcal D_{X/S}$-linear isomorphism compatible with the extension structures. It additionally respects the splittings if we choose $\mu$ in the following way:\\
Define an isomorphism of $\mathcal O_S$-modules $\psi: \mathcal O_S \oplus \mathcal H \rightarrow \mathcal O_S \oplus \mathcal H$ by the commutative diagram
\begin{equation*}
\begin{xy}
\xymatrix@C-0.3cm{
\epsilon^*\mathcal M\ar[r]^{\epsilon^*(f)} \ar[d]_{\varphi_M} & \epsilon^*\mathcal N \ar[d]^{\varphi_N} \\
\mathcal O_S \oplus \mathcal H \ar[r]^{\psi} & \mathcal O_S \oplus \mathcal H}
\end{xy}
\end{equation*}
and set $\mu$ to be the composition
\[\mu: \mathcal O_S \xrightarrow{\mathrm{can}} \mathcal O_S \oplus \mathcal H \xrightarrow{\ \ \psi^{-1}} \mathcal O_S \oplus \mathcal H \xrightarrow{\mathrm{can}} \mathcal H.\]
Then $f+j_N\circ \pi^*(\mu) \circ p_M$ indeed respects the splittings, as one can verify in a calculation.\\
\newline
\underline{Uniqueness}:\\
It clearly suffices to show that the following is true: if $(M, \varphi_M)$ is an extension with splitting as in the claim of the lemma, then any automorphism $h$ of this pair is the identity.\\
Note that by assumption $h: \mathcal M \xrightarrow{\sim} \mathcal M$ is a $\mathcal D_{X/S}$-linear isomorphism such that $\epsilon^*(h)$ is the identity.\\
That then already $h=\id$ holds follows by an argument already used in the proof of Lemma 1.4.1.:\\
Namely, we have to show the vanishing of the $\mathcal D_{X/S}$-linear arrow $h-\id:\mathcal M \rightarrow \mathcal M$ which we know to be zero after pullback via $\epsilon$. By the integrality of $X$ and the fact that $\mathcal M$ is a vector bundle one is reduced to consider from the beginning the situation $S=\Spec(k)$ with $k$ a field of characteristic zero. In this case \cite{Bert-Og}, §2, Prop. 2.16, yields that $h-\id$ is zero not only in the fiber, but already in the stalk of the zero point of $X$. As $X$ is integral and $\mathcal M$ is a vector bundle one can conclude from this that indeed $h-\id$ vanishes.
\end{proof}

\begin{lemma}
Suppose we are given two extensions of $\mathcal D_{X/\Q}$-modules
\begin{align*}
M&: \quad 0\rightarrow \mathcal H_X \xrightarrow{j_M} \mathcal M \xrightarrow{p_M} \mathcal  O_X \rightarrow 0\\
N&: \quad 0\rightarrow \mathcal H_X \xrightarrow{j_N} \mathcal N \xrightarrow{p_N} \mathcal  O_X \rightarrow 0
\end{align*}
with $\mathcal D_{S/\Q}$-linear splittings
\begin{align*}
\varphi_M&: \mathcal O_S \oplus \mathcal H \simeq \epsilon^*\mathcal M\\
\varphi_N&: \mathcal O_S \oplus \mathcal H \simeq \epsilon^*\mathcal N
\end{align*}
and the property that the classes of $M$ and $N$ in $\Ext^1_{\mathcal D_{X/\Q}}(\mathcal O_X,\mathcal H_X)$ are equal - e.g. if they both map to the identity under the upper projection of $(2.1.3)$.\\
Then there exists a unique isomorphism of $M$ and $N$ which respects the splittings.
\end{lemma}
\begin{proof}
The existence is shown completely analogously as in Lemma 2.1.1 by noting that the $\mathcal O_S$-linear map $\mu$ which was chosen there is now $\mathcal D_{S/\Q}$-linear and that hence the map $f+j_N\circ \pi^*(\mu) \circ p_M$ is $\mathcal D_{X/\Q}$-linear. For the uniqueness part it suffices to show that a pair $(M,\varphi_M)$ as in the claim has no nontrivial automorphisms. For this one restricts the $\Q$-connections to $S$-connections and then uses the same proof as in Lemma 2.1.1.
\end{proof}

\begin{remark}
It is easy to also extract from the preceding two proofs that if $(M,\varphi_M)$ is as in the claim of Lemma 2.1.1 resp. Lemma 2.1.2, then the extension $M$ has a nontrivial automorphism if and only if the sheaf $\mathcal H$ has a nonzero global section resp. a nonzero global horizontal section.
\end{remark}

\subsubsection{Reconstruction of the first logarithm sheaf from relative structures}

We can now come to the essential point of this preliminary discussion.

\begin{proposition}
Assume that we are given an exact sequence of $\mathcal D_{X/S}$-modules
\[\tag{\textbf{2.1.4}} 0\rightarrow \mathcal H_X \rightarrow \mathcal L_1' \rightarrow \mathcal O_X \rightarrow 0,\]
whose extension class maps to the identity under the lower projection of $(2.1.3)$, together with a $\mathcal O_S$-linear splitting
\[\tag{\textbf{2.1.5}} \varphi_1': \mathcal O_S \oplus \mathcal H \simeq \epsilon^*\mathcal L_1'\]
for its pullback along the zero section $\epsilon$. Then the following is true:\\
\newline
(i) The integrable $S$-connection on $\mathcal L_1'$ has a unique prolongation to an integrable $\Q$-connection $\nabla_1'$ on $\mathcal L_1'$ such that the following property holds:\\
If we endow $\mathcal L_1'$ with $\nabla_1'$, then the exact sequence in $(2.1.4)$ becomes $\mathcal D_{X/\Q}$-linear and the splitting $\varphi_1'$ in $(2.1.5)$ becomes $\mathcal D_{S/\Q}$-linear.\\
\newline
(ii) The class of the $\mathcal D_{X/\Q}$-linear extension given by $(2.1.4)$ via (i) projects to the identity and retracts to zero in the upper row of $(2.1.3)$.\\
\newline
In other words: $(\mathcal L_1',\nabla_1',\varphi_1')$ is the first logarithm sheaf of $X/S/\Q$.
\end{proposition}
\begin{proof} In view of our assumptions, part (i), and the commutativity of $(2.1.3)$ part (ii) is obvious. It hence remains to show part (i), i.e. the existence and uniqueness of the connection $\nabla_1'$.\\
\newline
\underline{Existence:}\\
Fix a triple $(\mathcal L_1,\nabla_1,\varphi_1)$ as in Def. 1.1.1 and consider all involved $\mathcal D_{X/\Q}$- resp. $\mathcal D_{S/\Q}$-linear structures restricted to $\mathcal D_{X/S}$- resp. $\mathcal O_S$-linear structures. The diagram $(2.1.3)$ and Lemma 2.1.1 then imply that there exists a $\mathcal D_{X/S}$-linear isomorphism
\[\eta: \mathcal L_1' \simeq \mathcal L_1\]
which respects the extension structures belonging to $\mathcal L_1'$ and $\mathcal L_1$ as well as the splittings $\varphi_1'$ and $\varphi_1$. Via $\eta$ and $\nabla_1$ we obtain an integrable $\Q$-connection $\nabla_1'$ on $\mathcal L_1'$ which is easily checked to prolong the integrable $S$-connection on $\mathcal L_1'$ and to satisfy the property formulated in (i).\\
\newline
\underline{Uniqueness:}\\
Assume that $\nabla_1'$ and $\widetilde{\nabla_1'}$ are two integrable $\Q$-connections on $\mathcal L_1'$ which prolong the integrable $S$-connection on $\mathcal L_1'$ and fulfill the property described in (i).\\
Endowing $\mathcal L_1'$ one time with $\nabla_1'$ and the other time with $\widetilde{\nabla_1'}$ we are given in $(2.1.4)$ two $\mathcal D_{X/\Q}$-linear extensions (by assumption). The images of their classes under the upper projection of $(2.1.3)$ are both times the identity: observe the diagram $(2.1.3)$ and that the image of $(2.1.4)$, considered as $\mathcal D_{X/S}$-linear extension, under the lower projection of $(2.1.3)$ is the identity, as was assumed from the beginning on. Furthermore, both of these classes retract to zero in $\Ext^1_{\mathcal D_{S/\Q}}(\mathcal O_S,\mathcal H)$ (by assumption). Altogether, we can thus conclude from the splitting of the upper row of $(2.1.3)$ that the two $\mathcal D_{X/\Q}$-linear extension classes we obtain from $(2.1.4)$ via $\nabla_1'$ and $\widetilde{\nabla_1'}$ are equal in $\Ext^1_{\mathcal D_{X/\Q}}(\mathcal O_X,\mathcal H)$.\\
By Lemma 2.1.2 we conclude that there exists a $\mathcal D_{X/\Q}$-linear isomorphism
\[\nu: (\mathcal L_1', \nabla_1') \simeq (\mathcal L_1',\widetilde{\nabla_1'})\]
which respects the extension structure of $\mathcal L_1'$ and the splitting $\varphi_1'$. Now restrict all involved $\mathcal D_{X/\Q}$-linear structures to $\mathcal D_{X/S}$-linear structures: then, as $\nabla_1'$ and $\widetilde{\nabla_1'}$ are equal when considered as $S$-connections, $\nu$ yields an automorphism of the $\mathcal D_{X/S}$-linear extension $(2.1.4)$ with its $\mathcal O_S$-linear splitting $\varphi_1'$. The uniqueness part of Lemma 2.1.1 then implies that $\nu=\id$. But $\nu$ is $\mathcal D_{X/\Q}$-linear , i.e. horizontal for the $\Q$-connections $\nabla_1'$ and $\widetilde{\nabla_1'}$. This shows $\nabla_1'=\widetilde{\nabla_1'}$.
\end{proof}

\newpage
\markright{\uppercase{The logarithm sheaves and the Poincaré bundle}}

\section{The Fourier-Mukai transformation}
\markright{\uppercase{The logarithm sheaves and the Poincaré bundle}}
Recall from 0.1.1 that we denote by $Y$ the dual abelian scheme of $X$ and by $(\mathcal P^0,r^0,s^0)$ the birigidified Poincaré bundle on $X\times_S Y$; the last means that we keep fixed a representative $(\mathcal P^0,r^0)$ for the universal isomorphism class in $\mathrm{Pic}^0(X\times_S Y/Y)$ and note the existence of a unique $X$-rigidification $s^0$ of $\mathcal P^0$ which is compatible with the $Y$-rigidification $r^0$ of $\mathcal P^0$.\\
\newline
Recall furthermore that $\Yr$ stands for the universal vectorial extension of $Y$ and $(\mathcal P,r,s,\nabla_\mathcal P)$ for the birigidified Poincaré bundle with universal integrable $\Yr$-connection on $X\times_S\Yr$; the last was obtained by first taking the pullback $(\mathcal P,r)$ of $(\mathcal P^0,r^0)$ along the canonical morphism $X\times_S\Yr\rightarrow X\times_S Y$, by then observing that there is a unique integrable $\Yr$-connection $\nabla_\mathcal P$ on $\mathcal P$ such that $(\mathcal P,r,\nabla_\mathcal P)$ represents the universal isomorphism class in $\mathrm{Pic}^\natural(X\times_S \Yr/\Yr)$ and by finally noting the existence of a unique trivialization $s$ for the pullback of $(\mathcal P,\nabla_\mathcal P)$ along the diagram
\begin{equation*}
\begin{xy}
\xymatrix{
X \ar[r]^{\mathrm{id}_X \times \epsilon^\natural \quad \ } \ar[d] & X \times_S \Yr \ar[d]\\
S \ar[r]^{\epsilon^\natural} & \Yr}
\end{xy}
\end{equation*}
which (on the level of modules) is compatible with the $\Yr$-rigidification $r$ of $\mathcal P$; in fact, the isomorphism $s$ is the one induced by $s^0$ in the natural way.\\
\newline
The focus of the present and subsequent sections will lie on the quadruple $(\mathcal P,r,s,\nabla_\mathcal P)$ and the information it contains infinitesimally around the zero section of the $X$-group scheme $X\times_S \Yr$.\\
The reason is that the named information comprises a construction of the logarithm sheaves of $X/S/\Q$ with their connections relative $S$. To ensure that this construction really produces the logarithm sheaves and to reconcile us with the fact that at first we only have access to the relative connections is the role of the preceding Prop. 2.1.4.\\
\newline
In the following, we will denote by $p: X \times_S \Yr \rightarrow X$ and $q:X\times_S \Yr \rightarrow \Yr$ the two projections:
\begin{equation*}
\begin{xy}
\xymatrix{
X \times_S \Yr \ar[r]^{\quad q} \ar[d]_{p} & \Yr \ar[d]^{\pi^\natural}\\
X \ar[r]^{\pi} & S}
\end{xy}
\end{equation*}
\textit{\textbf{The results of the whole present section hold for the general situation of an abelian scheme $X/S$ over an arbitrary noetherian base scheme $S$ of characteristic zero.}}

\subsection{The definition of the Fourier-Mukai transformation}
Our geometric construction of the first logarithm sheaf from the Poincaré bundle, at which we are ultimately aiming, naturally embeds in the formalism of the Fourier-Mukai transformation between sheaves on $\Yr$ and $\mathcal D_{X/S}$-modules on $X$, as introduced by Laumon in \cite{Lau}, 3. This transformation was studied independently by Rothstein (cf. \cite{Ro}) who uses an approach by rather explicit formulas, assuming $S$ as the spectrum of an algebraically closed field. For our purposes, the perspective adopted in \cite{Lau} is more profitable, and we begin by recalling and explicating the definition of the relevant Fourier-Mukai functor given there as well as the result that this functor is in fact an equivalence.\\
Of course, the origin for this circle of ideas stems back to Mukai's seminal work \cite{Mu}.\\
\newline
Recall from 0.2.1 (v) that $\mathcal O$-module pullback via $q$ induces a functor
\[q^*:\Mod(\mathcal O_{\Yr})\rightarrow \Mod(\mathcal D_{X\times_S\Yr/\Yr}),\]
defined by endowing such a pullback with its canonical integrable $\Yr$-connection. The left derivation $Lq^*$ is part of a commutative diagram (with lower horizontal arrow given by the usual left derivation on the level of $\mathcal O$-modules)
\begin{equation*}
\begin{xy}
\xymatrix{
D^-(\mathcal O_{\Yr}) \ar[r]^{Lq^* \quad} \ar[d]_{\id} &  D(\mathcal D_{X\times_S \Yr/\Yr}) \ar[d]^{\mathrm{can}}\\
D^-(\mathcal O_{\Yr}) \ar[r]^{Lq^* \ } & D(\mathcal O_{X\times_S \Yr})}
\end{xy}
\end{equation*}
and yields a triangulated functor
\[\tag{\textbf{2.2.1}} Lq^*: D^b_{\mathrm{qc}}(\mathcal O_{\Yr}) \rightarrow D^b_{\mathrm{qc}}(\mathcal D_{X\times_S \Yr/\Yr}).\]
Note that $Lq^*$ is actually given by termwise pullback of complexes as $q$ is flat and thus $q^*$ is exact.\\
\newline
Moreover, taking tensor product with the Poincaré bundle $\mathcal P$ and its universal integrable $\Yr$-connection gives a functor (cf. 0.2.1 (iv))
\[\mathcal P\otimes_{\mathcal O_{X\times_S \Yr}} (.): \Mod(\mathcal D_{X\times_S\Yr/\Yr}) \rightarrow \Mod(\mathcal D_{X\times_S\Yr/\Yr}).\]
Its left derivation $\mathcal P\otimes_{\mathcal O_{X\times_S \Yr}}^L (.)$ sits again in a commutative diagram (with lower horizontal arrow given by the usual left derivation on the level of $\mathcal O$-modules)
\begin{equation*}
\begin{xy}
\xymatrix{
D^-(\mathcal D_{X\times_S \Yr/\Yr}) \ar[rr]^{\mathcal P\otimes_{\mathcal O_{X\times_S \Yr}}^L (.)} \ar[d]_{\mathrm{can}} & &  D(\mathcal D_{X\times_S \Yr/\Yr}) \ar[d]^{\mathrm{can}}\\
D^-(\mathcal O_{X\times_S \Yr}) \ar[rr]^{\mathcal P\otimes_{\mathcal O_{X\times_S \Yr}}^L (.)} & & D(\mathcal O_{X\times_S \Yr})}
\end{xy}
\end{equation*}
and induces a triangulated functor
\[\tag{\textbf{2.2.2}} \mathcal P\otimes_{\mathcal O_{X\times_S \Yr}}^L (.): D^b_{\mathrm{qc}}(\mathcal D_{X\times_S \Yr/\Yr}) \rightarrow D^b_{\mathrm{qc}}(\mathcal D_{X\times_S \Yr/\Yr}),\]
given in fact by tensoring a complex termwise with $\mathcal P$ because the line bundle $\mathcal P$ is flat over $\mathcal O_{X\times_S \Yr}$ and thus $\mathcal P\otimes_{\mathcal O_{X\times_S \Yr}} (.)$ is exact.\\
\newline
Finally, recall from 0.2.1 (vi) that taking direct image of $\mathcal O$-modules gives rise to a functor
\[p_*:\Mod(\mathcal D_{X\times_S \Yr/\Yr}) \rightarrow \Mod(\mathcal D_{X/S}).\]
Its right derivation $Rp_*$ sits in a commutative diagram (with lower horizontal arrow given by the usual right derivation on the level of $\mathcal O$-modules)
\begin{equation*}
\begin{xy}
\xymatrix{
D^+(\mathcal D_{X\times_S \Yr/\Yr}) \ar[r]^{\quad \ \ Rp_*} \ar[d]_{\mathrm{can}} &  D(\mathcal D_{X/S}) \ar[d]^{\mathrm{can}}\\
D^+(\mathcal O_{X\times_S \Yr}) \ar[r]^{\quad \ Rp_*} & D(\mathcal O_X)}
\end{xy}
\end{equation*}
and is checked to induce a triangulated functor
\[\tag{\textbf{2.2.3}} Rp_*: D^b_{\mathrm{qc}}(\mathcal D_{X\times_S \Yr/\Yr}) \rightarrow D^b_{\mathrm{qc}}(\mathcal D_{X/S}). \vspace{0.1cm}\]
Composition of $(2.2.1)$, $(2.2.2)$ and $(2.2.3)$ yields the triangulated functor
\[\tag{\textbf{2.2.4}} \Phi_\mathcal P: D^b_{\mathrm{qc}}(\mathcal O_{\Yr}) \rightarrow D^b_{\mathrm{qc}}(\mathcal D_{X/S}), \quad \mathcal F^{\bullet} \mapsto Rp_*(\mathcal P\otimes_{\mathcal O_{X\times_S \Yr}}^L Lq^*\mathcal F^{\bullet})\]
which one can view as a "Fourier-Mukai transformation with kernel $(\mathcal P,\nabla_\mathcal P)$".\\
We record the following fundamental insight about the Fourier-Mukai functor $\Phi_\mathcal P$, proven (using different notations) in \cite{Lau}, $(3.2)$, by explicit exhibition of the quasi-inverse:
\begin{theorem}
The triangulated functor $\Phi_\mathcal P: D^b_{\mathrm{qc}}(\mathcal O_{\Yr}) \rightarrow D^b_{\mathrm{qc}}(\mathcal D_{X/S})$ of $(2.2.4)$ is an equivalence of categories.
\end{theorem}

Let us append that for a single sheaf the cohomology of its Fourier-Mukai transformation has the following description, as one sees from the above definitions and from what was said in 0.2.1 (vi):

\begin{remark}
If $\mathcal F$ is a quasi-coherent $\mathcal O_{\Yr}$-module, considered as object of $D^b_{\mathrm{qc}}(\mathcal O_{\Yr})$ in the natural way, the $i$-th cohomology sheaf $H^i(\Phi_\mathcal P(\mathcal F))$ of its Fourier-Mukai transformation is the quasi-coherent $\mathcal O_X$-module with integrable $S$-connection given as follows: Consider $\mathcal P \otimes_{\mathcal O_{X\times_S \Yr}} q^*\mathcal F$ which is endowed with the tensor product of $\nabla_\mathcal P$ with the canonical integrable $\Yr$-connection on $q^*\mathcal F$. Apply the $i$-th higher direct image functor $R^ip_*$ to this connection and identify
\[R^ip_*(\Omega^1_{X\times_S \Yr/\Yr}\otimes_{\mathcal O_{X\times_S\Yr}} (\mathcal P \otimes_{\mathcal O_{X\times_S \Yr}} q^*\mathcal F)) \simeq \Omega^1_{X/S}\otimes_{\mathcal O_X} R^ip_*(\mathcal P \otimes_{\mathcal O_{X\times_S \Yr}} q^*\mathcal F)\]
via the canonical isomorphism $\Omega^1_{X\times_S \Yr/\Yr} \simeq p^*\Omega^1_{X/S}$ together with the projection formula. In this way one obtains a homomorphism of abelian sheaves on $X$:
\[R^ip_*(\mathcal P \otimes_{\mathcal O_{X\times_S \Yr}} q^*\mathcal F) \rightarrow \Omega^1_{X/S}\otimes_{\mathcal O_X} R^ip_*(\mathcal P \otimes_{\mathcal O_{X\times_S \Yr}} q^*\mathcal F)\]
which defines an integrable $S$-connection. The cohomology sheaf $H^i(\Phi_\mathcal P(\mathcal F))$ is then given by the quasi-coherent $\mathcal O_X$-module $R^ip_*(\mathcal P \otimes_{\mathcal O_{X\times_S \Yr}} q^*\mathcal F)$, equipped with this integrable $S$-connection.
\end{remark}

\subsection{WIT-sheaves on the universal vectorial extension}
To make the Fourier-Mukai transformation fertile for our aims we need a way to leave the derived setting and work with honest sheaves. The most convenient way to do this is by force and can already be found in Mukai's classical notion of a WIT-sheaf (WIT = weak index theorem) on an abelian variety (cf. \cite{Mu}, Def. 2.3). We here define its analogue for our given Fourier-Mukai functor on the universal vectorial extension and then present a rather important class of WIT-sheaves of index $0$.
\begin{definition}
A quasi-coherent $\mathcal O_{\Yr}$-module $\mathcal F$ is called \underline{WIT-sheaf of index $i$} if
\[H^j(\Phi_\mathcal P(\mathcal F))= 0 \ \ \textrm{for all} \ \ j \neq i.\]
In this case we write
\[\widehat{\mathcal F}:= H^i(\Phi_\mathcal P(\mathcal F))\]
for the remaining cohomology sheaf in $\Mod_{\mathrm{qc}}(\mathcal D_{X/S})$ and call it the \underline{Fourier-Mukai transform} of $\mathcal F$.
\end{definition}
\begin{remark}
(i) If $\mathcal F$ is a WIT-sheaf of index $i$, then $\Phi_\mathcal P(\mathcal F)\simeq \widehat{\mathcal F}[-i]$ in $D^b_{\mathrm{qc}}(\mathcal D_{X/S})$.\vspace{1mm}\\
(ii) Taking the Fourier-Mukai transform defines a covariant functor from the full subcategory of $\Mod_{\mathrm{qc}}(\mathcal O_{\Yr})$ given by the WIT-sheaves of index $i$ into the category $\Mod_{\mathrm{qc}}(\mathcal D_{X/S})$. It is fully faithful, as one can easily deduce from part (i), Thm. 2.2.1 and the standard fact that $\Mod_{\mathrm{qc}}(\mathcal O_{\Yr})$ resp. $\Mod_{\mathrm{qc}}(\mathcal D_{X/S})$ are fully faithfully embedded in $D^b_{\mathrm{qc}}(\mathcal O_{\Yr})$ resp. $D^b_{\mathrm{qc}}(\mathcal D_{X/S})$.
\end{remark}
As is common for the classical Fourier-Mukai transformation we will frequently need the base change and projection formula. As in our situation integrable connections are involved we need an auxiliary statement ensuring horizontality of these identifications in an adapted sense. The precise form in which this will be used is recorded in the following lemma whose proof we postpone to 2.2.4.
\begin{lemma}
Let
\begin{equation*}
\begin{xy}
\xymatrix{
Z' \ar[r]^{\beta} \ar[d]_{\delta} & Z \ar[d]^{\gamma} \\
T' \ar[r]^{\alpha} & T}
\end{xy}
\end{equation*}
be a cartesian diagram of schemes with a closed immersion $\alpha$ and a smooth morphism $\gamma$.\\
If $\mathcal E\in \Mod_{\mathrm{qc}}(\mathcal D_{Z/T})$, $\mathcal F \in \Mod_{\mathrm{qc}}(\mathcal D_{Z'/T'})$ and $\mathcal G \in \Mod_{\mathrm{qc}}(\mathcal O_{T'})$, then:\\
\newline
(i) The base change isomorphism (cf. \cite{EGAI}, Ch. I, Cor. $(9.3.3)$)
\[\gamma^* \alpha_*\mathcal G \xrightarrow{\sim}\beta_*\delta^*\mathcal G\]
is horizontal if one considers both sides as objects in $\Mod_{\mathrm{qc}}(\mathcal D_{Z/T})$ in the natural way (via 0.2.1 (v) and (vi)).\\
\newline
(ii) The projection formula isomorphism (cf. \cite{EGAI}, Ch. I, Cor. $(9.3.9)$)
\[\mathcal E \otimes_{\mathcal O_Z} \beta_*\mathcal F \xrightarrow{\sim} \beta_*(\beta^*\mathcal E \otimes_{\mathcal O_{Z'}} \mathcal F)\]
is horizontal if one considers both sides as objects in $\Mod_{\mathrm{qc}}(\mathcal D_{Z/T})$ in the natural way (via 0.2.1 (vi), (iv) and (v)).
\end{lemma}

With this observation we can now compute the Fourier-Mukai transforms for an important class of WIT-sheaves of index $0$: these are those $\mathcal O_{\Yr}$-modules that come from quasi-coherent $\mathcal O_S$-modules via the zero section $\epsilon^\natural: S\rightarrow \Yr$.\\
One will note from its proof that the identification recorded in the following proposition crucially uses (apart from standard canonical isomorphisms) the $\mathcal D_{X/S}$-linear trivialization $s$ of $(\mathcal P,\nabla_\mathcal P)$. We nevertheless call it "canonical" because the quadruple $(\mathcal P,r,s,\nabla_\mathcal P)$ is always regarded as fixed.

\begin{proposition}
Let $\mathcal G$ be a quasi-coherent $\mathcal O_S$-module. Then $\epsilon^{\natural}_*\mathcal G$ is a WIT-sheaf of index $0$ for whose Fourier-Mukai transform we have a canonical  identification in $\Mod_{\mathrm{qc}}(\mathcal D_{X/S})$:
\[\widehat{(\epsilon^{\natural}_* \mathcal G)} \simeq \pi^*\mathcal G,\]
where $\pi^*\mathcal G$ is endowed with its canonical integrable $S$-connection.
\end{proposition}
\begin{proof}
Consider the following cartesian diagram, where we abbreviate $\beta:=\id_X \times \epsilon^{\natural}$.
\begin{equation*}
\begin{xy}
\xymatrix@C-0.3cm{
X \simeq X\times_S S\ar[r]^{\quad \beta} \ar[d]_{\pi} & X\times_S \Yr \ar[d]^{q} \\
S \ar[r]^{\epsilon^{\natural}} & \Yr}
\end{xy}
\end{equation*}
By Lemma 2.2.5 (i) we have a canonical isomorphism of $\mathcal D_{X\times_S \Yr/\Yr}$-modules:
\[q^*\epsilon^{\natural}_*\mathcal G \simeq \beta_* \pi^* \mathcal G.\]
Together with Lemma 2.2.5 (ii) and the $\mathcal D_{X/S}$-linear trivialization $s:\beta^* (\mathcal P,\nabla_{\mathcal P}) \simeq (\mathcal O_X,\mathrm{d})$ (belonging to the quadruple $(\mathcal P,r,s,\nabla_\mathcal P)$) we obtain the isomorphism of $\mathcal D_{X\times_S \Yr/\Yr}$-modules:
\[\mathcal P \otimes_{\mathcal O_{X\times_S \Yr}} q^*\epsilon^{\natural}_*\mathcal G \simeq \mathcal P \otimes_{\mathcal O_{X\times_S \Yr}} \beta_*\pi^*\mathcal G \simeq \beta_* \pi^*\mathcal G.\]
Computing the $0$-th cohomology sheaf of $\Phi_\mathcal P(\epsilon_*^\natural \mathcal G)$ then means applying to $\beta_* \pi^*\mathcal G$ the functor\\
$p_*: \Mod_{\mathrm{qc}}(\mathcal D_{X\times_S \Yr/\Yr})\rightarrow \Mod_{\mathrm{qc}}(\mathcal D_{X/S})$, such that in sum we get the $\mathcal D_{X/S}$-linear isomorphism
\[p_*(\mathcal P\otimes_{\mathcal O_{X\times_S \Yr}} q^*\epsilon^{\natural}_*\mathcal G) \simeq (p \circ \beta)_* \pi^*\mathcal G \simeq \pi^*\mathcal G.\]
We finally need to verify
\[H^j(\Phi_\mathcal P(\epsilon^{\natural}_*\mathcal G)) = 0 \ \ \textrm{for all} \ \ j \neq 0,\]
for which we show that the $j$-th higher direct image
\[R^jp_*(\mathcal P\otimes_{\mathcal O_{X\times_S \Yr}} q^*\epsilon^{\natural}_*\mathcal G) \simeq R^jp_*(\beta_*\pi^*\mathcal G)\]
vanishes for all $j\neq 0$. But as a closed immersion $\beta$ is affine, such that the Leray spectral sequence
\[E_2^{j,k}=R^jp_*R^k\beta_*(\pi^*\mathcal G) \Rightarrow E^{j+k}=R^{j+k} (\mathrm{id}_X)_* (\pi^*\mathcal G)\]
implies for each $j$ an isomorphism
\[R^jp_*(\beta_* \pi^*\mathcal G) \simeq R^j(\mathrm{id}_X)_*(\pi^*\mathcal G).\]
This shows the remaining claim.
\end{proof}
We finally append the following easy observation:
\begin{lemma}
If
\[0\rightarrow \mathcal F' \rightarrow \mathcal F \rightarrow \mathcal F{''} \rightarrow 0\]
is an exact sequence in $\Mod_{\mathrm{qc}}(\mathcal O_{\Yr})$ and $\mathcal F',\mathcal F^{''}$ are WIT-sheaves of index $i$, then the same holds for $\mathcal F$ and the sequence of Fourier-Mukai transforms in $\Mod_{\mathrm{qc}}(\mathcal D_{X/S})$
\[0\rightarrow \widehat{\mathcal F'} \rightarrow \widehat{\mathcal F} \rightarrow \widehat{\mathcal F{''}} \rightarrow 0\]
is exact.
\end{lemma}
\begin{proof}
By general theory of derived categories the given exact sequence naturally defines a distinguished triangle
\[\mathcal F' \rightarrow \mathcal F \rightarrow \mathcal F{''} \rightarrow \mathcal F'[1]\]
in $D^b_{\mathrm{qc}}(\mathcal O_{\Yr})$. Applying the triangulated functor $\Phi_\mathcal P$ to it and then taking the long exact sequence of cohomology yields both claims of the lemma.
\end{proof}

\subsection{Categories of unipotent sheaves}
The present geometric situation $X/S$ permits to define the notion of unipotency for vector bundles with integrable $S$-connection on $X$ and to collect such bundles, dependent on their length, in categories $U_n(X/S)$ - analogously as we did in 1.3.1 for integrable $\Q$-connections by considering a setting $X/S/\Q$. Our motivation to study unipotency now in a purely relative situation is twofold:\\
First, the category of unipotent vector bundles (without connections) on an abelian variety already appears in the study of the classical Fourier-Mukai transformation, where one proves its equivalence with the category of coherent modules on the dual variety which are supported in the zero point (cf. \cite{Mu}, Ex. 2.9). It is thus natural and of its own interest to ask for an analogue of this result for the present Fourier-Mukai transformation which involves connections and is defined over a base scheme. We will answer this problem by revealing the category of $\mathcal O_{\Yr}$-modules which corresponds to $U_n(X/S)$ under Fourier-Mukai transformation; the essential ingredients here are the observation of Prop. 2.2.6 and Laumons derived equivalence result of Thm. 2.1.1. We will remark that over a field this category coincides precisely with the category of coherent modules on $\Yr$ which are supported in the zero point, as one would expect from the mentioned result in the classical case.\\
Second, if we want to construct the logarithm sheaves for a setting $X/S/\Q$, then by the discussion in 2.1 it essentially suffices to work with $S$-connections, which means moving in the categories $U_n(X/S)$. But as was just explained we will prove that the objects of these categories are realizable as Fourier-Mukai transforms of sheaves on $\Yr$. Looking at the definition of the Fourier functor we thus get the guarantee that there is a way to obtain the logarithm sheaves from the geometry of the Poincaré bundle on $X\times_S \Yr$. The successive sections 2.3 and 2.4 will then explore this way in detail.\\
A final result of this subsection is concerned with describing for a Fourier-Mukai transform its pullback along the zero section $\epsilon$. By what we just said it is clear that we should address this question because a construction of the first logarithm sheaf always involves the choice of a splitting along $\epsilon$.

\subsubsection{Unipotency and the equivalence result}
Let us denote by $\textit{VIC}(X/S)$ the category whose objects are the vector bundles on $X$ with integrable $S$-connection and whose morphisms are the horizontal $\mathcal O_X$-module homomorphisms. We write $V(S)$ for the category of vector bundles on $S$. By endowing pullbacks via $X \xrightarrow{\pi} S$ with their canonical integrable $S$-connection (cf. 0.2.1 (v)) we obtain an exact functor $\pi^*: V(S) \rightarrow \textit{VIC}(X/S)$.

\begin{definition}
Let $n\geq 0$.\\
(i) An object $\mathcal U$ of $\textit{VIC}(X/S)$ is called \underline{unipotent of length $n$ for $X/S$} if there exists a filtration
\[\mathcal U=A^0\mathcal U \supseteq A^1\mathcal U \supseteq...\supseteq A^n\mathcal U \supseteq A^{n+1}\mathcal U=0\]
by subvector bundles stable under the connection of $\mathcal U$ such that for all $i=0,...,n$ there are objects $\mathcal Y_i$ of $V(S)$ and $\mathcal D_{X/S}$-linear isomorphisms
\[A^i\mathcal U/A^{i+1}\mathcal U \simeq \pi^*\mathcal Y_i.\]
(ii) We write $U_n(X/S)$ for the full subcategory of $\textit{VIC}(X/S)$ consisting of those $\mathcal U$ in $\textit{VIC}(X/S)$ which are unipotent of length $n$ for $X/S$.\\
(iii) We write $U(X/S)$ for the full subcategory of $\textit{VIC}(X/S)$ consisting of those $\mathcal U$ in $\textit{VIC}(X/S)$ which are unipotent of some length for $X/S$. In other words, $U(X/S)$ is the union of the $U_n(X/S)$ for the canonical embeddings
\[U_0(X/S) \hookrightarrow U_1(X/S) \hookrightarrow U_2(X/S) \hookrightarrow... \hookrightarrow \textit{VIC}(X/S).\] 
Note that the zero vector bundle on $X$ with its unique $S$-connection is an object of each $U_n(X/S)$ and that $U_0(X/S)$ is just the essential image of the functor $\pi^*: V(S) \rightarrow \textit{VIC}(X/S)$.
\end{definition}

\begin{remark}
By contrast to the categories $\textit{VIC}(X/\Q)$, $\textit{VIC}(S/\Q)$, $U_n(X/S/\Q)$ and $U(X/S/\Q)$, which were all recognized as abelian (cf. the comment preceding Def. 1.3.1 and the results in 1.3.4) the categories $\textit{VIC}(X/S)$, $V(S)$, $U_n(X/S)$ and $U(X/S)$ are in general not abelian.\\
One can demonstrate this easily by a rather general type of example for which one actually only needs that the morphism $\pi: X\rightarrow S$ is smooth and surjective and that $S$ is an integral scheme whose ring of global sections $\Gamma(S,\mathcal O_S)$ is not a field.\\
Namely, in this situation one first observes (by flatness and surjectivity of $\pi$) that the pullback $\pi^*\mathcal K$ of a coherent $\mathcal O_S$-module $\mathcal K$ which is not a vector bundle on $S$ neither is a vector bundle on $X$.\\
Now choose a non-zero non-unit element $\xi \in \Gamma(S,\mathcal O_S)$ and note (by integrality of $S$) that it induces an exact sequence
\[0 \rightarrow \mathcal O_S \xrightarrow{\cdot \xi} \mathcal O_S \rightarrow \mathcal K \rightarrow 0\]
with a non-zero coherent $\mathcal O_S$-module $\mathcal K$ which then can't be a vector bundle on $S$. Endowing pullbacks with their canonical integrable $S$-connection we obtain an exact $\mathcal D_{X/S}$-linear sequence
\[0 \rightarrow \mathcal O_X \rightarrow \mathcal O_X \rightarrow \pi^*\mathcal K \rightarrow 0,\]
where $\pi^*\mathcal K$ can't be a vector bundle on $X$, as was already remarked.\\
We have thus constructed a morphism in $V(S)$ resp. in $\textit{VIC}(X/S)$, $U_n(X/S)$, $U(X/S)$ whose cokernel in $\Mod(\mathcal O_S)$ resp. in $\Mod(\mathcal D_{X/S})$ does not belong to these categories.
\end{remark}

\begin{definition}
Let $n\geq 0$.\\
An object $\mathcal F$ of $\Mod_{\mathrm{qc}}(\mathcal O_{\Yr})$ is called \underline{unipotent of length n for $\Yr/S$} if there exists a filtration
\[\mathcal F =A^0\mathcal F \supseteq A^1\mathcal F \supseteq...\supseteq A^n\mathcal F \supseteq A^{n+1}\mathcal F=0\]
by quasi-coherent $\mathcal O_{\Yr}$-submodules such that for all $i=0,...,n$ there are objects $\mathcal Y_i$ of $V(S)$ and $\mathcal O_{\Yr}$-linear isomorphisms
\[A^i\mathcal F/A^{i+1}\mathcal F \simeq \epsilon^{\natural}_*\mathcal Y_i.\]
The categories $U_n(\Yr/S)$ and $U(\Yr/S)$ are defined analogously as in Def. 2.2.8 (ii) and (iii).\\
Note that the zero sheaf on $\Yr$ lies in all $U_n(\Yr/S)$ and that $U_0(\Yr/S)$ is the essential image of the functor $\epsilon^\natural_*: V(S) \rightarrow \Mod_{\mathrm{qc}}(\mathcal O_{\Yr})$.
\end{definition}

\begin{remark}
It is easy to see that a $\mathcal O_{\Yr}$-module $\mathcal F$ which belongs to $U(\Yr/S)$ is actually coherent.
\end{remark}

We write $\mathcal J \subseteq \mathcal O_{\Yr}$ for the ideal sheaf defined by the zero section $\epsilon^\natural:S \rightarrow \Yr$.

\begin{theorem}
Let $n \geq 0$. Then:\\
\newline
(i)  If $\mathcal F$ is an object of $U_n(\Yr/S)$ then it is a WIT-sheaf of index $0$ and satisfies $\mathcal J^{n+1}\cdot \mathcal F=0$. The Fourier-Mukai transform $\widehat{\mathcal F}$ is a vector bundle and belongs to $U_n(X/S)$.\\
\newline
(ii) Taking the Fourier-Mukai transform induces equivalences of categories
\begin{align*}
\widehat{(.)}&:U_n(\Yr/S)\xrightarrow{\sim} U_n(X/S),\\
\widehat{(.)}&:U(\Yr/S)\xrightarrow{\sim} U(X/S).
\end{align*}
\end{theorem}
\begin{proof}
Part (i) is easily shown by induction over $n$ and usage of Prop. 2.2.6 and Lemma 2.2.7.\\
As to part (ii) recall from Rem. 2.2.4 (ii) that the Fourier-Mukai transform defines a fully faithful functor from the category of WIT-sheaves of index $0$ into the category $\Mod_{\mathrm{qc}}(\mathcal D_{X/S})$. Together with (i) we can then conclude that the functors in (ii) are well-defined and fully faithful. It remains to show their essential surjectivity, and it is sufficient to do this for the functors $\widehat{(.)}:U_n(\Yr/S) \rightarrow U_n(X/S)$, where $n\geq 0$. The proof proceeds by induction over $n$ as follows:\\
For $n=0$ the claim clearly follows from Prop. 2.2.6.\\
Let now $n\geq 1$ and $\mathcal U$ be an object of $U_n(X/S)$. We find an exact sequence in $\textit{VIC}(X/S)$
\[\tag{\textbf{2.2.5}} 0 \rightarrow A^1\mathcal U \rightarrow \mathcal U \rightarrow \pi^*\mathcal Y\rightarrow 0\]
for some $\mathcal O_S$-vector bundle $\mathcal Y$ and $A^1\mathcal U$ an object of $U_{n-1}(X/S)$. By induction hypothesis there exists $\mathcal F$ in $U_{n-1}(\Yr/S)$ such that $\widehat{\mathcal F} \simeq A^1\mathcal U$, and together with Prop. 2.2.6 we obtain from $(2.2.5)$ an exact sequence in $\textit{VIC}(X/S)$:
\[0 \rightarrow \widehat{\mathcal F} \rightarrow \mathcal U \rightarrow \widehat{(\epsilon^\natural_*\mathcal Y)} \rightarrow 0.\]
This naturally provides a distinguished triangle in $D^b_{\mathrm{qc}}(\mathcal D_{X/S})$:
\[\tag{\textbf{2.2.6}} \widehat{\mathcal F} \rightarrow \mathcal U \rightarrow \widehat{(\epsilon^\natural_*\mathcal Y)} \rightarrow \widehat{\mathcal F} [1].\]
From part (i) and Prop. 2.2.6 we know that $\mathcal F$ and $\epsilon^\natural_*\mathcal Y$ are WIT-sheaves of index $0$, such that $(2.2.6)$ writes in view of Rem. 2.2.4 (i) as a distinguished triangle in $D^b_{\mathrm{qc}}(\mathcal D_{X/S})$:
\[\tag{\textbf{2.2.7}} \Phi_\mathcal P(\mathcal F) \rightarrow \mathcal U \rightarrow \Phi_\mathcal P(\epsilon^\natural_*\mathcal Y) \rightarrow \Phi_\mathcal P(\mathcal F) [1].\]
By Thm. 2.2.1 the functor $\Phi_\mathcal P: D^b_{\mathrm{qc}}(\mathcal O_{\Yr}) \rightarrow D^b_{\mathrm{qc}}(\mathcal D_{X/S})$ has a quasi-inverse $\Phi_\mathcal P^{-1}$. Applied to $(2.2.7)$ it yields a distinguished triangle in $D^b_{\mathrm{qc}}(\mathcal O_{\Yr})$:
\[\mathcal F \rightarrow \Phi_\mathcal P^{-1}(\mathcal U) \rightarrow \epsilon^\natural_*\mathcal Y \rightarrow \mathcal F [1].\]
Going into the long exact sequence of cohomology we obtain that $H^j(\Phi_\mathcal P^{-1}(\mathcal U))=0$ for $j\neq 0$ and an exact sequence in $\Mod_{\mathrm{qc}}(\mathcal O_{\Yr})$:
\[\tag{\textbf{2.2.8}} 0 \rightarrow \mathcal F \rightarrow H^0(\Phi_\mathcal P^{-1}(\mathcal U))\rightarrow \epsilon^\natural_*\mathcal Y \rightarrow 0.\]
As $\Phi_\mathcal P^{-1}(\mathcal U)$ has cohomology only in degree zero we have $H^0(\Phi_\mathcal P^{-1}(\mathcal U)) \simeq \Phi_\mathcal P^{-1}(\mathcal U)$ in $D^b_{\mathrm{qc}}(\mathcal O_{\Yr})$. In particular, we obtain $\Phi_\mathcal P(H^0(\Phi_\mathcal P^{-1}(\mathcal U))) \simeq \mathcal U$ in $D^b_{\mathrm{qc}}(\mathcal D_{X/S})$. If we know that $H^0(\Phi_\mathcal P^{-1}(\mathcal U))$ is an object of $U_n(\Yr/S)$ the preceding observation and Rem. 2.2.4 (i) imply that the Fourier-Mukai transform of $H^0(\Phi_\mathcal P^{-1}(\mathcal U))$ is isomorphic to $\mathcal U$ in $D^b_{\mathrm{qc}}(\mathcal D_{X/S})$ and then also in $\Mod_{\mathrm{qc}}(\mathcal D_{X/S})$ (by general theory of derived categories), which proves the claim. But as $\mathcal F$ is in $U_{n-1}(\Yr/S)$ the exact sequence $(2.2.8)$ clearly shows that $H^0(\Phi_\mathcal P^{-1}(\mathcal U))$ belongs to $U_n(\Yr/S)$.
\end{proof}

\begin{remark}
From part (i) of the previous theorem we see in particular that sheaves of $U(\Yr/S)$ are supported in the zero section of $\Yr$. If $S=\Spec(k)$, with $k$ a field of characteristic zero, we can show that a coherent $\mathcal O_{\Yr}$-module is in $U(\Yr/k)$ already if its support is concentrated in the zero point $e^{\natural}$ of $\Yr$. Hence:
\[U(\Yr/k)=(\textrm{Coherent} \ \mathcal O_{\Yr}\textrm{-modules with support in} \ e^{\natural}) \simeq (\mathcal O_{\Yr,e^{\natural}}\textrm{-modules of finite length}),\]
where the right equivalence is induced by taking the stalk in $e^{\natural}$. Herunder, the category $U_n(\Yr/k)$ corresponds to those $\mathcal O_{\Yr}$- resp. $\mathcal O_{\Yr,e^{\natural}}$-modules which are annihilated by $\mathcal J^{n+1}$ resp. by $\mathcal J^{n+1}_{e^{\natural}}$.\\
By part (ii) of the previous theorem these categories are equivalent to $U(X/k)$ (resp. to $U_n(X/k)$), and one can also check that the length of a $\mathcal O_{\Yr,e^{\natural}}$-module of finite length equals the rank of the corresponding unipotent vector bundle with integrable $k$-connection on $X$.\\
These are analogues of the results for the classical Fourier-Mukai transformation in \cite{Mu}, Ex. 2.9. 
\end{remark}

\subsubsection{The pullback along the zero section}

\begin{definition}
We keep denoting by $\mathcal J \subseteq \mathcal O_{\Yr}$ the (coherent) ideal sheaf of the zero section $\epsilon^\natural: S \rightarrow \Yr$.\\
For each $n\geq 0$ we write $\Yr_n$ for the closed subscheme of $\Yr$ defined by the (coherent) ideal sheaf $\mathcal J^{n+1}\subseteq \mathcal O_{\Yr}$.\\
The following two diagrams (whose hitherto undefined arrows are always the evident ones) introduce some relevant notation associated with $\Yr_n$.
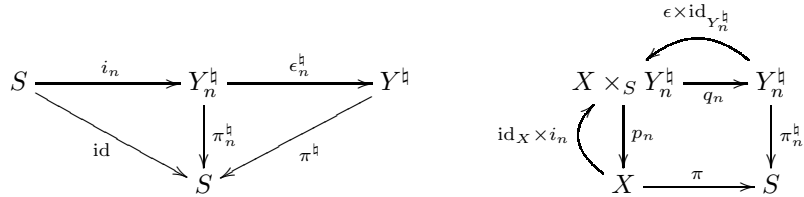
\begin{equation*}
\begin{split}
\begin{xy}
\xymatrix{
S \ar[drr]_{\id} \ar[rr]^{i_n} & & \Yr_n \ar[d]^{\pi_n^\natural} \ar[rr]^{\epsilon^\natural_n}& & \Yr \ar[dll]^{\pi^\natural}  & & X\times_S \Yr_n \ar[r]_{\quad \ q_n} \ar[d]^{p_n}& \Yr_n \ar[d]^{\pi^{\natural}_n} \ar@/_ 0.6cm/ [l]_{\epsilon \times \mathrm{id}_{\Yr_n}}\\
&  & S &  &   & & X \ar@ /^0.6cm/[u]^{\mathrm{id}_X \times i_n} \ar[r]^{\pi} & S}
\end{xy}
\end{split}
\end{equation*}
Finally, we denote by $(\mathcal P_n,\nabla_{\mathcal P_n})$ the pullback (cf. 0.2.1 (v)) of $(\mathcal P,\nabla_\mathcal P)$ along the diagram
\begin{equation*}
\begin{xy}
\xymatrix{
X\times_S \Yr_n \ar[rr]^{\mathrm{id}_X\times \epsilon^\natural_n} \ar[d]_{q_n} & & X\times_S \Yr \ar[d]^{q} \\
\Yr_n \ar[rr]^{\epsilon^\natural_n} & & \Yr}
\end{xy}
\end{equation*}
and by $r_n$ resp. $s_n$ the $\mathcal O_{\Yr_n}$-linear trivialization of $\mathcal P$ along $\epsilon \times \mathrm{id}_{\Yr_n}$ induced by $r$ resp. the $\mathcal D_{X/S}$-linear trivialization of $(\mathcal P_n,\nabla_{\mathcal P_n})$ along $\id_X \times i_n$ induced by $s$.\\
We collect these data in the quadruple $(\mathcal P_n,r_n,s_n,\nabla_{\mathcal P_n})$.
\end{definition}
\begin{remark}
One can check that $\pi^{\natural}_n:\Yr_n \rightarrow S$ is a finite locally free morphism.
\end{remark}
The canonical identification in the following proposition makes (apart from standard canonical isomorphisms) essentially use of the trivialization $r_n$ induced by $r$, as will be seen in its proof. Like in Prop. 2.2.6, where the analogous situation with the trivialization $s$ occurred, we speak of a "canonical" identification because the Poincaré quadruple $(\mathcal P,r,s,\nabla_\mathcal P)$ is always fixed.
\begin{proposition}
For $n\geq 0$ and $\mathcal F$ in $U_n(\Yr/S)$ there is a canonical functorial isomorphism of vector bundles on $S$:
\[\epsilon^*\widehat{\mathcal F} \simeq (\pn)_* \mathcal F,\]
where $\mathcal F$ is naturally considered as $\mathcal O_{\Yr_n}$-module due to $\mathcal J^{n+1} \cdot \mathcal F=0$ (cf. Thm. 2.2.12 (i)).
\end{proposition}
\begin{proof}
Observing the cartesian diagram
\begin{equation*}
\begin{xy}
\xymatrix@C-0.3cm{
X\times_S \Yr_n \ar[rr]^{\mathrm{id}_X\times \epsilon^\natural_n} \ar[d]_{q_n} & & X\times_S \Yr \ar[d]^{q} \\
\Yr_n \ar[rr]^{\epsilon^\natural_n} & & \Yr}
\end{xy}
\end{equation*}
and Lemma 2.2.5 one deduces a canonical $\mathcal D_{X/S}$-linear isomorphism
\[\tag{\textbf{2.2.9}} \widehat{\mathcal F} \simeq (p_n)_*(\mathcal P_n \otimes_{\mathcal O_{X\times_S \Yr_n}} q_n^*\mathcal F),\]
where on the right side the sheaf in brackets is endowed with the tensor product of $\nabla_{\mathcal P_n}$ with the canonical integrable $\Yr_n$-connection on $q^*_n\mathcal F$; taking direct image along the second (cartesian) diagram of Def. 2.2.14 gives the integrable $S$-connection on the right side of $(2.2.9)$ (cf. 0.2.1 (vi)).\\
Moreover, abbreviating $h_n:=\epsilon \times \mathrm{id}_{\Yr_n}$, the cartesian diagram of affine maps
\begin{equation*}
\begin{xy}
\xymatrix{
\Yr_n \ar[r]^{\pi^{\natural}_n} \ar[d]_{h_n} & S \ar[d]^{\epsilon} \\
X\times_S \Yr_n \ar[r]^{\quad \ p_n} & X}
\end{xy}
\end{equation*}
yields the canonical $\mathcal O_S$-linear base change isomorphism (cf. \cite{EGAI}, Ch. I, Cor. $(9.3.3)$):
\[\epsilon^*(p_n)_*(\mathcal P_n \otimes_{\mathcal O_{X\times_S \Yr_n}} q_n^*\mathcal F) \rightarrow (\pi^{\natural}_n)_*h_n^*(\mathcal P_n \otimes_{\mathcal O_{X\times_S \Yr_n}} q_n^*\mathcal F)\]
By combination with $(2.2.9)$ we obtain
\[\epsilon^*\widehat{\mathcal F} \simeq \epsilon^*(p_n)_*(\mathcal P_n \otimes_{\mathcal O_{X\times_S \Yr_n}} q_n^*\mathcal F)\simeq (\pi^{\natural}_n)_*h_n^*(\mathcal P_n \otimes_{\mathcal O_{X\times_S \Yr_n}} q_n^*\mathcal F) \simeq (\pi^{\natural}_n)_*\mathcal F,\]
using in the last step the trivialization $r_n:h_n^*\mathcal P_n \simeq \mathcal O_{\Yr_n}$ induced by the $\Yr$-rigidification $r$ of $\mathcal P$.
\end{proof}

\subsection{Proof of Lemma 2.2.5}
Recall that the lemma in question is occupied with a cartesian diagram of schemes
\begin{equation*}
\begin{xy}
\xymatrix{
Z' \ar[r]^{\beta} \ar[d]_{\delta} & Z \ar[d]^{\gamma} \\
T' \ar[r]^{\alpha} & T}
\end{xy}
\end{equation*}
with a closed immersion $\alpha$, a smooth morphism $\gamma$ and sheaves $\mathcal E\in \Mod_{\mathrm{qc}}(\mathcal D_{Z/T})$, $\mathcal F \in \Mod_{\mathrm{qc}}(\mathcal D_{Z'/T'})$ and $\mathcal G \in \Mod_{\mathrm{qc}}(\mathcal O_{T'})$.\\
\newline
\underline{Proof of (i)}:\\
Let $\varphi$ be the base change isomorphism for $\mathcal G$ on the level of $\mathcal O$-modules. We have to show that the diagram
\begin{equation*}
\begin{xy}
\xymatrix{
\gamma^*\alpha_*\mathcal G \ar[r] \ar[d]_{\varphi} &\Omega^1_{Z/T} \otimes_{\mathcal O_Z} \gamma^*\alpha_*\mathcal G \ar[d]^{\id \otimes \varphi} \\
\beta_*\delta^*\mathcal G \ar[r]&  \Omega^1_{Z/T} \otimes_{\mathcal O_Z} \beta_*\delta^*\mathcal G}
\end{xy}
\end{equation*}
is commutative, where the horizontal arrows are the integrable $T$-connections on $\gamma^*\alpha_*\mathcal G$ resp. $\beta_*\delta^*\mathcal G$.\\
We may assume that $Z$ is affine and thus also that $T$ is affine (by starting with an open affine covering of $T$, taking inverse images in $Z$ and covering these by open affines). But as $\alpha$ and $\beta$ are closed immersions (hence affine maps) we are reduced to the situation where all involved schemes are affine.\\
Our cartesian diagram of schemes then expresses as a diagram
\begin{equation*}
\begin{xy}
\xymatrix{
\Spec(B/\mathfrak a B) \ar[r]^{\beta} \ar[d]_{\delta} & \Spec(B) \ar[d]^{\gamma} \\
\Spec(A/\mathfrak a)\ar[r]^{\alpha}& \Spec(A)}
\end{xy}
\end{equation*}
with rings $A,B$, an ideal $\mathfrak a \subseteq A$ and with $\mathfrak a B$ the ideal generated by $\mathfrak a$ in $B$ via $A\rightarrow B$.\\
The quasi-coherent sheaf $\mathcal G$ corresponds to a $A/\mathfrak a$-module $G$.\\
One checks that $\delta^*\mathcal G$ with its canonical integrable $T'$-connection is given by the $A/\mathfrak a$-linear map
\begin{equation*}
B/\mathfrak a B\otimes _{A/\mathfrak a} G  \rightarrow \Omega^1_{(B/\mathfrak a B)/(A/\mathfrak a)} \otimes_{B/\mathfrak a B} (B/\mathfrak a B \otimes_{A/ \mathfrak a} G), \quad \bar{b}\otimes g \mapsto d\bar{b} \otimes 1 \otimes g.
\end{equation*}
The sheaf $\beta_*\delta^*\mathcal G$ with its integrable $T$-connection then arises from this by using the isomorphism $\Omega^1_{(B/\mathfrak a B)/(A/\mathfrak a)} \simeq \Omega^1_{B/A} \otimes_B B/\mathfrak a B$ and considering $ B/\mathfrak a B \otimes _{A/\mathfrak a} G$ as $B$-module. In brief, one checks that the integrable $T$-connection on $\beta_*\delta^*\mathcal G$ is given by the $A$-linear map
\begin{equation*}
\tag{\textbf{2.2.10}} B/\mathfrak a B \otimes_{A/\mathfrak a} G \rightarrow \Omega^1_{B/A} \otimes_B (B/\mathfrak a B \otimes_{A/\mathfrak a} G), \quad \bar{b} \otimes g \mapsto db\otimes 1 \otimes g,
\end{equation*}
where $B/\mathfrak a B \otimes_{A/\mathfrak a} G$ has the obvious $B$-structure (i.e. via the first factor) and where with $b$ we mean some representative of $\bar{b}$ in $B$. But as $B$-modules we have canonically \[B/\mathfrak a B \otimes_{A/\mathfrak a} G \simeq B \otimes_A G \quad \bar{b} \otimes g \mapsto b \otimes g,\]
which is precisely the base change isomorphism. Under this identification $(2.2.10)$ translates into
\begin{equation*}
B \otimes_A G \rightarrow \Omega^1_{B/A} \otimes_B (B \otimes_A G), \quad b\otimes g \mapsto db \otimes 1\otimes g,
\end{equation*}
which corresponds to $\gamma^*\alpha_*\mathcal G$ with its integrable $T$-connection. This shows (i).\\
\newline
\underline{Proof of (ii)}:\\
We may again assume that all schemes are affine and adapt notations as in the proof of (i).\\
Let $\mathcal F$ correspond to the $B/\mathfrak aB$-module $F$ with integrable $T'$-connection given by the $A/\mathfrak a$-linear map
\begin{equation*}
\nabla_F:F\rightarrow \Omega^1_{(B/\mathfrak a B)/(A/\mathfrak a)} \otimes_{B/\mathfrak a B} F
\end{equation*}
as well as $\mathcal E$ to the $B$-module $E$ with integrable $T$-connection given by the $A$-linear map
\begin{equation*}
\nabla_E:E\rightarrow \Omega^1_{B/A} \otimes_B E,
\end{equation*}
such that the integrable $T$-connection on $\beta_*\mathcal F$ is given by using the map $\nabla_F$ and the $B$-linear isomorphism $\Omega^1_{(B/\mathfrak a B)/(A/\mathfrak a)}\otimes_{B/\mathfrak a B} F \simeq \Omega^1_{B/A} \otimes_B F$.\\
\newline
If we have $e\in E$ and $f\in F$ we write
\[\nabla_E(e)=\sum_jr_j\cdot dt_j\otimes e_j, \quad \nabla_F(f)=\sum_i\bar{u}_i\cdot d\bar{v}_i\otimes f_i,\]
with $r_j, t_j \in B, e_j \in E$ and with $\bar{u}_i,\bar{v}_i \in B/\mathfrak a B, f_i\in F$. Note that we use the description of the differential module $\Omega^1_{B/A}$ as free $B$-module over the symbols $db \ (b\in B)$ modulo the usual equivalence relations, and similar for $\Omega^1_{(B/\mathfrak a B)/(A/\mathfrak a)}$.\\
\newline
Then the tensor product connection on $\mathcal E \otimes_{\mathcal O_Z} \beta_*\mathcal F$ is given by the $A$-linear map
\[\tag{\textbf{2.2.11}} \begin{split} E\otimes_B F &\rightarrow \Omega^1_{B/A} \otimes_B (E\otimes_B F)\\
e\otimes f & \mapsto \sum_j r_j\cdot dt_j \otimes e_j\otimes f + \sum_i u_i \cdot dv_i \otimes e \otimes f_i.\end{split}\]
On the other hand, $\beta^*\mathcal E$ corresponds to the $B/\mathfrak a B$-module $B/\mathfrak a B \otimes_B E$ with integrable $T'$-connection given by a $A/\mathfrak a$-linear map
\begin{equation*}
B/\mathfrak a B\otimes_B E \rightarrow \Omega^1_{(B/\mathfrak a B)/(A/\mathfrak a)}\otimes_{B/\mathfrak a B}(B/\mathfrak a B \otimes_B E).
\end{equation*}
Under the identifications $B/\mathfrak a B \simeq A/\mathfrak a \otimes_A B$ and $\Omega^1_{(B/\mathfrak a B)/(A/\mathfrak a)} \simeq \Omega^1_{B/A} \otimes_B B/\mathfrak a B$ this writes as
\[A/\mathfrak a \otimes_A E \rightarrow \Omega^1_{B/A}\otimes_B (A/\mathfrak a \otimes_A  E), \quad \bar{a} \otimes e \mapsto \sum_j r_j \cdot dt_j \otimes \bar{a} \otimes e_j.\]
Note that the $B/\mathfrak aB$-structure on $A/\mathfrak a \otimes_A E$ is given in the natural way, i.e. by multiplying representatives into the $B$-module $E$.\\
\newline
If we again use the identification $\Omega^1_{(B/\mathfrak a B)/(A/\mathfrak a)} \simeq B/\mathfrak aB \otimes_B \Omega^1_{B/A}$, then it is easily checked that the integrable $T'$-connection on $\beta^*\mathcal E \otimes_{\mathcal O_{Z'}} \mathcal F$ writes as the $A/\mathfrak a$-linear map
\[\tag{\textbf{2.2.12}} \begin{split} (A/\mathfrak a \otimes_A E) \otimes_{B/\mathfrak a B} F &\rightarrow \Omega^1_{B/A}\otimes_B ((A/\mathfrak a \otimes_A E) \otimes_{B/\mathfrak a B} F)),\\
\bar{a} \otimes e \otimes f &\mapsto \sum_j  r_j \cdot dt_j \otimes \bar{a} \otimes e_j \otimes f + \sum_i u_i \cdot dv_i \otimes \bar{a} \otimes e \otimes f_i. \end{split}\]
Now observe that the operation $\beta_*$ applied to $\beta^* \mathcal E \otimes_{\mathcal O_{Z'}}\mathcal F$ just means that we have to regard the modules in $(2.2.12)$ as $B$-modules. But we have the canonical isomorphism of $B$-modules
\[(A/\mathfrak a \otimes_A E) \otimes_{B/\mathfrak a B} F \simeq (A/\mathfrak a \otimes_A B \otimes_B E) \otimes_{(A/\mathfrak a \otimes_A B)} F \simeq E \otimes_B F.\]
As $\bar{1} \otimes e \otimes f$ corresponds to $e\otimes f$ under this identification we see that $(2.2.12)$ becomes
\[\tag{\textbf{2.2.13}} \begin{split} E\otimes_B F &\rightarrow  \Omega^1_{B/A} \otimes_B (E\otimes_B F),\\
e\otimes f &\mapsto \sum_j  r_j\cdot dt_j \otimes e_j \otimes f+ \sum_i u_i \cdot dv_i \otimes e \otimes f_i,\end{split}\]
which by construction is the integrable $T$-connection on $\mathcal E \otimes_{\mathcal O_Z} \beta_*\mathcal F$ obtained from the connection on $\beta_*(\beta^* \mathcal E \otimes_{\mathcal O_{Z'}} \mathcal F)$ via the projection formula. Comparing $(2.2.13)$ and $(2.2.11)$ proves (ii).\\
\qed

\newpage
\markright{\uppercase{The logarithm sheaves and the Poincaré bundle}}
\section{The first logarithm sheaf and the Poincaré bundle}
\markright{\uppercase{The logarithm sheaves and the Poincaré bundle}}

We return to the geometric situation $X/S/\Q$ fixed at the beginning of this chapter and recall the crucial observation made at the end of 2.1: Equip $\mathcal H_X$ with its canonical integrable $S$-connection and $\mathcal O_X$ with exterior $S$-derivation. Assume we have an exact sequence of $\mathcal D_{X/S}$-modules
\[0 \rightarrow \mathcal H_X \rightarrow \mathcal L_1' \rightarrow \mathcal O_X \rightarrow 0,\]
mapping to the identity under the lower projection in $(2.1.3)$:
\[\Ext^1_{\mathcal D_{X/S}}(\mathcal O_X,\mathcal H_X) \rightarrow \Hom_{\mathcal O_S}(\mathcal O_S,\mathcal H^\vee\otimes_{\mathcal O_S}\mathcal H),\]
together with a $\mathcal O_S$-linear splitting
\[\varphi_1':\mathcal O_S \oplus \mathcal H \simeq \epsilon^*\mathcal L_1'\]
for its pullback along the zero section $\epsilon: S \rightarrow X$. Then, by Prop. 2.1.4, the integrable $S$-connection on $\mathcal L_1'$ extends uniquely to an integrable $\Q$-connection such that the previous data become the first logarithm sheaf of $X/S/\Q$.\\
We are now prepared to construct such data from the first infinitesimal restriction $(\mathcal P_1,r_1,s_1,\nabla_{\mathcal P_1})$ of the quadruple $(\mathcal P,r,s,\nabla_\mathcal P)$ - for the definition of this restriction recall Def. 2.2.14.\\
It is also useful to dispose of an equivalent viewpoint provided by the more hermetic machinery of the Fourier-Mukai transformation, with which we begin.
\subsection{The construction of the fundamental data}

\subsubsection{Construction via the Fourier-Mukai formalism}

Consider the canonical exact sequence of $\mathcal O_{\Yr}$-modules:
\[\tag{\textbf{2.3.1}} 0 \rightarrow \mathcal J/\mathcal J^2 \rightarrow \mathcal O_{\Yr}/\mathcal J^2 \rightarrow \mathcal O_{\Yr}/\mathcal J \rightarrow 0.\]
Noting that $\mathcal J/\mathcal J^2$ is the conormal sheaf of the regularly embedded section $\epsilon^\natural: S \rightarrow \Yr$ we have a natural isomorphism of $\mathcal O_{\Yr}$-modules (cf. \cite{Fu-La}, Ch. IV, Lemma 3.8):
\[\mathcal J/\mathcal J^2 \simeq (\epsilon^{\natural})_*(\Lie(\Yr/S)^\vee).\]
Together with the fundamental canonical identification $\Lie(\Yr/S) \simeq H^1_{\mathrm{dR}}(X/S)$ of $(0.1.5)$ we get
\[\tag{\textbf{2.3.2}} \mathcal J/\mathcal J^2 \simeq (\epsilon^{\natural})_* \mathcal H.\]
Moreover, we clearly have
\[\mathcal O_{\Yr}/\mathcal J \simeq (\epsilon^{\natural})_*\mathcal O_S,\]
such that $(2.3.1)$ translates into the $\mathcal O_{\Yr}$-linear exact sequence
\[\tag{\textbf{2.3.3}} 0 \rightarrow (\epsilon^{\natural})_* \mathcal H \rightarrow \mathcal O_{\Yr}/\mathcal J^2 \rightarrow (\epsilon^{\natural})_*\mathcal O_S \rightarrow 0.\]
By Prop. 2.2.6 and Lemma 2.2.7 its terms are WIT-sheaves of index $0$ and the associated exact sequence of Fourier-Mukai transforms writes as a $\mathcal D_{X/S}$-linear sequence
\[\tag{\textbf{2.3.4}} 0 \rightarrow \mathcal H_X \rightarrow \widehat{\mathcal O_{\Yr}/\mathcal J^2} \rightarrow \mathcal O_X \rightarrow 0,\]
where $\mathcal H_X$ is equipped with its canonical integrable $S$-connection and $\mathcal O_X$ with exterior $S$-derivation.\\
Observe furthermore that by Prop. 2.2.16 the pullback of $(2.3.4)$ along $\epsilon: S \rightarrow X$ identifies with the exact sequence $(2.3.3)$, where now all its terms are to be considered as $\mathcal O_S$-modules:
\[\tag{\textbf{2.3.5}} 0 \rightarrow \mathcal H \rightarrow (\pi^\natural_1)_*\mathcal O_{\Yr_1} \rightarrow \mathcal O_S \rightarrow 0.\]
Note that the morphism $\pi^\natural_1: \Yr_1 \rightarrow S$ obviously provides a section for the sequence $(2.3.5)$ and thus a $\mathcal O_S$-linear splitting
\[\tag{\textbf{2.3.6}} \mathcal O_S \oplus \mathcal H \simeq (\pi^\natural_1)_*\mathcal O_{\Yr_1}.\]
Altogether, we obtain a $\mathcal O_S$-linear splitting for the pullback of $(2.3.4)$ along $\epsilon: S\rightarrow X$:
\[\tag{\textbf{2.3.7}} \mathcal O_S \oplus \mathcal H \simeq \epsilon^*(\widehat{\mathcal O_{\Yr}/\mathcal J^2}).\]
Let us now illustrate how the data $(2.3.4)$ and $(2.3.7)$ are explicitly induced by the quadruple $(\mathcal P_1,r_1,s_1,\nabla_{\mathcal P_1})$.

\subsubsection{Construction via the Poincaré bundle}
It is useful to have in sight the following diagram of cartesian squares:
\begin{equation*}
\begin{xy}
\xymatrix{
X \ar[r]^{\pi} \ar[d]_{\mathrm{id}_X \times i_1} & S \ar[d]^{i_1} \\
X \times_S \Yr_1 \ar[d]_{\mathrm{id}_X \times \epsilon^{\natural}_1} \ar[r]^{\quad \ q_1} & \Yr_1 \ar[d]^{\epsilon^{\natural}_1}\\
X \times_S \Yr \ar[r]^{\quad \ q} \ar[d]_{p} & \Yr \ar[d]^{\pi^\natural} \\
X \ar[r]^{\pi} & S
}
\end{xy}
\end{equation*}
Associated with the closed immersion $i_1$ is the canonical exact sequence of $\mathcal O_{\Yr_1}$-modules
\[\tag{\textbf{2.3.8}} 0 \rightarrow (i_1)_*\mathcal H \rightarrow \mathcal O_{\Yr_1} \rightarrow (i_1)_*\mathcal O_S \rightarrow 0,\]
where we have identified the ideal sheaf of $i_1$ with $(i_1)_*\mathcal H$ by $(0.1.5)$.\\
Pullback of $(2.3.8)$ along the (flat) projection $q_1$ gives the exact sequence
\[\tag{\textbf{2.3.9}} 0 \rightarrow (\mathrm{id}_X \times i_1)_*\mathcal H_X \rightarrow \mathcal O_{X\times_S \Yr_1} \rightarrow (\mathrm{id}_X \times i_1)_*\mathcal O_X \rightarrow 0\]
which identifies with the canonical exact sequence associated to the closed immersion $\id_X \times i_1$: note that the ideal sheaf of $\id_X \times i_1$ is the pullback of the ideal sheaf of $i_1$. The sequence $(2.3.9)$ is horizontal for the integrable $\Yr_1$-connections given by the direct image along the upper square of the canonical integrable $S$-connections on $\mathcal H_X$ resp. $\mathcal O_X$ and by exterior $\Yr_1$-derivation on $\mathcal O_{X\times_S \Yr_1}$.\\
Tensoring with $(\mathcal P_1,\nabla_{\mathcal P_1})$ and using Lemma 2.2.5 (ii) together with the $\mathcal D_{X/S}$-linear trivialization
\[s_1: (\mathrm{id}_X \times i_1)^*(\mathcal P_1,\nabla_{\mathcal P_1}) \simeq (\mathcal O_X,\mathrm{d})\]
yields the exact sequence of $\mathcal D_{X\times_S \Yr_1/\Yr_1}$-modules
\[\tag{\textbf{2.3.10}} 0 \rightarrow (\mathrm{id}_X \times i_1)_*\mathcal H_X \rightarrow \mathcal P_1 \rightarrow (\mathrm{id}_X \times i_1)_*\mathcal O_X \rightarrow 0;\]
the surjection comes alternatively from the adjunction $\mathcal P_1 \rightarrow (\mathrm{id}_X \times i_1)_*(\mathrm{id}_X \times i_1)^*\mathcal P_1$ and from $s_1$.\\
Finally, taking direct image of $(2.3.10)$ along the digram
\begin{equation*}
\begin{xy}
\xymatrix{
X\times_S \Yr_1 \ar[r]^{\quad \ q_1} \ar[d]_{p_1} & \Yr_1 \ar[d]^{\pi^\natural_1} \\
X \ar[r]^{\pi} & S}
\end{xy}
\end{equation*}
gives the exact sequence of $\mathcal D_{X/S}$-modules
\[\tag{\textbf{2.3.11}} 0 \rightarrow \mathcal H_X \rightarrow (p_1)_*\mathcal P_1 \rightarrow \mathcal O_X \rightarrow 0.\]
Note that we have a canonical isomorphism between the extensions $(2.3.4)$ and $(2.3.11)$:
\begin{equation*}\tag{\textbf{2.3.12}} \begin{split}
\begin{xy}
\xymatrix@C-0.3cm{
0 \ar[r] & \mathcal H_X \ar[d]_{\id} \ar[r] & \widehat{\mathcal O_{\Yr}/\mathcal J^2}\ar[r]\ar[d]_{\sim} & \mathcal O_X \ar[d]_{\id} \ar[r] & 0 \\
0 \ar[r] & \mathcal H_X \ar[r] & (p_1)_*\mathcal P_1 \ar[r] & \mathcal O_X \ar[r] & 0}
\end{xy}
\end{split}
\end{equation*}
which is induced by the chain of natural identifications
\[\begin{split} \widehat{\mathcal O_{\Yr}/\mathcal J^2}&=p_*(\mathcal P\otimes_{\mathcal O_{X\times_S \Yr}}q^*(\mathcal O_{\Yr}/\mathcal J^2)) \simeq p_*(\mathcal P\otimes_{\mathcal O_{X\times_S \Yr}}q^*(\epsilon^{\natural}_1)_*\mathcal O_{\Yr_1}) \\
&\simeq p_*(\mathcal P\otimes_{\mathcal O_{X\times_S \Yr}}(\mathrm{id}_X \times \epsilon^{\natural}_1)_*\mathcal O_{X\times_S \Yr_1}) \simeq p_*(\mathrm{id}_X\times \epsilon^\natural_1)_*\mathcal P_1 \simeq (p_1)_*\mathcal P_1 \end{split}\]
whose horizontality is guaranteed by Lemma 2.2.5.\\
The pullback of $(2.3.11)$ along $\epsilon:S \rightarrow X$ is $\mathcal O_S$-linearly split by
\[\tag{\textbf{2.3.13}} \epsilon^*(p_1)_*\mathcal P_1 \xrightarrow{\sim} (\pi^\natural_1)_*(\epsilon\times \mathrm{id}_{\Yr_1})^*\mathcal P_1 \simeq (\pi^\natural_1)_*\mathcal O_{\Yr_1} \simeq \mathcal O_S \oplus \mathcal H,\]
where the first map is the base change isomorphism along the cartesian diagram of affine maps
\begin{equation*}
\begin{xy}
\xymatrix{
\Yr_1 \ar[r]^{\pi^{\natural}_1} \ar[d]_{\epsilon \times \mathrm{id}_{\Yr_1}} & S \ar[d]^{\epsilon} \\
X\times_S \Yr_1 \ar[r]^{\quad \ p_1} & X}
\end{xy}
\end{equation*}
the second is induced by the trivialization
\[r_1: (\epsilon\times \mathrm{id}_{\Yr_1})^*\mathcal P_1 \simeq \mathcal O_{\Yr_1}\]
and the third is given by $(2.3.6)$. Under $(2.3.12)$ the splittings $(2.3.7)$ and $(2.3.13)$ correspond.

\subsection{The main result}
With the explanations at the outset of this section our goal must now clearly consist in proving
\begin{theorem}
The class of the extension $(2.3.11)$ maps to the identity under the lower projection in $(2.1.3)$:
\[\tag{\textbf{2.3.14}} \Ext^1_{\mathcal D_{X/S}}(\mathcal O_X,\mathcal H_X) \rightarrow \Hom_{\mathcal O_S}(\mathcal H, \mathcal H).\]
\end{theorem}
\begin{proof}
Let us write $\xi: \mathcal H \rightarrow \mathcal H$ for the image of the class of $(2.3.11)$ under $(2.3.14)$ and recall from $(2.3.12)$ that $\xi$ is also equal to the image of the class of $(2.3.4)$ under $(2.3.14)$.\\
\newline
We now proceed in several steps.\\
\newline
\underline{Step 1}: (Reduction to the case $S=\Spec(k)$)\\
\newline
As $\mathcal H$ is a vector bundle and $S$ is integral the map $\xi$ is the identity on $\mathcal H=H^1_{\mathrm{dR}}(X/S)^\vee$ already if for all points $s \in S$ its pullback along the canonical morphism $f_s: \Spec(k(s)) \rightarrow S$ is the identity on $H^1_{\mathrm{dR}}(X_s/k(s))^\vee$. Here, we set $X_s:=X\times_S \Spec(k(s))$, viewed as abelian variety over $k(s)$, and recall that $H^1_{\mathrm{dR}}(X/S)$ is compatible with arbitrary base change (cf. the beginning of Chapter 1).
\begin{equation*}
\begin{xy}
\xymatrix{
X_s \ar[r]^{\pi_s \qquad} \ar[d]_{g_s} & \Spec(k(s))\ar[d]^{f_s} \\
X \ar[r]^{\pi} & S}
\end{xy}
\end{equation*}
Now observe that clearly everything we said in 2.3.1 applies equally well for the situation of an abelian variety over a field of characteristic zero and that in this situation we also have the map $(2.3.14)$, defined in the same way as done in 2.1 under different assumptions on the base scheme $S$. This is important to note because in the following we will reduce to and work in this situation.\\
Indeed, in our previous consideration of the fibers $X_s$ over points $s\in S$ the respective maps $(2.3.14)$ are checked to fit into a commutative diagram
\begin{equation*}
\begin{xy}
\xymatrix@C-0.3cm{
\Ext^1_{\mathcal D_{X/S}}(\mathcal O_X,\pi^*H^1_{\mathrm{dR}}(X/S)^\vee) \ar[r] \ar[d]_{g_s^*} & \Hom_{\mathcal O_S}(H^1_{\mathrm{dR}}(X/S)^\vee, H^1_{\mathrm{dR}}(X/S)^\vee) \ar[d]^{f_s^*} \\
\Ext^1_{\mathcal D_{X_s/k(s)}}(\mathcal O_{X_s},\pi_s^*H^1_{\mathrm{dR}}(X_s/k(s))^\vee) \ar[r] & \Hom_{k(s)}(H^1_{\mathrm{dR}}(X_s/k(s))^\vee, H^1_{\mathrm{dR}}(X_s/k(s))^\vee) }
\end{xy}
\end{equation*}
Moreover, it is straightforward to see that under the left vertical arrow $(2.3.11)$ maps to the class of the extension obtained by performing the construction that led to $(2.3.11)$ with the Poincaré quadruple for $X_s \times \Yr_s$ naturally induced by $(\mathcal P,r,s,\nabla_\mathcal P)$. From this one easily concludes that from the beginning on one may assume $S=\Spec(k)$ with $k$ a field of characteristic zero, which we will henceforth do.\\
\newline
\underline{Step 2}: (The situation on the universal vectorial extension)\\
\newline
\underline{Claim}: There exists a canonical isomorphism of $k$-vector spaces
\[\tag{\textbf{2.3.15}} \Ext^1_{\mathcal O_{\Yr}}(\mathcal O_{\Yr}/\mathcal J, \mathcal J/\mathcal J^2) \xrightarrow{\sim} \Hom_{\mathcal O_{\Yr}}(\mathcal J/\mathcal J^2,\mathcal J/\mathcal J^2)\]under which the class of $(2.3.1)$ maps to the identity.\\
\underline{Proof of the claim}: From the natural exact sequence
\[0 \rightarrow \mathcal J \rightarrow \mathcal O_{\Yr} \rightarrow \mathcal O_{\Yr}/\mathcal J \rightarrow 0\]
we obtain\footnote{Note that
\[\Ext^1_{\mathcal O_{\Yr}}(\mathcal O_{\Yr},\mathcal O_{\Yr}/ \mathcal J) \simeq H^1(\Yr,\mathcal O_{\Yr}/\mathcal J)\simeq H^1(\Spec(k),\mathcal O_{\Spec(k)})=0\]
and that the map $\mathcal O_{\Yr} \rightarrow \mathcal O_{\Yr}/\mathcal J$ induces an isomorphism
\[\Hom_{\mathcal O_{\Yr}}(\mathcal O_{\Yr}/\mathcal J,\mathcal O_{\Yr}/\mathcal J) \xrightarrow{\sim} \Hom_{\mathcal O_{\Yr}}(\mathcal O_{\Yr}, \mathcal O_{\Yr}/\mathcal J).\] 
} from the long exact sequence for $\Hom_{\mathcal O_{\Yr}}(-,\mathcal O_{\Yr}/\mathcal J)$ an isomorphism
\[\Hom_{\mathcal O_{\Yr}}(\mathcal J,\mathcal O_{\Yr}/\mathcal J) \xrightarrow{\sim} \Ext^1_{\mathcal O_{\Yr}}(\mathcal O_{\Yr}/\mathcal J,\mathcal O_{\Yr}/ \mathcal J).\]
We precompose it with the isomorphism
\[\Hom_{\mathcal O_{\Yr}}(\mathcal J/\mathcal J^2,\mathcal O_{\Yr}/\mathcal J) \xrightarrow{\sim} \Hom_{\mathcal O_{\Yr}}(\mathcal J,\mathcal O_{\Yr}/\mathcal J),\]
induced by the natural map $\mathcal J \rightarrow \mathcal J/\mathcal J^2$, in order to obtain an isomorphism
\[\tag{\textbf{2.3.16}} \Hom_{\mathcal O_{\Yr}}(\mathcal J/\mathcal J^2,\mathcal O_{\Yr}/\mathcal J) \xrightarrow{\sim} \Ext^1_{\mathcal O_{\Yr}}(\mathcal O_{\Yr}/\mathcal J,\mathcal O_{\Yr}/ \mathcal J).\]
The desired isomorphism $(2.3.15)$ is then defined to be the composition
\[\begin{split} \Ext^1_{\mathcal O_{\Yr}}(\mathcal O_{\Yr}/\mathcal J, \mathcal J/\mathcal J^2)\simeq \Ext^1_{\mathcal O_{\Yr}}(\mathcal O_{\Yr}/\mathcal J,\mathcal O_{\Yr}/ \mathcal J) \otimes_k \mathcal J/\mathcal J^2\\
\simeq \Hom_{\mathcal O_{\Yr}}(\mathcal J/\mathcal J^2,\mathcal O_{\Yr}/\mathcal J) \otimes_k \mathcal J/\mathcal J^2 \simeq \Hom_{\mathcal O_{\Yr}}(\mathcal J/\mathcal J^2,\mathcal J/\mathcal J^2) \end{split}\]
in which the second identification is given by $(2.3.16)$ and the others are the canonical ones. Let us remark that here and henceforth in this proof $\mathcal J/\mathcal J^2$ and $\mathcal O_{\Yr}/\mathcal J$ are often freely viewed either as $\mathcal O_\Yr$-modules or as $k$-vector spaces. We have thus defined a canonical isomorphism as in the claim.\\
\newline
In terms of a $k$-basis $\{e_1,...,e_{2g}\}$ for $\mathcal J/\mathcal J^2$, if we are given an extension
\[0 \rightarrow \mathcal J/\mathcal J^2 \rightarrow \mathcal F \rightarrow \mathcal O_{\Yr}/\mathcal J \rightarrow 0,\]
then by pushout along the projections corresponding to the $e_i$:
\[p_i:\mathcal J/\mathcal J^2 \rightarrow \mathcal O_{\Yr}/\mathcal J, \quad i=1,...,2g\]
and by using the inverse of the isomorphism $(2.3.16)$ we get an element in
\[\bigoplus_{i=1}^{2g} \Hom_{\mathcal O_{\Yr}}(\mathcal J/\mathcal J^2,\mathcal O_{\Yr}/\mathcal J) \cdot e_i \simeq \Hom_{\mathcal O_{\Yr}}(\mathcal J/\mathcal J^2,\mathcal J/\mathcal J^2)\]
which is the image under $(2.3.15)$ of the class of the given extension.\\
\newline
Now consider the canonical exact sequence of $(2.3.1)$:
\[0 \rightarrow \mathcal J/\mathcal J^2 \rightarrow \mathcal O_{\Yr}/\mathcal J^2 \rightarrow \mathcal O_{\Yr}/\mathcal J \rightarrow 0\]
and write
\[\tag{\textbf{2.3.17}} 0 \rightarrow \mathcal O_{\Yr}/\mathcal J \rightarrow \mathcal F_i \rightarrow \mathcal O_{\Yr}/\mathcal J \rightarrow 0\]
for the pushout of this extension via the above projection $p_i$. The second assertion of the claim follows if we can verify that $p_i$ maps to $(2.3.17)$ under the isomorphism of $(2.3.16)$:
\[\Hom_{\mathcal O_{\Yr}}(\mathcal J/\mathcal J^2,\mathcal O_{\Yr}/\mathcal J) \xrightarrow{\sim} \Ext^1_{\mathcal O_{\Yr}}(\mathcal O_{\Yr}/\mathcal J,\mathcal O_{\Yr}/\mathcal J).\]
This in turn is easily seen as follows: observe the commutative diagram
\begin{equation*}\tag{\textbf{2.3.18}} \begin{split}
\begin{xy}
\xymatrix@C-0.3cm{
0 \ \ar[r] & \mathcal J \ar[r] \ar[d]^{\mathrm{can}} & \mathcal O_{\Yr} \ar[r] \ar[d]^{\mathrm{can}} & \mathcal O_{\Yr}/\mathcal J \ar[r] \ar[d]^{\id} &   0\\
0 \ \ar[r] & \mathcal J/\mathcal J^2 \ar[r] \ar[d]^{p_i} & \mathcal O_{\Yr}/\mathcal J^2 \ar[r] \ar[d] & \mathcal O_{\Yr}/\mathcal J \ar[r] \ar[d]^{\id} &   0\\
0 \ \ar[r] & \mathcal O_{\Yr}/\mathcal J \ar[r]  & \mathcal F_i \ar[r]  & \mathcal O_{\Yr}/\mathcal J \ar[r] &  0}
\end{xy}
\end{split}
\end{equation*}
in which the second row is the pushout of the first along $\mathcal J \xrightarrow{\mathrm{can}} \mathcal J/\mathcal J^2$ (obvious) and the third is the pushout of the second along $p_i$ (by definition) such that the third row is the pushout of the first along the map $p_i \circ \mathrm{can}$. Now recall that the isomorphism $(2.3.16)$ was defined as the composition
\[\Hom_{\mathcal O_{\Yr}}(\mathcal J/\mathcal J^2,\mathcal O_{\Yr}/\mathcal J) \xrightarrow{\sim} \Hom_{\mathcal O_{\Yr}}(\mathcal J,\mathcal O_{\Yr}/\mathcal J) \xrightarrow{\sim} \Ext^1_{\mathcal O_{\Yr}}(\mathcal O_{\Yr}/\mathcal J,\mathcal O_{\Yr}/ \mathcal J)\]
in which the first map is induced by the arrow $\mathcal J \xrightarrow{\mathrm{can}} \mathcal J/\mathcal J^2$ and the second is given by pushing out the top row of $(2.3.18)$ along morphisms $\mathcal J \rightarrow \mathcal O_{\Yr}/\mathcal J$. This clearly implies the remaining claim.\\
\newline
\underline{Step 3}: (A reduction step)\\
\newline
As $S=\Spec(k)$ the kernel $\Ext^1_{\mathcal O_S}(\mathcal O_S,\mathcal H)$ of the projection $(2.3.14)$ vanishes, and hence $(2.3.14)$ is in fact an isomorphism:
\[\Ext^1_{\mathcal D_{X/k}}(\mathcal O_X, \mathcal H_X) \xrightarrow{\sim} \mathrm{Hom}_k(\mathcal H, \mathcal H).\]
By Step 2 the claim of the theorem is proven if the following diagram commutes:
\begin{equation*}\tag{\textbf{2.3.19}}\begin{split}
\begin{xy}
\xymatrix@C-0.3cm{
\Ext^1_{\mathcal D_{X/k}}(\mathcal O_X, \mathcal H_X) \ar[r]^{ \ \ \sim} & \mathrm{Hom}_k(\mathcal H, \mathcal H) \ar@{-}[d]^{\sim} \\
\Ext^1_{\mathcal O_{\Yr}}(\mathcal O_{\Yr}/\mathcal J, \mathcal J/\mathcal J^2)\ar[u]^{\widehat{(.)}} \ar[r]^{\sim} & \Hom_{\mathcal O_{\Yr}}(\mathcal J/\mathcal J^2,\mathcal J/\mathcal J^2)}
\end{xy}
\end{split}
\end{equation*}
Here, the upper resp. lower horizontal arrow is $(2.3.14)$ resp. $(2.3.15)$, the right vertical map is induced by $(2.3.2)$ and the left vertical map is given as follows: analogously as explained in $(2.3.1)$-$(2.3.4)$ an extension of $\mathcal O_{\Yr}/\mathcal J$ by $\mathcal J/\mathcal J^2$ writes as exact sequence of $\mathcal O_{\Yr}$-modules
\[0 \rightarrow (\epsilon^{\natural})_*\mathcal H \rightarrow \mathcal F \rightarrow (\epsilon^{\natural})_*\mathcal O_{\Spec(k)} \rightarrow 0\]
whose terms are WIT-sheaves of index $0$, providing an exact sequence of Fourier-Mukai transforms
\[0 \rightarrow \mathcal H_X \rightarrow \widehat{\mathcal F} \rightarrow \mathcal O_X \rightarrow 0,\]
where $\mathcal H_X$ resp. $\mathcal O_X$ has the canonical integrable $k$-connection (cf. Prop. 2.2.6 and Lemma 2.2.7).\\
It is a routine task to deduce commutativity of $(2.3.19)$ if the following diagram is known to commute:
\begin{equation*}\tag{\textbf{2.3.20}}\begin{split}
\begin{xy}
\xymatrix@C-0.3cm{
\Ext^1_{\mathcal D_{X/k}}(\mathcal O_X, \mathcal O_X) \ar[r]^{\sim \qquad} & H^1_{\mathrm{dR}}(X/k)\ar@{-}[d]^{\sim} \\
\Ext^1_{\mathcal O_{\Yr}}(\mathcal O_{\Yr}/\mathcal J, \mathcal O_{\Yr}/\mathcal J)\ar[u]^{\widehat{(.)}} \ar[r]^{\sim \qquad \qquad \quad} & \Hom_{\mathcal O_{\Yr}}(\mathcal J/\mathcal J^2,\mathcal O_{\Yr}/\mathcal J) \simeq Hom_k(\mathcal J_{e^{\natural}}/\mathcal J^2_{e^{\natural}},k)}
\end{xy}
\end{split}
\end{equation*}
where with $e^{\natural}\in\Yr$ we denote the zero point of $\Yr$. The left vertical arrow is defined analogously as in $(2.3.19)$ by taking Fourier-Mukai transforms, the right vertical arrow is induced by $(2.3.2)$, the lower horizontal map is given by $(2.3.16)$ and the upper horizontal map is defined in the same way as $(2.3.14)$, with $\mathcal H_X$ replaced by $\mathcal O_X$.\\
It is clear that to prove commutativity of $(2.3.20)$ we will need a better understanding of its right vertical arrow and thus of the canonical isomorphism of $k$-vector spaces induced by $(2.3.2)$:
\[\tag{\textbf{2.3.21}} (\mathcal J_{e^{\natural}}/\mathcal J^2_{e^{\natural}})^\vee \simeq H^1_{\mathrm{dR}}(X/k).\]

\underline{Step 4}: (The identification $(2.3.21)$)\\
\newline
Abbreviating the scheme of dual numbers over $k$ with $D:=\Spec(k[\epsilon]/(\epsilon^2))$, the identification in question comes about as the composition
\[\tag{\textbf{2.3.22}} (\mathcal J_{e^{\natural}}/\mathcal J^2_{e^{\natural}})^\vee \simeq \ker(\Yr(D) \rightarrow \Yr(k)) \simeq H^1_{\mathrm{dR}}(X/k),\]
where the first identification is standard (cf. e.g. \cite{Gö-We}, Ch. 6, $(6.4)$) and the second is described in \cite{Maz-Mes}, Ch. I, § 4, as follows:\\
By definition (cf. 0.1.1) and Lemma 0.1.12 for a $k$-scheme $T$ the $T$-valued points of $\Yr$ are given by
\[\mathrm{Pic}^\natural(X_T/T)=\textrm{\{Isomorphism classes of triples} \ (\mathcal L, \alpha, \nabla_{\mathcal L})\},\]
where $\mathcal L$ is a line bundle on $X_T=X\times_kT$ with $T$-rigidification $\alpha$ and integrable $T$-connection $\nabla_{\mathcal L}$.\\
If $T$ has trivial Picard group, then the assignment $(\mathcal L, \alpha, \nabla_{\mathcal L}) \mapsto (\mathcal L, \nabla_{\mathcal L})$ is easily checked to induce an isomorphism of groups
\[\mathrm{Pic}^\natural(X_T/T) \simeq \textrm{\{Isomorphism classes of pairs} \ (\mathcal L, \nabla_{\mathcal L})\},\]
and (as in \cite{Maz-Mes}, Ch. I, Prop. $(4.1.2)$) we identify the last group with $\mathbb H^1(X_T,\Omega^*_{X_T/T})$, where
\[\Omega^*_{X_T/T}: \quad \lbrack \mathcal O_{X_T}^* \xrightarrow{\mathrm{dlog}} \Omega^1_{X_T/T} \xrightarrow{\mathrm{d}} \Omega^2_{X_T/T} \xrightarrow{\mathrm{d}} ... \rbrack \]
denotes the multiplicative de Rham complex for $X_T/T$ (starting in degree zero). We thus obtain:
\[\tag{\textbf{2.3.23}} \ker(\Yr(D) \rightarrow \Yr(k)) \simeq \ker(\mathbb H^1(X_D,\Omega^*_{X_D/D}) \rightarrow \mathbb H^1(X,\Omega^*_{X/k})).\]
The long exact sequence for hypercohomology associated with the canonically split natural short exact sequence of complexes of abelian sheaves on $X$
\[0 \rightarrow \Omega^\bullet_{X/k} \rightarrow \Omega^*_{X_D/D} \rightarrow \Omega^*_{X/k} \rightarrow 0\]
identifies $H^1_{\mathrm{dR}}(X/k)$ with $\ker(\mathbb H^1(X_D,\Omega^*_{X_D/D}) \rightarrow \mathbb H^1(X,\Omega^*_{X/k}))$. Combined with $(2.3.23)$ we obtain in sum an isomorphism of groups as in $(2.3.22)$ which is in fact $k$-linear.\\
A final remark:\\
With $\mathcal O_{X_D}=\mathcal O_X \oplus \epsilon \cdot \mathcal O_X$, $\mathcal O_{X_D}^*=\mathcal O_X^*\oplus \epsilon\cdot \mathcal O_X$ and $\Omega^i_{X_D/D} \simeq \Omega^i_{X/k} \otimes_{\mathcal O_X} \mathcal O_{X_D} \simeq \Omega^i_{X/k} \oplus \epsilon \cdot \Omega^i_{X/k}$ the previous split exact sequence of complexes writes as
\begin{equation*} \tag{\textbf{2.3.24}} \begin{split}
\begin{xy}
\xymatrix@C-0.3cm{
0 \ \ar[r] & \mathcal O_X \ar[r]^{\eta \qquad } \ar[d]^{\mathrm{d}} & \mathcal O_X^*\oplus \epsilon \cdot \mathcal O_X \ar[r] \ar[d]^{\theta} & \mathcal O_X^* \ar[r] \ar[d]^{\mathrm{dlog}} &   0\\
0 \ \ar[r] & \Omega^1_{X/k} \ar[r] \ar[d]^{\mathrm{d}} & \Omega^1_{X/k} \oplus \epsilon \cdot \Omega^1_{X/k} \ar[r] \ar[d]^{\mathrm{d}+ \epsilon \cdot \mathrm{d}} &  \Omega^1_{X/k}\ar[r] \ar[d]^{\mathrm{d}} &   0\\
0 \ \ar[r] & \Omega^2_{X/k} \ar[r] \ar[d]^{\mathrm{d}} & \Omega^2_{X/k} \oplus \epsilon \cdot \Omega^2_{X/k} \ar[r] \ar[d]^{\mathrm{d}+ \epsilon \cdot \mathrm{d}} & \Omega^2_{X/k}\ar[r] \ar[d]^{\mathrm{d}} &   0\\
& ... & ... & ... & &}
\end{xy}
\end{split}
\end{equation*}
where $\theta(u+ \epsilon \cdot v)= \mathrm{dlog}(u)+ \epsilon \cdot \big(\frac{\mathrm{d}v}{u}-\frac{v}{u}\cdot \mathrm{dlog}(u) \big)$, $\eta(v)=1+\epsilon \cdot v$, the other left horizontal arrows are given by multiplication with $\epsilon$ and inclusion into the second component, the right horizontal arrows by projection to the first component and the sections by inclusion into the first component.\\
\newline
\underline{Step 5}: (The map $KS$)\\
\newline
Denote by $\mathrm{pr}_X, \mathrm{pr}_D$ the projections of $X_D=X\times_kD$ and by $i_X: X \rightarrow X_D$ the nilpotent closed immersion of square zero induced by base change of the canonical closed immersion $i:\Spec(k) \rightarrow D$.
\begin{equation*} \tag{\textbf{2.3.25}} \begin{split}
\begin{xy}
\xymatrix{
X \ar[d]^{i_X} \ar[r]^{\pi \quad \ } \ar@ /_0.9cm/[dd]_{\id} & \Spec(k) \ar[d]_{i} \ar@ /^0.9cm/[dd]^{\id}\\
X_D \ar[d]^{\mathrm{pr}_X} \ar[r]^{\mathrm{pr}_D} & \ar[d] D \\
X \ar[r]^{\pi \quad \ } & \Spec(k)}
\end{xy}
\end{split}
\end{equation*}
Moreover, observe the $(k[\epsilon]/(\epsilon)^2)$-linear exact sequence
\[0 \rightarrow k \xrightarrow{\cdot \epsilon} k[\epsilon]/(\epsilon)^2 \rightarrow k \rightarrow 0\]
with maps $a \mapsto \epsilon \cdot a$ resp. $a +\epsilon \cdot b \mapsto a$ and $(k[\epsilon]/(\epsilon)^2)$-module structure of $k$ defined by the second map. The associated $\mathcal O_D$-linear exact sequence writes as
\[\tag{\textbf{2.3.26}} 0 \rightarrow i_*\mathcal O_{\Spec(k)} \xrightarrow{\cdot \epsilon} \mathcal O_D \rightarrow i_*\mathcal O_{\Spec(k)} \rightarrow 0,\]
the second arrow belonging to the closed immersion $i$.\\
Now assume that we are given an element $f: D \rightarrow \Yr$ of $\ker(\Yr(D) \rightarrow \Yr(k))$.\\
The $\mathcal O_{X_D}$-line bundle $(\mathrm{id}_X \times f)^*\mathcal P$ carries the integrable $D$-connection induced by $\nabla_\mathcal P$ by pullback along the lower square of the diagram
\begin{equation*} \tag{\textbf{2.3.27}} \begin{split}
\begin{xy}
\xymatrix{
X \ar[d]^{i_X} \ar[r] \ar@ /_0.9cm/[dd]_{\mathrm{id}_X \times \epsilon^{\natural} }& \Spec(k) \ar[d]_{i} \ar@ /^0.9cm/[dd]^{\epsilon^{\natural}}\\
X\times_k D \ar[d]^{\mathrm{id}_X \times f} \ar[r]^{ \quad \mathrm{pr}_D} & D \ar[d]_{f}\\
X\times_k \Yr \ar[r] & \Yr}
\end{xy}
\end{split}
\end{equation*}
Pull back the sequence $(2.3.26)$ along $\mathrm{pr}_D$ (the pullbacks endowed with the canonical integrable $D$-connections), then tensor with $(\mathrm{id}_X \times f)^*\mathcal P$ and finally push out along the lower square of $(2.3.25)$:
\[(\mathrm{pr}_X)_*((\mathrm{id}_X \times f)^*\mathcal P \otimes_{\mathcal O_{X_D}}\mathrm{pr}^*_D(.)).\]
Observing Lemma 2.2.5 and the $\mathcal D_{X/k}$-linear trivialization $s$ of $(\mathcal P,\nabla_\mathcal P)$ along $\id_X \times \epsilon^\natural$ this procedure transforms $(2.3.26)$ into an exact sequence of $\mathcal D_{X/k}$-modules:
\[0 \rightarrow \mathcal O_X \rightarrow (\mathrm{pr}_X)_*(\mathrm{id}_X \times f)^*\mathcal P \rightarrow \mathcal O_X \rightarrow 0,\]
where $\mathcal O_X$ is endowed with exterior $k$-derivation and $(\mathrm{pr}_X)_*(\mathrm{id}_X\times f)^*\mathcal P$ with the pushout connection along the lower square of $(2.3.25)$.\\
This assignment yields a $k$-linear arrow
\[\tag{\textbf{2.3.28}} KS: \ker(\Yr(D) \rightarrow \Yr(k)) \rightarrow \Ext^1_{\mathcal D_{X/k}}(\mathcal O_X, \mathcal O_X).\]
If we forget the integrable connections, then this is precisely the usual Kodaira-Spencer map at the zero point $e^\natural \in \Yr$ associated to the sheaf $\mathcal P$ (cf. e.g. \cite{Bri}, 3.4, for this notion).\\
It fits into a commutative diagram
\begin{equation*}\tag{\textbf{2.3.29}}\begin{split}
\begin{xy}
\xymatrix@C-0.3cm{
\Ext^1_{\mathcal D_{X/k}}(\mathcal O_X, \mathcal O_X) & \ker(\Yr(D) \rightarrow \Yr(k)) \ar[l]_{KS} \ar@{-}[d]^{\sim} \\
\Ext^1_{\mathcal O_{\Yr}}(\mathcal O_{\Yr}/\mathcal J, \mathcal O_{\Yr}/\mathcal J)\ar[u]^{\widehat{(.)}} & (\mathcal J_{e^{\natural}}/\mathcal J^2_{e^{\natural}})^\vee \ar@{-}[l]^{\qquad \quad \sim}}
\end{xy}
\end{split}
\end{equation*}
where the left vertical and lower horizontal arrow is as in $(2.3.20)$ and the right vertical identification is the canonical one, already mentioned in $(2.3.22)$.\\
For the commutativity of $(2.3.29)$ one first checks (rather straightforwardly from the definitions) that the composite of the right vertical and lower horizontal identification
\[\ker(\Yr(D) \rightarrow \Yr(k)) \xrightarrow{\sim} (\mathcal J_{e^{\natural}}/\mathcal J^2_{e^{\natural}})^\vee \xrightarrow{\sim} \Ext^1_{\mathcal O_{\Yr}}(\mathcal O_{\Yr}/\mathcal J, \mathcal O_{\Yr}/\mathcal J)\]
maps an arrow $f: D \rightarrow \Yr$ to the direct image of $(2.3.26)$ along the (affine) map $f$:
\[0 \rightarrow \mathcal O_{\Yr}/\mathcal J \rightarrow f_*\mathcal O_D \rightarrow \mathcal O_{\Yr}/\mathcal J \rightarrow 0.\]
That $KS(f)$ is isomorphic to the extension given by Fourier-Mukai transformation of the previous sequence is a simple consequence of the definitions together with Lemma 2.2.5. A slightly different argument for this compatibility of the Kodaira-Spencer map and the Fourier-Mukai transformation (without considering integrable connections) can be found in \cite{Bri}, Lemma 4.2.3.\\
\newline
\underline{Step 6}: (A further reduction step)\\
\newline
Recall from Step 3 that the theorem is proven as soon as $(2.3.20)$ is seen to commute. But as $(2.3.29)$ is commutative a brief reflection shows that this is the case as soon as we know commutativity of
\begin{equation*} \tag{\textbf{2.3.30}} \begin{split}
\begin{xy}
\xymatrix@C-0.3cm{
 & \Ext^1_{\mathcal D_{X/k}}(\mathcal O_X, \mathcal O_X) \ar[rr]^{\sim} & &  H^1_{\mathrm{dR}}(X/k) \\
 & & \ker(\Yr(D) \rightarrow \Yr(k)) \ar[ul]_{KS}   \ar[ur]^{\sim}&}
\end{xy}
\end{split}
\end{equation*}
where the upper arrow is as in $(2.3.20)$ and the identification on the right is as explained in Step 4.\\
Our remaining task thus consists in verifying that $(2.3.30)$ commutes.\\
\newline
\underline{Step 7}: (Some preparations for the final step of proof)\\
\newline
By Step 4 (and the definition of the Poincaré bundle on $X\times_S \Yr$) the right arrow in $(2.3.30)$ is given as follows: for a (henceforth fixed) element $f: D \rightarrow \Yr$ of $\ker(\Yr(D) \rightarrow \Yr(k))$ the $\mathcal O_{X_D}$-line bundle $(\mathrm{id}_X \times f)^*\mathcal P$ with its integrable $D$-connection (induced by $\nabla_\mathcal P$) defines a class in $\mathbb H^1(X_D,\Omega^*_{X_D/D})$ which maps to zero under the projection in the exact sequence
\[\tag{\textbf{2.3.31}} 0\rightarrow \mathbb H^1(X,\Omega_{X/k}^{\bullet}) \rightarrow \mathbb H^1(X_D,\Omega^*_{X_D/D}) \rightarrow \mathbb H^1(X,\Omega^*_{X/k}) \rightarrow 0\]
induced by the split short exact sequence of complexes of abelian sheaves on $X$ recorded in $(2.3.24)$.\\
The thus obtained class in $\mathbb H^1(X,\Omega_{X/k}^{\bullet})$ gives the desired image of $f$.\\
On the other hand and as explained in Step 5, applying $(\mathrm{pr}_X)_*((\mathrm{id}_X\times f)^*\mathcal P\otimes_{\mathcal O_{X_D}}\mathrm{pr}_D^*(.))$ to
\[0 \rightarrow i_*\mathcal O_{\Spec(k)} \xrightarrow{\cdot \epsilon} \mathcal O_D \rightarrow i_*\mathcal O_{\Spec(k)} \rightarrow 0\]
yields the $\mathcal D_{X/k}$-linear exact sequence
\[\tag{\textbf{2.3.32}} 0 \rightarrow \mathcal O_X \rightarrow (\mathrm{pr}_X)_*(\mathrm{id}_X\times f)^*\mathcal P \rightarrow \mathcal O_X \rightarrow 0.\]
The image of $1\in k$ under the map $k \rightarrow \mathbb H^1(X,\Omega^{\bullet}_{X/k})$ which appears in the long exact sequence of hypercohomology for the associated sequence of de Rham complexes is then the image of $f$ under the composition of the left with the upper arrow of $(2.3.30)$.\\
For the following, observe that the topological spaces of $X_D$ and $X$ are the same and that taking direct image of a $\mathcal O_{X_D}$-module along $\mathrm{pr}_X: X_D \rightarrow X$ just means considering it as $\mathcal O_X$-module via the morphism of ring sheaves $\mathcal O_X \rightarrow \mathcal O_{X_D}=\mathcal O_X \oplus \epsilon\cdot \mathcal O_X$ given by inclusion into the first component. As it will be clear which structure is meant we will thus henceforth leave away the notation $(\mathrm{pr}_X)_*$.\\
From now on let us abbreviate
\[\mathcal L:=(\mathrm{id}_X\times f)^*\mathcal P\]
which is an abelian sheaf on the topological space $X_D$ resp. $X$ and can be viewed as invertible $\mathcal O_{X_D}$-module with integrable $D$-connection resp. as $\mathcal O_X$-vector bundle with integrable $k$-connection. If we then write the exact $\mathcal D_{X/k}$-linear sequence $(2.3.32)$ as
\[\tag{\textbf{2.3.33}} 0 \rightarrow \mathcal O_X \xrightarrow{\psi} \mathcal L \xrightarrow{\varphi} \mathcal O_X \rightarrow 0\]
it is easy to see that the $\mathcal O_{X_D}=\mathcal O_X \oplus \epsilon\cdot \mathcal O_X$-structure of $\mathcal L$ is recovered from $(2.3.33)$ by the formula
\[\tag{\textbf{2.3.34}} (a+\epsilon \cdot b)\cdot l =a\cdot l +b\cdot (\psi \circ \varphi)(l).\]
The $D$-connection and the $k$-connection on $\mathcal L$ translate into each other via the canonical identification
\[\Omega^1_{X_D/D}\otimes_{\mathcal O_{X_D}}\mathcal L \simeq (\Omega^1_{X/k} \otimes_{\mathcal O_X}\mathcal O_{X_D}) \otimes_{\mathcal O_{X_D}}\mathcal L \simeq \Omega^1_{X/k}\otimes_{\mathcal O_X}\mathcal L,\]
and we will denote both of them with $\nabla_{\mathcal L}$.\\
Finally, one checks that if a local section $s: \mathcal O_U \rightarrow \mathcal L_{|U}$ of $\varphi$ over an open subset $U \subseteq X$ is given, then $\mathcal L$ as $\mathcal O_{X_D}$-module is trivialized on $U$ by means of the map
\[\tag{\textbf{2.3.35}} \mathcal L_{|U} \xrightarrow{\sim} \mathcal O_U \oplus \epsilon \cdot \mathcal O_U, \quad l \mapsto \varphi(l)+\epsilon \cdot (l-(s\circ \varphi)(l)),\]
where we tacitly view $\mathcal O_X$ as included into $\mathcal L$ via $\psi$.\\
\newline
\underline{Step 8}: (Conclusion of the proof)\\
\newline
Choose an open covering $\{U_i\}_{i\in I}$ of $X$ with local sections $s_i: \mathcal O_{U_i} \rightarrow \mathcal L_{|U_i}$ of $\varphi$ and write
\[l_i:=s_i(1) \in \Gamma(U_i,\mathcal L).\]
Consider
\[(l_j-l_i, \nabla_{\mathcal L}(l_i))\in \prod_{i,j} \Gamma (U_{ij},\mathcal L) \oplus \prod_{i}\Gamma(U_i, \Omega^1_{X/k}\otimes_{\mathcal O_X} \mathcal L).\]
Then there clearly exists
\[(u_{ij}, \omega_i)\in \prod_{i,j} \Gamma (U_{ij},\mathcal O_X) \oplus \prod_{i}\Gamma(U_i, \Omega^1_{X/k})\]
with $\psi(u_{ij})=l_j-l_i$ and $(\id\otimes \psi)(\omega_i)=\nabla_{\mathcal L}(l_i)$, where $\id\otimes \psi: \Omega^1_{X/k} \rightarrow \Omega^1_{X/k}\otimes_{\mathcal O_X}\mathcal L$.\footnote{For the final result of our \v{C}ech hypercohomology computations in this last step of proof it is irrelevant which sign convention is adopted for the hyperdifferentials resp. for the cocycle representing the class of $(\mathcal L, \nabla_{\mathcal L})$, as long as this is done consistently.}\\
Note that $(u_{ij}, \omega_i)$ is a cocycle which represents in $\mathbb H^1(X,\Omega^{\bullet}_{X/k})$ the image of $f$ under the composition of the left with the upper arrow of $(2.3.30)$ (cf. Step 7).\\
On the other hand, we now represent the class of $(\mathcal L,\nabla_{\mathcal L})$ in $\mathbb H^1(X_D,\Omega^*_{X_D/D})$ by a cocycle in
\[\prod_{i,j} \Gamma (U_{ij},\mathcal O_{X_D}^*) \oplus \prod_{i}\Gamma(U_i, \Omega^1_{X_D/D})=\prod_{i,j} \Gamma (U_{ij},\mathcal O_X^*\oplus \epsilon \cdot \mathcal O_X) \oplus \prod_{i}\Gamma(U_i, \Omega^1_{X/k}\oplus \epsilon \cdot \Omega^1_{X/k}).\]
Recall from $(2.3.35)$ that we have a trivialization on $U_i$ of $\mathcal L$ as $\mathcal O_{X_D}$-module, defined by
\[\tag{\textbf{2.3.36}} t_i: \mathcal L_{|U_i} \xrightarrow{\sim} \mathcal O_{U_i} \oplus \epsilon \cdot \mathcal O_{U_i}, \quad l \mapsto \varphi (l) + \epsilon \cdot (l-(s_i\circ \varphi) (l)).\]
A representative of $(\mathcal L,\nabla_{\mathcal L})$ in $\mathbb H^1(X_D,\Omega^*_{X_D/D})$ is then given by the expression
\[((t_i \circ t_j^{-1})(1), \eta_i) \in \prod_{i,j} \Gamma (U_{ij},\mathcal O_X^*\oplus \epsilon \cdot \mathcal O_X) \oplus \prod_{i}\Gamma(U_i, \Omega^1_{X/k}\oplus \epsilon \cdot \Omega^1_{X/k}),\]
where $\eta_i$ is the image of $1$ under the composition
\[\tag{\textbf{2.3.37}} \begin{split} &\mathcal O_{U_i} \oplus \epsilon \cdot \mathcal O_{U_i} \xrightarrow{t_i^{-1}} \mathcal L_{|U_i} \xrightarrow{\nabla_{\mathcal L}} (\Omega^1_{X_D/D} \otimes_{\mathcal O_{X_D}} \mathcal L)_{|U_i} \\
&\simeq ((\Omega^1_{X/k}\oplus \epsilon \cdot \Omega^1_{X/k}) \otimes_{\mathcal O_X \oplus \epsilon \cdot \mathcal O_X} \mathcal L)_{|U_i}\xrightarrow{\id \otimes t_i} \Omega^1_{U_i/k}\oplus \epsilon \cdot \Omega^1_{U_i/k}.\end{split}\]
Observe that the section $l_i$ corresponds to $1$ in $(2.3.36)$. We obtain
\[(t_i \circ t_j^{-1})(1)=\varphi (l_j) + \epsilon \cdot (l_j-(s_i\circ \varphi) (l_j))=1+\epsilon \cdot (s_j(1)-(s_i \circ \varphi \circ s_j)(1))= 1+\epsilon \cdot (l_j-l_i)=1+\epsilon \cdot u_{ij}.\]
Furthermore, as section of $\Omega^1_{X/k}\otimes_{\mathcal O_X}\mathcal L$ over $U_i$ we have
\[\nabla_{\mathcal L}(l_i)=(\id\otimes \psi)(\omega_i)=\omega_i \otimes \psi(1) =\omega_i\otimes (\psi\circ\varphi)(l_i)=\omega_i \otimes \epsilon \cdot l_i,\]
the first three equations by definition and the last by $(2.3.34)$. But under the identification
\[\Omega^1_{X/k}\otimes_{\mathcal O_X}\mathcal L \simeq \Omega^1_{X_D/D}\otimes_{\mathcal O_{X_D}}\mathcal L \simeq (\Omega^1_{X/k}\oplus \epsilon \cdot \Omega^1_{X/k}) \otimes_{\mathcal O_X \oplus \epsilon \cdot \mathcal O_X} \mathcal L\]
the section $\omega_i \otimes \epsilon \cdot l_i$ corresponds to $\epsilon\cdot \omega_i \otimes l_i$, as one can readily check. We conclude that the image of $1$ under the chain $(2.3.37)$ is the section $\epsilon \cdot \omega_i$ of $\Omega^1_{U_i/k}\oplus \epsilon \cdot \Omega^1_{U_i/k}$.\\
Altogether, we have shown that the class of $(\mathcal L,\nabla_{\mathcal L})$ in $\mathbb H^1(X_D, \Omega^*_{X_D/D})$ is represented by
\[((t_i \circ t_j^{-1})(1), \eta_i)=(1+\epsilon \cdot u_{ij}, \epsilon \cdot \omega_i).\]
From $(2.3.24)$ it then obviously follows that the preceding class comes from the class of the cocycle $(u_{ij},\omega_i)$ under the inclusion in $(2.3.31)$:
\[\mathbb H^1(X,\Omega^{\bullet}_{X/k})\rightarrow \mathbb H^1(X_D,\Omega^*_{X_D/D}).\]
Hence, $(u_{ij},\omega_i)$ represents the image of $f$ under the right arrow in $(2.3.30)$ (cf. Step 7). But - as remarked at the beginning of the present step of proof - it also represents the image of $f$ under the composition of the left with the upper arrow in $(2.3.30)$. The commutativity of $(2.3.30)$ is thus shown, which according to Step 6 concludes the proof of the theorem.
\end{proof}

\begin{corollary}
Consider the exact sequence of $\mathcal D_{X/S}$-modules constructed in $(2.3.11)$:
\[0 \rightarrow \mathcal H_X \rightarrow (p_1)_*\mathcal P_1 \rightarrow \mathcal O_X \rightarrow 0\]
together with the $\mathcal O_S$-linear splitting for its pullback along $\epsilon: S \rightarrow X$ constructed in $(2.3.13)$:
\[\mathcal O_S \oplus \mathcal H \simeq \epsilon^*(p_1)_*\mathcal P_1.\]
Then the integrable $S$-connection on $(p_1)_*\mathcal P_1$ has a unique prolongation to an integrable $\Q$-connection such that the previous data become the first logarithm sheaf of $X/S/\Q$ in the sense of Def. 1.1.1.
\end{corollary}
\begin{proof}
As already explained this follows by combining Prop. 2.1.4 and Thm. 2.3.1.
\end{proof}

\section{The higher logarithm sheaves and the Poincaré bundle}
\markright{\uppercase{The logarithm sheaves and the Poincaré bundle}}
\subsection{An equivalence of categories}
The following auxiliary result is analogous to \cite{Lau}, $(2.3)$. For our purposes, however, it is convenient to give a different and direct proof, without introducing an intermediate equivalent category.

\begin{lemma}
Let $\mathcal Y$ be a vector bundle on $S$ and denote by $(\mathcal Y_X, \nabla_{\mathcal Y_X})$ its pullback via $\pi:X \rightarrow S$ together with its canonical integrable $S$-connection. Then the following categories are equivalent:\\
\newline
(1) The category of $\mathcal D_{X/S}$-linear extensions of $(\mathcal Y_X, \nabla_{\mathcal Y_X})$ by $(\mathcal O_X,\mathrm{d})$:
\[0 \rightarrow (\mathcal O_X,\mathrm{d}) \rightarrow (\mathcal U, \nabla_{\mathcal U}) \rightarrow (\mathcal Y_X, \nabla_{\mathcal Y_X}) \rightarrow 0.\]
(2) The category consisting of (associative, commutative, unital) quasi-coherent $(\mathcal O_X,\mathrm{d})$-algebras with integrable $S$-connection $(\mathcal B, \nabla_{\mathcal B})$\footnote{This shall mean that $\mathcal B$ is a quasi-coherent $\mathcal O_X$-algebra such that the structure map $\mathcal O_X \rightarrow\mathcal B$ and the multiplication map $\mathcal B \otimes_{\mathcal O_X}\mathcal B \rightarrow\mathcal B$ are horizontal.} together with an exhaustive filtration by locally free $\mathcal O_X$-submodules of finite rank with integrable $S$-connection $(\mathcal B_n, \nabla_{\mathcal B_n}), n\geq  0$:
\[(0)\subseteq (\mathcal B_0, \nabla_{\mathcal B_0}) \subseteq (\mathcal B_1, \nabla_{\mathcal B_1}) \subseteq (\mathcal B_2, \nabla_{\mathcal B_2}) \subseteq...\subseteq 
(\mathcal B, \nabla_{\mathcal B})\]
such that $\mathcal B_n \cdot \mathcal B_m \subseteq \mathcal B_{n+m}$, and together with an isomorphism of graded $(\mathcal O_X,\mathrm{d})$-algebras with integrable $S$-connection
\[\psi_{\bullet}: \mathrm{gr}_{\bullet}(\mathcal B, \nabla_{\mathcal B})\simeq \mathrm{Sym}^{\bullet}_{\mathcal O_X}(\mathcal Y_X, \nabla_{\mathcal Y_X}).\]
\end{lemma}
\begin{proof}
We define functors $S: (1) \rightarrow (2), T: (2) \rightarrow (1)$ which will be quasi-inverse to each other.\\
\underline{Definition of $T$}:\\
Given an object in $(2)$ the exact sequence of $\mathcal D_{X/S}$-modules
\[0 \rightarrow (\mathcal B_0,\nabla_{\mathcal B_0}) \rightarrow (\mathcal B_1,\nabla_{\mathcal B_1}) \rightarrow (\mathcal B_1, \nabla_{\mathcal B_1})/(\mathcal B_0,\nabla_{\mathcal B_0}) \rightarrow 0\]
identifies by means of $\psi_0$ and $\psi_1$ with an exact sequence
\[0 \rightarrow (\mathcal O_X,\mathrm{d}) \rightarrow (\mathcal B_1,\nabla_{\mathcal B_1})\rightarrow (\mathcal Y_X, \nabla_{\mathcal Y_X}) \rightarrow 0\]
which shall be the associated object in $(1)$.\\
\underline{Definition of $S$}:\\
Let an extension
\[0 \rightarrow (\mathcal O_X,\mathrm{d}) \xrightarrow{j_1} (\mathcal U, \nabla_{\mathcal U}) \xrightarrow{p_1} (\mathcal Y_X, \nabla_{\mathcal Y_X}) \rightarrow 0\]
be given and note that $\mathcal U$ is a $\mathcal O_X$-vector bundle.\\
For $n\geq 0$ set $(\mathcal B_n,\nabla_{\mathcal B_n}):=\mathrm{Sym}^n_{\mathcal O_X}(\mathcal U, \nabla_{\mathcal U})$ with horizontal monomorphisms resp. epimorphisms
\[j_{n+1}:\mathcal B_n \rightarrow \mathcal B_{n+1}  \ \ \textrm{resp.} \ \ p_{n+1}:\mathcal B_{n+1} \rightarrow \mathrm{Sym}^{n+1}_{\mathcal O_X}(\mathcal Y_X),\]
given locally by
\[u_1 \cdot...\cdot u_n \mapsto u_1 \cdot...\cdot u_n\cdot j_1(1) \ \ \textrm{resp.} \ \  u_1 \cdot...\cdot u_{n+1} \mapsto p_1(u_1) \cdot...\cdot p_1(u_{n+1}), \qquad\]
from which one gets exact sequences
\[\tag{\textbf{2.4.1}} 0 \rightarrow (\mathcal B_n,\nabla_{\mathcal B_n}) \xrightarrow{j_{n+1}} (\mathcal B_{n+1},\nabla_{\mathcal B_{n+1}}) \xrightarrow{p_{n+1}} \mathrm{Sym}^{n+1}_{\mathcal O_X}(\mathcal Y_X, \nabla_{\mathcal Y_X})\rightarrow 0.\]
Define a quasi-coherent $\mathcal O_X$-module $\mathcal B:=\mathrm{Sym}^{\bullet}_{\mathcal O_X}(\mathcal U)/\left(1-j_1(1)\right)$. Note that $\mathcal B$ is the direct limit of the $\mathcal B_n$ for the above maps $j_{n+1}$. Endow it with the $\mathcal O_X$-algebra structure coming from $\mathrm{Sym}^{\bullet}_{\mathcal O_X}(\mathcal U)$ and with the induced integrable $S$-connection, denoted by $\nabla_{\mathcal B}$. Together with the filtration and the isomorphism $\psi_{\bullet}: \mathrm{gr}_{\bullet}(\mathcal B, \nabla_{\mathcal B})\simeq \mathrm{Sym}^{\bullet}_{\mathcal O_X}(\mathcal Y_X, \nabla_{\mathcal Y_X})$ induced by $(2.4.1)$ we obtain an object in $(2)$.\\
That $T \circ S \simeq \mathrm{id}_{(1)}$ is clear. To see that $S \circ T \simeq \mathrm{id}_{(2)}$ we let an object of $(2)$ be given.\\
The functor $T$ maps it to the extension
\[0 \rightarrow (\mathcal O_X,\mathrm{d}) \xrightarrow{j_1} (\mathcal B_1,\nabla_{\mathcal B_1}) \xrightarrow{p_1} (\mathcal Y_X,\nabla_{\mathcal Y_X}) \rightarrow 0,\]
where we have already taken into account the identifications of $\psi_0$ and $\psi_1$.\\
If we set $\widetilde{\mathcal B}:=\mathrm{Sym}^{\bullet}_{\mathcal O_X}(\mathcal B_1)/\left(1-j_1(1)\right)$, then the functor $S$ maps the preceding extension to $(\widetilde{\mathcal B},\nabla_{\widetilde{\mathcal B}})$ together with filtration
\[(0)\subseteq (\mathcal O_X,\mathrm{d}) \subseteq (\mathcal B_1, \nabla_{\mathcal B_1}) \subseteq \mathrm{Sym}^2_{\mathcal O_X}(\mathcal B_1, \nabla_{\mathcal B_1}) \subseteq...\subseteq (\widetilde{\mathcal B}, \nabla_{\widetilde{\mathcal B}})\]
and isomorphism $\widetilde{\psi_{\bullet}}: \mathrm{gr}_{\bullet}(\widetilde{\mathcal B}, \nabla_{\widetilde{\mathcal B}}) \simeq \mathrm{Sym}^{\bullet}_{\mathcal O_X}(\mathcal Y_X,\nabla_{\mathcal Y_X})$ induced by the exact sequences
\[0 \rightarrow \mathrm{Sym}^n_{\mathcal O_X}(\mathcal B_1,\nabla_{\mathcal B_1}) \xrightarrow{j_{n+1}} \mathrm{Sym}^{n+1}_{\mathcal O_X}(\mathcal B_1,\nabla_{\mathcal B_1}) \xrightarrow{p_{n+1}} \mathrm{Sym}^{n+1}_{\mathcal O_X}(\mathcal Y_X, \nabla_{\mathcal Y_X})\rightarrow 0\]
whose maps are defined by
\[j_{n+1} (b_1 \cdot...\cdot b_n)= b_1 \cdot...\cdot b_n \cdot j_1(1) \ \ \textrm{and} \ \ p_{n+1}(b_1\cdot...\cdot b_{n+1})= p_1(b_1)\cdot...\cdot p_1(b_{n+1}).\]
For all $n \geq 0$ we define $\mathcal O_X$-linear morphisms
\[\tag{\textbf{2.4.2}} \mathrm{Sym}^n_{\mathcal O_X}\mathcal B_1 \rightarrow \mathcal B_n\]
by $b_1 \cdot...\cdot b_n \mapsto b_1 \cdot ... \cdot b_n$, where on the right side multiplication in $\mathcal B$ is meant.\\
These morphisms are horizontal (because the multiplication map for $\mathcal B$ is horizontal) and fit into commutative diagrams
\begin{equation*}
\begin{xy}
\xymatrix{
0\ar[r] &\mathrm{Sym}^n_{\mathcal O_X}(\mathcal B_1, \nabla_{\mathcal B_1})  \ar[r]^{j_{n+1}} \ar[d]& \mathrm{Sym}^{n+1}_{\mathcal O_X}(\mathcal B_1, \nabla_{\mathcal B_1})  \ar[r]^{p_{n+1} \ } \ar[d] & \mathrm{Sym}^{n+1}_{\mathcal O_X}(\mathcal Y_X, \nabla_{\mathcal Y_X}) \ar@{=}[d]\ar[r] & 0 \\
0\ar[r] & (\mathcal B_n,\nabla_{\mathcal B_n}) \ar[r] & (\mathcal B_{n+1},\nabla_{\mathcal B_{n+1}})\ar[r] & \mathrm{Sym}^{n+1}_{\mathcal O_X}(\mathcal Y_X, \nabla_{\mathcal Y_X}) \ar[r] & 0}
\end{xy}
\end{equation*}
where we identify $(\mathcal B_{n+1},\nabla_{\mathcal B_{n+1}})/(\mathcal B_n,\nabla_{\mathcal B_n})$ with $\mathrm{Sym}^{n+1}_{\mathcal O_X}(\mathcal Y_X, \nabla_{\mathcal Y_X})$ by means of $\psi_{n+1}$.\\
The maps $(2.4.2)$ are isomorphisms: use induction over $n$ and the previous commutative diagrams.\\
We thus obtain isomorphisms between the exhaustive filtrations
\[(0)\subseteq (\mathcal O_X,\mathrm{d}) \subseteq (\mathcal B_1, \nabla_{\mathcal B_1}) \subseteq \mathrm{Sym}^2_{\mathcal O_X}(\mathcal B_1, \nabla_{\mathcal B_1})...\subseteq (\widetilde{\mathcal B}, \nabla_{\widetilde{\mathcal B}}),\]
\[(0) \subseteq (\mathcal O_X,\mathrm{d}) \subseteq (\mathcal B_1, \nabla_{\mathcal B_1}) \subseteq (\mathcal B_2, \nabla_{\mathcal B_2}) \subseteq...\subseteq (\mathcal B, \nabla_{\mathcal B}),\]
compatible with $\widetilde{\psi_{\bullet}}$ and $\psi_{\bullet}$, and thus also between the $(\mathcal O_X,\mathrm{d})$-algebras $(\widetilde{\mathcal B}, \nabla_{\widetilde{\mathcal B}})$ and $(\mathcal B, \nabla_{\mathcal B})$.
\end{proof}

The preceding proof shows:
\begin{corollary}
Suppose we are given data
\[(\mathcal B,\nabla_{\mathcal B}),\]
\[(0) \subseteq (\mathcal B_0, \nabla_{\mathcal B_0}) \subseteq (\mathcal B_1, \nabla_{\mathcal B_1}) \subseteq (\mathcal B_2, \nabla_{\mathcal B_2}) \subseteq...\subseteq 
(\mathcal B, \nabla_{\mathcal B}),\]
\[\psi_{\bullet}: \mathrm{gr}_{\bullet}(\mathcal B, \nabla_{\mathcal B})\simeq \mathrm{Sym}^{\bullet}_{\mathcal O_X}(\mathcal Y_X, \nabla_{\mathcal Y_X})\]
as in the definition of the category $(2)$ in Lemma 2.4.1. We then have horizontal isomorphisms
\[\mathrm{Sym}^n_{\mathcal O_X}\mathcal B_1 \xrightarrow{\sim} \mathcal B_n,\]
\[\mathrm{Sym}^{\bullet}_{\mathcal O_X}(\mathcal B_1)/\left(1-j_1(1)\right) \xrightarrow{\sim} \mathcal B,\]
induced by the rule $b_1 \cdot...\cdot b_n \mapsto b_1 \cdot ... \cdot b_n$, where on the right side multiplication in $\mathcal B$ is meant. Here, $j_1: \mathcal O_X \rightarrow \mathcal B_1$ is the inclusion coming from $\mathcal B_0\subseteq \mathcal B_1$ via the identification $\psi_0$. \qquad \qed
\end{corollary}
\begin{remark}
As $X$ is a $\Q$-scheme one has for each $\mathcal O_X$-vector bundle $\mathcal U$ and $n\geq 0$ a canonical isomorphism
\[\mathrm{Sym}^n_{\mathcal O_X}(\mathcal U^\vee) \xrightarrow{\sim} \mathrm{Sym}^n_{\mathcal O_X}(\mathcal U)^\vee\]
induced by
\[g_1 \cdot...\cdot g_n \mapsto \{f_1 \cdot...\cdot f_n \mapsto \sum_{\sigma \in \Sigma_n}g_{\sigma(1)}(f_1)\cdot ... \cdot g_{\sigma(n)}(f_n) \}.\]
If $\mathcal U$ carries an integrable $S$-connection this isomorphism is horizontal for the naturally induced connections on $\mathrm{Sym}^n_{\mathcal O_X}(\mathcal U^\vee)$ and $\mathrm{Sym}^n_{\mathcal O_X}(\mathcal U)^\vee$.
\end{remark}

\subsection{The construction of the higher logarithm sheaves}
By considering the first infinitesimal restriction to $X\times_S \Yr_1$ of the Poincaré quadruple $(\mathcal P,r,s,\nabla_\mathcal P)$ we have achieved in 2.3 a completely geometric construction of $\mathcal L_1$, $\varphi_1$ and the restriction of $\nabla_1$ relative $S$ (cf. Cor. 2.3.2). The higher logarithm sheaves $(\mathcal L_n, \nabla_n, \varphi_n)$ of $X/S/\Q$ are defined as the symmetric powers of $(\mathcal L_1,\nabla_1,\varphi_1)$ (cf. 1.1), and in all our later applications we will only need this approach. Nevertheless, the question naturally arises if one can recover also the data $\mathcal L_n$, $\varphi_n$ and the restriction of $\nabla_n$ relative $S$ from the Poincaré bundle on $X\times_S \Yr$, expectably from its $n$-th infinitesimal restriction $(\mathcal P_n, r_n, s_n, \nabla_{\mathcal P_n})$ on $X\times_S \Yr_n$.\\
The answer we will give subsequently can be summarized as follows. By means of the rigidification $r_n$ the zero fiber of the $\mathcal O_X$-vector bundle with integrable $S$-connection $(p_n)_*\mathcal P_n$ identifies with the structure sheaf of $\Yr_n$ as $\mathcal O_S$-module and thus disposes of a distinguished section given by $1$. We define a (unique) $\mathcal D_{X/S}$-linear isomorphism $(p_n)_*\mathcal P_n \xrightarrow{\sim} \mathcal L_n$ under which this section corresponds to $1^{(n)}=\frac{1}{n!}$ and from which we obtain the desired interpretation of the higher logarithm sheaves by the Poincaré bundle. The construction of the isomorphism uses an infinitesimal comultiplication which is in turn naturally induced by the $\natural$-$1$-structure of the $\mathbb G_{m,S}$-biextension $(\mathcal P,\nabla_\mathcal P)$.\\
The proper verification of these assertions involves the technicality of dualizing the situation, using the auxiliary results in 2.4.1 and then redualizing.\\
\newline
At first, let $n\geq 2$. As $\epsilon^{\natural}: S \rightarrow \Yr$ is a regular embedding the sheaf $\mathcal J^n/\mathcal J^{n+1}$ as $\mathcal O_S$-module naturally identifies with $\mathrm{Sym}^n_{\mathcal O_S}(\mathcal J/\mathcal J^2)$ (cf. \cite{Fu-La}, Ch. IV, Lemma 3.8 and Cor. 2.4) and hence by $(2.3.2)$ with $\mathrm{Sym}^n_{\mathcal O_S}\mathcal H$. Observing this one can apply an entirely analogous procedure as in $(2.3.8)$-$(2.3.11)$, now with the diagrams
\begin{equation*}
\begin{xy}
\xymatrix@C-0.3cm{
X \times_S \Yr_{n-1} \ar[r] \ar[d] & \Yr_{n-1} \ar[d] & & X\times_S \Yr_n \ar[r] \ar[d]^{p_n} & \Yr_n \ar[d] \\
X \times_S \Yr_n \ar[d] \ar[r] & \Yr_n \ar[d] & & X \ar[r] & S\\
X \times_S \Yr \ar[r] \ar[d] & \Yr \ar[d] & & \\
X \ar[r] & S
}
\end{xy}
\end{equation*}
and obtains an exact sequence of $\mathcal O_X$-vector bundles with integrable $S$-connection
\[\tag{\textbf{2.4.3}} 0 \rightarrow \mathrm{Sym}^n_{\mathcal O_X}\mathcal H_X \rightarrow (p_n)_*\mathcal P_n \rightarrow (p_{n-1})_*\mathcal P_{n-1}\rightarrow 0\]
which (with Lemma 2.2.5 and Prop. 2.2.6) can also be viewed as the Fourier-Mukai transformation of the following natural $\mathcal O_{\Yr}$-linear exact sequence whose terms are WIT-sheaves of index $0$:
\[0 \rightarrow \mathcal J^n/\mathcal J^{n+1} \rightarrow \mathcal O_{\Yr}/\mathcal J^{n+1} \rightarrow \mathcal O_{\Yr}/\mathcal J^n \rightarrow 0.\]
Note that the trivialization $s_n$ is essential for the construction of $(2.4.3)$.\\
Setting $\mathcal B_0:=\mathcal O_X$, equipped with exterior $S$-derivation, and for each $n\geq 1$:
\[\mathcal B_n:=((p_n)_*\mathcal P_n)^\vee,\]
equipped with the dual integrable $S$-connection of $(p_n)_*\mathcal P_n$, we get from dualizing $(2.3.11)$ and the sequences $(2.4.3)$ horizontal monomorphisms
\[\tag{\textbf{2.4.4}} \mathcal B_{n-1}:=(0) \subseteq \mathcal B_0\subseteq \mathcal B_1 \subseteq \mathcal B_2 \subseteq ...\]
and (with Rem. 2.4.3) for each $n\geq 0$ a horizontal isomorphism
\[\tag{\textbf{2.4.5}} \psi_n: \mathcal B_n / \mathcal B_{n-1} \xrightarrow{\sim} \mathrm{Sym}^n_{\mathcal O_X}(\mathcal H_X^\vee).\vspace{3mm}\]
Next, let an ordered tuple $(n,m)$ with $n,m \geq 0$ be fixed and write
\[p_{n,m;12}: X\times_S \Yr_n \times_S \Yr_m \rightarrow X\times_S \Yr_n \quad \textrm{resp.} \quad p_{n,m;13}: X\times_S \Yr_n \times_S \Yr_m \rightarrow X\times_S \Yr_m\]
for the projections and
\[\mu^{\natural}_{n,m}: \Yr_n \times_S \Yr_m \rightarrow \Yr_{n+m}\]
for the morphism naturally induced by the multiplication map $\mu^{\natural}: \Yr\times_S \Yr\rightarrow \Yr$.\\
The horizontal isomorphism of $(0.1.21)$ then induces by pullback via the canonical closed embedding
\[X\times_S \Yr_n \times_S \Yr_m \rightarrow X\times_S \Yr \times_S \Yr\]
a $\mathcal D_{X\times_S \Yr_n \times_S \Yr_m/\Yr_n \times_S \Yr_m}$-linear isomorphism
\[(\mathrm{id}_X \times \mu^{\natural}_{n,m})^*\mathcal P_{n+m} \xrightarrow{\sim} (p_{n,m;12})^*\mathcal P_n \otimes_{\mathcal O_{X\times_S \Yr_n\times_S \Yr_m}} (p_{n,m;13})^*\mathcal P_m.\]
Together with adjunction we obtain a $\mathcal D_{X\times_S \Yr_{n+m} / \Yr_{n+m}}$-linear morphism
\[\mathcal P_{n+m} \rightarrow (\mathrm{id}_X \times \mu^{\natural}_{n,m})_*\big[(p_{n,m;12})^*\mathcal P_n \otimes_{\mathcal O_{X\times_S \Yr_n\times_S \Yr_m}} (p_{n,m;13})^*\mathcal P_m \big].\]
Taking its direct image along $p_{n+m}: X\times_S \Yr_{n+m} \rightarrow X$, noting $p_{n+m} \circ (\mathrm{id}_X \times \mu^{\natural}_{n,m})=p_n \circ p_{n,m;12}$ as well as the projection formula and base change along the cartesian diagram of affine maps
\begin{equation*}
\begin{xy}
\xymatrix{
X\times_S \Yr_n\times_S \Yr_m \ar[r]^{\quad \ \ p_{n,m;13}} \ar[d]_{p_{n,m;12}} & X\times_S \Yr_m \ar[d]^{p_m} \\
X\times_S \Yr_n \ar[r]^{\qquad p_n} & X
}
\end{xy}
\end{equation*}
then yields a $\mathcal D_{X/S}$-linear morphism
\[\begin{split} (p_{n+m})_*\mathcal P_{n+m} \rightarrow (p_n)_*(p_{n,m;12})_*\big[(p_{n,m;12})^*\mathcal P_n \otimes_{\mathcal O_{X\times_S \Yr_n\times_S \Yr_m}} (p_{n,m;13})^*\mathcal P_m \big] \\
\xrightarrow{\sim} (p_n)_*\big[\mathcal P_n \otimes_{\mathcal O_{X\times_S \Yr_n}}(p_{n,m;12})_*(p_{n,m;13})^*\mathcal P_m\big] \xrightarrow{\sim} (p_n)_*\big[\mathcal P_n \otimes_{\mathcal O_{X\times_S \Yr_n}} (p_n)^*(p_m)_*\mathcal P_m\big]\\
\xrightarrow{\sim} (p_n)_*\mathcal P_n \otimes_{\mathcal O_X} (p_m)_*\mathcal P_m.\end{split}\]
We only remark that compatibility with the connections can be checked similarly as in Lemma 2.2.5 and that the resulting map of the previous chain is the same if one uses $p_{n+m}\circ(\mathrm{id}_X\times \mu^{\natural}_{n,m})=p_m\circ p_{n,m;13}$ instead of $p_{n+m} \circ (\mathrm{id}_X \times \mu^{\natural}_{n,m})=p_n \circ p_{n,m;12}$ and then proceeds analogously.\\
\newline
In sum, for a fixed ordered tuple $(n,m)$ with $n,m\geq 0$ we have constructed from the isomorphism $(0.1.21)$ - which is part of the $\natural$-$1$-structure on $\mathcal P$ (cf. 0.1.3) - in a canonical way a horizontal map
\[\tag{\textbf{2.4.6}} \xi_{(n,m)}:(p_{n+m})_*\mathcal P_{n+m} \rightarrow (p_n)_*\mathcal P_n \otimes_{\mathcal O_X} (p_m)_*\mathcal P_m.\]
The fact that $\mathcal P$ is a commutative $\mathbb G_{m,X}$-extension of $X\times_S \Yr$, i.e. the corresponding commutative diagrams $(0.1.8)$ and $(0.1.9)$, implies the commutativity of the diagrams
\begin{equation*}
\begin{xy}
\xymatrix{
(p_{n+m})_*\mathcal P_{n+m}\ar[rr]^{\xi_{(n,m)} \quad \ \ } \ar[d]_{\id} & & (p_n)_*\mathcal P_n \otimes_{\mathcal O_X} (p_m)_*\mathcal P_m \ar[d]^{\mathrm{can}}\\
(p_{m+n})_*\mathcal P_{m+n}\ar[rr]^{\xi_{(m,n)} \quad \ \ } & & (p_m)_*\mathcal P_m \otimes_{\mathcal O_X} (p_n)_*\mathcal P_n
}
\end{xy}
\end{equation*}
and
\begin{equation*}
\begin{xy}
\xymatrix{
(p_{n+m+l})_*\mathcal P_{n+m+l}\ar[rr]^{\xi_{(n+m,l)} \ \ \ \ } \ar[d]_{\xi_{(n,m+l)}} & & (p_{n+m})_*\mathcal P_{n+m} \otimes_{\mathcal O_X} (p_l)_*\mathcal P_l \ar[d]^{\xi_{(n,m)\otimes \id}}\\
(p_n)_*\mathcal P_n \otimes_{\mathcal O_X} (p_{m+l})_*\mathcal P_{m+l}\ar[rr]^{\id\otimes\xi_{(m,l)} \quad \ \ \ \ } & & (p_n)_*\mathcal P_n \otimes_{\mathcal O_X} (p_m)_*\mathcal P_m \otimes_{\mathcal O_X} (p_l)_*\mathcal P_l
}
\end{xy}
\end{equation*}
The morphisms in $(2.4.6)$ are compatible with the transition maps arising from the projections in $(2.4.3)$ and for $m=0$ become the identity on $(p_n)_*\mathcal P_n$ under the $\mathcal D_{X/S}$-linear trivialization $s$ of $\mathcal P_0$.\\
\newline
Let $\mathcal B$ be the quasi-coherent $\mathcal O_X$-module with integrable $S$-connection defined as the direct limit over the $\mathcal B_n=((p_n)_*\mathcal P_n)^\vee$ in $(2.4.4)$. By dualizing the maps of $(2.4.6)$ it becomes (according to the previous remarks) a well-defined associative, commutative and unital quasi-coherent $(\mathcal O_X,\mathrm{d})$-algebra with integrable $S$-connection.\\
Moreover, the arrow $(2.4.6)$ induces on the subsheaf $\mathrm{Sym}^{n+m}_{\mathcal O_X}\mathcal H_X$ of $(p_{n+m})_*\mathcal P_{n+m}$ a map
\[\mathrm{Sym}^{n+m}_{\mathcal O_X}\mathcal H_X \rightarrow \mathrm{Sym}^n_{\mathcal O_X}\mathcal H_X \otimes_{\mathcal O_X} \mathrm{Sym}^m_{\mathcal O_X}\mathcal H_X\]
which under the identification of Rem. 2.4.3 equals precisely the dual of the morphism
\[\mathrm{Sym}^n_{\mathcal O_X}(\mathcal H_X^\vee) \otimes_{\mathcal O_X} \mathrm{Sym}^m_{\mathcal O_X}(\mathcal H_X^\vee) \rightarrow \mathrm{Sym}^{n+m}_{\mathcal O_X}(\mathcal H_X^\vee)\]
given by multiplication in symmetric powers; one can verify this by calculating on the one hand the dual of the previous multiplication map under the isomorphism of Rem. 2.4.3 and by using on the other hand the definition of $(2.4.6)$ and the infinitesimal group law of $\Yr$ to check the claimed equality.\\
Altogether, we see that $\mathcal B$ together with the filtration
\[(0) \subseteq \mathcal B_0 \subseteq \mathcal B_1 \subseteq \mathcal B_2 \subseteq ... \subseteq \mathcal B\]
induced by $(2.4.4)$ and the isomorphism
\[\psi_{\bullet}:\mathrm{gr}_{\bullet}\mathcal B \xrightarrow{\sim} \mathrm{Sym}^{\bullet}_{\mathcal O_X}(\mathcal H_X^\vee)\]
induced by $(2.4.5)$ provides data as in the hypothesis of Cor. 2.4.2. We thus conclude that the map
\[\mathrm{Sym}^n_{\mathcal O_X}\mathcal B_1 \rightarrow \mathcal B_n\]
defined by multiplication in $\mathcal B$ is a $\mathcal D_{X/S}$-linear isomorphism. With the isomorphism obtained from the previous map by dualizing and using Rem. 2.4.3 one can verify that the following composition in which the first arrow comes from the maps $(2.4.6)$ is a $\mathcal D_{X/S}$-linear isomorphism:
\[\tag{\textbf{2.4.7}} (p_n)_*\mathcal P_n \rightarrow (p_1)_*\mathcal P_1 \otimes_{\mathcal O_X} ... \otimes_{\mathcal O_X} (p_1)_*\mathcal P_1 \xrightarrow{\frac{1}{n!}\cdot \mathrm{can}} \mathrm{Sym}^n_{\mathcal O_X}((p_1)_*\mathcal P_1).\]
Recall from Cor. 2.3.2 and Def. 1.1.1 that the right side with its integrable $S$-connection is the $n$-th logarithm sheaf $\mathcal L_n$ of $X/S/\Q$ with $\nabla_n$ restricted relative $S$. We have the $\mathcal O_S$-linear identification
\[\tag{\textbf{2.4.8}} \epsilon^*(p_n)_*\mathcal P_n \xrightarrow{\sim} (\pi^\natural_n)_*(\epsilon\times \mathrm{id}_{\Yr_n})^*\mathcal P_n \simeq (\pi^\natural_n)_*\mathcal O_{\Yr_n},\]
where the first map is the base change isomorphism along the cartesian diagram of affine maps
\begin{equation*}
\begin{xy}
\xymatrix{
\Yr_n \ar[r]^{\pi^{\natural}_n} \ar[d]_{\epsilon \times \mathrm{id}_{\Yr_n}} & S \ar[d]^{\epsilon} \\
X\times_S \Yr_n \ar[r]^{\quad \ p_n} & X}
\end{xy}
\end{equation*}
and the second is induced by the trivialization
\[r_n: (\epsilon\times \mathrm{id}_{\Yr_n})^*\mathcal P_n \simeq \mathcal O_{\Yr_n}.\]
With the identifications induced by $(2.4.8)$ the pullback of the isomorphism $(2.4.7)$ along $\epsilon$ is checked to translate into the $\mathcal O_S$-linear isomorphism
\[\tag{\textbf{2.4.9}} \mathcal O_{\Yr}/\mathcal J^{n+1}\rightarrow \mathcal O_{\Yr}/\mathcal J^2 \otimes_{\mathcal O_S} ... \otimes_{\mathcal O_S} \mathcal O_{\Yr}/\mathcal J^2 \xrightarrow{\frac{1}{n!}\cdot \mathrm{can}} \mathrm{Sym}^n_{\mathcal O_S}(\mathcal O_{\Yr}/\mathcal J^2),\]
where the first arrow is induced by $n$-fold multiplication $\Yr_1\times_S...\times_S \Yr_1 \rightarrow \Yr_n$.
The splitting $\varphi_n$ of $\epsilon^*\mathcal L_n$ is induced by identifying the $\mathcal O_S$-module $\mathcal O_{\Yr}/\mathcal J^2$ in $\mathrm{Sym}^n_{\mathcal O_S}(\mathcal O_{\Yr}/\mathcal J^2)$ with $\mathcal O_S \oplus \mathcal H$.

\begin{proposition}
The isomorphism $(2.4.7)$ is the unique $\mathcal D_{X/S}$-linear isomorphism
\[(p_n)_*\mathcal P_n \xrightarrow{\sim} \mathcal L_n\]
whose pullback along $\epsilon$ maps $1 \in \Gamma(\Yr_n,\mathcal O_{\Yr_n}) = \Gamma(S,(\pi^\natural_n)_*\mathcal O_{\Yr_n}) \simeq \Gamma(S,\epsilon^*(p_n)_*\mathcal P_n)$ to $1^{(n)}=\frac{1}{n!}$.\\
Under this identification the splitting
\[\varphi_n: \epsilon^*\mathcal L_n \simeq \prod_{k=0}^n\mathrm{Sym}^k_{\mathcal O_S}\mathcal H\]
corresponds to the composition of $(2.4.8)$ with $(2.4.9)$, where one observes $\mathcal O_{\Yr}/\mathcal J^2 \simeq \mathcal O_S \oplus \mathcal H$ as $\mathcal O_S$-modules (cf. $(2.3.6)$):
\[\epsilon^*(p_n)_*\mathcal P_n \simeq \mathcal O_{\Yr}/\mathcal J^{n+1}\xrightarrow{\sim} \mathrm{Sym}^n_{\mathcal O_S}(\mathcal O_{\Yr}/\mathcal J^2) \simeq \prod_{k=0}^n\mathrm{Sym}^k_{\mathcal O_S}\mathcal H.\]
Finally, the transition maps of the $\mathcal L_n$ correspond to those of the $(p_n)_*\mathcal P_n$ as given in $(2.4.3)$.
\end{proposition}
\begin{proof}
It only remains to check the uniqueness claim and the statement about the transition maps. The first follows from the observation that a $\mathcal D_{X/S}$-linear automorphism of $\mathcal L_n$ which in the zero fiber maps $1^{(n)}$ to itself is the identity because of the isomorphism $\Hom_{\mathcal D_{X/S}}(\mathcal L_n,\mathcal L_n) \xrightarrow{\sim} \Gamma(S,\epsilon^*\mathcal L_n)$ induced by $(1.3.5)$. We finally prove the assertion about the transition maps.\\
By what we already remarked about the morphisms $(2.4.6)$ we know that we have a commutative diagram
\begin{equation*}
\begin{xy}
\xymatrix@C-0.3cm{
(p_n)_*\mathcal P_n \ar[r] \ar[d] & (p_1)_*\mathcal P_1 \otimes_{\mathcal O_X} ... \otimes_{\mathcal O_X} (p_1)_*\mathcal P_1 \ar[d]\\
(p_{n-1})_*\mathcal P_{n-1} \ar[r] & (p_1)_*\mathcal P_1 \otimes_{\mathcal O_X} ... \otimes_{\mathcal O_X} (p_1)_*\mathcal P_1 }
\end{xy}
\end{equation*}
in which the left vertical arrow is the transition map in $(2.4.3)$, the horizontal arrows are given by $(2.4.6)$ and the right vertical arrow is induced by tensoring $(n-1)$-fold the identity map on $(p_1)_*\mathcal P_1$ with the transition morphism $\mathrm{pr}: (p_1)_*\mathcal P_1 \rightarrow \mathcal O_X$ in the $i$-th component, where $1\leq i \leq n$. This implies that the diagram is still commutative if we replace the right vertical arrow by the map
\[\frac{1}{n}\cdot (\mathrm{pr} \otimes \id \otimes ... \otimes \id +...+ \id \otimes ... \otimes \id \otimes \mathrm{pr}).\]
But this map obviously sits in a commutative diagram
\begin{equation*}
\begin{xy}
\xymatrix{
(p_1)_*\mathcal P_1 \otimes_{\mathcal O_X} ... \otimes_{\mathcal O_X} (p_1)_*\mathcal P_1 \ar[d] \ar[rr]^{\quad \frac{1}{n!}\cdot \mathrm{can}} & & \mathrm{Sym}^n_{\mathcal O_X} ((p_1)_*\mathcal P_1) =\mathcal L_n\ar[d] \\
(p_1)_*\mathcal P_1 \otimes_{\mathcal O_X} ... \otimes_{\mathcal O_X} (p_1)_*\mathcal P_1 \ar[rr]^{\quad \frac{1}{(n-1)!}\cdot \mathrm{can}} & & \mathrm{Sym}^{n-1}_{\mathcal O_X} ((p_1)_*\mathcal P_1)=\mathcal L_{n-1}}
\end{xy}
\end{equation*}
with right vertical arrow given by the usual transition morphism of the logarithm sheaves (cf. $(1.1.3)$). Recalling the definition of the maps $(2.4.7)$ the remaining claim obviously follows.
\end{proof}

Hence, if we prolong the integrable $S$-connection of $(p_n)_*\mathcal P_n$ to the integrable $\Q$-connection provided by $\nabla_n$ via the isomorphism in Prop. 2.4.4, the sheaf $(p_n)_*\mathcal P_n$ becomes an object of $U_n(X/S/\Q)$, where the filtration $A^i((p_n)_*\mathcal P_n):=\ker((p_n)_*\mathcal P_n \rightarrow (p_{i-1})_*\mathcal P_{i-1})$ corresponds to the canonical filtration on $\mathcal L_n$ (cf. also Rem. 1.1.6). Together with the section $1 \in \Gamma(\Yr_n,\mathcal O_{\Yr_n}) = \Gamma(S,(\pi^\natural_n)_*\mathcal O_{\Yr_n}) \simeq \Gamma(S,\epsilon^*(p_n)_*\mathcal P_n)$ it becomes the $n$-th logarithm sheaf of $X/S/\Q$ as characterized at the end of 1.3.2. The mentioned prolongation is the unique one with this property.\\
Via the $\mathcal O_S$-linear isomorphism induced by $(2.4.9)$ and the usual decomposition $\mathcal O_{\Yr}/\mathcal J^2 \simeq \mathcal O_S \oplus \mathcal H$:
\[\tag{\textbf{2.4.10}} \mathcal O_{\Yr}/\mathcal J^{n+1} \xrightarrow{\sim} \mathrm{Sym}^n_{\mathcal O_S}(\mathcal O_{\Yr}/\mathcal J^2) \simeq \prod_{k=0}^n \mathrm{Sym}^k_{\mathcal O_S} \mathcal H\]
we obtain on $\prod_{k=0}^n \mathrm{Sym}^k_{\mathcal O_S} \mathcal H$ the structure of a finite locally free $\mathcal O_S$-algebra: this structure is precisely the one introduced in $(1.3.9)$ and $(1.3.10)$, as one can check explicitly. The $\pi^*\big(\prod_{k=0}^n \mathrm{Sym}^k_{\mathcal O_S} \mathcal H\big)$-module structure on $(p_n)_*\mathcal P_n$, given by $\mathcal O_{X\times_S \Yr_n}$-multiplication on $\mathcal P_n$ and $(2.4.10)$, becomes under the isomorphism of Prop. 2.4.4 the $\pi^*\big(\prod_{k=0}^n \mathrm{Sym}^k_{\mathcal O_S} \mathcal H\big)$-module structure on $\mathcal L_n$ defined in 1.3.3.\\
In this interpretation the morphism in the pro-category of $U(X/S/\Q)$ induced by the above maps
\[\xi_{(n,m)}:(p_{n+m})_*\mathcal P_{n+m} \rightarrow (p_n)_*\mathcal P_n \otimes_{\mathcal O_X} (p_m)_*\mathcal P_m\]
is nothing else than the comultiplication on the projective system of the logarithm sheaves, defined as in \cite{Be-Le}, 1.2.10 (i), via their universal property $(1.3.5)$.

\markright{\uppercase{The logarithm sheaves and the Poincaré bundle}}
\section{The invariance under isogenies and the Poincaré bundle}
\markright{\uppercase{The logarithm sheaves and the Poincaré bundle}}
\subsection{The transpose endomorphism}
We begin with some general theory and recall how a homomorphism of abelian schemes induces (in the other direction) a homomorphism between the respective universal vectorial extensions of the duals. As in the end only the case of $N$-multiplication will be relevant for us we here focus on endomorphisms, which slightly reduces notation (for the analogous general case cf. Rem. 2.5.3 (i)).\\
\newline
Let $S$ be an arbitrary locally noetherian scheme and $X/S$ an abelian scheme, $\Yr/S$ the universal vectorial extension of its dual abelian scheme $Y/S$ and $(\mathcal P,r,\nabla_\mathcal P)$ the $\Yr$-rigidified Poincaré bundle on $X\times_S\Yr$ with its universal integrable $\Yr$-connection (cf. 0.1.1).\\
\newline
Let $u: X \rightarrow X$ be an endomorphism of $X$ over $S$.\\
Consider the endomorphism $u\times \mathrm{id}_{\Yr}$ of the abelian $\Yr$-scheme $X\times_S \Yr$ and the pullback line bundle $(u\times \mathrm{id}_{\Yr})^*\mathcal P$: it is algebraically equivalent to zero with respect to $\Yr$, as one easily deduces from the corresponding fact about $\mathcal P$. Moreover, equip this line bundle with the $\Yr$-rigidification $r_u$ naturally induced by $r$ and with the integrable $\Yr$-connection $(\nabla_\mathcal P)_u$ arising from $\nabla_\mathcal P$ by pullback.\\
\newline
Then, by definition of $\Yr$, there exists a unique $S$-morphism
\[\tag{\textbf{2.5.1}} u^\natural: \Yr\rightarrow \Yr\]
such that
\[\tag{\textbf{2.5.2}} (\mathrm{id}_X\times u^\natural)^*(\mathcal P,r,\nabla_\mathcal P) \simeq ((u\times \mathrm{id}_{\Yr})^*\mathcal P, r_u, (\nabla_\mathcal P)_u),\]
where the left side is equipped with the $\Yr$-rigidification and integrable $\Yr$-connection given by pullback of $(\mathcal P,r,\nabla_\mathcal P)$ along the cartesian diagram
\begin{equation*}
\begin{xy}
\xymatrix{
X\times_S \Yr \ar[r]^{ \ id_X \times u^\natural} \ar[d] & X\times_S \Yr \ar[d]_{q}\\
\Yr \ar[r]^{u^\natural} & \Yr}
\end{xy}
\end{equation*}
The $S$-morphism $u^\natural: \Yr \rightarrow \Yr$ of $(2.5.1)$ is an endomorphism and called the \underline{transpose endomorphism} of $u$. For an $S$-scheme $T$ it is given on $T$-rational points by sending the isomorphism class of a triple $(\mathcal L, \alpha, \nabla_{\mathcal L}) \in \mathrm{Pic}^\natural(X \times_S T/T)$ to the class of the pullback $(u\times \mathrm{id}_T)^*\mathcal L$ which is endowed with $T$-rigidification and integrable $T$-connection as performed above for $(u\times \mathrm{id}_{\Yr})^*\mathcal P$.\\
Recall also (cf. Lemma 0.1.5) that the isomorphism of $(2.5.2)$ is unique.\\
\newline
Working with $Y$ and $(\mathcal P^0,r^0)$ instead of $\Yr$ and $(\mathcal P,r,\nabla_\mathcal P)$ one defines in a completely analogous way the transpose endomorphism $Y \rightarrow Y$ associated with $u$.\\
\newline
To conclude this general part, we determine the transpose endomorphism in the important case of multiplication by integers. The result is what one would expect, its proof however not entirely trivial.

\begin{proposition}
For the $N$-multiplication endomorphism $[N]: X \rightarrow X$, where $N$ is an integer, the transpose endomorphism $[N]^\natural: \Yr\rightarrow \Yr$ equals the $N$-multiplication map of the $S$-group scheme $\Yr$. The same statement holds if $\Yr$ is replaced by $Y$.
\end{proposition}
\begin{proof}
At first, let $(\mathcal L,\alpha)$ be any $S$-rigidified line bundle on $X$ which is algebraically equivalent to zero with respect to $S$. The existence of $\alpha$ and the definition of algebraic equivalence to zero imply that there is some isomorphism of line bundles on $X\times_S X$:
\[\tag{\textbf{2.5.3}} \mu^*\mathcal L \otimes_{\mathcal O_{X\times_S X}} \mathrm{pr}_1^*\mathcal L^{-1} \otimes_{\mathcal O_{X\times_S X}} \mathrm{pr}_2^*\mathcal L^{-1} \simeq \mathcal O_{X\times_S X},\]
where $\mu, \mathrm{pr}_1, \mathrm{pr}_2: X\times_S X \rightarrow X$ denote the multiplication map resp. the two projections.\\
By restricting $(2.5.3)$ via the morphism $\id_X \times (-1)_X: X \rightarrow X\times_S X$, where $(-1)_X$ is the inverse map of $X/S$, and by taking into account the $S$-rigidification $\alpha$, one obtains an isomorphism
\[\tag{\textbf{2.5.4}} (-1)_X^*\mathcal L \simeq \mathcal L^{-1}.\]
In the following, if $\mathcal M$ and $\mathcal N$ are two line bundles on $X$ which are isomorphic up to a tensor factor given by the pullback of a line bundle on $S$, we will write $\mathcal M \simeq_S \mathcal N$.\\
Combining $(2.5.4)$ with the addition formula for $N$-multiplication:
\[[N]^*\mathcal L \simeq_S \mathcal L^{\otimes \frac{N^2+N}{2}}\otimes_{\mathcal O_X} (-1)_X^*\mathcal L^{\otimes \frac{N^2-N}{2}},\]
which is an easy consequence of the theorem of cube (cf. \cite{Ch-Fa}, Thm. 1.3), yields
\[\tag{\textbf{2.5.5}} [N]^*\mathcal L \simeq_S \mathcal L^{\otimes N}.\]
Equip each side with the $S$-rigidification naturally induced by $\alpha$. The existence of such rigidifications and $(2.5.5)$ imply that there is some isomorphism
\[\tag{\textbf{2.5.6}} [N]^*\mathcal L \simeq \mathcal L^{\otimes N}\]
which moreover can be chosen uniquely such that it respects these $S$-rigidifications.\footnote{The argument is the same as in footnote 6 of Chapter 0.}\\
This implies that the transpose endomorphism $Y \rightarrow Y$ of $[N]: X\rightarrow X$ evaluated in $S$-rational points equals the $N$-multiplication map $Y(S) \rightarrow Y(S)$. As the previous arguments work unalteredly in the situation $X\times_S T/T$ (with an $S$-scheme $T$) the second claim of the proposition follows.\\
Now assume in addition that $\mathcal L$ is endowed with an integrable $S$-connection $\nabla_{\mathcal L}$, i.e. $(\mathcal L,\alpha,\nabla_{\mathcal L})$ defines a class in $\mathrm{Pic}^\natural(X/S)$. Equip the line bundles in $(2.5.6)$ with the integrable $S$-connections given by pullback resp. by tensor product. If we can show that $(2.5.6)$ is horizontal for these connections, then we have verified that $[N]^\natural(S): \Yr(S)\rightarrow \Yr(S)$ is the $N$-multiplication map. As everything we do works equally in the situation $X\times_S T/T$ (with an $S$-scheme $T$) the remaining first claim of the proposition follows.\\
The horizontality of $(2.5.6)$ can be shown with the same trick as used in \cite{Maz-Mes}, proof of Prop. $(4.2.1)$ - in our situation applied to the isomorphism $(2.5.6)$. Only note the following:  to make the argument of \cite{Maz-Mes} work we need to see that also in our situation "$i(\nabla)$ depends only on $\mathcal L$ and not on the integrable connection $\nabla$ chosen" (cf. ibid., p. 49). As in \cite{Maz-Mes} we want to deduce this fact from ibid., Lemma $(3.2.6)$ - in the completely analogous way as is illustrated there after the proof of that lemma. The reasoning given there carries over to our situation, i.e. the equality "$\delta(\bar{\nabla})=\delta(\bar{\bar{\nabla}})$" (cf. ibid., p. 42) holds also in our case, because the morphism $[N]: X\rightarrow X$ acts as $N$-multiplication on invariant differential forms.
\end{proof}

\subsection{Interpretation of the invariance property}
We resume our familiar geometric setting $X/S/\Q$ and show that in our geometric interpretation of the logarithm sheaves via the Poincaré bundle their invariance under (endomorphic) isogenies is encoded in the symmetry isomorphism $(2.5.2)$. We explicate full details only for the case which will be concretely needed in the future, namely that of an isogeny endomorphism and the first logarithm sheaf. A brief outline of the (completely analogous) general case is included in a final remark.\\
\newline
Let $u: X\rightarrow X$ be an isogeny (cf. footnote 11 of Chapter 1) and $u^\natural: \Yr\rightarrow \Yr$ its transpose.\\
Write
\[u^\natural_1: \Yr_1 \rightarrow \Yr_1\]
for the morphism naturally induced by $u^\natural$, where as usual $\Yr_1$ is the first infinitesimal neighborhood of the zero section of $\Yr/S$ with associated closed immersion $\epsilon^{\natural}_1: \Yr_1 \rightarrow \Yr$ (cf. Def. 2.2.14).\\
We will abbreviate with $\iota_1: X\times_S \Yr_1\rightarrow X\times_S \Yr$ the closed immersion $\id_X \times \epsilon^{\natural}_1$.\\
\newline
Recall from Cor. 2.3.2 how the first logarithm sheaf $\mathcal L_1$ of $X/S/\Q$ was constructed as $(p_1)_*\mathcal P_1$, where $\mathcal P_1$ is the restriction of $\mathcal P$ via $\iota_1$ and $p_1: X \times_S \Yr_1 \rightarrow X$ the projection.\\
We have the following chain of canonical $\mathcal O_X$-linear isomorphisms:
\[u^*\mathcal L_1 = u^*(p_1)_*\mathcal P_1 \simeq (p_1)_*(u\times \mathrm{id}_{\Yr_1})^*\mathcal P_1\simeq (p_1)_*\iota_1^*(u\times \mathrm{id}_{\Yr})^*\mathcal P \simeq^{(!)} (p_1)_*\iota_1^*(\mathrm{id}_X\times u^\natural)^*\mathcal P \]
\[\simeq (p_1)_*(\mathrm{id}_X\times u^\natural_1)^*\mathcal P_1 \simeq (p_1)_*(\mathrm{id}_X\times u^\natural_1)_*(\mathrm{id}_X\times u^\natural_1)^*\mathcal P_1 \simeq (p_1)_*[(\mathrm{id}_X\times u^\natural_1)_* \mathcal O_{X\times_S \Yr_1} \otimes_{\mathcal O_{X\times_S \Yr_1}}\mathcal P_1].\]
The isomorphisms after the equality all are easy standard identifications coming from flat base change, commutative diagrams and the projection formula - except for the one decorated with the exclamation mark: this is the crucial identification of $(2.5.2)$.\\
The map on structure sheaves
\[\tag{\textbf{2.5.7}} \mathcal O_{X\times_S \Yr_1}\rightarrow (\mathrm{id}_X\times u^\natural_1)_* \mathcal O_{X\times_S \Yr_1}\]
induced by
\[\id_X\times u^\natural_1: X \times_S \Yr_1 \rightarrow X\times_S \Yr_1\]
is an isomorphism because the corresponding morphism of quasi-coherent $\mathcal O_X$-algebras is the map
\[\id \oplus (u_{{\mathrm{dR}}}^*)^\vee_X: \mathcal O_X \oplus \mathcal H_X \rightarrow \mathcal O_X \oplus \mathcal H_X\]
and thus an isomorphism. Here, we denote as in 1.4.2 with $u_{{\mathrm{dR}}}^*: H^1_{\mathrm{dR}}(X/S) \rightarrow H^1_{\mathrm{dR}}(X/S)$ the map on de Rham cohomology induced by $u$, which is an isomorphism (cf. the proof of Thm. 1.4.2), and with $(u_{{\mathrm{dR}}}^*)^\vee_X: \mathcal H_X \rightarrow \mathcal H_X$ the morphism obtained from it by dualizing and pullback to $X$.\\ 
With the identification $(2.5.7)$ we can conclude the above chain of isomorphisms and in sum receive
\[\tag{\textbf{2.5.8}} u^*\mathcal L_1=u^*(p_1)_*\mathcal P_1 \simeq......\simeq (p_1)_*\mathcal P_1 = \mathcal L_1.\]
One checks that $(2.5.8)$ respects the integrable $S$-connections of both sides, where $u^*\mathcal L_1$ carries the pullback connection (apart from standard verifications one observes in particular the horizontality of $(2.5.2)$). It is a further routine affair to see that under pullback via $\epsilon$ it maps the global $S$-section $1^{(1)}=1$ of $\epsilon^*u^*\mathcal L_1 \simeq \epsilon^*\mathcal L_1 \simeq \mathcal O_S \oplus \mathcal H$ to the global $S$-section $1^{(1)}=1$ of $\epsilon^*\mathcal L_1 \simeq \mathcal O_S \oplus \mathcal H$ (for this one has to recall how the splitting of $\mathcal L_1=(p_1)_*\mathcal P_1$ was defined explicitly in terms of the rigidification $r_1$ of $\mathcal P_1$, cf. $(2.3.13)$, and in particular observe the compatibility of $(2.5.2)$ with the occurring rigidifications). We have thus defined in $(2.5.8)$ a $\mathcal D_{X/S}$-linear isomorphism
\[\tag{\textbf{2.5.9}} \mathcal L_1 \xrightarrow{\sim} u^*\mathcal L_1\]
which in the zero fiber sends $1^{(1)}$ to $1^{(1)}$, and by Thm. 1.3.6 (ii) we know that it is in fact $\mathcal D_{X/\Q}$-linear. By comparison of the associated maps in the zero fiber we conclude:
\begin{proposition}
The $\mathcal D_{X/\Q}$-linear isomorphism constructed in $(2.5.9)$ is the invariance isomorphism of Thm. 1.4.2. \qquad \qed
\end{proposition}
Leaving away the unessential formal ballast in the chain of isomorphisms $(2.5.8)$ we thus see that after interpreting the logarithm sheaf $\mathcal L_1$ geometrically via the Poincaré bundle its invariance under (endomorphic) isogenies expresses as the isomorphism
\[(\mathrm{id}_X\times u^\natural)^*(\mathcal P,r,\nabla_\mathcal P) \simeq ((u\times \mathrm{id}_{\Yr})^*\mathcal P, r_u, (\nabla_\mathcal P)_u)\]
of $(2.5.2)$. This is the crucial insight of the above discussion.\\
\newline
Let us finally indicate some obvious generalizations of the preceding line of arguments:

\begin{remark}
(i) If $u: X\rightarrow X'$ is a homomorphism of abelian schemes over a locally noetherian base $S$ its transpose homomorphism $u^{\natural}: (Y')^{\natural} \rightarrow \Yr$ is characterized by the existence of a (unique) isomorphism
\[\tag{\textbf{2.5.10}} (\mathrm{id}_X \times u^{\natural})^*(\mathcal P,r,\nabla_\mathcal P) \simeq ((u\times \mathrm{id}_{(Y')^{\natural}})^*\mathcal P',r'_u,(\nabla_{\mathcal P'})_u),\]
where $r'_u$ resp. $(\nabla_{\mathcal P'})_u$ is the $(Y')^{\natural}$-rigidification resp. integrable $(Y')^{\natural}$-connection induced by $r'$ resp. $\nabla_{\mathcal P'}$. In $T$-rational points (with an $S$-scheme $T$) it is given by pulling back the class of a triple $(\mathcal L', \alpha', \nabla_{\mathcal L'})$ in $\mathrm{Pic}^{\natural}(X'\times_S T/T)$ via $u\times \id_T: X \times_S T \rightarrow X'\times_S T$.\vspace{1mm}\\
(ii) Resuming the usual situation $S/\Q$ one can then construct for an isogeny of abelian $S$-schemes $u: X \rightarrow X'$ in an entirely analogous way (but with more notational demand) a chain of isomorphisms as in $(2.5.8)$ whose crucial part is the identification $(2.5.10)$ and obtains a $\mathcal D_{X/\Q}$-linear isomorphism
\[\mathcal L_1 \xrightarrow{\sim} u^*\mathcal L_1'\]
which is precisely the (general) invariance isomorphism of Thm. 1.4.2.\vspace{1mm}\\
(iii) Finally, in the situation of (ii), when working with the $n$-th infinitesimal neighborhoods $\Yr_n$ resp. $(Y')^{\natural}_n$ for $n\geq 1$ one deduces from $(2.5.10)$ in the same way as before a $\mathcal D_{X/S}$-linear isomorphism
\[u^*(p'_n)_*\mathcal P'_n \simeq (p_n)_*\mathcal P_n\]
which in the zero fiber sends $1$ to $1$ in the induced isomorphism of $\mathcal O_S$-modules $\mathcal O_{(Y')^{\natural}_n}\simeq \mathcal O_{\Yr_n}$ (recall the identification $(2.4.8)$). If we identify $(p'_n)_*\mathcal P'_n \simeq \mathcal L'_n$ resp. $(p_n)_*\mathcal P_n \simeq \mathcal L_n$ as in Prop. 2.4.4, then we obtain a $\mathcal D_{X/S}$-linear isomorphism
\[\mathcal L_n \xrightarrow{\sim} u^*\mathcal L_n'\]
which sends $1^{(n)}=\frac{1}{n!}$ to $(1^{(n)})'=\frac{1}{n!}$ in the zero fiber, hence is $\mathcal D_{X/\Q}$-linear and coincides with the invariance isomorphism defined in Cor. 1.4.3.
\end{remark}

\section{Motivic description of the first logarithm extension class}
\markright{\uppercase{The logarithm sheaves and the Poincaré bundle}}
\subsection{Generalities on $1$-motives over a scheme}
We begin by recollecting fundamental definitions and facts concerning $1$-motives over a base scheme, fixing at the same time notation used in the further progress.\\
As main sources for our exposition we will use \cite{BaVi}, 2.2, \cite{An-BaVi}, 1.1-1.2 and \cite{Ber}, 2 and 4.\\
The theory goes back to the seminal work \cite{De2}, 10, which is focused on $1$-motives over an algebraically closed field but also contains the definition over an arbitrary scheme (cf. ibid., $(10.1.10)$). 
\subsubsection{Basic definitions}
Fix a locally noetherian base scheme $S$. The category of commutative $S$-group schemes will be tacitly viewed as a full subcategory of the category of abelian $fppf$-sheaves over $S$.\\
\newline
By definition, a (smooth) $1$-motive $M=[X \xrightarrow{u} G]$ over $S$ is the datum of commutative $S$-group schemes $X,G$ and a homomorphism $u$ between them, where we require that:\\
(i) $X$ is a lattice over $S$, i.e. étale locally on $S$ isomorphic to the constant $S$-group scheme defined by a finitely generated free abelian group;\\
(ii) $G$ is a semi-abelian scheme over $S$, i.e. an extension
\[0\rightarrow T \rightarrow G \rightarrow A \rightarrow 0\]
of an abelian scheme $A$ by a torus $T$.\footnote{By a torus $T$ we mean an $S$-group scheme which is étale locally on $S$ isomorphic to finitely many copies of $\mathbb G_{m,S}$.}\\
\newline
Note that any semi-abelian scheme (and a fortiori: any abelian scheme or torus) over $S$ naturally becomes a $1$-motive over $S$ by setting $X=0$; the analogous comment applies to a lattice over $S$.\\
\newline
A morphism $M_1 \rightarrow M_2$ between two $1$-motives over $S$ is a commutative diagram of homomorphisms
\begin{equation*}
\begin{xy}
\xymatrix{
X_1 \ar[r]^{u_1} \ar[d] & G_1 \ar[d]\\
X_2 \ar[r]^{u_2} & G_2}
\end{xy}
\end{equation*}
We remark that the so-defined morphisms respect the extension structures of $G_1$ and $G_2$: the reason for this is that any homomorphism of a torus into an abelian scheme is zero (cf. \cite{Berto}, Lemma 1.2.1).\\
\newline
The obtained category of $1$-motives over $S$ is additive but not abelian. By considering a $1$-motive $M=[X \xrightarrow{u} G]$ as complex in degree $-1$ and $0$ we may view the category of $1$-motives as a full subcategory of the bounded complexes of abelian $fppf$-sheaves over $S$, and we still obtain a full embedding when further passing to the bounded derived category $D^b(S_{fppf})$ of these sheaves (cf. \cite{BaVi}, Scholium 2.2.4). These viewpoints will be adopted several times in what follows.

\subsubsection{Universal vectorial extension and de Rham realization}
By an extension $[X \xrightarrow{v} E]$ of a $1$-motive $M=[X \xrightarrow{u} G]$ by a vector group $W$ over $S$ we mean an extension
\[0\rightarrow W \rightarrow E \rightarrow G \rightarrow 0\]
of $G$ by $W$ together with a homomorphism $v:X \rightarrow E$ which lifts $u$.\\
\newline
A universal (vectorial) extension $E(M)=[X \xrightarrow{v} E(M)_G]$ of $M$ is then defined to be an extension of $M$ by an $S$-vector group $V(M)$ which parametrizes all extensions of $M$ by $S$-vector groups $W$ via pushout along a unique vector group homomorphism $V(M) \rightarrow W$.\\
\newline
A universal extension exists for every $1$-motive $M$ and is determined up to canonical isomorphism (cf. \cite{An-BaVi}, 2.2-2.3, or \cite{Ber}, 2, for a proof).\\
We only record the following two special cases: if $M$ is an abelian scheme $A$ over $S$, then its universal extension is given by the $S$-group scheme $A^\natural$ introduced in Thm. 0.1.13. Further, if $M$ is a semi-abelian scheme $G$, extension of $A$ by $T$, then its universal extension identifies with $A^\natural \times_A G$, i.e. we take the pullback of the universal extension sequence $(0.1.3)$ for $A$ via the homomorphism $G \rightarrow A$.\\
\newline
Given a $1$-motive $M$ over $S$ with universal extension $E(M)=[X \xrightarrow{v} E(M)_G]$, its de Rham realization $\mathrm{T}_{\mathrm{dR}}(M)$ is defined to be the $\mathcal O_S$-vector bundle \[\mathrm{T}_{\mathrm{dR}}(M):=\Lie(E(M)_G/S),\] where the last means the Lie algebra relative $S$ of the (smooth) $S$-group scheme $E(M)_G$.\\
We recall from Thm. 0.1.13 that for $M=A$ an abelian scheme over $S$ we have canonically \[\mathrm{T}_{\mathrm{dR}}(A)\simeq H^1_{\mathrm{dR}}(A^\vee/S),\] where $A^\vee$ denotes the dual abelian scheme of $A$.\\
\newline
Taking universal extensions as well as de Rham realizations is covariant functorial and exact with respect to short exact sequences of $1$-motives over $S$ (cf. \cite{An-Ber}, Lemma 4.1).\footnote{A morphism between universal extensions is defined analogously as for $1$-motives; exactness for sequences of $1$-motives and universal extensions is to be understood on the level of complexes of abelian $fppf$-sheaves.}
\subsubsection{Cartier duality}
We briefly recall the construction of the Cartier dual for a $1$-motive. Full details can be found in \cite{Jo1}, 1.3 (where more general motives with torsion are admitted) and in \cite{BaVi}, 2.2.7.\\
\newline
If $M=[X \xrightarrow{u} G]$ is a $1$-motive over $S$ we define abelian $fppf$-sheaves over $S$
\[X^\vee:=\underline{\Hom}_{fppf}(T,\mathbb G_{m,S}), \quad G^\vee:=\underline{\Ext}^1_{fppf}([X\rightarrow A], \mathbb G_{m,S}), \quad T^\vee:=\underline{\Hom}_{fppf}(X,\mathbb G_{m,S}),\]
which are all recognized as being represented by commutative $S$-group schemes. $X^\vee$ resp. $T^\vee$ becomes a lattice resp. torus over $S$. By applying $R\underline{\Hom}_{fppf}(-,\mathbb G_{m,S})$ to the distinguished triangle
\[A \rightarrow [X\rightarrow A] \rightarrow [X\rightarrow 0] ,\]
going into cohomology and using the isomorphism of Barsotti-Rosenlicht-Weil (cf. Thm. 0.1.26).
\[A^\vee \simeq \underline{\Ext}^1_{fppf}(A, \mathbb G_{m,S})\]
yields an exact sequence $T^\vee \rightarrow G^\vee \rightarrow A^\vee$;  that this is in fact an extension of $A^\vee$ by $T^\vee$ follows from the vanishing of $\underline{\Hom}_{fppf}(A,\mathbb G_{m,S})$ and $\underline{\Ext}^1_{fppf}(X, \mathbb G_{m,S})$ (cf. \cite{SGA7-I}, exp. VIII, $(3.2.1)$ and \cite{Berto}, Lemma 1.1.4).\\
Similarly, the canonical distinguished triangle
\[T \rightarrow M \rightarrow [X \rightarrow A]\]
induces a homomorphism $u^\vee: X^\vee \rightarrow G^\vee$.\\
The $1$-motive $M^\vee=[X^\vee\xrightarrow{u^\vee} G^\vee]$ over $S$ obtained in this way is the Cartier dual of $M$. Its formation is contravariant functorial in $1$-motives over $S$ and satisfies a natural double duality isomorphism, such that it defines an antiautomorphism of the category of $1$-motives over $S$.\\
\newline
We remark that it is possible to interpret $M^\vee$ as representing the $\mathbb G_{m,S}$-biextension functor associated with $M$ on the category of $1$-motives over $S$ (cf. \cite{BaVi}, 2.2.24).\\
\newline
As an easy example, the Cartier dual of the $1$-motive $[\mathbb Z_S \rightarrow 0]$ is given by $\mathbb G_{m,S}$, where we write $\mathbb Z_S$ for the constant $S$-group scheme associated with the abstract group $\mathbb Z$.
\subsubsection{Deligne's pairing}
A canonical construction for Cartier dual $1$-motives over $S$
\[M=[X \xrightarrow{u} G], \quad M^\vee=[X^\vee \xrightarrow{u^\vee} G^\vee]\]
is the Deligne pairing
\[\Phi: \mathrm{T}_{\mathrm{dR}}(M) \otimes_{\mathcal O_S} \mathrm{T}_{\mathrm{dR}}(M^\vee) \rightarrow \mathcal O_S,\]
between their de Rham realizations, already explained at the end of 0.1.3 for the special case of $M=A$ an abelian scheme, where it writes as
\[\Phi: H^1_{\mathrm{dR}}(A^\vee/S) \otimes_{\mathcal O_S} H^1_{\mathrm{dR}}(A/S) \rightarrow \mathcal O_S.\]
We recall from \cite{Ber}, 4, that in the general motivic setting $\Phi$ is obtained as follows: The canonical $S$-connection of the Poincaré biextension on $E(M)_G \times_S E(M^\vee)_{G^\vee}$ induces via its curvature (which is an invariant $2$-form) an alternating $\mathcal O_S$-bilinear form on Lie algebras relative $S$, writing as
\[R: (\mathrm{T}_{\mathrm{dR}}(M) \oplus \mathrm{T}_{\mathrm{dR}}(M^\vee)) \oplus (\mathrm{T}_{\mathrm{dR}}(M) \oplus \mathrm{T}_{\mathrm{dR}}(M^\vee)) \rightarrow \mathcal O_S.\]
The Deligne pairing is then defined by
\[\Phi( v \otimes w):= R((v,0),(0,w)),\]
and it is a fundamental fact, proven in full generality in \cite{Ber}, Thm. 4.3, that $\Phi$ is perfect.
\subsection{The motivic Gauß-Manin connection and de Rham-Manin map}
The recent results in \cite{An-Ber} permit the definition of various realization maps on the sections of a $1$-motive. In the case of the de Rham realization, which we are interested in, this produces extensions of vector bundles with integrable connection and thus provides the suitable tool for a motivic interpretation of the logarithm extension. Let us hence review the relevant constructions of \cite{An-Ber}.
\subsubsection{Motivic Gauß-Manin connection}
Let $B$ be a locally noetherian scheme and $M=[X \xrightarrow{u}G]$ a $1$-motive over $B$ with universal extension $E(M)=[X \xrightarrow{v} E(M)_G]$. Further, let $q:B\rightarrow T$ be a smooth morphism with $T$ locally noetherian.\\
\newline
Then, following \cite{An-Ber}, 4.2, an integrable \underline{motivic Gauß-Manin connection for $M$ and $q:B\rightarrow T$}
\[\nabla_M: \mathrm{T}_{\mathrm{dR}}(M) \rightarrow \Omega^1_{B/T} \otimes_{\mathcal O_B} \mathrm{T}_{\mathrm{dR}}(M)\]
on the de Rham realization of $M$ can be constructed as follows:\\
If $\Delta^1_{B/T}$ denotes the first infinitesimal neighborhood of the diagonal immersion $B \rightarrow B\times_T B$, we have natural morphisms
\[B \xrightarrow{i} \Delta^1_{B/T} \rightarrow B\times_T B\]
composing to the diagonal, where $i$ is a nilpotent closed immersion of square zero.\\
Write $p_1,p_2: \Delta^1_{B/T} \rightarrow B$ for the maps induced by the projections of $B\times_T B$ and let $M_1, M_2$ be the $1$-motives over the (locally noetherian) scheme $\Delta^1_{B/T}$ obtained by base extension of $M$ via $p_1,p_2$.\\
Then, the crucial point for the construction of the connection on $\mathrm{T}_{\mathrm{dR}}(M)=\Lie(E(M)_G/B)$ is the existence of a canonical isomorphism of $\Delta^1_{B/T}$-group schemes
\[E_{\Delta^1_{B/T}}(\mathrm{id}_M):E(M_1)_{G_1} \simeq E(M_2)_{G_2}\]
which becomes the identity on $E(M)_G$ after further base change to $B$ via $i$. Passing to Lie algebras relative $\Delta^1_{B/T}$ in this isomorphism yields the desired connection on $\mathrm{T}_{\mathrm{dR}}(M)$ (cf. 0.2.1 (vii)). Its integrability is a consequence of the smoothness of $q: B \rightarrow T$ (cf. \cite{An-Ber}, 4.2 (d)).\\
\newline
The existence of $E_{\Delta^1_{B/T}}(\mathrm{id}_M)$ follows from a more general motivic deformation result for locally nilpotent $PD$-thickenings, proven by interpreting the universal extension in the framework of the crystalline site (cf. ibid., Thm. 2.1 and section 3).\\
\newline
The connection $\nabla_M$ enjoys a number of natural functorial properties and turns out to be the expectable connection in special cases: e.g. if $M=\mathbb G_{m,B}$ it equals the exterior derivation on $\mathcal O_B$, and for $M=A$ an abelian scheme it coincides with the usual Gauß-Manin connection on $H^1_{\mathrm{dR}}(A^\vee/B)$ (cf. ibid., Ex. 4.3 and Lemma 4.5).
\begin{remark}
When working with $\mathbb Z$-flat schemes any isomorphism of $\Delta^1_{B/T}$-group schemes
\[E(M_1)_{G_1} \simeq E(M_2)_{G_2}\]
restricting under $i$ to the identity on $E(M)_G$ is already the canonical isomorphism $E_{\Delta^1_{B/T}}(\mathrm{id}_M)$ of above (cf. the argument in \cite{An-Ber}, beginning of 6).
\end{remark}

\subsubsection{Motivic de Rham-Manin map}
As in \cite{An-Ber}, 7.1, we define the group $M(B)$ of $B$-rational points of $M$ by
\[M(B):=\Hom_{D^b(B_{fppf})}(\mathbb Z_B,M)\]
and the \underline{motivic de Rham-Manin map for $M$ and $q:B \rightarrow T$}
\[\mathcal M_{M,{\mathrm{dR}}}: M(B) \rightarrow \Ext^1_{\mathcal D_{B/T}}((\mathrm{T}_{\mathrm{dR}}(M^\vee), \nabla_{M^\vee}),(\mathcal O_B,\mathrm{d}))\]
as the following composition:
\[\Hom_{D^b(B_{fppf})}(\mathbb Z_B,M) \xrightarrow{\sim}  \Ext^1_{1-Mot}(M^\vee, \mathbb G_{m,B}) \rightarrow \Ext^1_{\mathcal D_{B/T}}((\mathrm{T}_{\mathrm{dR}}(M^\vee), \nabla_{M^\vee}),(\mathcal O_B,\mathrm{d})).\]
Here, the first arrow comes from the chain of natural identifications (cf. ibid., $(16)$)
\[\tag{\textbf{2.6.1}} \Hom_{D^b(B_{fppf})}(\mathbb Z_B,M)\simeq \Ext^1_{1-Mot}(\mathbb Z[1],M) \simeq  \Ext^1_{1-Mot}(M^\vee, \mathbb G_{m,B}),\]
and the second is given by passage to de Rham realizations, which are viewed as $\mathcal D_{B/T}$-modules via the associated motivic Gauß-Manin connection.\\
The map $\mathcal M_{M,{\mathrm{dR}}}$ is a group homomorphism behaving functorially in $M$ and also in $B$ when working with $\mathbb Z$-flat schemes (cf. ibid., Prop. 7.2).\\
\newline
Note that for $M=A$ an abelian scheme over $B$ the above group $A(B)$ identifies with the usual group of $B$-rational points $\mathrm{Hom}_B(B,A)$, which we will also denote by $A(B)$, not distinguishing between a section $s$ of $A/B$ and its associated homomorphism $\mathbb Z_B \xrightarrow{1\mapsto s} A$.\\
Observe further that the motivic de Rham-Manin map for $A$ and $q:B \rightarrow T$ writes as
\[\mathcal M_{A,{\mathrm{dR}}}: A(B) \rightarrow \Ext^1_{\mathcal D_{B/T}}(H^1_{\mathrm{dR}}(A/B),\mathcal O_B),\]
where (according to \cite{An-Ber}, Lemma 4.5) $H^1_{\mathrm{dR}}(A/B)$ carries its usual Gauß-Manin connection.\\
\newline
Let us point already now to the subtle appearance of a sign, which becomes important when interpreting the subsequent motivic results geometrically by the Poincaré bundle at the end of the section.
\begin{remark}
If $M=A$ is an abelian scheme over $B$ and $A^\vee$ denotes its dual, then the composite
\[A(B) \simeq \Ext^1_{1-Mot}(A^\vee,\mathbb G_{m,B}) \simeq \Ext^1_{fppf}(A^\vee,\mathbb G_{m,B})\]
of $(2.6.1)$ with the obvious identification differs from the Barsotti-Rosenlicht-Weil isomorphism in Thm. 0.1.26 - biduality of Thm. 0.1.10 tacitly implied - by a minus sign (cf. \cite{An-BaVi}, 1.2, 1)).\footnote{The cited sign change is based on the fact that in the chain $(2.6.1)$ one is actually using an intermediate identification
\[\Hom_{D^b(B_{fppf})}(\mathbb Z_B,M)\simeq \Ext^1_{D^b(B_{fppf})}(\mathbb Z,M[-1]) \simeq \Ext^1_{D^b(B_{fppf})}(\mathbb Z[1],M),\]
such that a shift of distinguished triangles appears (cf. also once again \cite{An-Ber}, $(16)$).}
\end{remark}
In view of this it might seem more natural to modify the motivic de Rham-Manin map and precompose it with the inversion automorphism of $M(B)$. The resulting homomorphism
\[\widetilde{\mathcal M}_{M,{\mathrm{dR}}}: M(B) \rightarrow \Ext^1_{\mathcal D_{B/T}}((\mathrm{T}_{\mathrm{dR}}(M^\vee), \nabla_{M^\vee}),(\mathcal O_B,\mathrm{d}))\]
then equals $-\mathcal M_{M,{\mathrm{dR}}}$ and in the case of $M=A$ an abelian scheme is now given by composing the Barsotti-Rosenlicht-Weil identification
\[A(B) \simeq \Ext^1_{fppf}(A^\vee,\mathbb G_{m,B}) \simeq \Ext^1_{1-Mot}(A^\vee,\mathbb G_{m,B})\]
with the de Rham realization map (motivic Gauß-Manin connections implied)
\[\Ext^1_{1-Mot}(A^\vee,\mathbb G_{m,B}) \rightarrow \Ext^1_{\mathcal D_{B/T}}(H^1_{\mathrm{dR}}(A/B),\mathcal O_B).\]

\subsection{The first logarithm extension and the motivic de Rham-Manin map}
In the following we show how the dual of the elliptic logarithm extension can be realized as image under a suitably chosen motivic de Rham-Manin map.\\
With our above preparations the precise statement is quickly given, but its detailed proof will require some work and occupy the whole of this subsection. The basic idea consists in first relating the logarithm extension to the so-called "classical Manin map" of \cite{Co2} and in then using a comparison theorem between this last map and the motivic de Rham-Manin map, established in \cite{An-Ber} by a translation of the whole situation into (log-) crystalline cohomology.\footnote{A different and very explicit proof for the equality of the two maps was found by Bertapelle and the author:\\
For this the realization sequence coming from the motivic map is described purely in terms of $\natural$-extension sheaves of various motives by $\mathbb G_a$, whereas the sequence induced by the classical map is interpreted (as in \cite{Co2}) via \v{C}ech hypercohomology. The required isomorphism is then obtained by using the cocycle data to construct in a natural way corresponding $\natural$-extensions.}
\subsubsection{Formulation of the main result and a first step towards its proof}
If $S$ is a connected scheme which is smooth, separated and of finite type over $\Spec(\Q)$ and if $E \xrightarrow{\pi} S$ is an elliptic curve, we may view $E \times_S E$ as abelian scheme relative $E$ via the second projection $\mathrm{pr}_2$. The motivic de Rham-Manin map for $E\times_S E$ and the smooth morphism $E \rightarrow \Spec(\Q)$ then expresses as a homomorphism
\[\mathcal M_{(E\times_S E), {\mathrm{dR}}}: (E\times_S E)(E) \rightarrow \Ext^1_{\mathcal D_{E/\Q}}(\pi^*H^1_{\mathrm{dR}}(E/S),\mathcal O_E),\]
where $\mathcal O_E$ carries the trivial connection and $\pi^*H^1_{\mathrm{dR}}(E/S)$ is endowed with the pullback of the Gauß-Manin connection on $H^1_{\mathrm{dR}}(E/S)$.\\
Writing $\bar{\Delta}_E \in (E\times_S E)(E)$ for the inverse of the diagonal section $\Delta_E$, i.e. for the section given in rational points by $x \mapsto (-x,x)$, our main goal is to show
\begin{theorem}
The dual extension of $\mathcal M_{(E\times_S E), {\mathrm{dR}}}(\bar{\Delta}_E)$ is equal to $\mathcal Log^1$ in $\Ext^1_{\mathcal D_{E/\Q}}(\mathcal O_E,\mathcal H_E)$.\\
Equivalently: The dual extension of $\widetilde{\mathcal M}_{(E\times_S E), {\mathrm{dR}}}(\Delta_E)$ is equal to $\mathcal Log^1$ in $\Ext^1_{\mathcal D_{E/\Q}}(\mathcal O_E,\mathcal H_E)$.
\end{theorem}
Let us write $\xi \in \Ext^1_{\mathcal D_{E/\Q}}(\mathcal O_E,\mathcal H_E)$ for the dual extension of $\mathcal M_{(E\times_S E), {\mathrm{dR}}}(\bar{\Delta}_E)$, such that we need to show the equality of extension classes $\xi =\mathcal Log^1$.\\
At first, we can prove without difficulty that $\xi$ splits on $S$:
\begin{lemma}
$\xi$ maps to zero under the retraction
\[ \Ext^1_{\mathcal D_{E/\Q}}(\mathcal O_E, \mathcal H_E) \xrightarrow{\epsilon^*} \Ext^1_{\mathcal D_{S/\Q}}(\mathcal O_S, \mathcal H)\]
of the exact sequence $(1.1.1)$.
\end{lemma}
\begin{proof}
Observe the cartesian diagram
\begin{equation*} \tag{\textbf{2.6.2}}\begin{split}
\begin{xy}
\xymatrix{
E \ar[r]^{\pi} \ar[d]_{\id \times \epsilon} & S \ar[d]^{\epsilon} \\
E \times_S E \ar[r]^{ \ \ \ \ \mathrm{pr}_2} & E}
\end{xy}
\end{split}
\end{equation*}
Functoriality of the motivic de Rham-Manin map with respect to $\mathbb Z$-flat base schemes (cf. \cite{An-Ber}, Prop. 7.2) yields the commutative diagram with horizontal motivic de Rham-Manin maps (the lower one is with respect to the abelian scheme $E\rightarrow S$ and the smooth morphism $S \rightarrow \Spec(\Q)$)
\begin{equation*}
\begin{xy}
\xymatrix{
(E \times_S E)(E) \ar[r] \ar[d]^{\mathrm{can}} & \Ext^1_{\mathcal D_{E/\Q}}(\pi^*H^1_{\mathrm{dR}}(E/S), \mathcal O_E) \ar[d]^{\epsilon^*} \\
E(S) \ar[r] & \Ext^1_{\mathcal D_{S/\Q}}(H^1_{\mathrm{dR}}(E/S), \mathcal O_S)}
\end{xy}
\end{equation*}
where $\mathrm{can}$ stands for the canonical arrow induced by $(2.6.2)$.\\
The image of $\bar{\Delta}_E$ under $\mathrm{can}$ is $\epsilon \in E(S)$, but the lower horizontal arrow is a homomorphism (cf. \cite{An-Ber}, Prop. 7.2), hence $\mathcal M_{(E\times_S E), {\mathrm{dR}}}(\bar{\Delta}_E)$ maps to zero under $\epsilon^*$. This suffices to conclude.
\end{proof}
The main effort to prove Thm. 2.6.3 consists in showing that the image of $\xi$ under the projection
\[ \Ext^1_{\mathcal D_{E/\Q}}(\mathcal O_E, \mathcal H_E) \rightarrow \Hom_{\mathcal D_{S/\Q}}(\mathcal H, \mathcal H)
\]
in $(1.1.1)$ is the identity. This is nothing at all clear if one recalls the construction of the motivic de Rham-Manin map. We will achieve the proof by relating the logarithm extension to Coleman's classical Manin map in \cite{Co2}, which is defined rather explicitly, and by then using a comparison result in \cite{An-Ber} between the motivic de Rham-Manin map and the classical Manin map.
\subsubsection{A reduction step}
Recall from $(2.1.3)$ the commutative diagram of split exact sequences with vertical forgetful arrows
\[\tag{\textbf{2.6.3}} \begin{split}
\begin{xy}
\xymatrix@C-0.3cm{
0 \ar[r] & \Ext^1_{\mathcal D_{S/\Q}}(\mathcal O_S, \mathcal H) \ar[d]^{\mathrm{can}} \ar[r]^{\pi^*} & \Ext^1_{\mathcal D_{E/\Q}}(\mathcal O_E, \mathcal H_E) \ar[d]^{\mathrm{can}} \ar[r]& \Hom_{\mathcal D_{S/\Q}}(\mathcal H, \mathcal H) \ar[d]^{\mathrm{can}} \ar[r] &0 \\
0 \ar[r] & \Ext^1_{\mathcal O_S}(\mathcal O_S, \mathcal H) \ar[r]^{\pi^* \ \ } & \Ext^1_{\mathcal D_{E/S}}(\mathcal O_E, \mathcal H_E) \ar[r] & \Hom_{\mathcal O_S}(\mathcal H, \mathcal H) \ar[r] & 0}
\end{xy}
\end{split}
\]
Then, by the previous lemma it obviously suffices to show that $\mathrm{can}(\xi)$ gives the identity when projected to $\Hom_{\mathcal O_S}(\mathcal H, \mathcal H)$. But an easy application of \cite{An-Ber}, Thm. 2.1 (iii), shows that the element $\mathrm{can}(\xi)$ is nothing else than the dual of the image of $\bar{\Delta}_E$ under the motivic de Rham-Manin map
\[\mathcal M_{(E\times_S E), {\mathrm{dR}}}^{res}: (E\times_S E)(E) \rightarrow \Ext^1_{\mathcal D_{E/S}}(\pi^*H^1_{\mathrm{dR}}(E/S),\mathcal O_E)\]
for the abelian scheme $E\times_S E \xrightarrow{\mathrm{pr}_2} E$ and the smooth morphism $E \xrightarrow{\pi} S$. Note the difference of $\mathcal M_{(E\times_S E), {\mathrm{dR}}}^{res}$ to the previous de Rham-Manin map
\[\mathcal M_{(E\times_S E), {\mathrm{dR}}}: (E\times_S E)(E) \rightarrow \Ext^1_{\mathcal D_{E/\Q}}(\pi^*H^1_{\mathrm{dR}}(E/S),\mathcal O_E)\]
and observe that now $\pi^*H^1_{\mathrm{dR}}(E/S)$ is equipped with the trivial $S$-connection (by ibid., 4.2 (c)).\\
Writing $\mathcal Log^{1,res}$ for the element in $\Ext^1_{\mathcal D_{E/S}}(\mathcal O_E, \mathcal H_E)$ which maps to zero under the retraction and to the identity under the projection in the lower row of $(2.6.3)$ we have shown:
\begin{Lemma}
In order to prove Thm. 2.6.3 it suffices to verify that
\[(\mathcal M_{(E\times_S E), {\mathrm{dR}}}^{res}(\bar{\Delta}_E))^\vee=\mathcal Log^{1,res} \ \ \textrm{in} \ \ \Ext^1_{\mathcal D_{E/S}}(\mathcal O_E, \mathcal H_E),\]
where $^\vee$ means the dual extension class.\\
(We have in fact shown that we only need to verify that $(\mathcal M_{(E\times_S E), {\mathrm{dR}}}^{res}(\bar{\Delta}_E))^\vee$ projects to the identity in the lower row of $(2.6.3)$, but it is convenient to formulate the lemma as we just did.)\qquad \qed
\end{Lemma}
But for an elliptic curve $E$ over \underline{any} connected, separated, noetherian, regular and finite-dimensional scheme $S$ of characteristic zero (not necessarily smooth and of finite type over $\Q$) we can define an element $\Ext^1_{\mathcal D_{E/S}}(\mathcal O_E, \mathcal H_E)$, characterized by mapping to zero resp. to the identity under the retraction resp. projection in the split exact sequence
\[0 \rightarrow \Ext^1_{\mathcal O_S}(\mathcal O_S, \mathcal H) \xrightarrow{\pi^*} \Ext^1_{\mathcal D_{E/S}}(\mathcal O_E, \mathcal H_E) \rightarrow \Hom_{\mathcal O_S}(\mathcal H, \mathcal H) \rightarrow 0.\]
Let us write $\mathfrak{Log}^1$ for this element. If $S$ happens to be smooth of finite type over $\Q$, then $\mathfrak{Log}^1$ is the class $\mathcal Log^{1,res}$ of Lemma 2.6.5. But we will from now on - until we have finished the proof of Thm. 2.6.3 - work with an elliptic curve $E$ over a general connected, separated, noetherian, regular and finite-dimensional $\Q$-scheme $S$.\\
We then have the motivic de Rham-Manin map for the abelian scheme $E\times_S E \xrightarrow{\mathrm{pr}_2} E$ and the smooth map $E \xrightarrow{\pi} S$
\[\mathcal M_{(E\times_S E), {\mathrm{dR}}}: (E\times_S E)(E) \rightarrow \Ext^1_{\mathcal D_{E/S}}(\pi^*H^1_{\mathrm{dR}}(E/S),\mathcal O_E),\]
which in the case of $\Q$-smooth $S$ is the map $\mathcal M_{(E\times_S E), {\mathrm{dR}}}^{res}$ of Lemma 2.6.5.\\
\newline
\textit{With these definitions, what we will do in the following (until Cor. 2.6.14) is to show that for an elliptic curve $E \xrightarrow{\pi} S$, with $S$ a connected, separated, noetherian, regular and finite-dimensional $\Q$-scheme, we have}
\[\tag{\textbf{2.6.4}} \mathfrak{Log}^1=(\mathcal M_{(E\times_S E), {\mathrm{dR}}}(\bar{\Delta}_E))^\vee \ \ in \ \ \Ext^1_{\mathcal D_{E/S}}(\mathcal O_E,\mathcal H_E).\]
With Lemma 2.6.5 and our remarks this will in particular prove Thm. 2.6.3.
\subsubsection{An auxiliary description for $\mathfrak{Log}^1$}
Write $U$ for the open complement of the zero section of $E$, $v: V \subset U \times_S U$ for the embedding of the open complement of the diagonal $\Delta_U: U \rightarrow U \times_S U$
and $p_2^U: U \times_S U \rightarrow U$ for the second projection.\\
\newline
We have the canonical distinguished triangle in $D^b_{\mathrm{qc}}(\mathcal D_{U \times_S U/S})$ (cf. $(0.2.5)$)
\[(\Delta_U)_+ \mathcal O_U[-1] \rightarrow \mathcal O_{U\times_S U} \rightarrow v_+ \mathcal O_V.\]
Applying $(p_2^U)_+$ gives the distinguished triangle in $D^b_{\mathrm{qc}}(\mathcal D_{U/S})$
\[(p_2^U\circ \Delta_U)_+\mathcal O_U[-1] \rightarrow (p_2^U)_+ \mathcal O_{U\times_S U} \rightarrow (p_2^U \circ v)_+ \mathcal O_V,\]
from which we obtain the exact sequence of vector bundles on $U$ with integrable $S$-connection
\[\tag{\textbf{2.6.5}} 0 \rightarrow H^1_{\mathrm{dR}}(U \times_S U/U) \rightarrow H^1_{\mathrm{dR}}(V/U) \rightarrow H^0_{\mathrm{dR}}(U/U) \rightarrow 0,\]
where we have used that $H^2_{\mathrm{dR}}(U\times_S U/U)=0$ (cf. Thm. 1.2.9 (ii)).\\
If we canonically identify $H^1_{\mathrm{dR}}(U \times_S U/U) \simeq H^1_{\mathrm{dR}}(E \times_S U /U) \simeq \mathcal H_U$ (cf. Thm. 1.2.9 (iii) and the beginning of Chapter 1) and 
observe that $p_2^U\circ \Delta_U=\id_U$, such that $H^0_{\mathrm{dR}}(U/U)$ equals $\mathcal O_U$, we arrive at
\[\tag{\textbf{2.6.6}} 0 \rightarrow \mathcal H_U \rightarrow H^1_{\mathrm{dR}}(V/U) \rightarrow \mathcal O_U \rightarrow 0.\]
Note that $\mathcal H_U$ is equipped with the trivial $S$-connection, $\mathcal O_U$ with the exterior derivative and $H^1_{\mathrm{dR}}(V/U)$ with the Gauß-Manin connection relative $S$, where $V$ is an $U$-scheme via $V \xrightarrow{v} U \times_S U \xrightarrow{p_2^U} U$.\footnote{This is a general fact concerning the functor $(-)_+$: Let $T$ be a noetherian, regular and finite-dimensional $\Q$-scheme, and let $X \xrightarrow{f} Y \rightarrow T$ be smooth arrows of finite type such that $X$ resp. $Y$ has relative dimension $d_{X/T}$ resp. $d_{Y/T}$ over $T$. Let $f_+: D^b_{\mathrm{qc}}(\mathcal D_{X/T}) \rightarrow D^b_{\mathrm{qc}}(\mathcal D_{Y/T})$ be the triangulated functor introduced in 0.2.3 and $\mathcal M \in \Mod_{\mathrm{qc}}(\mathcal D_{X/T})$. Then $H^i(f_+ \mathcal M) \simeq H_{\mathrm{dR}}^{i+d_{X/T}-d_{Y/T}}(X/Y, \mathcal M)$ canonically as $\mathcal D_{Y/T}$-modules, where the right side is equipped with its Gauß-Manin connection relative $T$ as in 0.2.2.. The proof is analogous to \cite{Dim-Ma-Sa-Sai}, Prop. 1.4.}
\begin{lemma}
The extension class of $(2.6.6)$ in $\Ext^1_{\mathcal D_{U/S}}(\mathcal O_U, \mathcal H_U)$ is the restriction of $\mathfrak{Log}^1$ to $U$, i.e. its image under the canonical map
\[\Ext^1_{\mathcal D_{E/S}}(\mathcal O_E, \mathcal H_E) \rightarrow  \Ext^1_{\mathcal D_{U/S}}(\mathcal O_U, \mathcal H_U).\]
\end{lemma}
\begin{proof}
The proof is entirely formal and (up to some minor supplements) a reproduction of the arguments in \cite{Ki4}, proof of Prop. 2.3.2, where the statement is shown in the case of $\ell$-adic sheaves.
\end{proof}
In fact, we will need the previous lemma only in the case where $S=\Spec(k)$ with $k$ a field of characteristic zero (namely for the proof of Cor. 2.6.11).
\begin{remark}
Assume that $S=\Spec(k)$, such that $U$ is affine and irreducible.\\
In \cite{Co2}, p. 404 (before and in the proof of Lemma 1.5.1), $\Gamma(U,H^0_{\mathrm{dR}}(U/U))$ is then interpreted as the group of divisors in $U \times_k U$ defined over $U$ and supported on $\Delta_U(U)\subseteq U\times_k U$. The identification $H^0_{\mathrm{dR}}(U/U) \simeq \mathcal O_U$ we made above then represents $H^0_{\mathrm{dR}}(U/U)$ as the free $\mathcal O_U$-module generated by the global section corresponding to the divisor $\Delta_U(U)$. We will need this trivial remark in the following when working with Coleman's classical Manin map.
\end{remark}

\subsubsection{The classical Manin map}
If $U$ is a smooth irreducible affine curve over a field $k$ of characteristic zero and $A$ is an abelian scheme over $U$, then by \underline{the classical Manin map} we mean the homomorphism
\[\mathcal M_A: A(U) \rightarrow \Ext^1_{\mathcal D_{U/k}}(H^1_{\mathrm{dR}}(A/U), \mathcal O_U)\]
defined in \cite{Co2}, 4. The investigation of its kernel is a major tool in Coleman's account of Manin's proof of the Mordell conjecture over function fields. We don't recapitulate its construction here (cf. ibid., 3 and 4), but instead record the following crucial comparison result of \cite{An-Ber}, Prop. 1.1:
\begin{theorem}
For all $s \in A(U)$ we have the equality
\[\mathcal M_A(s)=\mathcal M_{A,{\mathrm{dR}}}(s),\]
where $\mathcal M_{A,{\mathrm{dR}}}$ is the motivic de Rham-Manin map for $A$ and $U \rightarrow \Spec(k)$.
\end{theorem}
For the following, the strategy to prove Thm. 2.6.3 consists in relating the logarithm extension to the classical Manin map by a certain intermediate extension and the description given in Lemma 2.6.6; the preceding theorem will then give the desired realization via the motivic de Rham-Manin map.\\
When doing this, the restriction to irreducible affine curves as base schemes in the classical Manin map and in Thm. 2.6.8 will force us to consider the logarithm extension at first on a pointed single elliptic curve and to then, by appropriate techniques, haul the result up to the entire curve and families.
\subsubsection{An intermediate extension}
Assume that $U$ is as above and that $\tau: A=C \rightarrow U$ is a family of elliptic curves.\\
Let $s: U \rightarrow C$ be a section disjoint from the unit section $\epsilon$. Write $z:Z \rightarrow C$ for the closed subscheme of $C$ defined by $\epsilon \cup s$ and $\tilde{v}: V \rightarrow C$ for the open immersion of its complement $V$. As $Z$ is smooth over $U$ (and hence over $k$) we have the localization sequence in $D^b_{\mathrm{qc}}(\mathcal D_{C/k})$ (cf. $(0.2.5)$):
\[z_+ \mathcal O_Z[-1] \rightarrow \mathcal O_C \rightarrow \tilde{v}_+ \mathcal O_V.\]
By applying $\tau_+$ and taking cohomology we recover the exact $\mathcal D_{U/k}$-linear sequence of \cite{Co2}, $(5.2)$:
\[\tag{\textbf{2.6.7}} 0 \rightarrow H^1_{\mathrm{dR}}(C/U) \rightarrow H^1_{\mathrm{dR}}(V/U) \rightarrow H^0_{\mathrm{dR}}(Z/U) \rightarrow H^2_{\mathrm{dR}}(C/U) \rightarrow 0.\]
From the fact that $H^0_{\mathrm{dR}}(Z/U)$ consists of two copies of $\mathcal O_U$, generated by the divisors $s$ and $\epsilon$ of $C$ as in Coleman's interpretation of $H^0_{\mathrm{dR}}(Z/U)$ (cf. ibid., p. 404) and from the canonical identification $H^2_{\mathrm{dR}}(C/U) \simeq \mathcal O_U$ (cf. the beginning of 2) one obtains an exact $\mathcal D_{U/k}$-linear sequence
\[B_{s,\epsilon}: \quad 0 \rightarrow H^1_{\mathrm{dR}}(C/U) \rightarrow H^1_{\mathrm{dR}}(V/U) \rightarrow \mathcal O_U \rightarrow 0,\]
where $\mathcal O_U$ is identified with the free $\mathcal O_U$-module over the divisor $D=s-\epsilon$ (cf. ibid., p. 406).
\begin{lemma}
Under the Poincaré duality identification $H^1_{\mathrm{dR}}(C/U) \simeq H^1_{\mathrm{dR}}(C/U)^\vee$ the extension class of $B_{s,\epsilon}$ in $\Ext^1_{\mathcal D_{U/k}}(\mathcal O_U,H^1_{\mathrm{dR}}(C/U)^\vee)$ becomes equal to the opposite of the dual extension of $\mathcal M_C(s)$.
\end{lemma}
\begin{proof}
This follows from \cite{Co2}, Lemma 1.5.5 and Prop. 1.3.1, together with the definition of the classical Manin map (cf. ibid., p. 402).
\end{proof}
\subsubsection{Motivic description of the logarithm on a pointed single elliptic curve}
We consider a single elliptic curve $E \xrightarrow{\pi} \Spec(k)$ over a field $k$ of characteristic zero.\\
The open complement $U=E \backslash [0]$ of its zero point is a smooth irreducible affine curve over $k$.\\
We consider $E \times_k U$ as abelian scheme over $U$ via the second projection, and we define a section $s \in (E\times_k U)(U)$ to be the composition
\[s: U \xrightarrow{ \ \ \Delta_U} U\times_k U \xrightarrow{j} E \times_k U\]
with $j$ the canonical open immersion. Note that $s$ is disjoint from the unit section $\epsilon$ of $E\times_k U$.\\
Write $z:Z \rightarrow E\times_k U$ for the closed subscheme defined by $\epsilon \cup s$ and let $\tilde{v}:V\rightarrow E\times_k U$ be the open immersion of its complement. Note that the last coincides precisely with the open complement in $U\times_kU$ of the diagonal $\Delta_U: U \rightarrow U\times_k U$, embedded into $E\times_kU$ via $U \times_kU \rightarrow E\times_k U$.\\
If we apply the procedure in our previous exposition about the intermediate extension for the case $C:=E\times_kU$, then we obtain the exact sequence of $\mathcal D_{U/k}$-modules
\[B_{s, \epsilon}: \quad 0 \rightarrow H^1_{\mathrm{dR}}(E\times_k U /U) \rightarrow H^1_{\mathrm{dR}}(V/U) \rightarrow \mathcal O_U \rightarrow 0,\]
where again $\mathcal O_U$ actually is the free $\mathcal O_U$-module generated by the divisor $D=s-\epsilon$ of $E\times_k U$.\\
\newline
On the other hand, by Lemma 2.6.6 we can describe the restriction $\mathfrak{Log}^1_{|U}$ of $\mathfrak{Log}^1$ to $U$ by identifying $H^1_{\mathrm{dR}}(U \times_k U/U) \simeq H^1_{\mathrm{dR}}(E \times_k U /U) \simeq \mathcal H_U$ in the $\mathcal D_{U/k}$-linear exact sequence
\[\tag{\textbf{2.6.8}} 0 \rightarrow H^1_{\mathrm{dR}}(U \times_k U/U) \rightarrow H^1_{\mathrm{dR}}(V/U) \rightarrow \mathcal O_U \rightarrow 0,\]
which was won from the localization triangle for the closed immersion $\Delta_U: U \rightarrow U \times_k U$.\\
\newline
The following theorem relates the logarithm extension $\mathfrak{Log}^1_{|U}$ on $U$ with the classical Manin map $\mathcal M_{E\times_k U}$ for the elliptic curve $E\times_kU$ over $U$. The section $s\in (E\times_kU)(U)$ is as defined above.
\begin{theorem}
When identifying $H^1_{\mathrm{dR}}(E \times_k U/U) \simeq H^1_{\mathrm{dR}}(U \times_k U /U)$ the extensions $B_{s, \epsilon}$ and $(2.6.8)$ coincide. In particular, the class $\mathfrak{Log}^1_{|U}$ in $\Ext^1_{\mathcal D_{U/k}}(\mathcal O_U, \mathcal H_U)$ equals the opposite of the dual of $\mathcal M_{E\times_k U}(s)$.
\end{theorem}
\begin{proof}
We only need to show the first statement: the second then follows with Lemma 2.6.9.\\
We have commutative (in fact cartesian) squares
\begin{equation*}
\begin{xy}
\xymatrix{
U \ar[d]_{k} \ar[r]^{\Delta_U \quad}  & U \times_k U \ar[d]^{j} & V \ar[l]_{\quad \ \ v} \ar[d]^{\id}\\
Z \ar[r]^{z \quad} & E\times_k U & V \ar[l]_{\quad  \ \ \tilde{v}}}
\end{xy}
\end{equation*}
where $v$ is the open embedding of the complement of $\Delta_U$ and $k$ denotes the canonical open immersion (note that $j\circ\Delta_U=s$ and that $Z$ is the disjoint union of $\epsilon$ and $s$).\\
There are canonical (obvious) adjunction arrows\footnote{Observe \cite{Ho-Ta-Tan}, Ex. 1.5.22 and App. C, Prop. C. 2.4.} in $D^b_{\mathrm{qc}}(\mathcal D_{Z/k})$ resp. $D^b_{\mathrm{qc}}(\mathcal D_{E\times_k U/k})$ resp. $D^b_{\mathrm{qc}}(\mathcal D_{V/k})$: \[\mathcal O_Z \rightarrow k_+ \mathcal O_U, \quad \mathcal O_{E\times_k U} \rightarrow j_+\mathcal O_{U\times_k U}, \quad \mathcal O_V \rightarrow (\id)_+ \mathcal O_V, \] and it is easily seen that we obtain a morphism of distinguished triangles in $D^b_{\mathrm{qc}}(\mathcal D_{E\times_k U/k})$
\begin{equation*}
\begin{xy}
\xymatrix@C-0.3cm{
z_+ \mathcal O_Z[-1] \ar[d] \ar[r]  & \mathcal O_{E\times_kU} \ar[d] \ar[r] & \tilde{v}_+ \mathcal O_V \ar[d]\\
j_+ ((\Delta_U)_+ \mathcal O_U)[-1] \ar[r] & j_+ \mathcal O_{U\times_k U} \ar[r] & j_+(v_+ \mathcal O_V)}
\end{xy}
\end{equation*}
where the upper line is the localization triangle for $z: Z \rightarrow E\times_k U$, the lower one is $j_+$ applied to the localization triangle for $\Delta_U: U \rightarrow U \times_k U$ and the vertical arrows are induced by the three adjunction arrows and the above commutative diagram.\\
Applying the "lower plus" functor for the second projection of $E\times_k U$ and going into cohomology yields the commutative diagram of $\mathcal D_{U/k}$-modules with exact rows:
\begin{equation*}
\begin{xy}
\xymatrix@C-0.3cm{
0 \ar[r] & H^1_{\mathrm{dR}}(E \times_k U/U) \ar[d]_{\mathrm{can}}^{\sim} \ar[r]  & H^1_{\mathrm{dR}}(V/U) \ar[d]^{\id} \ar[r] & H^0_{\mathrm{dR}}(Z/U) \ar[r]^{\eta \quad \ } \ar[d] & H^2_{\mathrm{dR}}(E\times_k U/U) \ar[d] \ar[r] & 0 \\
0 \ar[r] & H^1_{\mathrm{dR}}(U \times_k U/U) \ar[r] & H^1_{\mathrm{dR}}(V/U) \ar[r] & H^0_{\mathrm{dR}}(U/U) \ar[r] & 0 \ar[r] & 0}
\end{xy}
\end{equation*}
Note that by construction the upper row is the sequence $(2.6.7)$ and the lower row is $(2.6.5)$.\\
It is clear that the arrow $H^0_{\mathrm{dR}}(Z/U) \rightarrow H^0_{\mathrm{dR}}(U/U)$ is given by restriction to the $s$-component of $Z$. Hence, the divisor $D= s - \epsilon$, which is the fixed generator for the kernel of $\eta$, hereunder maps to $\Delta_U$, which is the fixed $\mathcal O_U$-generator of $H^0_{\mathrm{dR}}(U/U)$ (cf. Rem. 2.6.7). This shows that the diagram
\begin{equation*}
\begin{xy}
\xymatrix@C-0.3cm{
0 \ar[r] & H^1_{\mathrm{dR}}(E \times_k U/U) \ar[d]^{\mathrm{can}}_{\sim} \ar[r]  & H^1_{\mathrm{dR}}(V/U) \ar[d]^{\id} \ar[r] & \mathcal O_U \ar[d]^{\id} \ar[r] &0 \\
0 \ar[r] & H^1_{\mathrm{dR}}(U \times_k U/U) \ar[r] & H^1_{\mathrm{dR}}(V/U) \ar[r] & \mathcal O_U \ar[r] & 0}
\end{xy}
\end{equation*}
is commutative, where the upper row is $B_{s, \epsilon}$ and the lower one is $(2.6.8)$; this shows the theorem.
\end{proof}

\begin{Korollar}
If $s \in (E\times_k U)(U)$ is the section $s: U \xrightarrow{\Delta_U} U \times_k U \xrightarrow{j} E\times_k U$, then the dual extension of $\mathcal M_{(E\times_k U), {\mathrm{dR}}} (s)$ is equal to $- \mathfrak{Log}^1_{|U}$ in $\Ext^1_{\mathcal D_{U/k}}(\mathcal O_U,\mathcal H_U)$.\\
Here,
\[\mathcal M_{(E\times_kU), {\mathrm{dR}}}: (E\times_k U)(U) \rightarrow \Ext^1_{\mathcal D_{U/k}}(\pi_U^*H^1_{\mathrm{dR}}(E/k), \mathcal O_U)\]
is the motivic de Rham-Manin map for the abelian $U$-scheme $E\times_k U$ and the map $\pi_U: U \rightarrow \Spec(k)$.
\end{Korollar}
\begin{proof}
Combine Thm. 2.6.8 and Thm. 2.6.10.
\end{proof}
\subsubsection{Motivic description of the logarithm on a single elliptic curve}
In Cor. 2.6.11 we had to work on the complement of the zero section of $E$ in order to apply the comparison result of Thm. 2.6.8. The next step consists in removing this restriction.\\
\newline
Let again $E \xrightarrow{\pi} \Spec(k)$ be an elliptic curve over a field $k$ of characteristic zero. We view $E \times_k E$ as abelian $E$-scheme via the second projection and consider the motivic de Rham-Manin map
\[\mathcal M_{(E\times_kE), {\mathrm{dR}}}: (E\times_k E)(E) \rightarrow \Ext^1_{\mathcal D_{E/k}}(\pi^*H^1_{\mathrm{dR}}(E/k), \mathcal O_E).\]
Write $\Delta_E \in (E\times_k E)(E)$ for the diagonal section.
\begin{Proposition}
The dual extension of $\mathcal M_{(E\times_kE), {\mathrm{dR}}}(\Delta_E)$ is equal to $- \mathfrak{Log}^1$ in $\Ext^1_{\mathcal D_{E/k}}(\mathcal O_E,\mathcal H_E)$.
\end{Proposition}
\begin{proof}
We have an obvious cartesian diagram
\begin{equation*}
\begin{xy}
\xymatrix{
E \times_k U \ar[r] \ar[d] & U \ar[d] \\
E \times_k E \ar[r] & E}
\end{xy}
\end{equation*}
As $E$ and $U$ are flat over $\mathbb Z$ we obtain from \cite{An-Ber}, Prop. 7.2, the commutativity of the diagram
\begin{equation*}
\begin{xy}
\xymatrix{
(E \times_k E)(E) \ar[r] \ar[d] & \Ext^1_{\mathcal D_{E/k}}(\pi^*H^1_{\mathrm{dR}}(E/k), \mathcal O_E) \ar[d] \\
(E \times_k U)(U) \ar[r] & \Ext^1_{\mathcal D_{U/k}}(\pi_U^*H^1_{\mathrm{dR}}(E/k), \mathcal O_U)}
\end{xy}
\end{equation*}
The vertical arrows are the canonical ones and the horizontal arrows are the respective motivic de Rham-Manin maps.\\
It is readily checked that the image of $\Delta_E$ under the left vertical arrow is precisely the section $s$ of Cor. 2.6.11. The same corollary and the commutativity of the preceding diagram imply that $-\mathfrak{Log}^1$ and the dual of $\mathcal M_{(E\times_kE), {\mathrm{dR}}}(\Delta_E)$ map to the same element under the restriction arrow
\[\Ext^1_{\mathcal D_{E/k}}(\mathcal O_E,\mathcal H_E) \rightarrow \Ext^1_{\mathcal D_{U/k}}(\mathcal O_U,\mathcal H_U).\]
But as we are working over a field this arrow is an isomorphism, as one easily sees by writing
\[\Ext^1_{\mathcal D_{E/k}}(\mathcal O_E,\mathcal H_E)\simeq H^1_{\mathrm{dR}}(E/k) \otimes_k \mathcal H \simeq H^1_{\mathrm{dR}}(U/k) \otimes_k \mathcal H\simeq \Ext^1_{\mathcal D_{U/k}}(\mathcal O_U, \mathcal H_U),\]
where we have used that the restriction map $H^1_{\mathrm{dR}}(E/k) \rightarrow H^1_{\mathrm{dR}}(U/k)$ is an isomorphism (cf. Thm. 1.2.9 (iii)). This proves the desired equality.
\end{proof}
\subsubsection{Motivic description of the logarithm for families of elliptic curves}
We now generalize Prop. 2.6.12 to relative elliptic curves.\\
Let $S$ be a connected, separated, noetherian, regular and finite-dimensional $\Q$-scheme and $E \xrightarrow{\pi} S$ an elliptic curve with zero section $\epsilon$. Again, we view $E \times_S E \rightarrow E$ as abelian $E$-scheme via the second projection and consider the motivic de Rham-Manin map for $E\times_S E$ and the morphism $E \xrightarrow{\pi} S$:
\[\mathcal M_{(E\times_S E), {\mathrm{dR}}}: (E\times_S E)(E) \rightarrow \Ext^1_{\mathcal D_{E/S}}(\pi^*H^1_{\mathrm{dR}}(E/S),\mathcal O_E).\]
\begin{theorem}
The dual extension of $\mathcal M_{(E\times_S E), {\mathrm{dR}}}(\Delta_E)$ is equal to $- \mathfrak{Log}^1$ in $\Ext^1_{\mathcal D_{E/S}}(\mathcal O_E,\mathcal H_E)$.
\end{theorem}
\begin{proof}
Let us write $\chi \in \Ext^1_{\mathcal D_{E/S}}(\mathcal O_E,\mathcal H_E)$ for the negative of the dual extension of $\mathcal M_{(E\times_S E), {\mathrm{dR}}}(\Delta_E)$, such that we need to show $\chi =\mathfrak{Log}^1$.\\
That the image of $\chi$ under the retraction
\[ \Ext^1_{\mathcal D_{E/S}}(\mathcal O_E, \mathcal H_E) \xrightarrow{\epsilon^*} \Ext^1_{\mathcal O_S}(\mathcal O_S, \mathcal H)\]
is zero follows by the completely analogous argument as applied for the proof of Lemma 2.6.4.\\
It remains to verify that $\chi$ maps to the identity under the projection
\[\tag{\textbf{2.6.9}} \mathrm{pr}:\Ext^1_{\mathcal D_{E/S}}(\mathcal O_E, \mathcal H_E) \rightarrow \Hom_{\mathcal O_S}(\mathcal H, \mathcal H).\]
By the integrality of $S$ (cf. footnote 5 of Chapter 1) one is reduced to show that for all $s \in S$ the image of $\mathrm{pr}(\chi)$ under the canonical arrow
\[\Hom_{\mathcal O_S}(\mathcal H, \mathcal H) \rightarrow \Hom_{k(s)}(H^1_{\mathrm{dR}}(E_s/k(s))^\vee, H^1_{\mathrm{dR}}(E_s/k(s))^\vee)\]
is the identity; here, we set $E_s:=E \times_S \Spec(k(s))$, considered as elliptic curve over $\Spec(k(s))$, and recall that $H^1_{\mathrm{dR}}(E/S)$ commutes with arbitrary base change (cf. the beginning of Chapter 1).\\
Consider the following diagram with horizontal motivic de Rham-Manin maps
\begin{equation*} \tag{\textbf{2.6.10}}\begin{split}
\begin{xy}
\xymatrix{
(E_s \times_{k(s)} E_s)(E_s) \ar[r] & \Ext^1_{\mathcal D_{E_s/k(s)}}(\pi_s^*H^1_{\mathrm{dR}}(E_s/k(s)), \mathcal O_{E_s}) \\
(E \times_S E)(E) \ar[r] \ar[u]^{\mathrm{can}} & \Ext^1_{\mathcal D_{E/S}}(\pi^*H^1_{\mathrm{dR}}(E/S), \mathcal O_E) \ar[u]^{\mathrm{can}}
}
\end{xy}
\end{split}
\end{equation*}
and observe that the right vertical arrow is understood to be given by pullback along
\begin{equation*}
\begin{xy}
\xymatrix{
E_s \ar[r] \ar[d]_{\pi_s} & E \ar[d]^{\pi}\\
\Spec(k(s)) \ar[r] & S
}
\end{xy}
\end{equation*}
as explained in 0.2.1 (v); the left vertical arrow is the obvious one.\\
Assuming that $(2.6.10)$ commutes the theorem clearly follows from our above remarks together with Prop. 2.6.12, noting that (2.6.9) respects the base change.\\
The commutativity of $(2.6.10)$ in turn is a straightforward application of the third functoriality statement of \cite{An-Ber}, Prop. 7.2, combined with ibid., Thm. 2.1 (iii).
\end{proof}
Recall that we write $\bar{\Delta}_E \in (E\times_S E)(E)$ for the antidiagonal section, given in rational points by $x \mapsto (-x,x)$. As it is the inverse of $\Delta_E$ in the group $(E\times_S E)(E)$ and $\mathcal M_{(E\times_S E), {\mathrm{dR}}}$ is a homomorphism (cf. \cite{An-Ber}, Prop. 7.2) the preceding theorem yields:
\begin{corollary}
The dual extension of $\mathcal M_{(E\times_S E), {\mathrm{dR}}}(\bar{\Delta}_E)$ is equal to $\mathfrak{Log}^1$ in $\Ext^1_{\mathcal D_{E/S}}(\mathcal O_E,\mathcal H_E)$. \qquad \qed
\end{corollary}
We have thus shown $(2.6.4)$ and hence, according to the explanations given there, also Thm. 2.6.3 finally is proven.
\begin{remark}
From the beginning on we have considered $E \times_S E$ as abelian $E$-scheme via the second projection. The results of Cor. 2.6.14 and Thm. 2.6.3 hold verbatim if one changes the convention and uses the first projection, such that $\bar{\Delta}_E$ then is given by $x\mapsto (x,-x)$ in points. A quick way to see this consists in applying the functoriality of the motivic de Rham-Manin map (cf. \cite{An-Ber}, Prop. 7.2) to the shift automorphism of $E \times_S E$.
\end{remark}

\subsection{Some corollaries}
We outline how our main theorem leads to a geometric interpretation of the logarithm extension by the Lie algebra of the Poincaré bundle. We finally give an equivalent approach using Deligne duality.

\subsubsection{Description via the Poincaré bundle}
Let $E\xrightarrow{\pi} S$ be an elliptic curve, where $S$ is a connected scheme which is smooth, separated and of finite type over $\Spec(\Q)$.\\
\newline
By Thm. 2.6.3 the logarithm extension class $\mathcal Log^1 \in \Ext^1_{\mathcal D_{E/\Q}}(\mathcal O_E,\mathcal H_E)$ is the dual of the image of the diagonal section $\Delta_E$ under the modified motivic de Rham-Manin map
\[\widetilde{\mathcal M}_{(E\times_S E), {\mathrm{dR}}}: (E\times_S E)(E) \rightarrow \Ext^1_{\mathcal D_{E/\Q}}(\pi^*H^1_{\mathrm{dR}}(E/S),\mathcal O_E).\]
Let us write $E^\vee$ for the dual abelian scheme of $E$ and recall that $\widetilde{\mathcal M}_{(E\times_S E), {\mathrm{dR}}}$ is the composite
\[(E \times_S E)(E) \xrightarrow{\sim} \Ext^1_{1-Mot}(E^\vee \times_S E, \mathbb G_{m,E}) \xrightarrow{\mathrm{T}_{\mathrm{dR}}(-)} \Ext^1_{\mathcal D_{E/\Q}}(\pi^*H^1_{\mathrm{dR}}(E/S),\mathcal O_E),\]
in which the first arrow comes from the Barsotti-Rosenlicht-Weil isomorphism (cf. Rem. 2.6.2).\\
\newline
It follows from the discussion in 0.1.3 that under this first arrow the section $\Delta_E$ maps to the $fppf$-extension $(0.1.19)$ induced by the Poincaré bundle $\mathcal P^0$ on $E \times_S E^\vee$:
\[\tag{\textbf{2.6.11}} 0 \rightarrow \mathbb G_{m,E} \rightarrow P^0 \rightarrow E^\vee \times_S E \rightarrow 0,\]
which we view as an extension of $1$-motives over $E$ in the obvious way. To compute $\widetilde{\mathcal M}_{(E\times_S E), {\mathrm{dR}}}(\Delta_E)$ we then apply the realization functor $\mathrm{T}_{\mathrm{dR}}(-)$ to this extension and equip the terms in the obtained sequence of $\mathcal O_E$-vector bundles with their motivic Gauß-Manin connections relative $\Q$.\\
\newline
The sequence of universal extensions related to $(2.6.11)$ is (e.g. by \cite{An-BaVi}, Lemma 2.2.1) precisely the extension of $E$-group schemes $(0.1.23)$ coming from the Poincaré bundle $\mathcal P$ on $E \times_S (E^\vee)^\natural$:
\[0 \rightarrow \mathbb G_{m,E} \rightarrow P \rightarrow (E^\vee)^\natural \times_S E \rightarrow 0.\]
The associated exact sequence of Lie algebras relative $E$ writes as
\[\tag{\textbf{2.6.12}} 0 \rightarrow \mathcal O_E \rightarrow \Lie(P/E) \rightarrow \pi^*H^1_{\mathrm{dR}}(E/S) \rightarrow 0,\]
where all terms are viewed as equipped with their motivic Gauß-Manin connection relative $\Q$; recall that for the outer terms this is just the exterior derivative resp. the pullback of the usual Gauß-Manin connection. As the obtained class in $\Ext^1_{\mathcal D_{E/\Q}}(\pi^*H^1_{\mathrm{dR}}(E/S),\mathcal O_E)$ is by construction equal to $\widetilde{\mathcal M}_{(E\times_S E), {\mathrm{dR}}}(\Delta_E)$ we conclude from Thm. 2.6.3:
\begin{corollary}
The class in $\Ext^1_{\mathcal D_{E/\Q}}(\mathcal O_E,\mathcal H_E)$ obtained from dualizing $(2.6.12)$ coincides with $\mathcal Log^1$. \qquad \qed
\end{corollary}
\begin{remark}
One can show that the restriction to an $S$-connection of the motivic Gauß-Manin connection relative $\Q$ on $\Lie(P/E)$ is induced by the universal $(E^\vee)^\natural$-connection on $\mathcal P$ in a natural way, as one would expect. The full (absolute) connection on $\Lie(P/E)$, however, does not have an intrinsic expression via the geometry of the Poincaré bundle.
\end{remark}

\subsubsection{Another description via Deligne's pairing}
Finally, consider another time the (non-modified) motivic de Rham-Manin map
\[\mathcal M_{(E\times_S E), {\mathrm{dR}}}: (E \times_S E)(E) \xrightarrow{\sim} \Ext^1_{1-Mot}(E^\vee \times_S E, \mathbb G_{m,E}) \rightarrow \Ext^1_{\mathcal D_{E/\Q}}(\pi^*H^1_{\mathrm{dR}}(E/S),\mathcal O_E).\]
Cartier dualizing the image of $\bar{\Delta}_E$ under the first arrow yields the extension of $1$-motives over $E$:
\[\tag{\textbf{2.6.13}} 0 \rightarrow E \times_S E \rightarrow [\mathbb Z_E \xrightarrow{1 \mapsto \bar{\Delta}_E} E \times_S E] \rightarrow [\mathbb Z_E \rightarrow 0] \rightarrow 0,\]
where $\mathbb Z_E$ denotes the constant $E$-group scheme and the maps of the sequence are the natural ones.\\
Taking de Rham realizations in $(2.6.13)$ yields an exact sequence of $\mathcal O_E$-vector bundles, which we equip with the motivic Gauß-Manin connections relative $\Q$. The realizations of the outer terms in $(2.6.13)$ are $\pi^* H^1_{\mathrm{dR}}(E^\vee/S)$, equipped with the pullback of the usual Gauß-Manin connection, resp. $\mathcal O_E$, equipped with exterior derivation (cf. \cite{An-Ber}, 4.2 (c), Ex. 4.2 and Ex. 4.4).\\
By the perfect Deligne pairing we have a canonical isomorphism of $\mathcal O_S$-vector bundles
\[H^1_{\mathrm{dR}}(E^\vee/S) \simeq H^1_{\mathrm{dR}}(E/S)^\vee,\]
which is in fact horizontal (cf. ibid., Cor. 2.8).\\
We thus obtain the natural identification of $\mathcal O_E$-vector bundles with integrable $\Q$-connection
\[\tag{\textbf{2.6.14}} \mathrm{T}_{\mathrm{dR}}(E \times_S E) \simeq \mathcal H_E.\]
Let us write $\zeta$ for the thus obtained $\mathcal D_{E/\Q}$-linear extension of $\mathcal O_E$ by $\mathcal H_E$, coming about by taking de Rham realizations with motivic Gauß-Manin connections in $(2.6.13)$ and using $(2.6.14)$.\\
\newline
From the fact that Deligne's duality behaves functorially\footnote{This means that for a morphism $t: M \rightarrow N$ of $1$-motives we have a commutative diagram
\begin{equation*}
\begin{xy}
\xymatrix{
\mathrm{T}_{\mathrm{dR}}(N^\vee) \ar[r]^{\sim} \ar[d]_{\mathrm{T}_{\mathrm{dR}}(t^\vee)} & \mathrm{T}_{\mathrm{dR}}(N)^\vee \ar[d]^{\mathrm{T}_{\mathrm{dR}}(t)^\vee} \\
\mathrm{T}_{\mathrm{dR}}(M^\vee)\ar[r]^{\sim} & \mathrm{T}_{\mathrm{dR}}(M)^\vee }
\end{xy}
\end{equation*}
with horizontal maps induced by the respective Deligne pairing. The statement can be viewed as a special case of \cite{An-Ber}, Rem. 2.7 (b) and (c) resp. the preceding Cor. 2.6.}, that it respects the motivic Gauß-Manin connections (cf. \cite{An-Ber}, Cor. 2.8) and that for the motive $[\mathbb Z_E \rightarrow 0]$ it just gives the identity on $\mathcal O_E$ (cf. \cite{Ber}, Ex. 4.4) we can hence conclude from Thm. 2.6.3:
\begin{corollary}
For the extension $\zeta \in \Ext^1_{\mathcal D_{E/\Q}}(\mathcal O_E,\mathcal H_E)$ induced by $(2.6.13)$ as described we have $\zeta= \mathcal Log^1$. \qquad \qed
\end{corollary}

\chapter{The explicit description on the universal elliptic curve}

\section{The birigidified Poincaré bundle for elliptic curves}
We explain how to obtain from the zero divisor of an elliptic curve an explicit construction of the birigidified Poincaré bundle $(\mathcal P^0,r^0,s^0)$ if one takes into account the self-duality of the curve.\\
\newline
Let $S$ be a locally noetherian scheme and $\pi: E\rightarrow S$ an elliptic curve with multiplication map $\mu : E\times_S E \rightarrow E$, zero section $\epsilon: S \rightarrow E$ and projections $\mathrm{pr}_1, \mathrm{pr}_2: E\times_S E \rightarrow E$.\\
Let $\widehat{E}$ denote the dual abelian scheme of $E$.\\
\newline
Recall that the $S$-scheme $\widehat{E}$ represents the dual functor of $E/S$ on the category of all $S$-schemes:
\[\tag{\textbf{3.1.1}} T \mapsto \mathrm{Pic}^0(E_T/T)=\textrm{\{Isomorphism classes of pairs} \ (\mathcal L, \alpha)\},\]
where $\mathcal L$ is a line bundle on $E_T=E\times_S T$ which is algebraically equivalent to zero and $\alpha$ a $T$-rigidification of $\mathcal L$ (cf. 0.1.1). The abelian group $\mathrm{Pic}^0(E_T/T)$ is canonically isomorphic to the group
\[\{[\mathcal L] \in \mathrm{Pic}(E_T)/\mathrm{Pic}(T)|\mathcal L \ \textrm{is algebraically equivalent to zero}\},\]
where $\mathrm{Pic}(T)$ becomes a subgroup of $\mathrm{Pic}(E_T)$ by pullback along the structure map $\pi_{E_T}: E_T \rightarrow T$ and $[\mathcal L]$ denotes the residue in the quotient group of the isomorphism class of a line bundle $\mathcal L$ on $E_T$; note that $[\mathcal L]=[\mathcal L']$ implies that $\mathcal L$ is algebraically equivalent to zero if and only if $\mathcal L'$ is.\\
The identification of $\mathrm{Pic}^0(E_T/T)$ with this group arises by mapping the class of $(\mathcal L, \alpha)$ to $[\mathcal L]$ and by conversely sending $[\mathcal L]$ to the class of $(\mathcal L \otimes_{\mathcal O_{E_T}} \pi_{E_T}^* \epsilon_{E_T}^* \mathcal L^{-1},\mathrm{can})$, where $\mathrm{can}$ denotes the canonical $T$-rigidification of $\mathcal L \otimes_{\mathcal O_{E_T}} \pi_{E_T}^* \epsilon_{E_T}^* \mathcal L^{-1}$. That one obtains well-defined homomorphisms which are inverse to each other is easy to check and indeed holds for $E$ replaced by any  abelian scheme over $S$.\\
\newline
We can now explain the self-duality of $E$: for this we define for each $S$-scheme $T$ a map
\[\tag{\textbf{3.1.2}} \begin{split} E(T) \xrightarrow{\sim} \{[\mathcal L] \in \mathrm{Pic}(E_T)/\mathrm{Pic}(T)|\mathcal L \ \textrm{is algebraically equivalent to zero}\},\\
Q \longmapsto \textrm{Class of} \ \mathcal O_{E_T}([-Q]-[0]). \qquad \qquad \qquad \qquad \qquad \qquad \qquad \ \ \end{split}\]
Here, note that for $Q \in E(T)$ its inverse $-Q \in E(T)$ defines a section $-Q$ of the abelian $T$-scheme $E_T$, which in turn induces an effective relative Cartier-Divisor $[-Q]$ of $E_T/T$ (cf. \cite{Kat-Maz}, Lemma 1.2.2.). Write $I([-Q])\subseteq \mathcal O_{E_T}$ for the associated invertible ideal sheaf and $\mathcal O_{E_T}([-Q]-[0])$ for $I([-Q])^{-1} \otimes_{\mathcal O_{E_T}} I([0])$, with $[0]$ the effective relative Cartier divisor of the zero section of $E_T/T$.\\
It is well-known that $(3.1.2)$ is an isomorphism of groups: this follows from \cite{Kat-Maz}, p. 64; cf. also \cite{Kat5}, p. 292. Note that our identification $(3.1.2)$ differs from the cited ones by the sign in $[-Q]$.\\
\newline
Combining $(3.1.1)$ and $(3.1.2)$ provides an isomorphism of $S$-group schemes
\[\tag{\textbf{3.1.3}} E \xrightarrow{\sim} \widehat{E}\]
which becomes the identification of \cite{Kat5}, p. 292, only after precomposition with the inverse map of $E$. As a side note, let us remark here that from a conceptual viewpoint $(3.1.3)$ is decisively the better self-duality isomorphism to use: we direct the interested reader to the discussion in \cite{Con2}, Ex. 2.5.\\
\newline
The identification $(3.1.3)$ will henceforth be fixed and referred to as the \underline{principal polarization} of $E/S$.\\
\newline
In the terminology of \cite{Ch-Fa}, p. 3-4, the morphism $(3.1.3)$ is the polarization $\lambda(\mathcal O_E([0])): E \rightarrow \widehat{E}$ associated with $\mathcal O_E([0])$. Under $\widehat{E}(E) \simeq \mathrm{Pic}^0(E\times_S E/E)$ it corresponds to the class of the pair
\[\tag{\textbf{3.1.4}} (\mathcal M \otimes_{\mathcal O_{E\times_S E}}(\pi \times \pi)^*\epsilon^*\mathcal O_E([0]), \mathrm{can}),\]
where $\mathcal M$ is the Mumford bundle for $\mathcal O_E([0])$ on $E\times_S E$:
\[\mathcal M:= \mu^* \mathcal O_E([0]) \otimes_{\mathcal O_{E\times_S E}} \mathrm{pr}_1^* \mathcal O_E([0])^{-1} \otimes_{\mathcal O_{E\times_S E}} \mathrm{pr}_2^*\mathcal O_E([0])^{-1}\]
and where $\mathrm{can}$ means the canonical rigidification along the second factor of $E \times_S E$.\\
Hence, now always identifying $E$ with $\widehat{E}$ via $(3.1.3)$, the pair $(3.1.4)$ represents the universal class in $\mathrm{Pic}^0(E\times_S E/E)$. Observe furthermore that $\mathcal M \otimes_{\mathcal O_{E\times_S E}}(\pi \times \pi)^*\epsilon^*\mathcal O_E([0])$ also carries a canonical rigidification along the first factor of $E\times_S E$, equally denoted by $\mathrm{can}$, which is obviously compatible with the canonical rigidification along the second factor after further restriction to $S$.\\
From the discussion in Rem. 0.1.8 we thus obtain that\\
\[\tag{\textbf{3.1.5}} (\mathcal M \otimes_{\mathcal O_{E\times_S E}}(\pi \times \pi)^*\epsilon^*\mathcal O_E([0]), \mathrm{can}, \mathrm{can})\]
is what we called in Def. 0.1.9 the birigidified Poincaré bundle on $E\times_S E$.

\markright{\uppercase{The explicit description on the universal elliptic curve}}
\section{Automorphy matrices for holomorphic vector bundles}
\markright{\uppercase{The explicit description on the universal elliptic curve}}
Throughout the subsequent sections we will use the yoga of automorphy matrices for vector bundles on complex manifolds. Such a matrix is obtained as soon as a trivialization for the pullback of the bundle to the universal covering exists and is chosen. One then has a convenient way of writing down the sections of the bundle as certain vectors of holomorphic functions on the universal covering.\\
We here give a brief self-contained account of the required techniques and fix conventions that will freely be used in the further proceeding. Despite slightly different priorities the material of this section is basically found in \cite{Ie}, 2 and 3, or obtained by generalizing \cite{Bi-Lan}, App. B, to vector bundles.
\newline

Let $X$ be a connected complex manifold with universal covering $p: \widetilde{X}\rightarrow X$ and write $\mathrm{Deck}(\widetilde{X}/X)$ for the group of deck transformations of $\widetilde{X}/X$.\\
\newline
(i) Assume we have a $\mathcal O_X$-vector bundle $\mathcal V$ of rank $n$ and a (henceforth fixed) trivialization
\[\mathcal O_{\widetilde{X}}^{\oplus n} \simeq p^*\mathcal V\]
with associated trivializing sections $\{e_1,...,e_n \} \in \Gamma(\widetilde{X}, p^*\mathcal V)$.\\
Let $\gamma \in \mathrm{Deck}(\widetilde{X}/X)$ and write for each $j=1,...,n$ the section $\gamma^*(e_j) \in \Gamma(\widetilde{X},p^*\mathcal V)$ as
\[\gamma^*(e_j)=\sum_{i=1}^n \varphi_{ij}^\gamma\cdot e_i\]
with uniquely determined $\varphi_{ij}^\gamma\in \Gamma(\widetilde{X},\mathcal O_{\widetilde{X}})$ which we collect in the matrix $\Big(\varphi_{ij}^\gamma\Big)_{\substack{i=1,...,n \\ j=1,...,n}}$.\\
Hence
\[(\gamma^*(e_1),...,\gamma^*(e_n))=(e_1,...,e_n)\cdot \Big(\varphi_{ij}^\gamma\Big)_{\substack{i=1,...,n \\ j=1,...,n}}.\]
We refer to the map
\[A: \mathrm{Deck}(\widetilde{X}/X) \times \widetilde{X} \rightarrow \mathrm{GL}_n(\C), \quad (\gamma, \widetilde{x}) \mapsto \Bigg(\Big(\varphi_{ij}^\gamma(\widetilde{x})\Big)_{\substack{i=1,...,n \\ j=1,...,n}}\Bigg)^{-1}\]
as the \underline{automorphy matrix} for $\mathcal V$ with respect to the (ordered) trivializing sections $\{e_1,...,e_n \}$ of $p^*\mathcal V$.\\
If $n=1$ we call the automorphy matrix the \underline{factor of automorphy} and use a small letter to denote it.\\
The automorphy matrix $A$ satisfies the relation
\[A(\gamma \cdot \gamma',\widetilde{x})=A(\gamma,\gamma'\cdot \widetilde{x})\cdot A(\gamma',\widetilde{x}) \quad \textrm{for all} \ \gamma, \gamma' \in \mathrm{Deck}(\widetilde{X}/X) \ \textrm{and} \ \widetilde{x} \in \widetilde{X}.\]
(ii) For an open subset $U$ of $X$ we will tacitly use the canonical identification
\[p^*:\Gamma(U,\mathcal V) \xrightarrow{\sim}\Gamma(p^{-1}(U),p^*\mathcal V)^{\mathrm{Deck}(\widetilde{X}/X)}:=\{s \in \Gamma(p^{-1}(U),p^*\mathcal V) \ | \ \gamma^*(s)=s \ \ \forall \ \gamma \in \mathrm{Deck}(\widetilde{X}/X) \}.\]
Then, if we express an element $s\in \Gamma(p^{-1}(U),p^*\mathcal V)$ in terms of the trivializing sections
\[s=\sum_{i=1}^n s_i\cdot e_i, \quad s_i\in \Gamma(p^{-1}(U),\mathcal O_{\widetilde{X}}),\]
it is invariant under $\mathrm{Deck}(\widetilde{X}/X)$ and hence defines a section of $\mathcal V$ over $U$ if and only if
\[\begin{pmatrix} (\gamma^*s_1)(\widetilde{x}) \\ \vdots \\ (\gamma^*s_n)(\widetilde{x}) \end{pmatrix}= A(\gamma,\widetilde{x})\cdot \begin{pmatrix} s_1(\widetilde{x}) \\ \vdots \\ s_n(\widetilde{x}) \end{pmatrix} \quad \textrm{for all} \ \gamma \in \mathrm{Deck}(\widetilde{X}/X), \ \widetilde{x} \in p^{-1}(U).\]
We will often directly write $s$ as the vector $\begin{pmatrix} s_1 \\ \vdots \\ s_n \end{pmatrix}$ if the underlying trivialization of $p^*\mathcal V$ is clear.\\
\newline
(iii) Consider the dual of the previously fixed trivialization of $p^*\mathcal V$:
\[\mathcal O_{\widetilde{X}}^{\oplus n} \simeq p^*\mathcal V^\vee\]
with associated trivializing sections $\{e_1^\vee,...,e_n^\vee \} \in \Gamma(\widetilde{X}, p^*\mathcal V^\vee)$. If
\[A: \mathrm{Deck}(\widetilde{X}/X) \times \widetilde{X} \rightarrow \mathrm{GL}_n(\C), \quad (\gamma, \widetilde{x}) \mapsto A(\gamma, \widetilde{x})\]
is the automorphy matrix for $\mathcal V$ with respect to $\{e_1,...,e_n \}$, then the automorphy matrix for $\mathcal V^\vee$ with respect to $\{e_1^\vee,...,e_n^\vee \}$ is given by
\[\big(A^t\big)^{-1}: \mathrm{Deck}(\widetilde{X}/X) \times \widetilde{X} \rightarrow \mathrm{GL}_n(\C), \quad (\gamma, \widetilde{x}) \mapsto \big(A(\gamma, \widetilde{x})^t\big)^{-1}.\]

(iv) Next, if we have a $\mathcal O_X$-vector bundle $\mathcal W$ of rank $m$ whose pullback to $\widetilde{X}$ is trivialized by the sections $\{e'_1,...,e'_m\}$, then the tensor product $\mathcal V\otimes_{\mathcal O_X} \mathcal W$ will be trivialized on $\widetilde{X}$ by the sections (in this order)
\[\{e_1\otimes e'_1, ...,e_1\otimes e'_m, e_2\otimes e'_1,...,e_2\otimes e'_m,......,e_n\otimes e'_1,...,e_n\otimes e'_m\}.\]
We will tacitly adopt this convention in the future. If the automorphy matrix of $\mathcal V$ resp. $\mathcal W$ is given by $A$ resp. $B$, then $\mathcal V\otimes_{\mathcal O_X} \mathcal W$ has the automorphy matrix $A \otimes B$, notation by which we mean the map
\[A\otimes B: \mathrm{Deck}(\widetilde{X}/X) \times \widetilde{X} \rightarrow \mathrm{GL}_{nm}(\C), \quad (\gamma,\widetilde{x})\mapsto A(\gamma,\widetilde{x})\otimes B(\gamma,\widetilde{x}),\]
where 
\[A(\gamma,\widetilde{x})\otimes B(\gamma,\widetilde{x}):=\begin{pmatrix}a_{11}(\gamma,\widetilde{x})\cdot B(\gamma,\widetilde{x}) & \dots & a_{1n}(\gamma,\widetilde{x})\cdot B(\gamma,\widetilde{x})\\
                             \vdots & \ddots & \vdots\\
														 a_{n1}(\gamma,\widetilde{x}) \cdot B(\gamma,\widetilde{x}) & \dots & a_{nn}(\gamma,\widetilde{x})\cdot B(\gamma,\widetilde{x})
 \end{pmatrix}\]
is the Kronecker product of $A(\gamma,\widetilde{x})=\Big(a_{ij}(\gamma,\widetilde{x})\Big)_{\substack{i=1,...,n \\ j=1,...,n}}$ with $B(\gamma,\widetilde{x})$.\\
\newline
(v) Let now $Y$ be another connected complex manifold with universal covering $q: \widetilde{Y}\rightarrow Y$.\\
Fix base points $\widetilde{y_0} \in  \widetilde{Y}$ resp. $\widetilde{x_0} \in \widetilde{X}$ with images $y_0$ resp. $x_0$ under $q$ resp. $p$.\\
Let $g: Y \rightarrow X$ be a holomorphic map with $g(y_0)=x_0$ and $\widetilde{g}: \widetilde{Y} \rightarrow \widetilde{X}$ the unique holomorphic map fitting into a commutative diagram
\begin{equation*}
\begin{xy}
\xymatrix{
(\widetilde{Y},\widetilde{y_0}) \ar[d]_{q} \ar[r]^{\widetilde{g}}& (\widetilde{X}, \widetilde{x_0}) \ar[d]^{p}\\
(Y,y_0)  \ar[r]^{g} & (X,x_0)}
\end{xy}
\end{equation*}
The standard isomorphisms $\mathrm{Deck}(\widetilde{X}/X) \simeq \pi_1(X,x_0), \mathrm{Deck}(\widetilde{Y}/Y) \simeq \pi_1(Y,y_0)$ and the canonical map $g_*:\pi_1(Y,y_0) \rightarrow \pi_1(X,x_0)$ provide a homomorphism $g_*:\mathrm{Deck}(\widetilde{Y}/Y) \rightarrow \mathrm{Deck}(\widetilde{X}/X)$ with
\[\widetilde{g}(\mu \cdot \widetilde{y})=g_*(\mu)\cdot \widetilde{g}(\widetilde{y}) \quad \textrm{for all} \ \mu \in \mathrm{Deck}(\widetilde{Y}/Y) \ \textrm{and} \ \widetilde{y}\in \widetilde{Y}.\]
By pullback along $\widetilde{g}$ we obtain from $\{e_1,...,e_n\}$ trivializing sections $\{f_1,...,f_n\}$ for $q^*g^*\mathcal V$.\\
If $\mu \in \mathrm{Deck}(\widetilde{Y}/Y), \gamma:=g_*(\mu)$ and for $j=1,...,n$:
\[\gamma^*(e_j)=\sum_{i=1}^n \varphi_{ij}^\gamma\cdot e_i,\]
then we have
\[\mu^*(f_j)=\sum_{i=1}^n(\widetilde{g})^*(\varphi_{ij}^\gamma)\cdot f_i.\]
In this way one obtains the automorphy matrix for $g^*\mathcal V$ with respect to $\{f_1,...f_n\}$ from the automorphy matrix for $\mathcal V$ with respect to $\{e_1,...,e_n\}$. If the underlying trivializing sections $\{e_1,...,e_n\}$ and the base points are fixed we will always work with the indicated trivialization of $q^*g^*\mathcal V$ and associated automorphy matrix for $g^*\mathcal V$.\\
If (as explained above) $s\in \Gamma(U, \mathcal V)$ is given with respect to $\{e_1,...,e_n\}$ by $\begin{pmatrix} s_1 \\ \vdots \\ s_n \end{pmatrix}, s_i \in \Gamma(p^{-1}(U),\mathcal O_{\widetilde{X}})$, satisfying
\[\begin{pmatrix} (\gamma^*s_1)(\widetilde{x}) \\ \vdots \\ (\gamma^*s_n)(\widetilde{x}) \end{pmatrix}= A(\gamma,\widetilde{x})\cdot \begin{pmatrix} s_1(\widetilde{x}) \\ \vdots \\ s_n(\widetilde{x}) \end{pmatrix} \quad \textrm{for all} \ \gamma \in \mathrm{Deck}(\widetilde{X}/X), \ \widetilde{x} \in p^{-1}(U).\]
then $g^*(s) \in \Gamma(g^{-1}(U),g^*\mathcal V)$ is given with respect to $\{f_1,...,f_n\}$ by $\begin{pmatrix} \widetilde{g}^*(s_1) \\ \vdots \\ \widetilde{g}^*(s_n) \end{pmatrix}$.

\section{The fundamental meromorphic Jacobi form and Eisenstein series}
\markright{\uppercase{The explicit description on the universal elliptic curve}}

\subsection{From canonical to classical theta functions}
A main goal of the future sections will consist in describing the analytified logarithm sheaves on the universal elliptic curve via automorphy matrices. The results of Chapter 2, the buildup of the Poincaré bundle $(3.1.5)$ and the techniques developed in 3.2 clearly suggest to carry this out first of all for the line bundle associated with the zero divisor of the curve. Trivializing this bundle on the universal covering corresponds to choosing a holomorphic function on the covering with the appropriate divisor, and hence the question arises which one we will want to choose.\\
The present subsection illustrates this question in the case of a single complex elliptic curve. Here, one possible trivialization is provided by the so-called canonical theta function associated to the zero divisor: In \cite{Ba-Ko-Ts}, 1, this function is the starting point for an explicit description of the polylogarithm on the curve minus its zero section; the induced trivializing section for the Poincaré bundle $(3.1.5)$ on the universal covering is the Kronecker theta function which is studied in depth in \cite{Ba-Ko}.\\
After a review of these functions in the context of 3.2 we point to the main problem of the canonical theta function: when varying the elliptic curve it no longer defines a holomorphic function and thus also no trivialization as before. As a substitute, we then present the so-called classical theta function for the zero divisor: it arises from the previous by multiplication with an exponential factor which erases the anti-holomorphic part but preserves a normalization property for its derivative by which it may also be characterized. Still working on a single curve we obtain the induced trivialization of the Poincaré bundle $(3.1.5)$ on the universal covering in the form of a meromorphic function $J$ which will turn out to play the main role in all subsequent sections. We finish this subsection with the fundamental observation that the function $J$, to which we were directed in the outlined natural way, coincides precisely with $2\pi i$-times the meromorphic Jacobi form introduced in \cite{Za2}, 3.

\subsubsection{The canonical theta function and factor of automorphy}

Fix a point $\tau \in \mathbb H$ in the upper half plane $\mathbb H$ of the complex numbers $\mathbb C$. Let $\Gamma_{\tau}:= \mathbb Z \tau \oplus \mathbb Z$ be the associated lattice in the complex plane and $E_{\tau}:=\mathbb C/\Gamma_{\tau}$ the complex torus defined by $\Gamma_{\tau}$. Moreover, set $A(\tau):=A(\Gamma_{\tau}):=\frac{1}{2\pi i}(\tau-\bar{\tau})$, which is the fundamental area of $\Gamma_{\tau}$ divided by $\pi$.\\
\newline
Consider the line bundle $\mathcal O_{E_{\tau}}([0])$ associated with the divisor given by the zero point $[0]$ of $E_{\tau}$. For its pullback along the canonical projection $p_{\tau}: \C \rightarrow E_{\tau}$ we have $p_{\tau}^*(\mathcal O_{E_{\tau}}([0])) \simeq \mathcal O_{\C}(\Gamma_{\tau})$.\\
Trivializing the previous line bundle amounts to giving a meromorphic function on $\C$ of divisor $-\Gamma_{\tau}$. Such a function is provided for example by $z\mapsto \frac{1}{\theta(z;\tau)}$ with
\[\tag{\textbf{3.3.1}} \theta(z;\tau):=\\\mathrm{exp}\bigg[-\frac{e_2^*(\tau)}{2}z^2\bigg]\cdot \sigma(z;\tau),\]
where $e_2^{*}(\tau):= \mathrm{lim}_{u \rightarrow 0^{+}} \sum_{\gamma \in \Gamma_{\tau} \backslash \{0\}}\gamma^{-2}|\gamma|^{-2u}$ is an Eisenstein-Kronecker number\footnote{Namely, in the notation of \cite{Ba-Ko}, Def. 1.5 resp. p. 238, it is the number $e^*_{0,2}(0,0;\Gamma_{\tau})$ resp. $e^*_{0,2}(\Gamma_{\tau})$.} and $\sigma(z;\tau)$ is the Weierstraß sigma function for the lattice $\Gamma_{\tau}$. Expressing the behaviour of $(3.3.1)$ under the deck transformations $\gamma \in \Gamma_{\tau}$ of $\C$ over $E_{\tau}$ as
\[\theta(z+\gamma;\tau)=a(\gamma, z) \cdot \theta(z;\tau),\]
then it is clear from the definition that
\[a: \Gamma_{\tau} \times \mathbb C \rightarrow \mathbb C^*\]
is the factor of automorphy for $\mathcal O_{E_{\tau}}([0])$ with respect to the trivializing section $z\mapsto \frac{1}{\theta(z;\tau)}$ of $p_{\tau}^*(\mathcal O_{E_{\tau}}([0]))$. This factor is well-known (cf. \cite{Ba-Ko}, Ex. 1.9), namely we have
\[\tag{\textbf{3.3.2}} a: \Gamma_{\tau} \times \mathbb C \rightarrow \mathbb C^*, \quad (\gamma,z)\mapsto \alpha(\gamma)\cdot \exp\bigg[\pi H(z,\gamma)+\frac{\pi}{2}H(\gamma,\gamma)\bigg],\]
where $H: \C\times \C \rightarrow \C$ is a hermitian form whose imaginary part is integral-valued on $\Gamma_{\tau} \times \Gamma_{\tau}$ and $\alpha: \Gamma_{\tau} \rightarrow \{z\in \C \ | \ |z|=1\}$ is a semicharacter for $H$, given explicitly by
\[H(z_1,z_2)=\frac{z_1 \bar{z}_2}{\pi A(\tau)} \quad \textrm{and} \quad \alpha (\gamma)= \begin{cases}
\ \ 1 & \text{for} \ \gamma \in 2  \Gamma_{\tau}\\
-1 & \ \textrm{otherwise.} \\
\end{cases}\]
\begin{remark}
The pair $(H,\alpha)$ is associated to $\mathcal O_{E_{\tau}}([0])$ via the Appell-Humbert theorem, the factor of automorphy $a$ in $(3.3.2)$ is the so-called \underline{canonical factor of automorphy} for $\mathcal O_{E_{\tau}}([0])$, and the function $z\mapsto \theta(z;\tau)$ in $(3.3.1)$ is a \underline{canonical theta function} for $\mathcal O_{E_{\tau}}([0])$. Details about this terminology can be found in \cite{Bi-Lan}, 2.2, 2.3 and 3.2, and also in \cite{Ba-Ko}, 1.2.\\
The function $z\mapsto \theta(z;\tau)$ is the unique holomorphic function on $\C$ with the property that its inverse defines a trivialization of $p_{\tau}^*(\mathcal O_{E_{\tau}}([0])) \simeq \mathcal O_{\C}(\Gamma_{\tau})$ which induces the canonical factor of automorphy for $\mathcal O_{E_{\tau}}([0])$ and such that its derivative at $z=0$ is normalized to $1$.
\end{remark}

Fix the base points $(0,0)\in \C\times \C$ and $0\in \C$. Then (according to the conventions recorded in 3.2 (iii)-(v)) we obtain a trivialization for the pullback of the Mumford bundle $\mathcal M_{\tau}$ of $\mathcal O_{E_{\tau}}([0])$ along the projection $\C\times \C \rightarrow E_{\tau}\times E_{\tau}=(\C\times \C)/(\Gamma_{\tau}\times \Gamma_{\tau})$ . With 3.2 (v) the associated factor of automorphy is straightforwardly computed from $(3.3.2)$ as
\[\tag{\textbf{3.3.3}} \Gamma_{\tau} \times \Gamma_{\tau} \times \C\times \C \rightarrow \C^*, \quad (\gamma, \mu, z,w) \mapsto \exp\bigg[\frac{\gamma\bar{\mu}+z\bar{\mu}+w\bar{\gamma}}{A(\tau)}\bigg].\]
One obtains the same factor for the Poincaré bundle $(3.1.5)$ on $E_{\tau}\times E_{\tau}$.\\
Note that $\mathcal M_{\tau}$ is given by $\mathcal  O_{E_{\tau} \times E_{\tau}}(\bar{\Delta}_{E_{\tau}}-([0]\times E_{\tau}) - (E_{\tau}\times [0]))$, where $\bar{\Delta}_{E_{\tau}}$ denotes the antidiagonal, that a trivialization of its pullback to $\C\times \C$ means giving a meromorphic function on $\C\times \C$ of divisor $-\{(z,w) \in \C\times \C \ | \ z+w\in \Gamma_{\tau} \}+(\Gamma_{\tau}\times \C) +(\C\times \Gamma_{\tau})$ and that the above trivialization then is the one defined by the function $(z,w) \mapsto \frac{1}{\Theta(z,w;\tau)}$, where
\[\tag{\textbf{3.3.4}} \Theta(z,w;\tau):=\frac{\theta(z+w;\tau)}{\theta(z;\tau)\theta(w;\tau)}\]
is the so-called \underline{Kronecker theta function} (cf. \cite{Ba-Ko}, 1.10).\\
Setting
\[F_1(z;\tau):=\mathrm{dlog}_z\theta(z;\tau)=\zeta(z;\tau)-e_2^*(\tau)\cdot z\]
we obtain from $\Theta(z,w;\tau)$ the function
\[\tag{\textbf{3.3.5}} \Xi(z,w;\tau):=\exp[-F_1(z;\tau)w] \cdot \Theta(z,w;\tau)\]
of \cite{Ba-Ko-Ts}, Def. 1.5. This function, in particular the coefficient functions obtained from it by Laurent expansion around $w=0$, lies at the center of the description given in \cite{Ba-Ko-Ts}, 1, for the polylogarithm on the single elliptic curve defined by $E_{\tau}$.\\
\newline
Observe that we have the equalities
\[\alpha(m\tau+n)=\exp[\pi i mn+\pi i m+\pi i n] \ \ \textrm{for all} \ m,n\in \Z \ \  \ \ \textrm{and} \  \quad e_2^*(\tau)=-\eta(1;\tau)-\frac{1}{A(\tau)},\]
where
\[\eta(m\tau +n ;\tau):=\zeta(z;\tau)-\zeta(z+m\tau+n;\tau) \quad \textrm{for all} \ m,n \in \mathbb Z\]
is the quasi-period defined via the Weierstraß zeta function $\zeta(z;\tau)$ for the lattice $\Gamma_{\tau}$; the formula for $e_2^*(\tau)$ may be deduced from the more general
\begin{lemma}
For each $\gamma \in \Gamma_{\tau}$ we have
\[e^*_2(\tau)\cdot \gamma=-\eta(\gamma;\tau)-\frac{\bar{\gamma}}{A(\tau)}.\]
\end{lemma}
\begin{proof}
By \cite{Ba-Ko-Ts}, 1.1, we have
\[F_1(z+\gamma;\tau)-F_1(z;\tau)=\frac{\bar{\gamma}}{A(\tau)},\]
from which the claim follows.
\end{proof}
With the formula for $\alpha(m\tau+n)$ and $\bar{\tau}=\tau-2\pi i \cdot A(\tau)$ resp. with the formula for $e_2^*(\tau)$ one verifies
\[a(m\tau+n,z)=\exp\bigg[\pi i m+\pi i n-\pi i m^2\tau+\frac{zm\tau}{A(\tau)}-2\pi i zm +\frac{nz}{A(\tau)}+\frac{(m\tau+n)^2}{2A(\tau)} \bigg]\]
resp.
\[\tag{\textbf{3.3.6}} \theta(z;\tau)=\exp\bigg[\frac{z^2}{2} \bigg (\frac{1}{A(\tau)}+\eta(1;\tau) \bigg)\bigg] \cdot \sigma(z;\tau),\]
such that moreover
\[\tag{\textbf{3.3.7}} \Theta(z,w;\tau)=\exp\bigg[zw \bigg (\frac{1}{A(\tau)}+\eta(1;\tau) \bigg)\bigg] \cdot \frac{\sigma(z+w;\tau)}{\sigma(z;\tau)\sigma(w;\tau)}\]
and
\[\tag{\textbf{3.3.8}} \Xi(z,w;\tau)=\exp[-\zeta(z;\tau)w] \cdot \frac{\sigma(z+w;\tau)}{\sigma(z;\tau)\sigma(w;\tau)}.\]
We see in particular that the function $\theta(z;\tau)$ of $(3.3.1)$, which we initially started with, does not vary holomorphically if the (so far fixed) parameter $\tau \in \H$ is moved: the obstruction comes from the antiholomorphic part of the normalized area function $A(\tau)$ in $(3.3.6)$. This presents very soon an overt problem when dealing with families of elliptic curves. Besides, also the fact that the naturally induced function $\Theta(z,w;\tau)$ needs to be altered in $(3.3.5)$ by an auxiliary exponential factor to obtain the central function $\Xi(z,w;\tau)$ in \cite{Ba-Ko-Ts}, 1, strongly hints at the following heuristic guideline:\\
\newline
\textit{\textbf{We should search from the beginning on a different trivialization for $p_{\tau}^*(\mathcal O_{E_{\tau}}([0])) \simeq \mathcal O_{\C}(\Gamma_{\tau})$.}}

\subsubsection{The classical theta function and factor of automorphy}
Instead of $\theta(z;\tau)$ let us consider - at first still for a fixed $\tau \in \H$ - the holomorphic function on $\C$:
\[\tag{\textbf{3.3.9}} z\mapsto \vartheta(z;\tau):=\exp\bigg[\frac{z^2}{2}\eta(1;\tau)\bigg]\cdot \sigma(z;\tau)=\exp\bigg[-\frac{z^2}{2A(\tau)}\bigg]\cdot \theta(z;\tau).\]
Its inverse $z \mapsto \frac{1}{\vartheta(z;\tau)}$ obviously provides again a trivialization for $p_{\tau}^*(\mathcal O_{E_{\tau}}([0])) \simeq \mathcal O_{\C}(\Gamma_{\tau})$, and its derivative in $z=0$ is again normalized to the value $1$. But now we have indeed defined a holomorphic function in both variables $(z,\tau) \in \C\times \H$; a more detailed analysis of $\vartheta$ will follow in 3.3.3.\\
\newline
The factor of automorphy for $\mathcal O_{E_{\tau}}([0])$
\[\widetilde{a}: \Gamma_{\tau} \times \mathbb C \rightarrow \mathbb C^*\]
induced by $z \mapsto \frac{1}{\vartheta(z;\tau)}$ is given by
\[\widetilde{a}(m\tau+n,z)=\frac{h(z)}{h(z+m\tau+n)}\cdot a(m\tau+n,z),\]
where $a$ is still the factor of automorphy defined by $z\mapsto \frac{1}{\theta(z;\tau)}$ and where
\[h(z):= \exp\bigg[\frac{z^2}{2A(\tau)}\bigg].\]
Explicitly, one computes
\[\frac{h(z)}{h(z+m\tau+n)}=\exp\bigg[-\frac{1}{2A(\tau)}\bigg((m\tau)^2+n^2+2zm\tau+2zn+2mn\tau\bigg)\bigg]\]
and hence
\[\tag{\textbf{3.3.10}} \widetilde{a}(m\tau+n,z)=\exp\bigg[\pi i m+\pi i n- 2 \pi i zm-\pi i m^2\tau \bigg].\]
\begin{remark}
(i) In the terms of \cite{Bi-Lan}, App. B, we have altered the cocycle in $\mathrm{Z}^1(\Gamma_{\tau},H^0(\mathcal O_{\C}^*))$ defined by the canonical factor of automorphy $a$ by the coboundary in $\mathrm{B}^1(\Gamma_{\tau},H^0(\mathcal O_{\C}^*))$ which is obtained from\\
$(m\tau+n,z)\mapsto \frac{h(z)}{h(z+m\tau+n)}$.\vspace{1mm}\\
(ii) The new factor of automorphy $\widetilde{a}$ in $(3.3.10)$ is exactly the so-called \underline{classical factor of automorphy} for the positive definite line bundle $\mathcal O_{E_{\tau}}([0])$ and its standard decomposition $\mathbb C=\mathbb R\cdot \tau\oplus \mathbb R \cdot 1$; the function $z\mapsto \vartheta(z;\tau)$ in $(3.3.9)$ then is a \underline{classical theta function} for $\mathcal O_{E_{\tau}}([0])$. Details about this terminology can be found in \cite{Bi-Lan}, 3.2.\\
The function $z\mapsto \vartheta(z;\tau)$ is the unique holomorphic function on $\C$ with the property that its inverse defines a trivialization of $p_{\tau}^*(\mathcal O_{E_{\tau}}([0])) \simeq \mathcal O_{\C}(\Gamma_{\tau})$ which induces the classical factor of automorphy for $\mathcal O_{E_{\tau}}([0])$ and such that its derivative at $z=0$ is normalized to $1$.
\end{remark}
Fixing as before $(0,0)\in \C\times \C$ and $0\in \C$ as base points we obtain (with the conventions of 3.2) a trivialization for the pullback to $\C\times \C$ of the Mumford bundle $\mathcal M_{\tau}$ associated with $\mathcal O_{E_{\tau}}([0])$ resp. for the pullback to $\C\times \C$ of the Poincaré bundle $(3.1.5)$ on $E_{\tau}\times E_{\tau}$. The associated factors of automorphy are computed (with 3.2 (v)) from $(3.3.10)$ both times as

\[{\small
\tag{\textbf{3.3.11}} \Gamma_{\tau} \times \Gamma_{\tau} \times \mathbb C \times \mathbb C \rightarrow \mathbb C^*, \  (m\tau+n,m'\tau+n',z,w)\mapsto \exp\bigg[-2\pi i mm'\tau -2\pi i m'z -2\pi i mw\bigg].
}
\]If we write again $\mathcal M_{\tau}$ as $\mathcal  O_{E_{\tau} \times E_{\tau}}(\bar{\Delta}_{E_{\tau}}-([0]\times E_{\tau}) - (E_{\tau}\times [0]))$, then the trivialization on $\C\times \C$ is defined by the function $(z,w) \mapsto \frac{1}{J(z,w;\tau)}$, where
\[\tag{\textbf{3.3.12}} J(z,w;\tau):=\frac{\vartheta(z+w;\tau)}{\vartheta(z;\tau)\vartheta(w;\tau)}=\exp[zw\cdot \eta(1;\tau)]\cdot \frac{\sigma(z+w;\tau)}{\sigma(z;\tau)\sigma(w;\tau)}.\]
\vspace{1mm}\\\textit{\textbf{The function $J$ in $(3.3.12)$ will be the key instrument for our explicit description of the $D$-variant of the elliptic polylogarithm in families. We have seen that it is the analogue for the Kronecker theta function $\Theta$ in $(3.3.4)$ when one performs the shift from the canonical to the classical factor of automorphy for $\mathcal O_{E_{\tau}}([0])$ (cf. Rem. 3.3.1 and Rem. 3.3.3). The classical theta function $\vartheta$ of $(3.3.9)$, as characterized in Rem. 3.3.3 (ii), instead of the canonical theta function $\theta$ of $(3.3.1)$ is the appropriate theta function to start with. The use of such a theta function is that it provides a factor of automorphy (and thus a way of writing down sections) for the Poincaré bundle $(3.1.5)$.}}\\
\newline
Let us conclude this more preliminary subsection by showing that with the function $J$ in $(3.3.12)$ we have not defined anything new to the existing literature.

\subsubsection{The relation to Zagier's meromorphic Jacobi form}
We keep working with a fixed $\tau \in \H$.\\
In \cite{Za2}, 3, Zagier introduces a meromorphic function $(z,w) \mapsto F_{\tau}(z,w)$ on $\C\times \C$ by analytic continuation of the function
\[(z,w) \mapsto \sum_{n\geq 0} \frac{\mathrm{e}^{-nw}}{\mathrm{e}^{-2\pi i n \tau}\cdot \mathrm{e}^z-1} - \sum_{m\geq 0} \frac{\mathrm{e}^{mz}\cdot \mathrm{e}^w}{\mathrm{e}^{-2\pi i m \tau}-\mathrm{e}^w}\]
which is defined on a certain domain of $\C \times \C$.\\
In the cited work the function $F_{\tau}(z,w)$ is shown to induce a generating function for the (adequately normalized) period polynomials of all Hecke eigenforms for the full modular group (cf. ibid., (17)). As explained in ibid., 3, Remark, $F_{\tau}(2\pi i z,2\pi i w)$ is a two-variable meromorphic Jacobi form.\\
\newline
A number of fundamental properties of the function $F_{\tau}(z,w)$, including a determination of its poles and residues and of its behaviour under modular transformations, are given in ibid., 3, Theorem.\\
The proof of the next lemma will make use of these results.
\begin{lemma}
With notation as above we have an equality of meromorphic functions in $(z,w) \in \C\times\C$:
\[\Theta(z,w;\tau)=2\pi i \cdot \exp\bigg[\frac{zw}{A(\tau)}\bigg] \cdot F_{\tau}(2\pi i z,2\pi i w).\]
\end{lemma}
\begin{proof}
As above let $\Gamma_{\tau}=\mathbb Z \tau \oplus \mathbb Z$ and write $2\pi i \Gamma_{\tau}:=\mathbb Z 2\pi i \tau \oplus \mathbb Z 2\pi i$. If $A(2\pi i \Gamma_{\tau})$ denotes the fundamental area divided by $\pi$ of the lattice $2\pi i \Gamma_{\tau}$, then we have
\begin{align*}		
		\tag{$*$} A(2\pi i \Gamma_{\tau}) &= 2\pi i \bar{\tau}-2\pi i \tau, \\
		\tag{$**$} A(2\pi i \Gamma_{\tau}) &= 4\pi^2 A(\tau).
\end{align*}
Let
\[\Theta(z,w;2\pi i \Gamma_{\tau}):=\frac{\theta(z+w;2\pi i \Gamma_{\tau})}{\theta(z;2\pi i\Gamma_{\tau})\theta(w;2\pi i \Gamma_{\tau})}\]
be the Kronecker theta function for the lattice $2\pi i \Gamma_{\tau}$, where the function $\theta(-;2\pi i \Gamma_{\tau})$ is defined as in $(3.3.1)$ by using the lattice $2\pi i \Gamma_{\tau}$ instead of $\Gamma_{\tau}$.\\
\newline
The canonical factor of automorphy for the Poincaré bundle on $(\mathbb C/2\pi i \Gamma_{\tau} \times \mathbb C /2\pi i \Gamma_{\tau})$ is given (by the same computation as in the deduction of $(3.3.3)$) as
\[2\pi i \Gamma_{\tau} \times 2\pi i \Gamma_{\tau} \times \C\times \C \rightarrow \C^*, \quad (\gamma, \mu, z,w) \mapsto \exp\bigg[\frac{\gamma\bar{\mu}+z\bar{\mu}+w\bar{\gamma}}{A(2\pi i \Gamma_{\tau})}\bigg].\]
Using $(*)$ in the numerator of this expression (and observing $\mathrm{e}^{2\pi i m'n}=1$) one calculates that it equals
\begin{align*}
&(2\pi i m\tau+2\pi in,2\pi i m'\tau+2\pi i n',z,w) \mapsto \exp\bigg[\frac{1}{A(2\pi i \Gamma_{\tau})}(4\pi^2mm'\tau^2+4\pi^2mn'\tau+4\pi^2m'n\tau\\
&+4\pi^2nn'-2\pi i m'z\tau-2\pi i mw\tau-2\pi i n'z-2\pi i nw) - 2\pi i mm'\tau-m'z-mw \bigg].
\end{align*}
Now set
\[g_{\tau}(z,w):= F_{\tau}(z,w)\cdot \exp\bigg[-\frac{zw}{A(2\pi i \Gamma_{\tau})}\bigg],\]
which is a meromorphic function in $(z,w)$ because this holds for $F_{\tau}(z,w)$ by \cite{Za2}, 3, Theorem, (ii).
Then
\[
\begin{split}
&g_{\tau}(z+2\pi i m\tau+2\pi i n, w+2\pi i m'\tau+2 \pi i n')\\
&=F_{\tau}(z+2\pi i m\tau+2\pi i n, w+2\pi i m'\tau+2 \pi i n') \cdot \exp\bigg[-\frac{zw}{A(2\pi i \Gamma_{\tau})}\bigg]\cdot \exp\bigg[\frac{1}{A(2\pi i \Gamma_{\tau})}(4\pi^2mm'\tau^2\\
&+4\pi^2mn'\tau+4\pi^2m'n\tau+4\pi^2nn'-2\pi i m'z\tau-2\pi i mw\tau-2\pi i n'z-2\pi i nw)\bigg],
\end{split}
\]
and the formula
\[F_{\tau}(z+2\pi i m\tau+2\pi i n, w+2\pi i m'\tau+2 \pi i n')=F_{\tau}(z,w)\cdot \exp\bigg[-2\pi i m m'\tau - m'z-mw\bigg]\]
thus implies that $(z,w) \mapsto g_{\tau}(z,w)$ defines a meromorphic section of the Poincaré bundle over $(\mathbb C/2\pi i \Gamma_{\tau} \times \mathbb C /2\pi i \Gamma_{\tau})$; for the preceding formula we use \cite{Za2}, 3, Theorem, (v).\\
\newline
The functions $g_{\tau}(z,w)$ and $\Theta(z,w; 2\pi i \Gamma_{\tau})$ moreover both have simple poles in $z=2\pi i m\tau+2\pi i n \ (m,n\in \mathbb Z)$ and in $w=2\pi i m'\tau+2\pi i n' \ (m',n'\in \mathbb Z)$ and are holomorphic elsewhere: for $F_{\tau}(z,w)$ - and hence for $g_{\tau}(z,w)$ - this is \cite{Za2}, 3, Theorem, (ii), and for $\Theta(z,w;2\pi i \Gamma_{\tau})$ this is clear by definition and by the zeroes of $\theta(-;2\pi i \Gamma_{\tau})$.\\
\newline
Let us calculate residues:\\
With \cite{Za2}, 3, Theorem, (ii), we obtain
\[\mathrm{Res}_{z=2\pi i m\tau+2\pi i n}\bigg \{\ g_{\tau}(z,w)\bigg \}\ =\exp[-wm]\cdot \exp\bigg[-\frac{(2\pi i m \tau+2\pi i n)w}{A(2\pi i \Gamma_{\tau})}\bigg].\]
On the other hand,
\begin{align*} \mathrm{Res}_{z=2\pi i m\tau+2\pi i n}\bigg \{\ \Theta(z,w; 2\pi i \Gamma_{\tau})\bigg \}\ &= \exp\bigg[\frac{(-2\pi im \bar{\tau}-2\pi i n)w}{A(2\pi i \Gamma_{\tau})}\bigg] \\
&=\exp\bigg[-\frac{2\pi i mw\tau +2\pi i nw}{A(2\pi i \Gamma_{\tau})}-mw\bigg],\end{align*}
where the first equation follows from \cite{Ba-Ko}, Lemma 1.15, and the second from $(*)$. In the same way one shows that the residues of $g_{\tau}(z,w)$ and $\Theta(z,w; 2\pi i \Gamma_{\tau})$ coincide also in $w=2\pi i m'\tau+2\pi i n'$.\\
\newline
In sum, we deduce that the difference of $g_{\tau}(z,w)$ and $\Theta(z,w;2\pi i \Gamma_{\tau})$ is a holomorphic section of the Poincaré bundle over $(\mathbb C/2\pi i \Gamma_{\tau} \times \mathbb C /2\pi i \Gamma_{\tau})$. By \cite{Ba-Ko}, Lemma 1.11, it must already be zero.\\
We hence get
\[\tag{$***$} \Theta(z,w;2\pi i \Gamma_{\tau})=g_{\tau}(z,w)=\exp\bigg[-\frac{zw}{A(2\pi i \Gamma_{\tau})}\bigg]\cdot F_{\tau}(z,w).\]
Now observe (e.g. by directly going into the definitions) the relation
\[\theta(z; 2\pi i \Gamma_{\tau})=2\pi i \cdot \theta \bigg(\frac{z}{2\pi i };\tau \bigg),\]
such that
\[\Theta(z,w;2\pi i \Gamma_{\tau})=\frac{1}{2\pi i} \cdot \Theta \bigg(\frac{z}{2\pi i },\frac{w}{2\pi i};\tau\bigg)\]
and hence
\[\begin{split} \Theta(z,w;\tau)&=2\pi i \cdot \Theta(2\pi i z,2\pi i w; 2\pi i \Gamma_{\tau})=2\pi i \cdot \exp\bigg[\frac{4\pi^2zw}{A(2\pi i \Gamma_{\tau})}\bigg]\cdot F_{\tau}(2 \pi i z,2\pi i w)\\
&=2\pi i \cdot \exp\bigg[\frac{zw}{A(\tau)}\bigg]\cdot F_{\tau}(2\pi i z,2\pi i w).\end{split}\]
Here, the second equality comes from $(***)$ and the third one is $(**)$. This establishes our claim.
\end{proof}

We now deduce the fundamental relation between the function $J(z,w;\tau)$ of $(3.3.12)$ and the meromorphic Jacobi form $F_{\tau}(2\pi iz, 2\pi i w)$:

\begin{proposition}
We have the equality
\[J(z,w;\tau)=2\pi i \cdot F_{\tau}(2\pi i z,2\pi i w).\]
\end{proposition}
\begin{proof}
As
\[\vartheta(z;\tau)=\exp\bigg[-\frac{z^2}{2A(\tau)}\bigg]\cdot \theta(z;\tau)\]
we have
\[J(z,w;\tau)=\exp\bigg[-\frac{zw}{A(\tau)}\bigg]\cdot \Theta(z,w;\tau).\]
Now use Lemma 3.3.4.
\end{proof}

\begin{remark}
We have the equality
\[J(z,w;\tau)=\exp[w\zeta(z;\tau)+zw\eta(1;\tau)]\cdot \Xi(z,w;\tau).\]
Namely, as $F_1(z;\tau)=\zeta(z;\tau)-e_2^*(\tau)\cdot z$ this follows by combining the last equation in the proof of Prop. 3.3.5 with Lemma 3.3.2.
\end{remark}

\subsection{Notations for some classical functions}

We briefly introduce notations for some well-known functions appearing in the complex theory of elliptic curves. The conventions adopted here will remain valid until the end of the work.\\
Details and basic properties concerning these functions can be found in \cite{Sil}, Ch. I, \cite{Kat3}, Ch. I, \cite{Kat4}, A 1.3, \cite{Po}, I, App. A, and \cite{Bi-Lan}, 8.5.\\
\newline
We use the abbreviations $q_z:=\mathrm{e}^{2\pi i z}, \ q_{\tau}:=\mathrm{e}^{2\pi i \tau}$.\\
\newline
(i) As in 3.3.1 we let
\[A(\tau):=\frac{1}{2\pi i}(\tau-\bar{\tau})\]
be the fundamental area of $\Gamma_{\tau}=\Z\tau\oplus \Z$ divided by $\pi$, now viewed as a function in $\tau \in \H$, as well as
\[\sigma(z,\tau):=z\cdot \prod_{\gamma \in \Gamma_{\tau} \backslash \{0\}}\bigg(1-\frac{z}{\gamma}\bigg)\cdot \exp\bigg[\frac{z}{\gamma}+\frac{z^2}{2\gamma^2}\bigg]\]
resp.
\[\zeta(z,\tau):=\frac{1}{z}+\sum_{\gamma \in \Gamma_{\tau}\backslash \{0\}} \bigg(\frac{1}{z-\gamma}+\frac{1}{\gamma}+\frac{z}{\gamma^2}  \bigg) \quad \quad \ \]
the Weierstraß sigma resp. zeta function, now viewed as functions in $(z,\tau)\in \C\times \H$.\\
\newline
(ii) For each $\tau \in \H$ the quasi-period $\eta(-,\tau): \Gamma_{\tau} \rightarrow \C$ for the lattice $\Gamma_{\tau}$ is again defined via
\[\eta(m\tau +n, \tau):=\zeta(z,\tau)-\zeta(z+m\tau+n,\tau).\]
We stress that our definition of the quasi-period follows the sign convention of \cite{Kat3}, 1.2.4 and \cite{Kat4}, A 1.3, which usually differs by the factor $(-1)$ from the definition in other sources.\\
The quasi-period in particular yields the functions in $\tau \in \H$ given by $\eta(1,\tau)$ and $\eta(\tau,\tau)$ which are connected via the Legendre relation (cf. \cite{Kat4}, A 1.3.4):
\[\eta(\tau,\tau)=2\pi i +\tau \cdot \eta(1,\tau).\]
We also remark the equality
\[-\eta(1,\tau)=G_2(\tau):=\sum_{n\neq 0}\frac{1}{n^2}+\sum_{m\neq 0}\sum_{n\in \Z}\frac{1}{(m\tau+n)^2} =\frac{\pi^2}{3}-8\pi^2\cdot \sum_{n \geq 1}\bigg(\sum_{\substack{
    d|n\\
    d>0}
  } d \bigg) q_{\tau}^n\]
of $-\eta(1,\tau)$ with the holomorphic Eisenstein series of weight two $G_2(\tau)$ (cf. \cite{Kat4}, Lemma A 1.3.9).\\
\newline
(iii) Furthermore, let
\[\eta(\tau):=\exp\bigg[\frac{2\pi i \tau}{24}\bigg]\cdot \prod_{n\geq 1}(1-q_{\tau}^n)\]
denote the Dedekind eta function in $\tau \in \H$.\\
\newline
(iv) We introduce the two-variable classical Riemann theta function of characteristic ${\begin{bmatrix}
   \frac{1}{2}\\ 
   \frac{1}{2}
\end{bmatrix}}$ by
\[\theta_{11}(z,\tau):=\sum_{n\in \mathbb Z}\exp\bigg[\pi i \bigg (n+\frac{1}{2} \bigg)^2\tau +2\pi i \bigg(n+\frac{1}{2} \bigg)\bigg(z+\frac{1}{2} \bigg)\bigg],\]
following in our notation \cite{Po}, I, App. A. It is also often denoted by $\vartheta \begin{bmatrix} 
   \frac{1}{2}\\ 
   \frac{1}{2} 
\end{bmatrix} (z,\tau)$, e.g. in \cite{Bi-Lan}, 8.5.\\
\newline
(v) Finally, we write
\[F(z,w,\tau):=F_{\tau}(z,w)\]
for the function defined in \cite{Za2}, 3, which we have already introduced prior to Lemma 3.3.4 and now view as a function in $(z,w,\tau) \in \C \times \C \times \H$.

\subsection{The elementary theta function and the fundamental meromorphic Jacobi form}
In 3.3.1 we explained at some length the motivation to consider the functions $\vartheta(z;\tau)$ and $J(z,w;\tau)$ which we now officially introduce with also the parameter $\tau \in \H$ varying. Before we explicitly apply them in the context of the universal elliptic curve we want to use the present and following subsection to record the most important analytic properties of these functions that will be needed later.

\subsubsection{The elementary theta function}

\begin{definition}
The \underline{elementary theta function} is the function in $(z,\tau) \in \mathbb C\times \mathbb H$ given by
\[\vartheta(z,\tau):=\exp\bigg[\frac{z^2}{2}\eta(1,\tau)\bigg]\cdot \sigma(z,\tau).\]
\end{definition}
\begin{remark}
As recorded in $(3.3.9)$ the relation between $\vartheta(z,\tau)$ and the function $\theta(z,\tau)$ of $(3.3.1)$ is given by
\[\tag{\textbf{3.3.13}} \vartheta(z,\tau)=\exp\bigg[-\frac{z^2}{2A(\tau)}\bigg]\cdot \theta(z,\tau).\]
\end{remark}
\begin{lemma}
We have the following two alternative expressions for the elementary theta function:
\[\tag{\textbf{3.3.14}} \vartheta(z,\tau)=\frac{1}{2\pi i }\bigg(q_z^{\frac{1}{2}}-q_z^{-\frac{1}{2}}\bigg) \cdot \prod_{n=1}^\infty \frac{(1-q_{\tau}^nq_z)(1-q_{\tau}^nq^{-1}_z)}{(1-q_{\tau}^n)^2}\]
and
\[\tag{\textbf{3.3.15}} \vartheta(z,\tau)=-\frac{\theta_{11}(z,\tau)}{2\pi \eta(\tau)^3}.\]
\end{lemma}
\begin{proof}
For $(3.3.14)$ use \cite{Sil}, I, §6, Thm. 6.4, and observe that the quasi-period used there differs from ours by a minus sign.\\
To derive $(3.3.15)$ use \cite{Po}, I, App. A, Thm. 3.9; note again the different sign of the quasi-period.
\end{proof}
As the functions $(z,\tau) \mapsto \theta_{11}(z,\tau)$ and $\tau \mapsto \eta(\tau)$ are holomorphic\footnote{The product in the Dedekind eta function converges absolutely and uniformly on compact subsets of $\mathbb H$, hence is holomorphic. For the holomorphicity of $\theta_{11}(z,\tau)$ in $\C\times \H$ set $c_1=c_2=\frac{1}{2}$ in \cite{Bi-Lan}, Prop. 8.5.4.} we see from $(3.3.15)$ that $\vartheta(z,\tau)$ varies holomorphically in $\C\times \H$. Its zeroes cut out the divisor $\{(m\tau+n,\tau)|\tau \in \mathbb H, m,n \in \mathbb Z \} \subseteq \C \times \H$, and its one-variable Taylor expansion around $z=0$ has the form
\[\tag{\textbf{3.3.16}} \vartheta(z,\tau)=z+ \ \textrm{higher terms}.\]
Finally, we will need the behaviour of $\vartheta(z,\tau)$ under modular transformations:
\begin{proposition}
For all $m,n \in \Z$ and $\begin{pmatrix} a & b \\ c & d \end{pmatrix} \in \mathrm{SL_2}(\Z)$ we have
\[\vartheta \bigg (\frac{z+m\tau+n}{c\tau+d},\frac{a\tau+b}{c\tau+d}\bigg)=\frac{1}{c\tau +d}\cdot \exp\bigg[\frac{\pi i c (z+m\tau+n)^2}{c\tau+d}+\pi i m+\pi i n - 2\pi i mz-\pi i m^2\tau\bigg] \cdot \vartheta(z,\tau).\]
\end{proposition}
\begin{proof}
Consider formula $(3.3.15)$:
\[\vartheta(z,\tau)=-\frac{\theta_{11}(z,\tau)}{2\pi \eta(\tau)^3}.\]
One checks that $-\theta_{11}(z,\tau)$ is equal to what in \cite{Si}, p. 30, is denoted by $\vartheta_1 (z,\tau)$. The transformation formulas for $\vartheta_1(z,\tau)$ and for $\eta(\tau)$ given in ibid., p. 34, then imply
\[\vartheta \bigg (\frac{z}{c\tau+d},\frac{a\tau+b}{c\tau+d}\bigg)=\frac{1}{c\tau +d}\cdot \exp\bigg[\frac{\pi i c z^2}{c\tau+d}\bigg] \cdot \vartheta(z,\tau).\]
For the proof of the proposition it thus remains to show that
\[\vartheta(z+m\tau+n,\tau)=\mathrm{e}^{\pi i m+ \pi  in -2\pi i m z-\pi i m^2\tau}\cdot \vartheta(z,\tau).\]
Taking into account Def. 3.3.7, the last formula is equivalent to the claim that $\sigma(z+m\tau+n,\tau)$ equals
\[\mathrm{e}^{\pi i m+\pi i n-2\pi i mz-\pi i m^2\tau-\eta(1,\tau)\cdot[\frac{1}{2}m^2\tau^2+\frac{1}{2}n^2+mz\tau+nz+mn\tau]}\cdot \sigma(z,\tau),\]
or, written differently, equals
\[
\begin{split}
&\mathrm{e}^{\pi i mn+ \pi  im +\pi i n}\cdot \mathrm{e}^{-\eta(1,\tau)\cdot [\frac{1}{2}m^2\tau^2+\frac{1}{2}n^2+mz\tau+nz+mn\tau]-2\pi i mz-\pi i m^2\tau-\pi i mn} \cdot \sigma(z,\tau)\\
&\qquad= \mathrm{e}^{\pi i mn+ \pi  im +\pi i n}\cdot \mathrm{e}^{[-\eta(1,\tau)\cdot(m\tau+n)-2\pi i m]\cdot [z+\frac{1}{2}n+\frac{1}{2}m\tau]} \cdot \sigma(z,\tau)\\
&\qquad=\mathrm{e}^{\pi i mn+ \pi  im +\pi i n}\cdot \mathrm{e}^{-\eta(m\tau+n,\tau)\cdot [z+\frac{1}{2}(m\tau+n)]} \cdot \sigma(z,\tau),
\end{split}
\]
where for the last line we have used the Legendre relation (cf. 3.3.2 (ii)):
\[\eta(\tau,\tau)=2\pi i + \tau\cdot \eta(1,\tau).\]
But that the last line is equal to $\sigma(z+m\tau+n, \tau)$ is well-known: this is \cite{Sil}, Ch. I, Prop. 5.4 (c): one only has to take into account that the quasi-period used there differs from ours by a minus sign and that the factor $\psi(m\tau+n)$ appearing there equals $\mathrm{e}^{\pi i mn+ \pi  im +\pi i n}$.
\end{proof}

\subsubsection{The fundamental meromorphic Jacobi form}

The holomorphic function $\vartheta(z,\tau)$ induces the following meromorphic function in three variables:
\begin{definition}
We define a meromorphic function in $(z,w,\tau) \in \C \times \C \times \H$ by setting
\[J(z,w,\tau):=\frac{\vartheta(z+w,\tau)}{\vartheta(z,\tau)\vartheta(w,\tau)}\]
and call it the \underline{fundamental meromorphic Jacobi form}.
\end{definition}

\begin{remark}
The relation between $J(z,w,\tau)$ and the function $\Theta(z,w,\tau)$ of $(3.3.4)$ is given by
\[\tag{\textbf{3.3.17}} J(z,w,\tau)=\exp\bigg[-\frac{zw}{A(\tau)}\bigg]\cdot \Theta(z,w,\tau).\]
\end{remark}
Let us also recall from $(3.3.12)$ that
\[\tag{\textbf{3.3.18}} J(z,w,\tau)=\exp[zw\cdot \eta(1,\tau)]\cdot \frac{\sigma(z+w,\tau)}{\sigma(z,\tau)\sigma(w,\tau)}.\]
In Prop. 3.3.5 we have already established the following relation between $J(z,w,\tau)$ and the meromorphic Jacobi form $F(2\pi i z,2\pi i w, \tau)$ of \cite{Za2}, 3:
\begin{Proposition}
We have the equality
\[J(z,w,\tau)=2\pi i \cdot F(2\pi i z,2\pi i w, \tau).\]
\qquad \qed
\end{Proposition}
In view of \cite{Za2}, 3, Theorem, (i)-(viii), this yields a number of non-trivial properties for $J(z,w,\tau)$.\\
\newline
The following transformation formula can be viewed as a corollary of Prop. 3.3.10 or be deduced from Prop. 3.3.13 and \cite{Za2}, 3, Theorem, (v) and (vi).
\begin{corollary}
For all $m,n,m',n' \in \Z$ and $\begin{pmatrix} a & b \\ c & d \end{pmatrix} \in \mathrm{SL_2}(\Z)$ we have
\begin{align*}
&J \bigg (\frac{z+m\tau+n}{c\tau+d}, \frac{w+m'\tau+n'}{c\tau+d},\frac{a\tau+b}{c\tau+d}\bigg)\cdot \big(J(z,w,\tau)\big)^{-1}\\
&=(c\tau +d)\cdot \exp\bigg[\frac{2\pi i c}{c\tau+d}\cdot(z+m\tau+n)(w+m'\tau+n')-2\pi i m'z-2\pi i mw-2\pi i mm'\tau\bigg].
\end{align*}
\qed
\end{corollary}
It follows e.g. from Prop. 3.3.13 and \cite{Za2}, 3, Theorem, (ii), that $J(z,w,\tau)$ has simple poles in
\[z=m\tau+n \ (m,n\in \Z, \tau \in \H)\]
with residue $\mathrm{e}^{-2\pi i mw}$ and simple poles in
\[w=m'\tau+n' \ (m',n'\in \Z, \tau \in \H)\]
with residue $\mathrm{e}^{-2\pi i m'z}$ and is holomorphic elsewhere.\\
We consider its Laurent expansion with respect to the variable $w$ around $w=0$:
\[\tag {\textbf{3.3.19}} J(z,w,\tau)=\frac{1}{w}+\sum_{k\geq 0} r_k(z,\tau)\cdot w^k.\]
From the mentioned knowledge about the residues of $J$ and comparison of expansions one easily deduces the following information about the meromorphic coefficient functions $r_k(z,\tau)$ for all $k\geq 0$, where we use the convention $0^0:=1$ and trivially remark that having residue $0$ at a pole of order at worst one means holomorphicity:

\[{\small
\tag{\textbf{3.3.20}}r_k \ \textrm{has at worst simple poles along} \ z=m\tau+n \ (m,n\in \Z, \tau \in \H), \textrm{with residue} \ \frac{(-1)^k(2\pi i m)^k}{k!}.
}
\]
From $(3.3.17)$, \cite{Ba-Ko-Ts}, p. 191, and Lemma 3.3.2 one finds
\[\tag{\textbf{3.3.21}} r_0(z,\tau)=\zeta(z,\tau)+\eta(1,\tau)\cdot z.\]
For later purposes we note as a consequence of $(3.3.19)$ and $(3.3.20)$ that for each $D\in \Z\backslash \{0\}$ we have a Laurent expansion around $w=0$:
\[\tag{\textbf{3.3.22}} D^2 \cdot J(z,-w,\tau)-D\cdot J\Big(Dz,-\frac{w}{D},\tau \Big)=s^D_0(z,\tau)+s^D_1(z,\tau)w+...,\]
where for all $k\geq 0$ the $s^D_k(z,\tau)$ are meromorphic functions on $\C \times \H$ with the property:

\[{\small
\tag{\textbf{3.3.23}}
\begin{split} &s^D_k \ \textrm{has at worst simple poles along} \ z=m\tau+n \ (m,n\in \Z, \tau \in \H), \textrm{with residue} \ (D^2-1)\cdot \frac{(2\pi i m)^k}{k!},\\
&\textrm{and along} \ z=\frac{m}{D}\tau+\frac{n}{D} \ (\textrm{with} \ D \ \textrm{not simultaneously dividing} \ m \ \textrm{and} \ n),  \textrm{with residue} \ -\frac{(2\pi i \frac{m}{D})^k}{k!}.
\end{split}
}
\]
From $(3.3.21)$ we see:
\[\tag{\textbf{3.3.24}} s^D_0(z,\tau)=D^2\cdot \zeta(z,\tau)-D\cdot \zeta(Dz,\tau).\]

\subsection{The fundamental meromorphic Jacobi form and Eisenstein series}
We continue the investigation of the function $J(z,w,\tau)$ by revealing that its Laurent expansion around $w=0$ involves as coefficients Eisenstein functions $e_k(z,\tau)$ which are obtained from the Eisenstein-Kronecker-Lerch series (resp. its analytic continuation) defined in \cite{Ba-Ko-Ts}, 2.1 or \cite{Ba-Ko}, 1.1.\\
It will turn out later that the analytified $D$-variant of the polylogarithm on the universal elliptic curve can be constructed from the $w$-expansion of $J(z,w,\tau)$. Hence, already in prevision of determining the specialization of the $D$-variant along torsion sections, we here evaluate the $e_k(z,\tau)$ at points $z=\frac{a}{N}\tau +\frac{b}{N}$. The result, which is also of independent interest, represents the obtained functions $\tau\mapsto e_k\Big(\frac{a}{N}\tau+\frac{b}{N},\tau\Big)$ by the modular forms $F^{(k)}_{\frac{a}{N},\frac{b}{N}}(\tau)$ defined and studied in \cite{Ka}, Ch. I.

\subsubsection{The two-variable Eisenstein functions}
For $\tau \in \mathbb H$ we will again write $\Gamma_{\tau}= \mathbb Z \tau \oplus \mathbb Z$ for the associated lattice in the complex plane and $A(\tau):=A(\Gamma_{\tau}):=\frac{1}{2\pi i}(\tau-\bar{\tau})$, which is the fundamental area of $\Gamma_{\tau}$ divided by $\pi$.\\
\newline
For each $k\geq 0$ we define a two-variable Eisenstein function by
\[e_k(z,\tau):=K^*_k(0,z,k;\tau),\]
where $K^*_k(-,-,-;\tau)$ denotes the Eisenstein-Kronecker-Lerch function (with asterisk) associated to the lattice $\Gamma_{\tau}$; for its definition and basic properties cf. \cite{Ba-Ko-Ts}, 2.1, or \cite{Ba-Ko}, 1.1. The $e_k(z,\tau)$ define $\mathscr{C}^{\infty}$-functions for $(z,\tau)$ in $\mathbb (\mathbb C\times \mathbb H) \backslash \{(m\tau+n,\tau)|\tau \in \mathbb H, m,n \in \mathbb Z \}$, see the argument below.\\
\newline
The Eisenstein functions have an important relation to the expansion of $J(z,w,\tau)$ in $w=0$; namely, we have the following equation (in which we don't write out the expansion of the exponential term):
\[\tag{\textbf{3.3.25}} J(z,w,\tau)=\exp \bigg[2\pi i \frac{\bar{z}w-zw}{\tau-\bar{\tau}} \bigg]\cdot \bigg( \frac{1}{w}+\sum_{k\geq 0 }(-1)^k \cdot e_{k+1}(z,\tau)  \cdot w^k\bigg).\]
This is deduced from the symmetry property $J(z,w,\tau)=J(w,z,\tau)$ together with the formulas
\begin{align*}
J(z,w,\tau) & =\exp \bigg[2\pi i \frac{z\bar{w}-zw}{\tau-\bar{\tau}} \bigg]\cdot K_1(z,w,1;\tau),\\
K_1(z,w,1;\tau) & =\bigg(\frac{1}{z}+\sum_{k\geq 0}(-1)^k\cdot e_{k+1}(w,\tau)\cdot z^k\bigg),
\end{align*}
where for the first equality one uses $(3.3.17)$ and \cite{Ba-Ko}, Thm. 1.13, and the second follows from the formula in \cite{Ba-Ko-Ts}, p. 226, noting ibid., Def. A.2 and Rem. A.5. Here, $K_1(z,w,1;\tau)$ is a certain Eisenstein-Kronecker-Lerch function (without asterisk) for $\Gamma_{\tau}$, as defined and studied in \cite{Ba-Ko}, 1.1.\\
\newline
From $(3.3.25)$ and the properties of the Laurent $w$-expansion of $J(z,w,\tau)$ (cf. $(3.3.19)$ and $(3.3.20)$) we see that the $e_k(z,\tau)$ are $\mathscr{C}^{\infty}$-functions for $(z,\tau)$ in $\mathbb (\mathbb C\times \mathbb H) \backslash \{(m\tau+n,\tau)|\tau \in \mathbb H, m,n \in \mathbb Z \}$.
\subsubsection{Specializing the Eisenstein functions to modular forms}

Fix $N\geq 3$ and two integers $a,b$ which are not simultaneously multiples of $N$, in other words: setting $\alpha:=\frac{a}{N}, \beta:=\frac{b}{N}$ we have $(\alpha,\beta)\neq (0,0)$ in $(\frac{1}{N}\mathbb Z/\mathbb Z)^2$.\\
\newline
We want to relate for each $k\geq 1$ the function $\tau \mapsto e_k\Big(\frac{a}{N}\tau + \frac{b}{N},\tau\Big)$ with the modular form $F^{(k)}_{\alpha,\beta}(\tau)$ of weight $k$ and level $N$ defined in \cite{Ka}, Ch. I, 3.6. For the reader's convenience, we quickly review the definition of these modular forms, for the present purpose understood in the classical analytic sense; for more details, in particular for their purely algebraic construction, we refer to ibid., Ch. I, 3.6-3.10. In fact, their algebraic origin won't play a role before 3.8.2, where we will also precisely describe the relation between the algebraic and the analytic approach (cf. Rem. 3.8.12 and Rem. 3.8.13). Until then, we will exclusively treat them as modular forms in the classical sense, working with their expressions as holomorphic functions in $\tau$ given in ibid., Ch. I, 3.8.\\
\newline
Namely, if we (now and in the following) write $\zeta_N:=\mathrm{e}^{\frac{2\pi i}{N}}$, then
\begin{align*}
F^{(k)}_{\alpha,\beta}(\tau) &:=N^{-k} \sum_{(x,y)\in (\mathbb Z/N\mathbb Z)^2} \zeta_N^{xb-ya}\cdot E^{(k)}_{\frac{x}{N},\frac{y}{N}}(\tau), \quad k\neq 2 \\
F^{(2)}_{\alpha,\beta}(\tau) &:=N^{-2} \sum_{(x,y)\in (\mathbb Z/N\mathbb Z)^2} \zeta_N^{xb-ya}\cdot \widetilde{E}^{(2)}_{\frac{x}{N},\frac{y}{N}}(\tau).
\end{align*}
The $E^{(k)}_{\frac{x}{N},\frac{y}{N}}(\tau)$ resp. $\widetilde{E}^{(2)}_{\frac{x}{N},\frac{y}{N}}(\tau)$, associated to $(\frac{x}{N},\frac{y}{N}) \in (\frac{1}{N}\mathbb Z/\mathbb Z)^2$, are modular forms of weight $k$ and level $N$ of algebraic origin. As functions in $\tau$ they are given as follows:\\
\newline
In the case $k\geq 3$ one has
\begin{align*}
E^{(k)}_{\frac{x}{N},\frac{y}{N}}(\tau) &=(-1)^k(k-1)!(2\pi i)^{-k}\sum_{(m,n)\in \mathbb Z^2}\frac{1}{(\frac{\widetilde{x}}{N}\tau+\frac{\widetilde{y}}{N}+m\tau+n)^k} \quad \textrm{for} \ (x,y) \neq (0,0) \ \textrm{in} \ (\mathbb Z/N\mathbb Z)^2,\\
E^{(k)}_{0,0}(\tau) &=(-1)^k(k-1)!(2\pi i)^{-k}\sum_{(m,n)\in \mathbb Z^2 \backslash \{(0,0)\}}\frac{1}{(m\tau+n)^k},
\end{align*}where $(\widetilde{x},\widetilde{y})$ represents $(x,y)$ in $(\mathbb Z/N\mathbb Z)^2$.\\
To obtain $E^{(k)}_{\frac{x}{N},\frac{y}{N}}(\tau)$, where $(x,y)\neq (0,0) \in (\mathbb Z/N\mathbb Z)^2,$ resp. $E^{(k)}_{0,0}(\tau)$ for the cases $k=1,2$ one proceeds by Hecke summation: this means that one takes the value at $s=0$ of the analytic continuation of the - in $\mathrm{Re}(s) >2-k$ absolutely convergent - series
\[ \ \ (-1)^k(k-1)!(2\pi i)^{-k}\sum_{(m,n)\in \mathbb Z^2}\frac{1}{(\frac{\widetilde{x}}{N}\tau+\frac{\widetilde{y}}{N}+m\tau+n)^k \cdot|\frac{\widetilde{x}}{N}\tau+\frac{\widetilde{y}}{N}+m\tau+n|^s}\]
resp.
\[(-1)^k(k-1)!(2\pi i)^{-k}\sum_{(m,n)\in \mathbb Z^2 \backslash \{(0,0)\}}\frac{1}{(m\tau+n)^k \cdot|m\tau+n|^s}. \ \ \quad \quad \qquad \qquad\]\\
\vspace{1mm}
$\widetilde{E}^{(2)}_{\frac{x}{N},\frac{y}{N}}(\tau)$ then is given by $E^{(2)}_{\frac{x}{N},\frac{y}{N}}(\tau)-E^{(2)}_{0,0}(\tau)$.\\
\newline
For more details and for proofs of tacitly assumed facts we refer the reader to \cite{Sc}, Ch. VII; see also the following remark.
\begin{remark}
Let $(x,y) \in (\mathbb Z/ N \Z)^2$. Then in \cite{Sc}, Ch. VII, for each $k\geq 1$ functions $G_{N,k,(x,y)}(\tau)$ are defined as the value at $s=0$ of the analytic continuation of the - in $\mathrm{Re}(s)>2-k$ absolutely convergent - series
\[s\mapsto \sum_{\substack{(m,n)\in \mathbb Z^2 \backslash \{(0,0)\}, \\ (m,n) \equiv (x,y) \ \textrm{mod} \ N}}\frac{1}{(m\tau+n)^k \cdot|m\tau+n|^s}.\]
One can check that the following relation with the above defined Eisenstein series holds:
\[G_{N,k,(x,y)}(\tau)= N^{-k}\frac{(-1)^k(2\pi i)^k}{(k-1)!} \cdot E^{(k)}_{\frac{x}{N},\frac{y}{N}}(\tau).\]
From \cite{Sc}, Ch. VII, (6) and (30), we also get the following formulas for any $(x,y)\in (\mathbb Z/ N\Z)^2 \backslash \{(0,0)\}$:
\begin{align*}
E^{(k)}_{\frac{x}{N},\frac{y}{N}}(\tau)& =(2\pi i)^{-k}\cdot \wp^{(k-2)} \Big(\frac{x}{N}\tau+\frac{y}{N},\tau \Big) \quad \textrm{if} \ k\geq 3, \\
\widetilde{E}^{(2)}_{\frac{x}{N},\frac{y}{N}}(\tau)& =(2\pi i)^{-2}\cdot \wp \Big(\frac{x}{N}\tau+\frac{y}{N},\tau \Big),
\end{align*}
with $\wp(z,\tau)$ resp. $\wp^{(k-2)}(z,\tau)$ the Weierstraß $\wp$-function resp. its $(k-2)$-th derivative in $z$-direction.
\end{remark}
We can now describe the specialization of the Eisenstein functions as follows:

\begin{theorem}
Let $N,a,b,\alpha,\beta$ be as above. Then for each $k\geq 1$ we have the equality
\[e_k \Big(\frac{a}{N}\tau+\frac{b}{N},\tau \Big)=\frac{(-1)^k(2\pi i)^k}{(k-1)!}\cdot F^{(k)}_{\alpha,\beta}(\tau).\]
\end{theorem}

\begin{proof}
We need to distinguish three cases.\\
\newline
Let us start with \underline{the case $k\geq 3$}:\\
\newline
Under this assumption we have the series representation

{\small
\[e_k(z,\tau)=\sum_{\gamma \in \Gamma_{\tau}\backslash \{0\}}\exp\bigg[2\pi i \frac{\bar{z}\gamma-z\bar{\gamma}}{\tau-\bar{\tau}}\bigg]\cdot \frac{1}{\gamma^k}=\sum_{(x,y)\in \mathbb Z^2\backslash \{(0,0)\}}\exp\bigg[2\pi i \frac{\bar{z}x\tau+\bar{z}y-zx\bar{\tau}-zy}{\tau-\bar{\tau}}\bigg]\cdot \frac{1}{(x\tau+y)^k},\]}following from the definition of $e_k$ and the fact that $\mathrm{Re}(k)>\frac{k}{2}+1$ holds for $k\geq 3$, such that we have an expression of $K^*_k(0,z,k;\tau)$ by the preceding absolute convergent series (cf. \cite{Ba-Ko-Ts}, Def. 2.1).\\
\newline
We get

{\small
\begin{align*}
& e_k \Big(\frac{a}{N}\tau + \frac{b}{N},\tau \Big) =\sum_{(x,y)\in \mathbb Z^2\backslash \{(0,0)\}}\zeta_N^{xb-ya}\cdot \frac{1}{(x\tau+y)^k}\\
& = \sum_{(x,y)\in (\mathbb Z/N\mathbb Z)^2\backslash \{(0,0)\}} \bigg \{\zeta_N^{xb-ya}\cdot \sum_{(m,n)\in \mathbb Z^2}\frac{1}{((\widetilde{x}+mN)\tau+\widetilde{y}+nN)^k}\bigg\}+\sum_{(m,n)\in \mathbb Z^2 \backslash \{(0,0)\}}\frac{1}{(mN\tau+nN)^k}\\
& = N^{-k} \bigg(\sum_{(x,y)\in (\mathbb Z/N\mathbb Z)^2\backslash \{(0,0)\}} \bigg \{\zeta_N^{xb-ya}\cdot \sum_{(m,n)\in \mathbb Z^2}\frac{1}{(\frac{\widetilde{x}}{N}\tau+\frac{\widetilde{y}}{N}+m\tau+n)^k}\bigg\}+\sum_{(m,n)\in \mathbb Z^2 \backslash \{(0,0)\}}\frac{1}{(m\tau+n)^k} \bigg).
\end{align*}
}

On the other hand, we have by definition

{\small
\begin{align*}
& F^{(k)}_{\alpha,\beta}(\tau)=N^{-k} \sum_{(x,y)\in (\mathbb Z/N\mathbb Z)^2} \zeta_N^{xb-ya}\cdot E^{(k)}_{\frac{x}{N},\frac{y}{N}}(\tau)\\
& =N^{-k} \bigg( \sum_{(x,y)\in (\mathbb Z/N\mathbb Z)^2\backslash \{(0,0)\}}\bigg \{\zeta_N^{xb-ya}\cdot (-1)^k(k-1)!(2\pi i)^{-k} \cdot \sum_{(m,n)\in \mathbb Z^2}\frac{1}{(\frac{\widetilde{x}}{N}\tau+\frac{\widetilde{y}}{N}+m\tau+n)^k} \bigg \} +E^{(k)}_{(0,0)}(\tau)\bigg)\\
&=\frac{(-1)^k(k-1)!}{(2\pi i)^k}N^{-k}\bigg(\sum_{(x,y)\in (\mathbb Z/N\mathbb Z)^2\backslash \{(0,0)\}}\bigg \{\zeta_N^{xb-ya}\cdot \sum_{(m,n)\in \mathbb Z^2}\frac{1}{(\frac{\widetilde{x}}{N}\tau+\frac{\widetilde{y}}{N}+m\tau+n)^k} \bigg \}\\
& +\sum_{(m,n)\in \mathbb Z^2 \backslash \{(0,0)\}}\frac{1}{(m\tau+n)^k} \bigg).
\end{align*}
}

Hence
\[e_k\Big(\frac{a}{N}\tau+\frac{b}{N},\tau \Big)=\frac{(-1)^k(2\pi i)^k}{(k-1)!}\cdot F^{(k)}_{\alpha,\beta}(\tau) \quad \textrm{for all} \ k\geq 3.\]\\
\vspace{1mm}
We now treat \underline{the case $k=1$}:\\
\newline
By definition and \cite{Ka}, p. 140, we have
{\small
\begin{align*}
& F^{(1)}_{\alpha,\beta}(\tau)=N^{-1} \bigg( \sum_{(x,y)\in (\mathbb Z/N\mathbb Z)^2\backslash \{(0,0)\}}\zeta_N^{bx-ay}\cdot E^{(1)}_{\frac{x}{N},\frac{y}{N}}(\tau)+E^{(1)}_{0,0}(\tau)\bigg)\\
&=-(2\pi iN)^{-1}\bigg(\sum_{(x,y)\in (\mathbb Z/N\mathbb Z)^2\backslash \{(0,0)\}}\zeta_N^{bx-ay}\cdot \bigg \{ \ \textrm{Anal. cont. of} \\
& \sum_{(m,n)\in \mathbb Z^2}\frac{1}{(\frac{\widetilde{x}}{N}\tau+\frac{\widetilde{y}}{N}+m\tau+n) \cdot |\frac{\widetilde{x}}{N}\tau+\frac{\widetilde{y}}{N}+m\tau+n|^s} \ \textrm{in} \ s=0 \bigg\} + \bigg \{ \ \textrm{Anal. cont. of} \\
&\sum_{(m,n)\in \mathbb Z^2 \backslash \{(0,0)\}}\frac{1}{(m\tau+n) \cdot |m\tau+n|^s}  \ \textrm{in} \ s=0 \bigg\} \bigg) =-(2\pi iN)^{-1} \bigg( \ \textrm{Anal. cont. of} \\
&\bigg\{\sum_{(x,y)\in (\mathbb Z/N\mathbb Z)^2\backslash \{(0,0)\}}\zeta_N^{bx-ay} \cdot \sum_{(m,n)\in \mathbb Z^2}\frac{1}{(\frac{\widetilde{x}}{N}\tau+\frac{\widetilde{y}}{N}+m\tau+n) \cdot |\frac{\widetilde{x}}{N}\tau+\frac{\widetilde{y}}{N}+m\tau+n|^s}\\
& + \sum_{(m,n)\in \mathbb Z^2 \backslash \{(0,0)\}}\frac{1}{(m\tau+n) \cdot |m\tau+n|^s}\bigg\} \ \textrm{in} \ s=0 \bigg).
\end{align*}
}
The series are absolutely convergent for $\mathrm{Re}(s)>1$, and there the expression in $\{...\}$ is equal to
{\small
\begin{align*}
& N^{1+s}\sum_{(x,y)\in (\mathbb Z/N\mathbb Z)^2\backslash \{(0,0)\}}\zeta_N^{bx-ay} \cdot \sum_{(m,n)\in \mathbb Z^2}\frac{1}{((\widetilde{x}+mN)\tau+\widetilde{y}+nN) \cdot |(\widetilde{x}+mN)\tau+\widetilde{y}+nN|^s} \\
& +N^{1+s}\sum_{(m,n)\in \mathbb Z^2 \backslash \{(0,0)\}}\frac{1}{(mN\tau+nN) \cdot |mN\tau+nN|^s}\\
& =N^{1+s}\sum_{(x,y)\in \mathbb Z^2\backslash \{(0,0)\}} \zeta_N^{bx-ay} \cdot \frac{1}{(x\tau+y)\cdot|x\tau+y|^s}=:f(s).
\end{align*}
}

On the other hand, we have for $\mathrm{Re}(s)>1$:

{\small
\[K_1^*\Big(0,\frac{a}{N}\tau+\frac{b}{N},1+\frac{s}{2}; \tau \Big)=\sum_{(x,y)\in \mathbb Z^2\backslash \{(0,0)\}} \zeta_N^{bx-ay} \cdot \frac{1}{(x\tau+y)\cdot|x\tau+y|^s},\]}

absolutely convergent (by applying the definition, cf. \cite{Ba-Ko-Ts}, Def. 2.1). We then have:
\[\tag{$*$} e_1(z,\tau)_{|z=\frac{a}{N}\tau+\frac{b}{N}}=\textrm{Anal. cont. of} \ K_1^* \Big(0,\frac{a}{N}\tau+\frac{b}{N},1+\frac{s}{2}; \tau \Big) \ \textrm{in} \ s=0,\]

again by definition.\\
The equality $N^{1+s} \cdot K_1^*(0,\frac{a}{N}\tau+\frac{b}{N},1+\frac{s}{2}; \tau)=f(s)$ in $\mathrm{Re}(s)>1$ implies the equality 
\[N \cdot K_1^*\Big(0,\frac{a}{N}\tau+\frac{b}{N},1;\tau \Big)=f(0)\]
for the analytic continuations in $s=0$. Now, $(*)$ together with the above computation for $F^{(1)}_{\alpha,\beta}(\tau)$ yields, as desired:
\[e_1 \Big(\frac{a}{N}\tau+\frac{b}{N},\tau \Big)=-2\pi i \cdot F^{(1)}_{\alpha,\beta}(\tau).\]
\\
\vspace{1mm}
We finally turn to \underline{the case $k=2$}:\\
\newline
By definition and \cite{Ka}, p. 140, we have:

{\small
\begin{align*}
F^{(2)}_{\alpha,\beta}(\tau) & =N^{-2} \bigg( \sum_{(x,y)\in (\mathbb Z/N\mathbb Z)^2}\zeta_N^{bx-ay}\cdot \widetilde{E}^{(2)}_{\frac{x}{N},\frac{y}{N}}(\tau)\bigg)\\
& =N^{-2} \bigg( \sum_{(x,y)\in (\mathbb Z/N\mathbb Z)^2\backslash \{(0,0)\}}\zeta_N^{bx-ay}\cdot \bigg[E^{(2)}_{\frac{x}{N},\frac{y}{N}}(\tau)-E^{(2)}_{0,0}(\tau)\bigg]\bigg)\\
& =-(4\pi^2N^2)^{-1}\bigg( \sum_{(x,y)\in (\mathbb Z/N\mathbb Z)^2\backslash \{(0,0)\}}\zeta_N^{bx-ay}\cdot \bigg[E\Big(2,\tau,\frac{\widetilde{x}}{N}\tau+\frac{\widetilde{y}}{N},0 \Big)-E_{0,0}(2,\tau,0)\bigg]       \bigg)\\
& =-(4\pi^2N^2)^{-1}\bigg( \sum_{(x,y)\in (\mathbb Z/N\mathbb Z)^2\backslash \{(0,0)\}}\zeta_N^{bx-ay}\cdot E\Big(2,\tau,\frac{\widetilde{x}}{N}\tau+\frac{\widetilde{y}}{N},0\Big)+E_{0,0}(2,\tau,0) \bigg)\\
& +(4\pi^2N^2)^{-1}\bigg( \sum_{(x,y)\in (\mathbb Z/N\mathbb Z)^2\backslash \{(0,0)\}}\zeta_N^{bx-ay}\cdot E_{0,0}(2,\tau,0)+E_{0,0}(2,\tau,0)\bigg).
\end{align*}
}

The first summand is by definition

{\small
\begin{align*}
& -(4\pi^2N^2)^{-1} \bigg( \ \textrm{Anal. cont. of} \ \bigg\{\sum_{(x,y)\in (\mathbb Z/N\mathbb Z)^2\backslash \{(0,0)\}}\zeta_N^{bx-ay} \cdot \sum_{(m,n)\in \mathbb Z^2}\frac{1}{(\frac{\widetilde{x}}{N}\tau+\frac{\widetilde{y}}{N}+m\tau+n)^2 \cdot |\frac{\widetilde{x}}{N}\tau+\frac{\widetilde{y}}{N}+m\tau+n|^s}\\
& +\sum_{(m,n)\in \mathbb Z^2 \backslash \{(0,0)\}}\frac{1}{(m\tau+n)^2 \cdot |m\tau+n|^s}\bigg\} \ \textrm{in} \ s=0 \bigg),
\end{align*}
}

and the term in $\{...\}$ is absolutely convergent in $\mathrm{Re}(s)>0$, where it equals
\[N^{2+s}\sum_{(x,y)\in \mathbb Z^2\backslash \{(0,0)\}}\zeta_N^{bx-ay}\cdot \frac{1}{(x\tau+y)^2\cdot|x\tau+y|^s}=:g(s).\]
On the other hand, by definition:
\[\tag{$**$}e_2(z,\tau)_{|z=\frac{a}{N}\tau+\frac{b}{N}}=\textrm{Anal. cont. of} \ K_2^* \Big (0,\frac{a}{N}\tau+\frac{b}{N},2+\frac{s}{2}; \tau \Big) \ \textrm{in} \ s=0.\]
For $\mathrm{Re}(s)>0$ we have
\[K_2^*\Big(0,\frac{a}{N}\tau+\frac{b}{N}, 2+\frac{s}{2};\tau \Big)=\sum_{(x,y)\in \mathbb Z^2\backslash \{(0,0)\}}\zeta_N^{bx-ay}\cdot \frac{1}{(x\tau+y)^2\cdot|x\tau+y|^s}=N^{-2-s}g(s),\]absolutely convergent (cf. \cite{Ba-Ko-Ts}, Def. 2.1).\\
Altogether, we obtain from $(**)$ and the above computation of $F^{(2)}_{\alpha,\beta}$:\\
\[\tag{$***$} e_2 \Big(\frac{a}{N}\tau+\frac{b}{N},\tau \Big)=-4\pi^2F^{(2)}_{\alpha,\beta}(\tau)+N^{-2}\sum_{(x,y)\in (\mathbb Z/N\mathbb Z)^2}\zeta_N^{bx-ay}\cdot E_{0,0}(2,\tau,0).\]
But
\[\sum_{(x,y)\in (\mathbb Z/N\mathbb Z)^2}\zeta_N^{bx-ay}=\bigg(\sum_{x\in \mathbb Z/N\mathbb Z}\exp\bigg[\frac{2\pi i}{N}bx\bigg]\bigg) \cdot \bigg(\sum_{y\in \mathbb Z/N\mathbb Z}\exp\bigg[-\frac{2\pi i}{N}ay\bigg] \bigg).\]
By hypothesis, we have that $N \nmid a$ or $N\nmid b$. If $N\nmid a$ we get
\[\sum_{y\in \mathbb Z/N\mathbb Z}\exp\bigg[-\frac{2\pi i}{N}ay\bigg]=\frac{1-\bigg(\exp\bigg[-\frac{2\pi i}{N}a \bigg]\bigg)^N}{1-\exp\bigg[-\frac{2\pi i}{N}a\bigg]} =0,\]
and if $N\nmid b$ we get
\[\sum_{x\in \mathbb Z/N\mathbb Z}\exp\bigg[\frac{2\pi i}{N}bx\bigg]=\frac{1-\bigg(\exp\bigg[\frac{2\pi i}{N}b \bigg]\bigg)^N}{1-\exp\bigg[\frac{2\pi i}{N}b\bigg]} =0.\]
With this, $(***)$ writes as
\[e_2 \Big (\frac{a}{N}\tau+\frac{b}{N},\tau \Big)=-4\pi^2F^{(2)}_{\alpha,\beta}(\tau),\]
which is what we wanted to show.
\end{proof}
Thm. 3.3.16 and $(3.3.25)$ will serve as an important computational tool when determining the specialization of the $D$-variant of the polylogarithm along torsion sections of the universal elliptic curve.\\
More precisely, we will see that this specialization expresses in terms of the modular forms $_DF^{(k)}_{\alpha,\beta}(\tau)$ defined in \cite{Ka}, Ch. I, p. 143; they are obtained from the above $F^{(k)}_{\alpha,\beta}(\tau)$ in the following way:
\begin{definition}
Let $N,a,b,\alpha,\beta$ be as above and let $D$ be an integer with $(D,N)=1$.\\
Then we define for each $k\geq 1$ the modular form $_DF^{(k)}_{\alpha,\beta}(\tau)$ of weight $k$ and level $N$ by setting
\[_DF^{(k)}_{\alpha,\beta}(\tau):=D^2F^{(k)}_{\alpha,\beta}(\tau)-D^{2-k}F^{(k)}_{D\alpha,D\beta}(\tau).\]
\end{definition}

\section{The analytic geometry of the basic objects}
\markright{\uppercase{The explicit description on the universal elliptic curve}}
The main goals for the rest of the work will be the following: to describe analytically the $D$-variant of the polylogarithm for the universal family of elliptic curves with level $N$-structure, and to determine its (algebraic) specialization along torsion sections.\\
In the following, we first review the definition of the universal elliptic curve over the modular curve of level $N$ (introduced as schemes over $\Q$) and the description of their associated complex manifolds. We also derive a suitable expression for the analytification of the universal vectorial extension of the dual elliptic curve. Then, as already illustrated in 3.3.1 for a single complex elliptic curve, the elementary theta function resp. the fundamental meromorphic Jacobi form is used to trivialize the pullback of the analytified Poincaré bundle to the universal covering, which provides us with a factor of automorphy. Expressing in this way sections of the Poincaré bundle via holomorphic functions on the universal covering, we finally deduce an explicit formula for its universal integrable connection.

\subsubsection{Conventions and notations}
Assume that $X$ is a scheme which is locally of finite type over $\C$ or over $\Q$.\\
We will denote by $X^{an}$ the complex analytic space associated with the $\C$-valued points $X(\C)=\mathrm{Hom}_{\C}(\mathrm{Spec(\C)}, X)$ resp. $(X\times_{\Q}\C)(\C)=\mathrm{Hom}_{\C}(\mathrm{Spec(\C)},X\times_{\Q}\C)=\mathrm{Hom}_{\Q}(\mathrm{Spec(\C)},X)$.\\
The assignment $X\mapsto X^{an}$ is functorial for morphisms in the respective category of schemes.\\
If $\mathcal F$ is a $\mathcal O_X$-module we write $\mathcal F^{an}$ for the $\mathcal O_{X^{an}}$-module given by pullback of $\mathcal F$ via the canonical map of locally ringed spaces $X^{an}\rightarrow X$; as this map is flat we get an exact functor $\mathcal F \mapsto \mathcal F^{an}$ which preserves coherence by Oka's theorem and \cite{EGAI}, Ch. 0, $(5.3.11)$.\\
For more details about the transition from the algebraic to the analytic category cf. \cite{SGA1}, exp. XII.\\
\newline
Furthermore, in view of 3.2 (v) we want to choose base points for the universal coverings of various manifolds (resp. their connected components) which will be encountered. Here is our convention:\\
\newline
\textbf{\textit{Until the end of the work we fix the base point $0$ for $\C$, the base point $i$ for $\H$ and the induced base points for factors, e.g. $(0,i)$ for $\C\times \H$ or $(0,0,0,i)$ for $\C \times \C^2\times \H$ etc.}}\\
\newline
From now on let an integer $N \geq 3$ be given.
\subsubsection{The universal elliptic curve with level $N$-structure}
Consider the (contravariant) set-valued functor on the category of $\Q$-schemes
\[\mathscr{S} \mapsto \textrm{\{Iso classes of pairs} \ (\mathscr{E}, \alpha) \ | \ \mathscr{E}/\mathscr{S} \ \textrm{elliptic curve,} \ \alpha:\big( \Z/N\Z \big)^2_{\mathscr{S}} \xrightarrow{\sim} \mathscr{E}[N] \ \textrm{iso of} \ \mathscr{S}\textrm{-groups}\}, \]
where $\big( \Z/N\Z \big)^2_{\mathscr{S}}$ means the constant $\mathscr{S}$-group scheme associated with the abstract group $\big( \Z/N\Z \big)^2$.\\
\newline
An isomorphism $\alpha$ as in the definition of the functor, called "level $N$-structure", is tantamount to give an ordered pair of $N$-torsion sections of $\mathscr{E}/\mathscr{S}$ inducing on each geometric fiber over $\mathscr{S}$ a basis for the (usual) $\big( \Z/N\Z \big)$-module of $N$-torsion points: we obtain this pair of sections as the images of $(1,0)$ and $(0,1)$ under the homomorphism $\big( \Z/N\Z \big)^2 \rightarrow \mathscr{E}[N](\mathscr{S})$ of abstract groups corresponding to $\alpha$ and call it the associated "Drinfeld basis" for $\mathscr{E}[N]$.\footnote{To see the equivalence between the isomorphy of $\alpha$ and the condition about the geometric fibers one may use \cite{Kat-Maz}, Prop. 1.10.12, (1.10.5) and Lemma 1.8.3.}\\
\newline
It is a well-known fundamental theorem that the above functor is representable by a $1$-dimensional affine scheme $S$ which is smooth, separated and of finite type over $\Spec(\Q)$, the (open) modular curve of level $N$ (cf. \cite{Kat-Maz}, Cor. 4.7.2, $(4.3)$, $(4.13)$ and $(1.2.1)$). One can check that $S$ is irreducible.\footnote{One may use e.g. the remarks at the beginning of \cite{Hi}, 2.9.3, together with \cite{Li}, Ch. 4, Prop. 3.8.}\\
We write
\[\pi: E \rightarrow S, \quad (e_1,e_2) \in E[N](S)\]
for the universal elliptic curve over the modular curve and the related Drinfeld basis for $E[N]$.\\
\newline
We want to have explicit expressions for the complex manifolds $S^{an}$ and $E^{an}$.\\
\newline
For this let $\Gamma(N):=\ker(\mathrm{SL_2}(\Z) \rightarrow \mathrm{SL_2}(\Z/N\Z))$ act as a group of automorphisms on $\H$ by linear fractional transformations. The action is properly discontinuous and as $N\geq 3$ it is also free; the orbit space $\Gamma(N) \backslash \H$, endowed with the quotient topology, naturally becomes a complex manifold.\\
Moreover, one checks that the set $\Z^2 \times \Gamma (N)$ obtains the structure of a group by
\[\Bigg(\begin{pmatrix} 
   m\\ 
   n 
\end{pmatrix}, \gamma\Bigg) \circ \Bigg(\begin{pmatrix} 
   m'\\ 
   n' 
\end{pmatrix}, \gamma'\Bigg):=\Bigg(\begin{pmatrix} 
   m'\\ 
   n'
\end{pmatrix}+(\gamma')^t \begin{pmatrix} 
   m\\ 
   n 
\end{pmatrix}, \gamma \gamma'\Bigg)\]
and that it then acts properly discontinuously and freely as a group of automorphisms on $\C\times \H$ via
\[\Bigg(\begin{pmatrix} 
   m\\ 
   n 
\end{pmatrix}, \begin{pmatrix} a & b \\ c & d \end{pmatrix}\Bigg) \cdot (z,\tau):=\bigg(\frac{z+m\tau+n}{c\tau+d},\frac{a\tau+b}{c\tau+d}\bigg).\]
We can thus again naturally form the quotient manifold $(\Z^2 \times \Gamma(N)) \backslash (\C\times \H)$.\\
\newline
In the following, we usually won't indicate it in the notation when actually working with orbits and will simply write down representatives, tacitly implying well-definedness of everything we do.\\
\newline
Because of $N\geq 3$ the natural projection
\[\tag{\textbf{3.4.1}} (\Z/N\Z)^*\times (\Z^2 \times \Gamma(N)) \backslash (\C \times \H) \xrightarrow{\mathrm{pr}} (\Z/N\Z)^*\times \Gamma(N) \backslash \H \]
with the section
\[\tag{\textbf{3.4.2}} \qquad \qquad \ \ (j,0,\tau) \mapsfrom (j,\tau)\]
defines an analytic family of elliptic curves.\footnote{By an analytic family of elliptic curves we understand a proper flat morphism of analytic spaces together with a section, satisfying that each fiber is a compact Riemann surface of genus $1$ which is then viewed as complex elliptic curve via the distinguished point induced by the section. In our situation, $(3.4.1)$ and $(3.4.2)$ define such a family: for every $(j,\tau) \in (\Z/N\Z)^*\times \Gamma(N) \backslash \H$ the fiber $\mathrm{pr}^{-1}((j,\tau))$ is (non-uniquely) isomorphic to the complex torus $\C/(\Z \tau \oplus \Z)$, e.g. by the map
\[\C/(\Z \tau \oplus \Z) \xrightarrow{\sim} \mathrm{pr}^{-1}((j,\tau)), \quad z \mapsto (j,z,\tau).\]
Here it is essential that $\Gamma(N)$ acts without fixed points on $\H$, i.e. that we have $N\geq 3$.} Furthermore, the two ordered sections
\[\tag{\textbf{3.4.3}} \Big(j,\frac{j\tau}{N},\tau \Big)         \mapsfrom   (j,\tau), \qquad \Big(j,\frac{1}{N},\tau \Big)  \mapsfrom    (j,\tau),\]
define an analytic Drinfeld basis for the $N$-torsion of this family.\footnote{For an analytic family of elliptic curves this means that the two sections fiberwise yield a basis for the $N$-torsion of the respective complex elliptic curve. For our situation this is clear by using the isomorphism of the previous footnote.}\\
For each $(j,\tau) \in (\Z/N\Z)^*\times \Gamma(N) \backslash \H$ the complex elliptic curve $\mathrm{pr}^{-1}((j,\tau))$, equipped with the level $N$-structure induced by this Drinfeld basis, identifies with $\Big(\C/(\Z \tau \oplus \Z), \frac{j\tau}{N}, \frac{1}{N}\Big)$ by
\[\C/(\Z \tau \oplus \Z) \xrightarrow{\sim} \mathrm{pr}^{-1}((j,\tau)),\quad z\mapsto (j,z,\tau),\]
and this is the only such isomorphism because of $N\geq 3$.\\
\newline
$(3.4.1)$-$(3.4.3)$ parametrizes all analytic families of elliptic curves with Drinfeld basis for the $N$-torsion (cf. \cite{Ha}, Thm. 5.2.29\footnote{One easily checks that our quotient object identifies canonically with the curve $\widetilde{\mathscr{E}}_N$ considered in the cited reference.}). In particular, applying this to the analytification of $(E/S, e_1,e_2)$ we obtain unique holomorphic maps constituting the vertical arrows of a commutative diagram
\begin{equation*} \tag{\textbf{3.4.4}} \begin{split}
\begin{xy}
\xymatrix{
E^{an} \ar[r]\ar[d]_{\pi^{an}} & (\Z/N\Z)^*\times (\Z^2\times \Gamma(N)) \backslash (\C \times \H) \ar[d]^{\mathrm{pr}} \\
S^{an} \ar[r]^{\Phi \qquad \qquad} & (\Z/N\Z)^* \times \Gamma(N) \backslash \H}
\end{xy}
\end{split}
\end{equation*}
such that for each $s \in S^{an}$ the induced map $(\pi^{an})^{-1}(s) \rightarrow \mathrm{pr}^{-1}(\Phi(s))$ is an isomorphism of elliptic curves with level $N$-structure. With this information\footnote{and with the well-known fact that bijective holomorphic maps of complex manifolds already are isomorphisms} one sees that the upper arrow of $(3.4.4)$ is an isomorphism if the lower one is; that the last in turn is true follows easily from the knowledge of the $\C$-valued points of $S$ by considering the functor which it represents.\\
\newline
In view of $(3.4.4)$ we will henceforth always identify the analytification of $(E/S, e_1, e_2)$ with the analytic family of elliptic curves with Drinfeld basis for the $N$-torsion defined by $(3.4.1)$-$(3.4.3)$.
\subsubsection{The universal vectorial extension}
In a next step we describe the analytification of the universal vectorial extension $\widehat{E}^\natural$ of the dual elliptic curve $\widehat{E}$ of $E$. For this we consider the following canonical exact sequence of abelian sheaves on $S^{an}$ (cf. \cite{Maz-Mes}, Ch. I, $(4.4)$):
\[\tag{\textbf{3.4.5}} 0 \rightarrow R^1\pi^{an}_* (2\pi i \Z) \rightarrow H^1_{\mathrm{dR}}(E^{an}/S^{an}) \rightarrow (\widehat{E}^\natural)^{an} \rightarrow 0.\]
We then trivialize the $\mathcal O_{S^{an}}$-vector bundle $H^1_{\mathrm{dR}}(E^{an}/S^{an})$ on the universal covering $\H$ of each connected component of $S^{an}$ by using the cartesian diagram
\begin{equation*} \tag{\textbf{3.4.6}} \begin{split}
\begin{xy}
\xymatrix@C-0.3cm{
\Z^2\backslash(\C\times \H) \ar[r]\ar[d] & (\Z^2\times \Gamma(N)) \backslash (\C \times \H) \ar[d]\\
\H \ar[r] & \Gamma(N) \backslash \H}
\end{xy}
\end{split}
\end{equation*}
and the $\mathcal O_{\H}$-basis $\{p(z,\tau)\mathrm{d}z,\mathrm{d}z\}$ for the de Rham cohomology of the left analytic family of elliptic curves in $(3.4.6)$. This trivialization induces the $2$-dimensional factor of automorphy
\[\tag{\textbf{3.4.7}} \Gamma(N) \times \H \rightarrow \mathrm{GL}_2(\C), \qquad \Bigg(\begin{pmatrix} a & b \\ c & d \end{pmatrix}, \tau \Bigg) \mapsto \begin{pmatrix} \frac{1}{c\tau+d} & 0 \\ 0 & c\tau +d \end{pmatrix}\]
for $H^1_{\mathrm{dR}}(E^{an}/S^{an})$ on each connected component $\Gamma(N) \backslash \H$ of $S^{an}$. The geometric vector bundle associated with $H^1_{\mathrm{dR}}(E^{an}/S^{an})$ is then given by
\begin{equation*} \tag{\textbf{3.4.8}} \begin{split}
\begin{xy}
\xymatrix{
(\Z/N\Z)^*\times \Gamma(N) \backslash(\C^2\times \H) \ar[d]^{\mathrm{can}} \\
(\Z/N\Z)^*\times \Gamma(N) \backslash \H}
\end{xy}
\end{split}
\end{equation*}
with $\Gamma(N)$ acting (properly discontinuously and freely) on $\C^2\times \H$ by the rule
\[\begin{pmatrix} a & b \\ c & d \end{pmatrix} \cdot (w,u,\tau)=\bigg(\frac{w}{c\tau+d}, (c\tau+d) u, \frac{a\tau+b}{c\tau+d} \bigg).\]
From $(3.4.5)$, $(3.4.8)$ and \cite{Kat5}, p. 301, one easily deduces the following description for $(\widehat{E}^\natural)^{an}/S^{an}$:
\begin{equation*} \tag{\textbf{3.4.9}} \begin{split}
\begin{xy}
\xymatrix{
(\Z/N\Z)^*\times (\Z^2\times \Gamma(N))\backslash (\C^2\times \H) \ar[d]^{\mathrm{can}} \\
(\Z/N\Z)^*\times \Gamma(N) \backslash \H}
\end{xy}
\end{split}
\end{equation*}
with $\Z^2\times \Gamma(N)$ acting (properly discontinuously and freely) on $\C^2\times \H$ by the rule
\[\Bigg(\begin{pmatrix} m' \\ n' \end{pmatrix}, \begin{pmatrix} a & b \\ c & d \end{pmatrix} \Bigg) \cdot (w,u,\tau):=\Bigg(\frac{w+m'\tau+n'}{c\tau+d}, (c\tau+d)(u-\eta(m'\tau+n',\tau)), \frac{a\tau +b}{c\tau+d}\Bigg).\]\\
\vspace{1mm}Under the principal polarization $E \simeq \widehat{E}$ in $(3.1.3)$ the analytification of the canonical arrow in $(0.1.3)$
\[\widehat{E}^\natural \rightarrow \widehat{E}\]
becomes the morphism $(\widehat{E}^\natural)^{an}\rightarrow E^{an}$ which is checked to be given as
\[\begin{split}
(\Z/N\Z)^*\times (\Z^2\times \Gamma(N))\backslash (\C^2\times \H) & \rightarrow  (\Z/N\Z)^*\times (\Z^2\times \Gamma(N))\backslash (\C\times \H),\\
(j,w,u,\tau) & \mapsto (j,-w,\tau).\end{split}\]

Taking the product with $\id_E$ we obtain the map $E\times_S \widehat{E}^\natural \rightarrow E \times_S \widehat{E} \simeq E \times_S E$ with analytification
\[\tag{\textbf{3.4.10}}
\begin{split}
(\Z/N\Z)^*\times (\Z^2\times \Z^2 \times \Gamma(N))\backslash (\C \times \C^2\times \H) & \rightarrow (\Z/N\Z)^*\times (\Z^2\times \Z^2 \times \Gamma(N))\backslash (\C \times \C \times \H),\\
(j,z,w,u,\tau) & \mapsto (j,z,-w,\tau),
\end{split}\]
where the action of $\Z^2\times \Z^2 \times \Gamma(N)$ on $\C\times \C^2\times \H$ resp. on $\C \times \C \times \H$ is given by

{\footnotesize
\[\Bigg(\begin{pmatrix}  m \\ n \end{pmatrix}, \begin{pmatrix} m' \\ n' \end{pmatrix}, \begin{pmatrix} a & b \\ c & d \end{pmatrix} \Bigg) \cdot (z, w,u,\tau):=\Bigg(\frac{z+m\tau+n}{c\tau+d}, \frac{w+m'\tau+n'}{c\tau+d}, (c\tau+d)(u-\eta(m'\tau+n',\tau)), \frac{a\tau +b}{c\tau+d}\Bigg)\]}

resp.
{\footnotesize
\[\Bigg( \begin{pmatrix}  m \\ n \end{pmatrix}, \begin{pmatrix} m' \\ n' \end{pmatrix}, \begin{pmatrix} a & b \\ c & d \end{pmatrix} \Bigg) \cdot (z,w,\tau):=\Bigg(\frac{z+m\tau+n}{c\tau+d}, \frac{w+m'\tau+n'}{c\tau+d}, \frac{a\tau +b}{c\tau+d}\Bigg).\]}

With $(3.4.10)$ we have a good analytic access to the objects over which the Poincaré bundle is defined. The next task is to find a suitable description for the analytification of this bundle and of its universal integrable connection. In view of its buildup (cf. $(3.1.5)$) it is clear that we first treat:

\subsubsection{The line bundle defined by the zero section}
Consider again the analytification of the universal elliptic curve with level $N$-structure (cf. $(3.4.4)$):
\begin{equation*} \begin{split}
\begin{xy}
\xymatrix{
(\Z/N\Z)^*\times (\Z^2\times \Gamma(N)) \backslash(\C\times \H) \ar[d]_{\pi^{an}} \\
(\Z/N\Z)^*\times \Gamma(N) \backslash \H \ar@/_ 0.3cm/[u]_{\epsilon^{an}}}
\end{xy}
\end{split}
\end{equation*}
where $\pi^{an}((j,z,\tau))=(j,\tau)$ and $\epsilon^{an}((j,\tau))=(j,0,\tau)$.\\
\newline
Fix $j\in (\Z/N\Z)^*$ and let $E^{an}_j$ resp. $S^{an}_j$ be the connected component of $E^{an}$ resp. $S^{an}$ belonging to $j$. The structure morphism resp. zero section for the analytic family of elliptic curves $E^{an}_j/S^{an}_j$ is given by $\pi^{an}_j((z,\tau))=\tau$ resp. $\epsilon^{an}_j(\tau)=(0,\tau)$. The line bundle $\mathcal O_{E^{an}}([0])$ on $E^{an}$ associated to the analytic effective relative Cartier divisor defined by the zero section of $E^{an}/S^{an}$ becomes the analogous line bundle $\mathcal O_{E^{an}_j}([0])$ for $E^{an}_j/S^{an}_j$ when restricted to the component $E^{an}_j$. Trivializing $\mathcal O_{E^{an}_j}([0])$ on the universal covering $\C\times \H$ of $E^{an}_j$ amounts to giving a meromorphic function on $\C\times \H$ of divisor $-\{(m\tau+n,\tau)|\tau \in \mathbb H, m,n \in \mathbb Z \}$. Such a function is given by $(z,\tau)\mapsto \frac{1}{\vartheta(z,\tau)}$ with $\vartheta(z,\tau)$ the elementary theta function (cf. 3.3.3). Fix this trivialization.\\
The associated factor of automorphy is then read off from Prop. 3.3.10 as
\[\tag{\textbf{3.4.11}} 
\begin{split} &\Z^2\times \Gamma(N) \times \C\times \H \rightarrow \C^*, \quad \Bigg( \begin{pmatrix}  m \\ n \end{pmatrix}, \begin{pmatrix} a & b \\ c & d \end{pmatrix}, (z,\tau) \Bigg)  \\ 
&\mapsto \frac{1}{c\tau +d}\cdot \exp\bigg[\frac{\pi i c (z+m\tau+n)^2}{c\tau+d}+\pi i m+\pi i n - 2\pi i mz-\pi i m^2\tau\bigg].\end{split}\]

\begin{remark}
The line bundle $\mathcal O_{E^{an}_j}(-[0])$ is trivialized on $\C \times \H$ by the function $\vartheta(z,\tau)$. With 3.2 (v) we obtain a trivialization on $\H$ of
\[(\epsilon^{an}_j)^*\mathcal O_{E^{an}_j}(-[0])\simeq (\epsilon^{an}_j)^*\Omega^1_{E^{an}_j/S^{an}_j}\simeq (\pi^{an}_j)_*\Omega^1_{E^{an}_j/S^{an}_j},\] i.e. of the co-Lie algebra $\omega_{E^{an}_j/S^{an}_j}$ of $E^{an}_j/S^{an}_j$.\\
On the other hand, there is a natural trivialization of $\omega_{E^{an}_j/S^{an}_j}$ on $\H$, namely the one induced by $(3.4.6)$ and the standard differential form $\mathrm{d}z$ on $\Z^2\backslash (\C\times \H)$; we remark that this trivialization (or rather: this trivialization multiplied by $2\pi i$) allows the interpretation of modular forms as sections in the tensor powers of the co-Lie algebra of $E^{an}/S^{an}$ (for details cf. \cite{Ka}, Ch. I, 3.8).\\
The two described trivializations of $\omega_{E^{an}_j/S^{an}_j}$ on $\H$, i.e. the one coming from $\vartheta(z,\tau)$ and the one induced by $\mathrm{d}z$, indeed coincide: to see this requires chasing through many identifications, but in the end boils down precisely to the fact that $\partial_z\vartheta(0,\tau)\mathrm{d}z=\mathrm{d}z$ (cf. $(3.3.16)$).
\end{remark}

\subsubsection{The Poincaré bundle I}

Recall from $(3.1.5)$ that under the fixed principal polarization $(3.1.3)$ the birigidified Poincaré bundle $(\mathcal P^0, r^0, s^0)$ on $E\times_S E$ is given by
\[(\mathcal M \otimes_{\mathcal O_{E \times_S E}} (\pi \times \pi)^*\epsilon^* \mathcal O_E([0]), \mathrm{can}, \mathrm{can}),\]
where
\[\mathcal M= \mu^* \mathcal O_E([0]) \otimes_{\mathcal O_{E \times_S E}} \mathrm{pr}_1^* \mathcal O_E([0])^{-1} \otimes_{\mathcal O_{E \times_S E}} \mathrm{pr}_2^*\mathcal O_E([0])^{-1}\]
and $\mathrm{can}$ is the canonical rigidification along the second resp. first factor.\\
\newline
By 3.2 (iii)-(v) we get from the above fixed trivialization of $\mathcal O_{E^{an}_j}([0])$ on $\C\times \H$ a trivialization of $(\mathcal P^0)^{an}_j$ on $\C\times \C \times \H$ whose associated factor of automorphy is computed from $(3.4.11)$ as
\[\tag{\textbf{3.4.12}}
\begin{split} & \Z^2 \times \Z^2\times \Gamma(N) \times \C \times \C \times \H \rightarrow \C^*, \quad \Bigg( \begin{pmatrix}  m \\ n \end{pmatrix}, \begin{pmatrix} m' \\ n' \end{pmatrix}, \begin{pmatrix} a & b \\ c & d \end{pmatrix}, (z,w,\tau) \Bigg)\\
& \mapsto \exp\bigg[\frac{2\pi i c}{c\tau+d}\cdot(z+m\tau+n)(w+m'\tau+n')-2\pi i m'z-2\pi i mw-2\pi i mm'\tau\bigg].\end{split}  \]
Explicitly, the trivialization of
\[(\mathcal P^0)^{an}_j \simeq \mathcal O_{E^{an}_j\times_{S^{an}_j}E^{an}_j}(\bar{\Delta}_{E^{an}_j}-[0]\times E^{an}_j - E^{an}_j\times [0]) \otimes_{\mathcal O_{E^{an}_j\times_{S^{an}_j}E^{an}_j}} (\pi^{an}_j \times \pi^{an}_j)^* (\omega_{E^{an}_j/S^{an}_j})^\vee\]on $\C\times \C \times \H$ is given by the section
\[\tag{\textbf{3.4.13}} t^0:=\frac{1}{J(z,w,\tau)} \otimes (\omega^0_{\mathrm{can}})^\vee,\]
with $J(z,w,\tau)$ the fundamental meromorphic Jacobi form (cf. 3.3.3) and $(\omega_{\mathrm{can}}^0)^\vee$ defined as follows:
\begin{definition}
We write $(\omega^0_{\mathrm{can}})^\vee$ for the trivializing section of $(\pi^{an}_j \times \pi^{an}_j)^*(\omega_{E^{an}_j/S^{an}_j})^\vee$ on $\C \times \C \times \H$ induced via
\begin{equation*}
\begin{xy}
\xymatrix{
\C\times \C \times \H \ar[d] \ar[rr]^{\qquad proj} & & \H \ar[d]\\
E^{an}_j\times_{S^{an}_j}E^{an}_j \ar[rr]^{\qquad \ \pi^{an}_j\times \pi^{an}_j} & & S^{an}_j
}
\end{xy}
\end{equation*}
from the trivialization of $(\omega_{E^{an}_j/S^{an}_j})^\vee$ on $\H$ given by the dual of the standard differential form - cf. Rem. 3.4.1.
\end{definition}

\subsubsection{The Poincaré bundle II}
Denoting as usual by $(\mathcal P,r,s)$ the birigidified Poincaré bundle on $E\times_S \widehat{E}^\natural$ we obtain an induced trivialization of $\mathcal P^{an}_j$ on $\C\times \C^2\times \H$ whose factor of automorphy is given in view of $(3.4.10)$ by
\[\tag{\textbf{3.4.14}} 
\begin{split} & \Z^2 \times \Z^2 \times \Gamma(N) \times \C \times \C^2 \times \H \rightarrow \C^*, \quad \Bigg( \begin{pmatrix}  m \\ n \end{pmatrix}, \begin{pmatrix} m' \\ n' \end{pmatrix}, \begin{pmatrix} a & b \\ c & d \end{pmatrix}, (z,w,u,\tau) \Bigg) \\
& \mapsto \exp\bigg[-\frac{2\pi i c}{c\tau+d}\cdot(z+m\tau+n)(w+m'\tau+n')+2\pi i m'z+2\pi i mw+2\pi i mm'\tau\bigg].\end{split}\]
Write
\[\mathcal P^{an}_j \simeq \mathcal O_{E^{an}_j\times_{S^{an}_j}(\widehat{E}^{\natural})^{an}_j}(K) \otimes_{\mathcal O_{E^{an}_j\times_{S^{an}_j}(\widehat{E}^{\natural})^{an}_j}} (\pi^{an}_j \times (\pi^{\natural})^{an}_j)^* (\omega_{E^{an}_j/S^{an}_j})^\vee\]
with $K$ the pullback divisor of $\bar{\Delta}_{E^{an}_j}-[0]\times E^{an}_j - E^{an}_j\times [0]$ via $E^{an}_j \times_{S^{an}_j} (\widehat{E}^\natural)^{an}_j \rightarrow E^{an}_j\times_{S^{an}_j}E^{an}_j$.\\
The trivialization of $\mathcal P^{an}_j$ on $\C\times \C^2\times \H$ then is given by the section
\[\tag{\textbf{3.4.15}} t:=\frac{1}{J(z,-w,\tau)} \otimes \omega^\vee_{\mathrm{can}},\]
where $(z,w,u,\tau) \mapsto J(z,-w,\tau)$ is considered as meromorphic function on $\C\times \C^2\times \H$ and:
\begin{definition}
We write $\omega^{\vee}_{\mathrm{can}}$ for the trivializing section of $(\pi^{an}_j \times (\pi^{\natural})^{an}_j)^*(\omega_{E^{an}_j/S^{an}_j})^\vee$ on $\C \times \C^2 \times \H$ induced via
\begin{equation*}
\begin{xy}
\xymatrix{
\C\times \C^2 \times \H \ar[d] \ar[rr]^{\qquad proj} & & \H \ar[d]\\
E^{an}_j\times_{S^{an}_j}(\widehat{E}^{\natural})^{an}_j \ar[rr]^{\qquad \quad \pi^{an}_j\times (\pi^{\natural})^{an}_j} & & S^{an}_j
}
\end{xy}
\end{equation*}
from the trivialization of $(\omega_{E^{an}_j/S^{an}_j})^\vee$ on $\H$ given by the dual of the standard differential form - cf. Rem. 3.4.1.
\end{definition}

\begin{remark}
Let us finally record how the birigidification $(r,s)$ of $\mathcal P$ expresses now analytically. An analogous comment applies to $(\mathcal P^0,r^0,s^0)$.\\
Assume we have a section of $\mathcal P^{an}_j$ defined on some open subset of $E^{an}_j\times_{S^{an}_j}(\widehat{E}^{\natural})^{an}_j$. It is given by a holomorphic function $f(z,w,u,\tau)$ which is defined on the inverse image of this subset under $\C\times \C^2\times \H \rightarrow E^{an}_j\times_{S^{an}_j}(\widehat{E}^{\natural})^{an}_j$ and transforms under the effect of $\Z^2\times \Z^2\times \Gamma(N)$ in the way
\begin{align*}
& f\Bigg(\frac{z+m\tau+n}{c\tau+d}, \frac{w+m'\tau+n'}{c\tau+d}, (c\tau+d)(u-\eta(m'\tau+n',\tau)), \frac{a\tau +b}{c\tau+d}\Bigg)\\
& =\exp\bigg[-\frac{2\pi i c}{c\tau+d}\cdot(z+m\tau+n)(w+m'\tau+n')+2\pi i m'z+2\pi i mw+2\pi i mm'\tau\bigg]\cdot f(z,w,u,\tau).
\end{align*}

The function $(w,u,\tau)\mapsto f(0,w,u,\tau)$ resp. $(z,\tau)\mapsto f(z,0,0,\tau)$ (defined on the appropriate open subset of $\C^2\times \H$ resp. $\C\times \H$) is then invariant under the effect of $\Z^2\times \Gamma(N)$ on $\C^2\times \H$ resp. $\C\times \H$ and represents the section of $\mathcal O_{(\widehat{E}^\natural)^{an}_j}$ resp. $\mathcal O_{E^{an}_j}$ given by pullback of the given section via $\epsilon^{an}_j\times \mathrm{id}_{(\widehat{E}^{\natural})^{an}_j}$ resp. $\mathrm{id}_{E^{an}_j}\times (\epsilon^{\natural})^{an}_j$ and using $r^{an}_j$ resp. $s^{an}_j$.\\
With the definition of $t^0$ (cf. $(3.4.13)$) the verification of these facts is easy if one first observes that $(r,s)$ is induced by $(r^0,s^0)$ (cf. 0.1.1), recalls how $(r^0,s^0)$ expresses for elliptic curves (cf. $(3.1.5)$) and then takes into account Rem. 3.4.1.
\end{remark}

\subsubsection{The universal integrable connection}
Recall that the Poincaré bundle $\mathcal P$ on $E\times_S \widehat{E}^\natural$ is equipped with the universal integrable connection
\[\nabla_\mathcal P: \mathcal P \rightarrow \Omega^1_{E\times_S \widehat{E}^\natural/\widehat{E}^\natural}\otimes_{\mathcal O_{E\times_S \widehat{E}^\natural}} \mathcal P.\]
We now want to give a formula for the analytification\footnote{The formally clean way to analytify connections is to interpret them equivalently in terms of $\mathcal O$-linear structures which one can then analytify as usual (cf. e.g. the viewpoint on connections outlined in 0.2.1 (v) or the approach via differential operators in \cite{Mal}, p. 152, whose formalism of course generalizes to a smooth relative situation of $\Q$-schemes).} of this connection:
\[\nabla_\mathcal P^{an}: \mathcal P^{an} \rightarrow \Omega^1_{E^{an}\times_{S^{an}} (\widehat{E}^\natural)^{an}/(\widehat{E}^\natural)^{an}}\otimes_{\mathcal O_{E^{an}\times_{S^{an}} (\widehat{E}^\natural)^{an}}} \mathcal P^{an},\]
which we do again by restricting to a fixed $j\in (\Z/N\Z)^*$:
\[(\nabla_\mathcal P^{an})_j: \mathcal P^{an}_j \rightarrow \Omega^1_{E^{an}_j\times_{S^{an}_j} (\widehat{E}^\natural)^{an}_j/(\widehat{E}^\natural)^{an}_j}\otimes_{\mathcal O_{E^{an}_j\times_{S^{an}_j} (\widehat{E}^\natural)^{an}_j}} \mathcal P^{an}_j.\]\\
\vspace{1mm}We trivialize the pullback of $\mathcal P^{an}_j$ resp. of $\Omega^1_{E^{an}_j\times_{S^{an}_j} (\widehat{E}^\natural)^{an}_j/(\widehat{E}^\natural)^{an}_j}$ to $\C\times \C^2 \times \H$ as in $(3.4.15)$ resp. by the standard differential form $\{ \mathrm{d}z \}$\footnote{Observe that the pullback of $\Omega^1_{E^{an}_j\times_{S^{an}_j} (\widehat{E}^\natural)^{an}_j/(\widehat{E}^\natural)^{an}_j}$ to $\C\times \C^2\times \H$ naturally identifies with $\Omega^1_{\C\times \C^2\times \H/\C^2\times \H}$.}, with associated factor of automorphy $(3.4.14)$ resp.
\[\Z^2 \times \Z^2 \times \Gamma(N) \times \C \times \C^2 \times \H \rightarrow \C^*, \quad \Bigg( \begin{pmatrix}  m \\ n \end{pmatrix}, \begin{pmatrix} m' \\ n' \end{pmatrix}, \begin{pmatrix} a & b \\ c & d \end{pmatrix}, (z,w,u,\tau) \Bigg) \mapsto\\ c\tau+d.\]

Assume then that a section of $\mathcal P^{an}_j$ on some open subset $U$ of $E^{an}_j\times_{S^{an}_j}(\widehat{E}^{\natural})^{an}_j$ is given: as usual, it is described by a holomorphic function $f(z,w,u,\tau)$, defined on the inverse image of $U$ under the canonical map $\mathrm{pr}:\C\times \C^2\times \H \rightarrow E^{an}_j\times_{S^{an}_j}(\widehat{E}^{\natural})^{an}_j$
and transforming under the effect of $\Z^2\times \Z^2\times \Gamma(N)$ with the factor of automorphy for $\mathcal P^{an}_j$ in $(3.4.14)$.\\
The section $(\nabla_\mathcal P^{an})_j(f(z,w,u,\tau))$ of $\Omega^1_{E^{an}_j\times_{S^{an}_j} (\widehat{E}^\natural)^{an}_j/(\widehat{E}^\natural)^{an}_j}\otimes_{\mathcal O_{E^{an}_j\times_{S^{an}_j} (\widehat{E}^\natural)^{an}_j}} \mathcal P^{an}_j$ is given by a holomorphic function on $\mathrm{pr}^{-1}(U)$ transforming with $(c\tau+d)$-times the factor in $(3.4.14)$.\\
By considering the situation fiberwise over points of $S^{an}_j$ and noting \cite{Kat5}, App. C, Thm. C. 6 (1), one straightforwardly derives the following formula for this function:\footnote{When verifying the computation in fibers the reader should observe the following equality:
\begin{align*}
&\partial_z \frac{f(z,w,u;\tau)}{J(z,-w;\tau)}\mathrm{d}z+ \frac{f(z,w,u;\tau)}{J(z,-w;\tau)} \cdot (\zeta(z-w;\tau)- \zeta(z;\tau)+u)\mathrm{d}z\\
&=\frac{1}{J(z,-w;\tau)}\bigg[\partial_z f(z,w,u;\tau)\mathrm{d}z-f(z,w,u;\tau)\cdot \frac{\partial_zJ(z,-w;\tau)}{J(z,-w;\tau)}\mathrm{d}z+f(z,w,u;\tau)(\zeta(z-w;\tau)-\zeta(z;\tau)+u)\mathrm{d}z \bigg]\\
&=\frac{1}{J(z,-w;\tau)}\bigg[\partial_zf(z,w,u;\tau)\mathrm{d}z -f(z,w,u;\tau)\cdot\big(-\eta(1;\tau)w+\zeta(z-w;\tau)-\zeta(z;\tau)\big)\mathrm{d}z\\
&+f(z,w,u;\tau)\cdot(\zeta(z-w;\tau)-\zeta(z;\tau)+u)\mathrm{d}z\bigg]\\
&=\frac{1}{J(z,-w;\tau)}\bigg[\partial_zf(z,w,u;\tau)+(\eta(1;\tau)w+u)\cdot f(z,w,u;\tau) \bigg]\mathrm{d}z,
\end{align*}
where the equation for the logarithmic $z$-derivative of $J(z,-w;\tau)$, used in the transition from the second to the third line, follows from $(3.3.18)$.}
\[\tag{\textbf{3.4.16}} (\nabla_\mathcal P^{an})_j(f(z,w,u,\tau))=\partial_z f(z,w,u,\tau) +(\eta(1, \tau)w+u)\cdot f(z,w,u,\tau).\]

\begin{remark}
The horizontality of $s^{an}_j: (\mathrm{id}_{E^{an}_j}\times (\epsilon^{\natural})^{an}_j)^*\mathcal P^{an}_j \simeq \mathcal O_{E^{an}_j}$ (cf. Rem. 0.1.17) is reflected by $(3.4.16)$ in the fact that this formula becomes (relative) exterior derivation on the $(\Z^2\times\Gamma(N))$-invariant function $(z,\tau) \mapsto f(z,0,0,\tau)$ after plugging in $w=0=u$ (cf. also Rem. 3.4.4).
\end{remark}

\markright{\uppercase{The explicit description on the universal elliptic curve}}

\section{The analytification of the logarithm sheaves}
\markright{\uppercase{The explicit description on the universal elliptic curve}}

\subsubsection{Preliminary remarks and conventions: generalities}

In the following, we will as usual write $\mathcal H$ for the dual of the $\mathcal O_S$-vector bundle $H^1_{\mathrm{dR}}(E/S)$.\\
The canonical map
\[H^1_{\mathrm{dR}}(E/S)^{an} \rightarrow H^1_{\mathrm{dR}}(E^{an}/S^{an})\]
is an isomorphism\footnote{One first considers the base extension $(E\times_{\Q}\C)/(S\times_{\Q}\C)$ of $E/S$, for which the relative de Rham cohomology commutes (by flatness of $S\times_{\Q}\C\rightarrow S$); then one applies \cite{De1}, Ch. II, Thm. 6.13 and Prop. 6.14, or \cite{Har}, Ch. IV, Prop. (4.1).}, and hence when writing $\mathcal H^{an}$ we will equally refer to the analytification of $\mathcal H$ or to the dual of the first (analytic) de Rham cohomology sheaf of $E^{an}/S^{an}$.\\
\newline
Moreover, we will have to consider the analytification of infinitesimal thickenings. The notational conventions of Def. 2.2.14 will remain fixed.\\
Hence, if $\mathcal J$ denotes again the ideal sheaf of the zero section $\epsilon^{\natural}: S \rightarrow \widehat{E}^\natural$, one can check that $(\epsilon^{\natural})^{an}: S^{an}\rightarrow (\widehat{E}^\natural)^{an}$ is a closed (analytic) immersion with ideal sheaf $\mathcal J^{an}$.\footnote{This is general and deduced from application of the (exact) analytification functor to the canoncial exact sequence defining the ideal sheaf and observing that the occurring direct image commutes with analytification: use  \cite{SGA1}, exp. XII, Thm. 4.2.}\\
One also checks that $(\widehat{E}^\natural_n)^{an}$ is identical with the $n$-th (analytic) infinitesimal neighborhood of $S^{an}$ in $(\widehat{E}^\natural)^{an}$, defined by $(\mathcal J^{an})^{n+1} \subseteq \mathcal O_{(\widehat{E}^\natural)^{an}}$.\\
\newline
Furthermore, we will freely use compatibility of analytification with fiber products, in the sense that the canonical map $(V \times_Z W)^{an} \rightarrow V^{an}\times_{Z^{an}}W^{an}$ is an isomorphism for schemes $V,W,Z$ locally of finite type over $\Spec(\C)$ providing a fiber product situation (cf. \cite{SGA1}, exp. XII, 1.2).\\
\newline
Finally, let us slightly ease notation and fix the following \underline{agreement}:\\
As in the previous section we work on a fixed connected component of $E^{an}/S^{an}$, and for quite a long time it still won't matter which $j \in (\Z/N\Z)^*$ is chosen - as an example, recall that the factor of automorphy for $\mathcal P^{an}_j$ in $(3.4.14)$ or the formula for the connection $(\nabla_\mathcal P)^{an}_j$ in $(3.4.16)$ is the same on each component. We therefore no longer use a subscript practice like $E^{an}_j, S^{an}_j, (\widehat{E}^\natural)^{an}_j, \mathcal P^{an}_j, \mathcal H^{an}_j,...$, but in abuse of notation leave away the $j$.\\
Hence, until explicitly said otherwise, $E^{an}, S^{an}, (\widehat{E}^\natural)^{an}, \mathcal P^{an}, \mathcal H^{an}$ etc. will denote what actually is some fixed component $(\Z^2\times\Gamma(N))\backslash (\C\times \H), \Gamma(N)\backslash \H, (\Z^2\times\Gamma(N))\backslash (\C^2\times \H)$, the restriction of $\mathcal P^{an}$ to $(\Z^2\times \Z^2 \times \Gamma(N))\backslash (\C \times \C^2 \times \H)$, the restriction of $\mathcal H^{an}$ to $\Gamma(N) \backslash \H$ resp. (by the above comment) the dual of the first (analytic) de Rham cohomology sheaf of $(\Z^2 \times \Gamma(N)) \backslash (\C \times \H) \rightarrow \Gamma(N) \backslash \H$ etc.\\
In the end, all of our constructions on the henceforth fixed component which apparently don't depend on the choice of the component are to be understood as carried out simultaneously on all components.

\subsubsection{Preliminary remarks and conventions: trivialization for de Rham cohomology}

As explained in the previous section (cf. $(3.4.6)$-$(3.4.8)$) we trivialize $H^1_{\mathrm{dR}}(E/S)^{an}\simeq \mathrm{H^1_{dR}}(E^{an}/S^{an})$ on the universal covering $\H$ of $S^{an}$ by $\{\eta,\omega\}:=\{p(z,\tau)\mathrm{d}z,\mathrm{d}z\}$ and then $\mathcal H^{an}$ by $\{\eta^\vee, \omega^\vee\}$.\\
Hence, denoting for a moment by $p: \H \rightarrow \Gamma(N)\backslash \H$ the canonical projection:
\[\tag{\textbf{3.5.1}} p^*\mathcal H^{an} \simeq \mathcal O_{\H}\cdot  \eta^\vee \oplus \mathcal O_{\H}\cdot \omega^\vee \simeq \mathcal O_{\H} \oplus \mathcal O_{\H}.\]
We obtain the following automorphy matrix for $\mathcal H^{an}$:
\[\tag{\textbf{3.5.2}} \Gamma(N) \times \H \rightarrow \mathrm{GL}_2(\C), \qquad \Bigg(\begin{pmatrix} a & b \\ c & d \end{pmatrix}, \tau \Bigg) \mapsto \begin{pmatrix} c\tau+d & 0 \\ 0 & \frac{1}{c\tau +d} \end{pmatrix}.\]
For each $n \geq 0$ we then trivialize the pullback of $\prod_{k=0}^n\mathrm{Sym}^k_{\mathcal O_{S^{an}}}\mathcal H^{an}$ to $\H$ by
\[p^*\prod_{k=0}^n\mathrm{Sym}^k_{\mathcal O_{S^{an}}}\mathcal H^{an} \simeq \bigoplus_{\substack{0\leq i+j \leq n\\
0 \leq i,j}} \mathcal O_{\H}\cdot \frac{(\eta^\vee)^i (\omega^\vee)^j}{(n-i-j)!} \simeq \bigoplus_{\substack{0\leq i+j \leq n \\ 0\leq i,j}} \mathcal O_{\H}.\]
As it will play a role for writing down automorphy matrices we fix - for now and for the future - the order of the basis elements in the preceding trivialization;  we do this by letting increase the number of $\omega^\vee$-factors and within this ordering rule letting increase the number of $\eta^\vee$-factors:

{\footnotesize
\[\bigg \{ \frac{1}{n!}, \frac{\eta^\vee}{(n-1)!}, \frac{(\eta^\vee)^2}{(n-2)!},..., \frac{(\eta^\vee)^n}{(n-n)!}, \frac{\omega^\vee}{(n-1)!}, \frac{\eta^\vee\omega^\vee}{(n-2)!},\frac{(\eta^\vee)^2\omega^\vee}{(n-3)!},..., \frac{(\eta^\vee)^{n-1}\omega^\vee}{(n-n)!}, \frac{(\omega^\vee)^2}{(n-2)!},.........,\frac{(\omega^\vee)^n}{(n-n)!} \bigg \}.\]}

If $r(n):=\mathrm{rk}_{\mathcal O_{S^{an}}}\bigg(\prod_{k=0}^n \mathrm{Sym}^k_{\mathcal O_{S^{an}}}\mathcal H^{an}\bigg)$ we obtain an automorphy matrix for $\prod_{k=0}^n \mathrm{Sym}^k_{\mathcal O_{S^{an}}}\mathcal H^{an}$ of size $r(n)\times r(n)$ whose entries lie only on the diagonal and are computed by the transformation property of the preceding basis elements:
\[\Gamma(N) \times \H \rightarrow \mathrm{GL}_{r(n)}(\C), \qquad \Bigg(\begin{pmatrix} a & b \\ c & d \end{pmatrix}, \tau \Bigg) \mapsto
\begin{pmatrix}
1	& 	&  & & \\
	& c\tau +d	&   & & \\
	&  	& (c\tau+d)^2 & \\
	&   &     & \ddots & \\
	&   &     &        & (c\tau +d)^{-n}
\end{pmatrix}
 \]
In general, the entry coming from $\frac{(\eta^\vee)^i (\omega^\vee)^j}{(n-i-j)!}$ is given by $(c\tau+d)^{i-j}$.

\subsection{The analytification of the first logarithm sheaf}
We consider the analytification of the first logarithm sheaf $\mathcal L_1$ of $E/S/\Q$ and (componentwise) fix a trivialization for its pullback to the universal covering $\C\times \H$. We then compute the associated automorphy matrix and, expressing in this way sections via holomorphic functions on $\C\times \H$, deduce a formula for the analytified integrable $S$-connection belonging to $\mathcal L_1$.\\
Our definition of the trivialization and the determination of the relative connection is based substantially on the construction of $\mathcal L_1$ in terms of the Poincaré bundle on $E\times_S \widehat{E}^\natural$ and its universal integrable $\widehat{E}^\natural$-connection (cf. Cor. 2.3.2).\\
We remark that the splitting $\varphi_1$ will be described analytically in 3.5.2, and that the problem of finding a formula for the analytification of the absolute connection $\nabla_1$ of $\mathcal L_1$ will be addressed in 3.5.4.

\subsubsection{A commutative diagram of ringed spaces}
Consider the closed analytic immersion
\[\mathrm{id}_{E^{an}}\times (\epsilon^{\natural}_1)^{an}: E^{an}\times_{S^{an}}(\widehat{E}^\natural_1)^{an} \rightarrow E^{an}\times_{S^{an}}(\widehat{E}^\natural)^{an}.\]
As a ringed space $ E^{an}\times_{S^{an}}(\widehat{E}^\natural_1)^{an}$ is the topological space $E^{an}$, endowed with a sheaf of rings which (by analytifying $(2.3.5)$-$(2.3.6)$) identifies with
\[\mathcal O_{E^{an}} \oplus \mathcal H_{E^{an}}^{an},\]
where multiplication is given by $(s_1, h_1) \cdot (s_2 , h_2)=(s_1s_2, s_1 h_2+s_2h_1)$; here, $\mathcal H_{E^{an}}^{an}$ means the pullback of $\mathcal H^{an}$ via $\pi^{an}: E^{an}\rightarrow S^{an}$; we could equally well have written $(\mathcal H_E)^{an}$.\\
\newline
We define a diagram of ringed spaces

\begin{equation*} \tag{\textbf{3.5.3}} 
\begin{xy}
\xymatrix@C-0.7cm{
(\C \times \H, \mathcal O_{\C\times \H}\oplus \mathcal O_{\C\times \H}\cdot \eta^\vee\oplus \mathcal O_{\C\times \H}\cdot \omega^\vee) \ar[r]^{\qquad  \qquad k_1}\ar[d]^{\tau_1} & (\C \times \C^2 \times \H,\mathcal O_{\C \times \C^2 \times \H}) \ar[d]_{\mathrm{can}}\\
((\Z^2 \times \Gamma(N))\backslash(\C\times \H), \mathcal O_{E^{an}} \oplus \mathcal H_{E^{an}}^{an}) \ar[r] & ((\Z^2 \times \Z^2 \times \Gamma(N))\backslash (\C \times \C^2 \times \H),\mathcal O_{E^{an}\times_{S^{an}}(\widehat{E}^\natural)^{an}})}
\end{xy}
\end{equation*}
as follows:\\
The right part consists of $E^{an}\times_{S^{an}}(\widehat{E}^\natural)^{an}$ and the projection from its universal covering. The lower part was explained above. The upper left space consists of $\C\times \H$ and the sheaf of rings $\mathcal O_{\C\times \H}\oplus\mathcal O_{\C\times \H}\cdot \eta^\vee\oplus \mathcal O_{\C\times \H}\cdot \omega^\vee$ whose addition works componentwise and whose multiplication is defined by linearity, multiplication in $\mathcal O_{\C\times \H}$ and letting $(\eta^\vee)^2=(\omega^\vee)^2=\eta^\vee \omega^\vee=\omega^\vee \eta^\vee=0$. The left vertical arrow $\tau_1$ topologically is the projection and on sheaf level is defined by pulling back sections of the $\mathcal O_{E^{an}}$-vector bundle $\mathcal O_{E^{an}} \oplus \mathcal H_{E^{an}}^{an}$ along the usual projection
\[\mathrm{pr}: (\C\times \H, \mathcal O_{\C\times \H}) \rightarrow E^{an}\]
and by then using the trivialization induced by $(3.5.1)$:
\[\mathrm{pr}^*(\mathcal O_{E^{an}} \oplus \mathcal H_{E^{an}}^{an}) \simeq \mathcal O_{\C\times \H} \oplus \mathcal O_{\C\times \H} \cdot \eta^\vee \oplus \mathcal O_{\C\times \H}\cdot \omega^\vee.\]
Finally, the upper horizontal arrow $k_1$ topologically is given by $(z,\tau)\mapsto (z,0,0,\tau)$ and on sheaves maps a section $g(z,w,u,\tau)$ to
\[g_{0,0}(z,\tau) + g_{1,0}(z,\tau) \cdot \eta^\vee + g_{0,1}(z,\tau) \cdot \omega^\vee,\]
where $g_{0,0}(z,\tau):=g(z,0,0,\tau)$, $g_{1,0}(z,\tau):=\partial_wg(z,0,0,\tau), \ g_{0,1}(z,\tau):=\partial_ug(z,0,0,\tau)$, i.e.
the $g_{i,j}(z,\tau)$ are defined by the Taylor expansion of $g(z,w,u,\tau)$ around $(w,u)=(0,0)$:
\[g(z,w,u,\tau)=g_{0,0}(z,\tau)+g_{1,0}(z,\tau) \cdot w + g_{0,1}(z,\tau)\cdot u + \ \textrm{higher terms}.\]
It is a straightforward exercise to check that the diagram of ringed spaces $(3.5.3)$ commutes.
\subsubsection{Trivializing the first logarithm sheaf on the universal covering}
Recall from $(3.4.15)$ that the pullback of $\mathcal P^{an}$ to the universal covering $\C\times \C^2 \times \H$ is trivialized by the section
\[t:=\frac{1}{J(z,-w,\tau)} \otimes \omega_{\mathrm{can}}^{\vee}.\]
Observing the above commutative diagram we get by further pullback via $k_1$ an induced isomorphism
\[\tau_1^*\mathcal P_1^{an} \simeq \mathcal O_{\C\times \H}\oplus \mathcal O_{\C\times \H}\cdot \eta^\vee\oplus \mathcal O_{\C\times \H}\cdot \omega^\vee.\]
As $\mathcal O_{\C\times \H}$-modules, the right side is $\mathcal O_{\C\times \H}^{\oplus 3}$ and the left side is $\mathrm{pr}^* (p^{an}_1)_*\mathcal P_1^{an}$, where
\[\mathrm{pr}: (\C\times \H, \mathcal O_{\C\times \H}) \rightarrow E^{an},\]
\[p_1^{an}: E^{an}\times_{S^{an}} (\widehat{E}^\natural_1)^{an} \rightarrow E^{an}.\]
We obtain an isomorphism
\[\mathrm{pr}^*(p_1^{an})_*\mathcal P_1^{an} \simeq \mathcal O_{\C \times \H}^{\oplus 3}.\]
Because $(p_1)_*\mathcal P_1$ is the first logarithm sheaf $\mathcal L_1$ of $E/S/\Q$ (cf. Cor. 2.3.2) we arrive at a trivialization
\[\mathrm{pr}^*\mathcal L_1^{an} \simeq \mathcal O_{\C \times \H}^{\oplus 3}\]
for the pullback of the $\mathcal O_{E^{an}}$-vector bundle $\mathcal L_1^{an}$ to the universal covering $\C\times \H$ of $E^{an}$.
\begin{definition}
Let
\[\mathrm{pr}: (\C\times \H, \mathcal O_{\C\times \H}) \rightarrow E^{an},\]
\[p_1^{an}: E^{an}\times_{S^{an}} (\widehat{E}^\natural_1)^{an} \rightarrow E^{an}\]
be the projections, let $\widetilde{\mathcal P}$ be the pullback of $\mathcal P^{an}$ to $\C \times \C^2 \times \H$, and let
\[\tag{\textbf{3.5.4}} \mathcal O_{\C \times \C^2 \times \H} \xrightarrow{\sim} \widetilde{\mathcal P} \quad 1\mapsto t:=\frac{1}{J(z,-w,\tau)}\otimes \omega_{\mathrm{can}}^\vee\]
be the trivialization of $(3.4.15)$.\\
We then define and henceforth fix a trivialization
\[\tag{\textbf{3.5.5}} \mathcal O_{\C \times \H}^{\oplus 3} \simeq \mathrm{pr}^*\mathcal L_1^{an}\]
by combining the $\mathcal O_{\C\times \H}$-linear isomorphisms
\[\tag{\textbf{3.5.6}} k_1^*\widetilde{\mathcal P}\simeq \tau_1^*\mathcal P_1^{an}\simeq \mathrm{pr}^*(p^{an}_1)_*\mathcal P_1^{an}=\mathrm{pr}^*\mathcal L_1^{an}\]
and
\[\tag{\textbf{3.5.7}}
\begin{split} \mathcal O_{\C\times \H}^{\oplus 3} & \xrightarrow{\sim} k_1^{-1}\widetilde{\mathcal P} \otimes_{k_1^{-1}\mathcal O_{\C\times \C^2 \times \H}}(\mathcal O_{\C\times \H} \oplus \mathcal O_{\C\times \H}\cdot \eta^\vee \oplus \mathcal O_{\C\times \H}\cdot \omega^\vee)=k_1^*\widetilde{\mathcal P},\\
(1,0,0) & \mapsto k_1^{-1}(t) \otimes 1,\\
(0,1,0) & \mapsto k_1^{-1}(t)\otimes \eta^\vee,\\
(0,0,1) & \mapsto k_1^{-1}(t)\otimes \omega^\vee.
\end{split}\]
Here, $\tau_1$ and $k_1$ are the morphisms of ringed spaces appearing in the commutative diagram $(3.5.3)$.
\end{definition}

We can formulate the content of Def. 3.5.1 differently as follows:\\
\newline
Analytify (and restrict to our fixed connected component) the exact sequence $(2.3.11)$:
\[0 \rightarrow \mathcal H^{an}_{E^{an}} \rightarrow (p^{an}_1)_*\mathcal P^{an}_1 \rightarrow \mathcal O_{E^{an}} \rightarrow 0.\]
Its pullback via $\mathrm{pr}: (\C\times \H, \mathcal O_{\C\times \H}) \rightarrow E^{an}$ reads as the $\mathcal O_{\C\times \H}$-linear exact sequence
\[\tag{\textbf{3.5.8}} 0 \rightarrow \mathcal O_{\C\times \H} \cdot \eta^\vee \oplus \mathcal O_{\C\times \H}\cdot \omega^\vee \rightarrow \mathrm{pr}^*(p^{an}_1)_*\mathcal P^{an}_1 \rightarrow \mathcal O_{\C\times \H} \rightarrow 0,\]
where we have used the identification induced by $(3.5.1)$. Via the natural isomorphism of $(3.5.6)$:
\[k_1^{-1}\widetilde{\mathcal P}\otimes_{k_1^{-1}\mathcal O_{\C\times \C^2 \times \H}}(\mathcal O_{\C\times \H} \oplus \mathcal O_{\C\times \H}\cdot \eta^\vee \oplus \mathcal O_{\C\times \H}\cdot \omega^\vee)=k_1^*\widetilde{\mathcal P}\simeq \tau_1^*\mathcal P_1^{an}\simeq \mathrm{pr}^*(p^{an}_1)_*\mathcal P_1^{an}\]
we obtain from the trivializing section $t$ of $\widetilde{\mathcal P}$ a global section
\[\tag{\textbf{3.5.9}} k_1^{-1}(t)\otimes 1 \in \Gamma(\C\times \H,\mathrm{pr}^*(p^{an}_1)_*\mathcal P_1^{an}).\]
It maps to $1$ under the projection in $(3.5.8)$, as is easily verified from the definition of $t^0$ (cf. $(3.4.13)$) by recalling how the projection is induced by the trivialization $s$ (cf. $(2.3.9)$-$(2.3.11)$), how $s$ is induced by $s^0$ (cf. 0.1.1), how $s^0$ expresses for elliptic curves (cf. $(3.1.5)$) and by noting Rem. 3.4.1.\\
In a similar way, by additionally observing how $(3.4.9)$ came about, one sees that the two remaining trivializing sections in $(3.5.7)$:
\begin{align*}
k_1^{-1}(t)\otimes \eta^\vee &\in \Gamma(\C\times \H,\mathrm{pr}^*(p^{an}_1)_*\mathcal P_1^{an})\\
k_1^{-1}(t)\otimes \omega^\vee &\in \Gamma(\C\times \H,\mathrm{pr}^*(p^{an}_1)_*\mathcal P_1^{an})
\end{align*}
are nothing else than the images of the sections $\eta^\vee$ and $\omega^\vee$ under the inclusion of $(3.5.8)$.\\
\newline
Let us summarize:

\begin{proposition}
If
\[\{ k_1^{-1}(t)\otimes 1, \eta^\vee, \omega^\vee \}\]
are the global sections of $\mathrm{pr}^*\mathcal L_1^{an}$ given by $(3.5.9)$ resp. via the inclusion in $(3.5.8)$, then
\[\mathrm{pr}^*\mathcal L_1^{an} = \mathcal O_{\C\times \H}\cdot \big(k_1^{-1}(t)\otimes 1 \big) \oplus \mathcal O_{\C\times \H} \cdot \eta^\vee \oplus \mathcal O_{\C\times \H}\cdot \omega^\vee.\]
Under this decomposition the maps in $(3.5.8)$ are given by the inclusion of the second and third component resp. by the projection to the coefficient function of the first component.\\
The induced trivialization
\[\mathrm{pr}^*\mathcal L_1^{an} \simeq \mathcal O_{\C \times \H}^{\oplus 3}\]
equals $(3.5.5)$, and the associated automorphy matrix for $\mathcal L_1^{an}$ is given by
\begin{align*}
A_1: & \ \Z^2 \times \Gamma(N) \times \C \times \H \rightarrow \mathrm{GL}_3(\C),\\
&\Bigg( \begin{pmatrix}  m \\ n \end{pmatrix}, \begin{pmatrix} a & b \\ c & d \end{pmatrix}, (z,\tau) \Bigg) \mapsto
\begin{pmatrix} 1 & 0 & 0\\
                        2\pi i (dm-cz-cn) & c\tau+d & 0\\
												0 & 0 & \frac{1}{c\tau+d}
\end{pmatrix}.
\end{align*}
\end{proposition}
\begin{proof}
It only remains to compute the automorphy matrix. For this let
\[\Bigg(\begin{pmatrix}  m \\ n \end{pmatrix}, \begin{pmatrix} a & b \\ c & d \end{pmatrix}\Bigg) \in \Z^2\times \Gamma(N),\]
and let
\[\sigma: \C\times \H \rightarrow \C\times \H, \quad (z,\tau) \mapsto \bigg(\frac{z+m\tau+n}{c\tau+d},\frac{a\tau+b}{c\tau+d}\bigg)\]
be the corresponding deck transformation of $\C\times \H$ over $E^{an}$.\\
Now observe that
\begin{align*}
\sigma^*(k_1^{-1}(t) \otimes 1)& =k_1^{-1}\bigg(\exp\bigg[\frac{2\pi i}{c\tau+d}\cdot (cz+cn-dm)\cdot w\bigg]\cdot t\bigg)\otimes 1\\
& =k_1^{-1}(t)\otimes \bigg(1\oplus \frac{2\pi i }{c\tau+d}(cz+cn-dm)\cdot \eta^\vee \bigg).
\end{align*}
Here, the first equality is deduced from the fact that the section $t=\frac{1}{J(z,-w,\tau)} \otimes \omega_{\mathrm{can}}^\vee$ of $\widetilde{\mathcal P}$ transforms under pullback along the deck transformation of $\C\times \C^2\times \H$ over $E^{an}\times_{S^{an}}(\widehat{E}^\natural)^{an}$ belonging to
\[\Bigg( \begin{pmatrix}  m \\ n \end{pmatrix}, \begin{pmatrix} 0 \\ 0 \end{pmatrix}, \begin{pmatrix} a & b \\ c & d \end{pmatrix} \Bigg) \in \Z^2\times \Z^2 \times \Gamma(N)\]
with the factor
\[\exp\bigg[\frac{2\pi ic}{c\tau+d}\cdot (z+m\tau+n)\cdot w-2\pi i m w\bigg]=\exp\bigg[\frac{2\pi i}{c\tau+d}\cdot (cz+cn-dm)\cdot w\bigg],\]
as one recollects by observing (the inverse of) the factor of automorphy in $(3.4.14)$. For the second equality one recalls how $\mathcal O_{\C\times \H} \oplus \mathcal O_{\C\times \H}\cdot \eta^\vee \oplus \mathcal O_{\C\times \H}\cdot \omega^\vee$ becomes a sheaf of rings over $k_1^{-1}\mathcal O_{\C\times \C^2 \times \H}$, i.e. the definition of $k_1$ as a morphism of ringed spaces (see the diagram $(3.5.3)$).\\
We have thus shown that
\[\sigma^*(k_1^{-1}(t) \otimes 1)=k_1^{-1}(t) \otimes 1+\frac{2\pi i}{c\tau+d} (cz+cn-dm)\cdot \big(k_1^{-1}(t) \otimes \eta^\vee\big),\]
and the section $ k_1^{-1}(t) \otimes \eta^\vee$ was already noticed to be the image of $\eta^\vee$ under $(3.5.8)$.\\
Moreover,
\begin{align*}
\sigma^*(\eta^\vee) &=\frac{1}{c\tau+d}\cdot \eta^\vee,\\
\sigma^*(\omega^\vee) &=(c\tau+d) \cdot \omega^\vee,
\end{align*}
such that altogether the matrix describing the change of basis is given by
\[\begin{pmatrix} 1 & 0 & 0\\
                        \frac{2\pi i}{c\tau+d} (cz+cn-dm) & \frac{1}{c\tau+d} & 0\\
												0 & 0 & c\tau+d
\end{pmatrix}.\]
By definition (cf. 3.2 (i)), the desired automorphy matrix evaluated at {\footnotesize$\Bigg( \begin{pmatrix}  m \\ n \end{pmatrix}, \begin{pmatrix} a & b \\ c & d \end{pmatrix}, (z,\tau) \Bigg) $} is the inverse of the previous matrix. But this inverse is precisely the matrix in the claim.
\end{proof}

\subsubsection{Description of the relative connection}
The first logarithm sheaf $\mathcal L_1=(p_1)_*\mathcal P_1$ of $E/S/\Q$ is equipped with an integrable $\Q$-connection
\[\nabla_1: \mathcal L_1\rightarrow \Omega^1_{E/\Q}\otimes_{\mathcal O_E} \mathcal L_1\]
whose restriction to an integrable $S$-connection
\[\nabla_1^{res}: \mathcal L_1\rightarrow \Omega^1_{E/S}\otimes_{\mathcal O_E} \mathcal L_1\]
is induced by the universal integrable $\widehat{E}^\natural$-connection $\nabla_\mathcal P$ of $\mathcal P$ in the way explained in $(2.3.10)$-$(2.3.11)$ (recall again Cor. 2.3.2). We want to describe explicitly the analytification (and restriction to our fixed connected component) of this relative connection:
\[(\nabla_1^{res})^{an}: \mathcal L_1^{an}\rightarrow \Omega^1_{E^{an}/S^{an}}\otimes_{\mathcal O_{E^{an}}} \mathcal L_1^{an}.\]
Of course, the essential ingredient to derive the desired formula for $(\nabla_1^{res})^{an}$ will be the expression for $(\nabla_\mathcal P)^{an}$ deduced in $(3.4.16)$.
\begin{proposition}
Trivialize $\mathcal L_1^{an}$ resp. $\Omega^1_{E^{an}/S^{an}}$ on the universal covering $\C\times \H$ of $E^{an}$ as in Prop. 3.5.2 resp. by $\{ \mathrm{d}z \}$.\footnote{Observe that we have a canonical identification $\mathrm{pr}^*\Omega^1_{E^{an}/S^{an}}\xrightarrow{\sim} \Omega^1_{\C\times \H/\H}$.} Then $(\nabla_1^{res})^{an}$ writes on sections as
\[\begin{pmatrix} \lambda(z,\tau) \\ \mu(z,\tau) \\ \nu(z,\tau) \end{pmatrix}\mapsto \begin{pmatrix} \partial_z \lambda(z,\tau) \\ \partial_z\mu(z,\tau)+\eta(1,\tau)\cdot \lambda(z,\tau) \\ \partial_z\nu(z,\tau)+\lambda(z,\tau) \end{pmatrix}.\]
\end{proposition}
\begin{proof}
At first consider the pullback of the analytified universal integrable connection
\[(\nabla_\mathcal P)^{an}: \mathcal P^{an} \rightarrow \Omega^1_{E^{an}\times_{S^{an}} (\widehat{E}^\natural)^{an}/(\widehat{E}^\natural)^{an}}\otimes_{\mathcal O_{E^{an}\times_{S^{an}} (\widehat{E}^\natural)^{an}}} \mathcal P^{an}\]
along the following commutative diagram in which all arrows are the canonical projections:
\begin{equation*}
\begin{xy}
\xymatrix@C-0.3cm{
\C \times \C^2 \times \H \ar[r]\ar[d] & E^{an}\times_{S^{an}}(\widehat{E}^\natural)^{an} \ar[d]\\
\C^2\times \H \ar[r] & (\widehat{E}^\natural)^{an}}
\end{xy}
\end{equation*}
Taking into account our fixed trivialization of $\widetilde{\mathcal P}$ (cf. $(3.5.4)$) the obtained integrable connection
\[\widetilde{\nabla}: \widetilde{\mathcal P} \rightarrow \Omega^1_{\C \times \C^2 \times \H/\C^2 \times \H } \otimes_{\mathcal O_{\C \times \C^2 \times \H}} \widetilde{\mathcal P}\]
writes as a map
\[ \ \ \ \ \widetilde{\nabla}: \mathcal O_{\C\times \C^2 \times\H} \cdot t \rightarrow \mathcal O_{\C\times \C^2 \times \H}\cdot (\mathrm{d}z\otimes t)\]
and is determined by the value at the section $t$. It is readily deduced from $(3.4.16)$ that this value is
\[\widetilde{\nabla}(t)=(\eta(1,\tau)w+u)\cdot (\mathrm{d}z\otimes t).\]
Altogether, we obtain
\[\tag{$*$} \begin{split}
\widetilde{\nabla}(t)& =(\eta(1,\tau)w+u)\cdot (\mathrm{d}z\otimes t),\\
\widetilde{\nabla}(wt)& =(\eta(1,\tau)w^2+wu)\cdot (\mathrm{d}z\otimes t),\\
\widetilde{\nabla}(ut)& =(\eta(1,\tau)wu+u^2)\cdot (\mathrm{d}z\otimes t).
\end{split}\]
Note that the sections $t,wt,ut$ of $\widetilde{\mathcal P}$ induce sections of the $\mathcal O_{\C\times\H}$-vector bundle
\[k_1^*\widetilde{\mathcal P} \simeq \tau_1^*\mathcal P_1^{an}\simeq \mathrm{pr}^*(p_1^{an})_*\mathcal P_1^{an}\] 
which (by definition of $k_1$) are just the three basic sections 
\[k_1^{-1}(t) \otimes 1, \quad k_1^{-1}(t)\otimes \eta^\vee, \quad k_1^{-1}(t)\otimes \omega^\vee\]
used for the trivialization of $\mathrm{pr}^*\mathcal L_1^{an}=\mathrm{pr}^*(p_1^{an})_*\mathcal P_1^{an}$ (cf. Def. 3.5.1 resp. Prop. 3.5.2).\\
Observing this and the equalities in $(*)$ one derives the following fact:\\
The pullback of the connection $(\nabla_1^{res})^{an}: \mathcal L_1^{an}\rightarrow \Omega^1_{E^{an}/S^{an}}\otimes_{\mathcal O_{E^{an}}} \mathcal L_1^{an}$ along
\begin{equation*}
\begin{xy}
\xymatrix{
\C \times \H \ar[r]^{ \ \ \mathrm{pr}}\ar[d]_{\mathrm{can}} & E^{an}\ar[d]^{\pi^{an}}\\
\H \ar[r]^{\mathrm{can}} & S^{an}}
\end{xy}
\end{equation*}
is given (and determined) by
\[\begin{split}
k_1^{-1}(t) \otimes 1& \mapsto \eta(1,\tau)\cdot (\mathrm{d}z\otimes (k_1^{-1}(t)\otimes \eta^\vee))+\mathrm{d}z\otimes (k_1^{-1}(t)\otimes \omega^\vee),\\
k_1^{-1}(t)\otimes \eta^\vee& \mapsto 0,\\
k_1^{-1}(t)\otimes \omega^\vee& \mapsto 0.
\end{split}\]
The claim of the proposition then easily follows.
\end{proof}

\subsection{The analytification of the higher logarithm sheaves}
We proceed to study the analytification of the higher logarithm sheaves $\mathcal L_n=\mathrm{Sym}^n_{\mathcal O_E}\mathcal L_1$ of $E/S/\Q$. Its (componentwise) pullback to the universal covering $\C\times \H$ receives a trivialization which is induced by the one constructed for the first logarithm sheaf in 3.5.1. Based on this trivialization we describe the resulting automorphy matrix, determine the formula for the analytified integrable $S$-connection of $\mathcal L_n$ and finally also make explicit the analytified splitting of $\epsilon^*\mathcal L_n$.

\subsubsection{Trivializing the higher logarithm sheaves on the universal covering}
\begin{definition}
(i) The trivializing sections
\[\{ k_1^{-1}(t)\otimes 1, \eta^\vee, \omega^\vee \}\]
for the $\mathcal O_{\C\times \H}$-vector bundle $\mathrm{pr}^*\mathcal L_1^{an}$ (cf. Prop. 3.5.2) will henceforth be abbreviated with
\[\{ e,f,g \}.\]
(ii) For each $n\geq 0$ we will denote by $\nabla_n^{res}$ the restriction to an integrable $S$-connection of the integrable $\Q$-connection $\nabla_n=\mathrm{Sym}^n\nabla_1$ belonging to $\mathcal L_n=\mathrm{Sym}^n_{\mathcal O_E}\mathcal L_1$.\\
(iii) For each $n\geq 0$ we fix the trivialization
\[\tag{\textbf{3.5.10}} \mathcal O_{\C\times \H}^{\oplus r(n)} \simeq \mathrm{pr}^*\mathcal L_n^{an}\]
defined by the ordered basic sections
{\footnotesize
\[\bigg \{\frac{e^n}{n!}, \frac{e^{n-1}f}{(n-1)!}, \frac{e^{n-2}f^2}{(n-2)!},..., \frac{f^n}{(n-n)!}, \frac{e^{n-1}g}{(n-1)!}, \frac{e^{n-2}fg}{(n-2)!},\frac{e^{n-3}f^2g}{(n-3)!},..., \frac{f^{n-1}g}{(n-n)!}, \frac{e^{n-2}g^2}{(n-2)!},.........,\frac{g^n}{(n-n)!} \bigg \},\]}where we let
\[r(n):=\mathrm{rk}_{\mathcal O_E}(\mathcal L_n)=\mathrm{rk}_{\mathcal O_{E^{an}}}(\mathcal L_n^{an})=\mathrm{rk}_{\mathcal O_S}\bigg(\prod_{k=0}^n \mathrm{Sym}^k_{\mathcal O_S}\mathcal H\bigg)=\mathrm{rk}_{\mathcal O_{S^{an}}}\bigg(\prod_{k=0}^n \mathrm{Sym}^k_{\mathcal O_{S^{an}}}\mathcal H^{an}\bigg).\footnote{This number equals $\frac{(n+1)(n+2)}{2}$.}\]
We denote by
\[A_n: \Z^2 \times \Gamma(N) \times \C \times \H \rightarrow \mathrm{GL}_{r(n)}(\C)\]
the associated automorphy matrix for $\mathcal L_n^{an}$.
\end{definition}

We now describe the structure of $A_n$.\\
\newline
Let
\[\Bigg(\begin{pmatrix}  \underline{m} \\ \underline{n} \end{pmatrix}, \begin{pmatrix} a & b \\ c & d \end{pmatrix}\Bigg) \in \Z^2\times \Gamma(N),\]
and let
\[\sigma: \C\times \H \rightarrow \C\times \H,  \quad (z,\tau) \mapsto \bigg(\frac{z+\underline{m}\tau+\underline{n}}{c\tau+d},\frac{a\tau+b}{c\tau+d}\bigg)\]
be the corresponding deck transformation of $\C\times \H$ over $E^{an}$ (we use the underlines so as not to double the notation $n$).\\
Recall from 3.2 (i) that the value
\[A_n \Bigg( \begin{pmatrix}  \underline{m} \\ \underline{n} \end{pmatrix}, \begin{pmatrix} a & b \\ c & d \end{pmatrix}, (z,\tau) \Bigg)\]
is defined by the equation
\begin{align*}
&\bigg(\frac{e^n}{n!}, \frac{e^{n-1}f}{(n-1)!},......,\frac{g^n}{(n-n)!}\bigg)\\
&=\bigg(\sigma^*\frac{e^n}{n!}, \sigma^*\frac{e^{n-1}f}{(n-1)!},......,\sigma^*\frac{g^n}{(n-n)!}\bigg) \cdot A_n \Bigg( \begin{pmatrix}  \underline{m} \\ \underline{n} \end{pmatrix}, \begin{pmatrix} a & b \\ c & d \end{pmatrix}, (z,\tau) \Bigg).
\end{align*}

In other words, the columns of
\[A_n \Bigg( \begin{pmatrix}  \underline{m} \\ \underline{n} \end{pmatrix}, \begin{pmatrix} a & b \\ c & d \end{pmatrix}, (z,\tau) \Bigg)\]
are obtained by expressing each basic section $\frac{e^{n-i-j}f^ig^j}{(n-i-j)!}$ in terms of the $\sigma$-transformed basic sections.\\
But from the already computed matrix
\[A_1 \Bigg( \begin{pmatrix}  \underline{m} \\ \underline{n} \end{pmatrix}, \begin{pmatrix} a & b \\ c & d \end{pmatrix}, (z,\tau) \Bigg)=\begin{pmatrix} 1 & 0 & 0\\
                        2\pi i (d\underline{m}-cz-c\underline{n}) & c\tau+d & 0\\
												0 & 0 & \frac{1}{c\tau+d}
\end{pmatrix}\]
(cf. Prop. 3.5.2) we obtain all the information we need:\\
\newline
Namely, for all $0 \leq i,j$ with $i+j\leq n$ we calculate:
\begin{align*}
&\frac{e^{n-i-j}f^ig^j}{(n-i-j)!}=\frac{(\sigma^*e+2\pi i (d\underline{m}-cz-c\underline{n})\sigma^*f)^{n-i-j}((c\tau+d)\sigma^*f)^i((c\tau+d)^{-1}\sigma^*g)^j}{(n-i-j)!}\\
&=\sum_{k=0}^{n-i-j}\frac{1}{(n-i-j-k)!\cdot k!}(\sigma^*e)^k(2\pi i (d\underline{m}-cz-c\underline{n}))^{n-i-j-k}(\sigma^*f)^{n-j-k}(\sigma^*g)^j (c\tau+d)^{i-j}\\
&=\sum_{k=0}^{n-i-j} \frac{1}{(n-i-j-k)!}\bigg(\frac{2\pi i (d\underline{m}-cz-c\underline{n})}{c\tau+d}\bigg)^{n-i-j-k}\cdot (c\tau+d)^{n-2j-k}\cdot \sigma^*\frac{e^kf^{n-j-k}g^j}{k!},
\end{align*}
and from the last line one can read off all entries of the respective column in our matrix.\\
\newline
Let us express this last line somewhat differently: we consider the factor of automorphy
\[\begin{split} a & :\Z^2 \times \Z^2 \times \Gamma(N) \times \C \times \C^2 \times \H \rightarrow \C^*, \quad \Bigg( \begin{pmatrix}  \underline{m} \\ \underline{n} \end{pmatrix}, \begin{pmatrix} m' \\ n' \end{pmatrix}, \begin{pmatrix} a & b \\ c & d \end{pmatrix}, (z,w,u,\tau) \Bigg) \mapsto\\
& \exp\bigg[-\frac{2\pi i c}{c\tau+d}\cdot(z+\underline{m}\tau+\underline{n})(w+m'\tau+n')+2\pi i m'z+2\pi i \underline{m}w+2\pi i \underline{m}m'\tau\bigg].\end{split}\]of $\mathcal P^{an}$ (cf. $(3.4.14)$) and make the following
\begin{definition}
Let
\[\Bigg(\begin{pmatrix}  \underline{m} \\ \underline{n} \end{pmatrix}, \begin{pmatrix} a & b \\ c & d \end{pmatrix}, (z,\tau)\Bigg) \in \Z^2 \times \Gamma(N)\times \C\times \H.\]
Then, for each $r\in \Z$ we define
\[a_r\Bigg(\begin{pmatrix}  \underline{m} \\ \underline{n} \end{pmatrix}, \begin{pmatrix} a & b \\ c & d \end{pmatrix}, (z,\tau) \Bigg)\]
to be the coefficient at $w^r$ in the expansion around $w=0$ of
\[a\Bigg( \begin{pmatrix}  \underline{m} \\ \underline{n} \end{pmatrix}, \begin{pmatrix} 0 \\ 0 \end{pmatrix}, \begin{pmatrix} a & b \\ c & d \end{pmatrix}, (z,w,u,\tau) \Bigg).\]
Hence, it is zero for all $r <0$ and
\[a_r\Bigg(\begin{pmatrix}  \underline{m} \\ \underline{n} \end{pmatrix}, \begin{pmatrix} a & b \\ c & d \end{pmatrix}, (z,\tau) \Bigg)=\frac{1}{r!}\bigg(\frac{2\pi i (d\underline{m}-cz-c\underline{n})}{c\tau+d}\bigg)^r\]
for all $r\geq 0$.
\end{definition}
Our previous transformation formula, determining $A_n$, may then conveniently be rewritten as
\[\frac{e^{n-i-j}f^ig^j}{(n-i-j)!}=\sum_{k=0}^n a_{n-i-j-k}\Bigg(\begin{pmatrix}  \underline{m} \\ \underline{n} \end{pmatrix}, \begin{pmatrix} a & b \\ c & d \end{pmatrix}, (z,\tau) \Bigg)\cdot (c\tau+d)^{n-2j-k}\cdot \sigma^*\frac{e^kf^{n-j-k}g^j}{k!}.\]With notation as above we have thus shown:
\begin{proposition}
The matrix
\[A_n \Bigg( \begin{pmatrix}  \underline{m} \\ \underline{n} \end{pmatrix}, \begin{pmatrix} a & b \\ c & d \end{pmatrix}, (z,\tau) \Bigg) \in \mathrm{GL}_{r(n)}(\C)\]
is determined by the equation
\[\frac{e^{n-i-j}f^ig^j}{(n-i-j)!}=\sum_{k=0}^n a_{n-i-j-k}\Bigg(\begin{pmatrix}  \underline{m} \\ \underline{n} \end{pmatrix}, \begin{pmatrix} a & b \\ c & d \end{pmatrix}, (z,\tau) \Bigg)\cdot (c\tau+d)^{n-2j-k}\cdot \sigma^*\frac{e^kf^{n-j-k}g^j}{k!}.\]Concretely, this means:\\
Given $1\leq r,s \leq r(n)$, the entry $A_n(r,s)$ of the $s$-th column and $r$-th line of the above matrix is obtained as follows:\\
Take the $s$-th of the sections
{\footnotesize
\[\bigg \{\frac{e^n}{n!}, \frac{e^{n-1}f}{(n-1)!}, \frac{e^{n-2}f^2}{(n-2)!},..., \frac{f^n}{(n-n)!}, \frac{e^{n-1}g}{(n-1)!}, \frac{e^{n-2}fg}{(n-2)!},\frac{e^{n-3}f^2g}{(n-3)!},..., \frac{f^{n-1}g}{(n-n)!}, \frac{e^{n-2}g^2}{(n-2)!},.........,\frac{g^n}{(n-n)!}\bigg \}.\]}It is of the form $\frac{e^{n-i-j}f^ig^j}{(n-i-j)!}$ for some $0\leq i,j$ with $i+j\leq n$.\\
Now consider the $r$-th of these sections: if it is of the form $\frac{e^kf^{n-j-k}g^j}{k!}$ for some $0\leq k \leq n$, then
\[A_n(r,s)=a_{n-i-j-k}\Bigg(\begin{pmatrix}  \underline{m} \\ \underline{n} \end{pmatrix}, \begin{pmatrix} a & b \\ c & d \end{pmatrix}, (z,\tau) \Bigg)\cdot (c\tau+d)^{n-2j-k}.\]
Otherwise we have $A_n(r,s)=0$. \qquad \qed
\end{proposition}
With the preceding proposition one can visualize the matrix $A_n \Bigg( \begin{pmatrix}  \underline{m} \\ \underline{n} \end{pmatrix}, \begin{pmatrix} a & b \\ c & d \end{pmatrix}, (z,\tau) \Bigg)$:\\
\newline
Namely, leaving away the argument $\Bigg( \begin{pmatrix}  \underline{m} \\ \underline{n} \end{pmatrix}, \begin{pmatrix} a & b \\ c & d \end{pmatrix}, (z,\tau) \Bigg)$ of $A_n$ and of the $a_r$, one figures out that:
\[\tag{\textbf{3.5.11}}
A_n=
\left( \begin{array}{cccc|cccc}
a_0	& 0 	& \dots  & 0& 0& \dots & \dots & 0 \\
	(c\tau+d)\cdot a_1& (c\tau+d)\cdot a_0	& \dots   &0 &0 & \dots & \dots & 0 \\
	\vdots & \vdots  	&   &  \vdots &   \vdots & & & \vdots \\
	(c\tau+d)^n \cdot a_n& (c\tau+d)^n\cdot a_{n-1}   & \dots     & (c\tau+d)^n\cdot a_0 & 0 & \dots & \dots &0 \\
    \hline
 	0 & 0   & \dots     & 0       & * & \dots & \dots & * \\
	\vdots &  \vdots   &     &\vdots   &\vdots & & & \vdots \\
	\vdots & \vdots & & \vdots & \vdots & & & \vdots \\
	0&   0& \dots    &    0    &* &\dots & \dots & *
\end{array}\right)
\]
The upper left block is a $(n+1)\times (n+1)$-matrix whose $(r,s)$-entry is given by $(c\tau+d)^{r-1}\cdot a_{r-s}$.\\
The upper right resp. lower left block is a $(n+1)\times (r(n)-n-1)$ resp. $(r(n)-n-1) \times (n+1)$-matrix consisting only of zeroes.\\
The lower right block is a $(r(n)-n-1)\times (r(n)-n-1)$-matrix whose detailed entries won't be of importance, but can nevertheless be computed as explained in Prop. 3.5.6.

\subsubsection{Description of the relative connection}
With the knowledge of $(\nabla_1^{res})^{an}$ (cf. Prop. 3.5.3) the computation of the connection
\[(\nabla_n^{res})^{an}: \mathcal L_n^{an}\rightarrow \Omega^1_{E^{an}/S^{an}}\otimes_{\mathcal O_{E^{an}}} \mathcal L_n^{an}\]
is a straightforward affair. At first, let us once again clarify how to understand the result:\\
Letting as usual
\[\mathrm{pr}: (\C\times \H,\mathcal O_{\C\times\H}) \rightarrow E^{an}\]
be the projection of the universal covering, we trivialize $\mathrm{pr}^*\mathcal L_n^{an}$ resp. $\mathrm{pr}^*\Omega^1_{E^{an}/S^{an}}\simeq \Omega^1_{\C\times \H/\H}$ as in $(3.5.10)$ resp. by $\{\mathrm{d}z\}$. Recall from 3.2 that then a section of $\mathcal L_n^{an}$ resp. of $\Omega^1_{E^{an}/S^{an}}\otimes_{\mathcal O_{E^{an}}} \mathcal L_n^{an}$ over an open subset $V \subseteq E^{an}$ is given by a vector (of length $r(n)$) of holomorphic functions on $\mathrm{pr}^{-1}(V)$:
\[\begin{pmatrix} l_{0,0}(z,\tau) \\ l_{1,0}(z,\tau) \\ l_{2,0}(z,\tau)\\ \vdots\\ \vdots\\l_{0,n}(z,\tau)\end{pmatrix} \quad \textrm{resp.} \quad \begin{pmatrix} \widetilde{l}_{0,0}(z,\tau) \\ \widetilde{l}_{1,0}(z,\tau) \\ \widetilde{l}_{2,0}(z,\tau)\\ \vdots\\ \vdots\\\widetilde{l}_{0,n}(z,\tau)\end{pmatrix},\]
transforming under the effect of $\Z^2\times \Gamma(N)$ with the matrix $(3.5.11)$ resp. with the matrix obtained from $(3.5.11)$ by multiplying each entry with $c\tau+d$.\\
The indices $i,j$ of the function $l_{i,j}(z,\tau)$ resp. $\widetilde{l}_{i,j}(z,\tau)$ refer to the basic section
\[\frac{e^{n-i-j}f^ig^j}{(n-i-j)!} \quad \textrm{resp.} \quad \mathrm{d}z\otimes \frac{e^{n-i-j}f^ig^j}{(n-i-j)!}\]
in the trivialization of $\mathrm{pr}^*\mathcal L_n^{an}$ resp. $\mathrm{pr}^*(\Omega^1_{E^{an}/S^{an}}\otimes_{\mathcal O_{E^{an}}} \mathcal L_n^{an})$ to which the function is multiplied.\\
\newline
With these explanations we have:
\begin{proposition}
$(\nabla_n^{res})^{an}$ is given on sections of $\mathcal L_n^{an}$ by the rule
\[\begin{pmatrix} l_{0,0}(z,\tau) \\ l_{1,0}(z,\tau) \\ l_{2,0}(z,\tau)\\ \vdots\\ \vdots\\l_{0,n}(z,\tau) \end{pmatrix}\mapsto \begin{pmatrix} \widetilde{l}_{0,0}(z,\tau) \\ \widetilde{l}_{1,0}(z,\tau) \\ \widetilde{l}_{2,0}(z,\tau)\\ \vdots\\ \vdots\\\widetilde{l}_{0,n}(z,\tau)\end{pmatrix},\]
where for each $0\leq i,j$ with $i+j\leq n$ we have
\[\widetilde{l}_{i,j}(z,\tau)=\partial_zl_{i,j}(z,\tau)+\eta(1,\tau)\cdot l_{i-1,j}(z,\tau)+l_{i,j-1}(z,\tau).\]
Here, we set $l_{-1,j}(z,\tau)=l_{i,-1}(z,\tau)\equiv 0$ for all $i,j$.
\end{proposition}
\begin{proof}
Use that the pullback of $(\nabla_1^{res})^{an}$ along
\begin{equation*}
\begin{xy}
\xymatrix{
\C \times \H \ar[r]^{ \ \ \mathrm{pr}}\ar[d]_{\mathrm{can}} & E^{an}\ar[d]^{\pi^{an}}\\
\H \ar[r]^{\mathrm{can}} & S^{an}}
\end{xy}
\end{equation*}
is given by
\[\begin{split}
e& \mapsto \eta(1,\tau)\cdot (\mathrm{d}z\otimes f)+\mathrm{d}z\otimes g,\\
f& \mapsto 0,\\
g& \mapsto 0.
\end{split}\]
(cf. the proof of Prop. 3.5.3) and that $(\nabla_n^{res})^{an}=\mathrm{Sym}^n(\nabla_1^{res})^{an}$. The rest is routine.
\end{proof}
\subsubsection{Description of the splitting}
Finally, we also want to express explicitly (on the fixed connected component) the analytification
\[\varphi_n^{an}: \prod_{k=0}^n \mathrm{Sym}^k_{\mathcal O_{S^{an}}}\mathcal H^{an} \simeq (\epsilon^{an})^*\mathcal L_n^{an}\]
of the splitting
\[\varphi_n: \prod_{k=0}^n \mathrm{Sym}^k_{\mathcal O_S}\mathcal H \simeq \epsilon^*\mathcal L_n\]
in the formalism of automorphy matrices.\\
\newline
Recall (from the preliminary remarks preceding 3.5.1) that the pullback of $\prod_{k=0}^n \mathrm{Sym}^k_{\mathcal O_{S^{an}}}\mathcal H^{an}$ to the universal covering $\H$ of $S^{an}$ is trivialized by the ordered basic sections

{\footnotesize
\[\bigg \{\frac{1}{n!}, \frac{\eta^\vee}{(n-1)!}, \frac{(\eta^\vee)^2}{(n-2)!},..., \frac{(\eta^\vee)^n}{(n-n)!}, \frac{\omega^\vee}{(n-1)!}, \frac{\eta^\vee\omega^\vee}{(n-2)!},\frac{(\eta^\vee)^2\omega^\vee}{(n-3)!},..., \frac{(\eta^\vee)^{n-1}\omega^\vee}{(n-n)!}, \frac{(\omega^\vee)^2}{(n-2)!},.........,\frac{(\omega^\vee)^n}{(n-n)!} \bigg \}.\]}

Recall furthermore (from Def. 3.5.4) that the pullback of $\mathcal L_n^{an}$ to the universal covering $\C\times \H$ of $E^{an}$ is trivialized by the ordered basic sections

{\footnotesize
\[\bigg \{\frac{e^n}{n!}, \frac{e^{n-1}f}{(n-1)!}, \frac{e^{n-2}f^2}{(n-2)!},..., \frac{f^n}{(n-n)!}, \frac{e^{n-1}g}{(n-1)!}, \frac{e^{n-2}fg}{(n-2)!},\frac{e^{n-3}f^2g}{(n-3)!},..., \frac{f^{n-1}g}{(n-n)!}, \frac{e^{n-2}g^2}{(n-2)!},.........,\frac{g^n}{(n-n)!} \bigg \}\]}

and that we obtain from it (as explained in 3.2 $(v)$ and the preliminary remarks of 3.4) an induced trivialization for the pullback of $(\epsilon^{an})^*\mathcal L_n^{an}$ to $\H$.\\
One deduces directly from Prop. 3.5.6 and the explicit formula for $a_r$ (cf. Def. 3.5.5) that the associated automorphy matrix for $(\epsilon^{an})^*\mathcal L_n^{an}$, evaluated at
\[\Bigg(\begin{pmatrix} a & b \\ c & d \end{pmatrix},\tau\Bigg) \in \Gamma(N)\times \H,\]
is a diagonal matrix of size $r(n)$ whose nontrivial entry in the $r$-th line is given as follows:\\
Take the $r$-th of the sections

{\footnotesize
\[\bigg \{\frac{e^n}{n!}, \frac{e^{n-1}f}{(n-1)!}, \frac{e^{n-2}f^2}{(n-2)!},..., \frac{f^n}{(n-n)!}, \frac{e^{n-1}g}{(n-1)!}, \frac{e^{n-2}fg}{(n-2)!},\frac{e^{n-3}f^2g}{(n-3)!},..., \frac{f^{n-1}g}{(n-n)!}, \frac{e^{n-2}g^2}{(n-2)!},.........,\frac{g^n}{(n-n)!} \bigg \}.\]}

It is of the form $\frac{e^{n-i-j}f^ig^j}{(n-i-j)!}$ for some $0\leq i,j$ with $i+j\leq n$, and the entry is then given by $(c\tau+d)^{i-j}$.\\
The thus obtained automorphy matrix for $(\epsilon^{an})^*\mathcal L_n^{an}$ hence coincides with the automorphy matrix for $\prod_{k=0}^n \mathrm{Sym}^k_{\mathcal O_{S^{an}}}\mathcal H^{an}$ (cf. again the remarks preceding 3.5.1).

\begin{proposition}
Trivializing the pullback of $\prod_{k=0}^n \mathrm{Sym}^k_{\mathcal O_{S^{an}}}\mathcal H^{an}$ resp. $(\epsilon^{an})^*\mathcal L_n^{an}$  to the universal covering $\H$ of $S^{an}$ as explained above, the associated automorphy matrices are equal. Expressing as usual the sections of these $\mathcal O_{S^{an}}$-vector bundles as vectors of holomorphic functions on (open subsets of) $\H$ which transform under the effect of $\Gamma(N)$ with this automorphy matrix, the isomorphism
\[\varphi_n^{an}: \prod_{k=0}^n \mathrm{Sym}^k_{\mathcal O_{S^{an}}}\mathcal H^{an} \simeq (\epsilon^{an})^*\mathcal L_n^{an}\]
is given as the identity on such vectors.\\
In particular, given a section
\[\begin{pmatrix} l_{0,0}(z,\tau) \\ l_{1,0}(z,\tau) \\ l_{2,0}(z,\tau)\\ \vdots\\ \vdots\\l_{0,n}(z,\tau)\end{pmatrix}\]
of $\mathcal L_n^{an}$ which is defined on some open subset of $E^{an}$ intersecting $\epsilon^{an}(S^{an})$, the induced section of $\prod_{k=0}^n \mathrm{Sym}^k_{\mathcal O_{S^{an}}}\mathcal H^{an}$ (via pullback along $\epsilon^{an}$ and using $\varphi_n^{an}$) is given by
\[\begin{pmatrix} l_{0,0}(0,\tau) \\ l_{1,0}(0,\tau) \\ l_{2,0}(0,\tau)\\ \vdots\\ \vdots\\l_{0,n}(0,\tau)\end{pmatrix}.\]
\end{proposition}
\begin{proof}
The first assertion was already shown, and the third follows from the second together with the explanations of 3.2 $(v)$. It thus remains to verify the second claim of the proposition.\\
For this let us write $p: \H\rightarrow S^{an}$ and $\mathrm{pr}: \C\times \H \rightarrow E^{an}$ for the projections and consider the pullback
\[p^*(\varphi_n^{an}): p^*(\epsilon^{an})^*\mathcal L_n^{an} \simeq p^*\prod_{k=0}^n \mathrm{Sym}^k_{\mathcal O_{S^{an}}}\mathcal H^{an}\]
of the isomorphism $\varphi_n^{an}$ along $p$. Combining it with the identification
\[(\widetilde{\epsilon})^*\mathrm{pr}^*\mathcal L_n^{an}\simeq p^*(\epsilon^{an})^*\mathcal L_n^{an},\]
induced by the commutative diagram
\begin{equation*}
\begin{xy}
\xymatrix{
\H \ar[r]^{\widetilde{\epsilon} \quad }\ar[d]_{p} & \C\times\H\ar[d]^{\mathrm{pr}}\\
S^{an} \ar[r]^{\epsilon^{an}} & E^{an}}
\end{xy}
\end{equation*}
with $\widetilde{\epsilon}(\tau):=(0,\tau)$, we obtain the isomorphism
\[(\widetilde{\epsilon})^*\mathrm{pr}^*\mathcal L_n^{an} \simeq p^*\prod_{k=0}^n \mathrm{Sym}^k_{\mathcal O_{S^{an}}}\mathcal H^{an}.\]
What we have to show is that hereunder the trivializing section
\[(\widetilde{\epsilon})^*\Big(\frac{e^{n-i-j}f^ig^j}{(n-i-j)!}\Big)\]
of the left side is sent to the following trivializing section of the right side:
\[\frac{(\eta^\vee)^i (\omega^\vee)^j}{(n-i-j)!}.\]
Because of $\varphi_n=\mathrm{Sym}^n\varphi_1$ one is quickly reduced to showing that
\[(\widetilde{\epsilon})^*\mathrm{pr}^*\mathcal L_1^{an} \simeq p^*(\mathcal O_{S^{an}}\oplus \mathcal H^{an}) \simeq \mathcal O_{\H}\oplus \mathcal O_{\H}\cdot \eta^\vee \oplus \mathcal O_{\H}\cdot \omega^\vee\]
sends $\{ (\widetilde{\epsilon})^*(e),(\widetilde{\epsilon})^*(f),(\widetilde{\epsilon})^*(g) \}$ to $\{ 1, \eta^\vee, \omega^\vee \}$.\\
For $(\widetilde{\epsilon})^*(f)$ and $(\widetilde{\epsilon})^*(g)$ this follows directly from the fact that $f$ and $g$ are the images of $\eta^\vee$ and $\omega^\vee$ under the inclusion of $(3.5.8)$ (cf. the comments preceding Prop. 3.5.2).\\
The claim for $(\widetilde{\epsilon})^*(e)$ follows by recalling how $\varphi_1$ was constructed via the rigidification $r$ of $\mathcal P$ (cf. Cor. 2.3.2 resp. $(2.3.13)$), that $r$ is induced by $r^0$ (cf. 0.1.1), how $r^0$ expresses for elliptic curves (cf. $(3.1.5)$) and by then observing the definition of $t^0$ (cf. $(3.4.13)$) and Rem. 3.4.1.
\end{proof}

\subsection{The pullback along torsion sections}
Given a torsion section $t: S\rightarrow E$ of $E/S$ we can pull back sections of $\Omega^1_{E/\Q}\otimes_{\mathcal O_E}\mathcal L_n$ along $t$ and use the invariance property $t^*\mathcal L_n \simeq \epsilon^*\mathcal L_n$ (cf. Lemma 1.4.5) together with the splitting $\varphi_n$ to obtain a section of $\Omega^1_{S/\Q}\otimes_{\mathcal O_S}\prod_{k=0}^n\mathrm{Sym}^k_{\mathcal O_S}\mathcal H$. This is basically the procedure to be carried out when one wants to specialize de Rham cohomology classes with coefficients in $\mathcal L_n$ along $t$.\\
On the analytic side, if we again (componentwise) trivialize the occurring vector bundles on the universal coverings, we can ask for the explicit formula of the preceding process. Its deduction is the main goal of this subsection, and it will present a key tool for our later (analytic and finally also algebraic) computation of the specialization of the $D$-variant polylogarithm classes along torsion sections.\\
To derive this formula we proceed in several steps. The first of these is the crucial one and consists in determining explicitly the effect of translating sections of $\mathcal L_1^{an}$ by a torsion section and subsequently using the invariance of $\mathcal L_1^{an}$ under such translations. To find the formula for this operation is the essential (and rather technical) task to be solved as it then immediately implies the corresponding result for general $n$, with which in turn one quickly arrives at the solution for the initial problem.\\
This subsection also contains (prior to Prop. 3.5.14) an insertion about how we trivialize the pullbacks of various differential modules to the universal coverings, recording in addition the associated automorphy matrices; these conventions will remain in force until the end of the work.

\subsubsection{Notations and conventions}

Recall from 3.4 that $(e_1,e_2) \in E[N](S)$ denotes the Drinfeld basis for $E[N]$ associated with $E/S$.\\
\newline
We let $a,b$ be two fixed integers and consider the associated $N$-torsion section
\[t_{a,b}:=ae_1+be_2: S \rightarrow E\]
which of course only depends on the class of $(a,b)$ in $(\Z/N\Z)^2$.\\
Moreover, we write
\[T_{a,b}: E \rightarrow E\]
for the translation by the torsion section $t_{a,b}$.\\
Recall from 1.4.2 that the invariance of $\mathcal L_n$ by $N$-multiplication
\[\mathcal L_n \simeq [N]^*\mathcal L_n\]
together with the equality $[N]\circ T_{a,b}=[N]$ yields a canonical isomorphism of $\mathcal D_{E/\Q}$-modules
\[T_{a,b}^*\mathcal L_n \simeq T_{a,b}^*[N]^*\mathcal L_n \simeq [N]^* \mathcal L_n \simeq  \mathcal L_n,\]
and that by further pullback via $\epsilon$ one obtains a canonical $\mathcal D_{S/\Q}$-linear identification
\[t_{a,b}^*\mathcal L_n  \simeq \epsilon^*\mathcal L_n.\]
Cf. Cor. 1.4.4 and Lemma 1.4.5.\\
\newline
Let us now consider the analytic situation.\\
We continue to work on a fixed connected component of $E^{an}/S^{an}$ but as usual suppress its index in the notation of the geometric objects, writing again e.g. $E^{an}, S^{an}, \pi^{an}, \epsilon^{an}, t_{a,b}^{an}, T_{a,b}^{an}, [N]^{an}, \mathcal L_n^{an}...$ for the restrictions of the honest analytifications to the chosen component.\\
With these conventions, note that by analytifying the algebraic isomorphisms recalled above and then restricting to the fixed component we obtain identifications
\[\tag{\textbf{3.5.12}} (T_{a,b}^{an})^*\mathcal L_n^{an} \simeq (T_{a,b}^{an})^*([N]^{an})^*\mathcal L_n^{an} \simeq ([N]^{an})^* \mathcal L_n^{an} \simeq  \mathcal L_n^{an}\]
resp.
\[\tag{\textbf{3.5.13}} (t_{a,b}^{an})^*\mathcal L_n^{an}  \simeq (\epsilon^{an})^*\mathcal L_n^{an}.\]
For the subsequent computations it is necessary to specify the chosen component.\\
Let us hence determine that the fixed connected component of $E^{an}/S^{an}$ is the one associated with the class of $j_0$ in $(\Z/N\Z)^*$, where $j_0$ is an integer with $(j_0,N)=1$.\\
\newline
The analytified torsion section $t_{a,b}^{an}: S^{an}\rightarrow E^{an}$ is then given by
\[t_{a,b}^{an}: \Gamma(N) \backslash \H \rightarrow (\Z^2\times \Gamma(N)) \backslash (\C\times \H),\quad \tau\mapsto \Big(\frac{aj_0\tau}{N}+\frac{b}{N},\tau\Big),\]
and the analytified translation $T_{a,b}^{an}: E^{an} \rightarrow E^{an}$ writes as
\[T_{a,b}^{an}: (\Z^2\times \Gamma(N)) \backslash (\C\times \H) \rightarrow (\Z^2\times \Gamma(N)) \backslash (\C\times \H), \quad (z,\tau)\mapsto \Big(z+\frac{aj_0\tau}{N}+\frac{b}{N},\tau \Big).\]
Cf. $(3.4.1)$-$(3.4.4)$.

\subsubsection{Explicit formulas I}

Consider $(3.5.12)$ for $n=1$ and assume a section of $\mathcal L_1^{an}$ over an open subset $V \subseteq E^{an}$ is given.\\
By pullback along $T_{a,b}^{an}$ and using $(3.5.12)$ we obtain a section of $\mathcal L_1^{an}$ defined over $(T^{an}_{a,b})^{-1}(V)$.\\
With respect to the fixed trivialization of $\mathrm{pr}^*\mathcal L_1^{an}$ on $\C\times \H$ (cf. Def. 3.5.1 resp. Prop. 3.5.2) we express the section we started with as usual by a vector
\[\begin{pmatrix}  l_{0,0}(z,\tau) \\ l_{1,0}(z,\tau) \\ l_{0,1}(z,\tau)\end{pmatrix}\]
of holomorphic functions on $\mathrm{pr}^{-1}(V)$ which under the effect of $\Z^2\times \Gamma(N)$ transforms with the automorphy matrix $A_1$ for $\mathcal L_1^{an}$.\\
The section of $\mathcal L_1^{an}$ obtained from it in the way just described writes similarly as such a vector of functions, defined on $\mathrm{pr}^{-1}((T^{an}_{a,b})^{-1}(V))$, and explicitly we have:
\begin{proposition}
The section in question is given by
\[\begin{pmatrix}  l_{0,0}(z+\frac{aj_0\tau}{N}+\frac{b}{N},\tau) \\ l_{1,0}(z+\frac{aj_0\tau}{N}+\frac{b}{N},\tau)-2\pi i \frac{aj_0}{N} \cdot l_{0,0}(z+\frac{aj_0\tau}{N}+\frac{b}{N},\tau)\\ l_{0,1}(z+\frac{aj_0\tau}{N}+\frac{b}{N},\tau)\end{pmatrix}.\]
\end{proposition}
\begin{proof}
We define maps on the universal covering $\C\times \H$ of $E^{an}$:
\begin{align*}
\widetilde{T}_{a,b}&:\C\times \H \rightarrow \C\times \H, \quad (z,\tau) \mapsto \Big(z+\frac{aj_0\tau}{N}+\frac{b}{N},\tau \Big)\\
\widetilde{[N]}&:\C\times \H \rightarrow \C\times \H, \quad (z,\tau) \mapsto (Nz,\tau)
\end{align*}
which then fit into commutative squares
\begin{equation*}
\begin{xy}
\xymatrix{
\C\times \H \ar[d]_{\mathrm{pr}} \ar[r]^{\widetilde{T}_{a,b}} & \C\times\H \ar[d]^{\mathrm{pr}} & & \C\times \H \ar[d]_{\mathrm{pr}} \ar[r]^{\widetilde{[N]}} & \C\times\H \ar[d]^{\mathrm{pr}}\\
E^{an} \ar[r]^{T_{a,b}^{an}} & E^{an} & & E^{an} \ar[r]^{[N]^{an}} & E^{an} }
\end{xy}
\end{equation*}
With these diagrams fixed, the pullback of $(3.5.12)$ via $\mathrm{pr}$ induces a chain of isomorphisms on $\C\times \H$:
\[\tag{$*$} \widetilde{T}_{a,b}^* \ \mathrm{pr}^*\mathcal L_1^{an} \simeq \widetilde{T}_{a,b}^*\widetilde{[N]}^*\mathrm{pr}^*\mathcal L_1^{an} \simeq \widetilde{[N]}^* \mathrm{pr}^*\mathcal L_1^{an} \simeq  \mathrm{pr}^*\mathcal L_1^{an}.\]
We now investigate what happens under $(*)$ with the following (global) sections of $\widetilde{T}_{a,b}^* \ \mathrm{pr}^*\mathcal L_1^{an}$:
\[\widetilde{T}_{a,b}^*(e), \ \ \widetilde{T}_{a,b}^*(f), \ \ \widetilde{T}_{a,b}^*(g),\]
where as usual $\{e,f,g\}$ are the fixed trivializing sections of $\mathrm{pr}^*\mathcal L_1^{an}$ (cf. Def. 3.5.4 resp. Prop. 3.5.2).\\
\newline
Let us begin with $e$. Recall from 3.5.1 that it was constructed as the section
\[e=k_1^{-1}(t)\otimes 1 \in \Gamma(\C\times \H,\mathrm{pr}^*\mathcal L_1^{an}),\]
where one uses the $\mathcal O_{\C\times \H}$-linear isomorphism of $(3.5.6)$:
\[k_1^{-1}\widetilde{\mathcal P}\otimes_{k_1^{-1}\mathcal O_{\C\times \C^2 \times \H}}(\mathcal O_{\C\times \H} \oplus \mathcal O_{\C\times \H}\cdot \eta^\vee \oplus \mathcal O_{\C\times \H}\cdot \omega^\vee)=k_1^*\widetilde{\mathcal P}\simeq \tau_1^*\mathcal P_1^{an}\simeq \mathrm{pr}^*(p^{an}_1)_*\mathcal P_1^{an}\]
and the trivializing section $t$ of $\widetilde{\mathcal P}$ on $\C\times \C^2\times \H$:
\[\mathcal O_{\C \times \C^2 \times \H} \xrightarrow{\sim} \widetilde{\mathcal P} \quad 1\mapsto t=\frac{1}{J(z,-w,\tau)}\otimes \omega_{\mathrm{can}}^\vee.\]
The isomorphism
\[\mathrm{pr}^*\mathcal L_1^{an} \xrightarrow{\sim} \widetilde{[N]}^*\mathrm{pr}^*\mathcal L_1^{an}\]
is given on the basic sections $\{ e,f,g \}$ of $\mathrm{pr}^*\mathcal L_1^{an}$ as follows:
\[\tag{$**$} \begin{split}
e & \mapsto \widetilde{[N]}^*(e)\\
f & \mapsto N\cdot \widetilde{[N]}^*(f)\\
g & \mapsto N\cdot \widetilde{[N]}^*(g).
\end{split} \]
We provide a detailed justification for this fact after we have finished the present proof.\\
This in particular implies that the image of $\widetilde{T}_{a,b}^*(e)$ in $\widetilde{T}_{a,b}^*\widetilde{[N]}^*\mathrm{pr}^*\mathcal L_1^{an}$ under the first isomorphism in $(*)$ is given by $(\widetilde{[N]} \circ \widetilde{T}_{a,b})^*(e)$. Because of
\[(\widetilde{[N]} \circ \widetilde{T}_{a,b})(z,\tau)=(Nz+aj_0\tau+b,\tau)\]
and the equality
\[\frac{1}{J(Nz+aj_0\tau+b,-w,\tau)}\otimes \omega^\vee_{\mathrm{can}}=\mathrm{exp}(-2\pi i aj_0 w)\cdot \frac{1}{J(Nz,-w,\tau)}\otimes \omega^\vee_{\mathrm{can}},\]
deduced from Cor. 3.3.14, the further image under the second arrow in $(*)$ is $\widetilde{[N]}^*(e-2\pi i aj_0 \cdot f)$: to see this recall the definition of $k_1$ and that the section $f$ equals $k^{-1}_1(wt)\otimes 1$.\\
Finally, under the last isomorphism $\widetilde{[N]}^*\mathrm{pr}^*\mathcal L_1^{an} \simeq \mathrm{pr}^*\mathcal L_1^{an}$ in $(*)$ this section $\widetilde{[N]}^*(e-2\pi i aj_0 \cdot f)$ goes to $e-2\pi i \frac{aj_0}{N}\cdot f$ by observing again $(**)$.\\
Altogether, we see that $\widetilde{T}_{a,b}^*(e)$ is indeed mapped to $e-2\pi i \frac{aj_0}{N}\cdot f$ under the chain $(*)$.\\
Observing that $f$ resp. $g$ comes via $k_1$ from $wt$ resp. from $ut$ one deduces with similar arguments that $\widetilde{T}_{a,b}^*(f)$ resp. $\widetilde{T}_{a,b}^*(g)$ maps to $f$ resp. to $g$ under $(*)$.\\
\newline
It is then clear that a section of $\mathrm{pr}^*\mathcal L_1^{an}$:
\[l_{0,0}(z,\tau)\cdot e +l_{1,0}(z,\tau)\cdot f +l_{0,1}(z,\tau)\cdot g\]
coming from a section of $\mathcal L_1^{an}$ goes under pullback via $\widetilde{T}_{a,b}$ and the $\mathcal O_{\C\times \H}$-linear chain $(*)$ to
\[\begin{split}
&l_{0,0}\Big(z+\frac{aj_0\tau}{N}+\frac{b}{N},\tau \Big)\cdot e\\
&+\bigg[l_{1,0}\Big(z+ \frac{aj_0\tau}{N}+\frac{b}{N},\tau\Big)-2\pi i \frac{aj_0}{N}\cdot l_{0,0}\Big(z+\frac{aj_0\tau}{N}+\frac{b}{N},\tau \Big)\bigg]\cdot f\\
&+l_{0,1}\Big(z+\frac{aj_0\tau}{N}+\frac{b}{N},\tau\Big)\cdot g.
\end{split}
\]
This readily implies the claim of the proposition.
\end{proof}

\subsubsection{Supplements to the proof of Prop. 3.5.9}
It remains to verify $(**)$ of the previous proof. This is what we have to do:\\
\newline
Consider the isomorphism
\[\mathrm{pr}^*\mathcal L_1^{an} \simeq \widetilde{[N]}^*\mathrm{pr}^*\mathcal L_1^{an},\]
obtained by pulling back along $\mathrm{pr}$ the analytification
\[\mathcal L_1^{an} \simeq ([N]^{an})^*\mathcal L_1^{an}\]
of the invariance isomorphism in Cor. 1.4.4. Combining it with the $\mathcal O_{\C\times \H}$-linear identification $\mathrm{pr}^*\mathcal L_1^{an}\simeq k_1^*\widetilde{\mathcal P}$ of $(3.5.6)$ we obtain
\[\tag{\textbf{3.5.14}} k_1^*\widetilde{\mathcal P} \simeq \widetilde{[N]}^* k_1^*\widetilde{\mathcal P},\]
where $k_1^*\widetilde{\mathcal P}$ is viewed as $\mathcal O_{\C\times \H}$-module. Then, the content of $(**)$ in the proof of Prop. 3.5.9 is that hereunder $k_1^{-1}(t)\otimes 1$ corresponds to $\widetilde{[N]}^*(k_1^{-1}(t)\otimes 1)$, where $t$ is the familiar section of $\widetilde{\mathcal P}$:
\[t=\frac{1}{J(z,-w,\tau)}\otimes\omega^\vee_{\mathrm{can}},\]
and that moreover $\frac{1}{N}\cdot (k_1^{-1}(t)\otimes \eta^\vee)$ resp. $\frac{1}{N}\cdot (k_1^{-1}(t)\otimes \omega^\vee)$ corresponds to $\widetilde{[N]}^*(k_1^{-1}(t)\otimes \eta^\vee)$ resp. $\widetilde{[N]}^*(k_1^{-1}(t)\otimes \omega^\vee)$.\\
\newline
We show this for $N$ replaced by an arbitrary integer $M\neq 0$.\\
\newline
Suppose for a moment that the following is true:
\begin{lemma}
Set
\[\quad \widetilde{\mathrm{id}_E\times [M]^\natural}: \C \times \C^2 \times \H \rightarrow \C\times \C^2 \times \H, \quad (z,w,u,\tau)\mapsto (z,Mw, Mu,\tau)\]
and
\[\widetilde{[M]\times \mathrm{id}_{\widehat{E}^\natural}}: \C \times \C^2 \times \H \rightarrow \C\times \C^2\times \H, \quad (z,w,u,\tau)\mapsto (Mz,w,u,\tau).\]
Analytifying (and restricting to the fixed component) the canonical isomorphism of $(2.5.2)$:
\[\tag{\textbf{3.5.15}} (\mathrm{id}_E\times [M]^\natural)^*(\mathcal P,r,\nabla_\mathcal P) \simeq (([M]\times \mathrm{id}_{\widehat{E}^\natural})^*\mathcal P,r_{[M]},(\nabla_\mathcal P)_{[M]}),\]
we obtain after pullback to $\C\times \C^2 \times \H$ an induced isomorphism
\[\tag{\textbf{3.5.16}} (\widetilde{\mathrm{id}_E \times [M]^\natural})^*\widetilde{\mathcal P} \simeq (\widetilde{[M]\times \mathrm{id}_{\widehat{E}^\natural}})^*\widetilde{\mathcal P}.\]
Then $(3.5.16)$ maps
\[(\widetilde{\mathrm{id}_E\times [M]^\natural})^*t=\frac{1}{J(z,-Mw,\tau)}\otimes\omega^\vee_{\mathrm{can}} \mapsto \frac{1}{J(Mz,-w,\tau)}\otimes\omega^\vee_{\mathrm{can}}=(\widetilde{[M]\times \mathrm{id}_{\widehat{E}^\natural}})^*t.\]
\end{lemma}
Assuming the content of this lemma we can deduce our above claims as follows.\\
We have a commutative diagram
\begin{equation*} 
\begin{xy}
\xymatrix{
(\C \times \H, \mathcal O_{\C\times \H}\oplus \mathcal O_{\C\times \H}\cdot \eta^\vee\oplus \mathcal O_{\C\times \H}\cdot \omega^\vee) \ar[r]^{\qquad  \qquad k_1}\ar[d] & (\C \times \C^2 \times \H,\mathcal O_{\C \times \C^2 \times \H}) \ar[d]^{\widetilde{\mathrm{id}_E\times [M]^\natural}}\\
(\C \times \H, \mathcal O_{\C\times \H}\oplus \mathcal O_{\C\times \H}\cdot \eta^\vee\oplus \mathcal O_{\C\times \H}\cdot \omega^\vee) \ar[r]^{\qquad  \qquad k_1} & (\C \times \C^2 \times \H,\mathcal O_{\C \times \C^2 \times \H})}
\end{xy}
\end{equation*}
in which the left arrow topologically is the identity and on ring sheaves is given by the isomorphism
\[\tag{\textbf{3.5.17}} \mathcal O_{\C\times \H}\oplus \mathcal O_{\C\times \H}\cdot \eta^\vee\oplus \mathcal O_{\C\times \H}\cdot \omega^\vee \xrightarrow{\sim} \mathcal O_{\C\times \H}\oplus \mathcal O_{\C\times \H}\cdot \eta^\vee\oplus \mathcal O_{\C\times \H}\cdot \omega^\vee\]
defined by
\[f(z,\tau)+g(z,\tau)\cdot \eta^\vee+h(z,\tau)\cdot \omega^\vee \mapsto f(z,\tau)+M g(z,\tau)\cdot \eta^\vee+M h(z,\tau)\cdot\omega^\vee.\]
Moreover, we have a second commutative diagram
\begin{equation*} 
\begin{xy}
\xymatrix{
(\C \times \H, \mathcal O_{\C\times \H}\oplus \mathcal O_{\C\times \H}\cdot \eta^\vee\oplus \mathcal O_{\C\times \H}\cdot \omega^\vee) \ar[r]^{\qquad  \qquad k_1}\ar[d] & (\C \times \C^2 \times \H,\mathcal O_{\C \times \C^2 \times \H}) \ar[d]^{\widetilde{[M]\times \mathrm{id}_{\widehat{E}^\natural}}}\\
(\C \times \H, \mathcal O_{\C\times \H}\oplus \mathcal O_{\C\times \H}\cdot \eta^\vee\oplus \mathcal O_{\C\times \H}\cdot \omega^\vee) \ar[r]^{\qquad  \qquad k_1} & (\C \times \C^2 \times \H,\mathcal O_{\C \times \C^2 \times \H})}
\end{xy}
\end{equation*}
in which the left arrow maps $(z,\tau)$ to $(Mz,\tau)$ and on ring sheaves is given by
\[f(z,\tau)+g(z,\tau)\cdot \eta^\vee+h(z,\tau)\cdot \omega^\vee \mapsto f(Mz,\tau)+g(Mz,\tau)\cdot \eta^\vee+h(Mz,\tau)\cdot \omega^\vee.\]
If we pull back the left side of $(3.5.16)$ along $k_1$ and use the first diagram together with the identification $(3.5.17)$ we obtain $k_1^*\widetilde{\mathcal P}$. If we pull back the right side of $(3.5.16)$ along $k_1$ and use the second diagram we obtain (as $\mathcal O_{\C\times\H}$-module) $\widetilde{[M]}^*k_1^*\widetilde{\mathcal P}$.\\
Hence, pulling back $(3.5.16)$ along $k_1$ yields a $\mathcal O_{\C\times \H}$-linear isomorphism
\[k_1^*\widetilde{\mathcal P} \simeq \widetilde{[M]}^*k_1^*\widetilde{\mathcal P}.\]
The point is that this is precisely the identification of $(3.5.14)$, which together with the assertion of Lemma 3.5.10 and a careful inspection implies the desired claims.\\
That the last isomorphism indeed coincides with $(3.5.14)$ basically follows by recalling from 2.5.2 how in terms of the Poincaré bundle the invariance isomorphism
\[\mathcal L_1 \simeq [M]^*\mathcal L_1\]
comes from the canonical isomorphism $(3.5.15)$ together with the identification
\[\tag{\textbf{3.5.18}} \mathcal O_{E\times_S \widehat{E}^\natural_1}\xrightarrow{\sim} (\mathrm{id}_E\times [M]^\natural_1)_* \mathcal O_{E\times_S \widehat{E}^\natural_1}=\mathcal O_{E\times_S \widehat{E}^\natural_1}\]
induced by the morphism
\[\id_E\times [M]^\natural_1: E \times_S \widehat{E}^\natural_1 \rightarrow E\times_S \widehat{E}^\natural_1.\]
Note that the map $\mathcal O_E \oplus \mathcal H_E \rightarrow \mathcal O_E \oplus \mathcal H_E$ associated with $(3.5.18)$ is given by $\id\oplus (\cdot M)$.\\
Cf. $(2.5.2)$, $(2.5.7)$-$(2.5.9)$, Prop. 2.5.1 and Prop. 2.5.2 for the details.\\
\newline
There is hence only one task left:\\
\newline
\textit{Proof of Lemma 3.5.10:}
Consider the canonical isomorphism of $\widehat{E}$-rigidified line bundles
\[(\mathrm{id}_E\times [M]^\wedge)^*(\mathcal P^0,r^0) \simeq (([M]\times \mathrm{id}_{\widehat{E}})^*\mathcal P^0, r^0_{[M]}),\]
constructed as $(2.5.2)$ by using the dual abelian scheme $\widehat{E}$ and $(\mathcal P^0,r^0)$ instead of $\widehat{E}^\natural$ and $(\mathcal P,r,\nabla_\mathcal P)$; note that by Prop. 2.5.1 the transpose endomorphism of $[M]$ is the $M$-multiplication map $[M]^\wedge$ on $\widehat{E}$.\\
With the principal polarization of $(3.1.3)$ this reads as isomorphism of $E$-rigidified line bundles
\[\tag{\textbf{3.5.19}} (\mathrm{id}_E\times [M])^*(\mathcal P^0,r^0) \simeq (([M]\times \id_E)^*\mathcal P^0, r^0_{[M]}).\]
Analytifying $(3.5.19)$ (and restricting to the fixed component) gives after pullback to $\C\times \C\times \H$ an induced isomorphism
\[\tag{\textbf{3.5.20}} (\widetilde{\mathrm{id}_E \times [M]})^*\widetilde{\mathcal P^0} \simeq (\widetilde{[M]\times \id_E})^*\widetilde{\mathcal P^0},\]
where $\widetilde{\mathcal P^0}$ denotes the pullback of $(\mathcal P^0)^{an}$ to $\C \times \C \times \H$ and where we set
\[\widetilde{\mathrm{id}_E \times [M]}: \C \times \C \times \H \rightarrow \C\times \C\times \H, \quad (z,w,\tau)\mapsto (z,Mw,\tau)\]
resp.
\[\widetilde{[M]\times \id_E}: \C \times \C \times \H \rightarrow \C\times \C\times \H, \quad (z,w,\tau)\mapsto (Mz,w,\tau).\]
Note that the trivializing section $t$ of $\widetilde{\mathcal P}$ is induced by pulling back the trivializing section
\[t^0=\frac{1}{J(z,w,\tau)}\otimes(\omega^0_{\mathrm{can}})^\vee\]
of $\widetilde{\mathcal P^0}$ along
\[\C\times \C^2\times \H \rightarrow \C\times \C\times \H, \quad (z,w,u,\tau)\mapsto (z,-w,\tau).\]
Cf. $(3.4.10)$, $(3.4.13)$ and $(3.4.15)$.\\
\newline
The goal for the rest of the proof is to show that $(3.5.20)$ maps
\[(\widetilde{\mathrm{id}_E\times [M]})^*t^0=\frac{1}{J(z,Mw,\tau)}\otimes(\omega^0_{\mathrm{can}})^\vee \mapsto \frac{1}{J(Mz,w,\tau)}\otimes(\omega^0_{\mathrm{can}})^\vee=(\widetilde{[M]\times \id_E})^*t^0;\]
because of compatibilities which are readily worked out this suffices to deduce the claim of the lemma.\\
\newline
Recall from $(3.1.5)$ that $(\mathcal P^0,r^0)$ writes as
\[(\mathcal M \otimes_{\mathcal O_{E\times_SE}} (\pi \times \pi)^*\epsilon^* \mathcal O_E([0]), \mathrm{can}),\]
where
\[\mathcal M= \mu^* \mathcal O_E([0]) \otimes_{\mathcal O_{E\times_SE}} \mathrm{pr}_1^* \mathcal O_E([0])^{-1} \otimes_{\mathcal O_{E\times_SE}} \mathrm{pr}_2^*\mathcal O_E([0])^{-1} \simeq \mathcal O_{E\times_S E}(\bar{\Delta}_E-[0]\times E - E\times [0])\]
and where $\mathrm{can}$ is the canonical rigidification along the inclusion of the second factor of $E\times_S E$.\\
Taking this into account, $(3.5.19)$ induces an isomorphism of line bundles on $E\times_S E$:

\[\tag{\textbf{3.5.21}} (\mathrm{id}_E\times [M])^*\mathcal O_{E\times_S E}(\bar{\Delta}_E-[0]\times E - E\times [0]) \simeq ([M]\times \id_E)^*\mathcal O_{E\times_S E}(\bar{\Delta}_E-[0]\times E - E\times [0])\]
whose restriction via $(\epsilon \times \epsilon): S \rightarrow E\times_S E$ is the identity on
\[\epsilon^* \mathcal O_E([0]) \otimes_{\mathcal O_S} \epsilon^* \mathcal O_E([0])^{-1} \otimes_{\mathcal O_S} \epsilon^* \mathcal O_E([0])^{-1};\]
note that there can only be one isomorphism of the form $(3.5.21)$ having this property: this is an easy consequence of Lemma 0.1.5.\\
Writing $\mathfrak D$ for the divisor $\bar{\Delta}_E-[0]\times E - E\times [0]$, the analytification of $(3.5.21)$ is an isomorphism
\[\tag{\textbf{3.5.22}} (\mathrm{id}_{E^{an}}\times [M]^{an})^*\mathcal O_{E^{an}\times_{S^{an}} E^{an}}(\mathfrak D^{an}) \simeq ([M]^{an}\times \mathrm{id}_{E^{an}})^*\mathcal O_{E^{an}\times_{S^{an}}E^{an}}(\mathfrak D^{an})\]
with the analogous restriction property along $(\epsilon^{an}\times \epsilon^{an}): S^{an}\rightarrow E^{an}\times_{S^{an}}E^{an}$, and there can only be one isomorphism of the form $(3.5.22)$ having this property: note that the canonical map $\mathcal O_{S^{an}}\rightarrow (\pi^{an}\times\pi^{an})_*\mathcal O_{E^{an}\times_{S^{an}}E^{an}}$ is an isomorphism\footnote{One can see this e.g. by taking the analytification of the respective (well-known) algebraic isomorphism and by noting that the appearing base change is an isomorphism by \cite{SGA1}, exp. XII, Thm. 4.2. In fact, the canonical morphism on structure sheaves is an isomorphism whenever one considers a proper smooth map of analytic spaces with connected fibers.}, such that one has the analytic analogue of Lemma 0.1.5 available (cf. its proof), from which the claimed uniqueness again easily follows.\\
Now comes the crucial point: the function
\[\tag{\textbf{3.5.23}} (z,w,\tau) \mapsto \frac{J(Mz,w,\tau)}{J(z,Mw,\tau)}\]
is invariant under the effect of $\Z^2 \times \Z^2 \times \Gamma(N)$ (use Cor. 3.3.14), and considering it hence as a meromorphic function on $E^{an}\times_{S^{an}}E^{an}$ its divisor is given by
\[([M]^{an}\times \mathrm{id}_{E^{an}})^*\mathfrak D^{an}-(\mathrm{id}_{E^{an}}\times [M]^{an})^*\mathfrak D^{an}.\]
The function $(3.5.23)$ thus defines an isomorphism between the line bundles in $(3.5.22)$ which is moreover checked to have the restriction property along $(\epsilon^{an}\times \epsilon^{an})$ of above, such that the mentioned uniqueness implies its equality with the isomorphism of $(3.5.22)$.\\
This implies straightforwardly the desired claim about $(3.5.20)$ and hence the statement of the lemma. \qquad \qed\\
\newline
As this finally finishes the proof of Prop. 3.5.9 in all relevant details we may return to the problem of finding explicit formulas for the translation of sections of the logarithm sheaves.

\subsubsection{Explicit formulas II}
With the knowledge of the case $n=1$ in Prop. 3.5.9 it is easy to compute the respective formula for the case of arbitrary $n\geq 0$, where the usual trivialization of $\mathrm{pr}^*\mathcal L_n^{an}$ in Def. 3.5.4 is fixed.\\
The result of the calculation is:
\begin{corollary}
Assume that a section
\[\begin{pmatrix} l_{0,0}(z,\tau) \\ l_{1,0}(z,\tau) \\ l_{2,0}(z,\tau)\\ \vdots\\ \vdots\\l_{0,n}(z,\tau) \end{pmatrix}\]
of $\mathcal L_n^{an}$ over an open subset $V\subseteq E^{an}$ is given. By pullback along $T_{a,b}^{an}$ and using $(3.5.12)$ we obtain a section of $\mathcal L_n^{an}$ over $(T^{an}_{a,b})^{-1}(V)$. It expresses as
\[\begin{pmatrix} \widehat{l}_{0,0}(z,\tau) \\ \widehat{l}_{1,0}(z,\tau) \\ \widehat{l}_{2,0}(z,\tau)\\ \vdots\\ \vdots\\\widehat{l}_{0,n}(z,\tau)\end{pmatrix},\]
where for each $0\leq i,j$ with $i+j\leq n$ the function $\widehat{l}_{i,j}(z,\tau)$ is given by
\[\widehat{l}_{i,j}(z,\tau)=\sum_{k=0}^i\frac{(-2\pi i \frac{aj_0}{N})^{i-k}}{(i-k)!}\cdot l_{k,j}\Big(z+\frac{aj_0\tau}{N}+\frac{b}{N},\tau\Big).\]
Equivalently, the function $\widehat{l}_{i,j}(z,\tau)$ is determined by the following equation, with free variables $w,u$ and to hold modulo expressions with factors $w^iu^j, \ i+j>n$:
\begin{align*}
&\widehat{l}_{0,0}(z,\tau)+\widehat{l}_{1,0}(z,\tau)w+\widehat{l}_{2,0}(z,\tau)w^2 +......+\widehat{l}_{0,n}(z,\tau)u^n=\mathrm{e}^{-2\pi i \frac{aj_0}{N}w}\cdot \bigg[l_{0,0}\Big(z+\frac{aj_0\tau}{N}+\frac{b}{N},\tau\Big)\\
&+l_{1,0}\Big(z+\frac{aj_0\tau}{N}+\frac{b}{N},\tau\Big)w+l_{2,0}\Big(z+\frac{aj_0\tau}{N}+\frac{b}{N},\tau\Big)w^2 +......+l_{0,n}\Big(z+\frac{aj_0\tau}{N}+\frac{b}{N},\tau\Big)u^n\bigg].
\end{align*}

\qquad \qed
\end{corollary}

\begin{remark}
To prevent possible confusion let us point to the following:\\
The integers $a$ and $b$ could very well be both divisible by $N$, such that $T_{a,b}=\id$ and $(3.5.12)$ coincides with the identity map on $\mathcal L_n^{an}$. The two vectors in Cor. 3.5.11 thus need to be equal, which at first sight they don't seem to be. But one should always recall that the functions $l_{i,j}(z,\tau)$, defining a section of $\mathcal L_n^{an}$, have a transformation property under the effect of $\Z^2\times \Gamma(N)$, namely by the matrix $A_n$ in Prop. 3.5.6. Taking this into account the two vectors indeed become equal.
\end{remark}
We can now explicitly describe the specialization of sections of $\mathcal L_n^{an}$ along $t^{an}_{a,b}$.\\
The trivialization for the pullback of $\prod_{k=0}^n \mathrm{Sym}^k_{\mathcal O_{S^{an}}}\mathcal H^{an}$ to the universal covering $\H$ is as fixed in the preliminary remarks preceding 3.5.1.

\begin{corollary}
Let
\[\begin{pmatrix} l_{0,0}(z,\tau) \\ l_{1,0}(z,\tau) \\ l_{2,0}(z,\tau)\\ \vdots\\ \vdots\\l_{0,n}(z,\tau)\end{pmatrix}\]
be a section of $\mathcal L_n^{an}$ which is defined on some open subset of $E^{an}$ intersecting $t_{a,b}^{an}(S^{an})$.\\
It induces a section of $\prod_{k=0}^n \mathrm{Sym}^k_{\mathcal O_{S^{an}}}\mathcal H^{an}$, obtained by its pullback along $t_{a,b}^{an}$ and using the isomorphism $(3.5.13)$ together with the splitting $\varphi_n^{an}$:
\[(t_{a,b}^{an})^*\mathcal L_n^{an}\simeq (\epsilon^{an})^*\mathcal L_n^{an} \simeq \prod_{k=0}^n \mathrm{Sym}^k_{\mathcal O_{S^{an}}}\mathcal H^{an}.\]
This induced section then expresses as
\[\begin{pmatrix} \widehat{l}_{0,0}(0,\tau) \\ \widehat{l}_{1,0}(0,\tau) \\ \widehat{l}_{2,0}(0,\tau)\\ \vdots\\ \vdots\\\widehat{l}_{0,n}(0,\tau)\end{pmatrix},\]
where for each $0\leq i,j$ with $i+j\leq n$ the function $\widehat{l}_{i,j}(0,\tau)$ is given by
\[\widehat{l}_{i,j}(0,\tau)=\sum_{k=0}^i\frac{(-2\pi i \frac{aj_0}{N})^{i-k}}{(i-k)!}\cdot l_{k,j}\Big(\frac{aj_0\tau}{N}+\frac{b}{N},\tau \Big).\]
\end{corollary}
\begin{proof}
Recall how $(3.5.13)$ is obtained from $(3.5.12)$. Then combine Cor. 3.5.11 and Prop. 3.5.8.
\end{proof}
In the end we want to specialize sections of $\Omega^1_{E^{an}}\otimes_{\mathcal O_{E^{an}}} \mathcal L_n^{an}$ along $t_{a,b}^{an}$. In order to reasonably formulate the respective result we need the following insertion.

\subsubsection{Conventions for the trivialization of differential modules}
We fix trivializations on the universal covering for modules of differential forms; they will tacitly remain valid until the end of the work.\\
For the various occurring manifolds we keep working on the component associated to the integer $j_0$ (where $(j_0,N)=1$) and thereby perform the usual abuse of notation, leaving away this index.\\
Here and henceforth, let us denote by $\Omega^i_{E^{an}}$ resp. $\Omega^i_{S^{an}}$ the sheaf of (absolute) differential $i$-forms on the complex manifolds $E^{an}$ resp. $S^{an}$.\\
\newline
$\bullet$ The pullback of $\Omega^1_{E^{an}}$ via
\[\mathrm{pr}: \C\times \H \rightarrow E^{an}\]
shall be trivialized by the ordered basic sections $\{\mathrm{d}z,\mathrm{d}\tau\}$.\\
The associated automorphy matrix for $\Omega^1_{E^{an}}$ is readily computed as
\[\tag{\textbf{3.5.24}} \begin{split}
&\Z^2 \times \Gamma(N) \times \C \times \H \rightarrow \mathrm{GL}_2(\C)\\
&\Bigg( \begin{pmatrix}  m \\ n \end{pmatrix}, \begin{pmatrix} a & b \\ c & d \end{pmatrix}, (z,\tau) \Bigg) \mapsto \begin{pmatrix}
                        c\tau+d & 0\\
											 (cz+cn-dm)(c\tau+d) &  (c\tau+d)^2\end{pmatrix}.
\end{split}\]

$\bullet$ The pullback of $\Omega^2_{E^{an}}$ to $\C\times \H$ will be trivialized by $\{\mathrm{d}z\wedge \mathrm{d}\tau\}$, which gives the factor of automorphy
\[\tag{\textbf{3.5.25}} \Z^2\times \Gamma(N) \times \C\times \H \rightarrow \C^*, \quad \Bigg( \begin{pmatrix}  m \\ n \end{pmatrix}, \begin{pmatrix} a & b \\ c & d \end{pmatrix}, (z,\tau)\Bigg)\mapsto (c\tau+d)^3.\]

The exterior derivation on $1$-forms then expresses as
\[\tag{\textbf{3.5.26}} \mathrm{d}: \Omega^1_{E^{an}}  \rightarrow \Omega^2_{E^{an}}, \quad \begin{pmatrix}f(z,\tau) \\ g(z,\tau) \end{pmatrix} \mapsto \partial_zg(z,\tau)-\partial_{\tau}f(z,\tau).\]

$\bullet$ The pullback of $\Omega^1_{S^{an}}$ via
\[p:\H\rightarrow S^{an}\]
shall be trivialized by $\{\mathrm{d}\tau\}$. The induced factor of automorphy is
\[\tag{\textbf{3.5.27}} \Gamma(N) \times \H \rightarrow \C^*, \quad \Bigg(\begin{pmatrix} a & b \\ c & d \end{pmatrix}, \tau\Bigg)\mapsto (c\tau+d)^2.\]
Assume then that we have a section
\[\begin{pmatrix}f(z,\tau) \\ g(z,\tau) \end{pmatrix}\]
of $\Omega^1_{E^{an}}$ which is defined on some open subset of $E^{an}$ intersecting $t_{a,b}^{an}(S^{an})$.\\
By pullback via $t_{a,b}^{an}$ and usage of the canonical map
\[(t_{a,b}^{an})^*\Omega^1_{E^{an}}\rightarrow \Omega^1_{S^{an}}\]
we obtain a section of $\Omega^1_{S^{an}}$ which expresses as follows:
\[\tag{\textbf{3.5.28}} \frac{aj_0}{N}\cdot f\Big(\frac{aj_0\tau}{N}+\frac{b}{N},\tau\Big)+g\Big(\frac{aj_0\tau}{N}+\frac{b}{N},\tau \Big).\]

$\bullet$ As already employed several times, the module $\Omega^1_{E^{an}/S^{an}}$ of relative differential forms is trivialized on $\C\times \H$ by $\{\mathrm{d}z\}$ and receives the associated factor of automorphy
\[\tag{\textbf{3.5.29}} \Z^2\times \Gamma(N)\times \C\times \H \rightarrow \C^*, \quad \Bigg( \begin{pmatrix}  m \\ n \end{pmatrix}, \begin{pmatrix} a & b \\ c & d \end{pmatrix}, (z,\tau)\Bigg)\mapsto c\tau+d.\]
The canonical map from absolute to relative differential $1$-forms then expresses as
\[\tag{\textbf{3.5.30}} \Omega^1_{E^{an}}\rightarrow \Omega^1_{E^{an}/S^{an}}, \quad \begin{pmatrix}f(z,\tau) \\ g(z,\tau) \end{pmatrix} \mapsto f(z,\tau).\]

\subsubsection{Explicit formulas III}
With the above conventions we can state the final result of this subsection, answering the question how to compute explicitly the pullback of sections of $\Omega^1_{E^{an}}\otimes_{\mathcal O_{E^{an}}} \mathcal L_n^{an}$ along analytified torsion sections. It will play a central auxiliary role later when we determine the specialization of the $D$-variant.\\
Its statement is easily deduced from Cor. 3.5.13 and $(3.5.28)$.
\begin{proposition}
Let
\[\begin{pmatrix} l_{0,0}(z,\tau) \\ l_{1,0}(z,\tau) \\ l_{2,0}(z,\tau)\\ \vdots\\ \vdots\\l_{0,n}(z,\tau)  \\  \lambda_{0,0}(z,\tau) \\ \lambda_{1,0}(z,\tau) \\ \lambda_{2,0}(z,\tau)\\ \vdots\\ \vdots\\ \lambda_{0,n}(z,\tau) \end{pmatrix}\]
be a section of $\Omega^1_{E^{an}}\otimes_{\mathcal O_{E^{an}}} \mathcal L_n^{an}$ which is defined on some open subset of $E^{an}$ intersecting $t_{a,b}^{an}(S^{an})$.\\
It induces a section of $\Omega^1_{S^{an}}\otimes_{\mathcal O_{S^{an}}} \prod_{k=0}^n \mathrm{Sym}^k_{\mathcal O_{S^{an}}}\mathcal H^{an}$, obtained from its pullback along $t_{a,b}^{an}$ and using the composition
\[(t_{a,b}^{an})^*(\Omega^1_{E^{an}}\otimes_{\mathcal O_{E^{an}}} \mathcal L_n^{an}) \simeq (t_{a,b}^{an})^*\Omega^1_{E^{an}}\otimes_{\mathcal O_{S^{an}}} \prod_{k=0}^n \mathrm{Sym}^k_{\mathcal O_{S^{an}}}\mathcal H^{an}  \xrightarrow{\mathrm{can}} \Omega^1_{S^{an}} \otimes_{\mathcal O_{S^{an}}} \prod_{k=0}^n \mathrm{Sym}^k_{\mathcal O_{S^{an}}}\mathcal H^{an},\]
where for the first map cf. Cor. 3.5.13.\\
This induced section then expresses as
\[\begin{pmatrix} \frac{aj_0}{N}\cdot \widehat{l}_{0,0}(0,\tau) + \widehat{\lambda}_{0,0}(0,\tau)\\ \frac{aj_0}{N}\cdot \widehat{l}_{1,0}(0,\tau) +\widehat{\lambda}_{1,0}(0,\tau) \\ \frac{aj_0}{N}\cdot \widehat{l}_{2,0}(0,\tau)+\widehat{\lambda}_{2,0}(0,\tau)\\ \vdots\\ \vdots\\ \frac{aj_0}{N}\cdot \widehat{l}_{0,n}(0,\tau)+\widehat{\lambda}_{0,n}(0,\tau)\end{pmatrix},\]
where for each $0\leq i,j$ with $i+j\leq n$ the function $\widehat{l}_{i,j}(0,\tau)$ resp. $\widehat{\lambda}_{i,j}(0,\tau)$ is given by
\[\widehat{l}_{i,j}(0,\tau)=\sum_{k=0}^i\frac{(-2\pi i \frac{aj_0}{N})^{i-k}}{(i-k)!}\cdot l_{k,j}\Big(\frac{aj_0\tau}{N}+\frac{b}{N},\tau \Big)\]
resp.
\[\widehat{\lambda}_{i,j}(0,\tau)=\sum_{k=0}^i\frac{(-2\pi i \frac{aj_0}{N})^{i-k}}{(i-k)!}\cdot \lambda_{k,j}\Big(\frac{aj_0\tau}{N}+\frac{b}{N},\tau \Big).\]
\qquad \qed
\end{proposition}

\subsection{The analytification of the absolute connection}
\subsubsection{Precise formulation of the problem}
The line bundle with integrable $\widehat{E}^\natural_1$-connection $(\mathcal P_1,\nabla_{\mathcal P_1})$ on $E\times_S \widehat{E}^\natural_1$ becomes, when considered as $\mathcal O_E$-module, a vector bundle $\mathcal L_1$ on $E$ with integrable $S$-connection $\nabla_1^{res}$. Essentially via the trivialization $s$ of $(\mathcal P,\nabla_\mathcal P)$ along $E$ it sits in an exact sequence of $\mathcal D_{E/S}$-modules:
\[\tag{\textbf{3.5.31}}0\rightarrow \mathcal H_E \rightarrow \mathcal L_1 \rightarrow \mathcal O_E \rightarrow 0.\]
The extension class of this sequence maps to the identity under the lower projection in $(2.1.3)$:
\[\Ext^1_{\mathcal D_{E/S}}(\mathcal O_E, \mathcal H_E)  \rightarrow \Hom_{\mathcal O_S}(\mathcal O_S, \mathcal H^\vee \otimes_{\mathcal O_S} \mathcal H).\]
Essentially via the rigidification $r$ of $\mathcal P$ along $\widehat{E}^\natural$ and the canonical $\mathcal O_S$-isomorphism $\mathcal O_S \oplus \mathcal H \simeq \mathcal O_{\widehat{E}^\natural_1}$ one constructs a $\mathcal O_S$-linear splitting $\varphi_1: \mathcal O_S \oplus \mathcal H \simeq \epsilon^*\mathcal L_1$ for the pullback of $(3.5.31)$ along $\epsilon$.\\
Cf. 2.3.1 and Thm. 2.3.1 for these constructions and facts.\\
\newline
Prop. 2.1.4 resp. Cor. 2.3.2 tells us that there is a unique prolongation of $\nabla_1^{res}$ to an integrable $\Q$-connection $\nabla_1$ that is characterized by the following property: If $\mathcal L_1$ is equipped with $\nabla_1$, then:
\[\tag{\textbf{3.5.32}} \begin{split}
&\textrm{(i) The $\mathcal D_{E/S}$-linear exact sequence $(3.5.31)$ becomes $\mathcal D_{E/\Q}$-linear.}\\
&\textrm{(ii) The $\mathcal O_S$-linear splitting $\varphi_1$ becomes $\mathcal D_{S/\Q}$-linear.}
\end{split}\]
The extension class of (i) then maps to the identity under the lower projection in $(2.1.3)$:
\[\Ext^1_{\mathcal D_{E/\Q}}(\mathcal O_E, \mathcal H_E)  \rightarrow \Hom_{\mathcal D_{S/\Q}}(\mathcal O_S, \mathcal H^\vee \otimes_{\mathcal O_S} \mathcal H).\]
Altogether, this means that $(\mathcal L_1,\nabla_1,\varphi_1)$ defines the first logarithm sheaf of $E/S/\Q$ in the sense of 1.1 (and that hence our notation is justified).\\
\newline
Analytically, with respect to our fixed (componentwise) trivializations on the universal coverings, we have already gained explicit knowledge about the following of the above data:\\
\newline
$\bullet$ The analytification
\[0\rightarrow \mathcal H_{E^{an}}^{an} \rightarrow \mathcal L_1^{an} \rightarrow \mathcal O_{E^{an}} \rightarrow 0\]
of $(3.5.31)$ expresses on sections as
\[\tag{\textbf{3.5.33}} \begin{pmatrix} \chi(z,\tau) \\ \xi(z,\tau) \end{pmatrix}  \mapsto \begin{pmatrix} 0 \\ \chi(z,\tau) \\ \xi(z,\tau) \end{pmatrix}, \quad \begin{pmatrix} \lambda (z,\tau) \\ \mu(z,\tau) \\ \nu(z,\tau) \end{pmatrix}  \mapsto \lambda(z,\tau),\]
as follows directly from the remarks preceding Prop. 3.5.2.\\
\newline
$\bullet$ The analytified splitting
\[\varphi_1^{an}: \mathcal O_{S^{an}} \oplus \mathcal H^{an} \simeq (\epsilon^{an})^*\mathcal L_1^{an}\]
expresses as the identity on sections written as vectors of holomorphic functions in $\tau$.\\
The pullback via $\epsilon^{an}$ of a section
\[\begin{pmatrix} \lambda (z,\tau) \\ \mu(z,\tau) \\ \nu(z,\tau) \end{pmatrix}\]
of $\mathcal L_1^{an}$, defined on some open subset of $E^{an}$ intersecting $\epsilon^{an}(S^{an})$, maps under $\varphi_1^{an}$ to the section
\[\begin{pmatrix} \lambda (0,\tau) \\ \mu(0,\tau) \\ \nu(0,\tau) \end{pmatrix}\]
of $\mathcal O_{S^{an}} \oplus \mathcal H^{an}$. Cf. Prop. 3.5.8 for these facts.\\
\newline
$\bullet$ Finally, recall from Prop. 3.5.3 that the formula for the analytification
\[(\nabla_1^{res})^{an}: \mathcal L_1^{an}\rightarrow \Omega^1_{E^{an}/S^{an}}\otimes_{\mathcal O_{E^{an}}} \mathcal L_1^{an}\]
of the relative connection
\[\nabla_1^{res}: \mathcal L_1\rightarrow \Omega^1_{E/S}\otimes_{\mathcal O_E} \mathcal L_1\]
is given by
\[\begin{pmatrix} \lambda(z,\tau) \\ \mu(z,\tau) \\ \nu(z,\tau) \end{pmatrix}\mapsto \begin{pmatrix} \partial_z \lambda(z,\tau) \\ \partial_z\mu(z,\tau)+\eta(1,\tau)\cdot \lambda(z,\tau) \\ \partial_z\nu(z,\tau)+\lambda(z,\tau) \end{pmatrix}.\]\\
\vspace{1mm}
The goal of the present subsection is to determine the formula for the analytification
\[\nabla_1^{an}: \mathcal L_1^{an}\rightarrow \Omega^1_{E^{an}}\otimes_{\mathcal O_{E^{an}}} \mathcal L_1^{an}\]
of the absolute connection
\[\nabla_1: \mathcal L_1\rightarrow \Omega^1_{E/\Q}\otimes_{\mathcal O_E} \mathcal L_1.\]
To possess such an explicit knowledge of $\nabla_1^{an}$ is not only interesting in its own right but also essential if one wants to give an analytic description of the $D$-variant polylogarithm system.\\
Let us hence proceed to develop the solution for this important problem. The already known analytic expressions of the other data, which we have repeated above, together with an analytic version for the characterizing property of $\nabla_1$ in $(3.5.32)$ will play the decisive auxiliary role.
\subsubsection{Description of the Gauß-Manin connection}
We continue to work on some fixed component of $E^{an}/S^{an}$ and to always consider the restrictions of all occurring objects to this component, omitting however the respective index in our notation.\\
\newline
It was already mentioned at the beginning of 3.5 that we freely identify
\[H^1_{\mathrm{dR}}(E/S)^{an}\simeq H^1_{\mathrm{dR}}(E^{an}/S^{an}) \quad \textrm{resp.} \quad \mathcal H^{an}\simeq H^1_{\mathrm{dR}}(E^{an}/S^{an})^\vee.\]
As explained in $(3.4.6)$ and $(3.4.7)$ we trivialize the pullback of these $\mathcal O_{S^{an}}$-vector bundles to the universal covering $\H$ of $S^{an}$ by the ordered basic sections $\{\eta,\omega\}=\{p(z,\tau)\mathrm{d}z,\mathrm{d}z\}$ resp. by $\{\eta^\vee, \omega^\vee\}$. Moreover, recall from $(3.5.27)$ that the pullback of $\Omega^1_{S^{an}}$ to $\H$ is trivialized by $\{\mathrm{d}\tau\}$.\\
\newline
The analytified Gauß-Manin connection on $\mathrm{H^1_{dR}}(E/S)$ expresses as the integrable connection
\[H^1_{\mathrm{dR}}(E^{an}/S^{an}) \rightarrow \Omega^1_{S^{an}}\otimes_{\mathcal O_{S^{an}}} H^1_{\mathrm{dR}}(E^{an}/S^{an})\]
which acts (with respect to our fixed trivializations) on sections over open subsets of $S^{an}$ as follows:
\[\begin{pmatrix} \chi(\tau) \\ \xi(\tau) \end{pmatrix} \mapsto \begin{pmatrix} \partial_{\tau}\chi(\tau)+\frac{1}{2\pi i}\eta(1,\tau)\cdot \chi(\tau)+\frac{1}{2\pi i } \xi(\tau)\\ \partial_{\tau}\xi(\tau)+\partial_{\tau}\eta(1,\tau)\cdot \chi(\tau)-\frac{1}{2\pi i }\eta(1,\tau)^2 \cdot \chi(\tau)-\frac{1}{2\pi i }\eta(1,\tau) \cdot \xi(\tau)\end{pmatrix}.\]
The formula is deduced from \cite{Kat4}, A 1.3.8.\\
Its dual connection
\[\mathcal H^{an} \rightarrow \Omega^1_{S^{an}}\otimes_{\mathcal O_{S^{an}}} \mathcal H^{an}\]
is then routinely computed as
\[\tag{\textbf{3.5.34}}  \begin{pmatrix} \chi(\tau) \\ \xi(\tau) \end{pmatrix} \mapsto \begin{pmatrix} \partial_{\tau}\chi(\tau)-\frac{1}{2\pi i}\eta(1,\tau)\cdot \chi(\tau)+\frac{1}{2\pi i}\eta(1,\tau)^2\cdot \xi(\tau)-\partial_{\tau}\eta(1,\tau)\cdot \xi(\tau)
\\ \partial_{\tau}\xi(\tau)-\frac{1}{2\pi i }\chi(\tau)+\frac{1}{2\pi i }\eta(1,\tau) \cdot \xi(\tau) \end{pmatrix}.\]

\subsubsection{The guess for the absolute connection}

Now let
\[\begin{pmatrix} \lambda (z,\tau) \\ \mu(z,\tau) \\ \nu(z,\tau) \end{pmatrix}\]
be a local section of $\mathcal L_1^{an}$ defined on some open subset of $E^{an}$ intersecting $\epsilon^{an}(S^{an})$. The image of its pullback along $\epsilon^{an}$ under the isomorphism $\varphi_1^{an}$ is the section
\[\begin{pmatrix} \lambda (0,\tau) \\ \mu(0,\tau) \\ \nu(0,\tau) \end{pmatrix}\]
of $\mathcal O_{S^{an}}\oplus \mathcal H^{an}$, and by $(3.5.34)$ the connection of $\mathcal O_{S^{an}} \oplus \mathcal H^{an}$ maps it to
\[\tag{\textbf{3.5.35}} \begin{pmatrix} \partial_{\tau}\lambda (0,\tau) \\ \partial_{\tau}\mu(0,\tau)-\frac{1}{2\pi i}\eta(1,\tau)\cdot \mu(0,\tau)+\frac{1}{2\pi i}\eta(1,\tau)^2\cdot \nu(0,\tau)-\partial_{\tau}\eta(1,\tau)\cdot \nu(0,\tau) \\ \partial_{\tau}\nu(0,\tau)-\frac{1}{2\pi i }\mu(0,\tau)+\frac{1}{2\pi i }\eta(1,\tau) \cdot \nu(0,\tau)  \end{pmatrix},\]
which is a section of $\Omega^1_{S^{an}}\oplus (\Omega^1_{S^{an}}\otimes_{\mathcal O_{S^{an}}} \mathcal H^{an})$.\\
On the other hand, the connection $\nabla_1^{an}$ applied to
\[\begin{pmatrix} \lambda (z,\tau) \\ \mu(z,\tau) \\ \nu(z,\tau) \end{pmatrix}\]
is a certain section of $\Omega^1_{E^{an}}\otimes_{\mathcal O_{E^{an}}} \mathcal L_1^{an}$, i.e. (recalling the agreements for the respective trivializations) given by a vector of functions
\[\nabla_1^{an}\Bigg(\begin{pmatrix} \lambda (z,\tau) \\ \mu(z,\tau) \\ \nu(z,\tau) \end{pmatrix} \Bigg) =\begin{pmatrix} \varrho_1 (z,\tau) \\ \varrho_2(z,\tau) \\ \varrho_3(z,\tau) \\ \varrho_4(z,\tau) \\ \varrho_5(z,\tau) \\ \varrho_6(z,\tau)\end{pmatrix}.\]
Our knowledge of the formula for $(\nabla_1^{res})^{an}$ implies that
\[\tag{\textbf{3.5.36}} \begin{pmatrix} \varrho_1(z,\tau) \\ \varrho_2(z,\tau) \\ \varrho_3(z,\tau) \end{pmatrix}= \begin{pmatrix} \partial_z \lambda(z,\tau) \\ \partial_z\mu(z,\tau)+\eta(1,\tau)\cdot \lambda(z,\tau) \\ \partial_z\nu(z,\tau)+\lambda(z,\tau) \end{pmatrix}.\]

Now consider the image of
\[(\epsilon^{an})^*\nabla_1^{an} \Bigg(\begin{pmatrix} \lambda (z,\tau) \\ \mu(z,\tau) \\ \nu(z,\tau) \end{pmatrix}\Bigg)\]
under the arrow
\[(\epsilon^{an})^*\Omega^1_{E^{an}}\otimes_{\mathcal O_{S^{an}}} (\epsilon^{an})^*\mathcal L_1^{an} \xrightarrow{\mathrm{can}} \Omega^1_{S^{an}}\otimes_{\mathcal O_{S^{an}}} (\epsilon^{an})^*\mathcal L_1^{an} \simeq \Omega^1_{S^{an}}\oplus (\Omega^1_{S^{an}}\otimes_{\mathcal O_{S^{an}}} \mathcal H^{an}),\]
where the isomorphism is induced by $\varphi_1^{an}$. It is given by
\[\begin{pmatrix}\varrho_4(0,\tau)\\ \varrho_5(0,\tau) \\ \varrho_6(0,\tau)  \end{pmatrix},\]
and if $\varphi_1^{an}: (\epsilon^{an})^*\mathcal L_1^{an} \simeq \mathcal O_{S^{an}}\oplus \mathcal H^{an}$ shall be horizontal, then this must equal $(3.5.35)$:
\[\begin{pmatrix}\varrho_4(0,\tau)\\ \varrho_5(0,\tau) \\ \varrho_6(0,\tau)  \end{pmatrix}=\begin{pmatrix} \partial_{\tau}\lambda (0,\tau) \\ \partial_{\tau}\mu(0,\tau)-\frac{1}{2\pi i}\eta(1,\tau)\cdot \mu(0,\tau)+\frac{1}{2\pi i}\eta(1,\tau)^2\cdot \nu(0,\tau)-\partial_{\tau}\eta(1,\tau)\cdot \nu(0,\tau) \\ \partial_{\tau}\nu(0,\tau)-\frac{1}{2\pi i }\mu(0,\tau)+\frac{1}{2\pi i }\eta(1,\tau) \cdot \nu(0,\tau)  \end{pmatrix}.\]
From the last equation and from $(3.5.36)$ we may arrive at the following
\begin{guess}
The absolute connection $\nabla_1^{an}: \mathcal L_1^{an} \rightarrow \Omega^1_{E^{an}}\otimes_{\mathcal O_{E^{an}}} \mathcal L_1^{an}$ is given on local sections by the formula
\[\begin{pmatrix} \lambda (z,\tau) \\ \mu(z,\tau) \\ \nu(z,\tau) \end{pmatrix} \mapsto \begin{pmatrix} \partial_z \lambda(z,\tau) \\ \partial_z\mu(z,\tau)+\eta(1,\tau)\cdot \lambda(z,\tau) \\ \partial_z\nu(z,\tau)+\lambda(z,\tau) \\ \partial_{\tau}\lambda (z,\tau) \\ \partial_{\tau}\mu(z,\tau)-\frac{1}{2\pi i}\eta(1,\tau)\cdot \mu(z,\tau)+\frac{1}{2\pi i}\eta(1,\tau)^2\cdot \nu(z,\tau)-\partial_{\tau}\eta(1,\tau)\cdot \nu(z,\tau) \\ \partial_{\tau}\nu(z,\tau)-\frac{1}{2\pi i }\mu(z,\tau)+\frac{1}{2\pi i }\eta(1,\tau) \cdot \nu(z,\tau)  \end{pmatrix}.\]
\end{guess}
\subsubsection{First properties of the guess}
\begin{lemma}
The formula given in Guess 3.5.15 defines an integrable connection on $\mathcal L_1^{an}$ which prolongs the relative connection $(\nabla_1^{res})^{an}$.
\end{lemma}
\begin{proof}
At first it is not even clear if the vector proposed in the guess is a section of $\Omega^1_{E^{an}}\otimes_{\mathcal O_{E^{an}}} \mathcal L_1^{an}$.\\
By $(3.5.24)$ and Prop. 3.5.2 we know that the automorphy matrix for $\Omega^1_{E^{an}}\otimes_{\mathcal O_{E^{an}}} \mathcal L_1^{an}$, evaluated at
\[\Bigg( \begin{pmatrix}  m \\ n \end{pmatrix}, \begin{pmatrix} a & b \\ c & d \end{pmatrix}, (z,\tau) \Bigg) \in \Z^2\times \Gamma(N)\times \C\times \H,\]is given (as explained in 3.2 (iv)) by

\[\begin{pmatrix}
                        c\tau+d & 0\\
											 (cz+cn-dm)(c\tau+d) &  (c\tau+d)^2
\end{pmatrix} \otimes \begin{pmatrix} 1 & 0 & 0\\
                        2\pi i (dm-cz-cn) & c\tau+d & 0\\
												0 & 0 & \frac{1}{c\tau+d}
\end{pmatrix},
\]
i.e. by the following $(6\times 6)$-matrix:
{\footnotesize
\[\begin{pmatrix} c\tau+d & 0 & 0     & 0 & 0 & 0    \\
                        2\pi i (dm-cz-cn)(c\tau+d) & (c\tau+d)^2 & 0  & 0 & 0 & 0           \\
												0 & 0 &  1 & 0 & 0 & 0\\
												(cz+cn-dm)(c\tau+d) & 0 & 0 & (c\tau+d)^2 & 0& 0\\
                        -2\pi i (cz+cn-dm)^2(c\tau+d) & (cz+cn-dm)(c\tau+d)^2 & 0 & 2\pi i (dm-cz-cn)(c\tau+d)^2 & (c\tau+d)^3 & 0\\
                        0 & 0 & cz+cn-dm & 0 & 0 & c\tau+d

\end{pmatrix}.
\]}Hence, what we first need to check is that the vector of length $6$ in Guess 3.5.15 transforms under the effect of $\Z^2\times \Gamma(N)$ with the preceding matrix.\\
By hypothesis, the vector
\[\begin{pmatrix}\lambda(z,\tau) \\ \mu(z,\tau) \\ \nu(z,\tau) \end{pmatrix}\]
transforms with the automorphy matrix of $\mathcal L_1^{an}$. Moreover, we have the equality 
\[\eta \bigg(1,\frac{a\tau+b}{c\tau+d}\bigg)=(c\tau+d)^2\eta(1,\tau)+2\pi i c (c\tau+d),\]
as one shows e.g. by using $\eta(1,\tau)=-G_2(\tau)$ (cf. 3.3.2 (ii)) and \cite{Di-Shu}, Ch. 1.2, formula $(1.4)$.\\
With these two pieces of information (and with the induced transformation formulas for derivatives and squares) one verifies by explicit computation that the vector proposed in the guess indeed transforms with the above $(6\times 6)$-matrix; this calculation is quite laborious, but it can be done and exactly gives the desired result.\\
Next, it is easily checked that the formula in the guess satisfies the Leibniz rule and hence defines a connection. To see its integrability one observes $(3.5.26)$ and then simply performs the required calculation; this is again a bit tedious, but in the end it precisely works out.\\
That its restriction to an $S^{an}$-connection equals $(\nabla_1^{res})^{an}$ is clear by construction, cf. $(3.5.36)$.
\end{proof}
The following two properties are easily verified in small calculations, the first using $(3.5.33)$ and $(3.5.34)$ and the second following essentially from how we arrived at the definition of our guess.
\begin{lemma}
If $\mathcal L_1^{an}$ is equipped with the integrable connection defined in Guess 3.5.15, then:\\
(i) The analytification
\[0 \rightarrow \mathcal H_{E^{an}}^{an}\rightarrow \mathcal L_1^{an} \rightarrow \mathcal O_{E^{an}} \rightarrow 0\]
of $(3.5.31)$ becomes $\mathcal D_{E^{an}}$-linear.\\
(ii) The analytification $\varphi_1^{an}: \mathcal O_{S^{an}}\oplus \mathcal H^{an} \simeq (\epsilon^{an})^*\mathcal L_1^{an}$ of $\varphi_1$ becomes $\mathcal D_{S^{an}}$-linear. \qquad \qed
\end{lemma}

\subsubsection{Technical preparations for the main theorem}
For a little while (until Thm. 3.5.21) we no longer restrict to a single connected component of $E^{an}/S^{an}$ but really work with the description as disconnected complex manifolds given in $(3.4.4)$.\\
\newline
As in the algebraic case we have a Leray spectral sequence for de Rham cohomology:
\[\tag{\textbf{3.5.37}} E_2^{p,q}=H^p_{\mathrm{dR}}(S^{an}, H^q_{\mathrm{dR}}(E^{an}/S^{an},\mathcal H^{an}_{E^{an}})) \Rightarrow E^{p+q}=H^{p+q}_{\mathrm{dR}}(E^{an},\mathcal H^{an}_{E^{an}})\]
(the remarks of \cite{Kat2}, $(3.1)$ and $(3.3)$ hold invariantly in the analytic situation).\\
Because of already cited identifications we have
\[H^0_{\mathrm{dR}}(E^{an}/S^{an})\simeq H^0_{\mathrm{dR}}(E/S)^{an} \simeq (\mathcal O_{S})^{an}\simeq \mathcal O_{S^{an}};\]
this together with the existence of the analytic zero section $\epsilon^{an}$ implies that the five term exact sequence for $(3.5.37)$ yields a short exact sequence of $\C$-vector spaces
\[0\rightarrow H^1_{\mathrm{dR}}(S^{an}, \mathcal H^{an}) \xrightarrow{(\pi^{an})^*} H^1_{\mathrm{dR}}(E^{an},\mathcal H^{an}_{E^{an}}) \rightarrow H^0_{\mathrm{dR}}(S^{an},(\mathcal H^{an})^\vee \otimes_{\mathcal O_{S^{an}}}\mathcal H^{an})\rightarrow 0\]
which is split by the retraction $(\epsilon^{an})^*$.\\
Using the canonical identifications of the occurring de Rham cohomology spaces with the respective $\Ext$-spaces, true also in the analytic setting (cf. \cite{Ho-Ta-Tan}, Prop. 4.2.1 and the subsequent discussion), the previous sequence writes as the familiar split exact sequence

\[ \tag{\textbf{3.5.38}}
\begin{split}
0 & \rightarrow \Ext^1_{\mathcal D_{S^{an}}}(\mathcal O_{S^{an}}, \mathcal H^{an}) \xrightarrow{(\pi^{an})^*} \Ext^1_{\mathcal D_{E^{an}}}(\mathcal O_{E^{an}}, \mathcal H^{an}_{E^{an}})\\
& \rightarrow \Hom_{\mathcal D_{S^{an}}}(\mathcal O_{S^{an}}, (\mathcal H^{an})^\vee \otimes_{\mathcal O_{S^{an}}} \mathcal H^{an}) \rightarrow 0,
\end{split}
\]
where the projection is given analogously as described in 1.1.\\

We have a Leray spectral sequence of hypercohomology (cf. \cite{Dim}, Thm. 1.3.19 (ii)) for the de Rham complex $\Omega^{\bullet}_{E^{an}/S^{an}}(\mathcal H^{an}_{E^{an}})$ of $\mathcal H^{an}_{E^{an}}$ relative $S^{an}$:
\begin{align*}
E^{p,q}_2&=H^p(S^{an}, H^q_{\mathrm{dR}}(E^{an}/S^{an})\otimes_{\mathcal O_{S^{an}}}\mathcal H^{an}) \\
\Rightarrow E^{p+q}&=\mathbb H^{p+q}(E^{an}, \Omega^{\bullet}_{E^{an}/S^{an}}(\mathcal H^{an}_{E^{an}})) \simeq \Ext^{p+q}_{\mathcal D_{E^{an}/S^{an}}}(\mathcal O_{E^{an}},\mathcal H^{an}_{E^{an}}).
\end{align*}
As happened algebraically in $(2.1.3)$ one obtains from the associated five term sequence a split short exact sequence forming the lower row of a commutative diagram whose upper row is $(3.5.38)$:

\[
{\small
\tag{\textbf{3.5.39}}
\begin{split}
\xymatrix@C-0.3cm{
0 \ar[r] & \Ext^1_{\mathcal D_{S^{an}}}(\mathcal O_{S^{an}}, \mathcal H^{an}) \ar[d]^{\mathrm{can}} \ar[r]^-{(\pi^{an})^*} & \Ext^1_{\mathcal D_{E^{an}}}(\mathcal O_{E^{an}}, \mathcal H^{an}_{E^{an}}) \ar[d]^{\mathrm{can}} \ar[r]& \Hom_{\mathcal D_{S^{an}}}(\mathcal O_{S^{an}}, (\mathcal H^{an})^\vee \otimes_{\mathcal O_{S^{an}}} \mathcal H^{an}) \ar[d]^{\mathrm{can}} \ar[r] &0 \\
0 \ar[r] & \Ext^1_{\mathcal O_{S^{an}}}(\mathcal O_{S^{an}}, \mathcal H^{an}) \ar[r]^-{(\pi^{an})^*} & \Ext^1_{\mathcal D_{E^{an}/S^{an}}}(\mathcal O_{E^{an}}, \mathcal H^{an}_{E^{an}}) \ar[r] & \Hom_{\mathcal O_{S^{an}}}(\mathcal O_{S^{an}}, (\mathcal H^{an})^\vee \otimes_{\mathcal O_{S^{an}}} \mathcal H^{an}) \ar[r] & 0}
\end{split}
}
\]

We will need the fact that also Lemma 2.1.1 and 2.1.2 are available in the analytic situation:

\begin{lemma}
Suppose we are given two extensions of $\mathcal D_{E^{an}/S^{an}}$-modules
\[M: \quad 0\rightarrow \mathcal H^{an}_{E^{an}} \xrightarrow{j_M} \mathcal M \xrightarrow{p_M} \mathcal  O_{E^{an}} \rightarrow 0\]
\[N: \quad 0\rightarrow \mathcal H^{an}_{E^{an}} \xrightarrow{j_N} \mathcal N \xrightarrow{p_N} \mathcal  O_{E^{an}} \rightarrow 0\]
with $\mathcal O_{S^{an}}$-linear splittings
\[\varphi_M: \mathcal O_{S^{an}} \oplus \mathcal H^{an} \simeq (\epsilon^{an})^*\mathcal M\]
\[\varphi_N: \mathcal O_{S^{an}} \oplus \mathcal H^{an} \simeq (\epsilon^{an})^*\mathcal N\]
and the property that the classes of $M$ and $N$ in $\Ext^1_{\mathcal D_{E^{an}/S^{an}}}(\mathcal O_{E^{an}},\mathcal H^{an}_{E^{an}})$ are equal - e.g. if they both map to the identity under the lower projection of $(3.5.39)$.\\
Then there exists a unique isomorphism of $M$ and $N$ which respects the splittings.
\end{lemma}
\begin{proof} \ \\
\underline{Existence}: Verbatim as in the proof of Lemma 2.1.1.\\
\newline
\underline{Uniqueness}:  It suffices to show: if $(M, \varphi_M)$ is an extension with splitting as in the claim of the lemma, then any automorphism $h$ of this pair is the identity.\\
The claim that $h: \mathcal M \rightarrow \mathcal M$ equals $\id$ can be checked on the fibers over the points of $S^{an}$; such a fiber writes as $\mathcal E^{an}$, the analytification of an elliptic curve $\mathcal E/\Spec(\C)$. Let $e^{an}$ resp. $e$ be the zero point of $\mathcal E^{an}$ resp. $\mathcal E$.\\
On $\mathcal E^{an}$ we are then given a horizontal automorphism $h: \mathcal M \simeq \mathcal M$ of a $\mathcal O_{\mathcal E^{an}}$-vector bundle $\mathcal M$ with integrable connection such that the induced map
\[\tag{$*$} h_{e^{an}}\otimes \id: \mathcal M_{e^{an}}\otimes_{\mathcal O_{\mathcal E^{an},e^{an}}} \C \simeq \mathcal M_{e^{an}}\otimes_{\mathcal O_{\mathcal E^{an},e^{an}}} \C\]
is the identity. We need to show that already $h$ is the identity.\\
By standard GAGA-results (cf. \cite{Mal}, §1) one knows that analytification induces an equivalence between the category of algebraic vector bundles with integrable connection on $\mathcal E$ and the respective analytic category on $\mathcal E^{an}$. Hence, the $\mathcal D_{\mathcal E^{an}}$-module $\mathcal M$ and the map $h$ come from the algebraic side. It is easily deduced from $(*)$ that the analogous "fiber-in-$e$-morphism" for the algebraic vector bundle automorphism is the identity. But then this automorphism is already the identity, by the same argument as in the uniqueness part in the proof of Lemma 2.1.1. This of course implies that the analytic automorphism $h$ is the identity, which is what we wanted to show.
\end{proof}
\begin{lemma}
Suppose we are given two extensions of $\mathcal D_{E^{an}}$-modules
\begin{align*}
M & : \quad 0\rightarrow \mathcal H^{an}_{E^{an}} \xrightarrow{j_M} \mathcal M \xrightarrow{p_M} \mathcal  O_{E^{an}} \rightarrow 0\\
N &: \quad 0\rightarrow \mathcal H^{an}_{E^{an}} \xrightarrow{j_N} \mathcal N \xrightarrow{p_N} \mathcal  O_{E^{an}} \rightarrow 0
\end{align*}
with $\mathcal D_{S^{an}}$-linear splittings
\begin{align*}
\varphi_M& : \mathcal O_{S^{an}} \oplus \mathcal H^{an} \simeq (\epsilon^{an})^*\mathcal M\\
\varphi_N& : \mathcal O_{S^{an}} \oplus \mathcal H^{an} \simeq (\epsilon^{an})^*\mathcal N
\end{align*}
and the property that the classes of $M$ and $N$ in $\Ext^1_{\mathcal D_{E^{an}}}(\mathcal O_{E^{an}},\mathcal H^{an}_{E^{an}})$ are equal - e.g. if they both map to the identity under the upper projection of $(3.5.39)$.
Then there exists a unique isomorphism of $M$ and $N$ which respects the splittings.
\end{lemma}
\begin{proof}
Same argument as in Lemma 2.1.2.
\end{proof}

It is formal to deduce from the preceding two lemmas the following analytic version of Prop. 2.1.4. We will only need the uniqueness statement, but also the existence of a prolongation with the mentioned properties can be shown analogously as in the algebraic case.

\begin{proposition}
Assume we are given a $\mathcal D_{E^{an}/S^{an}}$-linear extension
\[\tag{\textbf{3.5.40}} M: \quad 0\rightarrow \mathcal H^{an}_{E^{an}} \rightarrow \mathcal M \rightarrow \mathcal O_{E^{an}} \rightarrow 0,\]
whose class maps to the identity under the lower projection of $(3.5.39)$, together with a $\mathcal O_{S^{an}}$-linear splitting for its pullback along $\epsilon^{an}$:
\[\varphi_M: \mathcal O_{S^{an}} \oplus \mathcal H^{an} \simeq (\epsilon^{an})^*\mathcal M.\]
Then there exists at most one prolongation of the integrable $S^{an}$-connection on $\mathcal M$ to an absolute integrable connection such that the following holds:\\
If we endow $\mathcal M$ with this connection, then the $\mathcal D_{E^{an}/S^{an}}$-linear exact sequence $(3.5.40)$ becomes $\mathcal D_{E^{an}}$-linear and the $\mathcal O_{S^{an}}$-linear splitting $\varphi_M$ becomes $\mathcal D_{S^{an}}$-linear.\\
(A prolongation with these properties in fact exists.)
\end{proposition}
\begin{proof}
Let $\nabla_{\mathcal M}$ and $\widetilde{\nabla}_{\mathcal M}$ be two (absolute) integrable connections on $\mathcal M$ with the properties of the claim. Endowing $\mathcal M$ one time with $\nabla_{\mathcal M}$ and the other time with $\widetilde{\nabla}_{\mathcal M}$ we obtain (by assumption) from $(3.5.40)$ two $\mathcal D_{E^{an}}$-linear extensions. The images of their classes under the projection in the upper row of $(3.5.39)$ both times are the identity because of the commutativity of $(3.5.39)$ and the hypothesis about $(3.5.40)$. As by assumption both classes retract to zero in $\Ext^1_{\mathcal D_{S^{an}}}(\mathcal O_{S^{an}},\mathcal H^{an})$ we conclude from the splitting of the upper row of $(3.5.39)$ that the two extension classes we obtained from $(3.5.40)$ are equal in $\Ext^1_{\mathcal D_{E^{an}}}(\mathcal O_{E^{an}},\mathcal H^{an}_{E^{an}})$.\\
Then, by Lemma 3.5.19, there exists a $\mathcal D_{E^{an}}$-linear isomorphism
\[\nu: (\mathcal M, \nabla_{\mathcal M}) \simeq (\mathcal M,\widetilde{\nabla}_{\mathcal M})\]
respecting the extension structure of $\mathcal M$ and the splitting $\varphi_{\mathcal M}$. Now restrict absolute structures to $S^{an}$-structures: then, as $\nabla_{\mathcal M}$ and $\widetilde{\nabla}_{\mathcal M}$ are equal when considered as $S^{an}$-connections, $\nu$ yields an automorphism of the $\mathcal D_{E^{an}/S^{an}}$-linear extension $(3.5.40)$ with its $\mathcal O_{S^{an}}$-linear splitting $\varphi_M$. The uniqueness part of Lemma 3.5.18 then implies that $\nu=\id$. But as $\nu$ is horizontal for the absolute connections $\nabla_{\mathcal M}$ and $\widetilde{\nabla}_{\mathcal M}$ we conclude that $\nabla_{\mathcal M}=\widetilde{\nabla}_{\mathcal M}$.
\end{proof}

\subsubsection{Proof of the main theorem}
We are ready to show that the formula suggested in Guess 3.5.15 indeed describes the analytification $\nabla_1^{an}$ of the integrable $\Q$-connection $\nabla_1$.
\begin{theorem}
The connection $\nabla_1^{an}: \mathcal L_1^{an} \rightarrow \Omega^1_{E^{an}}\otimes_{\mathcal O_{E^{an}}} \mathcal L_1^{an}$ is given for each connected component on local sections by the formula
\[\begin{pmatrix} \lambda (z,\tau) \\ \mu(z,\tau) \\ \nu(z,\tau) \end{pmatrix} \mapsto \begin{pmatrix} \partial_z \lambda(z,\tau) \\ \partial_z\mu(z,\tau)+\eta(1,\tau)\cdot \lambda(z,\tau) \\ \partial_z\nu(z,\tau)+\lambda(z,\tau) \\ \partial_{\tau}\lambda (z,\tau) \\ \partial_{\tau}\mu(z,\tau)-\frac{1}{2\pi i}\eta(1,\tau)\cdot \mu(z,\tau)+\frac{1}{2\pi i}\eta(1,\tau)^2\cdot \nu(z,\tau)-\partial_{\tau}\eta(1,\tau)\cdot \nu(z,\tau) \\ \partial_{\tau}\nu(z,\tau)-\frac{1}{2\pi i }\mu(z,\tau)+\frac{1}{2\pi i }\eta(1,\tau) \cdot \nu(z,\tau)  \end{pmatrix}.\]
\end{theorem}

\begin{proof}
Consider the analytification
\[\tag{$*$} 0\rightarrow \mathcal H_{E^{an}}^{an} \rightarrow \mathcal L_1^{an} \rightarrow \mathcal O_{E^{an}} \rightarrow 0\]
of the $\mathcal D_{E/S}$-linear extension $(3.5.31)$ as well as of the splitting $\varphi_1$:
\[\tag{$**$}\varphi_1^{an}: \mathcal O_{S^{an}} \oplus \mathcal H^{an} \simeq (\epsilon^{an})^*\mathcal L_1^{an}.\]
The class of $(*)$ maps to the identity under the lower projection of $(3.5.39)$ because the respective fact holds algebraically for $(3.5.31)$.\\
\newline
The integrable $\Q$-connection $\nabla_1$ is characterized by the following property: it prolongs the integrable $S$-connection $\nabla_1^{res}$, and if $\mathcal L_1$ is equipped with $\nabla_1$, then:
\begin{align*}
&\textrm{(i) The $\mathcal D_{E/S}$-linear exact sequence $(3.5.31)$ becomes $\mathcal D_{E/\Q}$-linear.}\\
&\textrm{(ii) The $\mathcal O_S$-linear splitting $\varphi_1$ becomes $\mathcal D_{S/\Q}$-linear.}
\end{align*}
Cf. the recap at the beginning of this subsection.\\
\newline
The analytification $\nabla_1^{an}$ then satisfies the following property: it prolongs the integrable $S^{an}$-connection $(\nabla_1^{res})^{an}$, and if $\mathcal L_1^{an}$ is equipped with $\nabla_1^{an}$, then:
\begin{align*}
&\textrm{(i)$^{an}$ The exact sequence $(*)$ becomes $\mathcal D_{E^{an}}$-linear. \qquad \qquad \quad \ \ }\\
&\textrm{(ii)$^{an}$ The splitting $(**)$ becomes $\mathcal D_{S^{an}}$-linear. \qquad \qquad \quad \ \ }
\end{align*}
But by Lemma 3.5.16 and 3.5.17 the formula in the claim of the theorem likewise defines an (absolute) integrable connection on $\mathcal L_1^{an}$ which prolongs $(\nabla_1^{res})^{an}$ and satisfies (i)$^{an}$ and (ii)$^{an}$ if $\mathcal L_1^{an}$ is equipped with this connection.\\
As by Prop. 3.5.20 there can exist at most one such connection we conclude the claim.
\end{proof}

\subsubsection{A vanishing result}

With Thm. 3.5.21 we can prove the following observation which becomes relevant when one wants to define (absolute) de Rham cohomology classes with coefficients in $\mathcal L_n^{an}$. We will need it for our explicit construction of the analytified $D$-variant of the polylogarithm.\\
\newline
From now on we again work on a fixed connected component, applying the usual abuse of notation.

\begin{proposition}
Let $V$ be an open subset of $E^{an}$ and let $n\geq 0$. Assume that for $k=0,...,n+1$ we have holomorphic functions $s_k(z,\tau)$ on $\mathrm{pr}^{-1}(V)\subseteq \C\times \H$ such that the vector of length $2\cdot r(n)$ given by
\[\begin{pmatrix} s_0(z,\tau) \\ s_1(z,\tau) \\ \vdots\\s_n(z,\tau)  \\  0 \\ \vdots \\ \vdots \\ 0 \\-\frac{1}{2\pi i }s_1(z,\tau)\\ -\frac{2}{2\pi i }s_2(z,\tau) \\ \vdots\\ -\frac{n+1}{2\pi i }s_{n+1}(z,\tau) \\ 0 \\ \vdots \\ \vdots \\ 0\end{pmatrix}\]
defines a section of $\Omega^1_{E^{an}}\otimes_{\mathcal O_{E^{an}}} \mathcal L_n^{an}$ over $V$; here, after the first $n+1$ entries we fill up with zeroes until we reach $r(n)$ entries, and we perform the same procedure in the second half of the vector.\\
Assume furthermore that the functions $s_k(z,\tau)$ satisfy the differential equation
\[\partial_{\tau}s_k(z,\tau)=-\frac{k+1}{2\pi i }\cdot  \partial_zs_{k+1}(z,\tau) \quad \textrm{for all} \ k=0,...,n.\]
Then the above section goes to zero under the map
\[\Gamma(V,\Omega^1_{E^{an}}\otimes_{\mathcal O_{E^{an}}} \mathcal L_n^{an}) \rightarrow \Gamma(V,\Omega^2_{E^{an}}\otimes_{\mathcal O_{E^{an}}} \mathcal L_n^{an})\]
of the de Rham complex of $(\mathcal L_n^{an},\nabla_n^{an})$.
\end{proposition}

\begin{proof}
Recall that starting from the trivialization
\[\mathcal O_{\C\times \H}^{\oplus 3} \simeq \mathrm{pr}^*\mathcal L_1^{an}\]
by the ordered basic sections $\{e,f,g\}$ we have fixed for each $0\leq n$ the trivialization
\[\mathcal O_{\C\times \H}^{\oplus r(n)} \simeq \mathrm{pr}^*\mathcal L_n^{an}\]
by the $r(n)$ many ordered basic sections

{\footnotesize
\[\bigg \{\frac{e^n}{n!}, \frac{e^{n-1}f}{(n-1)!}, \frac{e^{n-2}f^2}{(n-2)!},..., \frac{f^n}{(n-n)!}, \frac{e^{n-1}g}{(n-1)!}, \frac{e^{n-2}fg}{(n-2)!},\frac{e^{n-3}f^2g}{(n-3)!},..., \frac{f^{n-1}g}{(n-n)!}, \frac{e^{n-2}g^2}{(n-2)!},.........,\frac{g^n}{(n-n)!} \bigg \}.\]}Recall moreover that the trivialization of $\mathrm{pr}^*\Omega^1_{E^{an}} \simeq \Omega^1_{\C\times \H}$ is given by $\{ \mathrm{d}z,\mathrm{d}\tau \}$.\\
Thm. 3.5.21 implies that the pullback connection $\mathrm{pr}^*\nabla_1^{an}$ acts on the basic sections $e,f,g$ as follows:
\begin{align*}
e&\mapsto \eta(1,\tau)\cdot \mathrm{d}z \otimes f + \mathrm{d}z \otimes g,\\
f&\mapsto -\frac{1}{2\pi i } \eta(1,\tau) \cdot \mathrm{d}\tau \otimes f -\frac{1}{2\pi i } \cdot \mathrm{d}\tau \otimes g,\\
g&\mapsto \bigg(\frac{1}{2\pi i } \eta(1,\tau)^2-\partial_{\tau}\eta(1,\tau)\bigg)\cdot \mathrm{d}\tau \otimes f+\frac{1}{2\pi i }\eta(1,\tau)\cdot d \tau \otimes g.
\end{align*}
Hence, for $0\leq n$ and $0\leq k \leq n$ one computes that $\mathrm{pr}^*\nabla_n^{an}$ acts on the basic section $\frac{e^{n-k}f^k}{(n-k)!}$ as
\begin{align*}
\frac{e^{n-k}f^k}{(n-k)!} \mapsto \ &\eta(1,\tau)\cdot \mathrm{d}z \otimes \frac{e^{n-k-1}f^{k+1}}{(n-k-1)!}+ \mathrm{d}z \otimes \frac{e^{n-k-1}f^kg}{(n-k-1)!}\\
&-\frac{k}{2\pi i } \eta(1,\tau)\cdot \mathrm{d}\tau \otimes \frac{e^{n-k}f^k}{(n-k)!}-\frac{k}{2\pi i }\cdot \mathrm{d}\tau \otimes \frac{e^{n-k}f^{k-1}g}{(n-k)!},
\end{align*}
where we make the convention that terms containing the expressions $e^{-1}$ or $f^{-1}$ shall be zero.\\
The section in the claim of the corollary is defined by the section of $\Omega^1_{\C \times \H} \otimes_{\mathcal O_{\C\times \H}} \mathrm{pr}^*\mathcal L_n^{an}$ given by
\[\sum_{k=0}^ns_k(z,\tau)\cdot \mathrm{d}z \otimes \frac{e^{n-k}f^k}{(n-k)!}-\sum_{k=0}^n \frac{(k+1)}{2\pi i }s_{k+1}(z,\tau)\cdot \mathrm{d}\tau \otimes \frac{e^{n-k}f^k}{(n-k)!}.\]
Observing $\mathrm{d}z\wedge \mathrm{d}z=0=\mathrm{d}\tau\wedge \mathrm{d}\tau$ we see that under the map of the de Rham complex this goes to
\begin{align*}
&\sum_{k=0}^n  \partial_{\tau}s_k(z,\tau)\cdot \mathrm{d}\tau\wedge \mathrm{d}z \otimes \frac{e^{n-k}f^k}{(n-k)!}-\sum_{k=0}^n \frac{(k+1)}{2\pi i }\partial_zs_{k+1}(z,\tau)\cdot \mathrm{d}z\wedge \mathrm{d}\tau \otimes \frac{e^{n-k}f^k}{(n-k)!}\\
- &\sum_{k=0}^n s_k(z,\tau)\cdot \mathrm{d}z \wedge\bigg(-\frac{k}{2\pi i } \eta(1,\tau)\cdot \mathrm{d}\tau \otimes \frac{e^{n-k}f^k}{(n-k)!}-\frac{k}{2\pi i }\cdot \mathrm{d}\tau \otimes \frac{e^{n-k}f^{k-1}g}{(n-k)!}\bigg)\\
+&\sum_{k=0}^n \frac{(k+1)}{2\pi i }s_{k+1}(z,\tau)\cdot \mathrm{d}\tau \wedge \bigg( \eta(1,\tau)\cdot \mathrm{d}z \otimes \frac{e^{n-k-1}f^{k+1}}{(n-k-1)!}+ \mathrm{d}z \otimes \frac{e^{n-k-1}f^kg}{(n-k-1)!} \bigg),
\end{align*}
which in turn writes out as $(\mathrm{d}z\wedge \mathrm{d}\tau)$- times the following expression:
\begin{align*}
-&\sum_{k=0}^n \partial_{\tau}s_k(z,\tau) \otimes  \frac{e^{n-k}f^k}{(n-k)!}-\sum_{k=0}^n\frac{(k+1)}{2\pi i }\partial_zs_{k+1}(z,\tau)  \otimes \frac{e^{n-k}f^k}{(n-k)!}\\
+ &\sum_{k=0}^n\frac{k}{2\pi i }\eta(1,\tau) s_k(z,\tau) \otimes \frac{e^{n-k}f^k}{(n-k)!} +\sum_{k=0}^n\frac{k}{2\pi i} s_k(z,\tau) \otimes\frac{e^{n-k}f^{k-1}g}{(n-k)!}\\
-&\sum_{k=0}^n\frac{(k+1)}{2\pi i } \eta(1,\tau) s_{k+1}(z,\tau) \otimes \frac{e^{n-k-1}f^{k+1}}{(n-k-1)!}-\sum_{k=0}^n\frac{(k+1)}{2\pi i } s_{k+1}(z,\tau) \otimes \frac{e^{n-k-1}f^kg}{(n-k-1)!}.
\end{align*}
The third and the fifth summand obviously cancel to zero; the same holds for the fourth and the sixth. Hence, what remains is $(\mathrm{d}z\wedge \mathrm{d}\tau)$-times the expression
\begin{align*}
-&\sum_{k=0}^n \partial_{\tau}s_k(z,\tau) \otimes  \frac{e^{n-k}f^k}{(n-k)!}-\sum_{k=0}^n\frac{(k+1)}{2\pi i }\partial_zs_{k+1}(z,\tau) \otimes \frac{e^{n-k}f^k}{(n-k)!}\\
=-&\sum_{k=0}^n\bigg(\partial_{\tau}s_k(z,\tau) +\frac{(k+1)}{2\pi i }\partial_zs_{k+1}(z,\tau) \bigg) \otimes \frac{e^{n-k}f^k}{(n-k)!}.
\end{align*}
It is now the differential equation
\[\partial_{\tau}s_k(z,\tau)=-\frac{k+1}{2\pi i }\cdot  \partial_zs_{k+1}(z,\tau) \quad \textrm{for all} \ k=0,...,n.\]
of the hypothesis which annihilates also this last expression.\\
What we have shown implies altogether that the image of the section in the claim under the map
\[\Omega^1_{E^{an}}\otimes_{\mathcal O_{E^{an}}} \mathcal L_n^{an} \rightarrow \Omega^2_{E^{an}}\otimes_{\mathcal O_{E^{an}}} \mathcal L_n^{an}\]
pulls back to the zero section of $\Omega^2_{\C\times \H}\otimes_{\mathcal O_{\C\times \H}} \mathrm{pr}^*\mathcal L_n^{an}$ under $\mathrm{pr}$. But then already the image in $\Omega^2_{E^{an}}\otimes_{\mathcal O_{E^{an}}} \mathcal L_n^{an}$ must be zero (cf. the isomorphism in 3.2 (ii)), and this was to prove.
\end{proof}

\markright{\uppercase{The explicit description on the universal elliptic curve}}
\section{The two fundamental systems of sections}
\markright{\uppercase{The explicit description on the universal elliptic curve}}
Fix an integer $D>1$ and let $U_{D}$ be the open complement of the $D$-torsion subscheme $E[D]$ of $E$.\\
\newline
When working with analytifications we restrict as usual to a fixed connected component whose index is suppressed in the notation (cf. the conventions at the beginning of 3.5 or 3.5.3). We then let again
\[\mathrm{pr}: \C\times \H\rightarrow E^{an}=(\Z^2\times \Gamma(N))\backslash (\C\times \H)\]
be the universal covering map.\\
\newline
Recall from $(3.3.22)$ and $(3.3.23)$ that we have a Laurent expansion
\[D^2 \cdot J(z,-w,\tau)-D\cdot J\Big(Dz,-\frac{w}{D},\tau \Big)=s^D_0(z,\tau)+s^D_1(z,\tau)w+...\]
with meromorphic functions $s^D_k(z,\tau), \ k\geq 0,$ on $\C \times \H$ which are holomorphic on $\mathrm{pr}^{-1}(U_{D}^{an})$.

\subsection{The construction}

By means of the coefficient functions $s_k^D(z,\tau)$ we construct two compatible systems of sections
\begin{align*}
\Big(q^D_n(z,\tau)\Big)_{n\geq 0} &\in \lim_{n\geq 0} \Gamma(U_D^{an},\Omega^1_{E^{an}/S^{an}}\otimes_{\mathcal O_{E^{an}}}\mathcal L_n^{an}),\\
\Big(p^D_n(z,\tau)\Big)_{n\geq 0} &\in \lim_{n\geq 0} \Gamma(U_D^{an},\Omega^1_{E^{an}}\otimes_{\mathcal O_{E^{an}}}\mathcal L_n^{an})
\end{align*}
and show that the $p_n^D(z,\tau)$ vanish under the morphism in the de Rham complex of $(\mathcal L_n^{an},\nabla_n^{an})$.\\
\newline
The future role of these systems will be the following:\\
The mentioned vanishing means that the second system induces an element of
\[\lim_{n\geq 0} H^1_{\mathrm{dR}}(U_D^{an},\mathcal L_n^{an}),\]
and a main result of the work (to be proven in 3.8.1) will show that it equals the analytification of
\[\Big(\varDn \Big)_{n\geq 0}.\]
In this context, we will view the first of the above systems as an element of
\[\lim_{n\geq 0} \Gamma(S^{an},H^1_{\mathrm{dR}}(U_D^{an}/S^{an},\mathcal L_n^{an}))\]
which we have to consider for the residue computation of the second system.

\subsubsection{The first section}

For each $n\geq 0$ we build the following vector of holomorphic functions on $\mathrm{pr}^{-1}(U_{D}^{an}) \subseteq \C\times \H$:
\[\tag{\textbf{3.6.1}} q_n^D(z,\tau):=\begin{pmatrix} s^D_0(z,\tau) \\ s^D_1(z,\tau) \\ \vdots\\s^D_n(z,\tau) \\ 0 \\ \vdots \\ \vdots \\ 0 \end{pmatrix},\]
where after the first $(n+1)$ entries we fill up with zeroes until we get a vector of length $r(n)$.

\begin{proposition}
With our fixed trivializations for the pullbacks of $\Omega^1_{E^{an}/S^{an}}$, $\mathcal L_n^{an}$ and $\Omega^1_{E^{an}/S^{an}}\otimes_{\mathcal O_{E^{an}}} \mathcal L_n^{an}$ to the universal covering the vector $(3.6.1)$ defines a section in $\Gamma(U_D^{an},\Omega^1_{E^{an}/S^{an}}\otimes_{\mathcal O_{E^{an}}} \mathcal L_n^{an})$.
\end{proposition}
\begin{proof}
Let
\[\widetilde{A}_n: \Z^2 \times \Gamma(N) \times \C \times \H \rightarrow \mathrm{GL}_{r(n)}(\C)\]
be the automorphy matrix for $\Omega^1_{E^{an}/S^{an}} \otimes_{\mathcal O_{E^{an}}} \mathcal L_n^{an}$. Its evaluation at
\[\Bigg(\begin{pmatrix}  \underline{m} \\ \underline{n} \end{pmatrix}, \begin{pmatrix} a & b \\ c & d \end{pmatrix}, (z,\tau) \Bigg) \in \Z^2\times \Gamma(N)\times \C \times \H\]
is given by
\[\widetilde{A}_n\Bigg(\begin{pmatrix}  \underline{m} \\ \underline{n} \end{pmatrix}, \begin{pmatrix} a & b \\ c & d \end{pmatrix},(z,\tau)\Bigg)=(c\tau+d)\cdot A_n\Bigg(\begin{pmatrix}  \underline{m} \\ \underline{n} \end{pmatrix}, \begin{pmatrix} a & b \\ c & d \end{pmatrix},(z,\tau)\Bigg),\]
which is the Kronecker product of the automorphy matrix for $\Omega^1_{E^{an}/S^{an}}$ with the automorphy matrix $A_n$ for $\mathcal L_n^{an}$ (cf. 3.2 (iv) and $(3.5.29)$).\\
\newline
It is of advantage to recall the shape of the matrix $A_n$, already recorded in $(3.5.11)$. For notational convenience we leave away the argument of $A_n$.

{\footnotesize
\[
A_n=
\left( \begin{array}{cccc|cccc}
a_0	& 0 	& \dots  & 0& 0& \dots & \dots & 0 \\
	(c\tau+d)\cdot a_1& (c\tau+d)\cdot a_0	& \dots   &0 &0 & \dots & \dots & 0 \\
	\vdots & \vdots  	&   &  \vdots &   \vdots & & & \vdots \\
	(c\tau+d)^n \cdot a_n& (c\tau+d)^n\cdot a_{n-1}   & \dots     & (c\tau+d)^n\cdot a_0 & 0 & \dots & \dots &0 \\
    \hline
 	0 & 0   & \dots     & 0       & * & \dots & \dots & * \\
	\vdots &  \vdots   &     &\vdots   &\vdots & & & \vdots \\
	\vdots & \vdots & & \vdots & \vdots & & & \vdots \\
	0&   0& \dots    &    0    &* &\dots & \dots & *
\end{array}\right)
\]}

Here, for all $r \geq 0$
\[a_r:=a_r\Bigg(\begin{pmatrix}  \underline{m} \\ \underline{n} \end{pmatrix}, \begin{pmatrix} a & b \\ c & d \end{pmatrix}, (z,\tau) \Bigg)=\frac{1}{r!}\bigg(\frac{2\pi i (d\underline{m}-cz-c\underline{n})}{c\tau+d}\bigg)^r\]
was defined as the coefficient at $w^r$ in the expansion around $w=0$ of
\[a\Bigg( \begin{pmatrix}  \underline{m} \\ \underline{n} \end{pmatrix}, \begin{pmatrix} 0 \\ 0 \end{pmatrix}, \begin{pmatrix} a & b \\ c & d \end{pmatrix}, (z,w,u,\tau) \Bigg)\]
with $a$ the factor of automorphy for $\mathcal P^{an}$ (cf. Def. 3.5.5).\\
\newline
What we need to show is (cf. 3.2 (ii)) that for each
\[\Bigg(\begin{pmatrix}  \underline{m} \\ \underline{n} \end{pmatrix}, \begin{pmatrix} a & b \\ c & d \end{pmatrix}\Bigg) \in \Z^2\times \Gamma(N)\]
we have the following equation of vectors of holomorphic functions on $\mathrm{pr}^{-1}(U_D^{an})$:
\[\begin{pmatrix} s^D_0(\frac{z+\underline{m}\tau+\underline{n}}{c\tau+d},\frac{a\tau+b}{c\tau+d}) \\ s^D_1(\frac{z+\underline{m}\tau+\underline{n}}{c\tau+d},\frac{a\tau+b}{c\tau+d}) \\ \vdots\\s^D_n(\frac{z+\underline{m}\tau+\underline{n}}{c\tau+d},\frac{a\tau+b}{c\tau+d}) \\ 0 \\ \vdots \\ \vdots \\ 0 \end{pmatrix}=\widetilde{A}_n\Bigg(\begin{pmatrix}  \underline{m} \\ \underline{n} \end{pmatrix}, \begin{pmatrix} a & b \\ c & d \end{pmatrix},(z,\tau)\Bigg)\cdot \begin{pmatrix} s^D_0(z,\tau) \\ s^D_1(z,\tau) \\ \vdots\\s^D_n(z,\tau) \\ 0 \\ \vdots \\ \vdots \\ 0 \end{pmatrix},\]
or equivalently (by the shape of $\widetilde{A}_n$) that for each $0\leq k \leq n$ we have
\[\tag{$*$} s^D_k\Big(\frac{z+\underline{m}\tau+\underline{n}}{c\tau+d},\frac{a\tau+b}{c\tau+d}\Big)=(c\tau+d)^{k+1}\cdot \Big(a_k \cdot s_0^D(z,\tau)+a_{k-1}\cdot s^D_1(z,\tau)+...+a_0\cdot s^D_k(z,\tau) \Big).\]
Now note at first that
\[\tag{$**$} \begin{split}
&J\bigg(\frac{z+\underline{m}\tau+\underline{n}}{c\tau+d},-\frac{w}{c\tau+d},\frac{a\tau+b}{c\tau+d}\bigg)\\
&=(c\tau+d) \cdot a\Bigg( \begin{pmatrix}  \underline{m} \\ \underline{n} \end{pmatrix}, \begin{pmatrix} 0 \\ 0 \end{pmatrix}, \begin{pmatrix} a & b \\ c & d \end{pmatrix}, (z,w,u,\tau) \Bigg) \cdot J(z,-w,\tau)\end{split}\]
by Cor. 3.3.14 and the formula $(3.4.14)$ for $a$.\\
Furthermore, the equality (implied by $(**)$)
\[\begin{split}
&J\bigg(\frac{D\cdot(z+\underline{m}\tau+\underline{n})}{c\tau+d},-\frac{\frac{w}{D}}{c\tau+d},\frac{a\tau+b}{c\tau+d}\bigg)\\
&= (c\tau+d)\cdot a\Bigg( \begin{pmatrix}  D\underline{m} \\ D\underline{n} \end{pmatrix}, \begin{pmatrix} 0 \\ 0 \end{pmatrix}, \begin{pmatrix} a & b \\ c & d \end{pmatrix}, (Dz,\frac{w}{D},u,\tau) \Bigg) \cdot J \Big(Dz,-\frac{w}{D},\tau \Big)
\end{split}\]
and the fact (following from the explicit formula $(3.4.14)$ for $a$) that
\[a\Bigg( \begin{pmatrix}  D\underline{m} \\ D \underline{n}\end{pmatrix}, \begin{pmatrix} 0 \\ 0 \end{pmatrix}, \begin{pmatrix} a & b \\ c & d \end{pmatrix}, \Big(Dz,\frac{w}{D},u,\tau \Big) \Bigg)=a\Bigg( \begin{pmatrix}  \underline{m} \\ \underline{n} \end{pmatrix}, \begin{pmatrix} 0 \\ 0 \end{pmatrix}, \begin{pmatrix} a & b \\ c & d \end{pmatrix}, (z,w,u,\tau) \Bigg)\]
yield
\[\tag{$***$}
\begin{split}
&J\bigg(\frac{D\cdot(z+\underline{m}\tau+\underline{n})}{c\tau+d},-\frac{\frac{w}{D}}{c\tau+d},\frac{a\tau+b}{c\tau+d}\bigg)\\
&=(c\tau+d) \cdot a\Bigg( \begin{pmatrix}  \underline{m} \\ \underline{n} \end{pmatrix}, \begin{pmatrix} 0 \\ 0 \end{pmatrix}, \begin{pmatrix} a & b \\ c & d \end{pmatrix}, (z,w,u,\tau) \Bigg) \cdot J \Big(Dz,-\frac{w}{D},\tau \Big).\end{split}
\]
From $(**)$ and $(***)$ we get
\begin{align*}
&D^2\cdot J\bigg(\frac{z+\underline{m}\tau+\underline{n}}{c\tau+d},-\frac{w}{c\tau+d},\frac{a\tau+b}{c\tau+d}\bigg)-D\cdot J\bigg(\frac{D\cdot(z+\underline{m}\tau+\underline{n})}{c\tau+d},-\frac{\frac{w}{D}}{c\tau+d},\frac{a\tau+b}{c\tau+d}\bigg)\\
&=(c\tau+d)\cdot a\Bigg( \begin{pmatrix}  \underline{m} \\ \underline{n} \end{pmatrix}, \begin{pmatrix} 0 \\ 0 \end{pmatrix}, \begin{pmatrix} a & b \\ c & d \end{pmatrix}, (z,w,u,\tau) \Bigg)\cdot \bigg[D^2\cdot J(z,-w,\tau)-D\cdot J\Big(Dz,-\frac{w}{D},\tau \Big) \bigg].
\end{align*}
Comparing coefficients at $w^k$ in the expansion around $w=0$ of this last equation yields $(*)$.
\end{proof}

\subsubsection{The second section}

Next, again by means of the coefficient functions $s^D_k(z,\tau)$, we define for each $n\geq 0$ as follows a vector of length $2\cdot r(n)$ consisting of holomorphic functions on $\mathrm{pr}^{-1}(U_D^{an})$:
\[\tag{\textbf{3.6.2}} p_n^D(z,\tau):=\begin{pmatrix} s^D_0(z,\tau) \\ s^D_1(z,\tau) \\ \vdots\\s^D_n(z,\tau)  \\  0 \\ \vdots \\ \vdots \\ 0 \\-\frac{1}{2\pi i }s^D_1(z,\tau)\\ -\frac{2}{2\pi i }s^D_2(z,\tau) \\ \vdots\\ -\frac{n+1}{2\pi i }s^D_{n+1}(z,\tau) \\ 0 \\ \vdots \\ \vdots \\ 0\end{pmatrix}\]
Here, note that after the first $(n+1)$ entries we fill up with zeroes until we have reached $r(n)$ entries, and we perform the same procedure in the second half of the vector.

\begin{theorem}
(i) With our fixed trivializations for the pullbacks of $\Omega^1_{E^{an}}$, $\mathcal L_n^{an}$ and $\Omega^1_{E^{an}}\otimes_{\mathcal O_{E^{an}}} \mathcal L_n^{an}$ to the universal covering the vector $(3.6.2)$ defines a section in $\Gamma(U_D^{an},\Omega^1_{E^{an}}\otimes_{\mathcal O_{E^{an}}} \mathcal L_n^{an})$.\\
(ii) The section in (i) goes to zero under the map
\[\Gamma(U_D^{an},\Omega^1_{E^{an}}\otimes_{\mathcal O_{E^{an}}} \mathcal L_n^{an}) \xrightarrow{(\nabla_n^{an})^1} \Gamma(U_D^{an},\Omega^2_{E^{an}}\otimes_{\mathcal O_{E^{an}}} \mathcal L_n^{an})\]of the de Rham complex of $(\mathcal L_n^{an},\nabla_n^{an})$.
\end{theorem}
\begin{proof}
(i): Recall from $(3.5.24)$ that the automorphy matrix of $\Omega^1_{E^{an}}$ is given by
\[\Bigg( \begin{pmatrix}  \underline{m} \\ \underline{m} \end{pmatrix}, \begin{pmatrix} a & b \\ c & d \end{pmatrix}, (z,\tau) \Bigg) \mapsto \begin{pmatrix}
                        c\tau+d & 0\\
											 (cz+c\underline{n}-d\underline{m})(c\tau+d) &  (c\tau+d)^2
\end{pmatrix}.
\]
The automorphy matrix for $\Omega^1_{E^{an}}\otimes_{\mathcal O_{E^{an}}} \mathcal L_n^{an}$, evaluated at
\[\Bigg(\begin{pmatrix}  \underline{m} \\ \underline{n} \end{pmatrix}, \begin{pmatrix} a & b \\ c & d \end{pmatrix},(z,\tau)\Bigg) \in \Z^2\times \Gamma(N) \times \C \times \H,\]
is then given (cf. 3.2 (iv)) by the Kronecker product
\[\begin{pmatrix}
                        c\tau+d & 0\\
											 (cz+c\underline{n}-d\underline{m})(c\tau+d) &  (c\tau+d)^2
\end{pmatrix} \otimes A_n=\begin{pmatrix}
                        (c\tau+d)\cdot A_n & \mathbf{0}\\
											 (cz+c\underline{n}-d\underline{m})(c\tau+d)\cdot A_n &  (c\tau+d)^2\cdot A_n
\end{pmatrix}\]
which is a $(2\cdot r(n) \times 2\cdot r(n))$-matrix; note that as usual we leave away the bulky argument of the automorphy matrix $A_n$ for $\mathcal L_n^{an}$.\\
\newline
What we need to show is (cf. 3.2 (ii)) that for each
\[\Bigg(\begin{pmatrix}  \underline{m} \\ \underline{n} \end{pmatrix}, \begin{pmatrix} a & b \\ c & d \end{pmatrix}\Bigg) \in \Z^2\times \Gamma(N)\]
we have the following equation of vectors of holomorphic functions on $\mathrm{pr}^{-1}(U_D^{an})$:
{\small
\[\begin{pmatrix} s^D_0(\frac{z+\underline{m}\tau+\underline{n}}{c\tau+d},\frac{a\tau+b}{c\tau+d}) \\ s^D_1(\frac{z+\underline{m}\tau+\underline{n}}{c\tau+d},\frac{a\tau+b}{c\tau+d}) \\ \vdots\\s^D_n(\frac{z+\underline{m}\tau+\underline{n}}{c\tau+d},\frac{a\tau+b}{c\tau+d}) \\ 0 \\ \vdots \\ \vdots \\ 0  \\ -\frac{1}{2\pi i }s^D_1(\frac{z+\underline{m}\tau+\underline{n}}{c\tau+d},\frac{a\tau+b}{c\tau+d})\\ -\frac{2}{2\pi i }s^D_2(\frac{z+\underline{m}\tau+\underline{n}}{c\tau+d},\frac{a\tau+b}{c\tau+d}) \\ \vdots\\ -\frac{n+1}{2\pi i }s^D_{n+1}(\frac{z+\underline{m}\tau+\underline{n}}{c\tau+d},\frac{a\tau+b}{c\tau+d}) \\ 0 \\ \vdots \\ \vdots \\ 0  \end{pmatrix}=\begin{pmatrix}
                        (c\tau+d)\cdot A_n & \mathbf{0}\\
											 (cz+c\underline{n}-d\underline{m})(c\tau+d)\cdot A_n &  (c\tau+d)^2\cdot A_n \end{pmatrix} \cdot \begin{pmatrix} s^D_0(z,\tau) \\ s^D_1(z,\tau) \\ \vdots\\s^D_n(z,\tau) \\ 0 \\ \vdots \\ \vdots \\ 0 \\ -\frac{1}{2\pi i }s^D_1(z,\tau)\\ -\frac{2}{2\pi i }s^D_2(z,\tau) \\ \vdots\\ -\frac{n+1}{2\pi i }s^D_{n+1}(z,\tau) \\ 0 \\ \vdots \\ \vdots \\ 0\end{pmatrix}.\]}
											
This seems quite an intimidating task to do, but one can reduce the problem step for step as follows.\\
We begin with the equality of the first $r(n)$ rows of the left vector with the corresponding rows of the right vector: this equality is equivalent to the following equation of vectors of length $r(n)$:
\[\begin{pmatrix} s^D_0(\frac{z+\underline{m}\tau+\underline{n}}{c\tau+d},\frac{a\tau+b}{c\tau+d}) \\ s^D_1(\frac{z+\underline{m}\tau+\underline{n}}{c\tau+d},\frac{a\tau+b}{c\tau+d}) \\ \vdots\\s^D_n(\frac{z+\underline{m}\tau+\underline{n}}{c\tau+d},\frac{a\tau+b}{c\tau+d}) \\ 0 \\ \vdots \\ \vdots \\ 0 \end{pmatrix}=((c\tau+d)\cdot A_n)\cdot \begin{pmatrix} s^D_0(z,\tau) \\ s^D_1(z,\tau) \\ \vdots\\s^D_n(z,\tau) \\ 0 \\ \vdots \\ \vdots \\ 0 \end{pmatrix},\]
and this was verified in the proof of Prop. 3.6.1.\\
Next, consider the second half of our big vectors. In this half, take a row where we filled the vector of the left side with a zero entry: such a row indeed equals the corresponding row at the right side because of the arrangement of zeroes in the matrix $A_n$ and in our vectors. This is immediate - for a convenient direct comparison here is again the matrix $A_n$:

{\footnotesize
\[
A_n=
\left( \begin{array}{cccc|cccc}
a_0	& 0 	& \dots  & 0& 0& \dots & \dots & 0 \\
	(c\tau+d)\cdot a_1& (c\tau+d)\cdot a_0	& \dots   &0 &0 & \dots & \dots & 0 \\
	\vdots & \vdots  	&   &  \vdots &   \vdots & & & \vdots \\
	(c\tau+d)^n \cdot a_n& (c\tau+d)^n\cdot a_{n-1}   & \dots     & (c\tau+d)^n\cdot a_0 & 0 & \dots & \dots &0 \\
    \hline
 	0 & 0   & \dots     & 0       & * & \dots & \dots & * \\
	\vdots &  \vdots   &     &\vdots   &\vdots & & & \vdots \\
	\vdots & \vdots & & \vdots & \vdots & & & \vdots \\
	0&   0& \dots    &    0    &* &\dots & \dots & *
\end{array}\right)
\]}What remains are the first $(n+1)$ rows in the second half of our big vectors. The equality of these rows of the left side with the corresponding ones of the right side is equivalent to the following claim, as one sees without difficulty:\\
Namely, that for each $1\leq k \leq n+1$ we have:
\begin{align*}
&-\frac{k}{2\pi i } \cdot s^D_k\bigg(\frac{z+\underline{m}\tau+\underline{n}}{c\tau+d},\frac{a\tau+b}{c\tau+d}\bigg)\\
&=(cz+c\underline{n}-d\underline{m})(c\tau+d)^k \bigg[a_{k-1}s^D_0(z,\tau)+a_{k-2}s^D_1(z,\tau)...+a_0s^D_{k-1}(z,\tau)\bigg]\\
&+(c\tau+d)^{k+1}\bigg[a_{k-1}\bigg(-\frac{1}{2\pi i}s_1^D(z,\tau)\bigg)+a_{k-2}\bigg(-\frac{2}{2\pi i}s_2^D(z,\tau)\bigg)+...+a_0\bigg(-\frac{k}{2\pi i }s_k^D(z,\tau)\bigg)\bigg].
\end{align*}
Recall (cf. $(3.5.11)$ and Def. 3.5.5) that $a_r$ abbreviates the expression
\[a_r\Bigg(\begin{pmatrix}  \underline{m} \\ \underline{n} \end{pmatrix}, \begin{pmatrix} a & b \\ c & d \end{pmatrix}, (z,\tau) \Bigg)=\frac{1}{r!}\bigg(\frac{2\pi i (d\underline{m}-cz-c\underline{n})}{c\tau+d}\bigg)^r=\frac{(-1)^r(2\pi i )^r}{r!}\bigg(\frac{cz+c\underline{n}-d\underline{m}}{c\tau+d}\bigg)^r.\]
Written differently, we have to show that for each $1\leq k \leq n+1$:
\begin{align*}
&s^D_k\bigg(\frac{z+\underline{m}\tau+\underline{n}}{c\tau+d},\frac{a\tau+b}{c\tau+d}\bigg)\\
&=s^D_0(z,\tau)\cdot \bigg[-\frac{2\pi i }{k}(cz+c\underline{n}-d\underline{m})(c\tau+d)^ka_{k-1}\bigg]\\
&+s^D_1(z,\tau)\cdot \bigg[-\frac{2\pi i }{k}(cz+c\underline{n}-d\underline{m})(c\tau+d)^ka_{k-2}+\frac{1}{k}(c\tau+d)^{k+1}a_{k-1}\bigg]\\
&+s^D_2(z,\tau)\cdot \bigg[-\frac{2\pi i }{k}(cz+c\underline{n}-d\underline{m})(c\tau+d)^ka_{k-3}+\frac{2}{k}(c\tau+d)^{k+1}a_{k-2}\bigg]+...\\
&+s^D_{k-1}(z,\tau)\cdot \bigg[-\frac{2\pi i }{k}(cz+c\underline{n}-d\underline{m})(c\tau+d)^ka_0+\frac{k-1}{k}(c\tau+d)^{k+1}a_1\bigg]\\
&+s_k^D(z,\tau)\cdot (c\tau+d)^{k+1}a_0.
\end{align*}
(Multiply the equation we want to show with $-\frac{2\pi i}{k}$ and then order the right side by the $s_r(z,\tau)$.)\\
Using the explicit formula for the $a_r$ we now look at the angled brackets on the right side of the previous equation. They are given as follows:\\
\newline
At $s_0^D(z,\tau)$:
\[\frac{(-1)^k(2\pi i )^k}{k!}(cz+c\underline{n}-d\underline{m})^k(c\tau+d).\]
At $s_j^D(z,\tau)$, for $1\leq j \leq k-1$:
\begin{align*}
&(-1)^{k-j}(2\pi i)^{k-j}(cz+c\underline{n}-d\underline{m})^{k-j}(c\tau+d)^{j+1}\cdot \bigg[\frac{1}{k\cdot (k-j-1)!}+\frac{j}{k\cdot (k-j)!}\bigg]\\
&=\frac{(-1)^{k-j}(2\pi i )^{k-j}}{(k-j)!}(cz+c\underline{n}-d\underline{m})^{k-j}(c\tau+d)^{j+1}.
\end{align*}
At $s_k(z,\tau)$:
\[(c\tau+d)^{k+1}.\]
The desired equality is thus equivalent with the following equation (where again $1\leq k \leq n+1$):
\begin{align*}
s^D_k\bigg(\frac{z+\underline{m}\tau+\underline{n}}{c\tau+d},\frac{a\tau+b}{c\tau+d}\bigg)&=\sum_{j=0}^k\frac{(-1)^{k-j}(2\pi i)^{k-j}}{(k-j)!}(cz+c\underline{n}-d\underline{m})^{k-j}(c\tau+d)^{j+1} \cdot s^D_j(z,\tau)\\
&=(c\tau+d)^{k+1}\sum_{j=0}^k\frac{(-1)^{k-j}(2\pi i)^{k-j}}{(k-j)!}\bigg(\frac{cz+c\underline{n}-d\underline{m}}{c\tau+d}\bigg)^{k-j} \cdot s^D_j(z,\tau)\\
&=(c\tau+d)^{k+1}\sum_{j=0}^ka_{k-j}\cdot s^D_j(z,\tau).
\end{align*}
This is true - cf. the equation $(*)$ in the proof of Prop. 3.6.1.\\
(ii): By Prop. 3.5.22 we only need to show that the functions $s^D_k(z,\tau)$ obey the differential equation
\[\partial_{\tau}s^D_k(z,\tau)=-\frac{k+1}{2\pi i }\cdot  \partial_zs^D_{k+1}(z,\tau) \qquad \textrm{for all} \ k=0,...,n.\]
But it is well-known (cf. \cite{Le-Ra}, $(14)$\footnote{Note that the function which they consider in this formula, written in the argument $(z,w,\tau)$, is equal to
\[2\pi i \cdot F(2\pi i z,2\pi i w,\tau),\]
where $F$ denotes again the function of \cite{Za2}, 3: to see this one may use the expression of $F$ in terms of a double series as given at the beginning of the proof of \cite{Za2}, 3, Theorem. By Prop. 3.3.13 the function considered in \cite{Le-Ra}, $(14)$, hence coincides precisely with what we call the fundamental meromorphic Jacobi form $J(z,w,\tau)$.}) resp. straightforwardly deduced from Prop. 3.3.13 and the series expression given in \cite{Za2}, 3, that the fundamental meromorphic Jacobi form satisfies the "mixed heat equation":
\[2\pi i \cdot \partial_{\tau} J(z,w,\tau)=\partial_z \partial_w J(z,w,\tau).\]From this we obtain the equality
\[\partial_{\tau}\Big[D^2J(z,w,\tau)-D\cdot J \Big(Dz,\frac{w}{D},\tau \Big) \Big]=\frac{1}{2\pi i} \partial_z\partial_w \Big[D^2J(z,w,\tau)-D\cdot J \Big(Dz,\frac{w}{D},\tau \Big) \Big],\]
which by consideration of Laurent expansions around $w=0$ yields
\[(-1)^k\cdot \partial_{\tau}s_k^D(z,\tau)=(-1)^{k+1} \cdot \frac{k+1}{2\pi i}\cdot \partial_zs^D_{k+1}(z,\tau)\]
(observe the definition of the $s_k(z,\tau)$ with the sign in front of $w$ and $\frac{w}{D}$). Dividing by $(-1)^k$ gives the desired differential equation.
\end{proof}

Let us finally record the following natural compatibilities.

\begin{lemma}
For all $n\geq 0$ the canonical arrow
\[\Gamma(U_D^{an},\Omega^1_{E^{an}}\otimes_{\mathcal O_{E^{an}}} \mathcal L_n^{an})\rightarrow \Gamma(U_D^{an},\Omega^1_{E^{an}/S^{an}}\otimes_{\mathcal O_{E^{an}}} \mathcal L_n^{an})\]maps the section $p_n^D(z,\tau)$ to the section $q_n^D(z,\tau)$.
\end{lemma}
\begin{proof}
This is immediately clear from how the occurring vector bundles are trivialized on the universal covering of $E^{an}$.
\end{proof}

\begin{lemma}
For all $n\geq 1$ the section $q_n^D(z,\tau)$ maps to the section $q_{n-1}^D(z,\tau)$ under the canonical arrow
\[\Gamma(U_D^{an},\Omega^1_{E^{an}/S^{an}}\otimes_{\mathcal O_{E^{an}}} \mathcal L_n^{an})\rightarrow \Gamma(U_D^{an},\Omega^1_{E^{an}/S^{an}}\otimes_{\mathcal O_{E^{an}}} \mathcal L_{n-1}^{an})\]
induced by the analytified transition map $\mathcal L_n \rightarrow \mathcal L_{n-1}$.\\
The analogous assertion holds for $p_n^D(z,\tau)$ and the canonical arrow
\[\Gamma(U_D^{an},\Omega^1_{E^{an}}\otimes_{\mathcal O_{E^{an}}} \mathcal L_n^{an})\rightarrow \Gamma(U_D^{an},\Omega^1_{E^{an}}\otimes_{\mathcal O_{E^{an}}} \mathcal L_{n-1}^{an}).\]
\end{lemma}
\begin{proof}
For each $n\geq 1$ the pullback of the analytified transition map $\mathcal L_n\rightarrow \mathcal L_{n-1}$ to $\C\times \H$ is given as the composition
\[\mathrm{Sym}^n_{\mathcal O_{\C\times \H}}(\mathrm{pr}^*\mathcal L_1^{an})\rightarrow \mathrm{Sym}^n_{\mathcal O_{\C\times \H}}(\mathrm{pr}^*\mathcal L_1^{an}\oplus \mathcal O_{\C\times \H}) \rightarrow \mathrm{Sym}^{n-1}_{\mathcal O_{\C\times \H}}(\mathrm{pr}^*\mathcal L_1^{an})\]
in which the first morphism is induced by the identity and the projection in $(3.5.8)$ and the second comes from the decomposition of the symmetric power of a direct sum.\\
For $0 \leq i,j$ with $i+j\leq n$ this composition maps the basic section
\[\frac{e^{n-i-j}f^ig^j}{(n-i-j)!}\]
of $\mathrm{Sym}^n_{\mathcal O_{\C\times \H}}(\mathrm{pr}^*\mathcal L_1^{an})$ to the basic section
\[\frac{e^{n-i-j-1}f^ig^j}{(n-i-j-1)!}\]
of $\mathrm{Sym}^{n-1}_{\mathcal O_{\C\times \H}}(\mathrm{pr}^*\mathcal L_1^{an})$, where terms with $e^{-1}$ are understood to be zero; this follows easily by recalling that the projection
\[\mathrm{pr}^*\mathcal L_1^{an} \rightarrow \mathcal O_{\C\times \H}\]
in $(3.5.8)$ sends the sections $ \{e,f,g \}$ to $\{ 1,0,0\}$ (cf. the explanations preceding Prop. 3.5.2).\\
With this the claims of the lemma are clear by the shape of the vectors $(3.6.1)$ and $(3.6.2)$ and by the definition of the trivialization for the pullback of the respective vector bundles to $\C\times \H$.
\end{proof}

\subsection{The pullback along torsion sections}
We now assume that the integer $D>1$ in addition satisfies $(D,N)=1$.\\
For two fixed integers $a,b$ which are not simultaneously divisible by $N$ the $N$-torsion section
\[t_{a,b}=ae_1+be_2: S \rightarrow E\]
by our assumptions on $D,a,b$ factors over the open subscheme $U_D=E - E[D]$:
\[t_{a,b}: S \rightarrow U_D \subseteq E.\]
On the analytic side we work as usual on a fixed connected component and again adopt the already familiar notational conventions formulated explicitly at the beginning of 3.5 or 3.5.3. We let $j_0$ be an integer with $(j_0,N)=1$ whose class in $(\Z/N\Z)^*$ determines the chosen component.\\
\newline
The analytification of $t_{a,b}$ then expresses (according to $(3.4.1)$-$(3.4.4)$) as the map
\[t_{a,b}^{an}: \Gamma(N) \backslash \H =S^{an} \rightarrow U_D^{an}\subseteq E^{an}=(\Z^2\times \Gamma(N)) \backslash (\C\times \H),\quad \tau\mapsto \Big(\frac{aj_0\tau}{N}+\frac{b}{N},\tau\Big).\]
An element in
\[\Gamma(U_D^{an},\Omega^1_{E^{an}}\otimes_{\mathcal O_{E^{an}}} \mathcal L_n^{an})\]
induces a section (its "specialization along $t_{a,b}^{an}$") in
\[\Gamma \Big(S^{an},\Omega^1_{S^{an}}\otimes_{\mathcal O_{S^{an}}} \prod_{k=0}^n \mathrm{Sym}^k_{\mathcal O_{S^{an}}}\mathcal H^{an} \Big),\]
obtained from its pullback along $t_{a,b}^{an}$ and using the composition
\[(t_{a,b}^{an})^*(\Omega^1_{E^{an}}\otimes_{\mathcal O_{E^{an}}} \mathcal L_n^{an}) \simeq (t_{a,b}^{an})^*\Omega^1_{E^{an}}\otimes_{\mathcal O_{S^{an}}} \prod_{k=0}^n \mathrm{Sym}^k_{\mathcal O_{S^{an}}}\mathcal H^{an}  \xrightarrow{\mathrm{can}} \Omega^1_{S^{an}} \otimes_{\mathcal O_{S^{an}}} \prod_{k=0}^n \mathrm{Sym}^k_{\mathcal O_{S^{an}}}\mathcal H^{an};\]
the isomorphism in the previous chain is induced by the identification $(3.5.13)$ and the splitting $\varphi_n^{an}$:
\[(t_{a,b}^{an})^*\mathcal L_n^{an}\simeq (\epsilon^{an})^*\mathcal L_n^{an} \simeq \prod_{k=0}^n \mathrm{Sym}^k_{\mathcal O_{S^{an}}}\mathcal H^{an}.\]
For the following result recall that the pullback of
\[\Omega^1_{E^{an}}\otimes_{\mathcal O_{E^{an}}} \mathcal L_n^{an}\]
resp. of
\[\Omega^1_{S^{an}}\otimes_{\mathcal O_{S^{an}}} \prod_{k=0}^n \mathrm{Sym}^k_{\mathcal O_{S^{an}}}\mathcal H^{an}\]
to $\C\times \H$ resp. to $\H$ has a fixed trivialization - as explained in 3.2 (iv) together with $(3.5.24)$ and Def. 3.5.4 resp. together with $(3.5.27)$ and the remarks preceding 3.5.1.

\begin{theorem}
Let $n\geq 0$. Specializing along $t_{a,b}^{an}$ the section in $\Gamma(U_D^{an},\Omega^1_{E^{an}}\otimes_{\mathcal O_{E^{an}}} \mathcal L_n^{an})$ defined by the vector $p^D_n(z,\tau)$ of $(3.6.2)$ (cf. Thm. 3.6.2 (i)) yields the section in $\Gamma(S^{an},\Omega^1_{S^{an}}\otimes_{\mathcal O_{S^{an}}} \prod_{k=0}^n \mathrm{Sym}^k_{\mathcal O_{S^{an}}}\mathcal H^{an})$ which is given by
\[\begin{pmatrix} -2\pi i \cdot \ _DF^{(2)}_{\frac{a j_0}{N},\frac{b}{N}}(\tau) \\ \vdots\\ \frac{(-1)^k(2\pi i )^k}{(k-1)!}\cdot \ _DF^{(k+1)}_{\frac{a j_0}{N},\frac{b}{N}}(\tau)\\ \vdots \\ \frac{(-1)^{n+1}(2\pi i)^{n+1}}{n!}\cdot \ _DF^{(n+2)}_{\frac{a j_0}{N},\frac{b}{N}}(\tau) \\ 0 \\ \vdots \\ \vdots \\ 0 \end{pmatrix},\]
i.e. by the vector of length $r(n)$ of holomorphic functions on $\H$ whose $k$-th entry consists of the modular form $\frac{(-1)^k(2\pi i )^k}{(k-1)!}\cdot \ _DF^{(k+1)}_{\frac{a j_0}{N},\frac{b}{N}}(\tau)$ if $1\leq k \leq n+1$ and whose further entries are zero.\\
Cf. Def. 3.3.17 for the definition of the function $_DF^{(k+1)}_{\frac{a j_0}{N},\frac{b}{N}}(\tau)$.
\end{theorem}

\begin{proof}
We apply Prop. 3.5.14. To connect the present situation with the notation used there we set

\[\tag{$*$}\begin{pmatrix} l_{0,0}(z,\tau) \\ l_{1,0}(z,\tau) \\ \vdots \\ l_{n,0}(z,\tau)\\ \vdots\\ \vdots\\l_{0,n}(z,\tau)  \\  \lambda_{0,0}(z,\tau) \\ \lambda_{1,0}(z,\tau) \\ \vdots \\ \lambda_{n,0}(z,\tau)\\ \vdots\\ \vdots\\ \lambda_{0,n}(z,\tau) \end{pmatrix} := \footnotesize{\begin{pmatrix}s^D_0(z,\tau) \\ s^D_1(z,\tau) \\ \vdots\\s^D_n(z,\tau)  \\  0 \\ \vdots \\ \vdots \\ 0 \\-\frac{1}{2\pi i }s^D_1(z,\tau)\\ -\frac{2}{2\pi i }s^D_2(z,\tau) \\ \vdots\\ -\frac{n+1}{2\pi i }s^D_{n+1}(z,\tau) \\ 0 \\ \vdots \\ \vdots \\ 0\end{pmatrix}}.\]
According to Prop. 3.5.14 the specialization of this section along $t_{a,b}^{an}$ is given by
\[\tag{$**$} \begin{pmatrix} \frac{aj_0}{N}\cdot \widehat{l}_{0,0}(0,\tau) + \widehat{\lambda}_{0,0}(0,\tau)\\ \frac{aj_0}{N}\cdot \widehat{l}_{1,0}(0,\tau) +\widehat{\lambda}_{1,0}(0,\tau) \\ \vdots \\ \frac{aj_0}{N}\cdot \widehat{l}_{n,0}(0,\tau)+\widehat{\lambda}_{n,0}(0,\tau)\\ \vdots\\ \vdots\\ \frac{aj_0}{N}\cdot \widehat{l}_{0,n}(0,\tau)+\widehat{\lambda}_{0,n}(0,\tau)\end{pmatrix},\]
where the functions $\widehat{l}_{i,j}(0,\tau)$ resp. $\widehat{\lambda}_{i,j}(0,\tau)$ are defined by the formulas
\[\widehat{l}_{i,j}(0,\tau)=\sum_{r=0}^i\frac{(-2\pi i \frac{aj_0}{N})^{i-r}}{(i-r)!}\cdot l_{r,j}\Big(\frac{aj_0\tau}{N}+\frac{b}{N},\tau \Big)\]
resp.
\[\widehat{\lambda}_{i,j}(0,\tau)=\sum_{r=0}^i\frac{(-2\pi i \frac{aj_0}{N})^{i-r}}{(i-r)!}\cdot \lambda_{r,j}\Big(\frac{aj_0\tau}{N}+\frac{b}{N},\tau \Big).\]
Having a look at $(*)$ it is obvious that
\[\widehat{l}_{i,j}(0,\tau)= 0 =\widehat{\lambda}_{i,j}(0,\tau) \qquad \textrm{for all} \ 0<j\leq n,\]
such that $(**)$ writes as
\[\begin{pmatrix} \frac{aj_0}{N}\cdot \widehat{l}_{0,0}(0,\tau) + \widehat{\lambda}_{0,0}(0,\tau)\\ \frac{aj_0}{N}\cdot \widehat{l}_{1,0}(0,\tau) +\widehat{\lambda}_{1,0}(0,\tau) \\ \vdots\\ \frac{aj_0}{N}\cdot \widehat{l}_{n,0}(0,\tau)+\widehat{\lambda}_{n,0}(0,\tau) \\ 0 \\ \vdots \\ \vdots \\ 0 \end{pmatrix}.\]

Now consider the $k$-th row of this vector for $1 \leq k \leq n+1$, i.e. the function
\[\frac{aj_0}{N}\cdot \widehat{l}_{k-1,0}(0,\tau)+\widehat{\lambda}_{k-1,0}(0,\tau).\]
The theorem is proven if we can show that this is equal to the function $\frac{(-1)^k(2\pi i )^k}{(k-1)!}\cdot \ _DF^{(k+1)}_{\frac{a j_0}{N},\frac{b}{N}}(\tau)$.\\
\newline
We have
\begin{align*}
&\frac{aj_0}{N}\cdot \widehat{l}_{k-1,0}(0,\tau)+\widehat{\lambda}_{k-1,0}(0,\tau)\\
&=\frac{aj_0}{N}\cdot \sum_{r=0}^{k-1} \frac{(-2\pi i \frac{aj_0}{N})^{k-1-r}}{(k-1-r)!}\cdot l_{r,0}\Big(\frac{aj_0\tau}{N}+\frac{b}{N},\tau\Big)+\sum_{r=0}^{k-1}\frac{(-2\pi i \frac{aj_0}{N})^{k-1-r}}{(k-1-r)!}\cdot \lambda_{r,0}\Big(\frac{aj_0\tau}{N}+\frac{b}{N},\tau\Big)\\
&=\frac{aj_0}{N}\cdot \sum_{r=0}^{k-1} \frac{(-2\pi i \frac{aj_0}{N})^{k-1-r}}{(k-1-r)!}\cdot s_r^D\Big(\frac{aj_0\tau}{N}+\frac{b}{N},\tau\Big)\\
&+\sum_{r=0}^{k-1}\frac{(-2\pi i \frac{aj_0}{N})^{k-1-r}}{(k-1-r)!}\cdot \bigg(-\frac{r+1}{2\pi i }\bigg)\cdot s^D_{r+1}\Big(\frac{aj_0\tau}{N}+\frac{b}{N},\tau\Big)\\
&=s^D_0\Big(\frac{aj_0\tau}{N}+\frac{b}{N},\tau\Big) \cdot \bigg[\frac{(-1)^{k-1}(2\pi i )^{k-1}}{(k-1)!}\bigg(\frac{aj_0}{N}\bigg)^k\bigg]\\
&+\sum_{r=1}^{k-1}s^D_r\Big(\frac{aj_0\tau}{N}+\frac{b}{N},\tau\Big) \cdot \bigg[\frac{(-1)^{k-1-r}(2\pi i )^{k-1-r}}{(k-1-r)!}\bigg(\frac{aj_0}{N}\bigg)^{k-r} +\frac{(-1)^{k-1-r}(2\pi i )^{k-1-r}}{(k-r)!}\bigg(\frac{aj_0}{N}\bigg)^{k-r} \cdot r \bigg]\\
&-s^D_k\Big(\frac{aj_0\tau}{N}+\frac{b}{N},\tau\Big)\cdot \frac{k}{2\pi i }\\
&=s^D_0\Big(\frac{aj_0\tau}{N}+\frac{b}{N},\tau\Big) \cdot \bigg[\frac{(-1)^{k-1}(2\pi i )^{k-1}}{(k-1)!}\bigg(\frac{aj_0}{N}\bigg)^k\bigg]\\
&+\sum_{r=1}^{k-1}s^D_r\Big(\frac{aj_0\tau}{N}+\frac{b}{N},\tau\Big) \cdot \bigg[\frac{(-1)^{k-1-r}(2\pi i )^{k-1-r}}{(k-r)!}\bigg(\frac{aj_0}{N}\bigg)^{k-r}\cdot k\bigg]-s^D_k\Big(\frac{aj_0\tau}{N}+\frac{b}{N},\tau\Big) \cdot \frac{k}{2\pi i }.
\end{align*}

Next, we combine the defining formula of the $s_r(z,\tau)$:
\[D^2 \cdot J(z,-w,\tau)-D\cdot J\Big(Dz,-\frac{w}{D},\tau \Big)=s^D_0(z,\tau)+s^D_1(z,\tau)w+...\]
with the $w$-expansion formula $(3.3.25)$:
\[J(z,w,\tau)=\exp \bigg[2\pi i \frac{\bar{z}w-zw}{\tau-\bar{\tau}} \bigg]\cdot \bigg( \frac{1}{w}+\sum_{k\geq 0 }(-1)^k \cdot e_{k+1}(z,\tau) \cdot w^k\bigg)\]
and obtain
\[s^D_0(z,\tau)+s^D_1(z,\tau)w+...=\exp \bigg[2\pi i \frac{zw-\bar{z}w}{\tau-\bar{\tau}} \bigg]\cdot \sum_{k\geq 0 }\bigg[D^2e_{k+1}(z,\tau)-D^{1-k}e_{k+1}(Dz,\tau)\bigg]\cdot w^k,\]
such that
\begin{align*}
&\sum_{k\geq 0 }\bigg[D^2e_{k+1}\Big(\frac{aj_0\tau}{N}+\frac{b}{N},\tau\Big)-D^{1-k}e_{k+1}\Big(D\frac{aj_0\tau}{N}+D\frac{b}{N},\tau\Big)\bigg]\cdot w^k\\
&=\exp \bigg[-2\pi i \frac{aj_0}{N}w\bigg] \cdot \bigg[s^D_0\Big(\frac{aj_0\tau}{N}+\frac{b}{N},\tau\Big)+s^D_1\Big(\frac{aj_0\tau}{N}+\frac{b}{N},\tau\Big)w+...\bigg].
\end{align*}
This implies that for each $k\geq 0$ we have
\begin{align*}
&D^2e_{k+1}\Big(\frac{aj_0\tau}{N}+\frac{b}{N},\tau\Big)-D^{1-k}e_{k+1}\Big(D\frac{aj_0\tau}{N}+D\frac{b}{N},\tau\Big)\\
&=s^D_0\Big(\frac{aj_0\tau}{N}+\frac{b}{N},\tau\Big)\cdot \frac{(-2\pi i \frac{aj_0}{N})^k}{k!}+s^D_1\Big(\frac{aj_0\tau}{N}+\frac{b}{N},\tau\Big)\cdot\frac{(-2\pi i \frac{aj_0}{N})^{k-1}}{(k-1)!}+...\\
&+s^D_{k-1}\Big(\frac{aj_0\tau}{N}+\frac{b}{N},\tau\Big)\cdot \Big(-2\pi i \frac{aj_0}{N}\Big)+s^D_k\Big(\frac{aj_0\tau}{N}+\frac{b}{N},\tau\Big),
\end{align*}

and hence

\begin{align*}
&-\frac{k}{2\pi i }\cdot \bigg[D^2e_{k+1}\Big(\frac{aj_0\tau}{N}+\frac{b}{N},\tau\Big)-D^{1-k}e_{k+1}\Big(D\frac{aj_0\tau}{N}+D\frac{b}{N},\tau\Big)\bigg] \\
&=s^D_0\Big(\frac{aj_0\tau}{N}+\frac{b}{N},\tau\Big) \cdot \bigg[\frac{(-1)^{k-1}(2\pi i )^{k-1}}{(k-1)!}\bigg(\frac{aj_0}{N}\bigg)^k\bigg]\\
&+\sum_{r=1}^{k-1}s^D_r\Big(\frac{aj_0\tau}{N}+\frac{b}{N},\tau\Big) \cdot \bigg[\frac{(-1)^{k-1-r}(2\pi i )^{k-1-r}}{(k-r)!}\bigg(\frac{aj_0}{N}\bigg)^{k-r}\cdot k\bigg]-s^D_k\Big(\frac{aj_0\tau}{N}+\frac{b}{N},\tau\Big) \cdot \frac{k}{2\pi i }.
\end{align*}

We have calculated above that this last expression for $1\leq k \leq n+1$ equals the function
\[\frac{aj_0}{N}\cdot \widehat{l}_{k-1,0}(0,\tau)+\widehat{\lambda}_{k-1,0}(0,\tau),\]
such that we have shown
\[\frac{aj_0}{N}\cdot \widehat{l}_{k-1,0}(0,\tau)+\widehat{\lambda}_{k-1,0}(0,\tau)=-\frac{k}{2\pi i }\cdot \bigg[D^2e_{k+1}\Big(\frac{aj_0\tau}{N}+\frac{b}{N},\tau\Big)-D^{1-k}e_{k+1}\Big(D\frac{aj_0\tau}{N}+D\frac{b}{N},\tau\Big)\bigg].\]
Thm. 3.3.16 tells us that
\begin{align*}
e_{k+1}\Big(\frac{aj_0}{N}\tau+\frac{b}{N},\tau\Big)&=\frac{(-1)^{k+1}(2\pi i)^{k+1}}{k!}\cdot F^{(k+1)}_{\frac{aj_0}{N},\frac{b}{N}}(\tau)\\
e_{k+1}\Big(D\frac{aj_0}{N}\tau+D\frac{b}{N},\tau\Big)&=\frac{(-1)^{k+1}(2\pi i)^{k+1}}{k!}\cdot F^{(k+1)}_{\frac{Daj_0}{N},\frac{Db}{N}}(\tau).
\end{align*}
Combining this with the previous equation we obtain
\begin{align*}
&\frac{aj_0}{N}\cdot \widehat{l}_{k-1,0}(0,\tau)+\widehat{\lambda}_{k-1,0}(0,\tau)=\frac{(-1)^k(2\pi i )^k}{(k-1)!}\cdot \bigg[D^2F^{(k+1)}_{\frac{aj_0}{N},\frac{b}{N}}(\tau)-D^{1-k}F^{(k+1)}_{\frac{Daj_0}{N},\frac{Db}{N}}(\tau)\bigg]\\
&=\frac{(-1)^k(2\pi i )^k}{(k-1)!}\cdot \ _DF^{(k+1)}_{\frac{aj_0}{N},\frac{b}{N}}(\tau),
\end{align*}
the last equality by definition. As was already observed this proves the theorem.
\end{proof}

\section{The analytic characterization result}
\markright{\uppercase{The explicit description on the universal elliptic curve}}
For the whole section fix an integer $D>1$.\\
\newline
As in 1.5.2 write $U_D$ for the open complement of the $D$-torsion subscheme $E[D]$ in $E$ as well as $i_D: E[D]\rightarrow E$ resp. $j_D: U_D \rightarrow E$ for the canonical closed resp. open immersion and $\pi_{E[D]}: E[D]\rightarrow S$ resp. $\pi_{U_D}: U_D \rightarrow S$ for the structure map of $E[D]$ resp. $U_D$ over $S$.
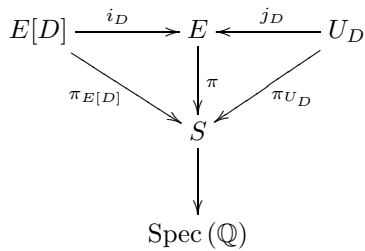
\begin{equation*}
\begin{xy}
\xymatrix{
 & E[D] \ar[dr]_{\pi_{E[D]}} \ar[r]^{i_D} & E  \ar[d]^{\pi} & U_D \ar[dl]^{\pi_{U_D}} \ar[l]_{j_D}\\
  &  & S \ar[d] & &\\
		&  & \Spec(\Q) & &}
\end{xy}
\end{equation*}
In 1.5.2 we saw that there exists an inverse system of cohomology classes
\[\tag{\textbf{3.7.1}} \varD = \Big(\varDn \Big)_{n\geq 0} \in \lim_{n\geq 0} H^1_{\mathrm{dR}}(U_D/\Q, \mathcal L_n)\]
satisfying a certain residue condition, called the $D$-variant of the elliptic polylogarithm for $E/S/\Q$.
\subsection{Preliminaries: Algebraic}

We recall the residue condition by which the system $(3.7.1)$ is uniquely determined. This can be done in slightly simpler terms than in 1.5.2, where more technicalities were needed to establish various isomorphisms by which we could ensure the existence of such a system.\\
Subsequently, we describe in more concrete terms the two arrows appearing in this characterization.

\subsubsection{Algebraic characterization of the $D$-variant of the polylogarithm}

Recall that for each $n\geq 0$ we have a chain of $\mathcal D_{E[D]/\Q}$-linear maps
\[\tag{\textbf{3.7.2}} \mathcal O_{E[D]} \hookrightarrow (\pi_{E[D]})^* \prod_{k=0}^n \mathrm{Sym}^k_{\mathcal O_S}\mathcal H \simeq (\pi_{E[D]})^*\epsilon^*\mathcal L_n \simeq i_D^*[D]^*\mathcal L_n \simeq i_D^*\mathcal L_n\]
in which the monomorphism is given by $\frac{1}{n!}$-times the obvious inclusion, the first isomorphism is the splitting of $\epsilon^*\mathcal L_n$, the second uses the commutative diagram defining $E[D]$ and the last comes from the invariance isomorphism of Cor. 1.4.4. Taking global sections induces an injection
\[\tag{\textbf{3.7.3}} H^0(E[D],\mathcal O_{E[D]}) \hookrightarrow \lim_{n \geq 0}H^0(E[D],i_D^*\mathcal L_n)\]
by which we may view the global section $D^2\cdot 1_{ \{\epsilon \}} - 1_{E[D]}$ of $\mathcal O_{E[D]}$ as an element of the right side.\\
\newline
Now observe that we have a natural map
\[\tag{\textbf{3.7.4}} \lim_{n\geq 0} H^1_{\mathrm{dR}}(U_D/\Q,\mathcal L_n) \rightarrow \lim_{n\geq 0} \Gamma(S,H^1_{\mathrm{dR}}(U_D/S,\mathcal L_n)),\]
given by taking the inverse limit in the composition
\[\begin{split}
&H^1_{\mathrm{dR}}(U_D/\Q,\mathcal L_n) \rightarrow H^0_{\mathrm{dR}}(S/\Q, H^1_{\mathrm{dR}}(U_D/S,\mathcal L_n))\\
&\simeq \Hom_{\mathcal D_{S/\Q}}(\mathcal O_S, H^1_{\mathrm{dR}}(U_D/S,\mathcal L_n)) \subseteq \Gamma(S,H^1_{\mathrm{dR}}(U_D/S,\mathcal L_n)).
\end{split}\]
The first arrow is an edge morphism in the obvious Leray spectral sequence and - as was observed in 1.5.2 - becomes bijective in the inverse limit; in particular, we see that $(3.7.4)$ is injective.\\
\newline
Moreover, the exact sequence of $\mathcal D_{S/\Q}$-modules (coming from the localization triangle for $\mathcal L_n$)

\[\tag{\textbf{3.7.5}} 0 \rightarrow H^1_{\mathrm{dR}}(E/S, \mathcal L_n) \xrightarrow{\mathrm{can}} H^1_{\mathrm{dR}}(U_D/S, \mathcal L_n) \xrightarrow{\mathrm{Res}^D_n} (\pi_{E[D]})_* i_D^*\mathcal L_n \xrightarrow{\sigma_n^D} H^2_{\mathrm{dR}}(E/S, \mathcal L_n) \rightarrow 0
\]
induces the arrow
\[\tag{\textbf{3.7.6}} \lim_{n\geq 0} \Gamma(S,H^1_{\mathrm{dR}}(U_D/S,\mathcal L_n)) \xrightarrow{\mathrm{Res}^D} \lim_{n\geq 0}H^0(E[D],i_D^*\mathcal L_n)\]
which with Thm. 1.2.1 (ii) is seen to be injective.\\
\newline
The composition of the injections $(3.7.4)$ and $(3.7.6)$ gives the injective $\Q$-linear map
\[\tag{\textbf{3.7.7}} \lim_{n \geq 0} H^1_{\mathrm{dR}}(U_D/\Q,\mathcal L_n) \rightarrow \lim_{n \geq 0} \Gamma(S,H^1_{\mathrm{dR}}(U_D/S,\mathcal L_n)) \xrightarrow{\mathrm{Res}^D} \lim_{n\geq 0}H^0(E[D],i_D^*\mathcal L_n),\]
which by definition of the $D$-variant sends $\varD$ to the element $D^2\cdot 1_{ \{\epsilon \}} - 1_{E[D]}$.\\
\newline
For all further proceeding it is essential to understand the two arrows in the composition $(3.7.7)$ as explicitly as possible. This is what we undertake next.
\subsubsection{Description of the first arrow}
The first map of $(3.7.7)$ was defined in $(3.7.4)$ using the a priori rather abstract edge morphism of the Leray spectral sequence. But it has an easy natural interpretation as follows:\\
\newline
Namely, as $\pi_{U_D}:U_D \rightarrow S$ is affine and the modular curve $S$ is affine the scheme $U_D$ is affine; as it is moreover $2$-dimensional and smooth over $\Q$ we see that the spectral sequence (cf. $(0.2.2)$)
\[E_1^{p,q}=H^q(U_D,\Omega^p_{U_D/\Q} \otimes_{\mathcal O_{U_D}} \mathcal L_n) \Rightarrow E^{p+q}=H^{p+q}_{\mathrm{dR}}(U_D/\Q,\mathcal L_n)\]
degenerates at $r=2$ and that the edge morphism $E_2^{1,0} \rightarrow E^1$ on the second sheet is an isomorphism:
\[\tag{\textbf{3.7.8}} \frac{\ker\bigg(\Gamma(U_D, \Omega^1_{U_D/\Q} \otimes_{\mathcal O_{U_D}} \mathcal L_n)\xrightarrow{\nabla_n^1}\Gamma(U_D, \Omega^2_{U_D/\Q} \otimes_{\mathcal O_{U_D}} \mathcal L_n)\bigg)}{\mathrm{im}\bigg(\Gamma(U_D,\mathcal L_n) \xrightarrow{\nabla_n} \Gamma(U_D,\Omega^1_{U_D/\Q} \otimes_{\mathcal O_{U_D}} \mathcal L_n)\bigg)} \xrightarrow{\sim} H^1_{\mathrm{dR}}(U_D/\Q,\mathcal L_n).\]
Similarly, as $\pi_{U_D}: U_D \rightarrow S$ is affine and smooth of relative dimension one the spectral sequence (cf. again $(0.2.2)$)
\[E_1^{p,q}=R^q(\pi_{U_D})_*(\Omega^p_{U_D/S} \otimes_{\mathcal O_{U_D}} \mathcal L_n) \Rightarrow E^{p+q}=H^{p+q}_{\mathrm{dR}}(U_D/S,\mathcal L_n)\]
degenerates at $r=2$ and the edge morphism $E_2^{1,0} \rightarrow E^1$ is an isomorphism; in global $S$-sections (recall that $S$ is affine) it reads as:
\[\tag{\textbf{3.7.9}} \frac{\Gamma(U_D, \Omega^1_{U_D/S} \otimes_{\mathcal O_{U_D}} \mathcal L_n)}{\mathrm{im}\bigg(\Gamma(U_D,\mathcal L_n) \xrightarrow{\nabla_n^{res}} \Gamma(U_D,\Omega^1_{U_D/S} \otimes_{\mathcal O_{U_D}} \mathcal L_n)\bigg)} \xrightarrow{\sim} \Gamma(S,H^1_{\mathrm{dR}}(U_D/S,\mathcal L_n)).\]
It is straightforward and purely formal to check that the arrow $(3.7.4)$:
\[\lim_{n\geq 0} H^1_{\mathrm{dR}}(U_D/\Q,\mathcal L_n) \rightarrow \lim_{n\geq 0} \Gamma(S,H^1_{\mathrm{dR}}(U_D/S,\mathcal L_n))\]
is given under the previous identifications $(3.7.8)$ resp. $(3.7.9)$ simply as the map induced by the natural restriction of differential forms:
\[\tag{\textbf{3.7.10}} \Gamma(U_D,\Omega^1_{U_D/\Q} \otimes_{\mathcal O_{U_D}} \mathcal L_n) \rightarrow \Gamma(U_D,\Omega^1_{U_D/S} \otimes_{\mathcal O_{U_D}} \mathcal L_n).\]
This is what we need to know about the first arrow in $(3.7.7)$.
\subsubsection{Description of the second arrow}
The second map $\mathrm{Res}^D$ of $(3.7.7)$ arises - by taking global $S$-sections and inverse limit - from the a priori rather abstract arrows (for all $n\geq 0$):
\[H^1_{\mathrm{dR}}(U_D/S, \mathcal L_n) \xrightarrow{\mathrm{Res}^D_n} H^0_{\mathrm{dR}}(E[D]/S,i_D^*\mathcal L_n) = (\pi_{E[D]})_* i_D^*\mathcal L_n,\]
defined by applying $\pi_+$ to the localization triangle in $D^b_{\mathrm{qc}}(\mathcal D_{E/\Q})$:
\[\mathcal L_n \rightarrow (j_D)_+ \mathcal L_{n|U_D} \rightarrow (i_D)_+ (i_D)^* \mathcal L_n\]
and then going into the long exact sequence of cohomology:
\[0 \rightarrow H^1_{\mathrm{dR}}(E/S, \mathcal L_n) \xrightarrow{\mathrm{can}} H^1_{\mathrm{dR}}(U_D/S, \mathcal L_n) \xrightarrow{\mathrm{Res}^D_n} H^0_{\mathrm{dR}}(E[D]/S,i_D^*\mathcal L_n) \rightarrow H^2_{\mathrm{dR}}(E/S, \mathcal L_n) \rightarrow 0.\]

Clearly, completely the same procedure can be applied to any vector bundle $\mathcal V$ on $E$ with integrable $\Q$-connection, yielding an exact sequence of vector bundles with integrable $\Q$-connection on $S$:

\[\tag{\textbf{3.7.11}} 0 \rightarrow H^1_{\mathrm{dR}}(E/S, \mathcal V) \xrightarrow{\mathrm{can}} H^1_{\mathrm{dR}}(U_D/S, \mathcal V) \xrightarrow{\mathrm{Res}} H^0_{\mathrm{dR}}(E[D]/S,i_D^*\mathcal V) \rightarrow H^2_{\mathrm{dR}}(E/S, \mathcal V) \rightarrow 0.\]

A different interpretation of the previous sequence $(3.7.11)$ together with an explicit description of the map $\mathrm{Res}$ is obtained by using logarithmic differential forms. Let us explain this:\\
\newline
$E[D]$ is an effective relative Cartier divisor of $E/S$ and smooth over the smooth variety $S$. In particular, it has (strict) relative normal crossings over $S$ and permits defining the sheaf $\Omega^1_{E/S}(\mathrm{log}(E[D]))$ of relative differential forms with logarithmic poles along $E[D]$ (cf. \cite{De1}, Ch. II, $(3.3.1)$).\\
We only give a brief reminder of how this sheaf looks like locally and refer to \cite{De1}, Ch. II, §3, and to \cite{Kat1}, 1.0, for more details concerning logarithmic differential forms.\\
Recall hence that we may cover $E$ by open subsets $V_i$ mapping étale to $\A^1_S$ via a coordinate $t^{(i)}$ such that $E[D]\cap V_i$ either is empty or - if it is not empty - is cut out by the equation $t^{(i)}=0$. The sheaf $\Omega^1_{E/S}(\mathrm{log}(E[D]))$, as a submodule of $(j_D)_*\Omega^1_{U_D/S}$, is then characterized by the fact that on such a $V_i$ it is freely generated over $\mathcal O_{V_i}$ by the basis element $\mathrm{d}t^{(i)}$ if $E[D]\cap V_i= \emptyset$ resp. by $\frac{\mathrm{d}t^{(i)}}{t^{(i)}}$ otherwise.\\
Observe moreover that due to the further smoothness of $S/\Q$ it is always possible to find an open covering $(V_i)_i$ of $E$ with a pair of coordinates $\{t^{(i)},s^{(i)}\}$ mapping $V_i$ étale to $\A^2_{\Q}$ and such that the coordinate $t^{(i)}$ has the property just mentioned.\\
\newline
Now let $\mathcal V$ be as above a vector bundle on $E$ with integrable $\Q$-connection.\\
Composing its relative de Rham complex
\[\Omega^{\bullet}_{E/S}(\mathcal V)=[\mathcal V \rightarrow \Omega^1_{E/S}\otimes_{\mathcal O_E} \mathcal V]\]
with the natural inclusion
\[\Omega^1_{E/S}\otimes_{\mathcal O_E}\mathcal V \rightarrow \Omega^1_{E/S}(\mathrm{log}(E[D]))\otimes_{\mathcal O_E}\mathcal V\]
defines the relative logarithmic de Rham complex of $\mathcal V$ along $E[D]$:
\[\Omega^{\bullet}_{E/S}(\mathrm{log}(E[D]))(\mathcal V):=[\mathcal V \rightarrow \Omega^1_{E/S}(\mathrm{log}(E[D]))\otimes_{\mathcal O_E}\mathcal V].\]
Enfolding the definitions it is easily seen that we have a well-defined and canonical short exact sequence of complexes
\[\tag{\textbf{3.7.12}} 0 \rightarrow \Omega^{\bullet}_{E/S}(\mathcal V) \xrightarrow{\mathrm{can}} \Omega^{\bullet}_{E/S}(\mathrm{log}(E[D]))(\mathcal V) \xrightarrow {\mathrm{res}} (i_D)_*i_D^*\mathcal V \ [-1] \rightarrow 0,\]
where the map $\mathrm{res}$ in degree $1$:
\[\Omega^1_{E/S}(\mathrm{log}(E[D]))\otimes_{\mathcal O_E} \mathcal V \rightarrow (i_D)_*i_D^*\mathcal V\]
is defined locally (using an open coordinate covering as above) by
\[\mathrm{d}t^{(i)} \otimes v \mapsto 0 \quad \textrm{resp.} \quad \frac{\mathrm{d}t^{(i)}}{t^{(i)}} \otimes v\mapsto i_D^*(v).\]
According to \cite{De1}, Ch. II, Cor. 3.14 (i) and Rem. 3.16, the canonical homomorphism
\[\tag{\textbf{3.7.13}} \Omega^{\bullet}_{E/S}(\mathrm{log}(E[D]))(\mathcal V) \rightarrow (j_D)_*(\Omega^{\bullet}_{U_D/S}(\mathcal V_{|U_D}))\]
is a quasi-isomorphism. The long exact sequence of hyperderived functors $\R^i \pi_*$ for $(3.7.12)$ combined with $(3.7.13)$ and the fact that $j_D$ is affine thus induces\footnote{Observe that because $j_D$ is affine one has
\[R^q(j_D)_*(\Omega^p_{U_D/S}\otimes_{\mathcal O_{U_D}}\mathcal V_{|U_D})=0 \quad \textrm{for all} \ q > 0,\]
which implies that the canonical arrow
\[(j_D)_*(\Omega^{\bullet}_{U_D/S}(\mathcal V_{|U_D})) \rightarrow \mathbb R(j_D)_*(\Omega^{\bullet}_{U_D/S}(\mathcal V_{|U_D}))\]
is a quasi-isomorphism (such a conclusion also appears in \cite{De1}, Ch. II, $(6.4.1)$ and $(6.4.3)$). The rest is then clear.} an exact sequence of vector bundles

\[\tag{\textbf{3.7.14}} 0 \rightarrow H^1_{\mathrm{dR}}(E/S, \mathcal V) \xrightarrow{\mathrm{can}} H^1_{\mathrm{dR}}(U_D/S, \mathcal V) \xrightarrow{\mathrm{res}} H^0_{\mathrm{dR}}(E[D]/S,i_D^*\mathcal V) \rightarrow H^2_{\mathrm{dR}}(E/S, \mathcal V) \rightarrow 0.\]

This is nothing else than the sequence $(3.7.11)$ deduced earlier by the machinery of the localization sequence.\footnote{A detailed verification of this fact consists in the tedious task of precisely retracing the construction of the localization sequence with all involved identifications and using the explicit local description of $(i_D)_+i_D^*\mathcal V$ as in \cite{Ho-Ta-Tan}, Ex. 1.5.23.} This observation provides a sufficiently explicit knowledge of the abstract residue map $\mathrm{Res}$ which appeared there, and hence in particular of the arrows $\mathrm{Res}^D_n$ in the case $\mathcal V=\mathcal L_n$.

\subsection{Preliminaries: Analytic}
We show that after base extension to $\Spec(\C)$ the logarithm sheaves are vector bundles with regular integrable connection. We then define the analogue of $(3.7.14)$ for the analytification of such a vector bundle and -  in the case that the bundle is moreover originally defined over $\Q$ - record compatibility of the algebraic and the analytic construction under comparison isomorphisms for de Rham cohomology.\\
The main technical input here is provided by the theory of regular integrable connections (cf. \cite{De1}).

\subsubsection{Regularity of the logarithm sheaves}
We denote by $E^{\C}$ resp. $S^{\C}$ the smooth complex algebraic varieties $E\times_{\Q}\C$ resp. $S\times_{\Q}\C$.\\
The maps $\pi^{\C}:=\pi \times \mathrm{id}$ resp. $\epsilon^{\C}:=\epsilon\times \mathrm{id}$ make $E^{\C}/S^{\C}$ into an elliptic curve.\\
For a vector bundle $\mathcal V$ on $E$ we write $\mathcal V^{\C}$ for its pullback along $E^{\C}\rightarrow E$; if $\mathcal V$ is equipped with an integrable $\Q$-connection $\nabla$, then $\mathcal V^{\C}$ carries a naturally induced integrable $\C$-connection $\nabla^{\C}$; the analogous remark holds for $E$ replaced by $S$.\\
\newline
For vector bundles with integrable connection on a smooth complex algebraic variety we have the well-known fundamental notion of regularity (cf. \cite{De1}, Ch. II, Def. 4.5). We then have the following basic result about the logarithm sheaves; the crucial input for its proof is the regularity of the Gauß-Manin connection, going back to Deligne, Griffiths and Katz.
\begin{proposition}
$(\mathcal L_n^{\C},\nabla_n^{\C})$ is regular for each $n\geq 0$.
\end{proposition}
\begin{proof}
As $(\mathcal L_0^{\C},\nabla_0^{\C})=(\mathcal O_{E^{\C}},\mathrm{d})$ the claim is certainly true for $n=0$.\\
Furthermore, for each $n\geq 0$ the exact sequence of $(1.1.3)$:
\[0 \rightarrow \mathrm{Sym}^{n+1}_{\mathcal O_E}\mathcal H_E \rightarrow \mathcal L_{n+1} \rightarrow \mathcal L_n \rightarrow 0\]
together with the canonical identifications\footnote{The first is obvious and the second follows from the compatibility of the de Rham cohomology of $E/S$ with base change (cf. the beginning of Chapter 1); the base change to be considered here (along $S^{\C}\rightarrow S$) is even flat.}
\[(\mathrm{Sym}^{n+1}_{\mathcal O_E}\mathcal H_E)^{\C} \simeq \mathrm{Sym}^{n+1}_{\mathcal O_{E^{\C}}}((\pi^{\C})^*\mathcal H^{\C}) \simeq \mathrm{Sym}^{n+1}_{\mathcal O_{E^{\C}}}((\pi^{\C})^*H^1_{\mathrm{dR}}(E^{\C}/S^{\C})^\vee)\]
induces the horizontal exact sequence
\[\tag{$*$} 0 \rightarrow \mathrm{Sym}^{n+1}_{\mathcal O_{E^{\C}}}((\pi^{\C})^*H^1_{\mathrm{dR}}(E^{\C}/S^{\C})^\vee) \rightarrow \mathcal L_{n+1}^{\C} \rightarrow \mathcal L_n^{\C} \rightarrow 0,\]
where $H^1_{\mathrm{dR}}(E^{\C}/S^{\C})$ is equipped with its Gauß-Manin connection relative $\Spec(\C)$.\\
The crucial point now is the regularity of $H^1_{\mathrm{dR}}(E^{\C}/S^{\C})$ (cf. \cite{De1}, Ch. II, Thm. 7.9 and Prop. 6.14), which then implies the regularity of $H^1_{\mathrm{dR}}(E^{\C}/S^{\C})^\vee$ (cf. ibid., Ch. II, Prop. 4.6 (ii)), hence also of $(\pi^{\C})^*H^1_{\mathrm{dR}}(E^{\C}/S^{\C})^\vee$ (cf. ibid., Ch. II, Prop. 4.6 (iii)) and then also of $\mathrm{Sym}^{n+1}_{\mathcal O_{E^{\C}}}((\pi^{\C})^*H^1_{\mathrm{dR}}(E^{\C}/S^{\C})^\vee)$ (by an analogous argument as used in ibid., Ch. II, Prop. 4.6 (ii), for the other tensor operations).\\
Using this and the fact that the property of regularity is closed under extensions (by ibid., Ch. II, Prop. 4.6. (i)) we see from $(*)$ and induction that regularity of $\mathcal L_1^{\C}$ will imply regularity of $\mathcal L_n^{\C}$ for all $n\geq 2$. For $n=0$ the sequence $(*)$ reads as
\[0 \rightarrow (\pi^{\C})^*H^1_{\mathrm{dR}}(E^{\C}/S^{\C})^\vee \rightarrow \mathcal L_1^{\C} \rightarrow \mathcal O_{E^{\C}} \rightarrow 0,\]
and from the regularity of the outer terms we deduce, as already noted, that $\mathcal L_1^{\C}$ is regular.
\end{proof}

The next goal is to assure that the algebraic constructions at the end of 3.7.1 by which we arrived at the sequence $(3.7.14)$ carry over to the analytic situation - provided that we assume regularity.

\subsubsection{Implications of regularity}
When working with analytifications we now no longer tacitly consider a fixed connected component (as we did in most of 3.5 and 3.6) but really mean the full objects with all components.\\
\newline
The complex submanifold $(E[D])^{an}=E^{an}[D]$ of $E^{an}$ defines an effective relative Cartier divisor in $E^{an}$ with normal crossings over $S^{an}$. The following diagram refreshes some associated notation.
\begin{equation*}
\begin{xy}
\xymatrix{
& E^{an}[D] \ar[dr]_{\pi_{E[D]}^{an}} \ar[r]^{ \ \ i_D^{an}} & E^{an}  \ar[d]_{\pi^{an}} & U_D^{an} \ar[dl]^{\pi_{U_D}^{an}} \ar[l]_{ \ \ j_D^{an}}\\
& & S^{an} & &}
\end{xy}
\end{equation*}
As in the algebraic setting we have the sheaf $\Omega^1_{E^{an}/S^{an}}(\mathrm{log}(E^{an}[D]))$ of relative differential forms with logarithmic poles along $E^{an}[D]$ (cf. \cite{De1}, Ch. II, $(3.3.1)$), with a local description as follows:\\
One can cover the $2$-dimensional complex manifold $E^{an}$ by open subsets $V_i$ with complex coordinate functions $\{t^{(i)},s^{(i)}\}$ such that $\mathrm{d}t^{(i)}$ freely generates $\Omega^1_{E^{an}/S^{an}}$ over $\mathcal O_{V_i}$ and $E^{an}[D]\cap V_i$ either is empty or - if it is not empty - is cut out by the equation $t^{(i)}=0$. The sheaf $\Omega^1_{E^{an}/S^{an}}(\mathrm{log}(E^{an}[D]))$, as a submodule of $(j^{an}_D)_*\Omega^1_{U_D^{an}/S^{an}}$, is then characterized by the fact that on such a $V_i$ it is freely generated over $\mathcal O_{V_i}$ by the basis element $\mathrm{d}t^{(i)}$ if $E^{an}[D]\cap V_i= \emptyset$ resp. by $\frac{\mathrm{d}t^{(i)}}{t^{(i)}}$ otherwise.\\
\newline
For each analytic vector bundle $\mathcal W$ on $E^{an}$ with (absolute) integrable connection the definition of its relative logarithmic de Rham complex along $E^{an}[D]$:
\[\Omega^{\bullet}_{E^{an}/S^{an}}(\mathrm{log}(E^{an}[D]))(\mathcal W):=[\mathcal W \rightarrow \Omega^1_{E^{an}/S^{an}}(\mathrm{log}(E^{an}[D]))\otimes_{\mathcal O_{E^{an}}}\mathcal W]\]
and of the canonical short exact sequence of complexes
\[\tag{\textbf{3.7.15}} 0 \rightarrow \Omega^{\bullet}_{E^{an}/S^{an}}(\mathcal W) \xrightarrow{\mathrm{can}} \Omega^{\bullet}_{E^{an}/S^{an}}(\mathrm{log}(E^{an}[D]))(\mathcal W) \xrightarrow {\mathrm{res}} (i_D^{an})_*(i_D^{an})^*\mathcal W \ [-1] \rightarrow 0\]
then proceeds analogously to the already treated algebraic counterpart in 3.7.1.\\
\newline
The problem now is the following: by contrast to the algebraic setting the hyperderivation $\mathbb R^i\pi^{an}_*$ of $\Omega^{\bullet}_{E^{an}/S^{an}}(\mathrm{log}(E^{an}[D]))(\mathcal W)$ does in general not compute the sheaf $H^i_{\mathrm{dR}}(U_D^{an}/S^{an},\mathcal W)$. Namely, by \cite{De1}, Ch. II, Prop. 3.13 and Rem. 3.16, it is instead true that the canonical homomorphism
\[\tag{\textbf{3.7.16}} \Omega^{\bullet}_{E^{an}/S^{an}}(\mathrm{log}(E^{an}[D]))(\mathcal W) \rightarrow (j^{an}_D)^{mer}_*(\Omega^{\bullet}_{U^{an}_D/S^{an}}(\mathcal W_{|U^{an}_D}))\]
is a quasi-isomorphism, where the superscript $"mer"$ on the right side means that we take the subcomplex of $(j^{an}_D)_*(\Omega^{\bullet}_{U^{an}_D/S^{an}}(\mathcal W_{|U^{an}_D}))$ of those sections which are meromorphic along $E^{an}[D]$.\\
\newline
It is the regularity condition which resolves this problem.\\
\newline
\textit{\textbf{Assume now that $\mathcal W=\mathcal Z^{an}$ for an algebraic vector bundle $\mathcal Z$ on $E^{\C}$ with regular integrable connection relative $\Spec(\C)$.}}\\
\newline
Of course, the special case we will be interested in is $\mathcal W=(\mathcal L_n^{\C})^{an}$, i.e. $\mathcal W=\mathcal L_n^{an}$ in our usual notation; Prop. 3.7.1 guarantees that this is indeed a special case.\\
\newline
Under this assumption, by \cite{De1}, Ch. II, Lemma 6.18\footnote{Observe the obvious typos in the cited lemma: one has to remove all bar-superscripts.}, the canonical inclusion
\[\tag{\textbf{3.7.17}} (j^{an}_D)^{mer}_*(\Omega^{\bullet}_{U^{an}_D/S^{an}}(\mathcal W_{|U^{an}_D})) \rightarrow (j^{an}_D)_*(\Omega^{\bullet}_{U^{an}_D/S^{an}}(\mathcal W_{|U^{an}_D}))\]
is a quasi-isomorphism. Combining this with $(3.7.16)$ we see that the natural homomorphism
\[\tag{\textbf{3.7.18}} \Omega^{\bullet}_{E^{an}/S^{an}}(\mathrm{log}(E^{an}[D]))(\mathcal W) \rightarrow (j^{an}_D)_*(\Omega^{\bullet}_{U^{an}_D/S^{an}}(\mathcal W_{|U^{an}_D}))\]
is in fact a quasi-isomorphism. Together with the fact that $j_D^{an}$ is Stein (because $j_D$ is affine) we conclude\footnote{Observe that by Cartan's Theorem B one has
\[R^q(j_D^{an})_*(\Omega^p_{U_D^{an}/S^{an}}\otimes_{\mathcal O_{U_D^{an}}}\mathcal W_{|U^{an}_D})=0 \quad \textrm{for all} \ q > 0,\]
which implies that the canonical arrow
\[(j_D^{an})_*(\Omega^{\bullet}_{U_D^{an}/S^{an}}(\mathcal W_{|U^{an}_D})) \rightarrow \mathbb R(j_D^{an})_*(\Omega^{\bullet}_{U_D^{an}/S^{an}}(\mathcal W_{|U^{an}_D}))\]
is a quasi-isomorphism (such a conclusion also appears in \cite{De1}, Ch. II, $(6.4.2)$ and $(6.4.4)$). The rest is then clear.} that for all $i$:
\[\tag{\textbf{3.7.19}} \R^i\pi^{an}_*(\Omega^{\bullet}_{E^{an}/S^{an}}(\mathrm{log}(E^{an}[D]))(\mathcal W)) \xrightarrow{\sim} H^i_{\mathrm{dR}}(U_D^{an}/S^{an},\mathcal W).\]
The long exact sequence of hyperderived functors $\R^i\pi^{an}_*$ for $(3.7.15)$ combined with $(3.7.19)$ then induces the following exact sequence of $\mathcal O_{S^{an}}$-vector bundles (for $H^2_{\mathrm{dR}}(U_D^{an}/S^{an}, \mathcal W)=0$ one may use that $\pi^{an}_{U_D}$ is Stein because $\pi_{U_D}$ is affine):
\[\tag{\textbf{3.7.20}}
\begin{split}
0 & \rightarrow H^1_{\mathrm{dR}}(E^{an}/S^{an}, \mathcal W) \xrightarrow{\mathrm{can}} H^1_{\mathrm{dR}}(U^{an}_D/S^{an}, \mathcal W) \xrightarrow{\mathrm{res}} H^0_{\mathrm{dR}}(E^{an}[D]/S^{an},(i^{an}_D)^*\mathcal W)\\
&\rightarrow H^2_{\mathrm{dR}}(E^{an}/S^{an}, \mathcal W) \rightarrow 0.
\end{split}\]

Note that all maps of the preceding exact sequence $(3.7.20)$ are defined in purely analytic terms; we only had to assume that the analytic bundle $\mathcal W$ on $E^{an}$ comes from a regular complex algebraic bundle $\mathcal Z$ on $E^{\C}$ in order to ensure that the (analytically defined) map $(3.7.18)$ is a quasi-isomorphism.

\subsubsection{Comparison and compatibility properties}
We conclude our general discussion with the following observations valid for any vector bundle $\mathcal V$ on $E$ with integrable $\Q$-connection such that the induced vector bundle $\mathcal V^{\C}$  on $E^{\C}$ with integrable $\C$-connection is regular, hence in particular for $\mathcal V=\mathcal L_n$ (cf. Prop. 3.7.1).\\
\newline
For each $i$ the canonical arrows
\begin{align*}
H^i_{\mathrm{dR}}(E/S,\mathcal V)\otimes_{\mathcal O_{S}}\mathcal O_{S^{\C}} & \rightarrow H^i_{\mathrm{dR}}(E^{\C}/S^{\C},\mathcal V^{\C})\\
H^i_{\mathrm{dR}}(U_D/S,\mathcal V)\otimes_{\mathcal O_{S}}\mathcal O_{S^{\C}} & \rightarrow H^i_{\mathrm{dR}}(U^{\C}_D/S^{\C},\mathcal V^{\C})\\
H^i_{\mathrm{dR}}(E[D]/S, i_D^*\mathcal V)\otimes_{\mathcal O_{S}}\mathcal O_{S^{\C}} & \rightarrow H^i_{\mathrm{dR}}(E^{\C}[D]/S^{\C},\mathcal (i^{\C}_D)^*\mathcal V^{\C})
\end{align*}
and
\begin{align*}
H^i_{\mathrm{dR}}(E^{\C}/S^{\C},\mathcal V^{\C})^{an} & \rightarrow H^i_{\mathrm{dR}}((E^{\C})^{an}/(S^{\C})^{an}, (\mathcal V^{\C})^{an})\\
H^i_{\mathrm{dR}}(U^{\C}_D/S^{\C},\mathcal V^{\C})^{an} & \rightarrow H^i_{\mathrm{dR}}((U_D^{\C})^{an}/(S^{\C})^{an},(\mathcal V^{\C})^{an})\\
H^i_{\mathrm{dR}}(E^{\C}[D]/S^{\C},\mathcal (i^{\C}_D)^*\mathcal V^{\C})^{an} &\rightarrow H^i_{\mathrm{dR}}((E^{\C}[D])^{an}/(S^{\C})^{an},((i_D^{\C})^{an})^*(\mathcal V^{\C})^{an})
\end{align*}
are all isomorphisms: the first triple because $S^{\C}\rightarrow S$ is flat and the second triple - due to our regularity assumption - by \cite{De1}, Ch. II, Prop. 6.14 and Prop. 4.6 (iii).\\
In sum, if $\mathcal V$ is as above we have for all $i$ canonical isomorphisms of $\mathcal O_{S^{an}}$-vector bundles:
\[\tag{\textbf{3.7.21}}
\begin{split}
H^i_{\mathrm{dR}}(E/S,\mathcal V)^{an} &\xrightarrow{\sim} H^i_{\mathrm{dR}}(E^{an}/S^{an},\mathcal V^{an})\\
H^i_{\mathrm{dR}}(U_D/S,\mathcal V)^{an} &\xrightarrow{\sim} H^i_{\mathrm{dR}}(U_D^{an}/S^{an},\mathcal V^{an})\\
H^i_{\mathrm{dR}}(E[D]/S, i_D^*\mathcal V)^{an} &\xrightarrow{\sim} H^i_{\mathrm{dR}}(E^{an}[D]/S^{an},(i_D^{an})^*\mathcal V^{an})
\end{split}
\]
which are moreover checked to be horizontal (for the connections induced by the algebraic Gauß-Manin connections on the left sides and the analytic Gauß-Manin connections on the right sides).\\
\newline
Finally, we then have a commutative diagram of $\mathcal O_{S^{an}}$-vector bundles with integrable connection:

\[\tag{\textbf{3.7.22}} \begin{split}
{\footnotesize
\begin{xy}
\xymatrix@C-0.3cm{
0 \ar[r] & H^1_{\mathrm{dR}}(E/S, \mathcal V)^{an} \ar[d]^{\sim} \ar[r]^-{\mathrm{can}^{an}} & H^1_{\mathrm{dR}}(U_D/S, \mathcal V)^{an} \ar[d]^{\sim} \ar[r]^-{\mathrm{Res}^{an}}& H^0_{\mathrm{dR}}(E[D]/S,i_D^*\mathcal V)^{an} \ar[d]^{\sim} \ar[r] & H^2_{\mathrm{dR}}(E/S, \mathcal V)^{an} \ar[d]^{\sim} \ar[r] &0 \\
0 \ar[r] & H^1_{\mathrm{dR}}(E^{an}/S^{an}, \mathcal V^{an}) \ar[r]^-{\mathrm{can}} & H^1_{\mathrm{dR}}(U^{an}_D/S^{an}, \mathcal V^{an}) \ar[r]^-{\mathrm{res}}& H^0_{\mathrm{dR}}(E^{an}[D]/S^{an},(i_D^{an})^*\mathcal V^{an}) \ar[r] & H^2_{\mathrm{dR}}(E^{an}/S^{an}, \mathcal V^{an}) \ar[r] &0}
\end{xy}}
\end{split}
\]
where the upper resp. lower row is given by the exact sequence $(3.7.11)^{an}=(3.7.14)^{an}$ resp. by $(3.7.20)$ and where the vertical arrows are the isomorphisms of $(3.7.21)$.\\
To check the commutativity of $(3.7.22)$ is indeed a routine formal procedure if for the upper row one uses the construction $(3.7.12)$-$(3.7.14)$.
 
\subsection{The construction of the fundamental commutative diagram}
As recalled in 3.7.1 the system $\varD$ is the unique element in $\lim_{n\geq 0} H^1_{\mathrm{dR}}(U_D/\Q, \mathcal L_n)$ mapping to $D^2\cdot 1_{ \{\epsilon \}} - 1_{E[D]}$ under the composition of injections $(3.7.7)$:
\[\lim_{n \geq 0} H^1_{\mathrm{dR}}(U_D/\Q,\mathcal L_n) \rightarrow \lim_{n \geq 0} \Gamma(S,H^1_{\mathrm{dR}}(U_D/S,\mathcal L_n)) \xrightarrow{\mathrm{Res}^D} \lim_{n\geq 0}H^0(E[D],i_D^*\mathcal L_n).\]
We will now construct a completely analogous composition of injections on the analytic side:
\[\lim_{n\geq 0} H^1_{\mathrm{dR}}(U_D^{an},\mathcal L_n^{an}) \rightarrow \lim_{n\geq 0} \Gamma(S^{an},H^1_{\mathrm{dR}}(U_D^{an}/S^{an},\mathcal L_n^{an})) \xrightarrow{ \ \mathrm{res}^D} \lim_{n\geq 0}H^0(E^{an}[D],(i_D^{an})^*\mathcal L_n^{an}),\]
which is possible because for $\mathcal L_n^{an}$ we have the analytic residue morphism in the lower row of $(3.7.22)$.\\
Via natural analytification maps we then place these two compositions into a big commutative diagram which in 3.7.4 will be used to formulate the main result of this section.

\subsubsection{An auxiliary result}

\begin{lemma}
(i) For each $n\geq 0$ the morphism
\[H^1_{\mathrm{dR}}(E^{an}/S^{an},\mathcal L_{n+1}^{an}) \rightarrow H^1_{\mathrm{dR}}(E^{an}/S^{an},\mathcal L_n^{an}),\]
induced by the analytified transition map $\mathcal L^{an}_{n+1} \rightarrow \mathcal L^{an}_n$, is zero.\\
(ii) The same holds for the morphism
\[H^0_{\mathrm{dR}}(U_D^{an}/S^{an},\mathcal L_{n+1}^{an}) \rightarrow H^0_{\mathrm{dR}}(U_D^{an}/S^{an},\mathcal L_n^{an}).\]
\end{lemma}
\begin{proof}
(i) We have a commutative diagram
\begin{equation*}
\begin{xy}
\xymatrix@C-0.3cm{
H^1_{\mathrm{dR}}(E/S,\mathcal L_{n+1})^{an} \ar[d]\ar[r] & H^1_{\mathrm{dR}}(E/S,\mathcal L_n)^{an} \ar[d]\\
H^1_{\mathrm{dR}}(E^{an}/S^{an},\mathcal L_{n+1}^{an})\ar[r] & H^1_{\mathrm{dR}}(E^{an}/S^{an},\mathcal L_n^{an})}
\end{xy}
\end{equation*}
Here, the upper horizontal arrow is the analytification of the morphism
\[\tag{$*$} H^1_{\mathrm{dR}}(E/S,\mathcal L_{n+1}) \rightarrow H^1_{\mathrm{dR}}(E/S,\mathcal L_n)\]
induced by the algebraic transition map, the lower horizontal arrow is the map of the claim and the vertical arrows are the canonical ones. From Thm. 1.2.1 (ii) we know that $(*)$ is zero, hence the same holds for $H^1_{\mathrm{dR}}(E/S,\mathcal L_{n+1})^{an} \rightarrow H^1_{\mathrm{dR}}(E/S,\mathcal L_n)^{an}$. Furthermore, due to $(3.7.21)$, the vertical maps of the diagram are isomorphisms - note that $(3.7.21)$ can be applied to the logarithm sheaves because of the regularity result in Prop. 3.7.1. This implies the claim of (i).\\
(ii) The assertion of (ii) follows in complete analogy if one knows that the transition map
\[H^0_{\mathrm{dR}}(U_D/S,\mathcal L_{n+1}) \rightarrow H^0_{\mathrm{dR}}(U_D/S,\mathcal L_n)\]
is zero for each $n\geq 0$. This is seen as follows: the canonical arrow
\[H^0_{\mathrm{dR}}(E/S,\mathcal L_n) \rightarrow H^0_{\mathrm{dR}}(U_D/S,\mathcal L_n)\]
is an isomorphism (prolong the exact sequence $(3.7.11)$ to the left); moreover, the transition map
\[H^0_{\mathrm{dR}}(E/S,\mathcal L_{n+1}) \rightarrow H^0_{\mathrm{dR}}(E/S,\mathcal L_n)\]
is zero, again by Thm. 1.2.1 (ii). This suffices to conclude.
\end{proof}

\subsubsection{The construction: first step}

The analytic residue maps of the lower row of $(3.7.22)$ are available for the logarithm sheaves because of Prop. 3.7.1 and induce - by taking global $S^{an}$-sections and inverse limit - an arrow
\[\tag{\textbf{3.7.23}} \begin{split} \lim_{n \geq 0} \Gamma(S^{an},H^1_{\mathrm{dR}}(U^{an}_D/S^{an},\mathcal L^{an}_n)) \xrightarrow{\mathrm{ \ res}^D} &\lim_{n\geq 0}\Gamma(S^{an},H^0_{\mathrm{dR}}(E^{an}[D]/S^{an},(i_D^{an})^*\mathcal L^{an}_n))\\
=&\lim_{n\geq 0}H^0(E^{an}[D],(i_D^{an})^*\mathcal L^{an}_n),\end{split}\]and Lemma 3.7.2 (i) implies that it is injective.\\
\newline
The injection $(3.7.23)$ fits into the following commutative square:
\begin{equation*} \tag{\textbf{3.7.24}} \begin{split}
\begin{xy}
\xymatrix{
\lim_{n\geq 0} \Gamma(S,H^1_{\mathrm{dR}}(U_D/S,\mathcal L_n)) \ar[r]^-{\mathrm{Res}^D} \ar[d]^{\mathrm{can}} & \lim_{n\geq 0}H^0(E[D],i_D^*\mathcal L_n) \ar[d]^{\mathrm{can}} \\
\lim_{n\geq 0}\Gamma(S^{an},H^1_{\mathrm{dR}}(U^{an}_D/S^{an},\mathcal L^{an}_n)) \ar[r]^-{\mathrm{res}^D} & \lim_{n\geq 0}H^0(E^{an}[D],(i_D^{an})^*\mathcal L^{an}_n)}
\end{xy}
\end{split}
\end{equation*}
where the upper horizontal arrow $\mathrm{Res}^D$ is the injection of $(3.7.6)$ and the vertical maps are the canonical ones: namely, note that by adjunction along $S^{an}\rightarrow S$ we have the canonical arrow
\[\Gamma(S,H^1_{\mathrm{dR}}(U_D/S,\mathcal L_n)) \rightarrow \Gamma(S^{an},H^1_{\mathrm{dR}}(U_D/S,\mathcal L_n)^{an})\]
which we may further compose with the natural map of $(3.7.21)$ in global $S^{an}$-sections
\[\Gamma(S^{an},H^1_{\mathrm{dR}}(U_D/S,\mathcal L_n)^{an}) \rightarrow \Gamma(S^{an},H^1_{\mathrm{dR}}(U_D^{an}/S^{an},\mathcal L_n^{an})).\]
The induced map in the inverse limit is the left vertical arrow of $(3.7.24)$, and the right vertical arrow comes from adjunction along $E^{an}[D]\rightarrow E[D]$ and going into the inverse limit.\\
That $(3.7.24)$ commutes is immediate if one recalls the commutativity of the middle square of $(3.7.22)$.

\subsubsection{The construction: second step}

We have the edge morphism
\[H^1_{\mathrm{dR}}(U_D^{an},\mathcal L_n^{an}) \rightarrow H^0_{\mathrm{dR}}(S^{an}, H^1_{\mathrm{dR}}(U_D^{an}/S^{an},\mathcal L_n^{an}))\]
of the Leray spectral sequence\footnote{The spectral sequence is available also for our analytic situation as is immediately clear by its construction in \cite{Kat2}, Rem. $(3.3)$, resp. \cite{Har}, Ch. III, §4, Thm. $(4.1)$.}
\[E_2^{p,q}=H^p_{\mathrm{dR}}(S^{an}, H^q_{\mathrm{dR}}(U_D^{an}/S^{an},\mathcal L_n^{an})) \Rightarrow E^{p+q}=H^{p+q}_{\mathrm{dR}}(U_D^{an},\mathcal L_n^{an}).\]
Its composition with
\[H^0_{\mathrm{dR}}(S^{an}, H^1_{\mathrm{dR}}(U_D^{an}/S^{an},\mathcal L_n^{an})) \subseteq \Gamma(S^{an},H^1_{\mathrm{dR}}(U_D^{an}/S^{an},\mathcal L_n^{an}))\]
yields in the inverse limit a morphism of $\C$-vector spaces
\[\tag{\textbf{3.7.25}} \lim_{n\geq 0} H^1_{\mathrm{dR}}(U_D^{an},\mathcal L_n^{an}) \rightarrow \lim_{n\geq 0} \Gamma(S^{an},H^1_{\mathrm{dR}}(U_D^{an}/S^{an},\mathcal L_n^{an})).\]
The above edge morphism becomes bijective in the limit as follows from the vanishing of the transition maps of $H^0_{\mathrm{dR}}(U_D^{an}/S^{an},\mathcal L_n^{an})$ (cf. Lemma 3.7.2 (ii)). This implies that $(3.7.25)$ is injective.\\
\newline
By composition of the injections $(3.7.25)$ and $(3.7.23)$ we obtain the injective $\C$-linear map

\[\tag{\textbf{3.7.26}} \lim_{n\geq 0} H^1_{\mathrm{dR}}(U_D^{an},\mathcal L_n^{an}) \rightarrow \lim_{n\geq 0} \Gamma(S^{an},H^1_{\mathrm{dR}}(U_D^{an}/S^{an},\mathcal L_n^{an})) \xrightarrow{ \ \mathrm{res}^D} \lim_{n\geq 0}H^0(E^{an}[D],(i_D^{an})^*\mathcal L_n^{an})\]
which becomes the lower row of the following crucial commutative diagram:
\begin{equation*} \tag{\textbf{3.7.27}} \begin{split}
\begin{xy}
\xymatrix@C-0.3cm{
\lim_{n\geq 0} H^1_{\mathrm{dR}}(U_D/\Q,\mathcal L_n) \ar[d]^{\mathrm{can}} \ar[r]  & \lim_{n\geq 0} \Gamma(S,H^1_{\mathrm{dR}}(U_D/S,\mathcal L_n)) \ar[r]^{\quad \mathrm{Res}^D} \ar[d]^{\mathrm{can}} & \lim_{n\geq 0}H^0(E[D],i_D^*\mathcal L_n) \ar[d]^{\mathrm{can}}\\
\lim_{n\geq 0} H^1_{\mathrm{dR}}(U_D^{an},\mathcal L_n^{an})\ar[r] & \lim_{n\geq 0} \Gamma(S^{an},H^1_{\mathrm{dR}}(U_D^{an}/S^{an},\mathcal L_n^{an})) \ar[r]^{\quad \ \mathrm{res}^D} & \lim_{n\geq 0}H^0(E^{an}[D],(i_D^{an})^*\mathcal L_n^{an})}
\end{xy}
\end{split}
\end{equation*}
Here, the upper row is the composition of injections given by $(3.7.7)$, the right square is $(3.7.24)$ and the left vertical arrow comes from the composition of canonical maps
\[\tag{\textbf{3.7.28}} H^1_{\mathrm{dR}}(U_D/\Q,\mathcal L_n) \rightarrow H^1_{\mathrm{dR}}(U_D/\Q,\mathcal L_n)\otimes_{\Q}\C \xrightarrow{\sim} H^1_{\mathrm{dR}}(U_D^{\C}/\C,\mathcal L_n^{\C})\xrightarrow{\sim} H^1_{\mathrm{dR}}(U_D^{an},\mathcal L_n^{an}).\]
Note that the first map in $(3.7.28)$ clearly is injective and that the second resp. third map is indeed an isomorphism by the flatness of $\Spec(\C)\rightarrow \Spec(\Q)$ resp. by the regularity of $\mathcal L_n^{\C}$ on $U_D^{\C}$ (use Prop. 3.7.1 and \cite{De1}, Ch. II, Prop. 4.6 (iii)) and the comparison result of \cite{De1}, Ch. II, Thm. 6.2.\\
In particular, we see that $(3.7.28)$ and hence also the left vertical arrow of $(3.7.27)$ is injective.

\subsubsection{Further remarks about the arrows of the diagram}
We have a sufficiently explicit description for the algebraic and analytic residue map $\mathrm{Res}^D$ and $\mathrm{res}^D$ in $(3.7.27)$ in terms of logarithmic differential forms (cf. $(3.7.12)$-$(3.7.14)$ and $(3.7.15)$-$(3.7.20)$).\\
Furthermore, the upper left horizontal arrow was expressed in $(3.7.8)$-$(3.7.10)$ via the transition from absolute to relative differential forms. There is an analogous description for its analytic counterpart beneath, but words need to be chosen more carefully. Let us explain this.\\
\newline
The edge morphism $E_2^{1,0} \rightarrow E^1$ at the second sheet of the spectral sequence of hyperderived functors
\[E_1^{p,q}=H^q(U_D^{an},\Omega^p_{U_D^{an}}\otimes_{\mathcal O_{U_D^{an}}} \mathcal L_n^{an})\Rightarrow E^{p+q}=H^{p+q}_{\mathrm{dR}}(U_D^{an},\mathcal L_n^{an})\]
writes as

\[\tag{\textbf{3.7.29}} \frac{\ker\bigg(\Gamma(U_D^{an}, \Omega^1_{U_D^{an}}\otimes_{\mathcal O_{U_D^{an}}} \mathcal L_n^{an})\xrightarrow{(\nabla_n^{an})^1}\Gamma(U_D^{an}, \Omega^2_{U_D^{an}} \otimes_{\mathcal O_{U_D^{an}}} \mathcal L_n^{an})\bigg)}{\mathrm{im}\bigg(\Gamma(U_D^{an},\mathcal L_n^{an}) \xrightarrow{\nabla_n^{an}} \Gamma(U_D^{an},\Omega^1_{U_D^{an}} \otimes_{\mathcal O_{U_D^{an}}}\mathcal L_n^{an})\bigg)} \rightarrow H^1_{\mathrm{dR}}(U_D^{an},\mathcal L_n^{an}).\]
Analogously, the second sheet of the spectral sequence
\[E_1^{p,q}=R^q(\pi_{U_D}^{an})_*(\Omega^p_{U_D^{an}/S^{an}} \otimes_{\mathcal O_{U_D^{an}}} \mathcal L_n^{an}) \Rightarrow E^{p+q}=H^{p+q}_{\mathrm{dR}}(U_D^{an}/S^{an},\mathcal L_n^{an})\]
gives an edge morphism
\[\tag{\textbf{3.7.30}} \frac{(\pi_{U_D}^{an})_*(\Omega^1_{U_D^{an}/S^{an}} \otimes_{\mathcal O_{U_D^{an}}} \mathcal L_n^{an})}{\mathrm{im}\bigg((\pi_{U_D}^{an})_*(\mathcal L_n^{an}) \xrightarrow{(\nabla_n^{res})^{an}} (\pi_{U_D}^{an})_*(\Omega^1_{U_D^{an}/S^{an}} \otimes_{\mathcal O_{U_D^{an}}} \mathcal L_n^{an})\bigg)} \rightarrow H^1_{\mathrm{dR}}(U_D^{an}/S^{an},\mathcal L_n^{an}).\]

The maps $(3.7.29)$ and $(3.7.30)$ are in fact isomorphisms (cf. Rem. 3.7.3 below.)\\
From $(3.7.30)$ we obtain - by taking global $S^{an}$-sections and using the canonical morphisms between presheaf and associated sheaf - a natural arrow (not necessarily an isomorphism; cf. Rem. 3.7.3)

\[\tag{\textbf{3.7.31}} \frac{\Gamma(U_D^{an}, \Omega^1_{U_D^{an}/S^{an}} \otimes_{\mathcal O_{U_D^{an}}} \mathcal L_n^{an})}{\mathrm{im}\bigg(\Gamma(U_D^{an},\mathcal L_n^{an}) \xrightarrow{(\nabla_n^{res})^{an}}\Gamma(U_D^{an}, \Omega^1_{U_D^{an}/S^{an}} \otimes_{\mathcal O_{U_D^{an}}} \mathcal L_n^{an})\bigg)} \rightarrow \Gamma(S^{an},H^1_{\mathrm{dR}}(U_D^{an}/S^{an},\mathcal L_n^{an})).\]The isomorphism $(3.7.29)$ and the map $(3.7.31)$ form a commutative diagram with the morphism
\[H^1_{\mathrm{dR}}(U_D^{an},\mathcal L_n^{an}) \rightarrow \Gamma(S^{an},H^1_{\mathrm{dR}}(U_D^{an}/S^{an},\mathcal L_n^{an}))\]
whose limit is $(3.7.25)$ and with the restriction of absolute to relative differential forms.\\
This finally expresses more concretely also the lower left horizontal map of $(3.7.27)$ (which is $(3.7.25)$).

\begin{remark}
Using the fact that the complex manifold $U_D^{an}$ is Stein (because $U_D^{\C}$ is an affine scheme) resp. that the morphism $\pi_{U_D}^{an}$ is Stein (because $\pi_{U_D}^{\C}$ is an affine morphism) together with Cartan's Theorem B one sees that both of the previous spectral sequences degenerate at $r=2$ and that the maps $(3.7.29)$ resp. $(3.7.30)$ are isomorphisms.\\
But although also $S^{an}$ is Stein one cannot conclude that $(3.7.31)$ is an isomorphism (by contrast to the algebraic situation in $(3.7.9)$, where this was possible): the reason is that Theorem B only holds for coherent analytic sheaves, and hence the global $S^{an}$-sections of the sheaves on the left side of $(3.7.30)$ a priori are not the same as the left side of $(3.7.31)$.
\end{remark}

\subsection{The characterization result}
Altogether, we have constructed the fundamental commutative diagram $(3.7.27)$:
\begin{equation*}
\begin{xy}
\xymatrix@C-0.3cm{
\lim_{n\geq 0}H^1_{\mathrm{dR}}(U_D/\Q,\mathcal L_n) \ar[d]^{\mathrm{can}} \ar[r]  &\lim_{n\geq 0} \Gamma(S,H^1_{\mathrm{dR}}(U_D/S,\mathcal L_n)) \ar[r]^{\quad \ \mathrm{Res}^D} \ar[d]^{\mathrm{can}} & \lim_{n\geq 0}H^0(E[D],i_D^*\mathcal L_n) \ar[d]^{\mathrm{can}}\\
\lim_{n\geq 0} H^1_{\mathrm{dR}}(U_D^{an},\mathcal L_n^{an})\ar[r] & \lim_{n\geq 0} \Gamma(S^{an},H^1_{\mathrm{dR}}(U_D^{an}/S^{an},\mathcal L_n^{an})) \ar[r]^{\ \quad \mathrm{res}^D} & \lim_{n\geq 0}H^0(E^{an}[D],(i_D^{an})^*\mathcal L_n^{an})}
\end{xy}
\end{equation*}
and deduced the information that its rows as well as its left vertical arrow are injective.\\
This finally enables us to characterize uniquely the analytification of the $D$-variant of the polylogarithm by an analytic residue condition which is completely analogous to the algebraic one in $(3.7.7)$.\\
\newline
We start with the algebraic system
\[\varD = \Big(\varDn \Big)_{n\geq 0} \in \lim_{n\geq 0} H^1_{\mathrm{dR}}(U_D/\Q, \mathcal L_n)\]
and make the following
\begin{definition}
(i) For each $n\geq 0$ we write
\[(\varDn)^{an}\]
for the image of $\varDn$ under the canonical map of $(3.7.28)$:
\[H^1_{\mathrm{dR}}(U_D/\Q,\mathcal L_n) \rightarrow H^1_{\mathrm{dR}}(U_D^{an},\mathcal L_n^{an}).\]
(ii) We set
\[(\varD)^{an}:=\Big(\varDn\Big)^{an}_{n\geq0} \in \lim_{n\geq 0} H^1_{\mathrm{dR}}(U_D^{an}, \mathcal L_n^{an}),\]
which is the image of $\varD$ under the left vertical map of $(3.7.27)$.
\end{definition}

For each $n\geq 0$ the analytification of $(3.7.2)$ yields the chain of $\mathcal D_{E^{an}[D]}$-linear maps:

\[\tag{\textbf{3.7.32}} \mathcal O_{E^{an}[D]} \hookrightarrow (\pi_{E[D]}^{an})^* \prod_{k=0}^n \mathrm{Sym}^k_{\mathcal O_{S^{an}}}\mathcal H^{an} \simeq (\pi_{E[D]}^{an})^*(\epsilon^{an})^*\mathcal L_n^{an} \simeq (i_D^{an})^*([D]^{an})^*\mathcal L_n^{an} \simeq (i_D^{an})^*\mathcal L_n^{an}\]
which induces an injection
\[\tag{\textbf{3.7.33}} H^0(E^{an}[D],\mathcal O_{E^{an}[D]}) \hookrightarrow \lim_{n\geq 0}H^0(E^{an}[D],(i_D^{an})^*\mathcal L_n^{an}).\]
The zero section $\epsilon: S \rightarrow E[D]$ of $E[D]/S$ induces by analytification the map $\epsilon^{an}: S^{an}\rightarrow E^{an}[D]$ which is an open and closed immersion. From the decomposition
\[E^{an}[D]=(E^{an}[D]\ - \{\epsilon^{an}\}) \amalg \{\epsilon^{an}\}\]
we can define precisely as in the algebraic situation the global section $D^2\cdot 1_{ \{\epsilon^{an} \}} - 1_{E^{an}[D]}$ of $\mathcal O_{E^{an}[D]}$ and may view it as an element of the right side of $(3.7.33)$.\\
\newline
The fact that the upper row of $(3.7.27)$ maps $\varD$ to $D^2\cdot 1_{ \{\epsilon \}} - 1_{E[D]}$ (cf. $(3.7.7)$) together with the commutativity of $(3.7.27)$ and of
\begin{equation*}
\begin{xy}
\xymatrix@C-0.3cm{
H^0(E[D],\mathcal O_{E[D]}) \ar[r] \ar[d]^{\mathrm{can}} &  \lim_{n\geq 0}H^0(E[D],i_D^*\mathcal L_n)\ar[d]^{\mathrm{can}}\\
H^0(E^{an}[D],\mathcal O_{E^{an}[D]}) \ar[r] & \lim_{n\geq 0}H^0(E^{an}[D],(i_D^{an})^*\mathcal L^{an}_n)}
\end{xy}
\end{equation*}
(with horizontal maps given by $(3.7.3)$ resp. $(3.7.33)$) obviously implies that
\[(\varD)^{an}\]
under the lower row of $(3.7.27)$ goes to
\[D^2\cdot 1_{ \{\epsilon^{an} \}} - 1_{E^{an}[D]}.\]
As we already explained (at $(3.7.26)$) that this row is injective we arrive at the following
\begin{theorem}
The analytification $(\varD)^{an}$ of the system $\varD$ (cf. Def. 3.7.4) is the unique element in $\displaystyle\lim_{n\geq 0} H^1_{\mathrm{dR}}(U_D^{an}, \mathcal L_n^{an})$ which under the lower row of $(3.7.27)$:
\[\lim_{n\geq 0} H^1_{\mathrm{dR}}(U_D^{an},\mathcal L_n^{an}) \rightarrow \lim_{n\geq 0} \Gamma(S^{an},H^1_{\mathrm{dR}}(U_D^{an}/S^{an},\mathcal L_n^{an})) \xrightarrow{ \ \mathrm{res}^D} \lim_{n\geq 0}H^0(E^{an}[D],(i_D^{an})^*\mathcal L_n^{an})\]
goes to $D^2\cdot 1_{ \{\epsilon^{an} \}} - 1_{E^{an}[D]}$. \qquad \qed
\end{theorem}

It is certainly worth mentioning that this characterization can also be used to determine the $D$-variant of the polylogarithm algebraically:

\begin{remark}
Combining the injectivity of the lower row of $(3.7.27)$ and of its left vertical arrow (cf. the explanations at $(3.7.26)$ and $(3.7.28)$) we obtain in addition:\\
The system $\varD$ is the unique element in $\displaystyle\lim_{n\geq 0} H^1_{\mathrm{dR}}(U_D/\Q, \mathcal L_n)$ whose analytification (i.e. image under the left vertical arrow of $(3.7.27)$) goes to $D^2\cdot 1_{ \{\epsilon^{an} \}} - 1_{E^{an}[D]}$ under the lower row of $(3.7.27)$.
\end{remark}

\markright{\uppercase{The explicit description on the universal elliptic curve}}
\section{The $D$-variant of the polylogarithm for the universal family}
\markright{\uppercase{The explicit description on the universal elliptic curve}}

\subsection{The description of the analytified $D$-variant of the polylogarithm}
Fix an integer $D>1$.\\
\newline
Essentially by recapitulating the results of 3.6.1 we explain how the sections $p^D_n(z,\tau)$ define a system
\[p^D \in \lim_{n\geq 0} H^1_{\mathrm{dR}}(U_D^{an},\mathcal L_n^{an}).\]
We then state one of the main theorems of this work, namely the equality
\[(\varD)^{an}=p^D.\]
Apart from laborious technical aspects - arising from the need to chase through the chain of identifications $(3.7.32)$ - there are two crucial ingredients for the proof of this identity: the first is Thm. 3.7.5 which characterizes the left side by a computable analytic residue condition, and the second is the precise information we have about the singularities and residues of the functions $s^D_k(z,\tau)$.\\
\newline
In our notation for analytifications we again no longer tacitly consider a fixed connected component but really mean the full objects with all components, provided that we don't explicitly say otherwise.\\
\newline
To begin with, recall that in 3.6.1 we used the coefficient functions appearing in the Laurent expansion
\[D^2 \cdot J(z,-w,\tau)-D\cdot J\Big(Dz,-\frac{w}{D},\tau \Big)=s^D_0(z,\tau)+s^D_1(z,\tau)w+...\]
to build for each $n\geq 0$ the following vectors of length $2\cdot r(n)$ resp. $r(n)$:
\[p_n^D(z,\tau)=\begin{pmatrix} s^D_0(z,\tau) \\ s^D_1(z,\tau) \\ \vdots\\s^D_n(z,\tau)  \\  0 \\ \vdots \\ \vdots \\ 0 \\-\frac{1}{2\pi i }s^D_1(z,\tau)\\ -\frac{2}{2\pi i }s^D_2(z,\tau) \\ \vdots\\ -\frac{n+1}{2\pi i }s^D_{n+1}(z,\tau) \\ 0 \\ \vdots \\ \vdots \\ 0\end{pmatrix} \qquad \textrm{resp.} \qquad q_n^D(z,\tau)=\begin{pmatrix} s^D_0(z,\tau) \\ s^D_1(z,\tau) \\ \vdots\\s^D_n(z,\tau) \\ 0 \\ \vdots \\ \vdots \\ 0 \end{pmatrix}.\]
Using Thm. 3.6.2 resp. Prop. 3.6.1 uniformly for each connected component of $E^{an}$, we obtain from the vectors $p_n^D(z,\tau)$ resp. $q_n^D(z,\tau)$ an element in
\[\ker\bigg(\Gamma(U_D^{an}, \Omega^1_{E^{an}}\otimes_{\mathcal O_{E^{an}}} \mathcal L_n^{an})\xrightarrow{(\nabla_n^{an})^1}\Gamma(U_D^{an}, \Omega^2_{E^{an}} \otimes_{\mathcal O_{E^{an}}} \mathcal L_n^{an})\bigg)\]
resp. in
\[\Gamma(U_D^{an},\Omega^1_{E^{an}/S^{an}}\otimes_{\mathcal O_{E^{an}}} \mathcal L_n^{an}).\]
\begin{definition}
(i) Via the canonical map of $(3.7.29)$ resp. $(3.7.31)$ we receive for each $n\geq 0$ an induced element in
\[H^1_{\mathrm{dR}}(U_D^{an},\mathcal L_n^{an}) \ \ \ \textrm{resp.} \ \ \ \Gamma(S^{an},H^1_{\mathrm{dR}}(U_D^{an}/S^{an},\mathcal L_n^{an}))\]
which we denote by $p^D_n$ resp. $q^D_n$.\\
(ii) We write $p^D$ resp. $q^D$ for the inverse system
\[p^D:=(p^D_n)_{n\geq 0} \in \lim_{n\geq 0} H^1_{\mathrm{dR}}(U_D^{an},\mathcal L_n^{an}) \ \ \ \textrm{resp.} \ \ \ q^D:=(q^D_n)_{n\geq 0} \in \lim_{n\geq 0} \Gamma(S^{an},H^1_{\mathrm{dR}}(U_D^{an}/S^{an},\mathcal L_n^{an})).\]
\end{definition}
For part (ii) of the previous definition one observes that it follows from Lemma 3.6.4 that the $p^D_n$ resp. $q^D_n$ are indeed compatible under the respective transition maps.

\begin{lemma}
Under the natural map of $(3.7.25)$:
\[\lim_{n\geq 0} H^1_{\mathrm{dR}}(U_D^{an},\mathcal L_n^{an}) \rightarrow \lim_{n\geq 0} \Gamma(S^{an},H^1_{\mathrm{dR}}(U_D^{an}/S^{an},\mathcal L_n^{an}))\]
the system $p^D$ goes to the system $q^D$.
\end{lemma}
\begin{proof}
This follows from Lemma 3.6.3 and the explanation subsequent to $(3.7.31)$.
\end{proof}
Recalling the notation
\[(\varD)^{an}=\Big(\varDn\Big)^{an}_{n\geq0} \in \lim_{n\geq 0} H^1_{\mathrm{dR}}(U_D^{an}, \mathcal L_n^{an})\]
for the analytification of the system $\varD$ (cf. Def. 3.7.4), we are now ready to formulate and to prove the first main result of this section.
\begin{theorem}
We have the following equality of inverse systems in $\displaystyle\lim_{n\geq 0} H^1_{\mathrm{dR}}(U_D^{an}, \mathcal L_n^{an})$:
\[(\varD)^{an}=p^D.\]
\end{theorem}
\begin{proof} The proof is long and involved. We therefore divide it into several steps.\\
\newline
\underline{Step 1}: (Reduction of the problem and recapitulation of the setup)\\
\newline
By the characterization of $(\varD)^{an}$ in Thm. 3.7.5 it suffices to show that $p^D$ is sent to $D^2\cdot 1_{ \{\epsilon^{an} \}} - 1_{E^{an}[D]}$ under the lower row of $(3.7.27)$:
\[\lim_{n\geq 0} H^1_{\mathrm{dR}}(U_D^{an},\mathcal L_n^{an}) \rightarrow \lim_{n\geq 0} \Gamma(S^{an},H^1_{\mathrm{dR}}(U_D^{an}/S^{an},\mathcal L_n^{an})) \xrightarrow{ \ \mathrm{res}^D} \lim_{n\geq 0} H^0(E^{an}[D],(i_D^{an})^*\mathcal L_n^{an}).\]
Because of Lemma 3.8.2 we thus need to show the following equation in $\displaystyle\lim_{n\geq 0}H^0(E^{an}[D],(i_D^{an})^*\mathcal L_n^{an})$:
\[\tag{\textbf{3.8.1}} \mathrm{res}^D(q^D)=D^2\cdot 1_{ \{\epsilon^{an} \}} - 1_{E^{an}[D]}.\]
For each $n\geq 0$ let us write
\[\begin{split}
\mathrm{res}^D_n: \Gamma(S^{an},H^1_{\mathrm{dR}}(U^{an}_D/S^{an},\mathcal L^{an}_n)) &\rightarrow \Gamma(S^{an},H^0_{\mathrm{dR}}(E^{an}[D]/S^{an},(i_D^{an})^*\mathcal L^{an}_n))\\
&=H^0(E^{an}[D],(i_D^{an})^*\mathcal L^{an}_n)
\end{split}\]
for the analytic residue map in global $S^{an}$-sections as defined in $(3.7.20)$, such that $\mathrm{res}^D=(\mathrm{res}^D_n)_{n\geq 0}$.\\
In order to obtain $(3.8.1)$ we then have to show that $\mathrm{res}^D_n(q^D_n)$ equals the image of the global $\mathcal O_{E^{an}[D]}$-section $D^2\cdot 1_{ \{\epsilon^{an} \}} - 1_{E^{an}[D]}$ under the injection
\[\tag{\textbf{3.8.2}} H^0(E^{an}[D],\mathcal O_{E^{an}[D]}) \hookrightarrow H^0(E^{an}[D],(i_D^{an})^*\mathcal L_n^{an})\]
which in turn is induced by taking global sections in the chain of maps $(3.7.32)$:

\[\tag{\textbf{3.8.3}} \mathcal O_{E^{an}[D]} \hookrightarrow (\pi_{E[D]}^{an})^* \prod_{k=0}^n \mathrm{Sym}^k_{\mathcal O_{S^{an}}}\mathcal H^{an} \simeq (\pi_{E[D]}^{an})^*(\epsilon^{an})^*\mathcal L_n^{an} \simeq (i_D^{an})^*([D]^{an})^*\mathcal L_n^{an} \simeq (i_D^{an})^*\mathcal L_n^{an};\]
observe that the monomorphism in $(3.8.3)$ is $\frac{1}{n!}$-times the canonical inclusion (cf. $(3.7.2)$).\\
\newline
Recall from 3.7.2 that the map $\mathrm{res}^D_n$ comes about as follows: one has the diagram of complexes
\begin{equation*}
\begin{xy}
\xymatrix@C-0.7cm{
\Omega^{\bullet}_{E^{an}/S^{an}}(\mathrm{log}(E^{an}[D]))(\mathcal L_n^{an}) \ar[rd] \ar[rr] &  & (j^{an}_D)_*(\Omega^{\bullet}_{U^{an}_D/S^{an}}((\mathcal L_n^{an})_{|U^{an}_D})) \\
&(i_D^{an})_*(i_D^{an})^*\mathcal L_n^{an} \ [-1] &}
\end{xy}
\end{equation*}
where the horizontal arrow is the natural homomorphism and the left map is given in degree one by
\[\tag{\textbf{3.8.4}} \begin{split}  
&\Omega^1_{E^{an}/S^{an}}(\mathrm{log}(E^{an}[D]))\otimes_{\mathcal O_{E^{an}}} \mathcal L_n^{an} \rightarrow (i^{an}_D)_*(i_D^{an})^* \mathcal L_n^{an},\\
&\mathrm{d}t^{(i)} \otimes l \mapsto 0 \quad \textrm{resp.} \quad \frac{\mathrm{d}t^{(i)}}{t^{(i)}} \otimes l\mapsto (i_D^{an})^*(l) \end{split}\]
whenever $(V_i)_i$ is an open covering of $E^{an}$ (in the analytic topology) with complex coordinate functions $\{t^{(i)},s^{(i)}\}$ such that $\mathrm{d}t^{(i)}$ freely generates $\Omega^1_{E^{an}/S^{an}}$ over $\mathcal O_{V_i}$ and $E^{an}[D]\cap V_i$ is empty resp. - if it is not empty - is cut out by $t^{(i)}=0$.\\
As the horizontal arrow is a quasi-isomorphism (cf. $(3.7.18)$) and $j_D^{an}$ is Stein, the diagram yields by applying $\R^1\pi^{an}_*$ a map
\[H^1_{\mathrm{dR}}(U_D^{an}/S^{an},\mathcal L_n^{an}) \xrightarrow{\sim} \R^1 \pi^{an}_*(\Omega^{\bullet}_{E^{an}/S^{an}}(\mathrm{log}(E^{an}[D]))(\mathcal L_n^{an})) \rightarrow (\pi_{E[D]}^{an})_*(i_D^{an})^*\mathcal L^{an}_n\]which in global $S^{an}$-sections gives $\mathrm{res}^D_n$.\\
\newline
It is immediately clear that in order to compute $\mathrm{res}^D_n(q^D_n)$ we can consider a fixed connected component of $E^{an}/S^{an}$; hence, from now on until the end of the proof we restrict to such a component and commit the familiar abuse of notation by leaving away its index. We then write again
\[\mathrm{pr}: \C\times \H \rightarrow E^{an}=(\Z^2\times \Gamma(N))\backslash (\C\times \H)\]
for the universal covering map.\\
\newline
\underline{Step 2}: (Construction of charts for the residue computation)\\
\newline
Recall that $\Z^2\times \Gamma(N)$ acts properly discontinuously and freely on $\C\times \H$, such that $E^{an}$ is a complex manifold in the well-known natural way: a complex atlas is constructed by local inverses of $\mathrm{pr}$.\\
As usual, charts of $E^{an}$ are always charts of the associated maximal atlas.\\
\newline
Let us now construct an open covering of $E^{an}$ together with local coordinate charts which satisfy the conditions recalled at the end of Step 1. For this we proceed as follows:\\
\newline
If $\mathrm{P}$ is a point in $E^{an}-E^{an}[D]$ we choose an open neighborhood $V$ of $\mathrm{P}$ such that $V\cap E^{an}[D]=\emptyset$, there is an open set $V'\subseteq \C\times \H$ projecting bijectively to $V$ and such that $\Omega^1_{E^{an}/S^{an}}$ is freely generated over $\mathcal O_V$ by the differential of the first coordinate function of $(\mathrm{pr}_{|V'})^{-1}$.\\
This is the coordinate system which we fix around $\mathrm{P}$; it satisfies the conditions recalled in Step 1.\\
\newline
Next, we need to choose an adapted chart for each point $\mathrm{P_0}$ in $E^{an}\cap E^{an}[D]$.\\
\newline
The first step is to select for every point of $E^{an}\cap E^{an}[D]$ a representative in $\C \times \H$ as follows:\\
For each point of $S^{an}$ we choose and keep fixed an orbit representative in $\H$. Then, if $\mathrm{P_0}$ is a point in $E^{an}\cap E^{an}[D]$ it is easy to see that it has the form $\mathrm{P_0}=[(\frac{m_0}{D}\tau_0+\frac{n_0}{D},\tau_0)]$ for some $m_0,n_0 \in \Z$, where $\tau_0$ is the representative of $\mathrm{Q_0}:=\pi^{an}(\mathrm{P_0})$ we have just fixed. One further checks that all other representatives of $\mathrm{P_0}$ with $\tau_0$ in the second component are obtained by replacing $m_0$ resp. $n_0$ by $m_0+m'_0D$ resp. $n_0+n'_0D$ with $m'_0,n'_0 \in \Z$. In the set of integers $m_0+\Z\cdot D$ resp. $n_0+\Z\cdot D$ we take the smallest one of those which are greater or equal to zero.\\
Altogether, this procedure fixes for each $\mathrm{P_0} \in E^{an}\cap E^{an}[D]$ uniquely a representative in $\C\times \H$.\\
It is worth noting that the $\mathrm{P_0}$ which are already contained in $\epsilon^{an}(S^{an}) \subseteq E^{an}[D]$ then are represented by $(0,\tau_0)$, where $\tau_0$ is associated with $\mathrm{Q_0}=\pi^{an}(\mathrm{P_0})$ as before.\\
\newline
We now come to construct the desired chart around each point $\mathrm{P_0} \in E^{an}\cap E^{an}[D]$.\\
For this we take the just fixed representative $(\frac{m_0}{D}\tau_0+\frac{n_0}{D},\tau_0)$ of $\mathrm{P_0}$ in $\C\times \H$ and choose an open set $V'_0\subseteq \C\times \H$ containing it and projecting bijectively to an open neighborhood $V_0$ of $\mathrm{P_0}$ in $E^{an}$.\\
The subsequent Lemma 3.8.4 claims the existence of a certain open neighborhood $U_0'$ in $\C\times \H$ of the point $(\frac{m_0}{D}\tau_0+\frac{n_0}{D},\tau_0)$, and we can assume (by shrinking $V'_0$ and $V_0$ if necessary) that $V_0$ and $V'_0$ are as just chosen and that moreover $V'_0$ is contained in this $U_0'$.\\
The composition of $(\mathrm{pr}_{|V'_0})^{-1}$ with
\[\psi_0: V'_0 \rightarrow \C\times \H, \quad (z,\tau) \mapsto \Big(z-\frac{m_0}{D}\tau-\frac{n_0}{D},\tau \Big)\]
gives a chart
\[(t_0,s_0):V_0 \xrightarrow{\sim} \mathrm{im}(\psi_0) \subseteq \C\times \H \subseteq \C\times \C.\]
This is not yet the desired one; however, before continuing the construction we first show
\[\tag{\textbf{3.8.5}} E^{an}[D] \cap V_0=\{\mathrm{Q} \in V_0| \ t_0(\mathrm{Q})=0\}.\]
The inclusion "$\supseteq$" is straightforward. On the other hand, if a point in $E^{an}[D] \cap V_0$ is given, then its inverse image under $(\mathrm{pr}_{|V'_0})^{-1}$ lies in $\mathrm{pr}^{-1}(E^{an}[D])\cap V'_0$, hence is of the form $(\frac{m}{D}\tau+ \frac{n}{D},\tau)$ for some $\tau \in \H$ and $m,n\in \Z$. But because $V'_0$ by choice is contained in $U_0'$ if follows that $m=m_0$ and $n=n_0$ (cf. Lemma 3.8.4). As $\psi_0(\frac{m_0}{D}\tau+ \frac{n_0}{D},\tau)=(0,\tau)$ we have shown the inclusion "$\subseteq$" in $(3.8.5)$.\\
\newline
To proceed with the construction of the chart around $\mathrm{P_0}$ we note that $\mathrm{im}(\psi_0)$ clearly is an open neighborhood of $(0,\tau_0)$ and that we may further shrink $V_0,V'_0$ to obtain a chart around $\mathrm{P_0}$ writing as
\[\tag{\textbf{3.8.6}} (t_0,s_0): V_0 \xrightarrow{\sim} V'_0 \xrightarrow{\sim} B_{\epsilon_0}(0)\times B_{\delta_0}(\tau_0)\footnote{For any $w_0 \in \C$ and $s >0$ we use the standard notation $B_s(w_0)$ for the open ball $B_s(w_0):=\{z\in \C| \ |z-w_0|<s\}$.}\]
for some $\epsilon_0, \delta_0 >0$, where again the first arrow is the inverse of the projection and the second is given by $(z,\tau) \mapsto (z-\frac{m_0}{D}\tau-\frac{n_0}{D},\tau)$. We can assume $\delta_0$ to be so small that $B_{\delta_0}(\tau_0) \subseteq \H$ projects bijectively to an open neighborhood of the point $\pi^{an}(\mathrm{P_0})=[\tau_0] \in S^{an}$.\\
Choose a map as in $(3.8.6)$ with all of the mentioned properties; note from $(3.8.5)$ that $E^{an}[D]\cap V_0$ hereunder corresponds to $(0) \times B_{\delta_0}(\tau_0)$.\\
This is the chart we fix around $\mathrm{P_0} \in E^{an} \cap E^{an}[D]$; it satisfies the conditions recalled in Step 1.\\
\newline
In sum, we have constructed an open covering of $E^{an}$ together with local coordinates adapted for the residue computation $(3.8.4)$.\\
\newline
\underline{Step 3}: (The residue computation)\\
\newline
Now consider the section
\[q_n^D(z,\tau)=\begin{pmatrix} s^D_0(z,\tau) \\ s^D_1(z,\tau) \\ \vdots\\s^D_n(z,\tau) \\ 0 \\ \vdots \\ \vdots \\ 0 \end{pmatrix} \in \Gamma(U_D^{an},\Omega^1_{E^{an}/S^{an}}\otimes_{\mathcal O_{E^{an}}} \mathcal L_n^{an})\]
which induces the element $q^D_n \in \Gamma(S^{an},H^1_{\mathrm{dR}}(U_D^{an}/S^{an},\mathcal L_n^{an}))$ (cf. Def. 3.8.1).\\
Crucial for the calculation of $\mathrm{res}^D_n(q^D_n)$ will be our precise knowledge of the singularity and residue behaviour of the functions $s^D_k(z,\tau),k\geq 0$, namely  (cf. $(3.3.23)$):

\[{\small
\tag{\textbf{3.8.7}} \begin{split} &s^D_k \ \textrm{has at worst simple poles along} \ z=x\tau+y \ (x,y \in \Z, \tau \in \H), \textrm{with residue} \ (D^2-1)\cdot \frac{(2\pi i x)^k}{k!},\\
&\textrm{and along} \ z=\frac{m}{D}\tau+\frac{n}{D} \ (\textrm{with} \ D \ \textrm{not simultaneously dividing} \ m\ \textrm{and} \ n),  \textrm{with residue} \ -\frac{(2\pi i \frac{m}{D})^k}{k!}.
\end{split}
}
\]
Using the coordinate charts of Step 2, the fact that $s_k^D(z,\tau)$ has at worst simple poles along each point $\frac{m_0}{D}\tau+\frac{n_0}{D}$ (with $\tau \in \H, m_0,n_0\in \Z$) implies that
$q_n^D(z,\tau)$ lies in the image of the canonical map
\[\tag{\textbf{3.8.8}} \begin{split}
\Gamma(E^{an},\Omega^1_{E^{an}/S^{an}}(\mathrm{log}(E^{an}[D]))\otimes_{\mathcal O_{E^{an}}} \mathcal L_n^{an}) & \rightarrow \Gamma(E^{an},(j^{an}_D)_*(\Omega^1_{U^{an}_D/S^{an}}\otimes_{U_D^{an}}(\mathcal L_n^{an})_{|U^{an}_D}))\\
&=\Gamma(U_D^{an},\Omega^1_{E^{an}/S^{an}}\otimes_{E^{an}}\mathcal L_n^{an}).
\end{split}\]
For this one needs to be aware of how the vector $q_n^D(z,\tau)$, which a priori is a section in
\[\Gamma(\mathrm{pr}^{-1}(U_D^{an}),\mathrm{pr}^*(\Omega^1_{E^{an}/S^{an}}\otimes_{\mathcal O_{E^{an}}} \mathcal L_n^{an})),\]
gives - by its invariance under deck transformations and our fixed trivialization for the pullback bundle $\mathrm{pr}^*(\Omega^1_{E^{an}/S^{an}}\otimes_{\mathcal O_{E^{an}}} \mathcal L_n^{an})$ - a corresponding section in
\[\Gamma(U_D^{an},\Omega^1_{E^{an}/S^{an}}\otimes_{\mathcal O_{E^{an}}} \mathcal L_n^{an}):\]
namely, by using the locally bijective projection $\mathrm{pr}$ and glueing; then, using the charts in Step 2, one enfolds the definition of logarithmic poles (cf. 3.7.2) and sees immediately that the claim follows from the mentioned property of the $s_k^D(z,\tau)$.\\
\newline
But $(3.8.8)$ implies (by recalling the diagramatic definition of $\mathrm{res}_n^D$ in Step 1 and by straightforward compatibilities) that $\mathrm{res}_n^D(q_n^D)$ is nothing else than the image of $q_n^D(z,\tau)$ under the map
\[\tag{\textbf{3.8.9}} \Gamma(E^{an},\Omega^1_{E^{an}/S^{an}}(\mathrm{log}(E^{an}[D]))\otimes_{\mathcal O_{E^{an}}} \mathcal L_n^{an}) \rightarrow H^0(E^{an}[D],(i_D^{an})^*\mathcal L_n^{an})\]
defined by $(3.8.4)$ in global $E^{an}$-sections.\\
\newline
The last is accessible by a direct computation, using the concrete coordinate charts of Step 2:\\
\newline
Namely, for each point $\mathrm{P_0}\in E^{an}[D]$ with its fixed representative $(\frac{m_0}{D}\tau_0+\frac{n_0}{D},\tau_0)$ in $\C\times \H$ let
\[(t_0,s_0): V_0 \xrightarrow{\sim} V'_0 \xrightarrow{\sim} B_{\epsilon_0}(0) \times B_{\delta_0}(\tau_0)\]
be the local coordinate of $E^{an}$ around $\mathrm{P_0}$, as defined in $(3.8.6)$. By the construction in Step 2 there are two possibilities for the representative of $\mathrm{P_0}$: it is given by $(0,\tau_0)$ iff $\mathrm{P_0} \in \epsilon^{an}(S^{an})\subseteq E^{an}[D]$ and by $(\frac{m_0}{D}\tau_0+\frac{n_0}{D},\tau_0)$, with $D$ not dividing both $m_0$ and $n_0$, otherwise.\\
Now, there exist unique sections $\widetilde{e}$ resp. $\widetilde{f}$ resp. $\widetilde{g}$ of $(\mathcal L_1^{an})_{|V_0}$ over $V_0$ corresponding under $V_0 \xrightarrow{\sim} V'_0$ to the trivializing sections $e_{|V'_0}$ resp. $f_{|V'_0}$ resp. $g_{|V'_0}$ of $(\mathrm{pr}^*\mathcal L_1^{an})_{|V'_0}$.\\
Then the residue formulas in $(3.8.7)$, the definition of the trivialization for $\mathrm{pr}^*(\Omega^1_{E^{an}/S^{an}}\otimes_{\mathcal O_{E^{an}}} \mathcal L_n^{an})$ as well as the definition and property $(3.8.5)$ of our chart imply:\\
$(3.8.9)$ maps $q_n^D(z,\tau)$ to a section in $H^0(E^{an}[D],(i_D^{an})^*\mathcal L_n^{an})$ whose restriction to $E^{an}[D]\cap V_0$ is

\[\begin{cases}
(D^2-1)\cdot \sum_{k=0}^n\frac{(2\pi i\cdot 0)^k}{k!}\cdot \frac{(\widetilde{e}_{|E^{an}[D]\cap V_0})^{n-k} \cdot (\widetilde{f}_{|E^{an}[D]\cap V_0})^k}{(n-k)!}=(D^2-1)\cdot \frac{(\widetilde{e}_{|E^{an}[D]\cap V_0})^n}{n!} \\
\ \ & \ \ resp.\\
-\sum_{k=0}^n\frac{(2\pi i \frac{m_0}{D})^k}{k!}\cdot \frac{(\widetilde{e}_{|E^{an}[D]\cap V_0})^{n-k} \cdot (\widetilde{f}_{|E^{an}[D]\cap V_0})^k}{(n-k)!}
\end{cases}\]
if
\[\begin{cases}
\mathrm{P_0} \in \epsilon^{an}(S^{an}) \subseteq E^{an}[D] \qquad \qquad \qquad \qquad \qquad \qquad \\
\qquad \qquad \qquad \qquad \qquad \qquad \qquad resp. & \ \ \qquad \qquad \qquad \qquad \qquad \qquad \qquad \\
\mathrm{P_0} \in E^{an}[D]-\epsilon^{an}(S^{an}). \qquad \qquad \qquad \qquad \qquad \qquad
\end{cases}\]
\\
As explained at the beginning, the theorem follows if we can show that the image of the preceding element of $H^0(E^{an}[D],(i_D^{an})^*\mathcal L_n^{an})$ under the chain of identifications
\[(i_D^{an})^*\mathcal L_n^{an} \simeq (i_D^{an})^*([D]^{an})^*\mathcal L_n^{an} \simeq  (\pi_{E[D]}^{an})^*(\epsilon^{an})^*\mathcal L_n^{an} \simeq   (\pi_{E[D]}^{an})^* \prod_{k=0}^n \mathrm{Sym}^k_{\mathcal O_{S^{an}}}\mathcal H^{an}\]
is equal to $\frac{1}{n!}\cdot ( D^2\cdot 1_{ \{\epsilon^{an} \}} - 1_{E^{an}[D]})$, the last naturally viewed as global $E^{an}[D]$-section of the outmost right term in the chain by considering the zeroth component of the product.\\
\newline
\underline{Step 4}: (Local considerations)\\
\newline
As we need to verify equality of two well-defined sections of a sheaf on $E^{an}[D]$, we may work locally on $E^{an}[D]$; it hence suffices to construct for each fixed point $\mathrm{P_0} \in E^{an}[D]$ an open neighborhood of $\mathrm{P_0}$ in $E^{an}[D]$ and to show that the two sections in question hereupon coincide.\\
\newline
At first, we have an open neighborhood $V_0$ of $\mathrm{P_0}$ in $E^{an}$, appearing in the definition of the coordinate chart $(3.8.6)$ around $\mathrm{P_0}$:
\[\tag{\textbf{3.8.10}} V_0 \xrightarrow{\sim} V'_0 \xrightarrow{\sim} B_{\epsilon_0}(0)\times B_{\delta_0}(\tau_0)\]
with $E^{an}[D]\cap V_0$ corresponding to $(0) \times B_{\delta_0}(\tau_0)$ and $B_{\delta_0}(\tau_0)\subseteq \H$ projecting bijectively to an open neighborhood of $\pi^{an}(\mathrm{P_0})=[\tau_0] \in S^{an}$. Here, $(\frac{m_0}{D}\tau_0+\frac{n_0}{D},\tau_0)$ is as usual the representative of $\mathrm{P_0}$ in $\C\times \H$, fixed as explained in Step 2.\\
\newline
We obtain the desired open neighborhood of $\mathrm{P_0}$ in $E^{an}[D]$ by possibly shrinking the $V_0$ of $(3.8.10)$ to a smaller open neighborhood of $\mathrm{P_0}$ in $E^{an}$ and then intersecting with $E^{an}[D]$. This goes as follows:\\
As the point $[D]^{an}(\mathrm{P_0})$ is in $\epsilon^{an}(S^{an})\subseteq E^{an}[D]$ it has a chart (defined as in $(3.8.6)$) of the form
\[\tag{\textbf{3.8.11}} W_0 \xrightarrow{\sim} W'_0 = B_{\mu_0}(0)\times B_{\nu_0}(\tau_0)\]
with $E^{an}[D]\cap W_0$ corresponding to $(0) \times B_{\nu_0}(\tau_0)$ and $B_{\nu_0}(\tau_0)\subseteq \H$ projecting bijectively to an open neighborhood of $\pi^{an}([D]^{an}(\mathrm{P_0}))=[\tau_0] \in S^{an}$. As $[D]^{an}(\mathrm{P_0})$ lies in the zero section we have to use (in the definition of the chart, cf. Step 2) its representative $(0,\tau_0)$ in $\C\times \H$, and hence the second arrow in $(3.8.11)$ really is just the identity, by definition of the charts in Step 2 (cf. $(3.8.6)$).\\
We can now clearly shrink $V_0$ to a possibly smaller open neighborhood $T_0$ of $\mathrm{P_0}$ in $E^{an}$ such that we again have a chain as in $(3.8.10)$:
\[\tag{\textbf{3.8.12}} T_0 \xrightarrow{\sim} T'_0 \xrightarrow{\sim} B_{r_0}(0)\times B_{s_0}(\tau_0)\]
with properties as before and additionally satisfying that the image of $B_{r_0}(0)\times B_{s_0}(\tau_0)$ under $(z,\tau) \mapsto (Dz,\tau)$ is contained in $B_{\mu_0}(0)\times B_{\nu_0}(\tau_0)$, i.e. $Dr_0 \leq \mu_0, s_0 \leq \nu_0$.\\
Then $E^{an}[D]\cap T_0$ is the open neighborhood of $\mathrm{P_0}$ in $E^{an}[D]$ which we wanted to construct.\\
Note that $T_0$ is contained in the coordinate region $V_0$ of $\mathrm{P_0}$ which was fixed in Step 2 and used for the residue computation in Step 3.\\
\newline
We use $(3.8.12)$ and its properties to get a chart for $\mathrm{P_0}$ in $E^{an}$ resp. in $E^{an}[D]$ as well as for $\pi^{an}(\mathrm{P_0})=[\tau_0]$ in $S^{an}$; we use $(3.8.11)$ to get a chart for $[D]^{an}(\mathrm{P_0})=\epsilon^{an}([\tau_0])=[(0,\tau_0)]$ in $E^{an}$.\\
\newline
To avoid confusion here: the residue computation in Step 3 was done with a fixed family of charts covering all of $E^{an}$, and this computation is over. The four charts introduced one moment ago serve a different purpose and are entirely tied to the fixed point $\mathrm{P_0} \in E^{an}[D]$.\\
\newline
One then verifies directly from the definitions that in these charts the maps
\begin{align*}
\pi^{an}&: E^{an}\rightarrow S^{an}\\
\epsilon^{an}&: S^{an}\rightarrow E^{an}\\
i_D^{an}&: E^{an}[D]\rightarrow E^{an}\\
\pi_{E[D]}^{an}&: E^{an}[D]\rightarrow S^{an}\\
[D]^{an}&: E^{an}\rightarrow E^{an}
\end{align*}
express as
\begin{align*}
\pi^{an}&: B_{r_0}(0)\times B_{s_0}(\tau_0) \rightarrow B_{s_0}(\tau_0), \ (z,\tau) \mapsto \tau\\
\epsilon^{an}&:B_{s_0}(\tau_0) \rightarrow B_{\mu_0}(0)\times B_{\nu_0}(\tau_0), \ \tau \mapsto (0,\tau)\\
i_D^{an}&:B_{s_0}(\tau_0) \rightarrow B_{r_0}(0)\times B_{s_0}(\tau_0), \ \tau \mapsto (0,\tau)\\
\pi_{E[D]}^{an}&:B_{s_0}(\tau_0) \rightarrow B_{s_0}(\tau_0), \ \tau \mapsto \tau\\
[D]^{an}&:B_{r_0}(0)\times B_{s_0}(\tau_0) \rightarrow B_{\mu_0}(0)\times B_{\nu_0}(\tau_0), \ (z,\tau) \mapsto (Dz,\tau)
\end{align*}
\\
\underline{Step 5}: (Coordinate transport of the situation)\\
\newline
As explained at the end of Step 3 and the beginning of Step 4 it now remains to show that
\[\begin{cases}
(D^2-1)\cdot \frac{(\widetilde{e}_{|E^{an}[D]\cap T_0})^n}{n!} \\
\ \ & \ \ resp.\\
-\sum_{k=0}^n\frac{(2\pi i \frac{m_0}{D})^k}{k!}\cdot \frac{(\widetilde{e}_{|E^{an}[D]\cap T_0})^{n-k} \cdot (\widetilde{f}_{|E^{an}[D]\cap T_0})^k}{(n-k)!}
\end{cases}\]
under the following chain of maps evaluated in $(E^{an}[D]\cap T_0)$-sections:
\[\tag{\textbf{3.8.13}} (i_D^{an})^*\mathcal L_n^{an} \simeq (i_D^{an})^*([D]^{an})^*\mathcal L_n^{an} \simeq  (\pi_{E[D]}^{an})^*(\epsilon^{an})^*\mathcal L_n^{an} \simeq   (\pi_{E[D]}^{an})^* \prod_{k=0}^n \mathrm{Sym}^k_{\mathcal O_{S^{an}}}\mathcal H^{an}\]
goes to $\frac{1}{n!}\cdot ( D^2\cdot 1_{ \{\epsilon^{an} \}\cap T_0} - 1_{E^{an}[D]\cap T_0})$; recall here that $T_0 \subseteq V_0$.\\
\newline
We transport the whole situation via the four charts introduced at the very end of Step 4.\\
\newline
Write
\[\mathsf L_n, \mathsf e, \mathsf f, \mathsf g \quad \textrm{resp.}\quad \mathsf H\]
to denote the pullback to $B_{r_0}(0)\times B_{s_0}(\tau_0)$ resp. to $B_{s_0}(\tau_0)$ of
\[(\mathcal L_n^{an})_{|T_0}, \widetilde{e}_{|T_0}, \widetilde{f}_{|T_0}, \widetilde{g}_{|T_0} \quad \textrm{resp. of (the adapted restriction of)} \ \mathcal H^{an}\]
via the chart around $\mathrm{P_0}$ in $E^{an}$ resp. around $\pi^{an}(\mathrm{P_0})=[\tau_0]$ in $S^{an}$ introduced in Step 4.\\
Recalling how $\widetilde{e}, \widetilde{f}, \widetilde{g}$ and $(3.8.12)$ were defined (cf. Steps 3 and 4) we can say this differently:\\
Namely, $\mathsf L_n, \mathsf e, \mathsf f, \mathsf g$ are the pullback of $(\mathrm{pr}^*\mathcal L_n^{an})_{|T'_0}, e_{|T'_0},f_{|T'_0},g_{|T'_0}$ along
\[B_{r_0}(0)\times B_{s_0}(\tau_0) \xrightarrow{\sim} T'_0, \quad (z,\tau)\mapsto \Big(z+\frac{m_0}{D}\tau+\frac{n_0}{D},\tau \Big).\]
Note that because the second arrow in $(3.8.11)$ is the identity we don't have to introduce extra notation when performing the analogous procedure with respect to the chart around $D\mathrm{P_0}$ in $E^{an}$.\\
For the induced maps between the various coordinate regions (recorded explicitly at the end of Step 4) we simply keep notations unchanged; note in particular that then $\pi_{E[D]}^{an}=\mathrm{id}_{B_{s_0}(\tau_0)}$.\\
\newline
The restriction of the chain $(3.8.13)$ to the coordinate region $E^{an}[D]\cap T_0 \subseteq E^{an}[D]$ then transforms into the following series of isomorphisms of $\mathcal O_{B_{s_0}(\tau_0)}$-modules:

\[\tag{\textbf{3.8.14}} \begin{split}
&(i_D^{an})^*\mathsf L_n \simeq (i_D^{an})^*([D]^{an})^*((\mathrm{pr}^*\mathcal L_n)_{|B_{\mu_0}(0)\times B_{\nu_0}(\tau_0)}) \simeq \\
&\simeq (\epsilon^{an})^*((\mathrm{pr}^*\mathcal L_n)_{|B_{\mu_0}(0)\times B_{\nu_0}(\tau_0)}) \simeq \prod_{k=0}^n \mathrm{Sym}^k_{\mathcal O_{B_{s_0}(\tau_0)}}\mathsf H.
\end{split}\]

We claim that under the first isomorphism
\[(i_D^{an})^*\mathsf L_1 \simeq (i_D^{an})^*([D]^{an})^*((\mathrm{pr}^*\mathcal L_1)_{|B_{\mu_0}(0)\times B_{\nu_0}(\tau_0)})\]
of $(3.8.14)$ for $n=1$ the section $(i_D^{an})^*(\mathsf e)$ resp. $(i_D^{an})^*(\mathsf f)$ of the left side corresponds to the section
\[(i_D^{an})^*([D]^{an})^*(e_{|B_{\mu_0}(0)\times B_{\nu_0}(\tau_0)}-2\pi i m_0 \cdot f_{|B_{\mu_0}(0)\times B_{\nu_0}(\tau_0)})\]
resp.
\[(i_D^{an})^*([D]^{an})^*(D\cdot f_{|B_{\mu_0}(0)\times B_{\nu_0}(\tau_0)})\]
of the right side.\\
We now insert a digression in which we give a detailed justification for the previous claim; it may be convenient for the stream of reading if one first jumps to Step 7 to see how the proof continues and later returns to the arguments of Step 6.\\
\newline
\underline{Step 6}: (Insertion I)\\
\newline
To show the claim at the end of Step 5 it suffices to see that under the isomorphism of $\mathcal O_{B_{r_0}(0)\times B_{s_0}(\tau_0)}$-modules
\[\mathsf L_1\simeq ([D]^{an})^*((\mathrm{pr}^*\mathcal L_1)_{|B_{\mu_0}(0)\times B_{\nu_0}(\tau_0)}),\]
naturally induced\footnote{This means that one restricts the isomorphism $\mathcal L_1^{an} \simeq ([D]^{an})^*\mathcal L_1^{an}$ to $T_0$, takes its pullback along the coordinate chart $T_0 \simeq B_{r_0}(0)\times B_{s_0}(\tau_0)$ of $\mathrm{P_0}$ in $E^{an}$ and uses the following commutative diagram with right vertical arrow given by the chart $(3.8.11)$ of $D\mathrm{P_0}$ in $E^{an}$; note that this last chart (by contrast to the first one) is just induced by the projection $\mathrm{pr}$.
\begin{equation*}
\begin{xy}
\xymatrix{
B_{r_0}(0)\times B_{s_0}(\tau_0) \ar[r]^{[D]^{an}} & B_{\mu_0}(0)\times B_{\nu_0}(\tau_0)\\
T_0 \ar[r]^{[D]^{an}}\ar@{-}[u]_{\sim} &  W_0 \ar@{-}[u]^{\sim}_{\mathrm{pr}}.}
\end{xy}
\end{equation*}
} by the invariance isomorphism
\[\tag{\textbf{3.8.15}} \mathcal L_1^{an} \simeq ([D]^{an})^*\mathcal L_1^{an},\]
the section $\mathsf e$ resp. $\mathsf f$ corresponds to
\[([D]^{an})^*(e_{|B_{\mu_0}(0)\times B_{\nu_0}(\tau_0)}-2\pi i m_0 \cdot f_{|B_{\mu_0}(0)\times B_{\nu_0}(\tau_0)})\]
resp.
\[([D]^{an})^*(D\cdot f_{|B_{\mu_0}(0)\times B_{\nu_0}(\tau_0)}).\]
But we have a big commutative diagram
\begin{equation*}
\begin{xy}
\xymatrix{
B_{r_0}(0)\times B_{s_0}(\tau_0) \ar[r]^{[D]^{an}} & B_{\mu_0}(0)\times B_{\nu_0}(\tau_0)\\
T'_0 \ar[r]^{\lambda  \ \ \qquad} \ar@{-}[u]_{\sim} & B_{\mu_0}(0)\times B_{\nu_0}(\tau_0) \ar@{-}[u]_{\id}\\
T_0 \ar[r]^{[D]^{an}}\ar@{-}[u]_{\sim}^{\mathrm{pr}} &  W_0 \ar@{-}[u]^{\sim}_{\mathrm{pr}}
}
\end{xy}
\end{equation*}
in which the vertical arrows are given by $(3.8.12)$ resp. $(3.8.11)$ and the map $\lambda$ is defined by
\[\lambda(z,\tau):=(Dz-m_0\tau-n_0,\tau).\]
It is hence enough to show that under the isomorphism
\[(\mathrm{pr}^*\mathcal L_1^{an})_{|T'_0} \simeq \lambda^*((\mathrm{pr}^*\mathcal L_1)_{|B_{\mu_0}(0)\times B_{\nu_0}(\tau_0)}),\]
induced by $(3.8.15)$ and the lower commutative square, the section $e_{|T'_0}$ resp. $f_{|T'_0}$ corresponds to
\[\lambda^*(e_{|B_{\mu_0}(0)\times B_{\nu_0}(\tau_0)}-2\pi i m_0 \cdot f_{|B_{\mu_0}(0)\times B_{\nu_0}(\tau_0)})\]
resp.
\[\lambda^*(D \cdot f_{|B_{\mu_0}(0)\times B_{\nu_0}(\tau_0)}),\]
because one can then apply pullback via the upper left vertical arrow of the diagram.\\
\newline
The last in turn follows if we can verify the following statement: under the isomorphism
\[\mathrm{pr}^*\mathcal L_1^{an} \simeq \sigma^*\mathrm{pr}^*\mathcal L_1^{an},\]
induced by $(3.8.15)$ and the commutative diagram
\begin{equation*}
\begin{xy}
\xymatrix{
\C\times \H \ar[d]_{\mathrm{pr}} \ar[r]^{\sigma} & \C\times \H \ar[d]^{\mathrm{pr}}\\
E^{an} \ar[r]^{[D]^{an}} &  E^{an}
}
\end{xy}
\end{equation*}
the section $e$ resp. $f$ corresponds to
\[\sigma^*(e-2\pi i m_0 \cdot f)\]
resp.
\[\sigma^*(D\cdot f);\]
here, we set
\[\sigma(z,\tau):=(Dz-m_0\tau-n_0,\tau).\]
We already encountered this kind of problem during the proof of Prop. 3.5.9, and using formula $(**)$ of that proof for $N$ replaced by $D$ (justified in the supplements subsequent to Prop. 3.5.9) our present claim is verified as follows:\\
\newline
First, one recalls the definition of the trivializing section $e$ resp. $f$ of $\mathrm{pr}^*\mathcal L_1^{an}$ over $\C\times \H$:
\[e=k_1^{-1}(t)\otimes 1 \quad \textrm{resp.} \quad f=k^{-1}_1(wt)\otimes 1,\]
where one uses the $\mathcal O_{\C\times \H}$-linear isomorphism of $(3.5.6)$:
\[k_1^{-1}\widetilde{\mathcal P}\otimes_{k_1^{-1}\mathcal O_{\C\times \C^2 \times \H}}(\mathcal O_{\C\times \H} \oplus \mathcal O_{\C\times \H}\cdot \eta^\vee \oplus \mathcal O_{\C\times \H}\cdot \omega^\vee)=k_1^*\widetilde{\mathcal P}\simeq \tau_1^*\mathcal P_1^{an}\simeq \mathrm{pr}^*\mathcal L_1^{an}\]
and the trivializing section $t$ of $\widetilde{\mathcal P}$ on $\C\times \C^2\times \H$:
\[\mathcal O_{\C \times \C^2 \times \H} \xrightarrow{\sim} \widetilde{\mathcal P} \quad 1\mapsto t=\frac{1}{J(z,-w,\tau)}\otimes \omega_{\mathrm{can}}^\vee.\]
Second, if one sets
\[\varphi(z,\tau):=(z-m_0\tau-n_0,\tau) \quad \textrm{and} \quad \widetilde{[D]}(z,\tau):=(Dz,\tau),\]
then $\sigma=\varphi \circ \widetilde{[D]}$ and the isomorphism
\[\sigma^*\mathrm{pr}^*\mathcal L_1^{an} \simeq \mathrm{pr}^*\mathcal L_1^{an}\]
of above equals the composition
\[\sigma^*\mathrm{pr}^*\mathcal L_1^{an} \simeq \widetilde{[D]}^*\varphi^*\mathrm{pr}^*\mathcal L_1^{an}\simeq \widetilde{[D]}^*\mathrm{pr}^*\mathcal L_1^{an} \simeq \mathrm{pr}^*\mathcal L_1^{an},\]
where the second identification is due to $\mathrm{pr}\circ \varphi = \mathrm{pr}$ and the third is induced by $(3.8.15)$ together with the commutative diagram
\begin{equation*}
\begin{xy}
\xymatrix{
\C\times \H \ar[d]_{\mathrm{pr}} \ar[r]^{\widetilde{[D]}} & \C\times \H \ar[d]^{\mathrm{pr}}\\
E^{an} \ar[r]^{[D]^{an}} &  E^{an}
}
\end{xy}
\end{equation*}
The section $\sigma^*(e)$ of $\sigma^*\mathrm{pr}^*\mathcal L_1^{an}$ maps to $\widetilde{[D]}^*(e+2\pi i m_0\cdot f)$ under the first two arrows of the previous chain because of the equality
\[\frac{1}{J(Dz-m_0\tau-n_0,-w,\tau)}\otimes \omega^\vee_{\mathrm{can}}=\mathrm{e}^{2\pi i m_0 w}\cdot \frac{1}{J(Dz,-w,\tau)}\otimes \omega^\vee_{\mathrm{can}}\]
and the definition of $k_1$. Furthermore, in $(**)$ of Prop. 3.5.9 (and the supplements after its proof) we have seen that $\widetilde{[D]}^*(e+2\pi i m_0\cdot f)$ goes to $e+2\pi i \frac{m_0}{D}\cdot f$ under the last arrow of the chain. Arguing analogously for $f$ one sees that $\sigma^*(f)$ goes to $\frac{1}{D}\cdot f$ under the chain.\\
\newline
Hence, under
\[\sigma^*\mathrm{pr}^*\mathcal L_1^{an} \xrightarrow{\sim} \mathrm{pr}^*\mathcal L_1^{an}\]
the section $\sigma^*(e)$ resp. $\sigma^*(f)$ goes to $e+2\pi i \frac{m_0}{D} \cdot f$ resp. $\frac{1}{D}\cdot f$. This clearly implies that under
\[\mathrm{pr}^*\mathcal L_1^{an} \xrightarrow{\sim} \sigma^*\mathrm{pr}^*\mathcal L_1^{an}\]
the section $e$ resp. $f$ goes to $\sigma^*(e-2\pi i m_0 \cdot f)$ resp. $\sigma^*(D\cdot f)$, as desired.\\
\newline
\underline{Step 7}: (The images under the final isomorphism)\\
\newline
Resuming the discussion prior to Step 6 we get that the image of $(i_D^{an})^*(\mathsf e)$ resp. $(i_D^{an})^*(\mathsf f)$ under
\[(i_D^{an})^*\mathsf L_1 \simeq (i_D^{an})^*([D]^{an})^*((\mathrm{pr}^*\mathcal L_1)_{|B_{\mu_0}(0)\times B_{\nu_0}(\tau_0)}) \simeq (\epsilon^{an})^*((\mathrm{pr}^*\mathcal L_1)_{|B_{\mu_0}(0)\times B_{\nu_0}(\tau_0)})\]
is given by
\[(\epsilon^{an})^*(e_{|B_{\mu_0}(0)\times B_{\nu_0}(\tau_0)}-2\pi i m_0 \cdot f_{|B_{\mu_0}(0)\times B_{\nu_0}(\tau_0)})\]
resp. by
\[(\epsilon^{an})^*(D \cdot f_{|B_{\mu_0}(0)\times B_{\nu_0}(\tau_0)}).\]
We claim that the images of these sections under the final isomorphism
\[\tag{\textbf{3.8.16}} (\epsilon^{an})^*((\mathrm{pr}^*\mathcal L_1)_{|B_{\mu_0}(0)\times B_{\nu_0}(\tau_0)}) \simeq \mathcal O_{B_{s_0}(\tau_0)} \oplus \mathsf H\]
of the chain $(3.8.14)$ are the sections $1_{B_{s_0}(\tau_0)}-2\pi i m_0 \cdot (\eta^{\vee})_{|B_{s_0}(\tau_0)}$ resp. $D\cdot (\eta^{\vee})_{|B_{s_0}(\tau_0)}$; note here that if we write
\[p: \H \rightarrow S^{an}=\Gamma(N)\backslash \H\]
for the projection, then $\mathsf H$ equals $(p^*\mathcal H^{an})_{|B_{s_0}(\tau_0)}$, and recall that $\eta^{\vee}$ is one of the two fixed trivializing sections $\{\eta^\vee, \omega^\vee\}$ for $p^*\mathcal H^{an}$ (cf. $(3.5.1)$).\\
We want to give a detailed verification of the preceding claim and thus include another intermediate step in which this is carried out; again, it might be reasonable to first proceed directly with the conclusion of the proof in Step 9 and to recur to Step 8 later.\\
\newline
\underline{Step 8}: (Insertion II)\\
\newline
To show the preceding claim note first that $(3.8.16)$ arises by construction from the isomorphisms
\[\tag{\textbf{3.8.17}} (\epsilon^{an})^*((\mathrm{pr}^*\mathcal L_1)_{|B_{\mu_0}(0)\times B_{\nu_0}(\tau_0)}) \simeq (p^*(\epsilon^{an})^*\mathcal L_1^{an})_{|B_{s_0}(\tau_0)}\]
and
\[\tag{\textbf{3.8.18}} (p^*(\epsilon^{an})^*\mathcal L_1^{an})_{|B_{s_0}(\tau_0)} \simeq (p^*(\mathcal O_{S^{an}}\oplus \mathcal H^{an}))_{|B_{s_0}(\tau_0)}.\]
Here, $(3.8.17)$ is the natural identification, in an obvious way due to the commutative diagram
\begin{equation*}
\begin{xy}
\xymatrix{
B_{s_0}(\tau_0) \ar[r]^{\epsilon^{an} \qquad} \ar[d]_{p}^{\sim} & B_{\mu_0}(0)\times B_{\nu_0}(\tau_0) \ar[d]^{\mathrm{pr}}_{\sim}\\
p(B_{s_0}(\tau_0)) \ar[r]^{\qquad \epsilon^{an}} &  W_0}
\end{xy}
\end{equation*}
and $(3.8.18)$ comes from pulling back via $p$ the splitting $(\epsilon^{an})^*\mathcal L_1^{an} \simeq \mathcal O_{S^{an}}\oplus \mathcal H^{an}$ and then restricting to $B_{s_0}(\tau_0)$. This and the fact that the upper arrow $\epsilon^{an}$ of the preceding diagram is $\epsilon^{an}(\tau)=(0,\tau)$ (by Step 4) directly implies that one obtains $(3.8.16)$ alternatively as follows: take the pullback of the splitting $(\epsilon^{an})^*\mathcal L_1^{an} \simeq \mathcal O_{S^{an}}\oplus \mathcal H^{an}$ along $p$ and use the commutative diagram
\begin{equation*}
\begin{xy}
\xymatrix{
\H \ar[r]^{\widetilde{\epsilon} \quad} \ar[d]_{p} & \C \times \H \ar[d]^{\mathrm{pr}}\\
S^{an} \ar[r]^{\epsilon^{an}} &  E^{an}}
\end{xy}
\end{equation*}
where $\widetilde{\epsilon}(\tau):=(0,\tau)$, to get
\[\tag{\textbf{3.8.19}} (\widetilde{\epsilon})^*\mathrm{pr}^*\mathcal L_1^{an} \simeq \mathcal O_{\H}\oplus p^*\mathcal H^{an}.\]
Restricting this to $B_{s_0}(\tau_0)$ then yields (after obvious canonical identifications which we don't write out) the isomorphism $(3.8.16)$.\\
Hence, to deduce our claim it clearly suffices to see that $(3.8.19)$ sends the section $(\widetilde{\epsilon})^*(e)$ resp. $(\widetilde{\epsilon})^*(f)$ to $1_{\H}$ resp. to $\eta^{\vee}$. This in turn was shown in the proof of Prop. 3.5.8.\\
\newline
\underline{Step 9}: (Conclusion of the proof)\\
\newline
Altogether, we have shown in the last four steps that the chain $(3.8.14)$ for $n=1$:
\[(i_D^{an})^*\mathsf L_1 \simeq (i_D^{an})^*([D]^{an})^*((\mathrm{pr}^*\mathcal L_1)_{|B_{\mu_0}(0)\times B_{\nu_0}(\tau_0)}) \simeq (\epsilon^{an})^*((\mathrm{pr}^*\mathcal L_1)_{|B_{\mu_0}(0)\times B_{\nu_0}(\tau_0)}) \simeq \mathcal O_{B_{s_0}(\tau_0)} \oplus \mathsf H\]
sends the section $(i_D^{an})^*(\mathsf e)$ resp. $(i_D^{an})^*(\mathsf f)$ to  $1_{B_{s_0}(\tau_0)}-2\pi i m_0 \cdot (\eta^{\vee})_{|B_{s_0}(\tau_0)}$ resp. to $D\cdot (\eta^{\vee})_{|B_{s_0}(\tau_0)}$.\\
\newline
Now observe that under our chart $E^{an}[D]\cap T_0 \simeq B_{s_0}(\tau_0)$ the sheaves $((i_D^{an})^*(\mathcal L_1^{an}))_{|E^{an}[D]\cap T_0}$ and $(i_D^{an})^*\mathsf L_1$ as well as their sections $\widetilde{e}_{|E^{an}[D]\cap T_0}, \widetilde{f}_{|E^{an}[D]\cap T_0}$ and $(i_D^{an})^*(\mathsf e), (i_D^{an})^*(\mathsf f)$ correspond to each other. Likewise, the sheaves $\mathcal O_{E^{an}[D]\cap T_0} \oplus (\mathcal H^{an}_{E^{an}[D]})_{|E^{an}[D]\cap T_0}$ and $\mathcal O_{B_{s_0}(\tau_0)} \oplus \mathsf H$ as well as their sections $1_{E^{an}[D] \cap T_0}, \upsilon$ and $1_{B_{s_0}(\tau_0)}, (\eta^{\vee})_{|B_{s_0}(\tau_0)}$ correspond.\\
Here, the auxiliary notation $\upsilon$ is used for the $(E^{an}[D]\cap T_0)$-section of $(\mathcal H^{an}_{E^{an}[D]})=(\pi_{E[D]}^{an})^*\mathcal H^{an}$ induced by the pullback of $(\eta^{\vee})_{|B_{s_0}(\tau_0)}$ via
\[E^{an}[D]\cap T_0 \xrightarrow{\pi_{E[D]}^{an}} p(B_{s_0}(\tau_0)) \xrightarrow{\sim} B_{s_0}(\tau_0).\]
All of this is checked directly from the definitions.\\
\newline
We thus conclude that the chain
\[(i_D^{an})^*\mathcal L_1^{an} \simeq (i_D^{an})^*([D]^{an})^*\mathcal L_1^{an} \simeq  (\pi_{E[D]}^{an})^*(\epsilon^{an})^*\mathcal L_1^{an} \simeq  \mathcal O_{E^{an}[D]}\oplus \mathcal H^{an}_{E^{an}[D]}\]
maps the $(E^{an}[D]\cap T_0)$-section $\widetilde{e}_{|E^{an}[D]\cap T_0}$ resp. $\widetilde{f}_{|E^{an}[D]\cap T_0}$ to $1_{E^{an}[D] \cap T_0}-2\pi i m_0 \cdot \upsilon$ resp. to $D\cdot \upsilon$.\\
\newline
We can now finally determine the image of
\[\begin{cases}
(D^2-1)\cdot \frac{(\widetilde{e}_{|E^{an}[D]\cap T_0})^n}{n!} \\
\ \ & \ \ resp.\\
-\sum_{k=0}^n\frac{(2\pi i \frac{m_0}{D})^k}{k!}\cdot \frac{(\widetilde{e}_{|E^{an}[D]\cap T_0})^{n-k} \cdot (\widetilde{f}_{|E^{an}[D]\cap T_0})^k}{(n-k)!}
\end{cases}\]
under the chain $(3.8.13)$ in $(E^{an}[D]\cap T_0)$-sections: namely, with what we have just deduced it is:
\[\begin{cases}
 (D^2-1)\cdot \frac{(1_{E^{an}[D] \cap T_0}-2\pi i m_0 \cdot \upsilon)^n}{n!} \\
\ \ & \ \ resp. \\
-\sum_{k=0}^n\frac{(2\pi i \frac{m_0}{D})^k}{k!}\cdot \frac{(1_{E^{an}[D] \cap T_0}-2\pi i m_0 \cdot \upsilon)^{n-k}\cdot (D\cdot \upsilon)^k}{(n-k)!}.
\end{cases}\]
Recall from Step 3 that the two cases distinguish the situation that $\mathrm{P_0} \in \epsilon^{an}(S^{an})$ resp. that\\
$\mathrm{P_0} \in E^{an}[D] - \epsilon^{an}(S^{an})$; recall moreover from Step 2 that in terms of the (permanently fixed) representative $(\frac{m_0}{D}\tau_0+\frac{n_0}{D},\tau_0)$ of $\mathrm{P_0}$ this is equivalent to $m_0=0=n_0$ resp. to $D$ not dividing simultaneously $m_0$ and $n_0$. Hence, the first entry of the preceding bracket is
\[(D^2-1)\cdot \frac{(1_{E^{an}[D] \cap T_0}-2\pi i \cdot 0 \cdot \upsilon)^n}{n!}=\frac{(D^2-1)}{n!} \cdot 1_{E^{an}[D] \cap T_0},\]
and the second entry is computed as
\begin{align*}
&-\sum_{k=0}^n\frac{(2\pi i \frac{m_0}{D})^k}{k!}\cdot \frac{(1_{E^{an}[D] \cap T_0}-2\pi i m_0 \cdot \upsilon)^{n-k}\cdot (D\cdot \upsilon)^k}{(n-k)!}\\
=&-\frac{1}{n!}\cdot \sum_{k=0}^n \binom{n}{k}\cdot (2\pi i m_0\cdot \upsilon)^k\cdot (1_{E^{an}[D] \cap T_0}-2\pi i m_0 \cdot \upsilon)^{n-k}\\
=&-\frac{1}{n!}\cdot 1_{E^{an}[D] \cap T_0}.\end{align*}
Observe furthermore that $E^{an}[D]\cap T_0$ is contained in $\epsilon^{an}(S^{an})$ resp. in $E^{an}[D] - \epsilon^{an}(S^{an})$.\\
After a short reflection we see that we have obtained precisely the section
\[\frac{1}{n!}\cdot ( D^2\cdot 1_{ \{\epsilon^{an} \}\cap T_0} - 1_{E^{an}[D]\cap T_0}),\]
which finishes the proof, as explained at the beginning of Step 5.
\end{proof}
During the proof of the preceding theorem (namely in Step 2) we used:
\begin{lemma}
Let $\tau_0 \in \H$ and $m_0,n_0 \in \Z$. Then there exists an open neighborhood $U_0'$ in $\C\times \H$ of the point $(\frac{m_0}{D}\tau_0+ \frac{n_0}{D},\tau_0)\in \C\times \H$ with the following property:\\
If $\tau \in \H$ and $m,n \in \Z$ such that $(\frac{m}{D}\tau+ \frac{n}{D},\tau)\in U_0'$, then $m=m_0$ and $n=n_0$.
\end{lemma}
\begin{proof}
Using that the $D$-torsion points of the real torus $\Z^2\backslash \R^2$ lie discrete, one easily finds an open neighborhood $S_0$ of $(\frac{m_0}{D},\frac{n_0}{D})$ in $\R^2$ such that
\[S_0\cap \frac{1}{D}\cdot \Z^2=  \Big\{ \Big(\frac{m_0}{D},\frac{n_0}{D} \Big) \Big\}.\]
Let $U_0'$ be the inverse image of $S_0\times \H \subseteq \R^2\times \H$ under the homeomorphism
\[\C\times \H \xrightarrow{\sim} \R^2\times \H, \quad (x\tau+y,\tau)\mapsto ((x,y),\tau).\]
It is readily checked that $(\frac{m_0}{D}\tau_0+ \frac{n_0}{D},\tau_0) \in U_0'$ and that $U_0'$ satisfies the property in the claim.
\end{proof}

\subsection{The specialization of the $D$-variant of the polylogarithm along torsion sections}
In this subsection we assume that the integer $D>1$ additionally satisfies $(D,N)=1$.\\
Furthermore, we fix two integers $a,b$ which are not simultaneously divisible by $N$.

\subsubsection{Formulation of the problem}
Using the Drinfeld basis $(e_1,e_2) \in E[N](S)$ for $E[N]$ we have the $N$-torsion section
\[t_{a,b}=ae_1+be_2: S \rightarrow E\]
which by our assumptions on $D,a,b$ factors over the open subscheme $U_D=E-E[D]$:
\[t_{a,b}: S \rightarrow U_D \subseteq E.\]
We can then "specialize" $\varD$ along $t_{a,b}$:
\begin{definition}
(i) For each $n\geq 0$ we let
\[t_{a,b}^*(\varDn) \in H^1_{\mathrm{dR}}\Big(S/\Q,\prod_{k=0}^n \mathrm{Sym}^k_{\mathcal O_S}\mathcal H \Big)\]
be the de Rham cohomology class obtained by pulling back
\[\varDn \in H^1_{\mathrm{dR}}(U_D/\Q, \mathcal L_n)\]
along $t_{a,b}$ and by using the (horizontal) identifications
\[t_{a,b}^*\mathcal L_n \simeq \epsilon^*\mathcal L_n \simeq \prod_{k=0}^n \mathrm{Sym}^k_{\mathcal O_S}\mathcal H,\]
where the first comes from Lemma 1.4.5 and the second is the splitting $\varphi_n$.\\
(ii) For each $n\geq 0$ we let
\[\big(t_{a,b}^*(\varDn)\big)^{(n)} \in H^1_{\mathrm{dR}}(S/\Q,\mathrm{Sym}^n_{\mathcal O_S}\mathcal H)\]
be the $n$-th component of $t_{a,b}^*(\varDn)$.
\end{definition}

\begin{remark}
Let $n \geq 0$ and $0 \leq k \leq n$. Then, if we more generally let $\big(t_{a,b}^*(\varDn)\big)^{(k)}$ be the $k$-th component of $t_{a,b}^*(\varDn)$, we have the equality
\[\big(t_{a,b}^*(\varDn)\big)^{(k)} = \frac{1}{(n-k)!}\cdot \big(t_{a,b}^*(\varDk)\big)^{(k)}\]
in $H^1_{\mathrm{dR}}(S/\Q,\mathrm{Sym}^k_{\mathcal O_S}\mathcal H)$. This is straightforwardly derived from Lemma 1.1.5.
\end{remark}

In order to relate the specialization of the $D$-variant to cohomology classes of modular forms it is necessary to make one further identification: namely, the
perfect alternating pairing
\[H^1_{\mathrm{dR}}(E/S)\otimes_{\mathcal O_S} H^1_{\mathrm{dR}}(E/S) \rightarrow \mathcal O_S,\]
coming from composition of the cup product and the trace isomorphism $\mathrm{tr}:H^2_{\mathrm{dR}}(E/S) \xrightarrow{\sim}\mathcal O_S$, furnishes a horizontal isomorphism (Poincaré duality):
\[\tag{\textbf{3.8.20}} H^1_{\mathrm{dR}}(E/S) \xrightarrow{\sim} \mathcal H, \quad x \mapsto \{\ y\mapsto \mathrm{tr}(x\cup y)\}.\]
For these facts cf. also the beginning of Chapter 1.\\
The identification $(3.8.20)$ and the induced isomorphism on symmetric powers will henceforth remain fixed; in particular, we will (without change of notation) consider the classes of Def. 3.8.5 as elements
\[t_{a,b}^*(\varDn) \in H^1_{\mathrm{dR}}\Big(S/\Q,\prod_{k=0}^n \mathrm{Sym}^k_{\mathcal O_S} H^1_{\mathrm{dR}}(E/S)\Big)\]
resp.
\[\big(t_{a,b}^*(\varDn)\big)^{(n)} \in H^1_{\mathrm{dR}}(S/\Q,\mathrm{Sym}^n_{\mathcal O_S} H^1_{\mathrm{dR}}(E/S)).\]
The goal of this final subsection then is to derive a concrete description for the preceding classes by showing that they come from certain algebraic modular forms which we will determine explicitly.\\
Note that in view of Rem. 3.8.6 it suffices to consider the classes $\big(t_{a,b}^*(\varDn)\big)^{(n)}$.

\subsubsection{Outline of the main result and strategy of proof}

The Hodge filtration and Kodaira-Spencer map induce a canonical homomorphism
\[\Gamma\big(S, \omega_{E/S}^{\otimes(n+2)}\big) \rightarrow H^1_{\mathrm{dR}}(S/\Q,\mathrm{Sym}^n_{\mathcal O_S}H^1_{\mathrm{dR}}(E/S))\]
which associates to a weakly holomorphic algebraic modular form of weight $n+2$ and level $N$ a class in the de Rham cohomology of the modular curve in the required coefficient sheaf.\\
Our main theorem then will assert that under this arrow the class $\big(t_{a,b}^*(\varDn)\big)^{(n)}$ is realized by the algebraic modular form
\[\begin{cases}
-{_D}F^{(2)}_{\frac{a}{N},\frac{b}{N}} \ \ &\textrm{if} \ \ n=0
\vspace{0.9mm}\\
\frac{(-1)^{n}}{n!} \cdot {_D}F^{(n+2)}_{\frac{a}{N},\frac{b}{N}} \ \ &\textrm{if} \ \ n> 0
\end{cases}
\]
introduced by Kato in \cite{Ka}, Ch. I, 4.2 resp. 3.6. There, it basically appears in the construction of the "zeta elements" in the space of modular forms; the periods of these elements are intimately related to values of operator-valued zeta functions (cf. ibid., Ch. II, 4, esp. 4.2, 4.5 and Thm. 4.6).\\
\newline
Our strategy to prove this result consists in first transporting the question to the analytic category by moving in the lower row of a natural commutative diagram
\begin{equation*}
\begin{xy}
\xymatrix@C-0.3cm{
\Gamma\big(S, \omega_{E/S}^{\otimes(n+2)}\big)\ar[r] \ar[d] & H^1_{\mathrm{dR}}(S/\Q,\mathrm{Sym}^n_{\mathcal O_S}H^1_{\mathrm{dR}}(E/S)) \ar[d]\\
\Gamma\big(S^{an}, \omega_{E^{an}/S^{an}}^{\otimes(n+2)}\big) \ar[r] &  H^1_{\mathrm{dR}}(S^{an},\mathrm{Sym}^n_{\mathcal O_{S^{an}}}H^1_{\mathrm{dR}}(E^{an}/S^{an}))}
\end{xy}
\end{equation*}
The key theorem on the analytic side will exhibit the class $\Big((t_{a,b}^{an})^*\big((\varDn)^{an}\big)\Big)^{(n)}$ as the image under the lower horizontal arrow of the collection of (classical) modular forms
\[\bigg(\frac{(-1)^{n+1}\cdot (2\pi i)^{n+2}}{n!}\cdot {_D}F^{(n+2)}_{\frac{aj}{N},\frac{b}{N}}(\tau)\bigg)_{j \in (\Z/N\Z)^*}\]
already introduced in Def. 3.3.17. The crucial input for its proof consists in Thm. 3.8.3, where we gave an explicit description of $(\varDn)^{an}$ via the section $p_n^D(z,\tau)$, and in Thm. 3.6.5, where we determined the specialization of this section along $t_{a,b}^{an}$.\\
As soon as the compatibility of the occurring algebraic and analytic specializations resp. modular forms under the vertical arrows of the diagram is ascertained\footnote{The hereby appearing factor $(2\pi i)^{n+2}$ is due to the fact that we will trivialize $\omega_{E^{an}/S^{an}}$ componentwise on the universal covering $\H$ by the basic section $ \{\omega=\mathrm{d}z \}$ and not by $\{2\pi i \mathrm{d}z \}$, though the last is the common trivialization when dealing with modular forms; the reason why we here don't follow this convention is that we want to keep consistent with the trivialization for $H^1_{\mathrm{dR}}(E^{an}/S^{an})$ on $\H$ by $\{\eta, \omega\}= \{p(z,\tau)\mathrm{d}z, \mathrm{d}z \}$ used throughout the previous sections.\\
Furthermore, the additional minus sign occurring in the case $n>0$ will be explained in Rem. 3.8.13.}, the desired algebraic specialization result then follows from the (already noted) injectivity of the right vertical arrow.\\
\newline
We finally apply the specialization theorem to derive for each $n\geq 0$ an explicit description for the $n$-th de Rham Eisenstein class evaluated at $t_{a,b}$: this is an element of $H^1_{\mathrm{dR}}(S/\Q,\mathrm{Sym}^n_{\mathcal O_S}\mathcal H)$ obtained essentially by specializing the polylogarithm class $\mathrm{pol}_{\mathrm{dR}}^{n+1}$ along $t_{a,b}$, and our result expresses it via the algebraic modular form $F^{(n+2)}_{\frac{a}{N},\frac{b}{N}}$ defined in \cite{Ka}, Ch. I, 3.6.\\
\newline
Let us now come to the detailed execution of the outlined strategy.

\subsubsection{Analytic Hodge filtration, Kodaira-Spencer map and Poincaré duality}
We stipulate that all of the following constructions in (i)-(iv) shall be performed equally on each connected component of $E^{an}/S^{an}$, so we don't need to add everywhere the word "componentwise".\\
\newline
(i) At first, in the same way as we trivialized in $(3.4.6)$ and $(3.4.7)$ the pullback of $H^1_{\mathrm{dR}}(E^{an}/S^{an})$ to $\H$ by the basic sections $\{\eta,\omega\}=\{p(z,\tau)\mathrm{d}z,\mathrm{d}z\}$ and received the automorphy matrix
\[\Gamma(N) \times \H \rightarrow \mathrm{GL}_2(\C), \qquad \Bigg(\begin{pmatrix} a & b \\ c & d \end{pmatrix}, \tau \Bigg) \mapsto \begin{pmatrix} \frac{1}{c\tau+d} & 0 \\ 0 & c\tau +d \end{pmatrix},\]
we do now for the line bundle (the co-Lie algebra of $E^{an}/S^{an}$)
\[\omega_{E^{an}/S^{an}}:=(\pi^{an})_*\Omega^1_{E^{an}/S^{an}} \simeq (\epsilon^{an})^*\Omega^1_{E^{an}/S^{an}}\]
by using the basic section $ \{\omega \} = \{\mathrm{d}z\}$ and obtain the factor of automorphy
\[\Gamma(N) \times \H \rightarrow \C^*, \qquad \Bigg(\begin{pmatrix} a & b \\ c & d \end{pmatrix}, \tau \Bigg) \mapsto c\tau +d.\]
We have a well-defined injection of $\mathcal O_{S^{an}}$-modules by setting
\[\tag{\textbf{3.8.21}} \omega_{E^{an}/S^{an}} \rightarrow H^1_{\mathrm{dR}}(E^{an}/S^{an}), \quad g(\tau) \mapsto \begin{pmatrix}0 \\ g(\tau) \end{pmatrix}.\]
Using for $\omega_{E^{an}/S^{an}}^{\otimes -1}$ the dual basic section $\{\omega^{\vee} \}$ and hence the factor of automorphy $\frac{1}{c\tau+d}$ (cf. 3.2. (iii)), we see that $(3.8.21)$ actually sits in a split short exact sequence ("analytic Hodge filtration")
\[\tag{\textbf{3.8.22}} 0 \rightarrow \omega_{E^{an}/S^{an}} \rightarrow H^1_{\mathrm{dR}}(E^{an}/S^{an}) \rightarrow \omega_{E^{an}/S^{an}}^{\otimes -1} \rightarrow 0\]
with projection resp. section given by
\[\begin{pmatrix} f(\tau) \\ g(\tau) \end{pmatrix} \mapsto f(\tau) \quad \textrm{resp.} \quad \begin{pmatrix} f(\tau) \\ 0 \end{pmatrix} \mapsfrom f(\tau).\]
\newline
(ii) Recall that according to our conventions (cf. 3.2 (iv)) the pullback of $\omega_{E^{an}/S^{an}}^{\otimes k}, \ k \in \Z$, to $\H$ is trivialized by the basic section $\{\omega^{\otimes k} \}=\{(\mathrm{d}z)^{\otimes k}\}$ and receives the factor of automorphy $(c\tau+d)^k$.\\
On the other hand, we had determined the pullback of $\Omega^1_{S^{an}}$ to $\H$ to be trivialized by the basic section $\{\mathrm{d}\tau\}$ and obtained the factor $(c\tau+d)^2$ (cf. $(3.5.27)$).\\
We hence get a well-defined isomorphism of $\mathcal O_{S^{an}}$-modules ("analytic Kodaira-Spencer map") by
\[\tag{\textbf{3.8.23}} \omega_{E^{an}/S^{an}}^{\otimes 2} \xrightarrow{\sim} \Omega^1_{S^{an}}, \quad g(\tau) \mapsto \frac{1}{2\pi i} \cdot g(\tau).\]
Let us remark that $(3.8.23)$ can equally be obtained as the composition
\[\omega_{E^{an}/S^{an}}^{\otimes 2} \rightarrow   \omega_{E^{an}/S^{an}}  \otimes_{\mathcal O_{S^{an}}}  H^1_{\mathrm{dR}}(E^{an}/S^{an}) \rightarrow \omega_{E^{an}/S^{an}}  \otimes_{\mathcal O_{S^{an}}} H^1_{\mathrm{dR}}(E^{an}/S^{an}) \otimes_{\mathcal O_{S^{an}}}  \Omega^1_{S^{an}} \rightarrow\]
\[\rightarrow \omega_{E^{an}/S^{an}} \otimes_{\mathcal O_{S^{an}}} \omega_{E^{an}/S^{an}}^{\otimes -1} \otimes_{\mathcal O_{S^{an}}} \Omega^1_{S^{an}}\rightarrow \Omega^1_{S^{an}}, \qquad \qquad \qquad \qquad \qquad\]
where the first and the third arrow come from $(3.8.22)$, the second from the Gauß-Manin connection on $H^1_{\mathrm{dR}}(E^{an}/S^{an})$ and the last is the canonical isomorphism; in view of the explicit formula for the appearing Gauß-Manin connection (cf. 3.5.4) the verification of this fact is trivial.\\
\newline
For each $n\geq 0$ we can now define a map
\[\tag{\textbf{3.8.24}} \omega_{E^{an}/S^{an}}^{\otimes(n+2)} \rightarrow \Omega^1_{S^{an}} \otimes_{\mathcal O_{S^{an}}} \mathrm{Sym}^n_{\mathcal O_{S^{an}}}H^1_{\mathrm{dR}}(E^{an}/S^{an})\]
as being the composition
\[\omega_{E^{an}/S^{an}}^{\otimes(n+2)} \xrightarrow{\sim} \omega_{E^{an}/S^{an}}^{\otimes 2} \otimes_{\mathcal O_{S^{an}}} \mathrm{Sym}^n_{\mathcal O_{S^{an}}} \omega_{E^{an}/S^{an}} \rightarrow \Omega^1_{S^{an}} \otimes_{\mathcal O_{S^{an}}} \mathrm{Sym}^n_{\mathcal O_{S^{an}}}H^1_{\mathrm{dR}}(E^{an}/S^{an}),\]
where the first arrow is the canonical isomorphism and the second is induced by $(3.8.23)$ and by taking the $n$-th symmetric power of $(3.8.21)$. It is easy to check that $(3.8.24)$ is actually injective.\\
\newline
(iii) The edge morphism $E^{1,0}_2\rightarrow E^1$ at the second sheet of the spectral sequence of hypercohomology
\[E^{p,q}_1=H^q(S^{an}, \Omega^p_{S^{an}}\otimes_{\mathcal O_{S^{an}}}\mathrm{Sym}_{\mathcal O_{S^{an}}}^nH^1_{\mathrm{dR}}(E^{an}/S^{an}))\Rightarrow E^{p+q}=H^{p+q}_{\mathrm{dR}}(S^{an}, \mathrm{Sym}_{\mathcal O_{S^{an}}}^nH^1_{\mathrm{dR}}(E^{an}/S^{an}))\]
furnishes a morphism

\[\tag{\textbf{3.8.25}} \Gamma(S^{an}, \Omega^1_{S^{an}}\otimes_{\mathcal O_{S^{an}}} \mathrm{Sym}^n_{\mathcal O_{S^{an}}}H^1_{\mathrm{dR}}(E^{an}/S^{an})) \rightarrow H^1_{\mathrm{dR}}(S^{an},\mathrm{Sym}^n_{\mathcal O_{S^{an}}}H^1_{\mathrm{dR}}(E^{an}/S^{an})).\]
Using that $S^{an}$ is Stein one can see that $(3.8.25)$ is surjective with kernel given by the image of the connection of $\mathrm{Sym}^n_{\mathcal O_{S^{an}}}H^1_{\mathrm{dR}}(E^{an}/S^{an})$ in global sections.\\
If we precompose $(3.8.25)$ with $(3.8.24)$ in global sections we get the $\C$-linear homomorphism
\[\tag{\textbf{3.8.26}} \Gamma \big(S^{an}, \omega_{E^{an}/S^{an}}^{\otimes(n+2)} \big) \rightarrow H^1_{\mathrm{dR}}(S^{an},\mathrm{Sym}^n_{\mathcal O_{S^{an}}}H^1_{\mathrm{dR}}(E^{an}/S^{an})).\]
(iv) Finally, recall from $(3.5.1)$ and $(3.5.2)$ that we trivialize the pullback of $\mathcal H^{an} \simeq H^1_{\mathrm{dR}}(E^{an}/S^{an})^{\vee}$ to $\H$ by the basic sections $\{\eta^{\vee}, \omega^{\vee} \}$, giving rise to the automorphy matrix  
\[\Gamma(N) \times \H \rightarrow \mathrm{GL}_2(\C), \qquad \Bigg(\begin{pmatrix} a & b \\ c & d \end{pmatrix}, \tau \Bigg) \mapsto \begin{pmatrix}  c\tau +d  & 0 \\ 0 & \frac{1}{c\tau+d}\end{pmatrix}.\]
By exchanging its diagonal entries we obtain just the automorphy matrix for $H^1_{\mathrm{dR}}(E^{an}/S^{an})$, hence
\[\tag{\textbf{3.8.27}} H^1_{\mathrm{dR}}(E^{an}/S^{an}) \xrightarrow{\sim} H^1_{\mathrm{dR}}(E^{an}/S^{an})^{\vee} \simeq \mathcal H^{an}, \quad \begin{pmatrix} f(\tau) \\ g(\tau) \end{pmatrix} \mapsto \begin{pmatrix} g(\tau) \\ -f(\tau) \end{pmatrix} \]
defines an isomorphism ("analytic Poincaré duality") which is easily checked to be horizontal; we will also need the induced identification in symmetric powers
\[\tag{\textbf{3.8.28}} \mathrm{Sym}^n_{\mathcal O_{S^{an}}}H^1_{\mathrm{dR}}(E^{an}/S^{an}) \simeq \mathrm{Sym}^n_{\mathcal O_{S^{an}}}\mathcal H^{an}, \quad n \geq 0.\]
With the identification $(3.8.27)$ we obtain yet another viewpoint on the Kodaira-Spencer map $(3.8.23)$: namely, it is also equal to the composition
\begin{align*}
\omega_{E^{an}/S^{an}}^{\otimes 2} & \rightarrow H^1_{\mathrm{dR}}(E^{an}/S^{an})  \otimes_{\mathcal O_{S^{an}}} H^1_{\mathrm{dR}}(E^{an}/S^{an}) \\
&\rightarrow H^1_{\mathrm{dR}}(E^{an}/S^{an})^\vee  \otimes_{\mathcal O_{S^{an}}} H^1_{\mathrm{dR}}(E^{an}/S^{an}) \otimes_{\mathcal O_{S^{an}}} \Omega^1_{S^{an}} \rightarrow \Omega^1_{S^{an}},
\end{align*}
where the first map comes from the inclusion $(3.8.21)$, the second from $(3.8.27)$ together with the Gauß-Manin connection and the last is given by evaluation. The verification of the claimed equality is again a simple computation, observing the formula for the Gauß-Manin connection in 3.5.4.

\subsubsection{Definition of the analytic specialization}
Recall from $(3.4.1)$-$(3.4.4)$ that the analytic picture for the $N$-torsion section $t_{a,b}$ is given by
\begin{equation*}
\begin{xy}
\xymatrix{
E^{an}=(\Z/N\Z)^*\times (\Z^2 \times \Gamma(N)) \backslash(\C\times \H) \ar[d]_{\pi^{an}} \\
S^{an}=(\Z/N\Z)^*\times \Gamma(N) \backslash \H \ar@/_ 0.3cm/[u]_{t_{a,b}^{an}}  }
\end{xy}
\end{equation*}
where
\[t_{a,b}^{an}: S^{an} \rightarrow U_D^{an} \subseteq E^{an}, \quad (j,\tau) \mapsto \Big(j,\frac{aj\tau}{N}+\frac{b}{N},\tau \Big).\]
We now consider the "specialization" along $t_{a,b}^{an}$ of the analytified $D$-variant (cf. Def. 3.7.4)
\[(\varD)^{an}=\Big(\varDn\Big)^{an}_{n\geq0} \in \lim_{n\geq 0} H^1_{\mathrm{dR}}(U_D^{an}, \mathcal L_n^{an}).\]
This means that we introduce for each $n\geq 0$ the analytic de Rham cohomology class
\[\Big((t_{a,b}^{an})^*\big((\varDn)^{an}\big)\Big)^{(n)} \in H^1_{\mathrm{dR}}(S^{an},\mathrm{Sym}^n_{\mathcal O_{S^{an}}} \mathcal H^{an}),\]
given by pulling back $(\varDn)^{an}$ via $t_{a,b}^{an}$, using the chain of analytified isomorphisms
\[\tag{\textbf{3.8.29}} (t_{a,b}^{an})^*\mathcal L_n^{an} \simeq (\epsilon^{an})^*\mathcal L_n^{an} \simeq \prod_{k=0}^n \mathrm{Sym}^k_{\mathcal O_{S^{an}}}\mathcal H^{an}\]
and by finally taking the $n$-th component.\\
By means of $(3.8.28)$ we will consider the preceding class (without change of notation) as an element
\[\tag{\textbf{3.8.30}} \Big((t_{a,b}^{an})^*\big((\varDn)^{an}\big)\Big)^{(n)} \in H^1_{\mathrm{dR}}(S^{an},\mathrm{Sym}^n_{\mathcal O_{S^{an}}} H^1_{\mathrm{dR}}(E^{an}/S^{an})).\]Before relating $(3.8.30)$ to modular forms we want to realize that it is nothing else than the analytification of the earlier $\big(t_{a,b}^*(\varDn)\big)^{(n)} \in H^1_{\mathrm{dR}}(S/\Q,\mathrm{Sym}^n_{\mathcal O_S} H^1_{\mathrm{dR}}(E/S))$.

\subsubsection{Insertion: compatibility of the algebraic and analytic specialization}
To formulate this in clean terms we need to consider\footnote{using the notation $(.)^{\C}$ as explained at the beginning of 3.7.2} the composition of canonical arrows
\[\tag{\textbf{3.8.31}} \begin{split} H^1_{\mathrm{dR}}(S/\Q,\mathrm{Sym}^n_{\mathcal O_S} H^1_{\mathrm{dR}}(E/S)) \rightarrow H^1_{\mathrm{dR}}(S^{\C}/\C,\mathrm{Sym}^n_{\mathcal O_{S^{\C}}} H^1_{\mathrm{dR}}(E/S)^{\C}) \\
\rightarrow H^1_{\mathrm{dR}}(S^{an},\mathrm{Sym}^n_{\mathcal O_{S^{an}}} H^1_{\mathrm{dR}}(E/S)^{an}) \rightarrow H^1_{\mathrm{dR}}(S^{an},\mathrm{Sym}^n_{\mathcal O_{S^{an}}} H^1_{\mathrm{dR}}(E^{an}/S^{an})). \end{split}\]
Here, the first map comes from the natural arrows
\begin{align*}
H^1_{\mathrm{dR}}(S/\Q,\mathrm{Sym}^n_{\mathcal O_S} H^1_{\mathrm{dR}}(E/S)) &\rightarrow H^1_{\mathrm{dR}}(S/\Q,\mathrm{Sym}^n_{\mathcal O_S} H^1_{\mathrm{dR}}(E/S)) \otimes_{\Q} \ \C \\
&\xrightarrow{\sim} H^1_{\mathrm{dR}}(S^{\C}/\C,\mathrm{Sym}^n_{\mathcal O_{S^{\C}}} H^1_{\mathrm{dR}}(E/S)^{\C})
\end{align*}
and in particular is injective (the second arrow is an isomorphism because $\Spec(\C)/ \Spec(\Q)$ is flat).\\
The second map is the standard one, but it is important to note that it is in fact an isomorphism: namely, the integrable $\C$-connection on $H^1_{\mathrm{dR}}(E/S)^{\C}$ - and hence also on $\mathrm{Sym}^n_{\mathcal O_{S^{\C}}} H^1_{\mathrm{dR}}(E/S)^{\C}$ - is regular: note that under the canonical identification (due to the flatness of $S^{\C}/S$)
\[H^1_{\mathrm{dR}}(E/S)^{\C} \xrightarrow{\sim} H^1_{\mathrm{dR}}(E^{\C}/S^{\C})\]
it identifies with the Gauß-Manin connection of $H^1_{\mathrm{dR}}(E^{\C}/S^{\C})$ relative $\Spec(\C)$ which is regular by \cite{De1}, Ch. II, Thm. 7.9 and Prop. 6.14; now use ibid., Ch. II, Thm. 6.2 for the claimed isomorphy.\\
The final map comes from the canonical (horizontal) isomorphism (cf. the beginning of 3.5)
\[H^1_{\mathrm{dR}}(E/S)^{an} \xrightarrow{\sim} H^1_{\mathrm{dR}}(E^{an}/S^{an})\]
and hence is an isomorphism.\\
In sum, we have explained how to understand the chain $(3.8.31)$ and seen that it is injective.
\begin{lemma}
For each $n\geq 0$ the image of $\big(t_{a,b}^*(\varDn)\big)^{(n)}$ under the injection $(3.8.31)$ coincides with $\Big((t_{a,b}^{an})^*\big((\varDn)^{an}\big)\Big)^{(n)}$ in $H^1_{\mathrm{dR}}(S^{an},\mathrm{Sym}^n_{\mathcal O_{S^{an}}} H^1_{\mathrm{dR}}(E^{an}/S^{an}))$.
\end{lemma}
\begin{proof}
This is, except for one point, a routine check enfolding the definitions and observing a couple of natural compatibilities: one writes down a large diagram of the various occurring de Rham cohomology spaces, furnishes commutativity of each square by an evident compatibility argument and deduces the claim from the resulting commutativity of the whole. As this is very straightforward we forbear from explicating it here, but mention that there is one piece of information, different from the routine compatibilities, which one thereby needs to make a clean conclusion:\\
Namely, that the analytification $H^1_{\mathrm{dR}}(E/S)^{an}\simeq \mathcal H^{an}$ of $(3.8.20)$ coincides with the isomorphism $H^1_{\mathrm{dR}}(E^{an}/S^{an})\simeq \mathcal H^{an}$ of $(3.8.27)$ when making the identification $H^1_{\mathrm{dR}}(E/S)^{an} \simeq H^1_{\mathrm{dR}}(E^{an}/S^{an})$ (for the last cf. the beginning of 3.5). As the dual of the preceding isomorphism was used in the definition of $(3.8.27)$ we see that we need to verify the following statement:\\
Under the natural identification $H^1_{\mathrm{dR}}(E/S)^{an} \simeq H^1_{\mathrm{dR}}(E^{an}/S^{an})$ and its dual $\mathcal H^{an} \simeq H^1_{\mathrm{dR}}(E^{an}/S^{an})^{\vee}$ the analytification
\[H^1_{\mathrm{dR}}(E/S)^{an}\xrightarrow{\sim} \mathcal H^{an}\]
of $(3.8.20)$ coincides with the isomorphism
\[H^1_{\mathrm{dR}}(E^{an}/S^{an}) \xrightarrow{\sim} H^1_{\mathrm{dR}}(E^{an}/S^{an})^{\vee}, \quad \begin{pmatrix} f(\tau) \\ g(\tau) \end{pmatrix} \mapsto \begin{pmatrix} g(\tau) \\ -f(\tau) \end{pmatrix},\]
where we recall that the (componentwise) pullbacks to $\H$ of the preceding two vector bundles are trivialized by the basic sections $\{\eta,\omega\}=\{p(z,\tau)\mathrm{d}z,\mathrm{d}z\}$ resp. by $\{\eta^{\vee}, \omega^{\vee}\}$.\\
This in turn can be verified fiberwise over points of $S^{an}$ and then boils down to the fact that $(3.8.20)$ for a single complex elliptic curve identifies the basis elements $\eta$ resp. $\omega$ with $-\omega^{\vee}$ resp. $\eta^{\vee}$, i.e. to the standard formula $\mathrm{tr}(\eta\cup \omega)=-1=-\mathrm{tr}(\omega\cup\eta)$ (for more details cf. \cite{Kat4}, Rem. A 1.3.13).
\end{proof}

\subsubsection{The analytic specialization result}
We will now describe for each $n\geq 0$ the cohomology class
\[\Big((t_{a,b}^{an})^*\big((\varDn)^{an}\big)\Big)^{(n)} \in H^1_{\mathrm{dR}}(S^{an},\mathrm{Sym}^n_{\mathcal O_{S^{an}}} H^1_{\mathrm{dR}}(E^{an}/S^{an}))\]
in terms of the modular forms $_DF^{(n+2)}_{\alpha,\beta}(\tau)$ introduced in Def. 3.3.17.\\
\newline
At first, according to our fixed trivialization by $\{\omega^{\otimes(n+2)} \}=\{ (\mathrm{d}z)^{\otimes (n+2)}\}$ for the pullback of $\omega_{E^{an}/S^{an}}^{\otimes(n+2)}$ to the universal covering $\H$ of each component of $S^{an}$, a section in $\Gamma \big(S^{an},\omega_{E^{an}/S^{an}}^{\otimes(n+2)} \big)$ is tantamount to a collection
\[\big(f_j(\tau)\big)_{j \in (\Z/N\Z)^*}\]
of holomorphic functions $f_j(\tau)$ on $\H$, each satisfying
\[f_j\bigg(\frac{a\tau+b}{c\tau+d}\bigg)=(c\tau+d)^{n+2}\cdot f_j(\tau) \quad \textrm{for all} \ \ \begin{pmatrix} a & b \\ c & d \end{pmatrix} \in \Gamma(N).\]
In the following, we will usually apply the same notation $j$ for a class in $(\Z/N\Z)^*$ or for a (arbitrarily chosen) representative for this class.\\
\newline
Recalling from 3.3.4 that for each integer $j$ with $(j,N)=1$ the function
\[\tau \mapsto {_D}F^{(n+2)}_{\frac{aj}{N},\frac{b}{N}}(\tau)\]
is a modular form of weight $n+2$ and level $N$, dependent only from the class of $j$ in $(\Z/N\Z)^*$, we see that the collection
\[\tag{\textbf{3.8.32}} \bigg(\frac{(-1)^{n+1}\cdot (2\pi i)^{n+2}}{n!}\cdot {_D}F^{(n+2)}_{\frac{aj}{N},\frac{b}{N}}(\tau)\bigg)_{j \in (\Z/N\Z)^*}\]
gives a well-defined element of $\Gamma \big(S^{an},\omega_{E^{an}/S^{an}}^{\otimes(n+2)} \big)$.\\
\newline
We can then formulate the main result about the specialization of the analytified $D$-variant of the polylogarithm.
\begin{theorem}
For each $n \geq 0$ the cohomology class $\Big((t_{a,b}^{an})^*\big((\varDn)^{an}\big)\Big)^{(n)}$ is equal to the image of $(3.8.32)$ under the map
\[\Gamma \big(S^{an}, \omega_{E^{an}/S^{an}}^{\otimes(n+2)} \big) \rightarrow H^1_{\mathrm{dR}}(S^{an},\mathrm{Sym}^n_{\mathcal O_{S^{an}}}H^1_{\mathrm{dR}}(E^{an}/S^{an}))\]
of $(3.8.26)$.
\end{theorem}

\begin{proof}
It is the assertion of Thm. 3.8.3 that $(\varDn)^{an} \in H^1_{\mathrm{dR}}(U_D^{an},\mathcal L_n^{an})$ is the image of $p_n^D(z,\tau)$ under the canonical arrow
\[\ker\bigg(\Gamma(U_D^{an}, \Omega^1_{E^{an}}\otimes_{\mathcal O_{E^{an}}} \mathcal L_n^{an})\xrightarrow{(\nabla_n^{an})^1}\Gamma(U_D^{an}, \Omega^2_{E^{an}} \otimes_{\mathcal O_{E^{an}}} \mathcal L_n^{an})\bigg) \rightarrow H^1_{\mathrm{dR}}(U_D^{an},\mathcal L_n^{an})\]
induced by $(3.7.29)$. Recall that $p_n^D(z,\tau)$ is the vector of functions defined in $(3.6.2)$.\\
Now Thm. 3.6.5 implies that if we pull back $p_n^D(z,\tau) \in \Gamma(U_D^{an}, \Omega^1_{E^{an}}\otimes_{\mathcal O_{E^{an}}} \mathcal L_n^{an})$ along $t_{a,b}^{an}$, use the chain of isomorphisms
\begin{align*}
(t_{a,b}^{an})^*(\Omega^1_{E^{an}}\otimes_{\mathcal O_{E^{an}}} \mathcal L_n^{an}) & \simeq (t_{a,b}^{an})^*\Omega^1_{E^{an}}\otimes_{\mathcal O_{S^{an}}} \prod_{k=0}^n \mathrm{Sym}^k_{\mathcal O_{S^{an}}}H^1_{\mathrm{dR}}(E^{an}/S^{an})^\vee \\
&\xrightarrow{\mathrm{can}} \Omega^1_{S^{an}} \otimes_{\mathcal O_{S^{an}}} \prod_{k=0}^n \mathrm{Sym}^k_{\mathcal O_{S^{an}}}H^1_{\mathrm{dR}}(E^{an}/S^{an})^\vee
\end{align*}
and finally take the $n$-part, then the obtained element in $\Gamma(S^{an},\Omega^1_{S^{an}}\otimes_{\mathcal O_{S^{an}}}  \mathrm{Sym}^n_{\mathcal O_{S^{an}}} H^1_{\mathrm{dR}}(E^{an}/S^{an})^\vee)$ is given by
\[\bigg(\frac{(-1)^{n+1}\cdot (2\pi i)^{n+1}}{n!}\cdot \ _DF^{(n+2)}_{\frac{aj}{N},\frac{b}{N}}(\tau) \cdot \mathrm{d}\tau \otimes \frac{(\eta^{\vee})^n}{(n-n)!}\bigg)_{j \in (\Z/N\Z)^*}.\]
Under our fixed identification
\[\mathrm{Sym}^n_{\mathcal O_{S^{an}}}H^1_{\mathrm{dR}}(E^{an}/S^{an})^\vee \simeq \mathrm{Sym}^n_{\mathcal O_{S^{an}}}H^1_{\mathrm{dR}}(E^{an}/S^{an})\]
of $(3.8.28)$ this becomes the element of $\Gamma(S^{an},\Omega^1_{S^{an}}\otimes_{\mathcal O_{S^{an}}}  \mathrm{Sym}^n_{\mathcal O_{S^{an}}} H^1_{\mathrm{dR}}(E^{an}/S^{an}))$ given by
\[\tag{$*$} \bigg(\frac{(-1)^{n+1}\cdot (2\pi i)^{n+1}}{n!}\cdot \ _DF^{(n+2)}_{\frac{aj}{N},\frac{b}{N}}(\tau) \cdot \mathrm{d}\tau \otimes \frac{\omega^n}{(n-n)!}\bigg)_{j \in (\Z/N\Z)^*}.\]
Putting everything together we obtain that
\[\Big((t_{a,b}^{an})^*\big((\varDn)^{an}\big)\Big)^{(n)}\]
is equal to the image of $(*)$ under the arrow
\[\Gamma(S^{an}, \Omega^1_{S^{an}}\otimes_{\mathcal O_{S^{an}}} \mathrm{Sym}^n_{\mathcal O_{S^{an}}}H^1_{\mathrm{dR}}(E^{an}/S^{an})) \rightarrow H^1_{\mathrm{dR}}(S^{an},\mathrm{Sym}^n_{\mathcal O_{S^{an}}}H^1_{\mathrm{dR}}(E^{an}/S^{an}))\]
of $(3.8.25)$. Now, the map of $(3.8.26)$
\[\Gamma \big(S^{an}, \omega_{E^{an}/S^{an}}^{\otimes(n+2)} \big) \rightarrow H^1_{\mathrm{dR}}(S^{an},\mathrm{Sym}^n_{\mathcal O_{S^{an}}}H^1_{\mathrm{dR}}(E^{an}/S^{an}))\]
in the claim of the proposition is the composition of the preceding arrow with
\[\tag{$**$} \Gamma \big(S^{an},\omega_{E^{an}/S^{an}}^{\otimes(n+2)}\big) \rightarrow \Gamma(S^{an},\Omega^1_{S^{an}} \otimes_{\mathcal O_{S^{an}}} \mathrm{Sym}^n_{\mathcal O_{S^{an}}}H^1_{\mathrm{dR}}(E^{an}/S^{an})),\]
where the last is $(3.8.24)$ in global $S^{an}$-sections. It hence remains to show that $(**)$ sends the section
\[\bigg(\frac{(-1)^{n+1}\cdot (2\pi i)^{n+2}}{n!}\cdot {_D}F^{(n+2)}_{\frac{aj}{N},\frac{b}{N}}(\tau) \cdot \omega^{\otimes (n+2)}\bigg)_{j \in (\Z/N\Z)^*} \in \Gamma \big(S^{an},\omega_{E^{an}/S^{an}}^{\otimes(n+2)}\big),\]
defined by $(3.8.32)$, to $(*)$. This is clear by the very definition of $(3.8.24)$ and the remark that the Kodaira-Spencer map appearing in it sends $\omega^{\otimes 2}$ to $\frac{1}{2\pi i}\cdot \mathrm{d}\tau$ on the universal covering (cf. 3.8.23).
\end{proof}

We can now approach the task of expressing the algebraic de Rham cohomology class
\[\big(t_{a,b}^*(\varDn)\big)^{(n)} \in H^1_{\mathrm{dR}}(S/\Q,\mathrm{Sym}^n_{\mathcal O_S} H^1_{\mathrm{dR}}(E/S))\]
as induced by an algebraic modular form of weight $n+2$ and level $N$.

\subsubsection{Hodge filtration and Kodaira-Spencer map}
(i) Writing $\omega_{E/S}$ for the line bundle (the co-Lie algebra of $E/S$)
\[\omega_{E/S}:=\pi_*\Omega^1_{E/S} \simeq \epsilon^*\Omega^1_{E/S},\]
the degeneration of the Hodge-de Rham spectral sequence at the first sheet (cf. the beginning of Chapter 1) furnishes a short exact sequence of $\mathcal O_S$-modules (Hodge filtration)
\[\tag{\textbf{3.8.33}} 0 \rightarrow \omega_{E/S} \rightarrow H^1_{\mathrm{dR}}(E/S) \rightarrow R^1\pi_*\mathcal O_E \rightarrow 0.\]
Identifying $R^1\pi_*\mathcal O_E \simeq \omega_{E/S}^{\otimes{-1}}$ by Serre duality\footnote{One considers the $\mathcal O_S$-linear map
\[\pi_*\Omega^1_{E/S}\otimes_{\mathcal O_S} R^1\pi_*\mathcal O_E \rightarrow R^1\pi_*\Omega^1_{E/S} \xrightarrow{\sim} \mathcal O_S,\]
where the first arrow comes from the natural bilinear map / cup product (cf. \cite{Li}, Ch. 6, $(4.11)$) and the second is the Grothendieck trace isomorphism (cf. \cite{Con1}, Ch. I, 1.1). It is a perfect pairing and the induced isomorphism (Serre duality)
\[R^1\pi_*\mathcal O_E \xrightarrow{\sim} \omega_{E/S}^{\otimes {-1}}\]
is the desired one.}, this reads as an exact sequence
\[\tag{\textbf{3.8.34}} 0 \rightarrow \omega_{E/S} \rightarrow H^1_{\mathrm{dR}}(E/S) \rightarrow \omega_{E/S}^{\otimes{-1}} \rightarrow 0\]
for which one can give an explicit local description after the choice of a local trivialization for $\omega_{E/S}$, showing in particular that $(3.8.34)$ splits in a canonical way (cf. \cite{Kat4}, A 1.2.1-A 1.2.6).\\
We remark that the arrows in $(3.8.33)$ are nothing else than those induced by application of $\R^1\pi_*$ to the maps in the canonical short exact sequence of complexes
\[0 \rightarrow \Omega^1_{E/S}[-1] \rightarrow \Omega_{E/S}^{\bullet} \rightarrow \mathcal O_E \rightarrow 0.\]
(ii) Next, let us briefly recall the Kodaira-Spencer map for our situation $E/S/\Q$ (cf. e.g. \cite{Kat-Maz}, $(10.13.10)$ or \cite{Kat4}, A 1.3.17). This is the $\mathcal O_S$-linear morphism 
\[\tag{\textbf{3.8.35}} \omega_{E/S}^{\otimes 2} \rightarrow \Omega^1_{S/\Q}\]
arising from the composition
\[\omega_{E/S}^{\otimes 2} \rightarrow H^1_{\mathrm{dR}}(E/S)  \otimes_{\mathcal O_S} H^1_{\mathrm{dR}}(E/S) \rightarrow H^1_{\mathrm{dR}}(E/S)^\vee  \otimes_{\mathcal O_S} H^1_{\mathrm{dR}}(E/S) \otimes_{\mathcal O_S} \Omega^1_{S/\Q} \rightarrow \Omega^1_{S/\Q},\]
where the first arrow comes from the injection in $(3.8.34)$, the second from Poincaré duality in $(3.8.20)$ together with the Gauß-Manin connection and the last is given by evaluation.\\
It is part of the theory of moduli schemes that $(3.8.35)$ is an isomorphism (cf. \cite{Kat-Maz}, $(10.13.10)$).
\begin{remark}
It is possible to give various equivalent definitions for $(3.8.35)$, e.g. as the composite map
\[\omega_{E/S}^{\otimes 2} \rightarrow   \omega_{E/S}  \otimes_{\mathcal O_S}  H^1_{\mathrm{dR}}(E/S) \rightarrow \omega_{E/S}  \otimes_{\mathcal O_S} H^1_{\mathrm{dR}}(E/S) \otimes_{\mathcal O_S}  \Omega^1_{S/\Q} \rightarrow \]
\[\rightarrow \omega_{E/S} \otimes_{\mathcal O_S} \omega_{E/S}^{\otimes -1} \otimes_{\mathcal O_S} \Omega^1_{S/\Q}\rightarrow \Omega^1_{S/\Q}, \qquad \qquad \qquad \qquad \ \ \]where the third arrow is induced by the projection in $(3.8.34)$ and the others should by now be clear.\\
One can also introduce it as the morphism obtained by tensoring with $\omega_{E/S}$ the map
\[\omega_{E/S}\rightarrow \omega_{E/S}^{\otimes{-1}} \otimes_{\mathcal O_S} \Omega^1_{S/\Q}\]
coming from the boundary $\omega_{E/S} =\pi_*\Omega^1_{E/S} \rightarrow R^1\pi_*\mathcal O_E\otimes_{\mathcal O_S} \Omega^1_{S/\Q}$ for the canonical exact sequence
\[0 \rightarrow \pi^*\Omega^1_{S/\Q}\rightarrow \Omega^1_{E/\Q}\rightarrow \Omega^1_{E/S} \rightarrow 0\]
and from the Serre duality isomorphism $R^1\pi_*\mathcal O_E \xrightarrow{\sim} \omega_{E/S}^{\otimes{-1}}$ (cf. footnote 29).\\
To show the equality of all these different approaches to $(3.8.35)$ is standard.\footnote{That the first two of the three approaches coincide follows from the well-known fact that Poincaré duality in $(3.8.20)$ and Serre duality in footnote 29 are compatible via the Hodge filtration of $(3.8.33)$ in the sense that the obtained diagram
\begin{equation*}
\begin{xy}
\xymatrix{
\omega_{E/S}  \ar[d] \ar[r]^{\sim \qquad} & (R^1\pi_*\mathcal O_E)^{\vee}\ar[d]\\
H^1_{\mathrm{dR}}(E/S)  \ar[r]^{\sim} & H^1_{\mathrm{dR}}(E/S)^\vee}
\end{xy}
\end{equation*}
commutes; for more details about this cf. the explanations in \cite{Ca}, 2.2, which hold invariantly over a base ring.\\
The equality of the last two definitions for the Kodaira-Spencer map requires to check that the composition
\[\omega_{E/S} \rightarrow H^1_{\mathrm{dR}}(E/S) \otimes_{\mathcal O_S} \Omega^1_{S/\Q} \rightarrow R^1\pi_*\mathcal O_E\otimes_{\mathcal O_S}\Omega^1_{S/\Q},\]
induced by the Hodge filtration and Gauß-Manin connection, is the same as the boundary map mentioned in Rem. 3.8.9. This in turn follows directly from the definition of the Gauß-Manin connection on $H^1_{\mathrm{dR}}(E/S)$.
}
\end{remark}
(iii) As in $(3.8.25)$ we obtain from the spectral sequence of hypercohomology a map
\[\tag{\textbf{3.8.36}} \Gamma(S, \Omega^1_{S/\Q}\otimes_{\mathcal O_S} \mathrm{Sym}^n_{\mathcal O_S}H^1_{\mathrm{dR}}(E/S)) \rightarrow H^1_{\mathrm{dR}}(S/\Q,\mathrm{Sym}^n_{\mathcal O_S} H^1_{\mathrm{dR}}(E/S))\]
which (as $S$ is affine) is surjective with kernel given by the image of the connection of $\mathrm{Sym}^n_{\mathcal O_S}H^1_{\mathrm{dR}}(E/S)$ in global sections.\\
Moreover, we define for each $n\geq 0$ a map
\[\tag{\textbf{3.8.37}} \omega_{E/S}^{\otimes(n+2)} \rightarrow \Omega^1_{S/\Q} \otimes_{\mathcal O_S} \mathrm{Sym}^n_{\mathcal O_S}H^1_{\mathrm{dR}}(E/S)\]
as the composition
\[\omega_{E/S}^{\otimes(n+2)} \xrightarrow{\sim} \omega_{E/S}^{\otimes 2} \otimes_{\mathcal O_S} \mathrm{Sym}^n_{\mathcal O_S} \omega_{E/S} \rightarrow \Omega^1_{S/\Q} \otimes_{\mathcal O_S} \mathrm{Sym}^n_{\mathcal O_S}H^1_{\mathrm{dR}}(E/S),\]
where the first arrow is the canonical isomorphism and the second is induced by $(3.8.35)$ and by taking the $n$-th symmetric power of the inclusion in $(3.8.34)$; the map $(3.8.37)$ is actually injective.\\
From $(3.8.37)$ and $(3.8.36)$ we obtain the $\Q$-linear homomorphism
\[\tag{\textbf{3.8.38}} \Gamma \big(S, \omega_{E/S}^{\otimes(n+2)} \big) \rightarrow H^1_{\mathrm{dR}}(S/\Q,\mathrm{Sym}^n_{\mathcal O_S}H^1_{\mathrm{dR}}(E/S)).\]
It is compatible with its analytic counterpart $(3.8.26)$:
\begin{lemma}
For each $n\geq 0$ the diagram
\begin{equation*}
\begin{xy}
\xymatrix@C-0.3cm{
\Gamma \big(S, \omega_{E/S}^{\otimes(n+2)}\big)\ar[r] \ar[d] & H^1_{\mathrm{dR}}(S/\Q,\mathrm{Sym}^n_{\mathcal O_S}H^1_{\mathrm{dR}}(E/S)) \ar[d]\\
\Gamma \big(S^{an}, \omega_{E^{an}/S^{an}}^{\otimes(n+2)}\big) \ar[r] &  H^1_{\mathrm{dR}}(S^{an},\mathrm{Sym}^n_{\mathcal O_{S^{an}}}H^1_{\mathrm{dR}}(E^{an}/S^{an}))}
\end{xy}
\end{equation*}
where the upper resp. lower horizontal arrow is given by $(3.8.38)$ resp. by $(3.8.26)$ and the left resp. right vertical arrow is given by the natural map resp. by $(3.8.31)$, is commutative.
\end{lemma}
\begin{proof}
According to the mere definitions of the maps we may split the diagram in question into the following two diagrams:\\
The first is
\begin{equation*}
\begin{xy}
\xymatrix@C-0.3cm{
\Gamma\big(S, \omega_{E/S}^{\otimes(n+2)}\big) \ar[r] \ar[d]&\Gamma(S, \Omega^1_{S/\Q}\otimes_{\mathcal O_S} \mathrm{Sym}^n_{\mathcal O_S}H^1_{\mathrm{dR}}(E/S))\ar[d] \\
\Gamma\big(S^{an}, \omega_{E^{an}/S^{an}}^{\otimes(n+2)}\big)\ar[r]&\Gamma(S^{an}, \Omega^1_{S^{an}}\otimes_{\mathcal O_{S^{an}}} \mathrm{Sym}^n_{\mathcal O_{S^{an}}}H^1_{\mathrm{dR}}(E^{an}/S^{an}))}
\end{xy}
\end{equation*}
Here, the horizontal arrows come from the algebraic resp. analytic Kodaira-Spencer map together with the algebraic resp. analytic Hodge filtration, and the right vertical arrow is the canonical one - always taking into account the identification $H^1_{\mathrm{dR}}(E/S)^{an}\xrightarrow{\sim} H^1_{\mathrm{dR}}(E^{an}/S^{an})$.\\
The second diagram is
\begin{equation*}
\begin{xy}
\xymatrix@C-0.3cm{
\Gamma(S, \Omega^1_{S/\Q}\otimes_{\mathcal O_S} \mathrm{Sym}^n_{\mathcal O_S}H^1_{\mathrm{dR}}(E/S))\ar[r] \ar[d] & H^1_{\mathrm{dR}}(S/\Q,\mathrm{Sym}^n_{\mathcal O_S}H^1_{\mathrm{dR}}(E/S)) \ar[d]\\
 \Gamma(S^{an}, \Omega^1_{S^{an}}\otimes_{\mathcal O_{S^{an}}} \mathrm{Sym}^n_{\mathcal O_{S^{an}}}H^1_{\mathrm{dR}}(E^{an}/S^{an}))\ar[r] &  H^1_{\mathrm{dR}}(S^{an},\mathrm{Sym}^n_{\mathcal O_{S^{an}}}H^1_{\mathrm{dR}}(E^{an}/S^{an}))}
\end{xy}
\end{equation*}
with horizontal arrows induced by the respective spectral sequence of hypercohomology for the de Rham complex of $\mathrm{Sym}^n_{\mathcal O_S}H^1_{\mathrm{dR}}(E/S)$ resp. $\mathrm{Sym}^n_{\mathcal O_{S^{an}}}H^1_{\mathrm{dR}}(E^{an}/S^{an})$, as in $(3.8.36)$ and $(3.8.25)$.\\
The first of the previous two diagrams commutes: the essential point is that all maps on the analytic side used for the definition of the lower horizontal arrow were constructed exactly in such a way that they are the analytifications of the corresponding algebraic maps used to define the upper horizontal arrow. To be more precise, always taking into account the usual canonical identification $H^1_{\mathrm{dR}}(E/S)^{an}\xrightarrow{\sim} H^1_{\mathrm{dR}}(E^{an}/S^{an})$, the analytification of $(3.8.20)$ resp. $(3.8.34)$ resp. $(3.8.35)$ becomes precisely the analytic Poincaré duality in $(3.8.27)$ resp. the analytic Hodge filtration in $(3.8.22)$ resp. the analytic Kodaira-Spencer map in $(3.8.23)$.\\
Everything of this is easy and can be seen in multiple ways; e.g., one may use the description of the analytic Kodaira-Spencer map via analytic Poincaré duality, explained subsequent to $(3.8.28)$, and as to Poincaré duality one may then recall the argument already given in the proof of Lemma 3.8.7.\\
With these remarks the commutativity of the first diagram becomes obvious.\\
The commutativity of the second diagram is due to the compatibility of the two spectral sequences used in its definition: the natural analytification maps between their initial and between their limit terms are part of a morphism of spectral sequences, as can be seen very easily in our special case but also holds in general, similar as in \cite{De1}, Ch. II, p. 101.
\end{proof}

\subsubsection{Preparations for the algebraic specialization result}

\begin{definition}
As in \cite{Ka}, Ch. I, 4.2, we define for each $k\geq 1$ an algebraic modular form of weight $k$ and level $N$:
\[_DF^{(k)}_{\frac{a}{N},\frac{b}{N}}:=D^2F^{(k)}_{\frac{a}{N},\frac{b}{N}}-D^{2-k}F^{(k)}_{\frac{Da}{N},\frac{Db}{N}}\]
with
\[F^{(k)}_{\frac{a}{N},\frac{b}{N}} \quad \textrm{resp.} \quad F^{(k)}_{\frac{Da}{N},\frac{Db}{N}}\]
the algebraic modular forms of weight $k$ and level $N$ introduced in ibid., Ch. I, 3.6.\\
For the construction and basic properties of the latter as well as for generalities about algebraic modular forms we refer to ibid., Ch. I, 3.
\end{definition}

\begin{remark}
An algebraic modular form of weight $k \in \Z$ and level $N$ defines (by restriction to the open modular curve $S$) a section in $\Gamma\big(S, \omega_{E/S}^{\otimes k}\big)$ and hence, via the canonical arrow
\[\tag{\textbf{3.8.39}} \Gamma \big(S, \omega_{E/S}^{\otimes k}\big) \rightarrow \Gamma\big(S^{an}, \omega_{E^{an}/S^{an}}^{\otimes k}\big),\]
an element of $\Gamma(S^{an}, \omega_{E^{an}/S^{an}}^{\otimes k})$.\\
As explained previous to Thm. 3.8.8, we view any such element
\[\Psi \in \Gamma\big(S^{an}, \omega_{E^{an}/S^{an}}^{\otimes k}\big)\]
as a collection
\[\big(f_j(\tau)\big)_{j \in (\Z/N\Z)^*}\]
of holomorphic functions $f_j(\tau)$ on $\H$ satisfying
\[f_j\bigg(\frac{a\tau+b}{c\tau+d}\bigg)=(c\tau+d)^k\cdot f_j(\tau) \ \ \textrm{for all} \ \ \begin{pmatrix} a & b \\ c & d \end{pmatrix} \in \Gamma(N).\]
Recall our conventions how to obtain from the abstract section $\Psi$ the associated functions $f_j(\tau)$:\\
Namely, for each $j\in (\Z/N\Z)^*$ consider the cartesian diagram 
\begin{equation*} \tag{\textbf{3.8.40}} \begin{split}
\begin{xy}
\xymatrix@C-0.3cm{
\mathbb E:=\Z^2\backslash(\C\times \H) \ar[r]\ar[d] & (\Z/N\Z)^*\times (\Z^2\times \Gamma(N)) \backslash (\C \times \H)=E^{an} \ar[d]\\
\H \ar[r]^{p_j \qquad \qquad \quad} & (\Z/N\Z)^* \times \Gamma(N) \backslash \H=S^{an}}
\end{xy}
\end{split}
\end{equation*}
where the lower arrow is given by the composition
\[\H \rightarrow \Gamma(N)\backslash \H \rightarrow (\Z/N\Z)^* \times \Gamma(N) \backslash \H\]
of the projection and inclusion into the $j$-component and where the upper arrow is the analogous chain
\[\Z^2\backslash(\C\times \H) \rightarrow (\Z^2\times \Gamma(N)) \backslash (\C \times \H) \rightarrow (\Z/N\Z)^*\times (\Z^2\times \Gamma(N)) \backslash (\C \times \H).\]
Then $\Psi$ naturally induces a section
\[p_j^*(\Psi) \in \Gamma\big(\H,\omega_{\mathbb E/\H}^{\otimes k}\big)\]
which - via the fixed trivialization of $\omega_{\mathbb E/\H}$ by the global section $\omega=\mathrm{d}z \in \Gamma(\mathbb E, \Omega^1_{\mathbb E/\H})=\Gamma(\H,\omega_{\mathbb E/\H})$ and the induced trivialization for $\omega_{\mathbb E/\H}^{\otimes k}$ by $\{\omega^{\otimes k} \}= \{ (\mathrm{d}z)^{\otimes k}\}$ - writes as
\[p_j^*(\Psi)=f_j(\tau) \cdot (\mathrm{d}z)^{\otimes k}\]
for a unique holomorphic function $f_j(\tau)$ on $\H$ (necessarily transforming as above).\\
Our deviance from the standard trivialization by $\{ (2\pi i\mathrm{d}z)^{\otimes k}\}$ was already pointed out and explained in footnote 27.
\end{remark}

\begin{remark}
In \cite{Ka}, Ch. I, p. 139, an algebraic modular form is assigned a classical modular form $f(\tau)$ by viewing it as an element of $\Gamma \big(S, \omega_{E/S}^{\otimes k} \big)$, taking its image under $(3.8.39)$ and then applying the process of Rem. 3.8.12 for $j=1$, but with the trivialization of $\omega_{\mathbb E/\H}^{\otimes k}$ by $\{ (2\pi i\mathrm{d}z)^{\otimes k}\}$.\\
In this sense, (i) of ibid., Ch. I, p. 140, states that for $k\geq 1, \ k\neq 2,$ the classical modular form assigned to the algebraic modular form $E^{(k)}_{\frac{a}{N},\frac{b}{N}}$defined in ibid., Ch. I, 3.3, is the $E^{(k)}_{\frac{a}{N},\frac{b}{N}}(\tau)$ we have introduced in 3.3.4; by definition, this would imply that for any such $k$ the function assigned to the algebraic modular form $F^{(k)}_{\frac{a}{N},\frac{b}{N}}$ resp. $_DF^{(k)}_{\frac{a}{N},\frac{b}{N}}$ is the $F^{(k)}_{\frac{a}{N},\frac{b}{N}}(\tau)$ resp. $_DF^{(k)}_{\frac{a}{N},\frac{b}{N}}(\tau)$ of 3.3.4.\\
But it is instead true that (still in Kato's convention and for $k\neq 2$) the function assigned to $E^{(k)}_{\frac{a}{N},\frac{b}{N}}$ is $-E^{(k)}_{\frac{a}{N},\frac{b}{N}}(\tau)$ and thus that $F^{(k)}_{\frac{a}{N},\frac{b}{N}}$ resp. $_DF^{(k)}_{\frac{a}{N},\frac{b}{N}}$ induces $-F^{(k)}_{\frac{a}{N},\frac{b}{N}}(\tau)$ resp. $-_DF^{(k)}_{\frac{a}{N},\frac{b}{N}}(\tau)$.\footnote{This is easily checked, and the reason for the erroneous claim in ibid., Ch. I, p. 140, is essentially a constant sign mistake in the logarithmic derivative of the Kato-Siegel function.}\\
On the other hand, the function assigned to the algebraic modular form $\widetilde{E}^{(2)}_{\frac{a}{N},\frac{b}{N}}$ defined in ibid., Ch. I, 3.4, is indeed the $\widetilde{E}^{(2)}_{\frac{a}{N},\frac{b}{N}}(\tau)$ we have introduced in 3.3.4, such that $F^{(2)}_{\frac{a}{N},\frac{b}{N}}$ resp. $_DF^{(2)}_{\frac{a}{N},\frac{b}{N}}$ induces the $F^{(2)}_{\frac{a}{N},\frac{b}{N}}(\tau)$ resp. $_DF^{(2)}_{\frac{a}{N},\frac{b}{N}}(\tau)$ of 3.3.4.\\
Hence, in \textit{our} conventions and notations of Rem. 3.8.12, if $k\geq 1$ and $\Psi_k$ denotes the image of $_DF^{(k)}_{\frac{a}{N},\frac{b}{N}}$ via $(3.8.39)$, then it follows altogether:
\[\tag{\textbf{3.8.41}}\begin{split}
p_1^*(\Psi_k)&=-(2\pi i)^k\cdot _DF^{(k)}_{\frac{a}{N},\frac{b}{N}}(\tau)\cdot (\mathrm{d}z)^{\otimes k}, \quad k\geq 1, \ k\neq 2,\\
p_1^*(\Psi_2)&=(2\pi i)^2\cdot _DF^{(2)}_{\frac{a}{N},\frac{b}{N}}(\tau)\cdot (\mathrm{d}z)^{\otimes 2}.
\end{split}\]
\end{remark}

\begin{lemma}
In the sense of Rem. 3.8.12, for each $k\geq 1, \ k \neq 2,$ the image of $_DF^{(k)}_{\frac{a}{N},\frac{b}{N}}$ under $(3.8.39)$ is given by
\[\bigg(-(2\pi i)^k \cdot {_D}F^{(k)}_{\frac{aj}{N},\frac{b}{N}}(\tau)\bigg)_{j \in (\Z/N\Z)^*},\]
and the image of $_DF^{(2)}_{\frac{a}{N},\frac{b}{N}}$ under $(3.8.39)$ is given by
\[\bigg((2\pi i)^2 \cdot {_D}F^{(2)}_{\frac{aj}{N},\frac{b}{N}}(\tau)\bigg)_{j \in (\Z/N\Z)^*},\]
where for each $k\geq 1$ the function $_DF^{(k)}_{\frac{aj}{N},\frac{b}{N}}(\tau)$ is the analytic modular form of Def. 3.3.17.
\end{lemma}
\begin{proof}
Let
\[\Psi_k \in \Gamma \big(S^{an}, \omega_{E^{an}/S^{an}}^{\otimes k} \big)\]
be the image of $_DF^{(k)}_{\frac{a}{N},\frac{b}{N}}$ under $(3.8.39)$ and for $j\in (\Z/N\Z)^*$ let
\[p_j: \H \rightarrow S^{an}\]
be as in $(3.8.40)$.\\
With the explanations in Rem. 3.8.12, what we have to show are the following equalities in $\Gamma\big(\H,\omega_{\mathbb E/\H}^{\otimes k}\big)$:
\[\tag{$*$} \begin{split}
p_j^*(\Psi_k)&= -(2\pi i)^k \cdot {_D}F^{(k)}_{\frac{aj}{N},\frac{b}{N}}(\tau) \cdot (\mathrm{d}z)^{\otimes k}, \quad k\geq 1, \ k\neq 2,\\
p_j^*(\Psi_2)&= (2\pi i)^2 \cdot {_D}F^{(2)}_{\frac{aj}{N},\frac{b}{N}}(\tau) \cdot (\mathrm{d}z)^{\otimes 2}.
\end{split}
\]
The case $j=1$:
\[\tag{$**$} \begin{split} p_1^*(\Psi_k)&=-(2\pi i)^k\cdot _DF^{(k)}_{\frac{a}{N},\frac{b}{N}}(\tau)\cdot (\mathrm{d}z)^{\otimes k}, \quad k\geq 1, \ k\neq 2,\\
p_1^*(\Psi_2)&= (2\pi i)^2 \cdot {_D}F^{(2)}_{\frac{a}{N},\frac{b}{N}}(\tau) \cdot (\mathrm{d}z)^{\otimes 2}
\end{split}\]
is known to be true (cf. $(3.8.41)$ in Rem. 3.8.13). For arbitrary $j\in (\Z/N\Z)^*$ let
\[\sigma_j:=\begin{pmatrix} j & 0 \\ 0 & 1 \end{pmatrix} \in \mathrm{GL}_2(\Z/N\Z)\]
and use the same notation
\[\sigma_j: (\Z/N\Z)^* \times \Gamma(N) \backslash \H \rightarrow (\Z/N\Z)^* \times \Gamma(N) \backslash \H\]
for the automorphism of $S^{an}$ sending $(l,\tau)$ to $(jl,\tau)$.\\
The diagram $(3.8.40)$ for our fixed $j$ then equals the frame diagram in
\begin{equation*}  \tag{\textbf{3.8.42}} \begin{split}
\begin{xy}
\xymatrix@C-0.3cm{
\mathbb E=\Z^2\backslash(\C\times \H) \ar[r]\ar[d] & (\Z/N\Z)^*\times (\Z^2\times \Gamma(N)) \backslash (\C \times \H) \ar[d] \ar[r] & (\Z/N\Z)^*\times (\Z^2\times \Gamma(N)) \backslash (\C \times \H) \ar[d]\\
\H \ar[r]^{p_1 \qquad \qquad} & (\Z/N\Z)^* \times \Gamma(N) \backslash \H \ar[r]^{\sigma_j} & (\Z/N\Z)^* \times \Gamma(N) \backslash \H}
\end{xy}
\end{split}
\end{equation*}
where the left square is the cartesian diagram $(3.8.40)$ for $j=1$ and the upper right horizontal arrow is given by $(l,z,\tau)\mapsto (jl,z,\tau)$. On the other hand, from the elliptic curve with level $N$-structure $(E/S, je_1,e_2)$ and the universality of $(E/S,e_1,e_2)$ we obtain a cartesian diagram
\begin{equation*} \tag{\textbf{3.8.43}} \begin{split}
\begin{xy}
\xymatrix{
E\ar[r]\ar[d] &E \ar[d]\\
S \ar[r] & S}
\end{xy}
\end{split}
\end{equation*}
whose analytification (under the identification $(3.4.4)$) becomes the right cartesian square in $(3.8.42)$.\\
Let $k\geq 1$ be arbitrary.\\
By \cite{Ka}, Ch. I, Lemma 3.7 (1) (iii), we know that $_DF^{(k)}_{\frac{a}{N},\frac{b}{N}}$ maps to $_DF^{(k)}_{\frac{aj}{N},\frac{b}{N}}$ under the arrow
\[\Gamma\big(S,\omega_{E/S}^{\otimes k}\big) \rightarrow \Gamma\big(S,\omega_{E/S}^{\otimes k}\big)\]
induced by $(3.8.43)$. We thus conclude that the image of $\Psi_k$ under the map
\[\sigma_j^*:\Gamma\big(S^{an},\omega_{E^{an}/S^{an}}^{\otimes k}\big) \rightarrow \Gamma\big(S^{an},\omega_{E^{an}/S^{an}}^{\otimes k}\big)\]
coming from the right square in $(3.8.42)$ is equal to $\Phi_k$, where $\Phi_k$ is the image of $_DF^{(k)}_{\frac{aj}{N},\frac{b}{N}}$ under $(3.8.39)$. Hence, we have the following equality in $\Gamma(\H, \omega_{\mathbb E/\H}^{\otimes k})$:
\[p_j^*(\Psi_k)=(\sigma_j\circ p_1)^*(\Psi_k)=p_1^* (\sigma_j^*(\Psi_k))=p_1^*(\Phi_k).\]
But considering $(**)$ with $\Psi_k$ replaced by $\Phi_k$ and $a$ replaced by $aj$ we have
\[\begin{split}
p_1^*(\Phi_k)&=-(2\pi i)^k\cdot {_D}F^{(k)}_{\frac{aj}{N},\frac{b}{N}}(\tau) \cdot (\mathrm{d}z)^{\otimes k}, \quad k\geq 1, \ k\neq 2,\\
p_1^*(\Phi_2)&=(2\pi i)^2\cdot {_D}F^{(2)}_{\frac{aj}{N},\frac{b}{N}}(\tau) \cdot (\mathrm{d}z)^{\otimes 2}.
\end{split}\]
The last two formulas show $(*)$ and hence the claim of the lemma.
\end{proof}

\subsubsection{The algebraic specialization result}

We can finally prove our main result about the specialization of the $D$-variant of the polylogarithm.

\begin{theorem}
For $n \geq 0$ the cohomology class $\big(t_{a,b}^*(\varDn)\big)^{(n)}$ is equal to the image of
\[\begin{cases}
-{_D}F^{(2)}_{\frac{a}{N},\frac{b}{N}} \ \ &\textrm{if} \ \ n=0
\vspace{0.9mm}\\
\frac{(-1)^{n}}{n!} \cdot {_D}F^{(n+2)}_{\frac{a}{N},\frac{b}{N}} \ \ &\textrm{if} \ \ n\neq 0
\end{cases}
\]
under the map $(3.8.38)$:
\[\Gamma\big(S, \omega_{E/S}^{\otimes(n+2)}\big) \rightarrow H^1_{\mathrm{dR}}(S/\Q,\mathrm{Sym}^n_{\mathcal O_S}H^1_{\mathrm{dR}}(E/S)).\]
\end{theorem}
\begin{proof}
Consider the commutative diagram of Lemma 3.8.10:
\begin{equation*}
\begin{xy}
\xymatrix@C-0.3cm{
\Gamma\big(S, \omega_{E/S}^{\otimes(n+2)}\big)\ar[r]^{\kappa\quad \qquad \ } \ar[d] & H^1_{\mathrm{dR}}(S/\Q,\mathrm{Sym}^n_{\mathcal O_S}H^1_{\mathrm{dR}}(E/S)) \ar[d]^{\iota}\\
\Gamma\big(S^{an}, \omega_{E^{an}/S^{an}}^{\otimes(n+2)}\big) \ar[r] &  H^1_{\mathrm{dR}}(S^{an},\mathrm{Sym}^n_{\mathcal O_{S^{an}}}H^1_{\mathrm{dR}}(E^{an}/S^{an}))}
\end{xy}
\end{equation*}
We use it to deduce the statement of the theorem as follows:\\
By Lemma 3.8.14 the left vertical arrow of this diagram maps the section in the claim
\[\begin{cases}
-{_D}F^{(2)}_{\frac{a}{N},\frac{b}{N}} \qquad \ \ \ \in \Gamma\big(S, \omega_{E/S}^{\otimes(2)}\big) \ \ &\textrm{if} \ \ n=0
\vspace{0.9mm}\\
\frac{(-1)^{n}}{n!} \cdot {_D}F^{(n+2)}_{\frac{a}{N},\frac{b}{N}} \in \Gamma\big(S, \omega_{E/S}^{\otimes(n+2)}\big) \ \ & \textrm{if} \ \ n\neq 0
\end{cases}
\]
to the section of $\Gamma\big(S^{an}, \omega_{E^{an}/S^{an}}^{\otimes(n+2)}\big)$ which is given - in the sense of Rem. 3.8.12 - by the collection
\[\bigg(\frac{(-1)^{n+1}\cdot (2\pi i)^{n+2}}{n!}\cdot {_D}F^{(n+2)}_{\frac{aj}{N},\frac{b}{N}}(\tau)\bigg)_{j \in (\Z/N\Z)^*}, \ \ \textrm{where} \ \ n\geq0.\]
According to Thm. 3.8.8 the previous section in turn is sent for each $n\geq 0$ to
\[\Big((t_{a,b}^{an})^*\big((\varDn)^{an}\big)\Big)^{(n)}\]
by the lower horizontal arrow of the diagram.\\
The commutativity of the diagram and Lemma 3.8.7 hence imply that
\[\begin{cases}
\kappa \big(-{_D}F^{(2)}_{\frac{a}{N},\frac{b}{N}}\big)
\vspace{0.9mm}\\
\kappa\Big(\frac{(-1)^{n}}{n!} \cdot {_D}F^{(n+2)}_{\frac{a}{N},\frac{b}{N}}\Big), \ \ \textrm{where} \ \ n\neq 0,
\end{cases}
\]
and
\[\hspace{1.3cm} \begin{cases}
\big(t_{a,b}^*(\mathrm{pol}^0_{\mathrm{dR},D^2\cdot 1_{ \{ \epsilon \} }-1_{E[D]}})\big)^{(0)}
\vspace{0.9mm}\\
\big(t_{a,b}^*(\varDn)\big)^{(n)}, \ \ \textrm{where} \ \ n\neq0,
\end{cases}
\]

have the same image under the right vertical arrow $\iota$.\\
But $\iota$ is injective, as we have seen directly after its definition in $(3.8.31)$. This proves the theorem.
\end{proof}
\begin{remark}
In degree $n=0$ the $D$-variant specializes to a cohomology class in $H^1_{\mathrm{dR}}(S/\Q)$; the preceding theorem allows to express this class in terms of the Siegel units on $S$, as defined in \cite{Ka}, Ch. I, 1.4.\\
Namely, we obtain that $\big(t_{a,b}^*(\mathrm{pol}^0_{\mathrm{dR},D^2\cdot 1_{ \{ \epsilon \} }-1_{E[D]}})\big)^{(0)} \in H^1_{\mathrm{dR}}(S/\Q)$ is the image under $(3.8.38)$ of
\[-{_D}F^{(2)}_{\frac{a}{N},\frac{b}{N}}=F^{(2)}_{\frac{Da}{N},\frac{Db}{N}} - D^2F^{(2)}_{\frac{a}{N},\frac{b}{N}} \in \Gamma \big(S, \omega_{E/S}^{\otimes 2} \big).\]
Using the Kodaira-Spencer isomorphism $(3.8.35)$, the previous section identifies (according to \cite{Ka}, Ch. I, Prop. 3.11 (2)) with
\[\tag{\textbf{3.8.44}} \mathrm{dlog} \bigg (\Big(g_{\frac{a}{N},\frac{b}{N}}\Big)^{D^2} \cdot \Big(g_{\frac{Da}{N},\frac{Db}{N}}\Big)^{-1} \bigg) \in \Gamma(S, \Omega^1_{{S/\Q}}),\]
where $g_{\frac{a}{N},\frac{b}{N}}$ resp. $g_{\frac{Da}{N},\frac{Db}{N}}$ are the mentioned Siegel unit elements in $\Gamma(S,\mathcal O_S)^*\otimes_{\Z}\Q$; the desired cohomology class then is the image of $(3.8.44)$ under the canonical map $\Gamma(S,\Omega^1_{S/\Q})\rightarrow H^1_{\mathrm{dR}}(S/\Q)$.\\
If we additionally assume $(D,6)=1$, then the last formula in ibid.,  Ch. I, 1.4, implies that the element of $(3.8.44)$ is equal to
\[\mathrm{dlog}(t_{a,b}^*(_D\theta_{E/S})) \in \Gamma(S,\Omega^1_{S/\Q}),\]
where $_D\theta_{E/S}\in \Gamma(U_D,\mathcal O_E)^*$ denotes the Kato-Siegel function associated with the elliptic curve $E/S$ and the integer $D$ (cf. ibid., Ch. I, Prop. 1.3 and 1.10, for its definition and construction).
\end{remark}
\subsubsection{Application to the de Rham Eisenstein classes}
The subsequent definition of the de Rham Eisenstein classes is modeled after \cite{Ki3}, 4.2, where the $\ell$-adic Eisenstein classes are introduced.\\
\newline
For each $n\geq 0$ we have a horizontal morphism
\[\mathcal H^\vee\otimes_{\mathcal O_S}\prod_{k=0}^{n+1}\mathrm{Sym}^k_{\mathcal O_S}\mathcal H \rightarrow \prod_{k=0}^n\mathrm{Sym}^k_{\mathcal O_S}\mathcal H\]
which is defined to be zero on $\mathcal H^\vee$ and for $1\leq k \leq n+1$ to be given on $\mathcal H^\vee\otimes_{\mathcal O_S}\mathrm{Sym}^k_{\mathcal O_S}\mathcal H $ by
\[\mathcal H^\vee\otimes_{\mathcal O_S}\mathrm{Sym}^k_{\mathcal O_S}\mathcal H \rightarrow \mathrm{Sym}^{k-1}_{\mathcal O_S}\mathcal H, \quad h^\vee\otimes h_1\cdot ... \cdot h_k \mapsto \frac{1}{k+1} \sum_{j=1}^kh^\vee(h_j)h_1\cdot ... \cdot \widehat{h_j}\cdot ... \cdot h_k.\]
Write
\[\mathrm{contr}_n: H^1_{\mathrm{dR}}\Big(S/\Q,\mathcal H^\vee\otimes_{\mathcal O_S}\prod_{k=0}^{n+1}\mathrm{Sym}^k_{\mathcal O_S}\mathcal H\Big) \rightarrow H^1_{\mathrm{dR}}\Big(S/\Q,\prod_{k=0}^n\mathrm{Sym}^k_{\mathcal O_S}\mathcal H\Big)\]
for the induced homomorphism (the $n$-th "contraction map").
\begin{definition}
For each $n\geq 0$ we set
\[\mathrm{Eis}^n(t_{a,b}):=-N^{n-1}\cdot\big(\mathrm{contr}_n(t_{a,b}^*\mathrm{pol}^{n+1}_{\mathrm{dR}})\big)^{(n)} \in H^1_{\mathrm{dR}}(S/\Q,\mathrm{Sym}^n_{\mathcal O_S}\mathcal H)\]
and call this element the \underline{$n$-th de Rham Eisenstein class evaluated at $t_{a,b}$}.
\end{definition}
As explained in \cite{Ki3}, Rem. 4.2.3, the factor $-N^{n-1}$ is owed to tradition.\\
\newline
Moreover, analogously as in \cite{Ki3}, Prop. 4.4.1, we have the following formula relating the specialization of the $D$-variant of the polylogarithm with the de Rham Eisenstein classes:
\begin{lemma}
For each $n\geq 0$ the following equality holds in $H^1_{\mathrm{dR}}(S/\Q,\mathrm{Sym}^n_{\mathcal O_S}\mathcal H)$:
\[\big(t_{a,b}^*(\varDn)\big)^{(n)}=-N^{1-n}\big(D^2\mathrm{Eis}^n(t_{a,b})-D^{-n}\mathrm{Eis}^n(t_{Da,Db})\big).\]
(For the left side we don't use Poincaré duality as we usually did before).\\
In particular, if $D\equiv 1 \ \mathrm{mod} \ N$ we obtain:
\[\big(t_{a,b}^*(\varDn)\big)^{(n)}=-N^{1-n}\frac{D^{n+2}-1}{D^n}\mathrm{Eis}^n(t_{a,b}).\]
\qed
\end{lemma}

Now, combining the preceding lemma with Thm. 3.8.15 implies that for $D\equiv 1 \ \mathrm{mod} \ N$ we have
\[\mathrm{Eis}^0(t_{a,b})=-N^{-1}\frac{1}{D^{2}-1} \cdot \big(-{_D}F^{(2)}_{\frac{a}{N},\frac{b}{N}}\big)\]
resp. for $n> 0$:
\[\mathrm{Eis}^n(t_{a,b})=-N^{n-1}\frac{D^n}{D^{n+2}-1} \cdot \frac{(-1)^{n}}{n!}\cdot {_D}F^{(n+2)}_{\frac{a}{N},\frac{b}{N}},\]
where here ${_D}F^{(2)}_{\frac{a}{N},\frac{b}{N}}$ resp. ${_D}F^{(n+2)}_{\frac{a}{N},\frac{b}{N}}$ means the element of $H^1_{\mathrm{dR}}(S/\Q)$ resp. $H^1_{\mathrm{dR}}(S/\Q,\mathrm{Sym}^n_{\mathcal O_S}\mathcal H)$ induced by the algebraic modular form ${_D}F^{(2)}_{\frac{a}{N},\frac{b}{N}}$ resp. ${_D}F^{(n+2)}_{\frac{a}{N},\frac{b}{N}}$ via the map $(3.8.38)$ resp. via the map $(3.8.38)$ together with the Poincaré duality identification $(3.8.20)$.\\
Taking furthermore into account that for $D\equiv 1 \ \mathrm{mod} \ N$:
\[{_D}F^{(n+2)}_{\frac{a}{N},\frac{b}{N}}=D^2F^{(n+2)}_{\frac{a}{N},\frac{b}{N}} - D^{-n}F^{(n+2)}_{\frac{Da}{N},\frac{Db}{N}}=\frac{D^{n+2}-1}{D^n}\cdot F^{(n+2)}_{\frac{a}{N},\frac{b}{N}}, \quad n\geq 0,\]
we obtain the subsequent explicit description for the de Rham Eisenstein classes; for a different approach to this problem (by computing residues at the cusps) see \cite{Ba-Ki2}, 3.

\begin{corollary}
We have the equalities
\[\begin{cases}
\mathrm{Eis}^0(t_{a,b})=-N^{-1}\cdot\big(-F^{(2)}_{\frac{a}{N},\frac{b}{N}}\big) \quad &\textrm{in} \ \ H^1_{\mathrm{dR}}(S/\Q)
\vspace{0.9mm}\\
\mathrm{Eis}^n(t_{a,b})=-N^{n-1} \frac{(-1)^{n}}{n!}\cdot F^{(n+2)}_{\frac{a}{N},\frac{b}{N}} \quad &\textrm{in} \ \ H^1_{\mathrm{dR}}(S/\Q,\mathrm{Sym}^n_{\mathcal O_S}\mathcal H), \ n>0, 
\end{cases}\]
where here $F^{(2)}_{\frac{a}{N},\frac{b}{N}}$ resp. $F^{(n+2)}_{\frac{a}{N},\frac{b}{N}}$ means the element of $H^1_{\mathrm{dR}}(S/\Q)$ resp. $H^1_{\mathrm{dR}}(S/\Q,\mathrm{Sym}^n_{\mathcal O_S}\mathcal H)$ induced by the algebraic modular form $F^{(2)}_{\frac{a}{N},\frac{b}{N}}$ resp. $F^{(n+2)}_{\frac{a}{N},\frac{b}{N}}$ via $(3.8.38)$ resp. via $(3.8.38)$ together with the identification $(3.8.20)$.\qquad \qed
\end{corollary}
\begin{remark}
(i) Observing again \cite{Ka}, Ch. I, Prop. 3.11 (2), the first equation of Cor. 3.8.19 translates into
\[\mathrm{Eis}^0(t_{a,b})=-N^{-1}\mathrm{dlog}\Big(g_{\frac{a}{N},\frac{b}{N}}\Big)\]
with the Siegel unit $g_{\frac{a}{N},\frac{b}{N}} \in \Gamma(S,\mathcal O_S)^*\otimes_{\Z}\Q$ of ibid., Ch. I, 1.4.\vspace{1mm}\\
(ii) Working with the identification
\[H^1_{\mathrm{dR}}(E/S) \xrightarrow{\sim} \mathcal H, \quad x \mapsto \{\ y\mapsto \mathrm{tr}(y\cup x)\}\]
instead of $(3.8.20)$:
\[H^1_{\mathrm{dR}}(E/S) \xrightarrow{\sim} \mathcal H, \quad x \mapsto \{\ y\mapsto \mathrm{tr}(x\cup y)\},\]
the second equation of Cor. 3.8.19 changes into
\[\mathrm{Eis}^n(t_{a,b})=-N^{n-1} \frac{1}{n!}\cdot F^{(n+2)}_{\frac{a}{N},\frac{b}{N}} \quad  \textrm{in} \ \ H^1_{\mathrm{dR}}(S/\Q,\mathrm{Sym}^n_{\mathcal O_S}\mathcal H), \ n>0.\]
\end{remark}

\end{document}